\documentclass[a4paper,11pt]{book}

\usepackage{amsfonts,qtree,stmaryrd,textcomp}
\usepackage{pifont,amsmath,amsthm,tabularx,graphicx}
\usepackage{tree-dvips,pictexwd,dcpic}
\usepackage{bbm}
\usepackage{bm}

\usepackage{stmaryrd}

\usepackage[T1]{fontenc}
\usepackage[latin1]{inputenc}
\usepackage[text={17cm,26cm},centering]{geometry}

\usepackage{hyperref}

\usepackage{etex}
\usepackage{fancyhdr}
\usepackage{graphicx}
\usepackage{enumitem}
\usepackage{amsmath}
\usepackage[bbgreekl]{mathbbol}
\usepackage{amssymb}
\usepackage{amsthm}
\usepackage[all]{xy}

\usepackage{epigraph}
\RequirePackage{bussproofs}

\usepackage{etex}
\usepackage{a4wide}
\usepackage{fancyhdr}
\usepackage{graphicx}

\usepackage{euscript}

\usepackage{mathabx}

\usepackage{scalerel,stackengine}

\usepackage{framed}
\usepackage{mathtools}
\usepackage{url}
\usepackage{proof}
\usepackage{xspace}
\usepackage{listings}
\usepackage{wrapfig}
\usepackage{tikz}
\usepackage{tikz-cd}

\usepackage{relsize}

\usetikzlibrary{positioning}
\usetikzlibrary{shapes}

\usetikzlibrary{arrows}
\usetikzlibrary{calc}
\usepackage{tikz-3dplot}
\usetikzlibrary{decorations.pathmorphing}
\usetikzlibrary{patterns}

\usetikzlibrary{decorations.pathreplacing}
\tikzset{axis/.style={&lt;-&gt;}}

\usepackage{amsfonts,amssymb,amsmath,qtree,stmaryrd,textcomp}
\usepackage{tree-dvips,pictexwd,dcpic,graphicx,latexsym}

\newtheorem{theorem}{Theorem}[section]
\newtheorem{proposition}[theorem]{Proposition}
\newtheorem{lemma}[theorem]{Lemma}
\newtheorem{corollary}[theorem]{Corollary}
\newtheorem{remark}[theorem]{Remark}

\newtheorem{definition}[theorem]{Definition}
\newtheorem{note}[theorem]{Note}

\usepackage[latin1]{inputenc}
\usepackage{appendix}
\usepackage{enumitem}
\usepackage{amsmath}
\usepackage[bbgreekl]{mathbbol}
\usepackage{amssymb}
\usepackage{amsthm}
\usepackage[all]{xy}

\usepackage{epigraph}

\usepackage{tipa}

\RequirePackage{bussproofs}

\newcommand{\B}{\boldsymbol}
\newcommand{\C}[1]{\mathcal{#1}}
\newcommand{\D}[1]{\mathbb{#1}}

\newcommand{\Nat}{{\mathbb N}}
\newcommand{\Real}{{\mathbb R}}

\newcommand{\BR}{{\Bic(\Real)}}

\newcommand{\Ap}{\mathrm{Ap}}

\newcommand{\Bic}{\mathrm{Bic}}

\newcommand{\BISH}{\mathrm{BISH}}

\newcommand{\CST}{\mathrm{CST}}

\newcommand{\Const}{\mathrm{Const}}

\newcommand{\dom}{\mathrm{dom}}

\newcommand{\ev}{\mathrm{ev}}

\newcommand{\id}{\mathrm{id}}

\newcommand{\Ind}{\mathrm{Ind}}

\newcommand{\PEM}{\mathrm{PEM}}

\newcommand{\Lim}{\mathrm{Lim}}

\newcommand{\Mor}{\mathrm{Mor}}

\newcommand{\Comp}{\mathrm{Comp}}

\newcommand{\LPO}{\mathrm{LPO}}

\newcommand{\ZF}{\mathrm{ZF}}

\newcommand{\HoTT}{\mathrm{HoTT}}

\newcommand{\Cat}{\mathrm{Cat}}

\newcommand{\Id}{\mathrm{Id}}

\newcommand{\BS}{\mathrm{BS}}

\newcommand{\Bis}{\mathrm{\mathbf{Bis}}}

\newcommand{\dotminus}{\mathbin{{{-\mkern -9.5mu%
\mathchoice{\raise 2pt\hbox{$\cdot$}\mkern6mu}%
{\raise 2pt\hbox{$\cdot$}\mkern6mu}%
{\raise 1.5pt\hbox{$\scriptstyle\cdot$}\mkern4mu}%
{\raise 1pt\hbox{$\scriptscriptstyle\cdot$}\mkern4mu}}}}}

\newcommand{\Fii}{{\mathbb F}}

\newcommand{\fXY}{f \colon X \to Y}

\newcommand{\pr}{\textnormal{\texttt{pr}}}

\newcommand{\idtoeqv}{\textnormal{\texttt{idtoEqv}}}

\newcommand{\refl}{\textnormal{\texttt{refl}}}

\newcommand{\UA}{\textnormal{\texttt{UA}}}

\newcommand{\MLTT}{\mathrm{MLTT}} 
\newcommand{\CZF}{\mathrm{CZF}}

\newcommand{\CSFT}{\mathrm{CSFT}}

\newcommand{\BCM}{\mathrm{BCM}}

\newcommand{\CM}{\mathrm{CM}}

\newcommand{\TOT}{\Leftrightarrow}

\newcommand{\To}{\Rightarrow}
\newcommand{\ot}{\leftarrow}

\newcommand{\oT}{\Leftarrow}

\newcommand{\sto}{\rightsquigarrow}

\newcommand{\eto}{\hookrightarrow}

\newcommand{\compl}{\textnormal{\texttt{compl}}}

\newcommand{\pto}{\rightharpoonup}

\newcommand{\Compl}{\textnormal{\texttt{Compl}}}

\newcommand{\Fam}{\textnormal{\texttt{Fam}}}

\newcommand{\LANC}{\mathrm{LANC}}

\newcommand{\ac}{\textnormal{\texttt{ac}}}

\newcommand{\prb}{\textnormal{\textbf{pr}}}

\newcommand{\BST}{\mathrm{BST}}

\newcommand{\LCST}{\mathrm{LCST}}

\newcommand{\CMT}{\mathrm{CMT}}

\newcommand{\cof}{\textnormal{\texttt{cof}}}

\newcommand{\mcup}{\mathsmaller{\cup}}

\newcommand{\oi}{\D I_{\mathsmaller{01}}}

\newcommand{\oic}{{\D I^{\mcup}_{\mathsmaller{01}}}}

\newcommand{\eqt}{\sim_{\mathsmaller{X}}^\mathsmaller{\tau}}

\newcommand{\Disj}{\B \rrbracket \B \llbracket}

\newcommand{\Hom}{\mathrm{Hom}}

\newcommand{\fib}{\textnormal{\texttt{fib}}}

\newcommand{\sw}{\textnormal{\texttt{sw}}}

\newcommand{\Map}{\textnormal{\texttt{Map}}}

\newcommand{\lt}{\preccurlyeq}
\newcommand{\mt}{\succcurlyeq}

\newcommand{\ltI}{D^{\lt}(I)}

\newcommand{\BCMT}{\mathrm{BCMT}}

\newcommand{\Prf}{\textnormal{\texttt{Prf}}}

\newcommand{\Eqq}{\textnormal{\texttt{PrfEql}}}
\newcommand{\Eq}{\Eqq_0}

\newcommand{\EwE}{\mathrm{EwE}}

\newcommand{\BHK}{\mathrm{BHK}}
\newcommand{\Set}{\textnormal{\texttt{Set}}}

\newcommand{\FunExt}{\mathrm{FunExt}}

\newcommand{\centr}{\textnormal{\texttt{centre}}}

\newcommand{\cofib}{\textnormal{\texttt{cofib}}}

\newcommand{\pcirc}{ \ \mathsmaller{\odot} \ }

\newcommand{\Eql}{\textnormal{\texttt{Eql}}}
\newcommand{\eql}{\textnormal{\texttt{eql}}}

\newcommand{\Dm}{\textnormal{\texttt{Dom}}}

\newcommand{\parl}{\textnormal{\texttt{par}}}

\newcommand{\Spec}{\textnormal{\texttt{Spec}}}

\newcommand{\Cnt}{\textnormal{\texttt{Cont}}}

\newcommand{\Cf}{\textnormal{\texttt{Cof}}}

\newcommand{\BMS}{\mathrm{BMS}}

\newcommand{\Borel}{\textnormal{\texttt{Borel}}}

\newcommand{\Baire}{\textnormal{\texttt{Baire}}}

\newcommand{\BMT}{\mathrm{BMT}}

\newcommand{\MS}{\mathrm{MS}}

\newcommand{\PMS}{\mathrm{PMS}}

\newcommand{\BCMS}{\mathrm{BCMS}}
\newcommand{\BCCMS}{\mathrm{BCCMS}}

\newcommand{\BCIS}{\mathrm{BCIS}}

\newcommand{\im}{\mathrm{im}}

\newcommand{\Even}{\textnormal{\texttt{Even}}}

\newcommand{\Odd}{\textnormal{\texttt{Odd}}}

\newcommand{\se}{\mathrm{se}}

\newcommand{\IS}{\mathrm{IS}}

\newcommand{\PIS}{\mathrm{PIS}}

\newcommand{\Morr}{\mathrm{mor}}

\newcommand{\Funct}{\mathrm{Funct}}

\newcommand{\mon}{\mathrm{mn}}

\newcommand{\Ob}{\mathrm{Ob}}

\usepackage{makeidx}
\makeindex

\begin{document}

\date{}

\title{\textbf{Families of Sets in Bishop Set Theory}}

\author{Iosif Petrakis\\	
Mathematics Institute, Ludwig-Maximilians-Universit\"{a}t M\"{u}nchen\\
petrakis@math.lmu.de}  

%





\maketitle

\tableofcontents

 \cleardoublepage

\chapter*{Abstract}

\markboth{Abstract}{Abstract}
%
%
%
%
%
%
%

We develop the theory of set-indexed families of sets within the informal Bishop Set Theory $(\BST)$, a
reconstruction of Bishop's theory of sets,. The latter is the informal theory of sets and functions underlying
Bishop-style constructive mathematics $(\BISH)$ and it is developed in Chapter 3 of Bishop's seminal 
book \textit{Foundations of Constructive Analysis}~\cite{Bi67} and in Chapter 
3 of \textit{Constructive Analysis}~\cite{BB85} that Bishop co-authored with Bridges. 

In the Introduction we briefly present the relation of Bishop's set theory to the set-theoretic
and type-theoretic foundations of mathematics, and we describe the features of $\BST$ that ``complete'' 
Bishop's theory of sets. These are the explicit use of the class ``universe of sets'', a clear distinction
between sets and classes, the explicit use of dependent operations, and the concrete formulation of 
various notions of families of sets. 

In Chapter~\ref{chapter: BST} we present the fundamentals of Bishop's theory of sets, extended with the features which form $\BST$. The universe $\D V_0$ of sets is implicit in Bishop's work, while the notion of a dependent operation
over a non-dependent assignment routine from a set to $\D V_0$ is explicitly mentioned, although in a rough way. 
These concepts are necessary to a concrete definition of a set-indexed family of sets, 
the main object of our study, which is only mentioned by Bishop.

In Chapter~\ref{chapter: familiesofsets} we develop the basic theory of set-indexed families of sets and of 
family-maps between them. We study the exterior union of a family of sets $\Lambda$, or the $\sum$-set of 
$\Lambda$, and the set
of dependent functions over $\Lambda$, or the $\prod$-set of $\Lambda$. We prove the distributivity of
$\prod$ over $\sum$ for families of sets indexed by a product of sets, which is the translation 
of the type-theoretic axiom of choice into $\BST$. Sets of sets are special set-indexed families of sets 
that allow ``lifting'' of functions on the index-set to functions on them. The direct families of sets and the 
set-relevant families of sets are introduced. The index-set of the former is a directed set, while the 
transport maps of the latter are more than one and appropriately indexed. With the use of the introduced universe
$\D V_0^{\im}$ of sets and impredicative sets we study families of families of sets, the next rung of the
ladder of set-like objects in $\D V_0^{\im}$.

In Chapter~\ref{chapter: familiesofsubsets} we develop the basic theory of set-indexed families of subsets
and of the corresponding family-maps between them. In contrast to set-indexed families of sets, the properties
of which are determined ``externally'' through their transport maps, the properties of a set-indexed family
$\Lambda(X)$ of subsets of a given set $X$ are determined ``internally'' through the embeddings of the
subsets of $\Lambda(X)$ to $X$. The interior union of $\Lambda(X)$ is the internal analogue to the $\sum$-set
of a set-indexed family of sets $\Lambda$, 
and the intersection of $\Lambda(X)$ is the internal analogue to the $\prod$-set of $\Lambda$.
Families of sets over products, sets of subsets, and direct families of subsets are the internal analogue to the 
corresponding notions for families of sets. Set-indexed families of partial functions
and set-indexed families of complemented subsets, together with their corresponding family-maps, are studied.

%
%

In Chapter~\ref{chapter: bspaces} we connect various notions and results from the theory of families of sets 
and subsets to the theory of Bishop spaces, a function-theoretic approach to constructive topology. Associating 
in an appropriate way to each set $\lambda_0(i)$ of an $I$-family of sets $\Lambda$ a Bishop topology $F_i$, a
spectrum $S(\Lambda)$ of Bishop spaces is generated. The $\sum$-set and the $\prod$-set of a spectrum $S(\Lambda)$ 
are equipped with canonical Bishop topologies. A direct spectrum of Bishop spaces is a family of Bishop spaces 
associated to a direct family of sets. The direct and inverse limits of direct spectra of Bishop spaces are studied.
Direct spectra of Bishop subspaces are also examined. Many Bishop topologies used in this chapter are defined inductively within the extension $\BISH^*$ of $\BISH$ with inductive definitions with rules of countably many premises.

In Chapter~\ref{chapter: measure} we study the Borel and Baire sets within Bishop spaces as a constructive 
counterpart to the study of Borel and Baire algebras within topological spaces. As we use the inductively 
defined least Bishop topology, and as the Borel and Baire sets over a family of $F$-complemented subsets are
defined inductively, we work again within $\BISH^*$. In contrast to the classical theory, we show that 
the Borel and the Baire sets of a Bishop space coincide. Finally, our reformulation within $\BST$ of the Bishop-Cheng 
definition of a measure space and of an integration space, based on the notions of families of complemented 
subsets and of families of partial functions, facilitates a predicative reconstruction of the originally 
impredicative Bishop-Cheng measure theory.

\chapter{Introduction}
\label{chapter: intro}

Bishop's theory of sets is Bishop's account of the informal theory of sets and functions that underlies 
Bishop-style constructive mathematics $\BISH$. We briefly present the relation of this theory to
the set-theoretic and type-theoretic foundations of mathematics. Bishop Set Theory $(\BST)$ is 
our ``completion'' of Bishop's theory of sets 
with a universe of sets, with a clear distinction between sets and classes, with an explicit use 
of dependent operations, and with a concrete formulation  of various notions of families of sets. 
We explain how the theory of families of sets within $\BST$ that is elaborated in this work is used,
in order to reveal proof-relevance in $\BISH$, to develop the theory of spectra of Bishop spaces, and 
to reformulate predicatively the fundamental notions of the impredicative Bishop-Cheng measure theory.

\section{Bishop's theory of sets}
\label{sec: bst}

The theory of sets underlying Bishop-style constructive mathematics $(\BISH)$ was only sketched in Chapter 3 of
Bishop's seminal book~\cite{Bi67}. Since Bishop's central aim in~\cite{Bi67} was to show 
that a large part of advanced mathematics can be done within a constructive and 
computational framework that does not contradict the classical practice, the inclusion of a 
detailed account of the set-theoretic foundations of $\BISH$ could possibly be against the effective 
delivery of his message.

The Bishop-Cheng measure theory, developed in~\cite{BC72}, was very
different from the measure theory of~\cite{Bi67}, and the inclusion of an enriched version of the former 
into~\cite{BB85}, the book on constructive analysis that Bishop co-authored with Bridges later, 
affected the corresponding Chapter 3 in two main respects. First, the inductively defined notion of the set 
of Borel sets generated by a given family of complemented subsets of a set $X$, with respect to a set 
of real-valued functions on $X$, was excluded, as unnecessary, and, second, the operations on the 
complemented subsets of a set $X$ were defined differently, and in accordance to the needs 
of the new measure theory. 
 
Yet, in both books many issues were left untouched, a fact that often was a source of confusion.
In many occasions, especially
in the measure theory of~\cite{BC72} and~\cite{BB85}, the powerset was treated as a set, while in the
measure theory of~\cite{Bi67}, Bishop generally avoided the powerset by using appropriate families of 
subsets instead. In later works of Bridges and Richman, like ~\cite{BR87} and~\cite{MRR88}, the powerset was
clearly used as a set, in contrast though, to the predicative spirit of~\cite{Bi67}. 

The concept of a family of
sets indexed by a (discrete) set, was asked to be defined in~\cite{Bi67} (Exercise 2, p.~72), and a definition, 
attributed to Richman, was given in~\cite{BB85} (Exercise 2, p.~78). An elaborate study though, of this concept 
within $\BISH$ is missing, despite its central character in the measure theory of~\cite{Bi67}, its 
extensive use in the theory of Bishop spaces~\cite{Pe15} and in abstract constructive algebra~\cite{MRR88}.
Actually, in~\cite{MRR88} Richman introduced the more general notion of a family of objects of 
a category indexed by some set, but the categorical component in the resulting mixture of Bishop's set theory and 
category theory was not explained in constructive terms\footnote{This was done e.g., in the
the formulation of category theory in homotopy type theory (Chapter 9 in~\cite{HoTT13}).}.

Contrary to the standard view on Bishop's relation to formalisation, Bishop was very interested in it. 
In~\cite{Bi70}, p.~60, he writes:
\begin{quote}
Another important foundational problem is to find a formal system that will efficiently
express existing predictive mathematics. I think we should keep the formalism as primitive as possible,
starting with a minimal system and enlarging it only if the enlargement serves a genuine mathematical
need. In this way the formalism and the mathematics will hopefully interact to the advantage
of both.
\end{quote}
Actually, in~\cite{Bi70} Bishop proposed $\Sigma$, a variant
of G\"odel's $T$, as a formal system for $\BISH$. In the last two pages of~\cite{Bi70} he sketched very 
briefly how $\Sigma$ can be presented as a functional programming language, like fortran and algol. 
In p.~72 he also added:
\begin{quote}
It would be interesting to take $\Sigma$ as the point of departure for a reasonable programming language,
and to write a compiler.
\end{quote}
Bishop's views on a full-scale program on the foundations of mathematics are realised in a more developed 
form in his, unfortunately, unpublished papers~\cite{Bi68a} and~\cite{Bi68b}. In the first, Bishop elaborated
a version of dependent type theory with one universe, in order to formalise $\BISH$. This was the first time that  some form of type theory is used to formalise constructive mathematics. 

As Martin-L\"of explains in~\cite{ML68}, p.~13, he got access to Bishop's book only shortly after his own book 
on constructive mathematics~\cite{ML68} was finished. Bishop's book~\cite{Bi67} also motivated his version of type theory.
Martin-L\"of opened his first published paper on type theory (\cite{ML75}, p.~73) as follows. 
\begin{quote}
The theory of types with which we shall be concerned is intended to be a full scale system for formalizing
intuitionistic mathematics as developed, for example, in the book of Bishop.
\end{quote}

The type-theoretic interpretation of Bishop's set theory into the theory of setoids (see especially the 
work of Palmgren~\cite{Pa05}-\cite{PW14}) has become nowadays the standard way to understand \textit{Bishop sets}
(as far as I know, this is a term due to Palmgren). A setoid is a type $A$ in a fixed universe $\C U$ equipped
with a term $\simeq \colon A \to A \to \C U$ that satisfies the properties of an equivalence relation. The identity
type of Martin-L\"of's intensional type theory ($\MLTT$) (see~\cite{ML98}), expresses, in a proof-relevant way,
the existence of the least reflexive relation on a type, a fact with no counterpart in Bishop's set theory. As
a consequence, the free setoid on a type is definable (see~\cite{Pa14}, p.~90), and the presentation axiom in
setoids is provable (see Note~\ref{not: onpresentationaxiom}). Moreover, in $\MLTT$ the families of types over a type $I$ is the type $I \to \C U$, 
which belongs to the successor universe $\C U{'}$ of $\C U$. In Bishop's set theory though, where only one
universe of sets is implicitly used, the set-character of the totality of all families of sets indexed by some set $I$
is questionable from the predicative point of view (see our comment after the Definition~\ref{def: map}). 

The quest $\B Q$\index{$\B Q$} of finding a formal system suitable for Bishop's system of informal constructive
mathematics $\BISH$ dominated the foundational studies of the 1970's. Myhill's system $\CST$,
introduced in~\cite{My75}, and later Aczel's $\CZF$ (see~\cite{AR10}),
Friedman's system $B$, developed in~\cite{Fr77}, and Feferman's system of explicit mathematics $T_0$ 
(see~\cite{Fe75} and~\cite{Fe79}), are some of the systems related to $\B Q$,  
but soon developed independently from it. These systems 
were influenced a lot from the classical Zermelo-Fraenkel set theory, and could be described as ``top-down'' 
approaches to the goal of $\B Q$, as they have many ``unexpected'' features with respect to 
$\BISH$. 
Using Feferman's terminology from~\cite{Fe79}, these formal systems are not completely
faithful\index{faithful formalisation} to $\BISH$. If $T$ is a formal theory of an informal body of mathematics $M$, 
Feferman gave in~\cite{Fe79} the following definitions.\\[1mm]
(i) $T$ is \textit{adequate} for $M$, if every concept, argument,\index{adequate formalisation}
and result of $M$ is represented by a (basic or defined) concept, proof, and a theorem, respectively, of $T$.\\[1mm]
(ii) $T$ is \textit{faithful} to $M$, if every basic concept of $T$ corresponds to a basic concept of $M$ and every 
axiom 
and rule of $T$ corresponds to or is implicit in the assumptions and reasoning followed in $M$ (i.e., $T$
does not go beyond $M$ conceptually or in principle).\\[1mm]
In~\cite{Be81}, p.~153, Beeson called $T$ \textit{suitable} to $M$, if $T$ is adequate for $M$ and faithful to $M$.\\[1mm]
Beeson's systems $S$ and $S_0$ in~\cite{Be81}, and Greenleaf's system
of liberal constructive set theory $\LCST$ in~\cite{Gr81} were dedicated to $\B Q$. Especially Beeson tried to 
find a faithful and adequate formalisation of $\BISH$, and, by including a serious amount of proof relevance, 
his systems stand in between the set-theoretic, proof-irrelevant point of view and the type-theoretic, 
proof-relevant point of view.

All aforementioned systems though, were not really ``tested'' with respect to $\BISH$. Only very small parts of $\BISH$
were actually implemented in them, and their adequacy for $\BISH$ was mainly a claim, rather than 
a shown fact. The implementation of Bishop's constructivism within a formal system for it was taken seriously 
in the type-theoretic formalisations of $\BISH$, and especially in the work of Coquand (see e.g.,~\cite{CP98} 
and~\cite{CS07}), Palmgren (see e.g.,~\cite{IP06} and the collaborative  work~\cite{CDPS05}), the
Nuprl research group of Constable (see e.g.,~\cite{Co86}), and of Sambin and Maietti
within the Minimalist Foundation (see~\cite{Sa19} and~\cite{Ma09}).

\section{Bishop Set Theory $(\BST)$ and Bishop's theory of sets}
\label{sec: BST}

Bishop set theory $(\BST)$ is an informal, constructive theory of totalities and assignment routines
that serves as a ``completion'' of Bishop's theory of sets. Its first aim is to fill in the ``gaps'', or 
highlight the fundamental notions that were suppressed by Bishop in his account of the set theory underlying $\BISH$. 
Its second aim is 
to serve as an intermediate step between Bishop's theory of sets and a suitable, in Beeson's sense, formalisation of 
$\BISH$. To assure faithfulness, we use concepts or principles that appear, explicitly or implicitly, in $\BISH$.
Next we describe briefly the features of $\BST$ that ``complete'' Bishop's theory of sets.\\[2mm]
\textbf{1. Explicit use of a universe of sets.} Bishop used a universe of sets only implicitly. E.g., he ``roughly'' 
describes in~\cite{Bi67}, p.~72, a set-indexed family of sets as 
\begin{quote}
$\ldots$ a rule which assigns to each $t$ in a discrete set $T$ a set $\lambda(t)$.
\end{quote}
Every other rule, or assignment routine mentioned by Bishop is from one given totality, the domain of the rule, 
to some other totality, its codomain. The only way to make the rule of a family of sets compatible with this 
pattern is to employ a totality of sets.
In~\cite{Bi68a} Bishop explicitly used a universe in his type theory.
Here we use the totality $\D V_0$ of sets, which is defined in an open-ended way,  and it contains the primitive 
set $\Nat$ and all defined sets. $\D V_0$ itself is not a set, but a class. It is a notion instrumental to the definition
of dependent operations, and of a set-indexed family of sets.\\[2mm]
\textbf{2. Clear distinction between sets and classes.} A class is a totality defined through a membership condition 
in which a quantification over $\D V_0$ occurs. The powerset $\C P(X)$ of a set $X$, the totality $\C P^{\Disj}(X)$ 
of complemented subsets of a set $X$, and the totality $\C F(X,Y)$ of partial functions from a set $X$ to a set $Y$ 
are characteristic examples of classes. A class is never used here as the domain of an assignment routine, only as
a codomain of an assignment routine. \\[2mm]
\textbf{3. Explicit use of dependent operations.} The standard view, even among practicioners of Bishop-style 
constructive mathematicians, is that dependency is not necessary to $\BISH$.
Dependent functions though, do appear explicitly in Bishop's definition of the intersection $\bigcap_{t \in T}
\lambda (t)$ of a family $\lambda$ of subsets of some set $X$ indexed by an inhabited set $T$
(see~\cite{Bi67}, p.~65, and ~\cite{BB85}, p.~70). We show that the elaboration of dependency within $\BISH$ is only fruitful 
to it. Dependent functions are not only necessary to the definition of products of families 
of sets indexed by an arbitrary set, but as we show throughout this work in many areas of constructive mathematics. 
Some form of dependency is also formulated in Bishop's type theory~\cite{Bi68a}. The somewhat ``silent'' role of 
dependency within Bishop's set theory is replaced by a central role within $\BST$.\\[2mm]
\textbf{4. Elaboration of the theory of families of sets.} With the use of the universe $\D V_0$, of the notion of 
a non-dependent assignment routine $\lambda_0$ from an index-set $I$ to $\D V_0$, and of a certain dependent operation $\lambda_1$, we 
define explicitly in Definition~\ref{def: famofsets} the notion of a family of sets indexed by $I$. 
Although an $I$-family of sets is a certain function-like object, it can be understood also as an object of 
a one level higher than that of a set. The corresponding notion of a ``function'' from an $I$-family $\Lambda$
to an $I$-family $M$ is that of a family-map. Operations between sets generate operations between families of 
sets and their family-maps. If the index-set $I$ is a directed set, the corresponding notion of a family of sets 
over it is that of a direct family of sets. The constructions for families of sets can be generalised 
appropriately for families of families of sets (see Section~\ref{sec: famoffam}). Families of subsets of a 
given set $X$ over an index-set $I$ are special $I$-families that deserve an independent treatment. Families 
of equivalence classes, families of partial functions, families of complemented subsets and direct families 
of subsets are some of the variations of set-indexed families of subsets that are studied here and have many 
applications in constructive mathematics.\\[2mm]
Here we apply the general theory of families of sets, in order:\\[2mm] 
\textbf{I. To reveal proof-relevance in $\BISH$.} Classical mathematics is proof-irrelevant, as it is indifferent 
to objects that ``witness'' a relation or a more complex formula. On the other extreme, Martin-L\"of type theory 
is proof-relevant, as every element of a type $A$ is a proof of the ``proposition'' $A$. Bishop's presentation
of $\BISH$ was on purpose closer to the proof-irrelevance of classical mathematics, although a form of proof-relevance
was evident in the use of several notions of moduli (of convergence, of uniform continuity, of uniform 
differentiability etc.). Focusing on membership and equality conditions for sets given by appropriate existential
formulas we define certain families of proof-sets that provide a $\BHK$-interpretation within $\BST$ of formulas 
that correspond to the standard atomic formulas of a first order theory. With the machinery of the general theory
of families of sets this $\BHK$-interpretation within $\BST$ is extended to complex formulas. Consequently, we can associate to many formulas $\phi$ of $\BISH$ a set $\Prf(\phi)$ of ``proofs'' or witnesses of $\phi$. Abstracting 
from several examples of totalities in $\BISH$ we define the notion of a set with a proof-relevant equality, 
and of a Martin-L\"of set, a special case of the former, the equality of which corresponds to the identity type of 
a type in intensional $\MLTT$. Through the concepts and results of $\BST$ notions and facts of $\MLTT$ and its 
extensions (either with the axiom of function extensionality, or with Vooevodsky's axiom of univalence) can be
translated into $\BISH$. While Bishop's theory of sets is standardly understood through its translation to 
$\MLTT$ (see e.g.,~\cite{CDPS05}), the development of $\BST$ offers a (partial) translation in the converse 
direction.\\[2mm]
\textbf{II. To develop the theory of spectra of Bishop spaces.} A Bishop space is a constructive, 
function-theoretic alternative to the notion of a topological space. A Bishop topology $F$ on a set $X$ is a 
subset of the real-valued function $\D F(X)$ on $X$ that includes the constant functions and it is closed 
under addition, composition with Bishop continuous functions $\BR$ from $\Real $ to $\Real$, and uniform limits. 
Hence, in contrast to topological spaces, continuity of real-valued functions is a primitive notion and a concept
of open set comes a posteriori. A Bishop topology on a set can be seen as an abstract and 
constructive approach to the ring of continuous functions $C(X)$ of a topological space $X$. 
Associating appropriately a Bishop topology to the set $\lambda_0(i)$ of a family of sets over a set $I$, 
for every $i \in I$, the notion of a spectrum of Bishop spaces is defined. If $I$ is a directed set, we get a 
direct spectrum. The theory of direct spectra of Bishop spaces and their limits is developed in 
Chapter~\ref{chapter: bspaces}, in analogy to the classical theory of spectra of topological spaces and their limits. 
The constructive theory of spectra of other structures, like groups, or rings, or modules, can be developed 
along the same lines.\\[2mm]
\textbf{III. To reformulate predicatively the basics of Bishop-Cheng measure theory.} The standard approach 
to measure theory (see e.g.,~\cite{Ta73},~\cite{Ha74}) is to take measure as a primitive notion, and to
define integration with respect to a given measure. An important alternative, and, as argued by 
Segal in~\cite{Se54} and~\cite{Se65}, a more natural approach to measure theory,
is to take the integral on a certain set of functions as a primitive notion, extend its definition to an appropriate, 
larger set of functions, and then define measure at a later stage. This is the idea of the Daniell integral, 
defined by Daniell in~\cite{Da18},
which was taken further by Weil, Kolmogoroff, and Carath\'eodory (see~\cite{We40},~\cite{Ko48}, and~\cite{Ca56},
respectively).

In the general framework of constructive-computable mathematics, there are many approaches 
to measure and probability theory. There is an extended literature in intuitionistic measure theory 
(see e.g.,~\cite{He56}), in measure theory within the computability framework of Type-2 Theory of Effectivity
(see e.g.,~\cite{Ed09}), in Russian constructivism (especially in the work of \v Sanin~\cite{Sa68}
and Demuth~\cite{BD91}), in type theory, where the main interest lies in the creation of probabilistic programming 
(see e.g.,~\cite{BAV12}), and recently also in homotopy type theory (see~\cite{FS16}), where 
homotopy type theory (see~\cite{HoTT13}) is applied to probabilistic programming. 

Within $\BISH$, measure and probability theory have taken two main directions.
The first direction, developed by Bishop and Cheng in~\cite{BC72} and by Chan in~\cite{Ch72}$-$\cite{Ch75},
is based on the notion of integration space, a constructive version of the Daniell integral, as
a starting point of constructive measure theory. Following the aforementioned 
spirit of classical algebraic integration theory, Bishop and Cheng defined first
the notion of an integrable function through the notion of an integration space, and afterwords
the measure of an integrable set. In their definition of integration space though, Bishop and Cheng used the 
impredicative concept $\C F(X)$ of all partial functions from a set $X$ to $\Real$. Such a notion makes the
extraction of the computational content of $\CMT$ and the implementation of $\CMT$ in some programming language impossible. 
The second direction to constructive measure theory, developed by Coquand, Palmgren and Spitters 
in~\cite{CP02},~\cite{Sp06} and~\cite{CS09}, is based on the recognition of the above problem of the 
Bishop-Cheng theory and of the advantages of working within the abstract, algebraic, and point-free framework of 
Boolean rings or of vector lattices. In analogy to Segal's notion of a probability algebra, the starting notion 
is a boolean ring equipped with 
an inequality and a measure function, which is called a measure ring, on which integrable and measurable 
functions can be defined. One can show that the integrable sets of Bishop-Cheng form a measure ring. In general, 
the second direction to constructive measure theory is considered technically and conceptually simpler.

In Chapter~\ref{chapter: measure} we reconstruct the Bishop-Cheng notion of measure space within $\BST$, where a 
set of measurable sets is not an appropriate set of complemented subsets, as it is usually understood, but an 
appropriate 
set-indexed family of complemented subsets. This fact is acknowledged by Bishop in~\cite{Bi70}, but it is 
completely suppressed later by him and his collaborators (Cheng and Chan). A similar indexing appears in a 
predicative formulation of the Bishop-Cheng notion of an integration space.\\[1mm]

The notions of a set-indexed family of sets and of a set-indexed family of subsets of a given set are shown here to be
important tools in the precise formulation of abstract notions in constructive mathematics. Avoiding them, 
makes the reading of constructive mathematics easier and very close to the reading of classical mathematics. 
Using them, makes the writing of constructive mathematics more precise, and seriously enriches its content. 

As the fundamental notion of a family of sets can be described both in categorical and type-theoretic terms,
many notions and constructions from category theory and dependent type theory are represented in $\BST$. 
While category theory and standard set-theory, or dependent type theory and standard set-theory do not match perfectly,
large parts of category theory and dependent type theory are reflected naturally in Bishop Set Theory (see also section~\ref{sec: typescats}).

\section{Notes}
\label{sec: notes0}

\begin{note}\label{not: onbishopapers}
\normalfont
Regarding the exact time that Bishop's unpublished papers~\cite{Bi68a} and~\cite{Bi68b} were written, 
it was difficult to find an answer. Bishop's scheme of presenting a formal system for $\BISH$ and of elaborating its implementation in some functional programming language is found both in~\cite{Bi70} and in Bishop's unpublished papers. 
The first is Bishop's contribution to the proceedings of the Buffalo meeting in 1968 that were 
published in~\cite{KMV70}. As Per Martin-L\"of informed me, Bishop was not present at the meeting. 
The presentation of the formal system $\Sigma$ and its presentation as a programming language in~\cite{Bi70} is 
very sketchy. Instead, the presentation of the type theory for $\BISH$ in~\cite{Bi68a}, and its presentation as
a programming language in~\cite{Bi68b} is an elaborated enterprise. I have heard a story of an unsuccessful 
effort of Bishop to publish~\cite{Bi68a}, due to some parallels between~\cite{Bi68a} and de Bruijn's work. 
According to that story, Bishop was unwilling to pursue the publication of his type-theoretic formalism after
that rejection. In any event, Bishop's unpublished papers must have been written between 1967 and 1970. Maybe, 
the period between 1968 and 1969 is a better estimation. In October 1970 Bishop and Cheng sent to the editors
of the Memoirs of the American Mathematical Society their short monograph~\cite{BC72}, a work that deviates a 
lot from the predicative character of~\cite{Bi67}. In my view, the papers~\cite{Bi68a} and~\cite{Bi68b} do not
fit to Bishop's period after 1970.
\end{note}

\begin{note}[The presentation axiom for setoids]\label{not: onpresentationaxiom}
\normalfont
If $A : \C U$, then, by Martin-L\"of's $J$-rule, $=_A$ is the least reflexive relation on $A$, and 
$\varepsilon A := (A, =_A)$ is the \textit{free setoid}\index{free setoid} on $A$. According to the universal 
property of a free setoid, for every setoid $\C B := (B, \sim_B)$ and every function $f : A \to B$, there is
a \textit{setoid-map} $\varepsilon f \colon A \to \C B$ such that the following left diagram commutes 
\begin{center}
\begin{tikzpicture}

\node (E) at (0,0) {$A$};
\node[above=of E] (F) {$A$};
\node[right=of F] (A) {$B$};
\node[right=of A] (K) {$A$};
\node[right=of K] (L) {$B$.};
\node[below=of K] (M) {$P$};

\draw[->,dashed] (E)--(A) node [midway,right] {$ \ \varepsilon f$};
\draw[->] (F)--(E) node [midway,left] {$\id_A$};
\draw[->] (F)--(A) node [midway,above] {$f$};

\draw[->,dashed] (M)--(K) node [midway,left] {$h$};
\draw[->>] (K)--(L) node [midway,above] {$f$};
\draw[->] (M)--(L) node [midway,right] {$\ g$};

\end{tikzpicture}
\end{center}
To show this, let $(\varepsilon f)(a) := f(a),$ and since $=_B$ is the least reflexive relation on $B$, we get
$a =_A a{'} \To (\varepsilon f)(a) =_B (\varepsilon f)(a{'})$, hence $f(a) \sim_B f(a{'})$.
A setoid $\C A$ is a \textit{choice setoid},\index{choice setoid}  if every $f: X \twoheadrightarrow A$, has a right 
inverse i.e., there $g \colon A \to X$ such that $f \circ g = \id_A$.
With the use of the type-theoretic axiom of choice (see~\cite{HoTT13}, section 1.6)  one can show that the free setoid 
$(A, =_A)$ is a choice setoid. Using the identity map, every setoid $\C A$ is the quotient of the free setoid on $A$, 
hence every setoid is the quotient of a choice setoid. 
If $\C C$ is a category, an object $P$ of $\C C$ is called \textit{projective},\index{projective object of a category} 
if for every objects $A, B$ of $\C C$ and every arrow $f : A \twoheadrightarrow B$ and $g \colon P \to B$, 
there is $h \colon P \to A$ such that the above right diagram commutes.
A category 
$\C C$ satisfies the presentation axiom,\index{presentation axiom} if for every object $C$ in $\C C$ there 
is $f: P \twoheadrightarrow C$, where
$P$ is projective. For the relation between the presentation axiom and various choice principles
see~\cite{Ra06}.
It is immediate to show that a projective setoid is a choice setoid. For the converse, and following~\cite{CDPS05}, 
p.~74, let $(P, \sim_P)$ be a choice setoid. To show that  it is a projective, we need to define a setoid-map $h$, 
given setoid maps $f$ and $g$ as above. Let 
\[Q := \sum_{(a,p) : A \times P}f(a) =_B g(p),\]
and let the projections $p_1 : Q \to A,$, where $p_1(a,p,e) := a$,
and $p_2 \colon Q \to P$, where $p_2(a,p,e) := p$. By the definition of $Q$ we get $f \circ p_1 = g \circ p_2$.
Since $p_2 \colon Q \twoheadrightarrow P$ and $P$ is a choice set, there is $k \colon P \to Q$ such that 
$p_2 \circ k = \id_P$. If $h := p_1 \circ k$, then
\begin{center}
\begin{tikzpicture}

\node (E) at (0,0) {$P$};
\node[above=of E] (F) {$A$};
\node[right=of F] (A) {$B$};
\node[left=of E] (B) {$ \ Q$};
\node[left=of B] (C) {$ \ P$};

\draw[->>] (F)--(A) node [midway,above] {$f$};
\draw[->] (E)--(A) node [midway,right] {$\ g$};
\draw[->,bend left] (C) to node [midway,left] {$h \ \ $} (F) ;
\draw[->>] (B)--(E) node [midway,below] {$p_2$};
\draw[->] (B)--(F) node [midway,left] {$p_1 \ $};
\draw[->] (C)--(B) node [midway,below] {$k$};

\end{tikzpicture}
\end{center}
$f \circ (p_1 \circ k) = (f \circ p_1) \circ k = (g \circ p_2) \circ k = g \circ (p_2 \circ k) = 
g \circ \id_P = g$. Consequently, every setoid is the surjective image of a choice setoid, hence of a projective setoid.
\end{note}

\begin{note}\label{not: onTYPES18paper}
\normalfont
A very first and short presentation of $\BST$  is found in~\cite{Pe19c}, where there we write 
$\CSFT$ instead of $\BST$. In~\cite{Pe19c} we also expressed
dependency through the universe of functions $\D V_1$ i.e., the totality of triplets 
$(A, B, f)$, where $A, B$ are sets and $f$ is a function from $A$ to $B$. Since dependent operations
are explicitly used by Bishop e.g., in the definition of the intersection $\bigcap_{t \in T}\lambda(t)$
of a $T$-family of subsets $(\lambda(t))_{t \in T}$ of a set $X$, while $\D V_1$ is neither explicitly, nor implicitly, mentioned, we 
use here the former concept.
\end{note}

\begin{note}\label{not: onfams(C)ZF}
\normalfont
As it is noted by Palmgren in~\cite{Pa12a}, p.~35, in $\ZF$, and also in its constructive version $\CZF$, a 
family of sets is represented by the fibers of a function $\lambda \colon B \to I$, where the fibers 
$\lambda_i := \{ b \in B \mid \lambda(b) = i\}$ of $\lambda$, for every $i \in I$, represent the sets of the
family. Hence the notion of a family of sets is reduced to that of a set. As this reduction rests on the 
replacement scheme, such a reduction is not possible neither in $\MLTT$ nor in $\BST$.
\end{note}

\chapter{Fundamentals of Bishop Set Theory}
\label{chapter: BST}

We present the basic elements of $\BST$, a reconstruction of Bishop's informal theory of sets, as this is 
developed in chapters 3 of~\cite{Bi67} and~\cite{BB85}. The main new features of $\BST$, with respect to Bishop's 
account, are the explicit use of the universe $\D V_0$ of sets and the elaboration of the study of dependent 
operations over a non-dependent assignment routines from a set to $\D V_0$. The first notion is implicit in Bishop's work,
while the second is explicitly mentioned, although in a rough way. These concepts are necessary to the concrete 
definition of a set-indexed family of sets, the main object of our study, which is only roughly mentioned by Bishop.
The various notions of families of sets introduced later, depend on the various notions of sets, subsets and assignment 
routines developed in this chapter.

\section{Primitives}
\label{sec: primitives}

The logical framework of $\BST$ is first-order intuitionistic logic with equality (see~\cite{SW12}, chapter 1). 
This primitive equality between terms is 
denoted by $s := t$\index{$s := t$}, and it is understood as a \textit{definitional}, or \textit{logical}, 
equality\index{definitional equality}\index{logical equality}. I.e., we read the equality $s := t$ as 
``the term $s$ is by definition equal to the term $t$''. If $\phi$ is an appropriate formula, for the standard axiom 
for equality $[a := b \ \& \ \phi(a)] \To \phi(b)$ we use the notation $[a := b \ \& \ \phi(a)] :\To \phi(b)$.
The equivalence notation $:\TOT$ is understood in the same way.
The set $(\Nat =_{\Nat}, \neq_{\Nat})$ of natural numbers, where its canonical equality is given by\index{$=_{\Nat}$} 
$m =_{\Nat} n :\TOT m := n$, and its canonical inequality by $m \neq_{\Nat} n :\TOT \neg(m =_{\Nat} n)$, is primitive.
The standard Peano-axioms are associated to $\Nat$. 

A global operation $( \cdot, \cdot)$ of pairing\index{pairing} is also considered primitive. I.e., if 
$s, t$ are terms, their pair $(s,t)$ is a new term. The corresponding equality axiom is  
$(s, t) := (s{'}, t{'}) :\TOT s := s{'} \ \& \ t := t{'}$. The $n$-tuples of given terms, for every $n$ larger
than $2$, are definable. The global projection routines $\prb_1(s, t) := s$ and $\prb_2(s, t) := t$ are also considered
primitive. The corresponding global projection routines for any $n$-tuples are definable.

An undefined notion of mathematical construction, or algorithm, or of finite routine is considered as primitive.
The main primitive objects of $\BST$ are totalities\index{totality} and assignment 
routines\index{assignment routine}. Sets are special totalities and 
functions are special assignment routines, where an assignment routine is a a special finite routine. 
All other equalities in $\BST$ are equalities on totalities defined though an equality condition.
A predicate\index{predicate on a set} on a set $X$ is a bounded formula $P(x)$ with $x$ a free variable ranging over $X$,
where a formula is bounded\index{bounded formula}, if every quantifier occurring in it is over a given set.

\section{Totalities}
\label{sec: totalities}

\begin{definition}\label{def: totalities}
 \normalfont (i) 
\itshape A \textit{primitive set}\index{primitive set} $\D A$ is a totality with a given membership 
 $x \in \D A$, and a given equality $x =_{\D A} y$, that satisfies axiomatically the properties of 
 an equivalence relation. The set $\Nat$ of natural numbers is the only primitive set considered here.\\[1mm]
 \normalfont (ii) 
\itshape A $($non-inductive$)$\textit{defined totality}\index{defined totality} $X$ is defined by a membership
condition $x \in X : \TOT \C M_X(x),$ where 
$\C M_X$ is a formula with $x$ as a free variable. 
If $X, Y$ are defined totalities with membership conditions $\C M_X $ and $\C M_Y$, respectively, we define 
$X := Y : \TOT \big[\C M_X (x) : \TOT \C M_Y (x)\big]$, and in this case 
we say that $X$ and $Y$ are \textit{definitionally equal}\index{definitionally equal defined totalities}.\\[1mm]
\normalfont (iii)
 \itshape There is a special ``open-ended'' defined totality $\D V_0$, which is called the universe of sets. $\D V_0$
 is not defined through a membership-condition, but in an open-ended way. When we say that a defined totality $X$ is
 considered to be a set we ``introduce'' $X$ as an element of $\D V_0$. We do not add the corresponding induction, 
 or elimination principle, as we want to leave open the possibility of adding new sets in $\D V_0$. \\[1mm]
 \normalfont (iv) 
\itshape A \textit{defined preset}\index{defined preset} $X$, or simply, a preset\index{preset}, is a defined totality 
$X$ the membership condition $\C M_X$ of which expresses a construction that can, in principle, 
be carried out in a finite time. Formally this is expressed by the requirement 
that no quantification over $\D V_0$ occurs in $\C M_X$.\\[1mm]
\normalfont (v) 
\itshape A \textit{defined totality $X$ with equality}\index{defined totality with equality}, or simply, a
totality $X$ with equality\index{totality with equality} is a defined totality $X$ equipped with an equality condition 
$x =_X y : \TOT \C E_X(x, y)$, where $\C E_X(x,y)$ is a formula with free variables $x$ and $y$ that 
satisfies the conditions of an equivalence relation i.e., $\C E_X(x, x)$ and
$\C E_X(x, y) \To \C E_X(y, x)$, and $[\C E_X(x, y) \ \& \ \C E_X(y, z)] \To \C E_X(x, y)$.
Two defined totalities with equality $(X,=_X)$ and $(Y, =_Y)$ are definitionally equal, if 
$\C M_X (x) : \TOT \C M_Y (x)$ and $\C E_X (x, y) : \TOT \C E_Y (x,y)$.\\[1mm]
\normalfont (vi) 
\itshape A \textit{defined set}\index{defined set} is a preset with a given equality.\\[1mm]
\normalfont (vii) 
\itshape A \textit{set}\index{set} is either a primitive set, or a defined set.\\[1mm]
\normalfont (viii) 
\itshape A totality is a \textit{class}\index{class}, if it is the universe $\D V_0$, or if 
quantification over $\D V_0$ occurs in its membership condition.
\end{definition}

\begin{definition}\label{ def: productofsets}
If $X, Y$ are sets, their \textit{product}\index{product of sets} $X \times Y$\index{$X \times Y$} is the 
defined totality
with equality
\[ (x, y) \in X \times Y : \TOT x \in A \ \& \ y \in B, \]
\[ z \in X \times Y :  \TOT \exists_{x \in A}\exists_{y \in B}\big(z := (x, y)\big). \]
\end{definition}

$X \times Y$ is considered to be a set, and its membership condition is written simpler as follows:
\[ (x, y) =_{X \times Y} (x{'}, y{'}) : \TOT x =_X x{'} \ \& \ y =_Y y{'}. \]

\begin{definition}\label{def: extsubset}
A bounded formula on a set $X$ is called an \textit{extensional property} on $X$\index{extensional property on a set}, if 
\[\forall_{x, y \in X}\big([x =_{X } y \ \& \ P(x)] \To P(y)\big). \]
The totality $X_P$\index{$P_X$} generated by $P(x)$ is defined by $x \in X_P : \TOT x \in X \ \& \ P(x)$, 
\[ x \in X_P : \TOT x \in X \ \& \ P(x), \]
and the equality 
of $X_P$ is inherited from the equality 
of $X$. We  also write $X_P := \{x \in X \mid P(x)\}$.
The totality $X_P$ is considered to be a set, and it is called the  
\textit{extensional subset}\index{extensional subset} of $X$ generated by $P$.
\end{definition}

Using the properties of an equivalence relation, it is immediate to show that an equality condition
$\C E_X(x,y)$ on a totality $X$ is an extensional property on the product $X \times X$ i.e., 
$[(x, y) =_{X \times Y} (x{'}, y{'}) \ \&\ x =_X y] \To x{'} =_X y{'}$. Let the following extensional subsets of 
$\Nat$: \index{$\D 1$} \index{$\D 2$}
 \[ \D 1 := \{x \in \Nat \mid x =_{\Nat} 0\} := \{0\}, \]
 \[ \D 2 := \{x \in \Nat \mid x =_{\Nat} 0 \ \vee x =_{\Nat} 1\} := \{0, 1\}. \]
Since $n =_{\Nat} m :\TOT n := m$, the property $P(x) :\TOT x =_{\Nat} 0 \ \vee x =_{\Nat} 1$ is extensional.

\begin{definition}\label{def: diagonal}
If $(X, =_X)$ is a set, its \textit{diagonal}\index{diagonal of a set} 
\index{$D(X, =_X)$} 
is the extensional subset of $X \times X$
\[ D(X, =_X) := \{(x, y) \in X \times X \mid x =_X y\}. \]
If $=_X$ is clear from the context, we just write $D(X)$\index{$D(X)$}. 
\end{definition}

\begin{definition}\label{def: apartness}
Let $X$ be a set. An \textit{inequality}\index{inequality} on $X$, or an 
\textit{apartness relation}\index{apartness relation} on $X$, is a relation $x \neq_X y$\index{$x \neq_X y$} such that 
the following conditions are satisfied:\\[1mm]
$(\Ap_1)$ $\forall_{x, y \in X}\big(x =_X y \ \& \ x \neq_X y \To \bot \big)$.\\
$(\Ap_2)$ $\forall_{x, y \in X}\big(x \neq_X y \To y \neq_X x\big)$.\\
$(\Ap_3)$ $\forall_{x, y \in X}\big(x \neq_X y \To \forall_{z \in X}(z \neq_X x \ \vee \ z \neq_X y)\big)$.\\[1mm]
We write $(X, =_X, \neq_X)$\index{$(X, =_X, \neq_X)$} to denote the equality-inequality structure of a
set\index{equality-inequality structure of a set} $X$, and for simplicity we refer the \textit{set}
$(X, =_X, \neq_X)$. The set $(X, =_X, \neq_X)$ is called \textit{discrete}\index{discrete}, if 
\[ \forall_{x, y \in X}\big(x =_X y \ \vee \ x \neq_X y\big). \]
An inequality $\neq_X$ on $X$ is called \textit{tight}\index{tight inequality}, 
if $\neg(x \neq_X y) \To x =_X y$, for every $x,y \in X$.
\end{definition}

\begin{remark}\label{rem: apartness1}
An inequality relation $x \neq_X y$ is extensional on $X \times X$.
\end{remark}

\begin{proof}
We show that if $x, y \in X$ such that $x \neq y$, and if $x{'}, y{'} \in X$ such that $x{'} =_X x$ 
and $y{'} =_X y$, then 
$x{'} \neq y{'}$. By $\Ap_3$ we have that $x{'} \neq x$, which is excluded from $\Ap_1$, or $x{'} \neq y$,
which has to be the case. Hence, $y{'} \neq x{'}$, or $y{'} \neq y$. Since the last option is excluded similarly,
we conclude that $y{'} \neq x{'}$, hence $x{'} \neq y{'}$.
\end{proof}

If $\neq_X$ is an inequality on $X$, and $P(x)$ is an extensional property on $X$, then $X_P$
inherits the inequality from $X$. Since $n \neq_{\Nat} m :\TOT \neg(n =_{\Nat} m)$,
the sets $\Nat$, $\D 1$, and $\D 2 $ are discrete. Clearly, if $(X, =_X, \neq_X)$ is discrete, then $\neq_X$ is tight.

\begin{remark}\label{rem: proddiscrete}
Let the sets $(X, =_X, \neq_X)$ and $(Y, =_Y, \neq_Y)$.\\[1mm]
\normalfont (i) 
\itshape The canonical inequality on $X \times Y$\index{canonical inequality on $X \times Y$} induced 
by $\neq_X$ and $\neq_Y$, which is 
defined by 
\[ (x, y) \neq_{X \times Y} (x{'}, y{'}) :\TOT x \neq_X x{'} \ \vee \ y \neq_Y y{'}, \]
for every $(x,y)$ and $(x{'}, y{'}) \in X \times Y$, is an inequality on $X \times Y$.\\[1mm]
\normalfont (ii) 
\itshape If $(X, =_X, \neq_X)$ and $(Y, =_Y, \neq_Y)$ are discrete, then 
$(X \times Y, =_{X \times Y}, \neq_{X \times Y})$ is discrete.
\end{remark}

\begin{proof}
The proof of (i) is immediate. To show (ii), let $(x, y), (x{'}, y{'}) \in X \times Y$. By our 
hypothesis $x =_X x{'} \ \vee \ x \neq_X x{'}$ and 
$y =_Y y{'} \ \vee \ y \neq_Y y{'}$. If $x =_X x{'}$ and $y =_Y y{'}$, then $(x, y) =_{X \times Y} (x{'}, y{'})$. 
In any other case we get $(x, y) \neq_{X \times Y} (x{'}, y{'})$.
\end{proof}

Uniqueness of an element of a set $X$ with respect to some property $P(x)$ on $X$ means that all elements 
of $X$ having this property are $=_X$-equal. We use the following abbreviation:
\[\exists_{!x \in X}P(x) :\TOT \exists_{x \in X}\big(P(x) \ \& \ \forall_{z \in X}\big(P(z) \To z =_X x\big)\big). \]

\begin{definition}\label{def: contractible}
Let $(X, =_X)$ be a set. \\[1mm]
\normalfont (i) 
\itshape $X$ is \textit{inhabited}\index{inhabited set}, if
$\exists_{x \in X}\big(x =_X x\big)$.\\[1mm]
\normalfont (ii) 
\itshape $X$ is a \textit{singleton}\index{singleton}, or \textit{contractible}\index{contractible}, or 
a $(-2)$-set\index{$(-2)$-set}, if 
$\exists_{x_0 \in X}\forall_{x \in X}\big(x_0 =_\mathsmaller{X} x\big)$. In this case,
$x_0$ is called a centre of contraction\index{centre of contraction} for $X$.\\[1mm]
\normalfont (iii) 
\itshape $X$ is a \textit{subsingleton}\index{subsingleton}, or a \textit{mere proposition}\index{mere proposition},
or a $(-1)$-set\index{$(-1)$-set}, if $\forall_{x, y \in X}\big(x =_\mathsmaller{X} y\big)$.\\[1mm]
\normalfont (iv) 
\itshape The \textit{truncation}\index{truncation} of $(X, =_X)$ is the 
set\index{$\mathsmaller{\mathsmaller{\lvert \lvert =_X \rvert \rvert}}$} 
$(X, \mathsmaller{\mathsmaller{\lvert \lvert =_X \rvert \rvert}})$, where
\[ x \ \mathsmaller{\mathsmaller{\lvert \lvert =_X \rvert \rvert}} \ y :\TOT x =_X x \ \& \ y =_X y. \]
We use the symbol $||X||$ to denote that the set $X$ is equipped with the truncated 
equality $\mathsmaller{\mathsmaller{\lvert \lvert =_X \rvert \rvert}}$.\index{$||X||$}
\end{definition}

Clearly, $ x \ \mathsmaller{\mathsmaller{\lvert \lvert =_X \rvert \rvert}} \ y$, for every $x, y \in X$, and
$(X, \mathsmaller{\mathsmaller{\lvert \lvert =_X \rvert \rvert}})$ is a subsingleton.

\section{Non-dependent assignment routines}
\label{sec: ndar}


\begin{definition}\label{def: ndar}
Let $X, Y$ be totalities. A \textit{non-dependent assignment routine}\index{non-dependent assignment routine}
$f$ from $X$ to 
$Y$, in symbols $f \colon X \sto Y$\index{$f \colon X \sto Y$}, is a finite routine that assigns an element $y$ of $Y$
to each given element $x$ of $X$. 
In this case we write $f(x) := y$. If $g \colon X \sto Y$, let 
\[f := g : \TOT \forall_{x \in X}\big(f(x) := g(x)\big). \]
If $f := g$, we say that  $f$ and $g$ are \textit{definitionally equal}\index{definitionally equal functions}.
If $(X, =_X)$ and $(Y, =_Y)$ are sets, an operation\index{operation} from $X$ to $Y$ is a non-dependent assignment routine 
from $X$ to $Y$, while a function\index{function} from $X$ to $Y$, in symbols $f \colon X \to Y$\index{$f \colon X \sto Y$},
is an operation from $X$ to $Y$ that respects equality i.e.,
\[\forall_{x, x{'} \in X}\big(x =_X x{'} \To f(x) =_Y f(x{'})\big). \]
If $f \colon X \sto Y$ is a function from $X$ to $Y$, we say
that $f$ is a function, without mentioning the expression ``from $X$ to $Y$''.
A function $\fXY$ is an \textit{embedding}\index{embedding}, in symbols 
$f \colon X \hookrightarrow Y$\index{$f \colon X \hookrightarrow Y$}, if 
\[\forall_{x, x{'} \in X}\big( f(x) =_Y f(x{'}) \To x =_X x{'}). \]
Let the sets $(X, =_X, \neq_X)$ and $(Y, =_Y, \neq_Y)$. A function $f \colon X \to Y$ is strongly extensional,
\index{strongly extensional function} if
\[ \forall_{x, x{'} \in X}\big(f(x) \neq_Y f(x{'}) \To x \neq_X x{'}\big). \]
If $\simeq_X$ is another equality on $X$, we use a new symbol e.g., $X^*$, for the same totality $X$. 
When we write $f \colon X^* \to Y$, then $f$ is a function from $X$, equipped with the equality 
$\simeq_X$, to $Y$.

\end{definition}

If $X$ is a set, the \textit{identity map}\index{identity map on a set} $\id_X$\index{$\id_X$} on $X$ is the operation
$\id_X \colon X \sto X$, defined by $\id_X (x) := x$, for every $x \in X$. Clearly, $\id_X$ is an embedding, 
which is strongly extensional, if $\neq_X$ is a given inequality on $X$. 
If $Y$ is also a set, the projection maps $\pr_X$ and 
$\pr_Y$\index{$\pr_Y$}\index{$\pr_X$} on $X$ and $Y$, respectively, are the operations
$\pr_X \colon X \times Y \sto X$ and $\pr_Y \colon X \times Y \sto Y$, where
\[ \pr_X (x, y) := \prb_1 (x, y) := x \ \ \& \ \ \pr_Y (x, y) := \prb_2 (x, y) := y; \ \ \ \ (x, y) \in X \times Y. \] 
Clearly, the operations $\pr_X$ and $\pr_Y$ are functions, which are strongly extensional, if $\neq_X, \neq_Y$ are
inequalities on $X, Y$, and $\neq_{X \times Y}$ is the canonical inequality on $X \times Y$ induced from them.
After introducing the universe $\D V_0$ of sets in section~\ref{sec: universe}, we shall define non-dependent
assignment routines from a set to a totality, like $\D V_0$, which is not considered to be a set.
In most of the cases the non-dependent assignment routines defined here have a set as a domain.There are cases though, 
see e.g., Definitions~\ref{def: unionofsubsets}~\ref{def: intofsubsets},~\ref{def: interiorunion}, 
and~\ref{def: intfamilyofsubsets}, where 
a non-dependent assignment routine is defined on a totality, before showing that this totality
is a set. \textit{We never define a non-dependent assignment routine from a class to a totality}.

Let the operation $m^* \colon \Real \sto \D Q$, defined by $m^*(a) := q_m$, where a real number $a$ 
is a regular sequence of rational numbers $(q_n)_n$ (see~\cite{BB85}, p.~18), and $q_m$ is the $m$-term of this sequence.
for some fixed $m$. The operation $m^*$ is an example of an operation, which is not a function, since unequal
real numbers, with respect to the definition of $=_{\Real}$ in~\cite{BB85}, p.~18, may have equal $m$-terms in $\D Q$.  
To define a function $\fXY$, \textit{first} we define the operation $f \colon X \sto Y$, and \textit{afterwords}
we prove that $f$ is a function (from $X$ to $Y$).

The \textit{composition}\index{composition of operations} $g \circ f$\index{$g \circ f$} of the operations $f \colon X \sto Y$ 
and $g \colon Y \sto Z$ is the operation $g \circ f \colon X \sto Z$, defined by $(g \circ f)(x) := g(f(x))$, for every 
$x \in X$. Clearly, $g \circ f$ is a function, if $f$ and $g$ are functions. If $h \colon Z \sto W$, notice the
following definitional equalities
\[ f \circ \id_X := f, \ \ \ \ \id_Y \circ f := f, \ \ \ \ h \circ (g \circ f) := (h \circ g) \circ f. \]

A diagram commutes always with respect to the equalities of the related sets.
E.g., the commutativity of the following diagram is the equality $e (f(x)) =_W g(h(x))$, for every $x \in X$. 
\begin{center}
\begin{tikzpicture}

\node (E) at (0,0) {$Z$};
\node[right=of E] (F) {$W$.};
\node[above=of F] (A) {$Y$};
\node [above=of E] (D) {$X$};

\draw[->] (E)--(F) node [midway,below] {$g$};
\draw[->] (D)--(A) node [midway,above] {$f $};
\draw[->] (D)--(E) node [midway,left] {$h$};
\draw[->] (A)--(F) node [midway,right] {$e$};

\end{tikzpicture}
\end{center}

\begin{definition}\label{def: setsoffunctions}
 Let $X, Y$ be sets, and $\neq_Y$ an inequality on $Y$. The totality $\D O(X, Y)$\index{$\D O(X, Y)$} of
 operations from $X$ to $Y$ 
 is equipped with the following canonical equality and inequality: 
 \vspace{-2mm}
 \[ f =_{\D O(X, Y)} g : \TOT \forall_{x \in X}\big(f(x) =_Y f(x)\big), \]
 \[ f \neq_{\D O(X, Y)} g : \TOT \exists_{x \in X}\big(f(x) \neq_Y f(x)\big). \]
 The totality $\D O(X, Y)$ is considered to be a set. The set $\D F(X, Y)$\index{$\D F(X, Y)$} of functions
 from $X$ to $Y$
 is defined by separation on $\D O(X, Y)$ through the extensional property $P(f) :\TOT
 \forall_{x, x{'} \in X}\big(x =_X x{'} \To f(x) =_Y f(x{'})\big)$. The equality $=_{\D F(X, Y)}$ and 
 the inequality $\neq_{\D F(X, Y)}$ are inherited from $=_{\D O(X, Y)}$ and $\neq_{\D O(X, Y)}$, respectively.
\end{definition}

\begin{remark}\label{rem: apartness2}
Let the sets $(X =_X)$ and $(Y, =_Y, \neq_Y)$. If $\fXY$, let   
$x_1 \neq_X^f x_2 : \TOT f(x_1) \neq_Y f(x_2)$, for every $x_1, x_2 \in X$.\\[1mm] 
\normalfont (i) 
\itshape $x_1 \neq_X^f x_2$ is an inequality on $X$.\\[1mm]
\normalfont (ii) 
\itshape If $(Y, =_Y, \neq_Y)$ is discrete, then $(X =_X,  \neq_X^f)$ is discrete if and only if $f$ is an 
embedding.\\[1mm]
\normalfont (ii) 
\itshape If $\neq_Y$ is tight, then $\neq_X^f$ is tight if and only if $f$ is an embedding.
\end{remark}

\begin{proof}
(i) Conditions $(\Ap_1)$-$(\Ap_3)$ for $\neq_X^f$ are reduced to conditions $(\Ap_1)$-$(\Ap_3)$ for $\neq_Y$.\\
(ii) If $(X =_X,  \neq_X^f)$ is discrete, let $f(x_1) =_X f(x_2)$, for some $x_1, x_2 \in X$. Since the possibility
$x_1 \neq_X^f x_2 : \TOT f(x_1) \neq_Y f(x_2)$ is impossible, we conclude that $x_1 =_X x_2$. If $f$ is an embedding, 
and since $f(x_1) =_X f(x_2)$ or $f(x_1) \neq_Y f(x_2)$, either $x_1 =_X x_2$, or $x_1 \neq_X^f x_2$.\\
(iii) If $\neq_X^f$ is tight, and  $f(x_1) =_X f(x_2)$, 
then $\neg(x_1 \neq_X^f x_2)$,
hence $x_1 =_X x_2$. If $f$ is an embedding and $\neg(x_1 \neq_X^f x_2) \TOT \neg\big(f(x_1) \neq_Y f(x_2)\big)$,
then $f(x_1) =_X f(x_2)$, and $x_1 =_X x_2$. 
\end{proof}

 \begin{definition}\label{def: surjective}
 A function $f \colon X \to Y$ is called \textit{surjective}\index{surjective function}, if 
 $\forall_{y \in Y}\exists_{x \in X}\big(f(x) =_Y y\big)$. A function $g \colon Y \to X$ is called a modulus of
 surjectivity for $f$\index{modulus of surjectivity}, if the following diagram commutes
\begin{center}
\begin{tikzpicture}

\node (E) at (0,0) {$Y$};
\node[right=of E] (F) {$X$};
\node[right=of F] (A) {$Y$.};

\draw[->] (E)--(F) node [midway,above] {$g$};
\draw[->] (F)--(A) node [midway,above] {$f$};
\draw[->,bend right] (E) to node [midway,below] {$\id_Y$} (A);

\end{tikzpicture}
\end{center}
If $g$ is a modulus of surjectivity for $f$, we also say that $f$ is a retraction\index{retraction}
and $Y$ is a retract\index{retract} of $X$.
If $y \in Y$, the fiber\index{fiber} $\fib^f(y)$\index{$\fib^f(y)$} of $f$ at $y$ is the following 
extensional subset of $X$
\[ \fib^f(y) := \{x \in X \mid f(x) =_Y y \}. \]
A function $\fXY$ is contractible\index{contractible function}, if $\fib^f(y)$ is contractible, for every $y \in Y$.
If $\neq_Y$ is an inequality on $Y$,  the cofiber\index{cofiber} $\cofib^f(y)$\index{$\cofib^f(y)$} of $f$ at $y$ 
is the following extensional subset of $X$
\[ \cofib^f(y) := \{x \in X \mid f(x) \neq_Y y \}. \]
\end{definition}

\section{The universe of sets}
\label{sec: universe}

The totality of all sets is the universe $\D V_0$ of sets\index{universe of sets}\index{$\D V_0$},
equipped with the canonical equality 
\[ X =_{\D V_0} Y :\TOT \exists_{f \in \D F(X,Y)}\exists_{g \in \D F(Y,X)}\big(g \circ f = \id_X \ \& \ 
f \circ g = \id_Y\big) \]
\begin{center}
\begin{tikzpicture}

\node (E) at (0,0) {$X$};
\node[right=of E] (F) {$Y$};
\node[right=of F] (A) {$X$};
\node[right=of A] (B) {$Y$.};

\draw[->] (E)--(F) node [midway,above] {$f$};
\draw[->] (F)--(A) node [midway,above] {$\ \ \ g$};
\draw[->] (A)--(B) node [midway,below] {$f$};
\draw[->,bend right] (E) to node [midway,below] {$\id_X$} (A) ;
\draw[->,bend left] (F) to node [midway,above] {$\id_Y$} (B) ;

\end{tikzpicture}
\end{center}
In this case we write $(f, g) : X =_{\D V_0} Y$. If $X, Y \in \D V_0$ such that $X =_{\D V_0} Y$, we  
define the set
\[ \Eq(X, Y) := \big\{(f, g) \in \D F(X, Y) \times \D F(Y, X) \mid (f, g) : X =_{\D V_0} Y\big\} \]
of all objects that ``witness'', or ``realise'', or prove the equality $X =_{\D V_0} Y$. The equality of 
$\Eq(X, Y)$ is the canonical one i.e., $(f, g) =_{\Eq(X, Y)} (f{'}, g{'}) :\TOT f =_{\D F(X, Y)} f{'} \ \&
\ g =_{\D F(Y, X)}g{'}$. Notice that, in general, not all elements of $\Eq(X, Y)$ are equal. As in~\cite{HoTT13},
Example 3.1.9, if $X := Y := \D 2 := \{0, 1\}$, then $(\id_{\D 2}, \id_{\D 2}) \in \Eq(\D 2, \D 2)$, and if
$\sw_{\D 2} : \D 2 \to \D 2$ maps $0$ to $1$ and $1$ to $0$, then $(\sw_{\D 2}, \sw_{\D 2}) \in \Eq(\D 2, \D 2)$, while
$\sw_{\D 2} \neq \id_{\D 2}$. 

It is expected that the proof-terms in $\Eq(X, Y)$ are compatible with the properties of the equivalence relation
$X =_{\D V_0} Y$. This means that we can define a distinguished proof-term
$\refl(X) \in \Eq(X, X)$ that proves the reflexivity of $X =_{\D V_0} Y$, an operation $^{-1}$, such that
if $(f, g) : X =_{\D V_0} Y$, then $(f, g)^{-1} : Y =_{\D V_0} X$, and an operation of 
``composition'' $\ast$ of proof-terms,
such that if $(f, g) : X =_{\D V_0} Y$ and $(h, k) : Y =_{\D V_0} Z$, then $(f, g) \ast (h, k) : X =_{\D V_0} Z$.
If $h \in \D F(Y, W)$ and $k \in \D F(W, Y)$, let
$$\refl(X) := \big(\id_X, \id_X\big) \ \ \& \ \ (f, g)^{-1} := (g, f) \ \ \& \ \ 
(f, g) \ast (h, k) := (h \circ f, g \circ k).$$
It is immediate to see that these operations satisfy the \textit{groupoid laws}\index{groupoid laws}:\\[1mm]
(i) $\refl (X) \ast (f, g) =_{\Eq(X, Y)} (f, g)$ and $(f, g) \ast \refl (Y) =_{\Eq(X, Y)} (f, g)$.\\[1mm]
(ii) $(f, g) \ast (f, g)^{-1} =_{\Eq(X, X)} \refl (X)$ and $(f, g)^{-1} \ast (f, g) =_{\Eq(Y, Y)} \refl (Y)$.\\[1mm]
(iii) $\big((f, g) \ast (h, k)\big) \ast (s, t) =_{\Eq(X, W)} (f, g) \ast \big((h, k) \ast (s, t)\big)$.\\[1mm]
Moreover, the following \textit{compatibility condition}\index{compatibility condition} is satisfied:\\[1mm]
(iv) If $(f, g), (f{'}, g{'}) \in \Eq(X, Y)$ and $(h, k), (h{'}, k{'}) \in \Eq(Y, Z)$, then
if $(f, g) =_{\Eq(X, Y)} (f{'}, g{'})$ and $(h, k) =_{\Eq(Y, Z)} (h{'}, k{'})$, then  
$(f, g) \ast (h, k) =_{\Eq(X, Z)} (f{'}, g{'}) \ast (h{'}, k{'})$.

\begin{proposition}\label{prp: fiber1}
Let $X, Y$ be sets, $f \in \D F(X, Y)$ and $g \in \D F(Y,X)$. If $(f, g) \colon X =_{\D V_0} Y$,
then the set $\fib^f(y)$ is contractible, for every $y \in Y$. 
\end{proposition}

\begin{proof}
If $y \in Y$, then $g(y) \in \fib^f(y)$, as $f(g(y)) =_Y \id_Y(y) := y$. 
If $x \in X$, 
$x \in \fib^f(y) :\TOT f(x) =_Y y$, and 
$x =_X g(f(x)) =_X g(y)$ i.e., $g(y)$ is a centre of contraction for $\fib^f(y)$.
\end{proof}

 
\begin{definition}\label{def: evaluation}
Let $X, Y$ be sets. The \textit{evaluation map}\index{evaluation map}\index{$\ev_{X,Y}$} 
$\ev_{X,Y} : \D F(X, Y) \times X 
\sto Y$ is defined by $\ev_{X,Y} (f, x) := f(x)$, for every $f \in \D F(X, Y)$ and $x \in X$.
\end{definition}

\begin{proposition}\label{prp: ccc}
Let $X, Y, Z$ be sets.\\[1mm]
\normalfont (i) 
\itshape The evaluation map $\ev_{X,Y}$ is a function from $\D F(X, Y) \times X$ to $Y$.\\[1mm]
\normalfont (ii) 
\itshape For every function $h : Z \times X \to Y$, there is a unique function $\hat{h} : Z \to \D F(X, Y)$ 
such that for every $z \in Z$ and $x \in X$
$\ev_{X,Y}\big(\hat h (z), x\big) =_Y h(z, x).$ 

\end{proposition}

\begin{proof}
(i) By definition $(f, x) =_{\D F(X, Y) \times X} (f{'}, x{'})$ if and only if $f =_{\D F(X, Y)} f{'}$ and $x =_X x{'}$. 
Hence
$\ev_{X,Y}(f, x) := f(x) =_Y f{'}(x) =_Y f{'}(x{'}) := \ev_{X,Y}(f{'}, x{'})$.\\
(ii) For every $z \in Z$, we define the assignment routine $\hat h$ from $Z$ to $\D F(X, Y)$ by
$z \mapsto \hat h (z)$, where
$\hat h(z)$ is the assignment routine from $X$ to $Y$, defined by
$\hat h(z)(x) := h(z, x),$
for every $x \in X$. First we show that $\hat h (z)$ is a function from $X$ to $Y$; if $x =_X x{'}$, 
then $(z, x) =_{Z \times X} (z, x{'})$, hence $\hat h(z)(x) := h(z, x) =_Y h(z, x{'}) := \hat h(z)(x{'}).$
Next we show that the assignment routine $\hat h$ is a function from $Z$ to $\D F(X, Y)$; if $z =_Z z{'}$, 
then, if $x \in X$, and since then $(z, x) =_{Z \times X} (z{'}, x)$, we have that 
$\hat h(z)(x) := h(z, x) =_Y h(z{'}, x) := \hat h(z{'})(x).$ Since $x \in X$ is arbitrary, we conclude 
that $\hat h (z) =_{\D F(X, Y)} \hat h(z{'})$. Since
$\ev_{X,Y}\big(\hat h (z), x\big) := \hat h(z)(x)  := h(z, x),$
we get the strong from of the required equality $\ev_{X,Y} \circ (\hat h \times 1_X) := h$. 
If $g : Z \to \D F(X, Y)$ satisfying the required equality,
and if $z \in Z$, then, for every $x \in X$ we have that 
$g(z)(x) := \ev_{X,Y}\big(g(z), x\big) =_Y h(z, x) =_Y \ev_{X,Y}\big(\hat h (z), x\big) := \hat h(z)(x),$
hence $g(z) =_{\D F(X, Y)} \hat h(z)$.
\end{proof}

\section{Dependent operations}
\label{sec: doperations}

\begin{definition}\label{def: dependops}
Let $I$ be a set and $\lambda_0 \colon I \sto \D V_0$\index{$\lambda_0 \colon I \sto \D V_0$} a non-dependent assignment 
routine from $I$ to $\D V_0$. 
A \textit{dependent operation}\index{dependent operation} $\Phi$ over $\lambda_0$, in symbols
\[ \Phi \colon \bigcurlywedge_{i \in I} \lambda_0 (i), \]
is an assignment routine that assigns to each element $i$ in $I$ an element $\Phi(i)$ in the set
$\lambda_0(i)$. If $i \in I$, we call $\Phi(i)$ the $i$-\textit{component}\index{component of a dependent operation}
of $\Phi$, and we also use the notation $\Phi_i := \Phi(i)$\index{$\Phi(i)$}\index{$\Phi_i$}.
An assignment routine\index{assignment routine} is either a non-dependent assignment routine, or a dependent operation 
over some non-dependent assignment routine from a set to the universe.
If $\Psi \colon \bigcurlywedge_{i \in I} \lambda_0 (i)$, 
let $\Phi := \Psi :\TOT \forall_{i \in I}\big(\Phi_i := \Psi_i\big).$
If $\Phi := \Psi$, we say that $\Phi$ and $\Psi$ are definitionally equal\index{definitionally
equal dependent operations}. 

\end{definition}

Let the non-dependent assignment routines $\lambda_0 \colon I \sto \D V_0, \mu_0 \colon I \sto \D V_0, 
\nu_0 \colon I \sto \D V_0$ and $\kappa_0 \colon I \sto \D V_0$. Let 
$\D F(\lambda_0, \mu_0) \colon
I \sto \D V_0$ be defined by $\D F (\lambda_0, \mu_0)(i) := \D F(\lambda_0(i), \mu_0(i)$, for every $i \in I$. The 
identity operation\index{identity operation} $\Id_{\lambda_0}$\index{$\Id_{\lambda_0}$} over $\lambda_0$ is the 
dependent operation 
\[ \Id_{\lambda_0} \colon \bigcurlywedge_{i \in I}\D F(\lambda_0(i), \mu_0(i)) \ \ \ \ \Id_{\lambda_0}(i) := 
\id_{\lambda_0(i)}; \ \ \ \ i \in I.\]
Let 
$\Psi \colon \bigcurlywedge_{i \in I}\D F(\mu_0(i), \nu_0(i))$ and 
$\Phi \colon \bigcurlywedge_{i \in I}\D F(\lambda_0(i), \mu_0(i))$. Their
\textit{composition}\index{composition of dependent operations} $\Psi \circ \Phi$\index{$\Psi \circ \Phi$}
is defined by 
\[ \Psi \circ \Phi \colon \bigcurlywedge_{i \in I}\D F(\lambda_0(i), \nu_0(i)) \ \ \ \ 
 (\Psi \circ \Phi)_i := \Psi_i \circ \Phi_i; \ \ \ \ i \in I. \]
If $\Xi \colon \bigcurlywedge_{i \in I}\D F(\nu_0(i), \kappa_0(i))$, notice the
following definitional equalities
\[ \Phi \circ \Id_{\lambda_0} := \Phi, \ \ \ \ \Id_{\mu_0} \circ \Phi := \Phi, \ \ \ \ \Xi \circ (\Psi \circ \Phi) 
:= (\Xi \circ \Psi) \circ \Phi. \]

\begin{definition}\label{def: setofdependops}
If $I$ is a set and $\lambda_0 : I \sto \D V_0$, let 
$\D A(I, \lambda_0)$\index{$\D A(I, \lambda_0)$} be the totality of dependent operations over $\lambda_0$,
equipped with the 
canonical equality:
\[ \Phi =_{\D A(I, \lambda_0)} \Psi :\TOT \forall_{i \in I}\big(\Phi_i =_{\lambda_0(i)} \Psi_i\big). \]
The totality $\D A(I, \lambda_0)$ is considered to be a set. If $\neq_{\lambda_0(i)}$ is an inequality on $\lambda_0(i)$,
for every $i \in I$, the canonical inequality $\neq_{\D A(I, \lambda_0)}$ on $\D A(I, \lambda_0)$ is defined by
$\Phi \neq_{\D A(I, \lambda_0)} \Psi :\TOT \exists_{i \in I}\big(\Phi_i \neq_{\lambda_0(i)} \Psi_i\big)$.
\end{definition}

Clearly, $\Phi =_{\D A(I, \lambda_0)} \Psi$ is an equivalence relation, and 
$\Phi \neq_{\D A(I, \lambda_0)} \Psi$ is an inequality relation. If $i \in I$, the
$i$-\textit{projection map}\index{$i$-projection map} on $\D A(I, \lambda_0)$ is the operation
$\pr_i^{\lambda_0} \colon \D A(I, \lambda_0) \sto \lambda_0(i)$, defined by $\pr_i^{\lambda_0}(\Phi) := \Phi_i$,
for every $i \in I$.
The operation $\pr_i^{\lambda_0}$ is a function. If $\Phi \colon \bigcurlywedge_{i \in I}\D F(\lambda_0(i),
\mu_0(i))$, a \textit{modulus of 
 surjectivity} for $\Phi$\index{modulus of surjectivity for a dependent operation} is a 
 dependent operation  $\Psi \colon \bigcurlywedge_{i \in I}\D F(\mu_0(i), \lambda_0(i)$ such that 
 $\Phi \circ \Psi =_{\D A(I,  \D F(\lambda_0, \lambda_0)} \Id_{\lambda_0}$. In this case, $\Psi_i$ is a modulus of 
 surjectivity for $\Phi$, for  every $i \in I$.
If $\fXY$, let $\fib^f \colon Y \sto \D V_0$ be defined by $y \mapsto \fib^f(y)$, for every $y \in Y$. If $f$ is 
contractible,
then by Definition~\ref{def: surjective} every fiber $\fib^f(y)$ of $f$ is contractible. A \textit{modulus of centres of 
contraction}\index{modulus of centres of contraction} for a contractible function $f$ is a dependent operation
$\centr^f \colon \bigcurlywedge_{y \in Y}\fib^f(y)$,
such that $\centr^f_y := \centr^f(y)$ is a centre of contraction for $f$.

\section{Subsets}
\label{sec: subsets}

\begin{definition}\label{def: subset}
Let $X$ be a set. A subset\index{subset} of $X$ is a pair $(A, i_A^X)$\index{$(A, i_A^X)$}, where $A$ is a set and 
$\i_A^X \colon A \hookrightarrow X$ is an embedding of $A$ into $X$.
If $(A, i_A^X)$ and $(B, i_B^X)$ are subsets of $X$, then $A$ is a \textit{subset}\index{subset} of $B$, in symbols 
$(A, i_A^X) \subseteq (B, i_B^X)$\index{$(A, i_A^X) \subseteq (B, i_B^X)$}, or simpler 
$A \subseteq B$\index{$A \subseteq B$},
if there is $f \colon A \to B$ such that the following diagram commutes
\begin{center}
\begin{tikzpicture}

\node (E) at (0,0) {$A$};
\node[right=of E] (B) {};
\node[right=of B] (F) {$B$};
\node[below=of B] (A) {$X$.};

\draw[->] (E)--(F) node [midway,above] {$f$};
\draw[right hook->] (E)--(A) node [midway,left] {$i_A^X \ $};
\draw[left hook->] (F)--(A) node [midway,right] {$\ i_B^X$};

\end{tikzpicture}
\end{center} 
In this case we use the notation $f \colon A \subseteq B$. Usually we write $A$ instead of $(A, i_A^X)$.
The totality of the subsets of $X$ is the \textit{powerset}\index{powerset} $\C P(X)$\index{$\C P(X)$} of $X$,
and it is equipped with the equality
\[ (A, i_A^X) =_{\C P(X)} (B, i_B^X) :\TOT A \subseteq B \ \& \ B \subseteq A. \]
If $f \colon A \subseteq B$ and $g \colon B \subseteq A$, we write $(f, g) \colon A =_{\C P(X)} B$.
\end{definition}

Since the membership condition for $\C P(X)$ requires quantification over $\D V_0$, the totality
$\C P(X)$ is a class. Clearly, $(X, \id_X) \subseteq X$. If $X_P$ is an extensional subset of $X$ (see Definition~\ref{def: 
extsubset}), then $(X_P, i_P^X) \subseteq X$, where $i_P^X \colon X_P \sto X$ is defined by $i_P^X(x) := x$, for every
$x \in X_P$.


\begin{proposition}\label{prp: subset1}
If $A, B \subseteq X$, and $f, g : A \subseteq B$, then $f$ is an embedding,
and $f =_{\D F(A, B)} h$
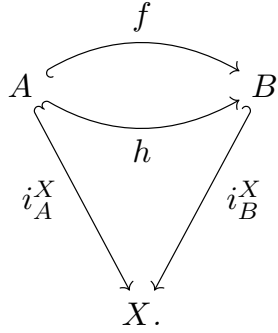
\begin{center}
\resizebox{4cm}{!}{%
\begin{tikzpicture}

\node (E) at (0,0) {$A$};
\node[right=of E] (B) {};
\node[right=of B] (F) {$B$};
\node[below=of B] (C) {};
\node[below=of C] (A) {$X$.};

\draw[left hook->,bend left] (E) to node [midway,above] {$f$} (F);
\draw[left hook->,bend right] (E) to node [midway,below] {$h$} (F);
\draw[right hook->] (E)--(A) node [midway,left] {$i_A^X \ $};
\draw[left hook->] (F)--(A) node [midway,right] {$ \ i_B^X$};

\end{tikzpicture}
}
\end{center} 
\end{proposition}

\begin{proof}
If $a, a{'} \in A$ such that $f(a) =_B f(a{'})$, then 
$i_B^X(f(a)) =_X i_B^X(f(a{'})) \TOT i_A^X(a) =_X i_A^X (a{'})$, which implies $a =_A a{'}$. Moreover, if
$i_B^X(f(a)) =_X i_A^X(a) =_X i_B^X(h(a))$, then $f(a) = h(a)$.
\end{proof}

The ``internal'' equality of subsets implies their ``external'' equality as sets  i.e., 
$(f, g) : A =_{\C P(X)} B \To (f, g) : A =_{\D V_0} B$. If
$a \in A$, then $i_A^X(g(f(a))) =_X i_B^X(f(a)) = i_A^X (a)$, hence $g(f(a)) =_A a$, and then
$g \circ f =_{\D F(A, A)} \id_A$. Similarly we get $f \circ g =_{\D F(B, B)} \id_B$.
Let the set
\[ \Eq(A, B) := \big\{(f, g) \in \D F(A, B) \times \D F(B, A) \mid f \colon A \subseteq B \ \& \ g 
\colon B \subseteq A\big\}, \]
equipped with the canonical equality of pairs as in the case of $\Eq(X, Y)$. 
Because of the Proposition~\ref{prp: subset1}, 
the set $\Eq(A, B)$ is a subsingleton i.e., 
\[ (f, g) \colon A =_{\C P(X)} B \ \& \ (f{'}, g{'}) \colon A =_{\C P(X)} B \To (f, g) = (f{'}, g{'}). \]
If $f \in \D F(A, B), g \in \D F(B, A), h \in \D F(B, C)$, and $k \in \D F(C, B)$, let 
$\refl(A) := \big(\id_A, \id_A\big)$ and $(f, g)^{-1} := (g, f)$, and $(f, g) \ast (h, k) := (h \circ f, g \circ k)$,
and the properties (i)-(iv) for $\Eq(A, B)$ hold by the equality of all their elements.

\begin{corollary}\label{cor: corapartness2}
Let the set $(X, =_X, \neq_X)$ and $\big(A, =_A, i_A^X, \neq_{A}^{\mathsmaller{i_A^X}}\big) \subseteq X$,
where the canonical inequality\index{canonical inequality on a subset} $\neq_{A}^{\mathsmaller{i_A^X}}$ on $A$ is
given by 
$a \neq_{A}^{\mathsmaller{i_A^X}} a{'} :\TOT i_A^X(a) \neq_X i_A^X(a{'})$, for every $a, a{'} \in A$.
If $(X, =_X, \neq_X)$ is discrete, then $\big(A, =_A, i_A^X, \neq_{A}^{\mathsmaller{i_A^X}}\big)$ is discrete,
and if $\neq_X$ is tight, $\neq_{A}^{\mathsmaller{i_A^X}}$ is tight.
\end{corollary}

\begin{proof} 
Since $i_A^X$ is an embedding, it follows immediately from Remark~\ref{rem: apartness2}. 
\end{proof}

\begin{remark}\label{rem: subset2}
If $P, Q$ are extensional properties on the set $X$, then
\[ X_P =_{\C P (X)} X_Q \TOT \forall_{x \in X}\big(P(x) \TOT Q(x)\big). \]
\end{remark}

\begin{proof}
 The implication $(\oT)$ is immediate to show, since the corresponding identity maps witness 
 the equality $X_P =_{\C P (X)} X_Q$.
 For the converse implication, let $(f, g) : X_P =_{\C P (X)} X_Q$. Let $x \in X$ such that $P(x)$. By
 the commutativity of the  following outer diagram 
 \begin{center}
 \resizebox{4cm}{!}
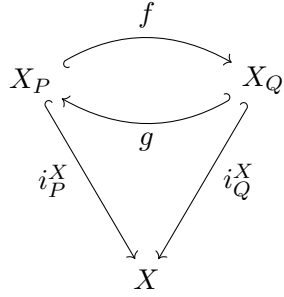
{%
\begin{tikzpicture}

\node (E) at (0,0) {$X_P$};
\node[right=of E] (B) {};
\node[right=of B] (F) {$X_Q$};
\node[below=of B] (C) {};
\node[below=of C] (A) {$X$};

\draw[left hook->,bend left] (E) to node [midway,above] {$f$} (F);
\draw[left hook->,bend left] (F) to node [midway,below] {$g$} (E);
\draw[right hook->] (E)--(A) node [midway,left] {$i_P^X \ $};
\draw[left hook->] (F)--(A) node [midway,right] {$ \ i_Q^X$};

\end{tikzpicture}
}
\end{center} 
we get $ f(x) := i_Q^X(f(x)) =_X  i_P^X(x) := x$, and by the extensionality of
 $Q$ and the fact that $Q(f(x))$ holds we get $Q(x)$. By the commutativity of the above inner diagram and the 
 extensionality of $P$ we get similarly the inverse implication.
\end{proof}

\begin{definition}\label{def: unionofsubsets}
If $(A, i_A^X), (B,i_B^X) \subseteq X$, their \textit{union}\index{union of subsets} 
$A \cup B$\index{$A \cup B$}
is the totality defined by
\[ z \in A \cup B : \TOT z \in A \ \vee \ z \in B, \]
equipped with the non-dependent assignment routine\footnote{Here we define a non-dependent assignment routine 
on the totality $A \cup B$, without knowing beforehand that $A \cup B$ is a set. It turns out that $A \cup B$
is set, but for that we need to define $i_{A \cup B}^X$ first.}
$i_{A \cup B}^X : A \cup B \sto X$, defined by 
\[ i_{A \cup B}^X(z) := \left\{ \begin{array}{lll}
                 i_A^X (z)    &\mbox{, $z \in A$}\\
                 {} \\
                 i_B^X (z)             &\mbox{, $z \in B$.}
                 \end{array}
          \right.\]           
If $z, w \in A \cup B$, we define $z =_{A \cup B} w :\Leftrightarrow i_{A \cup B}^X(z) =_X i_{A \cup B}^X(w).$
\end{definition}

Clearly, $=_{A \cup B}$ is an equality on $A \cup B$, which is considered to be a set,  
$i_{A \cup B}^X$ is an embedding of $A \cup B$ into $X$, and the pair 
$\big(A \cup B, i_{A \cup B}^X\big)$ is a subset of $X$. Note that if $P, Q$ are extensional properties on $X$, then 
$X_P \cup X_Q := X_{P \vee Q},$ since 
$z \in X_{P \vee Q} : \TOT (P \vee Q)(z) : \TOT P(z) \ \mbox{or} \ Q(z) \TOT : z \in X_P \cup X_Q,$
and the inclusion map $i : X_P \cup X_Q \eto X$ is the identity, as it is for $X_{P \vee Q}$ 
(see Definition~\ref{def: totalities}). If $\neq_X$ is a given inequality on $X$, the canonical inequality 
on $A \cup B$ is determined in Corollary~\ref{cor: 
corapartness2}.

\begin{definition}\label{def: intofsubsets}
If $(A, i_A^X), (B,i_B^X) \subseteq X$, their \textit{intersection}\index{intersection of subsets} 
$A \cap B$\index{$A \cap B$} is the totality defined by separation on $A \times B$ as follows:
\[ A \cap B := \{(a, b) \in A \times B \mid i_A^X (a) =_X i_B^X (b)\}. \]
Let the non-dependent assignment routine $i_{A \cap B}^X : A \cap B \sto X$, defined by 
$i_{A \cap B}^X(a, b) := i_A^X (a)$, for every $(a, b) \in A \cap B$.
If $(a, b)$ and $(a{'}, b{'})$ are in $A\cap B$, let
\[ (a, b) =_{A \cap B} (a{'}, b{'}) :\TOT i_{A \cap B}^X(a, b) =_X i_{A \cap B}^X(a{'}, b{'}) :\TOT 
i_A^X (a) =_X i_A^X (a{'}). \]
We write $A \between B$\index{$A \between B$} to denote that the intersection $A \cap B$ is inhabited.

\end{definition}

Clearly, $=_{A \cap B}$ is an equality on $A \cap B$, which is considered to be a set,  
$i_{A \cap B}^X$ is an embedding of $A \cap B$ into $X$, and 
$\big(A \cap B, i_{A \cap B}^X\big)$ is a subset of $X$.
If $\neq_X$ is a given inequality on $X$, the canonical inequality on $A \cap B$ is determined in Corollary~\ref{cor: 
corapartness2}.
If $P, Q$ are extensional properties on $X$, then $X_P \cap X_Q$ has elements in
$X \times X$, while $X_{P \wedge Q}$ has elements in $X$, hence the two subsets are not definitionally equal. 
Next we show that they are ``externally'' equal i.e., equal in $\D V_0$.

\begin{remark}\label{rem: extintersectionofsubsets}
If $P, Q$ are extensional properties on the set $X$, then $X_{P \wedge Q} =_{\D V_0} X_P \cap X_Q$.
\end{remark}

\begin{proof}
Since the inclusion maps corresponding to $X_P$ and $X_Q$ are the identities, let $f : X_{P \wedge Q}
\to X_P \cap X_Q$ with $f(z) := (z, z)$, for every $z \in X_{P \wedge Q}$, and let $g : X_P \cap X_Q
\to X_{P \wedge Q}$ with $g(a, b) := a$, for every $(a, b) \in X_P \cap X_Q$.  Hence, $f(g(a, b)) 
:= f(a) := (a, a)$, and since $(a, b) \in X_P \cap X_Q$, we have by definition that $P(a), Q(b)$ 
and $a =_X b$, hence $(a, a) =_{X \times X} (a, b)$. If $z \in X_{P \wedge Q}$, then $g(f(z)) := g(z, z) := z$. 
\end{proof}

Clearly, $X \cap X =_{\C P(X)} X$, while $\pr_A \colon (A \cap B, i_{A \cap B}^X) \subseteq (A, i_A)$ and the identity map 
$e_A \colon A \to A \cup B$ witnesses the inequality $(A, i_A^X) \subseteq (A \cup B, i_{A \cup B}^X)$
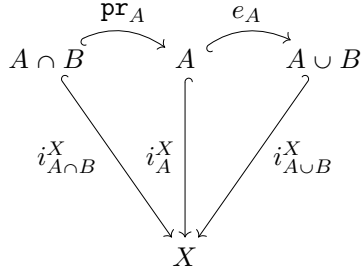
\begin{center}
 \resizebox{5cm}{!}{%
\begin{tikzpicture}

\node (E) at (0,0) {$A \cap B$};
\node[right=of E] (F) {$A$};
\node[below=of F] (C) {};
\node[below=of C] (A) {$X$};
\node[right=of F] (K) {$A \cup B$};

\draw[left hook->,bend left] (E) to node [midway,above] {$\pr_A$} (F);
\draw[right hook->] (F) to node [midway,left] {$i_A^X$} (A);
\draw[right hook->] (E)--(A) node [midway,left] {$i_{A \cap B}^X \ $};
\draw[left hook->] (K)--(A) node [midway,right] {$ \ i_{A \cup B}^X$};
\draw[left hook->,bend left] (F) to node [midway,above] {$e_A$} (K);

\end{tikzpicture}
}
\end{center} 

The following properties of the union and intersection of subsets are easy to show.

\begin{proposition}\label{prp: subset3}
Let $A, B$ and $C$ be subsets of the set $X$.\\[1mm]
\normalfont (i) 
\itshape $A \cup B =_{\C P(X)} B \cup A$ and $A \cap B =_{\C P(X)} B \cap A$.\\[1mm]
\normalfont (ii) 
\itshape $A \cup (B \cup C) =_{\C P(X)} (A \cup B) \cup C$ and $A \cap (B \cap C) =_{\C P(X)} (A \cap B) \cap C$.\\[1mm]
\normalfont (iii) 
\itshape $A \cap (B \cup C) =_{\C P(X)} (A \cap B) \cup (A \cap C)$ and $A \cup (B \cap C) =_{\C P(X)}
(A \cup B) \cap (A \cup C)$.
\end{proposition}

\begin{definition}\label{def: image}
Let $X, Y$ be sets, $(A, i_A^X) (C, i_C^X) \subseteq X$, $e \colon (A, i_A^X) \subseteq (C, i_C^X)$, $f \colon C \to Y$,
and $(B, i_B^Y) \subseteq Y$. The \textit{restriction}\index{restriction} $f_{|_{A}}$\index{$f_{|_{A}}$} of $f$ 
to $A$ is the function 
$f_{_{A}} := f \circ e$
\begin{center}
\begin{tikzpicture}

\node (E) at (0,0) {$A$};
\node[right=of E] (F) {$C$};
\node[right=of F] (A) {$Y$.};

\draw[right hook->] (E)--(F) node [midway,above] {$e$};
\draw[->] (F)--(A) node [midway,above] {$f$};
\draw[->,bend right] (E) to node [midway,below] {$f_{|_{A}}$} (A);

\end{tikzpicture}
\end{center}
The \textit{image} \index{image} $f(A)$\index{$f(A)$} of $A$ under $f$ is the pair $f(A) := (A, f_{_{A}}),$
where $A$ is equipped with the equality
$a =_{f(A)} a{'} : \TOT f_{|_{A}}(a) =_Y f_{|_{A}}(a{'}),$
for every $a, a{'} \in A$. We denote
$\{f(a) \mid a \in A\} := f(A)$\index{$\{f(a) \mid a \in A\}$}.
The \textit{pre-image} \index{pre-image} $f^{-1}(B)$\index{$f^{-1}(B)$} of $B$ under $f$ is the set
\[ f^{-1}(B) := \{(c, b) \in C \times B \mid f(c) =_Y i_B^Y (b)\}. \]
Let $i_{\mathsmaller{f^{-1}(B)}}^C \colon f^{-1}(B) \eto C$, defined by
$i_{\mathsmaller{f^{-1}(B)}}^C(c, b) := c$, for every $(c,b) \in f^{-1}(B)$. 
The equality of the extensional subset $f^{-1}(B)$ of $C \times B$ is
inherited from the equality of $C \times B$. 

\end{definition}

Clearly, the restriction $f_{|_{A}}$ of $\fXY$ to $(A, \i_A) \subseteq X$ is the function 
$f_{_{A}} := f \circ i_A^X$. It is immediate to show that $f(A) \subseteq Y$ and $f^{-1}(B) \subseteq C$.
Notice that 
\[ (A, i_A^X) =_{\C P(X)} (B, i_B^X) \To i_A^X (A) =_{\C P(X)} i_B^X (B), \]
since, if $(f, g) : (A, i_A^X) =_{\C P(X)} (B, i_B^X)$, then $i_A^X (a) =_X i_B^X (f(a))$ and $i_B^X (b) 
=_X i_A^X (g(b))$, for every $a \in A$ and $b \in B$, respectively. 
If $\neq_Y$ is a given inequality on $Y$, the canonical inequality on $f(A)$ is determined in 
Corollary~\ref{cor: corapartness2}. Similarly, if $\neq_X$ is an inequality on $X$, $f \colon X \to Y$,
and $(B, i_B^Y) \subseteq Y$, the canonical inequality on $f^{-1}(B)$ is given by  
$(x, b) \neq_{f^{-1}(B)} (x{'}, b{'}) :\TOT x \neq_X x{'}$, and not by the canonical inequality on $X \times B$.


\begin{proposition}\label{prp: subset4}
Let $X, Y$ be sets, $A, B$ subsets of $X$, $C, D$ subsets of $Y$, and $f : X \to Y$.\\[1mm]
\normalfont (i) 
\itshape $f^{-1}(C \cup D) =_{\C P(X)} f^{-1}(C) \cup f^{-1}(D)$.\\[1mm]
\normalfont (ii) 
\itshape $f^{-1}(C \cap D) =_{\C P(X)} f^{-1}(C) \cap f^{-1}(D)$.\\[1mm]
\normalfont (iii) 
\itshape $f(A \cup B) =_{\C P(Y)} f(A) \cup f(B)$.\\[1mm]
\normalfont (iv) 
\itshape $f(A \cap B) =_{\C P(Y)} f(A) \cap f(B)$.\\[1mm]
\normalfont (v) 
\itshape $A \subseteq f^{-1}(f(A))$.\\[1mm]
\normalfont (vi) 
\itshape $f(f^{-1}(C) \cap A) =_{\C P(Y)} C \cap f(A)$, and $f(f^{-1}(C)) =_{\C P(Y)} C \cap f(X)$.
\end{proposition}

\begin{proposition}\label{prp: subset15}
Let 
$(A, i_A^X), (B, i_B^X), (A{'}, i_{A{'}}^X) , (B, i_{B{'}}^X) \subseteq X$, such that 
$A =_{\C P(X)} A{'}$ and $B =_{\C P(X)} B{'}$. Let also $(C, i_C^Y), (C{'}, i_{C{'}}^Y), (D, i_D^Y) \subseteq Y$, such that 
$C =_{\C P(Y)} C{'}$, and let $\fXY$. \\[1mm]
\normalfont (i) 
\itshape  $A \cap B =_{\C P(X)} A{'} \cap B{'}$, and $A \cup B =_{\C P(X)} A{'} \cup B{'}$.\\[1mm]
\normalfont (ii) 
\itshape  $f(A) =_{\C P(Y)} f(A{'})$, and $f^{-1}(C) =_{\C P(X)} f^{-1}(C{'})$.\\[1mm]
\normalfont (iii) 
\itshape  $(A \times C, i_A^X \times i_C^Y) \subseteq X \times Y$,
where the map $i_A^X \times i_C^Y \colon A \times C \hookrightarrow X \times Y$ is defined by 
\[ (i_A^X \times i_C^Y)(a, c) := \big(i_A^X(a), i_C^Y(c)\big); \ \ \ \ (a, c) \in A \times C. \]
\normalfont (iv) 
\itshape  $A \times C =_{\C P(X \times Y)} A{'} \times C{'}$.\\[1mm]
\normalfont (v) 
\itshape  $A \times (C \cup D) =_{\C P(X \times Y)} (A \times C) \cup (A \times D)$.\\[1mm]
\normalfont (vi) 
\itshape  $A \times (C \cap D) =_{\C P(X \times Y)} (A \times C) \cap (A \cap D)$.

\end{proposition}
 
 \begin{proof}
All cases are straightforward to show.
 \end{proof}

\section{Partial functions}
\label{sec: partial}

\begin{definition}\label{def: partialfunction}
Let $X, Y$ be sets. A partial function\index{partial function} from $X$ to $Y$ is a triplet
$(A, i_A^X, f_A^Y)$\index{$(A, i_A^X, f_A^Y)$}, where $(A, i_A^X) \subseteq X$, and $f_A^Y \in \D F(A, Y)$.
Often, we use only the symbol $f_A^Y$ instead of the triplet $(A, i_A^X, f_A^Y)$, and we also write
$f_A^Y \colon X \pto Y$. If $(A, i_A^X, f_A^Y)$ and $(B, i_B^X, f_B^Y)$ are partial functions from $X$ to $Y$, 
we call $f_A^Y$ a \textit{subfunction}\index{sufunction} of $f_B^Y$, in symbols 
$(A, i_A^X, f_A^Y) \leq (B, i_B^X, f_B^Y)$\index{$(A, i_A^X, f_A^Y) \leq (B, i_B^X, f_B^Y)$}, or simpler 
$f_A^Y \leq f_B^Y$\index{$f_A^Y \leq f_B^Y$},
if there is $e_{AB} \colon A \to B$ such that the following inner diagrams commute
\begin{center}
\begin{tikzpicture}

\node (E) at (0,0) {$A$};
\node[right=of E] (B) {};
\node[right=of B] (F) {$B$};
\node[below=of B] (A) {$X$};
\node[below=of A] (C) {$ \  Y$.};

\draw[->] (E)--(F) node [midway,above] {$e_{AB}$};
\draw[->,bend right=40] (E) to node [midway,left] {$f_A^Y$} (C);
\draw[->,bend left=40] (F) to node [midway,right] {$f_B^Y$} (C);
\draw[right hook->] (E)--(A) node [midway,left] {$i_A^X \ $};
\draw[left hook->] (F)--(A) node [midway,right] {$\ i_B^X  $};

\end{tikzpicture}
\end{center} 
In this case we use the notation $e_{AB} \colon f_A^Y \leq f_B^Y$. 
The totality of partial functions from $X$ to $Y$ is the \textit{partial function space}
\index{partial function space}$\C F(X,Y)$\index{$\C F(X,Y)$},
and it is equipped with the equality
\[ (A, i_A^X, f_A^Y) =_{\C F(X,Y)} (B, i_B^X, f_B^Y) :\TOT f_A^Y \leq f_B^Y \ \& \ f_B^Y \leq f_A^Y. \]
If $e_{AB} \colon f_A^Y \leq f_B^Y$ and $e_{BA} \colon f_B^Y \leq f_A^Y$, we write $(e_{AB}, e_{BA}) : f_A^Y 
=_{\C F(X,Y)} f_B^Y$.
\end{definition}

Since the membership condition for $\C F(X,Y)$ requires quantification over $\D V_0$, the totality
$\C F(X,Y)$ is a class. Clearly, if $\fXY$, then $(X, \id_X, f) \in \C F(X,Y)$. 
If $(e_{AB}, e_{BA}) : f_A^Y =_{\C F(X,Y)} f_B^Y$, then $(e_{AB}, e_{BA}) : A =_{\C P(X)} B$, and 
$(e_{AB}, e_{BA}) : A =_{\D V_0} B$. 
Let the set
\[ \Eq(f_A^Y, f_B^Y) := \big\{(f, g) \in \D F(A, B) \times \D F(B, A) \mid f \colon f_A^Y \leq f_B^Y 
\ \& \ g \colon f_B^Y \leq f_A^Y\big\}, \]
equipped with the canonical equality of the product. 
All the elements of $\Eq(f_A^Y, f_B^Y)$ are equal to each other.
If  $f \in \D F(A, B), g \in \D F(B, A), h \in \D F(B, C)$, and $k \in \D F(C, B)$, let 
\[ \refl(f_A^Y) := \big(\id_A, \id_A\big) \ \ \& \ \ (f, g)^{-1} := (g, f) \ \ \& \ \ 
(f, g) \ast (h, k) := (h \circ f, g \circ k), \]
and the groupoid-properties for $\Eq(f_A^Y, f_B^Y)$ hold by the equality of its elements.

\begin{proposition}\label{prp: compositionpf}
Let $(A, i_A^X, f_A^Y) \in \C F(X, Y)$ and $(B, i_B^Y, g_B^Z) \in \C F(Y, Z)$. Their composition\index{composition
of partial functions}\index{$g_B^Z \pcirc f_A^Y$} 
\[ g_B^Z \mathsmaller{\pcirc} f_A^Y := \bigg(\big(f_A^Y\big)^{-1}(B), \ i_A^X \circ e_{\mathsmaller{(f_A^Y)^{-1}(B)}}^A, 
\ \big(g_B^Z \circ f_A^Y\big)^Z_{\mathsmaller{(f_A^Y)^{-1}(B)}}\bigg), \ \ \ \mbox{where}
 \]
\[ \big(f_A^Y\big)^{-1}(B) := \big\{(a, b) \in A \times B \mid f_A^Y(a) =_Y i_B^Y(b)\big\}, \]
\[  e_{\mathsmaller{(f_A^Y)^{-1}(B)}}^A \colon \big(f_A^Y\big)^{-1}(B) \eto A, \ \ \ \ (a, b) \mapsto a; \ \ \ \ (a,b) 
 \in \big(f_A^Y\big)^{-1}(B), \]
\[ \big(g_B^Z \circ f_A^Y\big)^Z_{\mathsmaller{(f_A^Y)^{-1}(B)}}(a, b) := g_B^Z(b); \ \ \ \  (a,b) 
 \in \big(f_A^Y\big)^{-1}(B), \]
is a partial function that belongs to $\C F(X, Z)$. If $(A, i_A^X, i_A^X) \in \C F(X, X), (B, i_B^Y, i_B^Y) \in \C F(Y, Y)$, and
$(C, i_C^Z, h_C^W) \in \C F(Z, W)$, the following properties hold:\\[1mm]
\normalfont (i) 
\itshape $f_A^Y \pcirc i_A^X =_{\C F(X,Y)} f_A^Y$ and $i_B^Y \pcirc f_A^Y =_{\C F(X,Y)} f_A^Y$.\\[1mm]
\normalfont (ii) 
\itshape $\big(h_C^W \pcirc g_B^Z\big) \pcirc f_A^X =_{\C F(X,Z)} h_C^W \pcirc \big(g_B^Z \pcirc f_A^X\big)$.
\end{proposition}

\begin{proof}
(i) We show only the first equality and for the second we work similarly. By definition
\[ f_A^Y \mathsmaller{\pcirc} i_A^X := \bigg(\big(i_A^X\big)^{-1}(A), \ i_A^X \circ e_{\mathsmaller{(i_A^X)^{-1}(A)}}^A, 
\ \big(f_A^Y \circ i_A^X\big)^Y_{\mathsmaller{(i_A^X)^{-1}(A)}}\bigg), \ \ \ \mbox{where}
 \]
\[ \big(i_A^X\big)^{-1}(A) := \big\{(a, a{'}) \in A \times A \mid i_A^X(a) =_X i_A^X(a{'})\big\}, \]
\[  e_{\mathsmaller{(i_A^X)^{-1}(A)}}^A \colon \big(i_A^X\big)^{-1}(A) \eto A, \ \ \ \ (a, a{'}) \mapsto a; \ \ \ \ (a,a{'}) 
 \in \big(i_A^X\big)^{-1}(A), \]
\[ \big(f_A^Y \circ i_A^X\big)(a, a{'}) := f_A^Y(a{'}); \ \ \ \  (a,a{'}) \in \big(i_A^X\big)^{-1}(A). \]
Let the operations $\phi \colon A \sto \big(i_A^X\big)^{-1}(A)$, defined by $\phi(a) := (a, a)$, for every $a \in A$, 
and $\theta \colon \big(i_A^X\big)^{-1}(A) \sto A$, defined by $\theta(a, a{'}) := a$, for every $(a, a{'}) \in
\big(i_A^X\big)^{-1}(A)$. It is immediate to show that $\phi$ and $\theta$ are well-defined functions. It is 
straightforward to show the commutativity of the following inner diagrams 
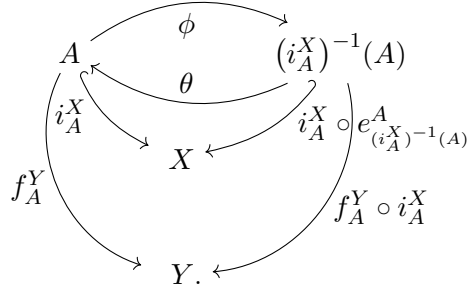
\begin{center}
\begin{tikzpicture}

\node (E) at (0,0) {$A$};
\node[right=of E] (B) {};
\node[right=of B] (F) {$\big(i_A^X\big)^{-1}(A)$};
\node[below=of B] (A) {$ \ X$};
\node[below=of A] (C) {$ \ \ Y$.};

\draw[->,bend left] (E) to node [midway,below] {$\phi$} (F);
\draw[->,bend left] (F) to node [midway,above] {$\theta$} (E);
\draw[->,bend right=50] (E) to node [midway,left] {$f_A^Y$} (C);
\draw[->,bend left=50] (F) to node [midway,right] {$f_A^Y \circ i_A^X$} (C);
\draw[right hook->,bend right=20] (E) to node [midway,left] {$ \ i_A^X \ $} (A);
\draw[left hook->,bend left=20] (F) to node [midway,right] {$\ \ i_A^X \circ e_{\mathsmaller{(i_A^X)^{-1}(A)}}^A \ $} (A);

\end{tikzpicture}
\end{center} 
(ii) We have that 
$ h_C^W \mathsmaller{\pcirc} g_B^Z := \bigg(\big(g_B^Z\big)^{-1}(C), \ i_B^Y \circ e_{\mathsmaller{(g_B^Z)^{-1}(C)}}^B, 
\ \big(h_C^W \pcirc g_B^Z\big)^W_{\mathsmaller{(g_B^Z)^{-1}(C)}}\bigg)$, where
\[ \big(g_B^Z\big)^{-1}(C) := \big\{(b, c) \in B \times C \mid g_B^Z(b) =_Z i_C^Z(c)\big\}, \]
\[  e_{\mathsmaller{(g_B^Z)^{-1}(C)}}^B \colon \big(g_B^Z\big)^{-1}(C) \eto B, \ \ \ \ (b, c) \mapsto b; \ \ \ \ (b,c) 
 \in \big(g_B^Z\big)^{-1}(C), \]
\[ \big(h_C^W \circ g_B^Z\big)(b, c) := h_C^W(c); \ \ \ \  (b,c) \in \big(g_B^Z\big)^{-1}(C). \] 
Hence, 
$\big(h_C^W \mathsmaller{\pcirc} g_B^Z\big) \mathsmaller{\pcirc} f_A^X := \bigg(D, i_A^X \circ e_D^A, 
\big[\big(h_C^W \circ g_B^Z\big) \circ f_Z^X\big]_D^W\bigg)$, where
\[ 
 D :=  \big(f_A^Y\big)^{-1}\big[\big(g_B^Z\big)^{-1}(C)\big] 
  := \bigg\{(a, d) \in A \times \big[\big(g_B^Z\big)^{-1}(C)\big] \mid f_A^Y(a) =_Y \big(i_B^Y \circ
  e_{\mathsmaller{(g_B^Z)^{-1}(C)}}^B\big)(d) \bigg\},
 \]
with $d:= (b,c) \in B \times C$ such that $g_B^Z(b) =_Z i_C^Z(c)$. The map $e_D^A \colon D \eto A$ is defined by 
the rule $(a, d) \mapsto a$, for every $(a,d) \in D$, and
\[ \big[\big(h_C^W \circ g_B^Z\big) \circ f_A^Y\big](a, d) := \big(h_C^W \circ g_B^Z\big)(d) := h_C^W(c); \ \ \ \ 
(a,d) := (a, (b,c)) \in D. \]
Moreover, 
$h_C^W \mathsmaller{\pcirc} \big(g_B^Z \mathsmaller{\pcirc} f_A^Y\big) := \bigg(E, \ 
i_A^X \circ e_{(f_A^Y)^{-1}(B)}^A \circ e_E^{\mathsmaller{(f_A^Y)^{-1}(B)}}, \  
\big[h_C^W \circ \big(g_B^Z \circ f_Z^X\big)\big]_E^W\bigg)$, where
\[ E :=  \bigg[\big(g_B^Z \circ f_A^Y\big)_{\mathsmaller{(f_A^Y)^{-1}(B)}}^Z\bigg]^{-1}(C) 
 := \bigg\{(u, c) \in \big[\big(f_A^Y\big)^{-1}(B)\big] \times C \mid (g_B^Z \circ f_A^Y)(u) =_Z i_C^Z(c) \bigg\},
\]
$e_E^{\mathsmaller{(f_A^Y)^{-1}(B)}} \colon E \eto \big((f_A^Y\big))^{-1}(B)$ is defined by the rule
$(u, c) \mapsto u$, for every $(u,c) \in E$, and 
\[ \big[h_C^W \circ \big(g_B^Z \circ f_A^Y\big)\big](u, c) := h_C^W(c); \ \ \ \ 
(u,c) \in E. \]
Let the operations $\phi \colon D \sto E$, defined by $\phi(a, (b,c)) := ((a,b), c)$, for every $(a, (b,c)) \in D$, 
and $\theta \colon E \sto D$, defined by $\theta((a,b),c) := (a, (b,c))$, for every $((a,b),c) \in 
E$. It is straightforward to show that $\phi$ and $\theta$ are well-defined functions, and that the 
following inner diagrams commute
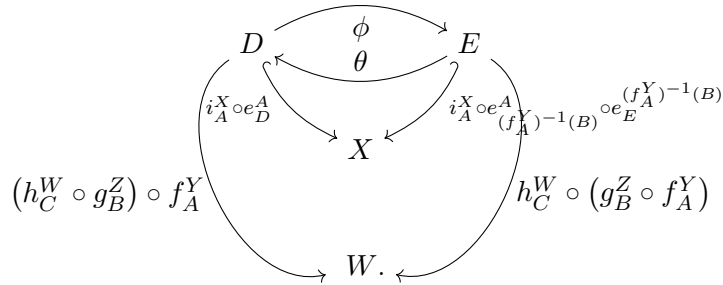
\begin{center}
\begin{tikzpicture}

\node (E) at (0,0) {$D$};
\node[right=of E] (B) {};
\node[right=of B] (F) {$E$};
\node[below=of B] (A) {$X$};
\node[below=of A] (C) {$ \ W$.};

\draw[->,bend left] (E) to node [midway,below] {$\phi$} (F);
\draw[->,bend left] (F) to node [midway,above] {$\theta$} (E);
\draw[->,bend right=80] (E) to node [midway,left] {$\big(h_C^W \circ g_B^Z\big) \circ f_A^Y$} (C);
\draw[->,bend left=80] (F) to node [midway,right] {$h_C^W \circ \big(g_B^Z \circ f_A^Y\big)$} (C);
\draw[right hook->,bend right=20] (E) to node [midway,left] {$ \ \mathsmaller{i_A^X \circ e_D^A} \ $} (A);
\draw[left hook->,bend left=20] (F) to node [midway,right] {$\ \mathsmaller{i_A^X \circ e_{(f_A^Y)^{-1}(B)}^A 
\circ e_E^{\mathsmaller{(f_A^Y)^{-1}(B)}}} \ $} (A);

\end{tikzpicture}
\end{center} 
\end{proof}

%
%
%
%

The next proposition is straightforward to show.

\begin{proposition}\label{prp: intunionpartial}
Let $(A, i_A^X, f_A^Y), (B, i_B^X, f_B^X) \in \C F(X, Y)$
\begin{center}
\begin{tikzpicture}

\node (E) at (0,0) {$A$};
\node[right=of E] (B) {$X$};
\node[right=of B] (F) {$B$};
\node[below=of B] (C) {$Y$.};

\draw[right hook->] (E)--(B) node [midway,above] {$i_A^X$};
\draw[left hook->] (F)--(B) node [midway,above] {$i_B^X$};
\draw[->] (E)--(C) node [midway,left] {$f_A^Y \ \ $};
\draw[->] (F)--(C) node [midway,right] {$ \ f_B^Y$};

\end{tikzpicture}
\end{center} 
Their left $f_A^Y \cap_l f_B^Y$ and right 
intersection\index{left intersection of partial functions}\index{$f_A^Y \cap_l f_B^Y$} 
$f_A^Y \cap_r f_B^Y$ are the partial functions\index{right intersection of partial functions}\index{$f_A^Y \cap_r f_B^Y$}
\[ f_A^Y \cap_l f_B^Y := \bigg(A \cap B, \ i_{A \cap B}^X, \ \big(f_A^Y \cap_l f_B^Y\big)_{A \cap B}^Y\bigg), \ \ \ \mbox{where}
 \]
\[ \big(f_A^Y \cap_l f_B^Y\big)_{A \cap B}^Y(a, b) := f_A^Y(a); \ \ \ \  (a,b) \in A \cap B, \ \ \ \ \mbox{and} \]
\[ f_A^Y \cap_r f_B^Y := \bigg(A \cap B, \ i_{A \cap B}^X, \ \big(f_A^Y \cap_r f_B^Y\big)_{A \cap B}^Y\bigg), \ \ \ \mbox{where}
 \]
\[ \big(f_A^Y \cap_r f_B^Y\big)_{A \cap B}^Y(a, b) := f_B^Y(b); \ \ \ \  (a,b) \in A \cap B. \]
Their union $f_A^Y \cup f_B^Y$\index{union of partial functions}\index{$f_A^Y \cup f_B^Y$} is the partial function
\[ f_A^Y \cup f_B^Y := \bigg(A \cup B, \ i_{A \cup B}^X, \ \big(f_A^Y \cup f_B^Y\big)_{A \cup B}^Y\bigg), \ \ \ \mbox{where}
 \]
\[ \big(f_A^Y \cup f_B^Y\big)_{A \cup B}^Y(z) := \left\{ \begin{array}{ll}
                 f_A^Y(z)   &\mbox{, $z \in A$}\\
                 f_B^Y(z)            &\mbox{, $z \in B$.}
                 \end{array}
                 \right.
\] 
\normalfont (i) 
\itshape $f_A^Y \cap_l f_B^Y \leq f_A^Y$ and $f_A^Y \cap_r f_B^Y \leq f_B^Y$.\\[1mm]
\normalfont (ii) 
\itshape If $f_A^Y(a) =_Y f_B^Y(b)$, for every $(a,b) \in A \cap B$, then $f_A^Y \cap_l f_B^Y =_{\C F(X,Y)}
f_A^Y \cap_r f_B^Y$.\\[1mm]
\normalfont (iii) 
\itshape $f_A^Y \leq f_A^Y \cup f_B^Y$ and $f_B^Y \leq f_A^Y \cup f_B^Y$.\\[1mm]
\normalfont (iv) 
\itshape  $f_A^Y \cup f_B^Y =_{\C F(X,Y)} f_B^Y \cup f_A^Y$.
\end{proposition}

\begin{definition}\label{def: multin2}
Let the operation of multiplication on $\D 2$, defined by $0 \cdot 1 := 1 \cdot 0 := 0 \cdot 0 := 0$ and $1 \cdot 1 := 1$.
If $(A, i_A^X, f_A^{\D 2}), (B, i_B^X, g_B^{\D 2}) \in \C F(X, \D 2)$, let 
\[ f_A \cdot g_B := \big(A \cap B, i_{A \cap B}^X, (f_A \cdot g_B)_{A \cap B}^{\D 2}\big), \]
where $(f_A \cdot g_B)_{A \cap B}^{\D 2} \colon A \cap B \to \D 2$ is defined, for every $(a, b) \in A \cap B$, 
by\index{$f_A \cdot g_B$}
\[ (f_A \cdot g_B)_{A \cap B}^{\D 2}(a, b) := f_A^{\D 2}(a) \cdot g_B^{\D 2}(b). \]
\end{definition}

By the equality of the product on $A \cap B$, it is immediate to show that the operation
$(f_A \cdot g_B)_{A \cap B}^{\D 2}$ is a function. More generally, operations on $Y$ induce
operations on $\C F(X, Y)$. The above example with $Y := \D 2$ is useful to the next section.

\section{Complemented subsets}
\label{sec: complemented}

An inequality on a set $X$ induces a positively defined notion of disjointness of subsets of $X$.

\begin{definition}\label{def: apartsubsets}
Let $(X, =_X, \neq_X)$ be a set, and $(A, i_A^X), (B, i_B^X) \subseteq X$. We say that 
$A$ and $B$ are disjoint with respect to $\neq_X$\index{disjoint subsets}, in symbols 
$A \Disj_{\mathsmaller{\neq_X}} B$\index{$A \Disj_{\mathsmaller{\neq}} B$}, if
\[ A \underset{\mathsmaller{\mathsmaller{\mathsmaller{\neq_X}}}} \Disj B : 
\TOT \forall_{a \in A}\forall_{b \in B}\big(i_A^X(a) \neq_X i_B^X(b) \big). \]
If $\neq_X$ is clear from the context, we only write $A \Disj B$\index{$A \Disj B$}.
\end{definition}

Clearly, if $A \Disj B$, then $A \cap B$ is not inhabited. 
The positive disjointness of subsets of $X$ induces the notion of a complemented subset of $X$, and
the negative notion of the complement of a set is avoided. We use bold letters to denote a complemented subset of a set.

\begin{definition}\label{def: complementedsubset}
A \textit{complemented subset}\index{complemented subset} of a set $(X, =_X, \neq_X)$ is a pair
$\B A := (A^1, A^0)$\index{$\B A := (A^1, A^0)$}, 
where $(A^1, i_{A^1}^X)$ and $(A^0, i_{A^0}^X)$ are subsets of $X$ such that $A^1 \Disj A^0$. We
call $A^1$ the $1$-component of $\B A$\index{$1$-component of a complemented subset} and $A^0$ the
$0$-component of $\B A$\index{$0$-component of a complemented subset}. If 
$\Dm(\B A) := A^1 \cup A^0$\index{$\Dm(\B A)$} is the domain of $\B A$\index{domain of a complemented subset},
the 
\textit{indicator function}\index{indicator function of a complemented subset}, or
\textit{characteristic function}\index{characteristic function of a complemented subset}, of
$\B A$ is the operation $\chi_{\B A} : \Dm(\B A) \sto \D 2$ defined by
\[ \chi_{\B A}(x) := \left\{ \begin{array}{ll}
                 1   &\mbox{, $x \in A^1$}\\
                 0             &\mbox{, $x \in A^0$.}
                 \end{array}
          \right. \] 
Let $x \in \B A :\TOT x \in A^1$ and $x \notin \B A :\TOT x \in A^0$. If $\B A, \B B$ are complemented subsets of $X$, let
\[ \B A \subseteq \B B : \TOT A^1 \subseteq B^1 \ \& \ B^0 \subseteq A^0.  \]
Let $\C P^{\Disj}(X)$\index{$\C P^{\Disj}(X)$}\index{$\B A \subseteq \B B$} be their totality, 
equipped with the equality 
$\B A =_{\C P^{\mathsmaller{\Disj}} (X)} \B B : \TOT \B A \subseteq \B B \ \& \ \B B \subseteq \B A$.
Let $\Eq(\B A, \B B) := \Eq(A^1, B^1) \times \Eq(A^0, B^0)$. A map $\B f \colon \B A \to \B B$ from $\B A$ 
to $\B B$ is a pair $(f^1, f^0)$, where $f^1 \colon A^1 \to B^1$ and
$f^0 \colon A^0 \to B^0$\index{map between complemented subsets}.

\end{definition}

Clearly, $ \B A =_{\C P^{\mathsmaller{\Disj}} (X)} \B B \TOT A^1 =_{\C P(X)} B^1 \ \& \ A^0 =_{\C P(X)} B^0$, 
and $\Eq(\B A, \B B)$ is a subsingleton, as the product of subsingletons.
Since the membership condition for  $\C P^{\mathsmaller{\Disj}} (X)$ requires quantification over $\D V_0$, the 
totality  $\C P^{\mathsmaller{\Disj}} (X)$ is a class. The operation
$\chi_{\B A}$ is a function, actually, $\chi_{\B A}$ is a partial function in $\C F(X, \D 2)$. 
Let $z, w \in A^1 \cup A^0$ such that $z =_{A^1 \cup A^0} w$ i.e., 
\[  \left.
\begin{array}{lll}
                 i_{A^1}^X(z)   &\mbox{, $z \in A^1$}\\
                 {}\\
                  i_{A^0}^X(z)             &\mbox{, $z \in A^0$}
                 \end{array}
          \right\}
     := i_{A^1 \cup A^0}^X(z) 
      =_X 
i_{A^1 \cup A^0}^X(w) := \left\{ \begin{array}{lll}
                  i_{A^1}^X(w)   &\mbox{, $w \in A^1$}\\
                 {}\\
                  i_{A^0}^X(w)           &\mbox{, $w \in A^0$.}
                 \end{array}
          \right.
\]
Let $z \in A^1$. If $w \in A^0$, then $i_{A^1}^X(z) := i_{A^1 \cup A^0}^X(z) =_X i_{A^1 \cup A^0}^X(w) := i_{A^0}^X(w)$
i.e., $(z, w) \in A^1 \cap A^0$, which contradicts the hypothesis $A^1 \Disj A^0$. Hence $w \in A^1$, and 
$\chi_{\B A}(z) = \chi_{\B A}(w)$. If $z \in A^0$, we proceed similarly.

\begin{definition}\label{def: detachable}
If $(X, =_X)$ is a set, let the inequality on $X$ defined by
\[ x \neq_X^{\mathsmaller{\D F(X, \D 2)}} x{'} : \TOT \exists_{f \in \D F(X, \D 2)}\big(f(x) 
=_{\mathsmaller{\D 2}} 1 \ \& \ f(x{'}) =_{\mathsmaller{\D 2}} 0 \big) 
\]
If $f \in \D F(X, \D 2)$, the following extensional subsets of $X$ 
\[ \delta_0^1(f) := \{x \in X \mid f(x) =_{\mathsmaller{\D 2}} 1\}, \]
\[ \delta_0^0(f) := \{x \in X \mid f(x) =_{\mathsmaller{\D 2}} 0\}, \]
are called detachable\index{detachable subset}, or free\index{free subset} subsets of $X$. Let also
their pair $\B \delta (f) := \big(\delta_0^1(f), \delta_0^0(f)\big)$.
\end{definition}

Clearly, 
$  x \neq_X^{\mathsmaller{\D F(X, \D 2)}} x{'} \TOT \exists_{f \in \D F(X, \D 2)}\big(f(x) \neq_{\D 2} f(x{'})\big)$,
and $\B \delta (f)$ is a complemented subset of $X$ with respect to the inequality $\neq_X^{\mathsmaller{\D F(X, \D 2)}}$.
The characteristic function $\chi_{\B \delta(f)}$ of $\B \delta(f)$ is definitionally equal to $f$ (recall that 
$f(x) =_{\mathsmaller{\D 2}} 1 :\TOT f(x) := 1$), and $\delta_0^1(f) \cup \delta_0^0(f) = X$.

\begin{definition}\label{def: operationscomplemented1}
If $\B A, \B B \in \C P^{\Disj}(X)$ and $\B C \in \C P^{\Disj}(Y)$, 
let\index{$\B A \cup \B B$}\index{$\B A \cap \B B$}\index{$- \B A$} 
\index{$\B A - \B B$}\index{$\B A \times \B C$}
\[ \B A \cup \B B  := (A^1 \cup B^1, A^0 \cap B^0), \]
\[ \B A \cap \B B := (A^1 \cap B^1, A^0 \cup B^0), \]
\[- \B A := (A^0, A^1), \]
\[ \B A - \B B := (A^1 \cap B^0, A^0 \cup B^1), \]
\[ \B A \times \B C := \big(A^1 \times C^1, \ [A^0 \times Y] \cup [X \times C^0]\big). \]
\end{definition}

The following diagrams depict $\B A \cup \B B, \B A \cap \B B$, $\B A - \B B$, and $\B A \times \B C$, respectively.

\begin{center}
\begin{tikzpicture}\begin{scope}[scale=0.8]

\node (AO) at (1,3) {};
\node (BO) at (6,5) {};
\node (AZ) at (10,3) {};
\node (BZ) at (5,0) {};

\def\AOne{plot [smooth cycle,tension=0.75] coordinates {(AO) (3,4) (4,3) (4,1) (2,1)}}
\def\BOne{plot [smooth cycle,tension=0.75] coordinates {(0,4) (2,5) (BO) (9,4) (5,3)}}
\def\AZero{plot [smooth cycle,tension=0.75] coordinates {(6,3) (8,4) (AZ) (9,1) (7,1)}}
\def\BZero{plot [smooth cycle,tension=0.75] coordinates {(0,1) (2,2) (8,2) (8,0) (BZ) (2,0)}}

\begin{scope}
	\clip\BZero;
	\draw[fill=red,opacity=0.3]\AZero;
\end{scope}
\draw[fill=blue,opacity=0.3]\AOne;
\draw[fill=white,opacity=1]\BOne;
\draw[fill=blue,opacity=0.3]\BOne;

\draw[thick,color=black]\AZero;
\draw[thick,color=black]\BZero;
\draw[thick,color=black]\BOne;
\draw[thick,color=black]\AOne;

\node [scale=1,anchor=east] at (AO) {$A^1$};
\node [scale=1,anchor=south] at (BO) {$B^1$};
\node [scale=1,anchor=west] at (AZ) {$A^0$};
\node [scale=1,anchor=north] at (BZ) {$B^0$};
\end{scope}
\end{tikzpicture}
\end{center}

\begin{center}

\begin{tikzpicture}
\begin{scope}[scale=0.8]

\node (AO) at (1,3) {};
\node (BO) at (6,5) {};
\node (AZ) at (8,3) {};
\node (BZ) at (5,0) {};

\def\AOne{plot [smooth cycle,tension=0.75] coordinates {(AO) (2,4) (3.5,3.5) (3,1) (1,1)}}
\def\BOne{plot [smooth cycle,tension=0.75] coordinates {(0.5,4) (2,5) (BO) (8,4) (5,3)}}
\def\AZero{plot [smooth cycle,tension=0.75] coordinates {(5,3) (6,4) (AZ) (7,1) (5,1)}}
\def\BZero{plot [smooth cycle,tension=0.75] coordinates {(0,1) (2,2) (7,2) (8,0) (BZ) (2,0)}}

\begin{scope}
\clip\BOne;
\draw[fill=blue,opacity=0.3]\AOne;
\end{scope}
\draw[fill=red,opacity=0.3]\AZero;
\draw[fill=white,opacity=1]\BZero;
\draw[fill=red,opacity=0.3]\BZero;

\draw[thick,color=black]\AZero;
\draw[thick,color=black]\BZero;
\draw[thick,color=black]\BOne;
\draw[thick,color=black]\AOne;

\node [scale=1,anchor=east] at (AO) {$A^1$};
\node [scale=1,anchor=south] at (BO) {$B^1$};
\node [scale=1,anchor=west] at (AZ) {$A^0$};
\node [scale=1,anchor=north] at (BZ) {$B^0$};
\end{scope}
\end{tikzpicture}
\end{center}

\begin{center}
\begin{tikzpicture}
\begin{scope}[scale=0.8]

\node (AO) at (0.5,3) {};
\node (BO) at (6,5) {};
\node (AZ) at (8,3) {};
\node (BZ) at (5,0) {};

\def\AOne{plot [smooth cycle,tension=0.75] coordinates {(AO) (2,4) (3.5,3.5) (3,1) (1,1)}}
\def\BOne{plot [smooth cycle,tension=0.75] coordinates {(0.5,4) (2,5) (BO) (8,4) (5,3) (2.5,3.5)}}
\def\AZero{plot [smooth cycle,tension=0.75] coordinates {(5,3) (6,4) (AZ) (7,1.5) (6,1)}}
\def\BZero{plot [smooth cycle,tension=0.75] coordinates {(0,1) (2,2) (6,2) (7,0) (BZ) (2,0)}}

\begin{scope}
\clip\BZero;
\draw[fill=blue,opacity=0.3]\AOne;
\end{scope}
\draw[fill=red,opacity=0.3]\BOne;
\draw[fill=white,opacity=1]\AZero;
\draw[fill=red,opacity=0.3]\AZero;

\draw[thick,color=black]\AZero;
\draw[thick,color=black]\BZero;
\draw[thick,color=black]\BOne;
\draw[thick,color=black]\AOne;

\node [scale=1,anchor=east] at (AO) {$A^1$};
\node [scale=1,anchor=south] at (BO) {$B^1$};
\node [scale=1,anchor=west] at (AZ) {$A^0$};
\node [scale=1,anchor=north] at (BZ) {$B^0$};
\end{scope}
\end{tikzpicture}
\end{center}

\begin{center}

\begin{tikzpicture}

\begin{scope}[scale=2]

\node (O) at (0,0) {};
\node (AOI) at (0.75,0) {};
\node (AOII) at (2.25,0) {};
\node (BOI) at (0,0.5) {};
\node (BOII) at (0,1.5) {};
\node (AZI) at (2.75,0) {};
\node (AZII) at (3.5,0) {};
\node (BZI) at (0,2) {};
\node (BZII) at (0,2.5) {};

\coordinate (y) at (0,3);
\coordinate (x) at (4,0);
\draw (y) node[above left] {$Y$} -- (0,0) --  (x) node[below right] {$X$};

\fill[color=blue,opacity=0.3] (AOI |- BOI) rectangle (AOII |- BOII);
\fill[color=red,opacity=0.3] (AZI |- O) rectangle (AZII |- BZI);
\fill[color=red,opacity=0.3] (O |- BZI) rectangle (AZII |- BZII);
\fill[color=red,opacity=0.3] (AZI |- BZII) rectangle (AZII |- y);
\fill[color=red,opacity=0.3] (AZII |- BZI) rectangle (x |- BZII);

\draw[color=black,thick] (AOI)++(0,-0.05) -- (AOI |- 0,0.05);
\draw[color=black,thick] (AOII)++(0,-0.05) -- (AOII |- 0,0.05);
\draw[color=black,thick] (BOI)++(-0.05,0) -- (0.05,0 |- BOI);
\draw[color=black,thick] (BOII)++(-0.05,0) -- (0.05,0 |- BOII);
\draw[color=black,thick] (AZI)++(0,-0.05) -- (AZI |- 0,0.05);
\draw[color=black,thick] (AZII)++(0,-0.05) -- (AZII |- 0,0.05);
\draw[color=black,thick] (BZI)++(-0.05,0) -- (0.05,0 |- BZI);
\draw[color=black,thick] (BZII)++(-0.05,0) -- (0.05,0 |- BZII);

\draw[color=black,dashed] (AOI) -- (AOI |- BOII);
\draw[color=black,dashed] (AOII) -- (AOII |- BOII);
\draw[color=black,dashed] (BOI) -- (AOII |- BOI);
\draw[color=black,dashed] (BOII) -- (AOII |- BOII);
\draw[color=black,dashed] (AZI) -- (AZI |- y);
\draw[color=black,dashed] (AZII) -- (AZII |- y);
\draw[color=black,dashed] (BZI) -- (x |- BZI);
\draw[color=black,dashed] (BZII) -- (x |- BZII);

\draw [decorate,decoration={brace,amplitude=10pt}]
	(AOII)++(0,-0.1) -- (AOI |- 0,-0.1) node [black,midway,yshift=-0.7cm] {$A^1$};
\draw [decorate,decoration={brace,amplitude=10pt}]
(AZII)++(0,-0.1) -- (AZI |- 0,-0.1) node [black,midway,yshift=-0.7cm] {$A^0$};

\draw [decorate,decoration={brace,amplitude=10pt}]
	(BOI)++(-0.1,0) -- (-0.1,0 |- BOII) node [black,midway,xshift=-0.7cm] {$C^1$};
\draw [decorate,decoration={brace,amplitude=10pt}]
(BZI)++(-0.1,0) -- (-0.1,0 |- BZII) node [black,midway,xshift=-0.7cm] {$C^0$};

\end{scope}
\end{tikzpicture}
\end{center}

\begin{remark}\label{rem: opercomp1}
If $\B A, \B B \in \C P^{\mathsmaller{\Disj}}(X)$ and $\B C \in \C P^{\mathsmaller{\Disj}}(Y)$, then
$\B A \cup \B B$, $\B A \cap \B B$, $- \B A$, and $\B A - \B B$ are in $\C P^{\mathsmaller{\Disj}}(X)$ 
and $\B A \times \B C$ is in $\C P^{\mathsmaller{\Disj}} (X \times Y)$.

\end{remark}

\begin{proof}
We show only the last membership. If $(a_1, b_1) \in A^1 \times B^1$ and $(a_0, b_0) \in A^0 \times B^0$, 
then $i_{A^1}^X(a_1) \neq_X i_{A^0}^X(a_0)$ and $i_{B^1}^Y(b_1) \neq_Y i_{B^0}^Y(b_0)$. By definition 
\[ i_{A^1 \times B^1}^{X \times Y}(a_1, b_1) := \big(i_{A^1}^X(a_1), i_{B^1}^Y(b_1)\big). \]
If $(a_0, y) \in A^0 \times Y$, then $(i_{A^0}^X \times \id_Y)(a_0, y) := (i_{A^0}^X(a_0), y)$, and if
$(x, b_0) \in X \times B^0$, then $(\id_X \times i_{B^0}^Y)(x, b_0) := (x, i_{B^0}^Y(b_0))$. In both 
cases we get the required inequality.
\end{proof}

\begin{remark}\label{rem: complemented2}
Let $\B A, \B B$ and $\B C$ be in $\C P^{\mathsmaller{\Disj}} (X)$. The following hold:\\[1mm]
\normalfont (i) 
\itshape  $- (- \B A) := \B A$.\\[1mm]
\normalfont (ii) 
\itshape $- (\B A \cup \B B) := (- \B A) \cap (- \B B)$.\\[1mm]
\normalfont (iii) 
\itshape $- (\B A \cap \B B) := (- \B A) \cup (- \B B)$.\\[1mm]
\normalfont (iv) 
\itshape $\B A \cup (\B B \cap \B C) =_{\C P^{\mathsmaller{\Disj}} (X)} (\B A \cup \B B) \cap (\B A \cup \B C)$.\\[1mm]
\normalfont (v) 
\itshape $\B A \cap (\B B \cup \B C) =_{\C P^{\mathsmaller{\Disj}} (X)} (\B A \cap \B B) \cup (\B A \cap \B C)$.\\[1mm]
\normalfont (vi) 
\itshape $\B A - \B B := \B A \cap (- \B B)$.\\[1mm]
\normalfont (vii) 
\itshape $\B A \subseteq \B B \TOT (\B A \cap \B B) =_{\C P^{\mathsmaller{\Disj}}(X)} \B A$.\\[1mm]
\normalfont (viii)
\itshape $\B A \subseteq \B B \TOT - \B B \subseteq - \B A$.\\[1mm]
\normalfont (ix)
\itshape If $\B A \subseteq \B B$ and $\B B \subseteq \B C$, then $\B A \subseteq \B C$.
\end{remark}

\begin{proposition}\label{prp: complemented4}
Let $\B A \in \C P^{\mathsmaller{\Disj}}(X)$ and $B, C \in \C P^{\mathsmaller{\Disj}}(Y)$.\\[1mm]
\normalfont (i)
\itshape $\B A \times (\B B \cup \B C) =_{\C P^{\mathsmaller{\Disj}}(X \times Y)} (\B A \times \B B) \cup
(\B A \times \B C)$.\\[1mm]
\normalfont (ii)
\itshape $\B A \times (\B B \cap \B C) =_{\C P^{\mathsmaller{\Disj}}(X \times Y)} (\B A \times \B B) \cap
(\B A \times \B C)$.\end{proposition}

\begin{proof}
We prove only (i). We have that
\begin{align*}
\B A \times (\B B \cup \B C) & := (A^1, A^0) \times (B^1 \cup C^1, B^0 \cap C^0)\\
& := \big(A^1 \times (B^1 \cup C^1), (A^0 \times Y) \cup [X \times (B^0 \cap C^0)]\big)\\
& =_{\C P^{\mathsmaller{\Disj}}(X \times Y)}  \big((A^1 \times B^1) \cup (A^1 \times C^1), [(A^0 \times Y)
\cup (X \times B^0)] \ \cap\\
& \ \ \ \ \ \ \ \ \ \ \ \ \ \ \ \  \cap  [(A^0 \times Y) \cup (X \times C^0)]\big)\\
& := (\B A \times \B B) \times (\B A \times \B C).\qedhere
\end{align*}
\end{proof}

\begin{proposition}\label{prp: complemented3}
Let the sets $(X, =_X, \neq_X^f)$ and $(Y, =_Y, \neq_Y)$, where $f \colon X \to Y$ $($see 
Remark~\ref{rem: apartness2}$)$.
Let also $\B A := (A^1, A^0)$ and $\B B := (B^1, B^0)$ in $\C P^{\mathsmaller{\Disj}} (Y)$.\\[1mm]  
\normalfont (i) 
\itshape $f^{-1}(\B A) := \big(f^{-1}(A^1), f^{-1}(A^0)\big) \in \C P^{\Disj} (X)$.\\[1mm]
\normalfont (ii) 
\itshape $f^{-1}(\B A \cup \B B) =_{\mathsmaller{\C P^{\Disj} (X)}}
f^{-1}(\B A) \cup f^{-1}(\B B)$.\\[1mm]
\normalfont (iii) 
\itshape $f^{-1}(\B A \cap \B B) =_{\mathsmaller{\C P^{\Disj} (X)}} f^{-1}(\B A) \cap f^{-1}(\B B)$.\\[1mm]
\normalfont (iv)
\itshape $f^{-1}(- \B A) =_{\mathsmaller{\C P^{\Disj} (X)}} - f^{-1}(\B A)$.\\[1mm]
\normalfont (v)
\itshape $f^{-1}(\B A - \B B) =_{\mathsmaller{\C P^{\Disj} (X)}} f^{-1}(\B A) - f^{-1}(\B B)$.
\end{proposition}

\begin{proof}
(i) By Definition~\ref{def: image} we have that
\[ f^{-1}(A^1) := \{(x, a_1) \in X \times A^1 \mid f(x) =_Y i_{A^1}^X (a_1)\}, \ \ \ \ 
i_{\mathsmaller{f^{-1}(A^1)}}^X(x, a_1) := x, \]
\[ f^{-1}(A^0) := \big(\{(x, a_0) \in X \times A^0 \mid f(x) =_Y i_{A^0}^X (a_0)\}, \ \ \ \ 
i_{\mathsmaller{f^{-1}(A^0)}}^X(x, a_0) := x. \]
Let $(x, a_1) \in f^{-1}(A^1)$ and $(z, a_0) \in f^{-1}(A^0)$. By the extensionality of $\neq_Y$ we have that 
\[ i_{\mathsmaller{f^{-1}(A^1)}}^X(x, a_1) \neq_X^f i_{\mathsmaller{f^{-1}(A^0)}}^X(z, a_0)
:\TOT x \neq_X^f z :\TOT f(x) \neq_Y f(z) \TOT i_{A^1}^X (a_1) \neq_Y i_{A^0}^X (a_0), \]
and the last inequality holds by the hypothesis $\B A \in \C P^{\mathsmaller{\Disj}} (Y)$.
Next we show only (ii):
\begin{align*}
f^{-1}(\B A \cup \B B) & := f^{-1}\big(A^1 \cup B^1, A^0 \cap B^0 \big)\\
& := \big(f^{-1}(A^1 \cup B^1), f^{-1}(A^0 \cap B^0)\big)\\
& = \big(f^{-1}(A^1) \cup f^{-1}(B^1), f^{-1}(A^0) \cap f^{-1}(B^0)\big)\\
& := f^{-1}(\B A) \cup f^{-1}(\B B).\qedhere
\end{align*}
\end{proof}

Alternatively, one can define the following operations between complemented subsets. 

\begin{definition}\label{def: operationscomplemented2}
If $\B A, \B B \in \C P^{\Disj}(X)$ and $\B C \in \C P^{\Disj}(Y)$, let\index{$\B A \vee \B B$}\index{$\B A \wedge \B B$}
\index{$\B A \ominus \B B$}\index{$\B A \otimes \B C$}
\[ \B A \vee \B B := \big([A^1 \cap B^1] \cup [A^1 \cap B^0] \cup [A^0 \cap B^1], \ A^0 \cap B^0\big), \] 
\[ \B A \wedge \B B := \big(A^1 \cap B^1, \ [A^1 \cap B^0] \cup [A^0 \cap B^1] \cup [A^0 \cap B^0]\big), \]
\[ \B A \ominus \B B := \B A \wedge (- \B B), \]
\[ \B A \otimes \B C := \big(A^1 \times C^1, \ [A^1 \times C^0] \cup [A^0 \times C^1] \cup [A^0 \times C^0]\big), \]

\end{definition}

The following diagrams depict $\B A \vee \B B, \B A \wedge \B B$, $\B A \ominus \B B$, and $\B A \otimes \B C$, 
respectively.

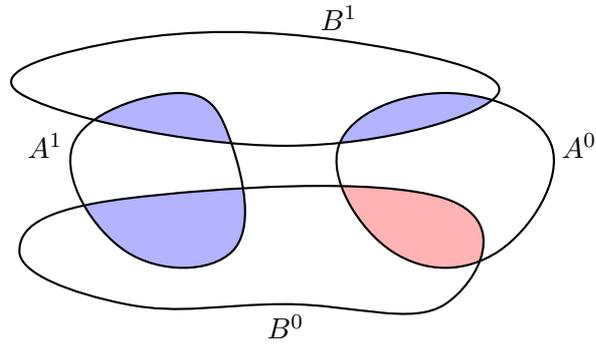
\begin{center}
\begin{tikzpicture}\begin{scope}[scale=0.7]

\node (AO) at (1,3) {};
\node (BO) at (6,5) {};
\node (AZ) at (10,3) {};
\node (BZ) at (5,0) {};

\def\AOne{plot [smooth cycle,tension=0.75] coordinates {(AO) (3,4) (4,3) (4,1) (2,1)}}
\def\BOne{plot [smooth cycle,tension=0.75] coordinates {(0,4) (2,5) (BO) (9,4) (5,3)}}
\def\AZero{plot [smooth cycle,tension=0.75] coordinates {(6,3) (8,4) (AZ) (9,1) (7,1)}}
\def\BZero{plot [smooth cycle,tension=0.75] coordinates {(0,1) (2,2) (8,2) (8,0) (BZ) (2,0)}}

\begin{scope}
	\clip\BZero;
	\draw[fill=red,opacity=0.3]\AZero;
	\draw[fill=blue,opacity=0.3]\AOne;
\end{scope}
\begin{scope}
	\clip\BOne;
	\draw[fill=blue,opacity=0.3]\AOne;
	\draw[fill=blue,opacity=0.3]\AZero;
\end{scope}

\draw[thick,color=black]\AZero;
\draw[thick,color=black]\BZero;
\draw[thick,color=black]\BOne;
\draw[thick,color=black]\AOne;

\node [scale=1,anchor=east] at (AO) {$A^1$};
\node [scale=1,anchor=south] at (BO) {$B^1$};
\node [scale=1,anchor=west] at (AZ) {$A^0$};
\node [scale=1,anchor=north] at (BZ) {$B^0$};
\end{scope}
\end{tikzpicture}
\end{center}

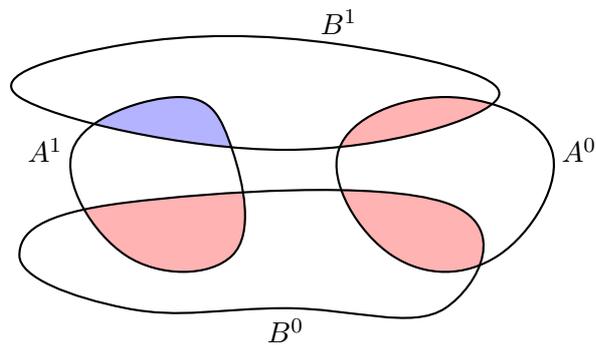
\begin{center}
\begin{tikzpicture}\begin{scope}[scale=0.7]

\node (AO) at (1,3) {};
\node (BO) at (6,5) {};
\node (AZ) at (10,3) {};
\node (BZ) at (5,0) {};

\def\AOne{plot [smooth cycle,tension=0.75] coordinates {(AO) (3,4) (4,3) (4,1) (2,1)}}
\def\BOne{plot [smooth cycle,tension=0.75] coordinates {(0,4) (2,5) (BO) (9,4) (5,3)}}
\def\AZero{plot [smooth cycle,tension=0.75] coordinates {(6,3) (8,4) (AZ) (9,1) (7,1)}}
\def\BZero{plot [smooth cycle,tension=0.75] coordinates {(0,1) (2,2) (8,2) (8,0) (BZ) (2,0)}}

\begin{scope}
	\clip\BZero;
	\draw[fill=red,opacity=0.3]\AZero;
	\draw[fill=red,opacity=0.3]\AOne;
\end{scope}
\begin{scope}
	\clip\BOne;
	\draw[fill=blue,opacity=0.3]\AOne;
	\draw[fill=red,opacity=0.3]\AZero;
\end{scope}

\draw[thick,color=black]\AZero;
\draw[thick,color=black]\BZero;
\draw[thick,color=black]\BOne;
\draw[thick,color=black]\AOne;

\node [scale=1,anchor=east] at (AO) {$A^1$};
\node [scale=1,anchor=south] at (BO) {$B^1$};
\node [scale=1,anchor=west] at (AZ) {$A^0$};
\node [scale=1,anchor=north] at (BZ) {$B^0$};
\end{scope}
\end{tikzpicture}
\end{center}

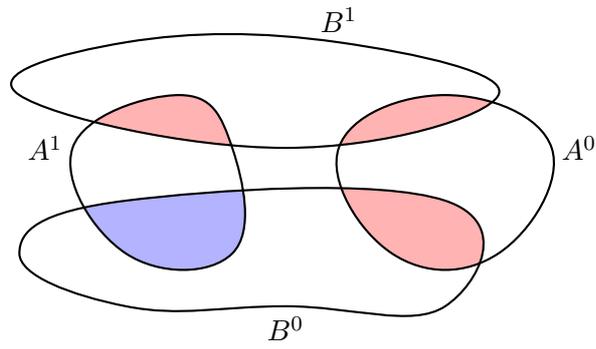
\begin{center}
\begin{tikzpicture}\begin{scope}[scale=0.7]

\node (AO) at (1,3) {};
\node (BO) at (6,5) {};
\node (AZ) at (10,3) {};
\node (BZ) at (5,0) {};

\def\AOne{plot [smooth cycle,tension=0.75] coordinates {(AO) (3,4) (4,3) (4,1) (2,1)}}
\def\BOne{plot [smooth cycle,tension=0.75] coordinates {(0,4) (2,5) (BO) (9,4) (5,3)}}
\def\AZero{plot [smooth cycle,tension=0.75] coordinates {(6,3) (8,4) (AZ) (9,1) (7,1)}}
\def\BZero{plot [smooth cycle,tension=0.75] coordinates {(0,1) (2,2) (8,2) (8,0) (BZ) (2,0)}}

\begin{scope}
	\clip\BZero;
	\draw[fill=red,opacity=0.3]\AZero;
	\draw[fill=blue,opacity=0.3]\AOne;
\end{scope}
\begin{scope}
	\clip\BOne;
	\draw[fill=red,opacity=0.3]\AOne;
	\draw[fill=red,opacity=0.3]\AZero;
\end{scope}

\draw[thick,color=black]\AZero;
\draw[thick,color=black]\BZero;
\draw[thick,color=black]\BOne;
\draw[thick,color=black]\AOne;

\node [scale=1,anchor=east] at (AO) {$A^1$};
\node [scale=1,anchor=south] at (BO) {$B^1$};
\node [scale=1,anchor=west] at (AZ) {$A^0$};
\node [scale=1,anchor=north] at (BZ) {$B^0$};
\end{scope}
\end{tikzpicture}
\end{center}

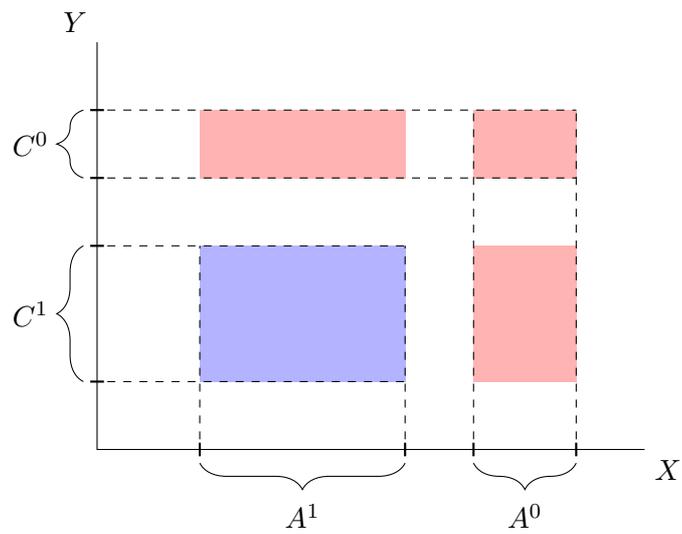
\begin{center}

\begin{tikzpicture}

\begin{scope}[scale=1.8]

\node (O) at (0,0) {};
\node (AOI) at (0.75,0) {};
\node (AOII) at (2.25,0) {};
\node (BOI) at (0,0.5) {};
\node (BOII) at (0,1.5) {};
\node (AZI) at (2.75,0) {};
\node (AZII) at (3.5,0) {};
\node (BZI) at (0,2) {};
\node (BZII) at (0,2.5) {};
\coordinate (y) at (0,3);
\coordinate (x) at (4,0);

\draw (y) node[above left] {$Y$} -- (0,0) --  (x) node[below right] {$X$};

\fill[color=blue,opacity=0.3] (AOI |- BOI) rectangle (AOII |- BOII);

\fill[color=red,opacity=0.3] (AZI |- BOI) rectangle (AZII |- BOII);
\fill[color=red,opacity=0.3] (AOI |- BZI) rectangle (AOII |- BZII);
\fill[color=red,opacity=0.3] (AZI |- BZII) rectangle (AZII |- BZI);

\draw[color=black,thick] (AOI)++(0,-0.05) -- (AOI |- 0,0.05);
\draw[color=black,thick] (AOII)++(0,-0.05) -- (AOII |- 0,0.05);
\draw[color=black,thick] (BOI)++(-0.05,0) -- (0.05,0 |- BOI);
\draw[color=black,thick] (BOII)++(-0.05,0) -- (0.05,0 |- BOII);
\draw[color=black,thick] (AZI)++(0,-0.05) -- (AZI |- 0,0.05);
\draw[color=black,thick] (AZII)++(0,-0.05) -- (AZII |- 0,0.05);
\draw[color=black,thick] (BZI)++(-0.05,0) -- (0.05,0 |- BZI);
\draw[color=black,thick] (BZII)++(-0.05,0) -- (0.05,0 |- BZII);

\draw[color=black,dashed] (AOI) -- (AOI |- BOII);
\draw[color=black,dashed] (AOII) -- (AOII |- BOII);
\draw[color=black,dashed] (BOI) -- (AOII |- BOI);
\draw[color=black,dashed] (BOII) -- (AOII |- BOII);
\draw[color=black,dashed] (AZI) -- (AZI |- BZII);
\draw[color=black,dashed] (AZII) -- (AZII |- BZII);
\draw[color=black,dashed] (BZI) -- (AZII |- BZI);
\draw[color=black,dashed] (BZII) -- (AZII |- BZII);

\draw [decorate,decoration={brace,amplitude=10pt}]
(AOII)++(0,-0.1) -- (AOI |- 0,-0.1) node [black,midway,yshift=-0.7cm] {$A^1$};
\draw [decorate,decoration={brace,amplitude=10pt}]
(AZII)++(0,-0.1) -- (AZI |- 0,-0.1) node [black,midway,yshift=-0.7cm] {$A^0$};

\draw [decorate,decoration={brace,amplitude=10pt}]
(BOI)++(-0.1,0) -- (-0.1,0 |- BOII) node [black,midway,xshift=-0.7cm] {$C^1$};
\draw [decorate,decoration={brace,amplitude=10pt}]
(BZI)++(-0.1,0) -- (-0.1,0 |- BZII) node [black,midway,xshift=-0.7cm] {$C^0$};

\end{scope}
\end{tikzpicture}

\end{center}

With the previous definitions the corresponding characteristic functions are expressed through the characteristic 
functions of
$\B A$ and $\B B$.

\begin{remark}\label{rem: complemented1}
If $\B A, \B B$ are complemented subsets of $X$, then $\B A \vee \B B, \B A \wedge \B B, \B A - \B B$
and $- \B A$ are complemented subsets of $X$ with characteristic functions
\[ \chi_{\B A \vee \B B} =_{\C F(X, \D 2)} \chi_{\B A} \vee \chi_{\B B}, \  \ \chi_{\B A \wedge \B B} =_{\C F(X, 
\D 2)} \chi_{\B A} \cdot \chi_{\B B}, \ \ 
\chi_{\B A - \B B} =_{\C F(X, \D 2)} \chi_{\B A}(1 - \chi_{\B B}), \ \ \ \] 
\[ \chi_{\B A \otimes \B B}(x, y) =_{\C F(X \times X, \D 2)} \chi_{\B A}(x) \cdot \chi_{\B B}(y), \ \ \ 
\chi_{- \B A} =_{\C F(X, \D 2)} 1 - \chi_{\B A}. 
 \] 
\end{remark}

\begin{proof}
We show only the equality $\chi_{\B A \wedge \B B} =_{\C F(X, \D 2)} \chi_{\B A} 
\cdot \chi_{\B B}$. By Definition~\ref{def: multin2} the multiplication of the partial maps $\chi_{\B A} \colon
\Dm(\B A) \to \D 2$ and $\chi_{\B B} \colon \Dm(\B B) \to \D 2$ is 
the partial function 
\[ \chi_{\B A} \cdot \chi_{\B B} := \big(\Dm(\B A) \cap \Dm(\B B), i_{\Dm(\B A) \cap \Dm(\B B)}^X, (\chi_{\B A} 
\cdot \chi_{\B B})_{\Dm(\B A) \cap \Dm(\B B)}^{\D 2} \big), \]
\[ (\chi_{\B A} \cdot \chi_{\B B})_{\Dm(\B A) \cap \Dm(\B B)}^{\D 2}(u, w) := \chi_{\B A}(u) \cdot \chi_{\B B}(w), \]
for every $(u, w) \in \Dm(\B A) \cap \Dm(\B B)$. The partial function $\chi_{\B A \wedge \B B}$ is the triplet
\[ \chi_{\B A \wedge \B B} := \big(\Dm(\B A \wedge \B B), i_{\Dm(\B A \wedge \B B)}^X, (\chi_{\B A 
\wedge \B B})_{\Dm(\B A \wedge \B B)}^{\D 2} \big). \]
Since $\Dm(\B A \wedge \B B) =_{\C P(X)} \Dm(\B A) \cap \Dm(\B B)$, and if $(f, g) \colon \Dm(\B A \wedge \B B)
=_{\C P(X)} \Dm(\B A) \cap \Dm(\B B)$, it is straightforward to show that also the following outer diagram commutes  
\begin{center}
\resizebox{12cm}{!}{%
\begin{tikzpicture}

\node (E) at (0,0) {$\Dm(\B A \wedge \B B)$};
\node[right=of E] (L) {};
\node[right=of L] (M) {};
\node[right=of M] (B) {};
\node[right=of B] (F) {$\Dm(\B A) \cap \Dm(\B B)$};
\node[below=of M] (N) {};
\node[below=of N] (A) {$X$};
\node[below=of A] (C) {$\D 2$};

\draw[->,bend left] (E) to node [midway,below] {$f$} (F);
\draw[->,bend left] (F) to node [midway,above] {$g$} (E);
\draw[->,bend right=50] (E) to node [midway,left] {$(\chi_{\B A \wedge \B B})_{\Dm(\B A \wedge \B B)}^{\D 2} \ \ $} (C);
\draw[->,bend left=50] (F) to node [midway,right] {$ \ \ (\chi_{\B A} \cdot \chi_{\B B})_{\Dm(\B A) \cap \Dm(\B B)}^{\D 2}$} (C);
\draw[right hook->,bend right=40] (E) to node [midway,right] {$ \ i_{\Dm(\B A \wedge \B B)}^X$} (A);
\draw[left hook->,bend left=40] (F) to node [midway,left] {$i_{\Dm(\B A) \cap \Dm(\B B)}^X \ $} (A);

\end{tikzpicture}
}
\end{center} 
and hence the two partial functions are equal in $\C F(X, \D 2)$. 
\end{proof}

\section{Notes}
\label{sec: notes2}

\begin{note}\label{not: primitive}
\normalfont
In~\cite{Gr81} Greenleaf introduced predicates on objects through the totality $\Omega$ of propositions 
and then he defined $\C P(X)$ as $\D F(X, \Omega)$. A similar treatment of the powerset $\C P(X)$ is found 
in~\cite{Sh18}. For us a predicate on a set $X$ is a bounded formula $P(x)$ with $x$ as a free variable.
In order to define new objects from $X$ through $P$ we ask $P$ to be extensional.

\end{note}

\begin{note}\label{not: cantor}
\normalfont
In~\cite{Ca82}, pp. 114-5, Cantor described a set as follows:
\begin{quote}
A manifold (a sum, a set) of elements belonging to some conceptual
sphere is called well-defined if, on the basis of its definition and in
accordance with the logical principle of the excluded third, it must
be regarded as internally determined, both whether any object of
that conceptual sphere belongs as an element to the mentioned set,
and also whether two objects belonging to the set, in spite of formal
differences in the mode of givenness, are equal to each other or not.
\end{quote}
Bishop's intuitive notion of set is similar to Cantor's, except that he does not invoke
the principle of the excluded middle ($\PEM$). As it was pointed to me by W.~Sieg, Dedekind's
primitive notions in~\cite{De88}  were ``systems'' and ``transformations of systems''.
Notice that here we study defined totalities that are not defined inductively. 
The inductively defined sets are expected to be studied in a future work within an extension $\BST^*$ of $\BST$.

\end{note}

\begin{note}\label{not: otherprimitive}
\normalfont
Although $\Nat$ is the only primitive set considered in $\BST$, one 
could, in principle, add more primitive sets. E.g., a primitive set of Booleans, of integers, and, 
more interestingly, a primitive continuous interval, or a primitive real line (see~\cite{Br99} for an axiomatic
treatment of the set $\Real$ of reals within $\BISH$).  
\end{note}

\begin{note}\label{not: defproduct}
\normalfont
In Martin-L\"of type theory the definitional, or judgemental equality $a := b$, where $a, b$ are terms
of some type $A$, is never used in a formula. We permit the use of the 
definitional equality $:=$ for membership conditions only.
In the membership condition for the product we use the primitive notion of a pair.
The membership condition for an extensional subset $X_P$ of $X$ implies that an object $x$ ``has not unique typing'',
as it can be an element of more than one sets.

\end{note}

\begin{note}\label{not: defdiscrete}
\normalfont
The positively defined notion of discrete set used here comes from~\cite{MRR88}, p.~9. There it is also mentioned that 
a set without a specified inequality i.e., a pair $(X, =_X)$, is discrete, if 
$\forall_{x,y \in X}\big(x =_X y \ \vee \ \neg(x =_X y)\big)$.
In~\cite{Pa13} it is mentioned that the above discreteness of $\D F(\Nat, \Nat)$ implies 
the non-constructive principle ``weak $\LPO$''
\[ \forall_{f \in \D F(\Nat, \Nat)}\bigg(\forall_{n \in \Nat}\big(f(n) =_{\Nat} 0\big) \ \vee \ \neg 
\forall_{n \in \Nat}\big(f(n) =_{\Nat} 0\big)\bigg). \]
Because of a result of Bauer and Swan in~\cite{BS18}, we cannot show in $\BISH$ the existence of an
uncountable separable metric space, hence, using the discrete metric, the existence of an uncountable
discrete set. Note that in~\cite{Bi67}, p.~66, a set $S$ is called discrete, if the set 
$D := \{(s, t) \in S \times S \mid s=_S t\}$ is a free\index{free subset}, or a detachable\index{detachable subset} subset of $S \times S$.
In Definition~\ref{def: diagonal} we use the symbol $D(S)$ for $D$ and we call it the diagonal of $S$. We employ here the diagonal of a set in the fundamental definition of a set-indexed family of sets (Definition~\ref{def: famofsets}).
\end{note}

\begin{note}\label{not: Onnegation}
\normalfont
In~\cite{Bi67} and~\cite{BB85}, the negation $\neg \phi$ of a formula $\phi$ is not mentioned explicitly. E.g., the exact
writing of condition $(\Ap_1)$ in Definition~\ref{def: apartness} is ``if $x =_X y$ and $x \neq_X y$, then $0 =_{\Nat} 1$''.
Similarly, the condition of tightness in Definition~\ref{def: apartness} is written as follows:
``if $x \neq_X y$ entails $0 = 1$, then $x =_X y$''. hence, if $\neq_X$ is tight, the implication $x \neq_X y \To 
0 =_{\Nat} 1$ is logically equivalent to the (positively defined, if $X$ is a defined totality) equality $x =_X y$.
Within intuitionistic logic one defines $\neg \phi := \phi \To \bot$. 
\end{note}

\begin{note}\label{not: nsets}
\normalfont
The definitions of $(-2)$-sets and $(-1)$-sets are proof-irrelevant translations of the corresponding notions in $\HoTT$,
which were introduced by Voevodsky (see~\cite{HoTT13}). The definition of a $0$-\textit{set} requires to determine a set 
$\Eq^X(x, y)$ of witnesses of the equality $x =_X y$. This is done in a universal way in $\MLTT$, while in $\BST$ 
in a ``local'' way, and by definition (see Definition~\ref{def: 0set}).
\end{note}

\begin{note}\label{not: operation}
\normalfont
In the literature of constructive mathematics (see e.g.,~\cite{Be85}, pp.~34--35)
the term \textit{preset} is used for a totality. 
Also, the term \textit{operation} is used for a non-dependent assignment routine from a totality $X$ to a totality $Y$
(see~\cite{Be85}, p.~44), while we use it only for a non-dependent assignment routine from a set $X$ to a set $Y$.
\end{note}

\begin{note}\label{not: LANC}
\normalfont
The notion of uniqueness associated to the definition of a function is \emph{local}, in the following sense:
if $f \colon X \to Y$, it is immediate to show that $\forall_{x \in X}\exists!_{y \in Y}\big(f(x) =_Y y\big)$.
The converse is the local version of Myhill's axiom of non-choice $(\LANC)$. Let $P(x, y)$ be an extensional property on
$X \times Y$ i.e., $\forall_{x, x{'}, y, y{'} \in X}\big([x =_X x{'} \ \& \ y =_Y y{'} \ \& \ P(x, y)]
\To P(x{'}, y{'})\big)$. The principle 
$(\LANC)$\index{local version of Myhill's axiom of non-choice}\index{$(\LANC)$} is the formula
\[ \forall_{x \in X}\exists!_{y \in Y}P(x, y) \To \exists_{f \in \D F(X, Y)}\forall_{x \in X}\big(P(x, 
f(x))\big). \]
Notice that $\LANC$ provides the existence of a function for which we only know how its outputs behave with 
respect to the equality of $Y$, and it gives no information on how $f$ behaves definitionally.
If we define $Q_x  (y) : = P(x, y)$, then if we suppose $Q_x (f(x))$ and $Q_x (g(x))$, for some $f, g \in
\D F(X, Y)$, we get $f(x) =_Y y =_Y g(x)$, and then $(\LANC)$ implies 
\[ \forall_{x \in X}\exists!_{y \in Y}P(x, y) \To \exists!_{f \in \D F(X, Y)}\forall_{x \in X}\big(P(x, 
f(x))\big). \]
We can use $(\LANC)$ to view an arbitrary subset $(A, i_A^X)$ of $X$ as an extensional subset of $X$. 
If $(A, i_A^X) \in \C P(X)$, then the property $P_A$ on $X$ defined by 
$P_A(x) :=  \exists_{a \in A}\big(i_A^X (a) =_X x\big)$,
is extensional, and $(i_A^X, j_A^{X}) : X_{P_A} =_{\C P(X)} (A, i_A^X)$, for some function $j_A^{X} \colon X_{P_A} \to A$.
To show this, let $x, y \in X$ such that $P_A(x)$ and $x =_X y$. By transitivity of $=_X$, if $i_A^{X} (a) =_X x$, then 
$i_A^{X} (a) =_X y$. If $x \in X$ and $a, b \in A$ such that $i_A^{X} (a) =_X x =_X i_A^{X} (b)$, then $a =_A b$ i.e., 
$\forall_{x \in X_{P_A}}\exists!_{a  \in A}\big(i_A^{X} (a) =_X x\big)$,
and since the property $Q(x, a) : \TOT i_A^{X} (a) =_X x$ is extensional on $X_{P_A} \times A$, by $(\LANC)$
there is  a (unique) function $j_A^{X} : X_{P_A} \to A$, such that for every $x \in X_P$ we have that 
$i_A^X(j_A ^{X}(x)) =_X x$, and the required diagram commutes.
The principle $(\LANC)$, which is also considered in~\cite{Be81}, is included in Myhill's system $\CST$ (see~\cite{My75})
as a principle of generating functions. This is in contrast to Bishop's algorithmic approach to the concept of
function.
\end{note}

\begin{note}\label{not: deffunction}
 \normalfont
In~\cite{BB85}, p.~67, a function $f : A \to B$ is defined as a finite routine which,
applied to \textit{any} element of $A$, produces an element $b \equiv f(a)$ of $B$, such that $f(a) =_B f(a{'})$,
whenever $a =_A a{'}$. In~\cite{BB85}, p.~15, we read that $f$ ``affords an explicit, finite mechanical reduction
of the procedure for 
constructing $f(a)$ to the procedure for constructing $a$''. 
The pattern of defining a function $\fXY$ by first defining an operation $f \colon X \sto Y$, and then proving that 
$f$ is a function, is implicit in the more elementary parts of~\cite{Bi67} and~\cite{BB85}, and more explicit 
in the later parts of the books. E.g., in~\cite{BB85}, p.~199, an inhabited subset $U$ of 
$\D C$ \textit{has the maximal extent property}, if there
is an operation $\mu$ from $U$ to $\Real^+$ satisfying certain properties. One can show \textit{afterwords} 
that $U$ is open and $\mu$ is a function on $U$. This property is used in Bishop's proof 
of the Riemann mapping theorem (see~\cite{BB85}, pp.~209--210).

\end{note}

\begin{note}\label{not: functionspace}
\normalfont
Regarding the set-character of $\D F(X, Y)$, Bishop, in~\cite{BB85}, p.~67, writes:
\begin{quote}
When $X$ is not countable, the set $\D F(X, Y)$ seems to have little practical interest, because 
to get a hold on its structure is too hard. For instance, it has been asserted by Brouwer that all
functions in $\D F(\Real, \Real)$ are continuous, but no acceptable proof of this assertion is known.
\end{quote}
Similar problems occur though, in function spaces where the domain of the functions is 
a countable set. E.g., we cannot accept constructively (i.e., in the sense of Bishop) that the Cantor
space $\D F(\Nat, 2)$ satisfies Markov's principle, but no one that we know of has doubted the
set-character of $\D F(\Nat, 2)$. The possibility of doubting the set-character of the Baire space $\D F(\Nat, \Nat)$ 
is discussed by Beeson in~\cite{Be85}, p.~46.
\end{note}

\begin{note}\label{not: funext}
\normalfont
In intensional Martin-L\"of Type Theory the type 
\[\bigg(\prod_{x \colon X}f(x) = g(x)\bigg) \to f = g \]
is not provable (inhabited), and its inhabitance is known as the 
\textit{axiom of function extensionality}\index{axiom of function extensionality} $(\FunExt)$\index{$(\FunExt)$}.
In $\BST$ this axiom is part of the canonical definition of the function space $\D F(X, Y)$. Because of this,
many results in $\MLTT + \FunExt$ are translatable in $\BST$ (see Chapter~\ref{chapter: proofrelevance}
).
\end{note}

\begin{note}\label{not: univalence}
\normalfont
The totality $\D V_0$ is not mentioned by Bishop, although it is necessary, if we want to formulate the 
fundamental notion of a set-indexed family of sets. 
The defined 
equality on the universe $\D V_0$ expresses that $\D V_0$ is \textit{univalent}, as isomorphic sets are
equal in $\D V_0$. In univalent type theory, which is $\MLTT$ extended with Voevodsky's 
\textit{axiom of univalence}\index{axiom of univalence} $\UA$\index{$\UA$}
(see~\cite{HoTT13}), the existence of a pair of quasi-inverses between types $A$ and $B$ implies
that they are equivalent in Voevodsky's sense, and by the univalence axiom, also propositionally equal. 
The axiom $\UA$ is partially translated in $\BST$ as the canonical definition of $\D V_0$. 
Because of this, results in $\MLTT + \UA$ that do not raise the level of the universe are translatable in $\BST$.
For example, Proposition~\ref{prp: univalence1} is lemma 4.9.2 in book HoTT~\cite{HoTT13}, where $\UA$ is used 
in its proof: if $e \colon X \simeq Y$, then $Z \to X \simeq Z \to Y$, and by
$\UA$ we get $e = \idtoeqv (p)$, for some $p \colon X =_{\C U} Y$. Notice that in the 
formulation of this lemma the universe-level is not raised.
\end{note}

\begin{note}\label{not: dar}
\normalfont
The notion of a dependent operation is explicitly mentioned by Bishop in~\cite{Bi67}, p.~65, and
repeated in~\cite{BB85}, p.~70, in the definition of the intersection of a family of subsets of a set 
indexed by some set 
$T$:
\begin{quote}
an element $u$ of $\bigcap_{t \in T}\lambda(t)$ is a finite routine which
associates an element $x_t$ of $\lambda(t)$ with each element $t$ of $T$, such that
$i_t(x_t) = i_{t{'}}(x_{t{'}})$ whenever $t, t{'} \in T$.
\end{quote}
This definition corresponds to our Definition~\ref{def: intfamilyofsubsets}.
\end{note}

\begin{note}\label{not: ondefsubset}
 \normalfont
 Bishop's definition of a subset of a set is related to the notion of a subobject in Category Theory
 (see~\cite{Aw10}, p.~89, and~\cite{Go84}, p.~75). In practice the subsets of a set $X$ are defined through an
extensional property on $X$. In~\cite{BR87}, p.~7, this approach to the notion of a subset is considered 
as its \textit{definition}. Note
that there the implication $x =_{X} y \To (P(y) \To P(x))$ is also included in 
the definition of an extensional property, something which follows though, from the symmetry of $=_X$.
Such a form of separation axiom is used implicitly in~\cite{Bi67} and in~\cite{BB85}. Myhill used in 
his system $\CST$ the axiom of bounded separation to implement the notion of an extensional subset of $X$. 
This axiom is also
included in Aczel's system $\CZF$ (see~\cite{AR10}, p.~26).
\end{note}

\begin{note}\label{not: defsubset}
\normalfont
One could have defined the equality $=_{A \cup B}$ without relying on the non-dependent assignment routine 
$i_{A \cup B}^X$. If we define first 
$$z =_{A \cup B} w  : \TOT \left\{ \begin{array}{lllllll}
                 i_A^X (z)  =_X i_A (w)  &\mbox{, $z, w \in A$}\\
                 {}\\
                 i_A^X (z)  =_X i_B (w)  &\mbox{, $z \in A \ \& \ w \in B$}\\
                 {}\\
                 i_B^X (z)  =_X i_B (w)  &\mbox{, $z, w \in B$}\\
                 {}\\
                 i_B^X (z)  =_X i_A (w)  &\mbox{, $z \in B \ \& \ w \in A$,}
                 \end{array}
          \right.$$
we can define afterwords the \textit{operation} $i_{A \cup B}^X : A \cup B \to X$ as in 
Definition~\ref{def: unionofsubsets}.
In this way the non-dependent assignment routine $i_{A \cup B}^X$ is defined on a set, and it is an operation. 
Bishop avoids this definition, probably because this pattern cannot be extended to the definition of a union of a family of 
subsets (see Definition~\ref{def: interiorunion}). In that case,
we cannot write down the corresponding case distinction for $z =_{A \cup B} w$. Moreover, 
the proof of $\big(A \cup B, i_{A \cup B}^X\big) \subseteq X$ is immediate, if one
uses Definition~\ref{def: unionofsubsets}. 
\end{note}

\begin{note}\label{not: emptyset}
\normalfont
The definition of the empty subset $\emptyset_X$ of a set $X$, given in~\cite{Bi67}, p.~65, can be 
formulated as follows. Let $X$ be a set and $x_0 \in X$. The totality $\emptyset_X$ is defined by
$z \in \emptyset_X : \TOT x_0 \in X \ \& \ 0 =_{\Nat} 1$.
Let $i_{\emptyset}^X \colon \emptyset_X \sto X$ be the non-dependent assignment routine, defined by $i(z) := x_0$, 
for every $z \in \emptyset_X$, and let $z =_{\emptyset_X} w : \TOT i(z) =_X i(w) : \TOT x_0 =_X x_0.$
The pair $(\emptyset_X, i_{\emptyset}^X)$ is the \textit{empty subset}\index{empty subset of a set} of $X$.
One can show that $=_{\emptyset_X}$ is an equality on $\emptyset_X$, and hence $\emptyset_X$ can be 
considered to be a set. The assignment routine $i_{\emptyset}^X$ is an embedding of $\emptyset_X$ into $X$, and hence 
$(\emptyset_X, i_{\emptyset}^X)$ is a subset of $X$. As Bishop himself writes in~\cite{Bi67}, p.~65, ``the definition
of $\emptyset$ is negativistic, and we prefer to mention the void set as seldom as possible''.
In~\cite{BB85}, p.~69, Bishop and Bridges define two subsets $A, B$ of $X$ to be disjoint, when 
$A \cap B$ ``is the void subset of $X$''. 
Clearly, this ``is'' cannot be $A \cap B := \emptyset_X$. If we interpret it as
$A \cap B =_{\C P(X)} \emptyset_X$, we need the existence of certain functions 
from $\emptyset_X$ to $A \cap B$ and from $A \cap B$ to $\emptyset_X$. The latter approach is 
followed in $\MLTT$ for the empty type. Following Bishop, we refrain from elaborating this negatively defined notion.

\end{note}

\begin{note}\label{not: defextimage}
\normalfont
If $(A, i_A^X) \subseteq A$, $(B, i_B^Y) \subseteq Y$, and $\fXY$, the 
\textit{extensional image}\index{extensional image} $f[A]$\index{$f[A]$} of 
$A$ under $f$ is defined through the extensional property
$P(y) := \exists_{a \in A}\big(f(i_A(a)) =_Y y\big)$.
Similarly, the \textit{extensional pre-image}\index{extensional pre-image} $f^{-1}[B]$\index{$f^{-1}[B]$} of 
$B$ under $f$ is defined through the extensional property
$Q(x) := \exists_{b \in B}\big(f(x) =_Y i_B (b)\big)$.
The subset $f(A)$ of $Y$ contains exactly the outputs $f(i_A^X(a))$ of $f$, for every $a \in A$, while 
the subset $f[A]$ of $Y$ contains all the elements of $Y$ that are $=_Y$-equal to some output  $f(i_A(a))$ 
of $f$, for every $a \in A$. It is useful to keep the ``distinction'' between the subsets  $f(A)$, $f[A]$, and
$f^{-1}(B)$, $f^{-1}[B]$.
We need the equality in $\C P(X)$ of a subset of $X$ to its extensional version (see Note~\ref{not: LANC}),
hence the principle $\LANC$, to get
$f(A) =_{\C P(Y)} f[A]$ and $f^{-1}(B) =_{\C P(X)} f^{-1}[B]$.
\end{note}

\begin{note}\label{not: powersetisset}
\normalfont
There are instances in Bishop's work indicating that the powerset of a set is treated as a set.
In~\cite{Bi67}, p.~68, and in~\cite{BB85}, p.~74, 
the following ``function'' is defined  
\[j : \C P^{\Disj} (X) \to \C P(X), \ \ \ \ (A^1, A^0) \mapsto A^1. \]
This is in complete contrast to our interpretation of a function as an operation between sets. Of course,
such a rule is an exception in~\cite{Bi67} and~\cite{BB85}. 
In the definition of an integration space, see~\cite{BB85}, p.~216, the ``set'' $\C F(X,Y)$ of all
strongly extensional partial functions from $X$ to $Y$ requires quantification over $\D V_0$.
Such a quantification is also implicit in the definition of a measure space given in~\cite{BB85}, p.~282, 
and in the definition of a complete measure space
in~\cite{BB85}, p.~289. These definitions appeared first in~\cite{BC72}, p.~47, and p.~55, respectively.
The powerset is repeatedly used as a set in~\cite{BR87} and~\cite{MRR88}. It is not known if the treatment 
of the powerset as a set implies some constructively unacceptable principle. 
\end{note}

\begin{note}\label{not: powersetnotset}
\normalfont
There are instances in Bishop's work indicating that the powerset of a set is \textit{not} treated as a set.
See e.g., the definition of a set-indexed family of sets in~\cite{BB85}, p.~78 (our Definition~\ref{def: famofsets}). 
Similarly, in the definition of a family of subsets of a set $A$ indexed 
by some set $T$ (see~\cite{BB85}, p.~69), the notion of a finite routine that assigns a 
subset of $A$ to an element of $T$ is used, and not the notion of a function from $T$ to $\C P (A)$. 
In the definition of a measure space in~\cite{Bi67}, p.~183, a subfamily of a given family of 
complemented sets is considered in order to avoid quantification over the class of all complemented subsets in the 
formulations of the definitional clauses of a measure space (see Note~\ref{not: Bms}).
The powerset axiom is also avoided in Myhill's formalization~\cite{My75} of $\BISH$ and in Aczel's subsequent system
$\CZF$ of constructive set theory (see~\cite{AR10}). Although, as we said, it is not known if the use
of the powerset as a set implies some constructively unacceptable principle, it is not accepted in
any \textit{predicative} development of constructive mathematics. 
\end{note}

\begin{note}\label{not: onpartialfunctions}
\normalfont
The notion of a partial function was introduced by Bishop and Cheng in~\cite{BC72}, p.~1, and this definition, 
together with the introduced term ``partial function'', was also included in Chapter 3 
of~\cite{BB85}, p.~71. The totality of partial functions $\C F(X)$ from a set $X$ to $\Real$ is crucial to the 
definition of an integration space in the new measure theory developed in~\cite{BC72}, and seriously 
extended in~\cite{BB85}. Only the basic algebraic operations on $\C F(X)$ were defined
in~\cite{BB85}, p.~71. The composition of partial functions is mentioned in~\cite{CDPS05},
pp.~66--67.
A notion of a partial dependent operation can be defined as follows. If $A, I$ are sets, a partial dependent 
operation\index{partial dependent operation} is a triplet $(A, i_A^I, \Phi_A^{\lambda_0})$, where
$(A, i_A) \subseteq I$, $\lambda_0 \colon A \sto \D V_0$, and $\Phi_A^{\lambda_0} \colon 
\bigcurlywedge_{a \in A}\lambda_0(a)$. If $\lambda_0(a) := Y$, for every $a \in A$, then the corresponding 
partial dependent operation is reduced to a partial function in $\C F(I, Y)$.

\end{note}

\begin{note}\label{not: notioncomplemented}
\normalfont
In the study of various subsets of a set $X$ we avoided to define the complement of a subset,
since this requires a negative definition. Recall that the negatively defined notion of empty subset of a set
is not really used.  In~\cite{Bi67} Bishop introduced a 
positive notion of the complement of a subset of a set $X$, the notion of a complemented subset of $X$. For
its definition we need a notion of a fixed inequality on $X$, which is compatible with the given equality of 
$X$. In this way we can express the disjointness of two subsets $A, B$ of a set $X$ in a positive way. Usually, 
$A, B$ are called \textit{disjoint}\index{disjoint subsets}, if $A \cap B$ is not inhabited. 
It is computationally more informative though, if a positive way is found to express disjointness of subsets. 
In~\cite{BV11} a positive notion of \textit{apartness} is used as a foundation of constructive topology.

\end{note}

\begin{note}\label{not: onoperationscomplemented67}
\normalfont
The definitions of $\B A \cap \B B, \B A \cup \B B$ and $\B A - \B B$ appear in~\cite{Bi67}, p.~66, 
where $\B A \cup \B B$ and $\B A \cap \B B$ are special cases of the complemented subsets 
$\bigcup_{i \in I}\B \lambda_0(i)$ and $\bigcap_{i \in I}\B \lambda_0(i)$, respectively (see Proposition~\ref{prp: unioncompl}). There the inequality on $X$ is induced by an inhabited set of functions from $X$ to $\Real$. The
definition of $\B A \times \B C$ appears in~\cite{Bi67}, p.~206, in the section of the product measures.
One can motivate these definitions applying a ``classical'' thinking. If $x \in X$, recall the definitions
\[ x \in \B A :\TOT x \in A^1 \ \ \ \& \ \ \ x \notin \B A :\TOT x \in A^0. \]
Interpreting the connectives in a classical way, we get 
\[ x \in \B A \cup \B B \TOT x \in \B A  \ \vee  \ x \in \B B :\TOT x \in A^1 \ \vee  \ x \in B^1 
:\TOT x \in A^1 \cup B^1, \]
\[ x \notin \B A \cup \B B \TOT x \notin \B A \ \& \ x \notin \B B :\TOT x \in A^0 \ \& \ x \in B^1 
:\TOT x \in A^1 \cap B^1, \]
\[ x \in \B A \cap \B B \TOT x \in \B A  \ \&  \ x \in \B B :\TOT x \in A^1 \ \&  \ x \in B^1 
:\TOT x \in A^1 \cap B^1, \]
\[ x \notin \B A \cap \B B \TOT x \notin \B A \ \vee \ x \notin \B B :\TOT x \in A^0 \ \vee \ x \in B^1 
:\TOT x \in A^1 \cup B^1, \]
\[ x \in - \B A \TOT x \notin \B A :\TOT x \in A^0 \ \ \ \ \& \ \ \ \ x \notin - \B A \TOT x \in \B A 
:\TOT x \in A^1, \]
\[ (x,y) \in \B A \times \B C \TOT x \in \B A \ \& \ y \in \B C :\TOT x \in A^1 \ \& \ y \in B^1 :\TOT 
(x,y) \in A^1 \times B^1,\]
\[ (x,y) \notin \B A \times \B C \TOT x \notin \B A \ \vee \ y \notin \B C :\TOT x \in A^0 \ \vee \ y 
\in B^0 :\TOT (x,y) \in (A^0 \times Y) \cup (X \times B^0).\]

\end{note}

\begin{note}\label{not: onoperationscomplemented72}
\normalfont
In~\cite{BC72}, pp.~16--17, and in~\cite{BB85}, p.~73, 
the operations between the complemented subsets of a set $X$ follow Definition~\ref{def: operationscomplemented2}
in order to employ the good behaviour of the corresponding characteristic functions in the new measure theory.
In the measure theory of~\cite{Bi67}, where the characteristic functions of complemented subsets are not crucial,
the operations between complemented subsets are defined according to Definition~\ref{def: operationscomplemented1}. 
Bishop and Cheng use the notation $\B A \times \B B$ instead of $\B A \otimes \B B$.
As it is evident from the previous figures, the $1$- and $0$-components of the complemented subsets in the Bishop-Cheng definition are subsets of the corresponding $1$- and $0$-components of the complemented subsets in the Bishop definition from~\cite{Bi67}.  Actually, the definitions of the operations of complemented subsets 
in~\cite{Bi67} associate to the $1$-component of the complemented subset a maximal complement. The two sets of operations though, share the same algebraic and set-theoretic properties. They only behave differently
with respect to their characteristic functions.
Based on the work~\cite{Sh18} of Shulman, we can motivate the second set of operations in a way similar to
the motivation provided for the first set of operations in Note~\ref{not: onoperationscomplemented67}. 
Keeping the definitions of $x \in \B A$ and $x \notin \B B$, we can apply a ``linear'' 
interpretation of the connectives $\vee$ and $\&$. As it is mentioned in~\cite{Sh18}, p.~2, the multiplicative
version $P \ \parl \ Q$\index{$P \ \parl \ Q$} of $P \vee Q$ in linear logic  represents the pattern
``if not $P$, then $Q$; and if not $Q$, then $P$''. Let 
\[ x \in \B A \vee \B B  :\TOT [x \notin \B A \To x \in \B B] \ \& \ [x \notin \B B \To x \in \B A]. \]
With the use of Ex falsum quodlibet the implication $x \notin \B A \To x \in \B B$ holds if $x \in \B A :\TOT x \in A^1$, or 
if $x \notin \B A :\TOT x \in A^0$ and $x \in \B B :\TOT x \in B^1$ i.e., if $x \in A^0 \cap B^1$. Hence, 
the first implication holds if $x \in A^1 \cup (A^0 \cap B^1)$. Similarly, the second holds if 
$x \in B^1 \cup (B^0 \cap A^1)$. Thus
\[ x \in \B A \vee \B B \TOT x \in [ A^1 \cup (A^0 \cap B^1)] \cap [B^1 \cup (B^0 \cap A^1)], \]
and the last intersection is equal to $\Dm(\B A \vee \B B)$! One then can define $x \notin \B A \vee \B B 
:\TOT x \notin \B A \ \& \ x \notin \B B$, and $x \in \B A \wedge \B B :\TOT x \in \B A \ \& \ x \in \B B$, 
and $x \notin \B A \wedge \B B :\TOT x \in (- \B A) \vee (- \B B)$. \\[1mm]
For the relation of complemented subsets to the Chu construction see~\cite{Pe21c}.
\end{note}

\chapter{Families of sets}
\label{chapter: familiesofsets}

We develop the basic theory of set-indexed families of sets and of family-maps between them. 
We study the exterior union of a family of sets $\Lambda$, or the $\sum$-set of $\Lambda$, and the set
of dependent functions over $\Lambda$, or the $\prod$-set of $\Lambda$. We prove the distributivity of
$\prod$ over $\sum$ for families of sets indexed by a product of sets, which is the translation 
of the type-theoretic axiom of choice into $\BST$. Sets of sets are special set-indexed families of sets 
that allow ``lifting'' of functions on the index-set to functions on them. The direct families of sets and the 
set-relevant families of sets are introduced. The index-set of the former is a directed set, while the 
transport maps of the latter are more than one and appropriately indexed. With the use of the introduced 
universe $\D V_0^{\im}$ of sets and impredicative sets we study families of families of sets.

\section{Set-indexed families of sets}
\label{sec: famsets}

Roughly speaking, a family of sets indexed by some set $I$ is an
assignment routine $\lambda_0 : I \sto \D V_0$ that 
behaves like a function i.e., if $i =_I j$, then $\lambda_0(i) =_{\D V_0} \lambda_0 (j)$. Next follows  
an exact formulation of this description that reveals the witnesses of the 
equality $\lambda_0(i) =_{\D V_0} \lambda_0 (j)$.

\begin{definition}\label{def: famofsets}
If $I$ is a set,
a \textit{family of sets}\index{family of sets} indexed by $I$, or an $I$-\textit{family 
of sets}\index{$I$-family of sets}, is a pair $\Lambda := (\lambda_0, \lambda_1)$, where\index{$\lambda_0$}
$\lambda_0 \colon I \sto \D V_0$, and\index{$\lambda_1$} $\lambda_1$, a
\textit{modulus of function-likeness for}\index{modulus of function-likeness} $\lambda_0$, is given by
\[ \lambda_1 \colon \bigcurlywedge_{(i, j) \in D(I)}\D F\big(\lambda_0(i), \lambda_0(j)\big), 
\ \ \ \lambda_1(i, j) := \lambda_{ij}, \ \ \ (i, j) \in D(I), \]
such that the \textit{transport maps}\index{transport map of a family of sets} $\lambda_{ij}$
\index{transport map of a family of sets}\index{$\lambda_{ij}$}
of $\Lambda$ satisfy the following conditions:\\[1mm]
\normalfont (a) 
\itshape For every $i \in I$, we have that $\lambda_{ii} := \id_{\lambda_0(i)}$.\\[1mm]
\normalfont (b) 
\itshape If $i =_I j$ and $j =_I k$, the following diagram commutes
\begin{center}
\begin{tikzpicture}

\node (E) at (0,0) {$\lambda_0(j)$};
\node[right=of E] (F) {$\lambda_0(k).$};
\node [above=of E] (D) {$\lambda_0(i)$};

\draw[->] (E)--(F) node [midway,below] {$\lambda_{jk}$};
\draw[->] (D)--(E) node [midway,left] {$\lambda_{ij}$};
\draw[->] (D)--(F) node [midway,right] {$\ \lambda_{ik}$};

\end{tikzpicture}
\end{center}
$I$ is the \textit{index-set}\index{index-set} of the family $\Lambda$.
If $X$ is a set, the \textit{constant} $I$-\textit{family of sets}\index{constant family of sets}
$X$\index{$C^X$} is the pair 
$C^X := (\lambda_0^X, \lambda_1^X)$, where $\lambda_0 (i) := X$, for every 
$i \in I$, and $\lambda_1 (i, j) := \id_X$, for every $(i, j) \in D(I)$ $($see the left diagram in
Definition~\ref{def: basicfamilies}$)$.
\end{definition}

The dependent operation $\lambda_1$ should have been written as follows 
\[ \lambda_1 \colon \bigcurlywedge_{z \in D(I)}\D F\big(\lambda_0(\prb_1(z)), \lambda_0(\prb_2(z))\big), \]
but, for simplicity, we avoid the use of the primitive projections $\prb_1, \prb_2$. Condition (a) of
Definition~\ref{def: famofsets} could have been written as $\lambda_{ii} =_{\mathsmaller{\D F(\lambda_0(i), 
\lambda_0(i))}} \id_{\lambda_0(i)}$. If $i =_I j$, then by conditions (b) and (a) of Definition~\ref{def: famofsets} 
we get 
$id_{\lambda_0(i)} := \lambda_{ii} = \lambda_{ji} \circ \lambda_{ij}$ and 
$ \id_{\lambda_0(j)} := \lambda_{jj} = \lambda_{ij} \circ \lambda_{ji} $
i.e., $(\lambda_{ij}, \lambda_{ji}) \colon \lambda_0 (i) =_{\D V_0} \lambda_0 (j)$. In this sense $\lambda_1$ is 
a modulus of function-likeness for $\lambda_0$.

\begin{definition}\label{def: basicfamilies}
The pair \index{$\Lambda^{\D 2}$}$\Lambda^{\D 2} := (\lambda_0^{\D 2}, \lambda_1^{\D 2})$, where $\lambda_0^{\D 2} \colon 
\D 2 \sto \D V_0$ with $\lambda_0^{\D 2} (0) := X$, $\lambda_0^{\D 2} (1) := Y$, and
$\lambda_1^{\D 2} (0, 0) := \id_X$ and $\lambda_1^{\D 2} (1, 1) := \id_Y$,
is the $\D 2$-family of $X$ and $Y$\index{$\D 2$-family of $X$ and $Y$}
\begin{center}
\begin{tikzpicture}

\node (E) at (0,0) {$X$};
\node[right=of E] (F) {$X$};
\node [above=of E] (D) {$X$};
\node[right=of F] (K) {$Y$};
\node [above=of K] (L) {$Y$};
\node [right=of K] (M) {$Y$.};

\draw[->] (E)--(F) node [midway,below] {$\id_X$};
\draw[->] (D)--(E) node [midway,left] {$\id_X$};
\draw[->] (D)--(F) node [midway,right] {$\ \id_X$};
\draw[->] (K)--(M) node [midway,below] {$\id_Y$};
\draw[->] (L)--(K) node [midway,left] {$\id_Y$};
\draw[->] (L)--(M) node [midway,right] {$\ \id_X$};

\end{tikzpicture}
\end{center}
The \index{$\Lambda^{\D n}$}$\D n$-family $\Lambda^{\D n}$\index{$\Lambda^{\D n}$} of the sets $X_1, \ldots X_n$, 
where $n \geq 1$,\index{$\D n$-family of $X_1, \ldots, X_n$} 
and the $\Nat$-family\index{$\Lambda^{\Nat}$} $\Lambda^{\Nat} := (\lambda_0^{\Nat}, \lambda_1^{\Nat})$
of the sets $(X_n)_{n \in \Nat}$\index{$\Nat$-family of $(X_n)_{n \in \Nat}$}
are defined similarly\footnote{It is immediate to show that $\Lambda^{\D n}$ is an $\D n$-family, and $\Lambda^{\Nat}$ 
is an $\Nat$-family.}.
\end{definition}

\begin{definition}\label{def: map}
Let $\Lambda := (\lambda_0, \lambda_1)$ and $M := (\mu_0, \mu_1)$ be $I$-families of sets.
A \textit{family-map} from $\Lambda$ to $M$\index{family-map}, in symbols $\Psi \colon \Lambda
\To M$\index{$\Psi \colon \Lambda \To M$} is a dependent operation
$\Psi \colon \bigcurlywedge_{i \in I}\D F \big(\lambda_0(i), \mu_0(i)\big)$ 
such that for every $(i, j) \in D(I)$ the following diagram commutes
\begin{center}
\begin{tikzpicture}

\node (E) at (0,0) {$\mu_0(i)$};
\node[right=of E] (F) {$\mu_0(j)$.};
\node[above=of F] (A) {$\lambda_0(j)$};
\node [above=of E] (D) {$\lambda_0(i)$};

\draw[->] (E)--(F) node [midway,below]{$\mu_{ij}$};
\draw[->] (D)--(A) node [midway,above] {$\lambda_{ij}$};
\draw[->] (D)--(E) node [midway,left] {$\Psi_i$};
\draw[->] (A)--(F) node [midway,right] {$\Psi_j$};

\end{tikzpicture}
\end{center}
Let $\Map_I(\Lambda, M)$\index{$\Map_I(\Lambda, M)$} be the totality of family-maps\index{totality of family-maps} 
from $\Lambda$ to $M$, which is equipped with the equality 
\[ \Psi =_{\Map_I(\Lambda, M)} \Xi :\TOT \forall_{i \in I}\big(\Psi_i =_{\D F(\lambda_0(i), \mu_0(i))} 
\Xi_i\big). \]
If $\Xi : M \To N$, the \textit{composition family-map}\index{composition family-map} 
$\Xi \circ \Psi \colon \Lambda \To N$ is defined, for every $i \in I$, by 
$(\Xi \circ \Psi)_i := \Xi_i \circ \Psi_i$
\begin{center}
\begin{tikzpicture}

\node (E) at (0,0) {$\lambda_0(i)$};
\node[right=of E] (F) {$\lambda_0(j)$};
\node[below=of F] (A) {$\mu_0(j)$};
\node[below=of E] (B) {$\mu_0(i)$};
\node[below=of B] (K) {$\nu_0(i)$};
\node[below=of A] (L) {$\nu_0(j)$.};

\draw[->] (E)--(F) node [midway,above] {$\lambda_{ij}$};
\draw[->] (F)--(A) node [midway,left] {$\Psi_j$};
\draw[->] (B)--(A) node [midway,below] {$\mu_{ij}$};
\draw[->] (E) to node [midway,right] {$\Psi_i$} (B);
\draw[->] (B) to node [midway,right] {$\Xi_i$} (K);
\draw[->] (A)--(L) node [midway,left] {$\Xi_j$};
\draw[->] (K)--(L) node [midway,below] {$\nu_{ij}$};
\draw[->,bend right] (E) to node [midway,left] {$(\Xi \circ \Psi)_i  \ $} (K);
\draw[->,bend left] (F) to node [midway,right] {$\  (\Xi \circ \Psi)_j$} (L);

\end{tikzpicture}
\end{center} 
The \textit{identity family-map}\index{identity family-map}\index{$\Id_{\Lambda}$} 
$Id_{\Lambda} \colon \bigcurlywedge_{i \in I}\D F 
\big(\lambda_0(i), \lambda_0(i)\big)$on $\Lambda$, is defined by $\Id_{\Lambda}(i) := \id_{\lambda_0(i)}$, 
for every $i \in I$.
Let $\Fam(I)$\index{$\Fam(I)$} be the totality of $I$-families\index{totality of $I$-families},
equipped with the canonical equality
\[ \Lambda =_{\Fam(I)} M :\TOT 
\exists_{\Phi \in \Map_I(\Lambda, M)}\exists_{\Xi \in \Map_I(M, \Lambda)}\big((\Phi, \Xi) \colon \Lambda =_{\Fam(I)} M\big), \]
\[ (\Phi, \Xi) \colon \Lambda =_{\Fam(I)} M :\TOT \big(\Phi \circ \Xi = \id_M \ \& \
\Xi \circ \Phi = \id_{\Lambda}\big). \]

\end{definition}

It is straightforward to show that the composition family-map $\Xi \circ \Psi$ is a family-map from $\Lambda$ to $N$, 
and that the equalities on $\Map_I(\Lambda, M)$ and $\Fam(I)$ satisfy the conditions of an 
equivalence relation. It is natural to accept the totality $\Map(\Lambda, M)$ as a set. If $\Fam(I)$ was a set
though, the constant $I$-family with value $\Fam(I)$ would be
defined though a totality in which it belongs to. From a predicative point of view, this cannot be accepted. The 
membership condition of the totality
$\Fam(I)$ though, does not depend on the universe $\D V_0$, therefore it is also natural not to consider $\Fam(I)$ to be a class. Hence, $\Fam(I)$ is a totality ``between'' a (predicative) set and a class. For this reason, we say that $\Fam(I)$ is an 
\textit{impredicative set}\index{impredicative set}. Next follows an obvious generalisation of a family-map.

\begin{definition}\label{def: prfeqfamofsets}
If $\Lambda, M \in \Fam(I)$, such that $\Lambda =_{\Fam(I)} M$, we define the set 
\[ \Eq(\Lambda, M) :=  \big\{(\Phi, \Psi) \in \Map_I(\Lambda, M) \times \Map_I(M, \Lambda) \mid (\Phi, \Psi) : 
\Lambda =_{\Fam(I)} M\big\}, \]
equipped with the equality of the product of sets. If $\Phi \in \Map_I(\Lambda, M), \Psi \in \Map_I(M, \Lambda), \Phi{'} \in \Map_I(M, N)$ and $\Psi{'} \in \Map_I(N, M)$, let
$\refl(\Lambda) := \big(\Id_{\Lambda}, \Id_{\Lambda}\big)$ and $(\Phi, \Psi)^{-1} := (\Psi, \Phi)$ and  
$(\Phi, \Psi) \ast (\Phi{'}, \Psi{'}) := (\Phi{'} \circ \Phi, \Psi \circ \Psi{'})$.

\end{definition}

As in the case of $\D V_0$ and the corresponding set $\Eq(X, Y)$,
in general, not all elements of $\Eq(\Lambda, M)$ are equal. If $I := \D 1 := \{0\}$, and $\lambda_0(0) := \D 2$, 
and if $\Phi_0
:= \id_{\D 2}$ and $\Psi_0 := \sw_{\D 2}$, 
then $(\Phi, \Phi) \in \Eq(\Lambda, \Lambda)$ and
$(\Psi, \Psi) \in \Eq(\Lambda, \Lambda)$, while $\Phi \neq_{\Map_I(\Lambda, \Lambda)} \Psi$, since 
$\Phi_0 \neq_{\D F(\D 2, \D 2)} \Psi_0$. 
It is immediate to show the groupoid-properties (i)-(iv) for the equality of the totality $\Fam(I)$.

\begin{definition}\label{def: hmap}
Let $I, J, K$ be sets, $h \in \D F(J, I), g \in \D F(K, J)$, 
$\Lambda := (\lambda_0, \lambda_1) \in \Fam(I)$, $M := (\mu_0, \mu_1) \in \Fam(J)$, and 
$N := (\nu_0, \nu_1) \in \Fam(K)$. 
A \textit{family-map} from $M$ to $\Lambda$ \textit{over} $h$\index{family-map over a function}
is a dependent operation $\Psi \colon \bigcurlywedge_{j \in J}\D F \big(\mu_0(j), \lambda_0(h(j))\big)$,
such that for every $(j, j{'}) \in D(J)$ the following diagram commutes
\begin{center}
\begin{tikzpicture}

\node (E) at (0,0) {$\lambda_0(h(j))$};
\node[right=of E] (F) {$\lambda_0(h(j{'}))$,};
\node[above=of F] (A) {$\mu_0(j{'})$};
\node [above=of E] (D) {$\mu_0(j)$};

\draw[->] (E)--(F) node [midway,below]{$\mathsmaller{\lambda_{h(j)h(j{'})}}$};
\draw[->] (D)--(A) node [midway,above] {$\mu_{jj{'}}$};
\draw[->] (D)--(E) node [midway,left] {$\Psi_j$};
\draw[->] (A)--(F) node [midway,right] {$\Psi_{j{'}}$};

\end{tikzpicture}
\end{center}
where $\Psi_j := \Psi(j)$ is the $j$-\textit{component} of $\Psi$, for every $j \in J$. We write $\Psi \colon M 
\stackrel{h} \To \Lambda$ for such a family-map. If $\Psi \colon M \stackrel{h} \To \Lambda$ and 
$\Xi \colon N \stackrel{g} \To M$, the \textit{composition family-map} $\Psi \circ \Xi \colon N \stackrel{h \circ g}
\Longrightarrow \Lambda$ \textit{over} $h \circ g$\index{composition family-map over
a composite function}  is defined, for every $k \in K$, by 
$(\Psi \circ \Xi)_k := \Psi_{g(k)} \circ \Xi_k$
\begin{center}
\begin{tikzpicture}

\node (E) at (0,0) {$\nu_0(k)$};
\node[right=of E] (X) {};
\node[right=of X] (F) {$\nu_0(k{'})$};
\node[below=of F] (A) {$\mathsmaller{\mu_0(g(k))}$};
\node[below=of E] (B) {$\mathsmaller{\mu_0(g(k{'}))}$};
\node[below=of B] (K) {$\mathsmaller{\lambda_0(h(g(k)))}$};
\node[below=of A] (L) {$\mathsmaller{\lambda_0(h(g(k{'})))}$.};

\draw[->] (E)--(F) node [midway,above] {$\nu_{kk{'}}$};
\draw[->] (F)--(A) node [midway,left] {$\Xi_{k{'}}$};
\draw[->] (B)--(A) node [midway,below] {$\mathsmaller{\mu_{g(k)g(k{'})}}$};
\draw[->] (E) to node [midway,right] {$\Xi_k$} (B);
\draw[->] (B) to node [midway,right] {$\Psi_{g(k)}$} (K);
\draw[->] (A)--(L) node [midway,left] {$\Psi_{g(k{'})}$};
\draw[->] (K)--(L) node [midway,below] {$\mathsmaller{\lambda_{h(g(k))h(g(k{'}))}}$};
\draw[->,bend right=50] (E) to node [midway,left] {$(\Psi \circ \Xi)_k  \ $} (K);
\draw[->,bend left=50] (F) to node [midway,right] {$\  (\Psi \circ \Xi)_{k{'}}$} (L);

\end{tikzpicture}
\end{center} 

\end{definition}

\begin{definition}\label{def: newfamsofsets1}
Let $\Lambda := (\lambda_0, \lambda_1), M := (\mu_0, \mu_1)$ be $I$-families of sets.\\[1mm] 
\normalfont (i) 
\itshape
The \textit{product family}\index{product of $I$-families} of $\Lambda$ and $M$ is the pair\index{$\Lambda \times M$} 
$\Lambda \times M := (\lambda_0 \times \mu_0, \lambda_1 \times \mu_1)$, where 
\[ (\lambda_0 \times \mu_0)(i) := \lambda_0 (i) \times \mu_0 (i); \ \ \ \ i \in I, \]
\[ \big(\lambda_1 \times \mu_1\big)_{ij} \colon \lambda_0 (i) \times \mu_0 (i) \to \lambda_0 (j) \times \mu_0 (j);
\ \ \ \ 
(i, j) \in D(I), \]
\[ \big(\lambda_1 \times \mu_1\big)_{ij}\big(x, y\big) := \big(\lambda_{ij}(x), \mu_{ij}(y)\big); \ \ \ \ 
x \in \lambda_0 (i) \ \& \ y \in \mu_0 (i). \]
\normalfont (ii) 
\itshape
The \textit{function space family}\index{function space
family} from $\Lambda$ to $M$ is the pair\index{$\D F(\Lambda, M)$}
$\D F(\Lambda, M) := \big(\D F(\lambda_0, \mu_0), \D F(\lambda_1, \mu_1)\big)$
where
\[ \big[\D F(\lambda_0, \mu_0)\big](i) := \D F\big(\lambda_0(i), \mu_0(i)\big); \ \ \ \ i \in I, \]
\[ \D F(\lambda_1, \mu_1) \ \colon \bigcurlywedge_{(i, j) \in D(I)}\D F\bigg(\D F\big(\lambda_0(i), \mu_0(i)\big), 
 \D F\big(\lambda_0(j), \mu_0(j)\big)\bigg) \]
 \[ \D F(\lambda_1, \mu_1)_{ij} := \D F(\lambda_1, \mu_1)(i, j) \colon 
 \D F\big(\lambda_0(i), \mu_0(i)\big) \to \D F\big(\lambda_0(j), \mu_0(j)\big); \ \ \ \ (i, j) \in D(I), \]
\[ \D F(\lambda_1, \mu_1)_{ij}(f) := \mu_{ij} \circ f \circ \lambda_{ji} \]
\begin{center}
\begin{tikzpicture}

\node (E) at (0,0) {$\lambda_0(j)$};
\node[right=of E] (H) {};
\node[right=of H] (F) {$\mu_0(j)$.};
\node[above=of F] (A) {$\mu_0(i)$};
\node [above=of E] (D) {$\lambda_0(i)$};

\draw[->] (E)--(F) node [midway,below]{$\mathsmaller{\D F(\lambda_1, \mu_1)_{ij}(f)}$};
\draw[->] (D)--(A) node [midway,above] {$f$};
\draw[->] (E)--(D) node [midway,left] {$\lambda_{ji}$};
\draw[->] (A)--(F) node [midway,right] {$\mu_{ij}$};

\end{tikzpicture}
\end{center}
\normalfont (iii) 
\itshape
If $K$ is a set, $\Sigma := (\sigma_0, \sigma_1)$ is a $K$-family of sets and $h : I \to K$, 
the \textit{composition family}\index{composition family}
of $\Sigma$ with $h$ is the pair\index{$\Sigma \circ h$} $\Sigma \circ h :=
(\sigma_0 \circ h, \sigma_1 \circ h)$, where
\[ (\sigma_0 \circ h)(i) := \sigma_0 (h(i)); \ \ \ \ i \in I, \] 
\[ (\sigma_1 \circ h)_{ij} := (\sigma_1 \circ h)(i, j) \colon \sigma_0 (h(i)) \to \sigma_0 (h(j)); \ \ \ \ 
(i, j) \in D(I), \]
\[ (\sigma_1 \circ h)_{ij} := \sigma_{h(i)h(j)}. \]
\end{definition}

It is straightforward to show that $\Lambda \times M$, $\D F(\Lambda, M)$, and $\Sigma \circ h$ are $I$-families.
E.g., for $\D F(\Lambda, M)$,
and if $i, j, k \in I$ and $i =_I j =_I k$, we have that
\[\D F(\lambda_1, \mu_1)_{ii}(f) := \mu_{ii} \circ f \circ \lambda_{ii} := \id_{\mu_0(i)} \circ f \circ 
\id_{\lambda_0(i)} := f, \]
\begin{align*}
 \D F(\lambda_1, \mu_1)_{jk}\bigg(\D F(\lambda_1, \mu_1)_{ij}(f)\bigg) & := 
 \mu_{jk} \circ \big[\mu_{ij} \circ f \circ \lambda_{ji}\big] \circ \lambda_{kj}\\
 & := \big[\mu_{jk} \circ \mu_{ij}\big] \circ f \circ \big[\lambda_{ji} \circ \lambda_{kj}\big]\\
 & = \mu_{ik} \circ f \circ \lambda_{ki}\\
 & := \D F(\lambda_1, \mu_1)_{ik}(f).
\end{align*}

\begin{proposition}\label{prp: constequal7}
Let $X, Y, I$ be sets and $C^X, C^Y, C^{X \times Y}, C^{\D F(X, Y)}$ the constant $I$-families $X, Y, X \times Y$, and
$\D F(X, Y)$, respectively. \\[1mm]
\normalfont (i) 
\itshape $C^X \times C^Y =_{\Fam(I)} C^{X \times Y}$.\\[1mm]
\normalfont (ii) 
\itshape $\D F(C^X, C^Y) =_{\Fam(I)} C^{\D F(X, Y)}$.
\end{proposition}

\begin{proof}
 (i) Let $\Phi \colon C^X \times C^Y \To  C^{X \times Y}$ and $\Psi \colon  C^{X \times Y} \To C^X \times C^Y$
 be defined by
 $\Phi_i := X \times Y := \Psi_i$, for every $i \in I$, then by the commutativity of the following left diagram
 
 \begin{center}
\begin{tikzpicture}

\node (E) at (0,0) {$X \times Y$};
\node[right=of E] (X) {};
\node[right=of X] (F) {$X \times Y$};
\node[above=of F] (A) {$X \times Y$};
\node [above=of E] (D) {$X \times Y$};
\node [right=of F] (P) {};
\node [right=of P] (G) {$\D F(X, Y)$};
\node [right=of G] (Z) {};
\node [right=of Z] (H) {$\D F(X, Y)$,};
\node[above=of H] (K) {$\D F(X, Y)$};
\node [above=of G] (L) {$\D F(X, Y)$};

\draw[->] (E)--(F) node [midway,below]{$\lambda_{ij}$};
\draw[->] (D)--(A) node [midway,above] {$(\lambda^X_1 \times \mu_1^Y)_{ij}$};
\draw[->] (D)--(E) node [midway,left] {$\id_{X \times Y}$};
\draw[->] (A)--(F) node [midway,right] {$\id_{X \times Y}$};
\draw[->] (G)--(H) node [midway,below] {$\mu_{ij}$};
\draw[->] (L)--(K) node [midway,above] {$\D F(\lambda_1^X, \mu_1^Y)_{ij}$};
\draw[->] (L)--(G) node [midway,left] {$\id_{\D F(X, Y)}$};
\draw[->] (K)--(H) node [midway,right] {$\id_{\D F(X, Y)}$};

\end{tikzpicture}
\end{center}
$\Phi, \Psi$ are well-defined family-maps  and $(\Phi, \Psi) \colon C^X \times C^Y =_{\Fam(I)} C^{X \times Y}$.\\
(ii) Let $\Phi \colon \D F(C^X, C^Y) \To  C^{\D F(X, Y)}$ and $\Psi \colon  C^{\D F(X, Y)} \To \D F(C^X, C^Y)$
be defined by
 $\Phi_i := \D F(X, Y) := \Psi_i$, for every $i \in I$, then by the commutativity of the above right diagram
 $\Phi, \Psi$ are well-defined family-maps  and $(\Phi, \Psi) \colon \D F(C^X, C^Y) =_{\Fam(I)} C^{\D F(X, Y)}$.
\end{proof}

The operations on families of sets generate operations on family-maps.

\begin{proposition}\label{prp: newfamilymaps1}
Let $\Lambda := (\lambda_0, \lambda_1), M := (\mu_0, \mu_1), N := (\nu_0, \nu_1)$, $K := (\kappa_0, \kappa_1) 
\in \Fam(I)$.\\[1mm]
\normalfont (i) 
\itshape If $\Phi \colon N \To \Lambda$ and $\Psi \colon N \To M$, then $\Phi \times \Psi \colon N \To \Lambda \times M$
is the \textit{product family-map}\index{product family-map}\index{$\Phi \times \Psi$} of $\Phi$ and $\Psi$, 
where, for every $i \in I$, the map $(\Phi \times \Psi)_i \colon \nu_0(i) \to \lambda_0(i) \times \mu_0(i)$ is defined by 
\[ (\Phi \times \Psi)_i(z) := \big(\Phi_i(z), \Psi_i(z)\big); \ \ \ \ z \in \nu_0(i). \]
\normalfont (ii) If $\Phi \colon N \To \Lambda$ and $\Psi \colon K \To M$, then $\Phi \times \Psi \colon N \times K
\To \Lambda \times M$ is the \textit{product family-map}\index{product family-map}\index{$\Phi \times \Psi$} of 
$\Phi$ and $\Psi$, where, for every $i \in I$, the map $(\Phi \times \Psi)_i \colon \nu_0(i) \times \kappa_0(i) 
\to \lambda_0(i) \times \mu_0(i)$ is defined by 
\[ (\Phi \times \Psi)_i(x, y) := \big(\Phi_i(x), \Psi_i(y)\big); \ \ \ \ (x, y) \in \nu_0(i) \times \kappa_0(i). \]
\itshape 
\normalfont (iii) 
\itshape If $\Phi \colon N \To \Lambda$, then $\D F(\Phi)^c \colon \D F(\Lambda, M) \To \D F(N, M)$, where,
for every $i \in I$, the function $\D F(\Phi)^c_i \colon \D F\big(\lambda_0(i), \mu_0(i)\big) \to 
\D F\big(\nu_0(i), \mu_0(i)\big)$ is defined by
\[ \D F(\Phi)^c_i(f) := f \circ \Phi_i; \ \ \ \ f \in \D F\big(\lambda_0(i), \mu_0(i)\big) \]
\begin{center}
\begin{tikzpicture}

\node (E) at (0,0) {$\nu_0(j)$};
\node[right=of E] (F) {$\lambda_0(i)$};
\node[right=of F] (A) {$\mu_0(i)$};
\node[right=of A] (B) {$\mu_0(j)$};
\node[right=of B] (C) {$\lambda_0(i)$};
\node[right=of C] (D) {$\nu_0(i)$.};

\draw[->] (E)--(F) node [midway,above]{$\Phi_i$};
\draw[->] (F)--(A) node [midway,above] {$f$};
\draw[->,bend right] (E) to node [midway,below] {$\D F(\Phi)^c_i(f)$} (A);
\draw[->] (B)--(C) node [midway,above]{$f$};
\draw[->] (C)--(D) node [midway,above] {$\Phi_i$};
\draw[->,bend right] (B) to node [midway,below] {$\D F(\Phi)^d_i(f)$} (D);

\end{tikzpicture}
\end{center}
If $\Phi \colon \Lambda \To N$, then $\D F(\Phi)^d \colon \D F(M, \Lambda) \To \D F(M, N)$, where,
for every $i \in I$ and $f \in \D F\big(\mu_0(i), \lambda_0(i)\big)$, the function $\D F(\Phi)^d_i 
\colon \D F\big(\mu_0(i), \lambda_0(i)\big) \to 
\D F\big(\mu_0(i), \nu_0(i)\big)$ is defined by
$\D F(\Phi)^d_i(f) := \Phi_i \circ f$.\\
\normalfont (iv) 
\itshape If $\Phi \colon N \To \Lambda$ and $\Psi \colon M \To K$, then $\D F(\Phi, \Psi) \colon \D F(\Lambda, M) \To 
\D F(N, K)$\index{$\D F(\Phi, \Psi)$}, where for every $i \in I$, the map $\D F(\Phi, \Psi)_i \colon 
\D F\big(\lambda_0(i), \mu_0(i) \to \D F\big(\nu_0(i), \kappa_0(i)$ is defined by
\[ \D F(\Phi, \Psi)_i(f) := \Psi_i \circ f \circ \Phi_i; \ \ \ \ f \in \D F\big(\lambda_0(i), \mu_0(i)\big) \]
\begin{center}
\begin{tikzpicture}

\node (E) at (0,0) {$\nu_0(i)$};
\node[right=of E] (H) {};
\node[right=of H] (F) {$\kappa_0(i)$.};
\node[above=of F] (A) {$\mu_0(i)$};
\node [above=of E] (D) {$\lambda_0(i)$};

\draw[->] (E)--(F) node [midway,below]{$\D F(\Phi, \Psi)_i(f)$};
\draw[->] (D)--(A) node [midway,above] {$f$};
\draw[->] (E)--(D) node [midway,left] {$\Phi_{i}$};
\draw[->] (A)--(F) node [midway,right] {$\Psi_{i}$};

\end{tikzpicture}
\end{center}
\end{proposition}

\begin{proof}
We prove (i) and (iii), as the proofs of (ii), (iv) are similar to that of (i), (iii), respectively. 
(i) If $i =_I j$, the following diagram is commutative
\begin{center}
\begin{tikzpicture}

\node (E) at (0,0) {$\mu_0(i) \times \lambda_0(i)$};
\node[right=of E] (H) {};
\node[right=of H] (F) {$\mu_0(j) \times \lambda_0(j)$,};
\node[above=of F] (A) {$\nu_0(j)$};
\node [above=of E] (D) {$\nu_0(i)$};

\draw[->] (E)--(F) node [midway,below]{$(\lambda_1 \times \mu_1)_{ij}$};
\draw[->] (D)--(A) node [midway,above] {$\nu_{ij}$};
\draw[->] (D)--(E) node [midway,left] {$(\Phi \times \Psi)_{i}$};
\draw[->] (A)--(F) node [midway,right] {$(\Phi \times \Psi)_{j}$};

\end{tikzpicture}
\end{center}
since by the commutativity of the following two diagrams
\begin{center}
\begin{tikzpicture}

\node (E) at (0,0) {$\lambda_0(i)$};
\node[right=of E] (F) {$\lambda_0(j)$};
\node[above=of F] (A) {$\nu_0(j)$};
\node [above=of E] (D) {$\nu_0(i)$};
\node [right=of F] (G) {$\mu_0(i)$};
\node [right=of G] (H) {$\mu_0(j)$,};
\node[above=of H] (K) {$\nu_0(j)$};
\node [above=of G] (L) {$\nu_0(i)$};

\draw[->] (E)--(F) node [midway,below]{$\lambda_{ij}$};
\draw[->] (D)--(A) node [midway,above] {$\nu_{ij}$};
\draw[->] (D)--(E) node [midway,left] {$\Phi_i$};
\draw[->] (A)--(F) node [midway,right] {$\Phi_j$};
\draw[->] (G)--(H) node [midway,below] {$\mu_{ij}$};
\draw[->] (L)--(K) node [midway,above] {$\nu_{ij}$};
\draw[->] (L)--(G) node [midway,left] {$\Psi_i$};
\draw[->] (K)--(H) node [midway,right] {$\Psi_j$};

\end{tikzpicture}
\end{center}
\begin{align*}
(\Phi \times \Psi)_{j}\big(\nu_{ij}(z)\big) & := \big(\Phi_j(\nu_{ij}(z)), \Psi_j(\nu_{ij}(z))\big)\\
& = \big(\lambda_{ij}(\Phi_i(z)), \mu_{ij}(\Psi_i(z))\big)\\
& := (\lambda_1 \times \mu_1)_{ij}\big(\Phi_i(z), \Psi_i(z)\big)\\
& := (\lambda_1 \times \mu_1)_{ij}\big((\Phi \times \Psi)_i(z)\big); \ \ \ \ z \in \nu_0(i).
\end{align*}
(ii) If $i =_I j$, the following diagram is commutative
\begin{center}
\begin{tikzpicture}

\node (E) at (0,0) {$\D F\big(\nu_0(i), \mu_0(i)\big)$};
\node[right=of E] (H) {};
\node[right=of H] (F) {$\D F\big(\nu_0(j), \mu_0(j)\big)$,};
\node[above=of F] (A) {$\D F\big(\lambda_0(j), \mu_0(j)\big)$};
\node [above=of E] (D) {$\D F\big(\lambda_0(i), \mu_0(i)\big)$};

\draw[->] (E)--(F) node [midway,below]{$\mathsmaller{\D F(\nu_1, \mu_1)_{ij}}$};
\draw[->] (D)--(A) node [midway,above] {$\mathsmaller{\D F(\lambda_1, \mu_1)_{ij}}$};
\draw[->] (D)--(E) node [midway,left] {$\D F(\Phi)^c_{i}$};
\draw[->] (A)--(F) node [midway,right] {$\D F(\Phi)^c_{j}$};

\end{tikzpicture}
\end{center}
\begin{align*}
\D F(\Phi)^c_{j}\big(\D F(\lambda_1, \mu_1)_{ij}(f)\big) & := \D F(\Phi)^c_{j}(\mu_{ij} \circ f \circ \lambda_{ji})\\
& := (\mu_{ij} \circ f \circ \lambda_{ji}) \circ \Phi_j\\
& := \mu_{ij} \circ f \circ (\lambda_{ji} \circ \Phi_j)\\
& := \mu_{ij} \circ f \circ (\Phi_i \circ \nu_{ji})\\
& := \mu_{ij} \circ (f \circ \Phi_i) \circ \nu_{ji}\\
& := \D F(\nu_1, \mu_1)_{ij}(f \circ \Phi_i)\\
& := \D F(\nu_1, \mu_1)_{ij}\big(\D F(\Phi)^c_{i}(f)\big); \ \ \ \ f \in \D F\big(\lambda_0(i), \mu_0(i)\big).
\end{align*}
The equality $\lambda_{ji} \circ \Phi_j = \Phi_i \circ \nu_{ji}$ used above follows from the definition of 
$\Phi \colon N \To \Lambda$ on $(j, i) \in D(I)$. The proof of $\D F(\Phi)^d \colon \D F(M, \Lambda) \To \D F(M, N)$ 
is similar.
%
%
%
%
\end{proof}

\section{The exterior union of a family of sets}
\label{sec: sigmaset}

\begin{definition}\label{def: sigmaset}
Let $\Lambda := (\lambda_0, \lambda_1)$ be an $I$-family of sets. The 
\textit{exterior union}\index{exterior union}, or 
\textit{disjoint union}\index{disjoint union}, or the $\sum$-\textit{set}\index{$\sum$-set of $\Lambda$}
$\sum_{i \in I}\lambda_0 (i)$\index{$\sum_{i \in I}\lambda_0 (i)$} of $\Lambda$, and its canonical equality
are defined by
\[ w \in \sum_{i \in I}\lambda_0 (i) : 
\TOT \exists_{i \in I}\exists_{x \in \lambda_0 (i)}\big(w := (i, x)\big), \]
\[ (i, x) =_{\mathsmaller{\sum_{i \in I}\lambda_0 (i)}} (j, y) : \TOT i =_I j \ \& \ \lambda_{ij} (x) 
=_{\lambda_0 (j)} y. \]
The $\sum$-set of the $\D 2$-family $\Lambda^{\D 2}$ of the sets $X$ and $Y$
is the \textit{coproduct}\index{coproduct}\index{$X + Y$} of $X$ and $Y$, and we write 
\[ X + Y := \sum_{i \in \D 2}\lambda_0^{\D 2}(i). \]
\end{definition}

\begin{proposition}\label{prp: sigmaset1}
\normalfont (i) 
\itshape The equality on $\sum_{i \in I}\lambda_0 (i)$ satisfies the conditions of an equivalence relation.\\[1mm]
\normalfont (ii) 
\itshape Let $(I, =_I, \neq_I)$ be a discrete set and $\neq_{\lambda_0(i)}$ an inequality on $\lambda_0(i)$, for every 
$i \in I$. If the transport map $\lambda_{ij}$ is strongly extensional, for every $(i, j) \in D(I)$, then the relation
\[(i,x) \neq_{\mathsmaller{\sum_{i \in I}\lambda_0 (i)}} (j, y) :\TOT i \neq_I j \ \vee \ \big(i =_I j \ \& \ 
\lambda_{ij}(x) \neq_{\lambda_0(j)} y \big) \]
is an inequality on $\sum_{i \in I}\lambda_0 (i)$. If $(\lambda_0(i), =_{\lambda_0(i)}, 
\neq_{\lambda_0(i)})$ is a discrete set, for every $i \in I$, then 
$\big(\sum_{i \in I}\lambda_0 (i), =_{\mathsmaller{\sum_{i \in I}\lambda_0 (i)}}, 
\neq_{\mathsmaller{\sum_{i \in I}\lambda_0 (i)}}\big)$ is discrete. Moreover, if $\neq_I$ is tight, and if, for every 
$i \in I$, the inequality $\neq_{\lambda_0(i)}$ is tight, then the inequality 
$\neq_{\mathsmaller{\sum_{i \in I}\lambda_0 (i)}} $ is tight.
 
\end{proposition}

\begin{proof}
(i)  Let 
$(i, x), (j, y)$, $(k, z) \in \sum_{i \in I}\lambda_0 (i)$. Since $i =_I i$ and $\lambda_{ii} 
:= \id_{\lambda_0 (i)}$, we get $(i, x) =_{\mathsmaller{\sum_{i \in I}\lambda_0 (i)}} (i, x)$. If $(i, x) 
=_{\mathsmaller{\sum_{i \in I}\lambda_0 (i)}} (j, y)$, then $j =_I i$ and 
$\lambda_{ji}(y) = \lambda_{ji}(\lambda_{ij}(x)) = \lambda_{ii}(x) := \id_{\lambda_0 (i)}(x) := x,$
hence $(j, y) =_{\mathsmaller{\sum_{i \in I}\lambda_0 (i)}} (i, x)$. If $(i, x)
=_{\mathsmaller{\sum_{i \in I}\lambda_0 (i)}} (j, y)$ and
$(j, y) =_{\mathsmaller{\sum_{i \in I}\lambda_0 (i)}} (k, z)$, then $i =_I j \ \& \ j =_I k \To i =_I k$, and 
\[ \lambda_{ik}(x) =_{\mathsmaller{\lambda_0(k)}} (\lambda_{jk} \circ \lambda_{ij})(x) := \lambda_{jk}(\lambda_{ij}(x)) 
=_{\mathsmaller{\lambda_0(k)}} \lambda_{jk}(y) =_{\mathsmaller{\lambda_0(k)}} z. \]
(ii) The condition $(\Ap_1)$ of Definition~\ref{def: apartness} is trivially satisfied. To show condition
$(\Ap_2)$, we suppose
first that $i \neq_I j$, hence by the corresponding condition of $\neq_I$ we get $(j, y) 
\neq_{\mathsmaller{\sum_{i \in I}\lambda_0 (i)}} (i, x)$. If $i =_I j \ \& \ 
\lambda_{ij}(x) \neq_{\lambda_0(j)} y$, we show that $\lambda_{ji}(y) \neq_{\lambda_0(i)} x$. By the extensionality of 
$\neq_{\lambda_0(j)}$ (Remark~\ref{rem: apartness1}) 
the inequality $\lambda_{ij}(x) \neq_{\lambda_0(j)} y$ implies the inequality
$\lambda_{ij}(x) \neq_{\lambda_0(j)} \lambda_{ij}\big(\lambda_{ji}(y)\big)$, and since 
$\lambda_{ij}$ is strongly extensional, we get $x \neq_{\lambda_0(i)} \lambda_{ji}(y)$. To show 
condition $(\Ap_3)$, let $(i,x) \neq_{\mathsmaller{\sum_{i \in I}\lambda_0 (i)}} (j, y)$, and let 
$(k, z) \in \sum_{i \in I}\lambda_0 (i)$. If $i \neq_I j$, then by condition $(\Ap_3)$ of $\neq_I$ we get
$k \neq_I i$, or $k \neq_I j$, hence $(k,z) \neq_{\mathsmaller{\sum_{i \in I}\lambda_0 (i)}} (i, x)$, or
$(k,z) \neq_{\mathsmaller{\sum_{i \in I}\lambda_0 (i)}} (j, y)$. Suppose next $i =_I j \ \& \ 
\lambda_{ij}(x) \neq_{\lambda_0(j)} y$. Since the set $(I, =_I, \neq_I)$ is discrete, $k \neq_I i$, or $k =_I i =_I j$.
If $k \neq_I i$, then what we want to show follows immediately. If $k =_I i =_I j$, then by the extensionality of 
$\neq_{\lambda_0(j)}$ and the strong extensionality of the transport map
$\lambda_{kj}$ we have that
\[ \lambda_{ij}(x) \neq_{\lambda_0(j)} y \To \lambda_{kj}\big(\lambda_{ik}(x)\big) \neq_{\lambda_0(j)} 
 \lambda_{kj}\big(\lambda_{jk}(y)\big) \To \lambda_{ik}(x) \neq_{\lambda_0(k)} \lambda_{jk}(y). \]
Hence, by condition $(\Ap_3)$ of $\neq_{\lambda_0(k)}$ we get $\lambda_{ik}(x) \neq_{\lambda_0(k)} z$, or
$\lambda_{jk}(y) \neq_{\lambda_0(k)} z$, hence $(i,x) \neq_{\mathsmaller{\sum_{i \in I}\lambda_0 (i)}} (k, z)$, or
$(j,y) \neq_{\mathsmaller{\sum_{i \in I}\lambda_0 (i)}} (k, z)$. Suppose next that 
$(\lambda_0(i), =_{\lambda_0(i)}, \neq_{\lambda_0(i)})$ is a discrete set, for every $i \in I$. We show that 
$(i,x) =_{\mathsmaller{\sum_{i \in I}\lambda_0 (i)}} (j, y)$ i.e., $i =_I j$ and $\lambda_{ij}(x) =_{\lambda_0(i)} y$, or 
$(i,x) \neq_{\mathsmaller{\sum_{i \in I}\lambda_0 (i)}} (j, y)$ i.e., 
$i \neq_I j$ or $i =_I j \ \& \ \lambda_{ij}(x) \neq_{\lambda_0(j)} y$. Since $(I, =_I, \neq_I)$ is discrete, 
$i =_I j$, or $i \neq_I j$. In the first case, and since $(\lambda_0(j), =_{\lambda_0(j)}, 
\neq_{\lambda_0(j)})$ is discrete, we get $\lambda_{ij}(x) =_{\lambda_0(j)} y$ or 
$\lambda_{ij}(x) \neq_{\lambda_0(j)} y$, and what we want follows immediately. If $i \neq_I j$, we get
$(i,x) \neq_{\mathsmaller{\sum_{i \in I}\lambda_0 (i)}} (j, y)$.
Finally, we suppose that $\neq_I$ is tight, and that $\neq_{\lambda_0(i)}$ is tight, for every 
$i \in I$. Let $\neg\big[(i,x) \neq_{\mathsmaller{\sum_{i \in I}\lambda_0 (i)}} (j, y)\big]$ i.e., 
\[ \big[i \neq_I j \ \vee \ \big(i =_I j \ \& \ \lambda_{ij}(x) \neq_{\lambda_0(j)} y \big)\big] \To \bot. \]
From this hypothesis we get the conjunction\footnote{Here we use the logical implication
$\big((\phi \vee  \psi) \To \bot\big) \To [(\phi \To \bot) \ \& \ (\psi \To \bot)]$.}
\[ \big[i \neq_I j \To \bot\big] \ \ \& \ \ \big[\big(i =_I j \ \& \ \lambda_{ij}(x) \neq_{\lambda_0(j)} y \big) \To
\bot\big]. \]
By the tightness of $\neq_I$ we get $i =_I j$. The implication 
$\big(i =_I j \ \& \ \lambda_{ij}(x) \neq_{\lambda_0(j)} y \big) \To
\bot$ logically implies the implication $(i =_I j) \To  \big(\lambda_{ij}(x) \neq_{\lambda_0(j)} y \To \bot\big)$,
and since 
its premiss $i =_I j$ is derived by the tightness of $\neq_I$, by Modus Ponens we get  
$\lambda_{ij}(x) \neq_{\lambda_0(j)} y \To \bot$. Since $\neq_{\lambda_0(i)}$ is tight, we conclude that 
$\lambda_{ij}(x) =_{\lambda_0(j)} y$, hence $(i, x) =_{\mathsmaller{\sum_{i \in I}\lambda_0 (i)}} (j, y)$.
\end{proof}

The totality $\sum_{i \in I}\lambda_0 (i)$ is considered to be a set.  
By the definition of $X + Y$ 
\begin{align*}
w \in X + Y & \TOT \exists_{i \in \D 2}\exists_{x \in \lambda_0^{\D 2}(i)}\big(w := (i, x)\big)\\
& \TOT \exists_{x \in X}\big(w := (0, x)\big) \ \vee \ \exists_{y \in Y}\big(w := (1, y)\big),
\end{align*}
\[ (i, x) =_{X + Y} (i{'}, x{'}) \TOT (i =_{\D 2} i{'} =_{\D 2} 0 \ \& \ x =_X x{'}) \ \vee \
(i =_{\D 2} i{'} =_{\D 2} 1 \ \& \ x =_Y x{'}). \]
One could have defined $X + Y$ independently from $\Lambda^{\D 2}$, and then prove 
$X + Y =_{\D V_0} \sum_{i \in \D 2}\lambda_0^{\D 2}(i)$.

\begin{corollary}\label{cor: coproddiscrete}
If $(X, =_X, \neq_X)$, $(Y, =_Y, \neq_Y)$ are discrete, 
 $(X + Y, =_{X + Y}, \neq_{X + Y})$ is
 discrete.
\end{corollary}

\begin{proof}
Since $(\D 2, =_{\D 2}, \neq_{\D 2})$ is a discrete set, we use Proposition~\ref{prp: sigmaset1}(ii).
\end{proof}

\begin{definition}\label{def: newfamsofsets2}
Let $\Lambda := (\lambda_0, \lambda_1), M := (\mu_0, \mu_1)$ be $I$-families of sets. 
The \textit{coproduct family}\index{coproduct of $I$-families} of $\Lambda$ and $M$ is the pair\index{$\Lambda + M$} 
$\Lambda + M := (\lambda_0 + \mu_0, \lambda_1 + \mu_1)$, where $(\lambda_0 + \mu_0)(i) := \lambda_0 (i) + \mu_0 (i)$, 
for every $i \in I$, and the map $\big(\lambda_1 + \mu_1\big)_{ij} \colon \lambda_0 (i) + 
\mu_0 (i) \to \lambda_0 (j) + \mu_0 (j)$ is defined by
\[ \big(\lambda_1 + \mu_1\big)_{ij}(w) := \left\{ \begin{array}{ll}
                 \big(0, \lambda_{ij}(x)\big)   &\mbox{, $w := (0, x)$}\\
                 \big(1, \mu_{ij}(y)\big)             &\mbox{, $w := (1, y)$}
                 \end{array}
          \right.
; \ \ \ \ w \in \lambda_0 (i) + \mu_0 (i). \]

\end{definition}

It is straightforward to show that $\Lambda + M$ is an $I$-family of sets.

\begin{proposition}\label{prp: newfamilymaps2}
Let $\Lambda := (\lambda_0, \lambda_1), M := (\mu_0, \mu_1)$, and $N := (\nu_0, \nu_1)$ be $I$-families of sets.
If $\Phi \colon \Lambda \To N$ and $\Psi \colon M \To N$, then $\Phi + \Psi \colon \Lambda + M \To N$ is
\index{$\Phi + \Psi $} the \textit{coproduct family-map}\index{coproduct family-map} of $\Phi$ and 
$\Psi$, where,
for every $i \in I$, the map $(\Phi + \Psi)_i \colon \lambda_0(i) + \mu_0(i) \to \nu_0(i)$ is defined by 
\[ (\Phi + \Psi)_i(w) := \left\{ \begin{array}{ll}
                 \Phi_i(x)   &\mbox{, $w := (0, x)$}\\
                 \Psi_i(y)             &\mbox{, $w := (1, y)$}
                 \end{array}
          \right.
; \ \ \ \ w \in \lambda_0 (i) + \mu_0 (i). \]
\end{proposition}

\begin{proof}
If $i =_I j$, the following diagram is commutative
\begin{center}
\begin{tikzpicture}

\node (E) at (0,0) {$\nu_0(i)$};
\node[right=of E] (H) {};
\node[right=of H] (F) {$\nu_0(j)$,};
\node[above=of F] (A) {$\mathsmaller{\lambda_0(j) + \mu_0(j)}$};
\node [above=of E] (D) {$\mathsmaller{\lambda_0(i) + \mu_0(i)}$};

\draw[->] (E)--(F) node [midway,below]{$\nu_{ij}$};
\draw[->] (D)--(A) node [midway,above] {$\mathsmaller{(\lambda_1 + \mu_1)_{ij}}$};
\draw[->] (D)--(E) node [midway,left] {$(\Phi + \Psi)_{i}$};
\draw[->] (A)--(F) node [midway,right] {$(\Phi + \Psi)_{j}$};

\end{tikzpicture}
\end{center}
since by the commutativity of the following left diagram
\begin{center}
\begin{tikzpicture}

\node (E) at (0,0) {$\nu_0(i)$};
\node[right=of E] (F) {$\nu_0(j)$};
\node[above=of F] (A) {$\lambda_0(j)$};
\node [above=of E] (D) {$\lambda_0(i)$};
\node [right=of F] (G) {$\nu_0(i)$};
\node [right=of G] (H) {$\nu_0(j)$,};
\node[above=of H] (K) {$\mu_0(j)$};
\node [above=of G] (L) {$\mu_0(i)$};

\draw[->] (E)--(F) node [midway,below]{$\nu_{ij}$};
\draw[->] (D)--(A) node [midway,above] {$\lambda_{ij}$};
\draw[->] (D)--(E) node [midway,left] {$\Phi_i$};
\draw[->] (A)--(F) node [midway,right] {$\Phi_j$};
\draw[->] (G)--(H) node [midway,below] {$\nu_{ij}$};
\draw[->] (L)--(K) node [midway,above] {$\mu_{ij}$};
\draw[->] (L)--(G) node [midway,left] {$\Psi_i$};
\draw[->] (K)--(H) node [midway,right] {$\Psi_j$};

\end{tikzpicture}
\end{center}
\begin{align*}
(\Phi + \Psi)_{j}\big((\lambda_1 + \mu_1)_{ij}(0, x)\big) & := (\Phi + \Psi)_{j}\big(0, \lambda_{ij}(x)\big)\\
& = \Phi_j\big(\lambda_{ij}(x)\big)\\
& = \nu_{ij}\big(\Phi_i(x)\big)\\
& := \nu_{ij}\big((\Phi + \Psi)_i(0, x)\big); \ \ \ \ x \in \lambda_0(i).
\end{align*}
By the commutativity of the right diagram, 
$(\Phi + \Psi)_{j}\big((\lambda_1 + \mu_1)_{ij}(1, y)\big) = \nu_{ij}\big((\Phi + \Psi)_i(1, y)\big)$.
\end{proof}

\begin{proposition}\label{prp: sumplus}
If $\Lambda := (\lambda_0, \lambda_1)$, $M := (\mu_0, \mu_1) \in \Fam(I)$, then
\[ \sum_{i \in I}\big(\lambda_0(i) + \mu_0(i)\big) =_{\D V_0} \bigg(\sum_{i \in I}\lambda_0(i)\bigg) + 
\bigg(\sum_{i \in I}\mu_0(i)\bigg). \]
\end{proposition}

\begin{proof}
Let $f \colon \sum_{i \in I}\big(\lambda_0(i) + \mu_0(i)\big) \sto \sum_{i \in I}\lambda_0(i) + 
\sum_{i \in I}\mu_0(i)$ be defined by 
\[ f(i, w) := \left\{ \begin{array}{ll}
                 \big(0, (i, x)\big)   &\mbox{, $w := (0, x)$}\\
                 \big(1, (i, y)\big)   &\mbox{, $w := (1, y)$}
                 \end{array}
          \right.
; \ \ \ \ i \in I, \ w \in \lambda_0 (i) + \mu_0 (i). \]
Clearly, $f$ is a well-defined operation. To show that $f$ is a function, we suppose that
\[ (i, w) =_{\mathsmaller{\sum_{i \in I}(\lambda_0(i) + \mu_0(i))}} (j, u) :\TOT i =_I j \ \& \ 
(\lambda_1 + \mu_1)_{ij}(w) =_{\mathsmaller{}\lambda_0(j) + \mu_0(j)} u, 
\]
and we show that $f(i, w) =
f(j, u)$. The equality $(\lambda_1 + \mu_1)_{ij}(w) =_{\mathsmaller{}\lambda_0(j) + \mu_0(j)} u$ amounts to 
$\lambda_{ij}(x) =_{\lambda_0(j)} x{'}$, if $w := (0, x)$ and $u := (0, x{'})$, or to   
$\mu_{ij}(y) =_{\mu_0(j)} y{'}$, if $w := (1, y)$ and $u := (1, y{'})$. With the use of these equalities and the 
definition of the canonical equality on the coproduct it
is straightforward to show that $\big(0, (i, x)\big) = \big(0, (j, x{'})\big)$, or 
$\big(1, (i, y)\big) = \big(1, (j, y{'})\big)$, hence $f(i, w) = f(j, u)$. Let $g \colon \sum_{i \in I}\lambda_0(i) + 
\sum_{i \in I}\mu_0(i) \sto \sum_{i \in I}\big(\lambda_0(i) + \mu_0(i)\big)$ be defined by 
\[ g(U) := \left\{ \begin{array}{ll}
                 \big(i, (0, x)\big)   &\mbox{, $U := \big(0, (i, x)\big)$}\\
                 \big(i, (1, y)\big)   &\mbox{, $U := \big(1, (i, y)\big)$}
                 \end{array}
          \right.
; \ \ \ \ U \in \sum_{i \in I}\lambda_0(i) + \sum_{i \in I}\mu_0(i). \]
Proceeding similarly, we show that the operation $g$ is a function. It is straightforward to show that $(f, g) \colon 
\sum_{i \in I}\big(\lambda_0(i) + \mu_0(i)\big) =_{\D V_0} \sum_{i \in I}\lambda_0(i) + 
\sum_{i \in I}\mu_0(i)$.
\end{proof}

\begin{proposition}\label{prp: map1}
Let $\Lambda := (\lambda_0, \lambda_1)$, $M := (\mu_0, \mu_1) \in \Fam(I)$, and $\Psi : \Lambda \To M$.\\[1mm]
\normalfont (i) 
\itshape For every $i \in I$ the operation
$e_i^{\Lambda} : \lambda_0(i) \sto \sum_{i \in I}\lambda_0(i)$, defined by 
$ e_i^{\Lambda}(x) := (i, x)$, for every $x \in \lambda_0(i)$,
is an embedding.\\[1mm]
\normalfont (ii) 
\itshape The operation
$\Sigma \Psi : \sum_{i \in I}\lambda_0(i) \sto \sum_{i \in I}\mu_0(i)$, defined by
\[ \Sigma \Psi (i, x) := (i, \Psi_i (x)); \ \ \ \ (i, x) \in \sum_{i \in I}\lambda_0(i), \] 
is a function,
such that for every
$i \in I$ the following diagram commutes
\begin{center}
\begin{tikzpicture}

\node (E) at (0,0) {$\sum_{i \in I}\lambda_0(i)$};
\node[right=of E] (F) {$\sum_{i \in I}\mu_0(i)$.};
\node[above=of F] (A) {$\mu_0(i)$};
\node [above=of E] (D) {$\lambda_0(i)$};

\draw[->] (E)--(F) node [midway,below]{$\Sigma \Psi$};
\draw[->] (D)--(A) node [midway,above] {$\Psi_{i}$};
\draw[->] (D)--(E) node [midway,left] {$e_i^{\Lambda}$};
\draw[->] (A)--(F) node [midway,right] {$e_i^M$};

\end{tikzpicture}
\end{center}
\normalfont (iii) 
\itshape If $\Psi_i$ is an embedding, for every $i \in I$, then $\Sigma \Psi$ is an
embedding.\\[1mm]
\normalfont (iv) 
\itshape If $\Psi_i$ is a surjection, for every $i \in I$, then $\Sigma \Psi$ is an
surjection.\\[1mm]
\normalfont (v) 
\itshape If $\Phi : M \To \Lambda$, where $\Phi_i$ is a modulus of surjectivity for $\Psi_i$, for every $i \in I$,
then a modulus of surjectivity for $\Sigma \Psi$ is the operation 
$\sigma \Psi \colon \sum_{i \in I}\mu_0(i) \sto \sum_{i \in I}\lambda_0(i)$,
defined by
\[ \sigma \Psi (i, y) := (i, \Phi_i (y)); \ \ \ \ (i, y) \in \sum_{i \in I}\mu_0(i). \]

\end{proposition}

\begin{proof}
(i)  If $x, y \in \lambda_0(i)$, then $e_i^{\Lambda}(x) =_{\mathsmaller{\sum_{i \in I}\lambda_0(i)}} e_i^{\Lambda}(y)$
if and only if $(i, x) =_{\mathsmaller{\sum_{i \in I}\lambda_0(i)}} (i, y) $, which is equivalent to 
$\lambda_{ii}(x) =_{\lambda_0(i)} y :\TOT x =_{\lambda_0(i)} y$.\\
(ii) If $(i, x) =_{\sum_{i \in I}\lambda_0(i)} (j, y)$ i.e., $ i =_I j$ and $\lambda_{ij}(x) 
=_{\lambda_0(j)} y$, we show that $(i, \Psi_i(x)) =_{\mathsmaller{\sum_{i \in I}\mu_0(i)}} (j, \Psi_j(y))$ 
i.e., $ i =_I j$ and $\mu_{ij}(\Psi_i(x)) =_{\mu_0(j)} \Psi_j(y)$. Since $\Psi \colon \Lambda \To M$, 
we get $\mu_{ij}(\Psi_i(x)) =_{\mu_0(j)} \Psi_j(\lambda_{ij}(x)) 
=_{\mu_0(j)} \Psi_j(y)$. The required commutativity of the diagram is immediate to show.\\
(iii) Since $\Psi$ is a family-map from $\Lambda$ to $M$, we have that
\begin{align*}
\Sigma \Psi (i, x) =_{\mathsmaller{\sum_{i \in I}\mu_0(i)}} \Sigma \Psi (j, y) & : \TOT (i, \Psi_i(x)) 
=_{\mathsmaller{\sum_{i \in I}\mu_0(i)}} (j, \Psi_j(y)) \\
& : \TOT i =_I j \ \& \ \mu_{ij}(\Psi(x)) =_{\mu_0(j)} \Psi_j(y)\\
& \TOT  i =_I j \ \& \ \Psi_{j}(\lambda_{ij}(x)) =_{\mu_0(j)} \Psi_j(y)\\
& \To i =_I j \ \& \ \lambda_{ij}(x) =_{\lambda_0(j)} y \\
& : \TOT (i, x) =_{\mathsmaller{\sum_{i \in I}\lambda_0(i)}} (j, y).
\end{align*} 
(iv) Let $(i, y) \in \sum_{i \in I}\mu_0(i)$. Since $\Psi_i$ is a surjection, there is $x \in \lambda_0(i)$ such that
$\Psi_i(x) = y$. Hence $\Sigma \Psi (i, x) := (i, \Psi_i (x)) =_{\mathsmaller{\sum_{i \in I}\mu_0(i)}} (i, y)$, 
since $\mu_{ii}\big(\Psi_i (x)\big) := \Psi_i (x) =_{\mu_0(i)} y$.\\
(v) If $y \in \mu_0(i)$, then $\Psi_i(\Phi_i(y)) =_{\mu_0(i)} y$. To show that the operation 
 $\sigma \Psi$ is a function, we suppose $(i, y) =_{\mathsmaller{\sum_{i \in I}\mu_0(i)}} (j, z) :\TOT
 i =_I j \ \& \ \mu_{ij}(y) =_{\mu_0(j)} z$ and we show that 
$(i, \sigma_i(y)) =_{\mathsmaller{\sum_{i \in I}\lambda_0(i)}} (j, \sigma_j(z)) :\TOT
 i =_I j \ \& \ \lambda_{ij}(\Phi_i(y)) =_{\lambda_0(j)} \Phi_j(z)$. Since $\Phi_j \colon \mu_0(j) \to \lambda_0(j)$,
 we have that $\mu_{ij}(y) =_{\mu_0(j)} z \To \Phi_j\big(\mu_{ij}(y)\big) =_{\lambda_0(j)} \Phi_j(z)$.
 By the commutativity of the diagram 
 \begin{center}
\begin{tikzpicture}

\node (E) at (0,0) {$\lambda_0(i)$};
\node[right=of E] (F) {$\lambda_0(j)$,};
\node[above=of F] (A) {$\mu_0(j)$};
\node [above=of E] (D) {$\mu_0(i)$};

\draw[->] (E)--(F) node [midway,below]{$\lambda_{ij}$};
\draw[->] (D)--(A) node [midway,above] {$\mu_{ij}$};
\draw[->] (D)--(E) node [midway,left] {$\Phi_i$};
\draw[->] (A)--(F) node [midway,right] {$\Phi_j$};

\end{tikzpicture}
\end{center}
$\Phi_j(z) =_{\lambda_0(j)}  \Phi_j\big(\mu_{ij}(y)\big)  =_{\lambda_0(j)} \lambda_{ij}(\Phi_i(y))$. 
Since $\mu_{ii}\big(\Psi_i(\Phi_i(y))\big) := \Psi_i(\Phi_i(y))\big) =_{\mu_0(i)} y$, 
\[ \Sigma \Psi\big(\sigma \Psi (i, y)\big) := \Sigma \Psi\big((i, \Phi_i (y))\big) := \big(i, \Psi_i(\Phi_i(y))\big)
 =_{\mathsmaller{\sum_{i \in I}\mu_0(i)}} (i, y).\qedhere \]
\end{proof}

\begin{definition}\label{def: projection1}
Let $\Lambda := (\lambda_0, \lambda_1)$ be an $I$-family of sets. The
\textit{first projection}\index{first projection} on
$\sum_{i \in I}\lambda_0 (i)$ is the operation\index{$\pr_1^{\Lambda}$}
$\pr_1^{\Lambda} \colon \sum_{i \in I}\lambda_0 (i) \sto I$, 
defined by $pr_1^{\Lambda} (i, x) : = \prb_1 (i, x) := i$, for every $(i, x) \in \sum_{i \in I}\lambda_0 (i)$.
We may only write $\pr_1$, if $\Lambda$ is clearly understood from the context.
\end{definition}

By the definition of the canonical equality on $\sum_{i \in I}\lambda_0 (i)$ we get  
that $\pr_1^{\Lambda}$ is a function.

\begin{definition}\label{def: forproj2}
Let $\Lambda := (\lambda_0, \lambda_1)$ be an $I$-family of sets. 
The $\sum$-\textit{indexing}\index{$\sum$-indexing of a family} of $\Lambda$ is the pair\index{$\Sigma^{\Lambda}$} 
$\Sigma^{\Lambda} := (\sigma_0^{\Lambda}, \sigma_1^{\Lambda})$, where $\sigma_0^{\Lambda}
\colon \sum_{i \in I}\lambda_0 (i) \sto \D V_0$ is defined by
$\sigma_0^{\Lambda}(i, x) := \lambda_0(i)$, for every $(i, x) \in \sum_{i \in I}\lambda_0 (i)$, and
$\sigma_1^{\Lambda}\big((i,x),(j,y)\big) := \lambda_{ij}$, for every 
$\big((i,x),(j,y)\big) \in 
D\big(\sum_{i \in I}\lambda_0 (i)\big)$.
%

\end{definition}

Clearly, $\Sigma^{\Lambda}$ is a family of sets over $\sum_{i \in I}\lambda_0 (i)$, and 
$\Sigma \colon \Sigma^{\Lambda} \stackrel{\mathsmaller{\pr_1^{\Lambda}}} \Longrightarrow \Lambda$ (see 
Definition~\ref{def: hmap}), where, if $w := (i, x) \in \sum_{i \in I}\lambda_0 (i)$, we define 
$\Sigma_w : \lambda_0(i) \to \lambda_0(\mathsmaller{\pr_1^{\Lambda}}(w)) $ to be the identity $\id_{\lambda_0(i)}$.

\begin{definition}\label{def: projection2}
Let $\Lambda := (\lambda_0, \lambda_1)$ be an $I$-family of sets. The
\textit{second projection}\index{second projection} on
$\sum_{i \in I}\lambda_0 (i)$ is the dependent operation\index{$\pr_2^{\Lambda}$}
$\pr_2^{\Lambda} \colon \bigcurlywedge_{(i,x) \in \sum_{i \in I}\lambda_0 (i)}\lambda_0(i)$, defined by
$\pr_2^{\Lambda}(i,x) := \prb_2(i,x) := x$, for every $(i, x) \in \sum_{i \in I}\lambda_0 (i)$.
We may only write $\pr_2$, when the family of sets $\Lambda$ is clearly understood from the context.
\end{definition}

In Remark~\ref{rem: proj2} we show that $\pr_2^{\Lambda}$ is a dependent function over the family
$\Sigma^{\Lambda}$.

\section{Dependent functions over a family of sets}
\label{sec: piset}

\begin{definition}\label{def: piset}
Let $\Lambda := (\lambda_0, \lambda_1)$ be an $I$-family of sets. 
The totality $\prod_{i \in I}\lambda_0(i)$\index{$\prod_{i \in I}\lambda_0(i)$} of \textit{dependent functions over}
\index{dependent functions over a family of sets} $\Lambda$, or the 
$\prod$-set of $\Lambda$\index{$\prod$-set of $\Lambda$}, is defined by
$$\Theta \in \prod_{i \in I}\lambda_0(i) :\TOT \Theta \in \D A(I, \lambda_0) \ \& \ \forall_{(i,j) \in D(I)}\big(\Theta_j 
=_{\lambda_0(j)} \lambda_{ij}(\Theta_i)\big),$$
and it is equipped with the canonical equality and the canonical inequality 
of the set $\D A(I, \lambda_0)$. 
If $X$ is a set and $\Lambda^X$ is the constant $I$-family $X$ $($see Definition~\ref{def: famofsets}$)$,
we use the notation\index{$X^I$}
\[ X^I := \prod_{i \in I}X. \]
\end{definition}

Clearly, the property $P(\Phi) :\TOT \forall_{(i,j) \in D(I)}\big(\Theta_j 
=_{\lambda_0(j)} \lambda_{ij}(\Theta_i)\big)$
is extensional on $\D A(I, \lambda_0)$,  
the equality on $\prod_{i \in I}\lambda_0(i)$ is an equivalence relation.   
$\prod_{i \in I}\lambda_0(i)$ is considered to be a set.

\begin{remark}\label{rem: proj2}
If $\Lambda := (\lambda_0, \lambda_1)$ is an $I$-family of sets and 
$\Sigma^{\Lambda} := (\sigma_0^{\Lambda}, \sigma_1^{\Lambda})$ is the $\sum$-indexing of $\Lambda$, then
$\pr_2^{\Lambda}$ is a dependent function over $\Sigma^{\Lambda}$.
 
\end{remark}

\begin{proof}
By Definition~\ref{def: projection2} the second projection $\pr_2^{\Lambda}$ of $\Lambda$ is the dependent
assignment 
$\pr_2^{\Lambda} \colon \bigcurlywedge_{(i,x) \in \sum_{i \in I}\lambda_0 (i)}\lambda_0(i),$
such that $\pr_2^{\Lambda}(i,x) := x$, for every $(i, x) \in  \sum_{i \in I}\lambda_0 (i)$. It suffices to show that
if $(i, x) =_{\mathsmaller{\sum_{i \in I}\lambda_0 (i)}} (j, y) : \TOT i =_I j \ \& \ \lambda_{ij} (x) 
=_{\lambda_0 (j)} y$, then 
\[ \pr_2^{\Lambda}(j,y) := y
=_{\lambda_0 (j)} \lambda_{ij} (x)
:= \sigma_1^{\Lambda}\big((i,x),(j,y)\big)\big(\pr_2^{\Lambda}(i,x)\big).\qedhere
\]
\end{proof}

\begin{remark}\label{rem: proddep}
\normalfont (i) 
\itshape If $\Lambda^{\D 2}$ is the $\D 2$-family of the sets $X$ and $Y$, then
$\prod_{i \in \D 2}\lambda_0^{\D 2} (i) =_{\D V_0} X \times Y$.\\[1mm] 
\normalfont (ii) 
\itshape
If $I, A$ are sets, and $\Lambda := (\lambda_0^A, \lambda_1)$ is the constant $I$-family $A$, then
$A^I =_{\D V_0} \D F(I, A)$.
\end{remark}

\begin{proof}
(i) Let $f \colon \prod_{i \in \D 2}\lambda_0^{\D 2} (i) \sto X \times Y$ be defined by $f(\Phi) := (\Phi_0, \Phi_1)$, 
for every $\Phi \in \prod_{i \in \D 2}\lambda_0^{\D 2} (i)$. Let 
$g \colon X \times Y \sto \prod_{i \in \D 2}\lambda_0^{\D 2} (i)$ be defined by $g(x, y) := \Phi_{(x, y)}$, 
for every $(x, y) \in X \times Y$. It is easy to show that $f, g$ are well-defined functions and $(f, g) \colon 
\prod_{i \in \D 2}\lambda_0^{\D 2} (i) =_{\D V_0} X \times Y$.\\
(ii)  Let $h \colon A^I \sto \D F(I, A)$ be defined by $h(\Phi) := h_{\Phi} \colon I \to A$, where 
$ h_{\Phi}(i) := \Phi_i$, for every $\Phi \in A^I$ and $i \in I$. Let 
$k \colon \D F(I, A) \sto A^I$ be defined by $k(e) := \Phi_e$, where $[\Phi_e]_i := e(i)$,
for every $e \in \D F(I, A)$ and $i \in I$. Then $h, k$ 
are well-defined functions and $(h, k) \colon A^I =_{\D V_0} \D F(I, A)$.
\end{proof}

\begin{corollary}\label{cor: mapdependent}
If $\Lambda, M \in \Fam(I)$ and $\Psi \colon \bigcurlywedge_{i \in I}\D F\big(\lambda_0(i), \mu_0(i)\big)$, 
the following are equivalent:\\[1mm]
\normalfont (i)
\itshape $\Psi \colon \Lambda \To M$.\\[1mm]
\normalfont (ii) 
\itshape $\Psi \in \prod_{i \in I}\big[\D F(\lambda_0, \mu_0)\big](i)$.

\end{corollary}

\begin{proof}
If $i =_I j$, the commutativity of the following left diagram
\begin{center}
\begin{tikzpicture}

\node (E) at (0,0) {$\mu_0(i)$};
\node[right=of E] (F) {$\mu_0(j)$};
\node[above=of F] (A) {$\lambda_0(j)$};
\node [above=of E] (D) {$\lambda_0(i)$};
\node [right=of F] (G) {$\mu_0(i)$};
\node [right=of G] (H) {$\mu_0(j)$,};
\node[above=of H] (K) {$\lambda_0(j)$};
\node [above=of G] (L) {$\lambda_0(i)$};

\draw[->] (E)--(F) node [midway,below]{$\mu_{ij}$};
\draw[->] (D)--(A) node [midway,above] {$\lambda_{ij}$};
\draw[->] (D)--(E) node [midway,left] {$\Phi_i$};
\draw[->] (A)--(F) node [midway,right] {$\Phi_j$};
\draw[->] (G)--(H) node [midway,below] {$\mu_{ij}$};
\draw[->] (K)--(L) node [midway,above] {$\lambda_{ji}$};
\draw[->] (L)--(G) node [midway,left] {$\Phi_i$};
\draw[->] (K)--(H) node [midway,right] {$\Phi_j$};

\end{tikzpicture}
\end{center}
is equivalent to the commutativity of the above right one, hence the defining condition for $\Psi \in  \Map(\Lambda, M)$
is equivalent to the defining condition $\Psi_j = \D F(\lambda_1, \mu_1)_{ij}(\Psi_i) := \mu_{ij} \circ \Psi_i 
\circ \lambda_{ji}$ for $\Psi \in \prod_{i \in I}\big(\lambda_0(i) \times \mu_0(i)\big)$ 
(see Definition~\ref{def: newfamsofsets1}(ii)).
 \end{proof}

%
%

\begin{proposition}\label{prp: map2}
Let $\Lambda := (\lambda_0, \lambda_1)$, $M := (\mu_0, \mu_1) \in \Fam(I)$, and $\Psi : \Lambda \To M$.\\[1mm]
\normalfont (i) 
\itshape If $i \in I$, the operation
$\pi_i^{\Lambda} : \prod_{i \in I}\lambda_0(i) \sto \lambda_0(i)$, defined by
$\Theta \mapsto \Theta_i,$ is a function.\\[1mm]
\normalfont (ii) 
\itshape The operation
$\Pi \Psi : \prod_{i \in I}\lambda_0(i) \sto \prod_{i \in I}\mu_0(i)$, defined by
\[ [\Pi \Psi (\Theta)]_i := \Psi_i (\Theta_i); \ \ \ \ i \in I, \]
is a function, 
such that for every $i \in I$ the following diagram commutes
\begin{center}
\begin{tikzpicture}

\node (E) at (0,0) {$\prod_{i \in I}\lambda_0(i)$};
\node[right=of E] (F) {$\prod_{i \in I}\mu_0(i)$.};
\node[above=of F] (A) {$\mu_0(i)$};
\node [above=of E] (D) {$\lambda_0(i)$};

\draw[->] (E)--(F) node [midway,below]{$\Pi \Psi$};
\draw[->] (D)--(A) node [midway,above] {$\Psi_{i}$};
\draw[->] (E)--(D) node [midway,left] {$\pi_i^{\Lambda}$};
\draw[->] (F)--(A) node [midway,right] {$\pi_i^M$};

\end{tikzpicture}
\end{center}
\normalfont (iii) 
\itshape If $\Psi_i$ is an embedding, for every $i \in I$, then $\Pi \Psi$ is an 
embedding.\\[1mm]
\normalfont (iv) 
\itshape If $\Phi : M \To \Lambda$ such that $\Phi_i$ is a modulus of surjectivity for $\Psi_i$, for every $i \in I$,
the operation $\pi \Psi \colon \prod_{i \in I}\mu_0(i) \sto \prod_{i \in I}\lambda_0(i)$ is a modulus of surjectivity for
$\prod \Psi$, where
\[ \big[\pi \Psi (\Omega)\big]_i := \Phi_i (\Omega_i); \ \ \ \ \Omega \in \prod_{i \in I}\mu_0(i), \ i \in I. \]
\end{proposition}

\begin{proof}
(i) This follows immediately from the definition of equality on $\prod_{i \in I}\lambda_0(i)$.\\
(ii) First we show that $\Pi \Psi$ is well-defined i.e., $\Pi \Psi (\Theta) \in
\prod_{i \in I}\mu_0(i)$. If $i =_I j$, then by the commutativity of the following left diagram from
the definition of a family-map 
\begin{center}
\begin{tikzpicture}

\node (E) at (0,0) {$\mu_0(i)$};
\node[right=of E] (F) {$\mu_0(j)$};
\node[above=of F] (A) {$\lambda_0(j)$};
\node [above=of E] (D) {$\lambda_0(i)$};
\node [right=of F] (G) {$\lambda_0(i)$};
\node[right=of G] (H) {$\lambda_0(j)$,};
\node[above=of H] (K) {$\mu_0(j)$};
\node [above=of G] (L) {$\mu_0(i)$};

\draw[->] (E)--(F) node [midway,below]{$\mu_{ij}$};
\draw[->] (D)--(A) node [midway,above] {$\lambda_{ij}$};
\draw[->] (D)--(E) node [midway,left] {$\Psi_i$};
\draw[->] (A)--(F) node [midway,right] {$\Psi_j$};
\draw[->] (G)--(H) node [midway,below]{$\lambda_{ij}$};
\draw[->] (L)--(K) node [midway,above] {$\mu_{ij}$};
\draw[->] (L)--(G) node [midway,left] {$\Phi_i$};
\draw[->] (K)--(H) node [midway,right] {$\Phi_j$};

\end{tikzpicture}
\end{center}
\[ [\Pi \Psi (\Theta)]_j := \Psi_j (\Theta_j) = \Psi_j \big(\lambda_{ij}(\Theta_i)\big) = 
 \mu_{ij}\big(\Psi_i (\Theta_i)\big) := \mu_{ij}\big([\Pi \Psi (\Theta)]_i\big).
\]
It is immediate to show that $\Pi \Psi$ is a function and that the required diagram commutes.\\
(iii) If $\Theta, \Theta{'} \in \prod_{i \in I}\lambda_0(i)$, then
\begin{align*}
\Pi \Psi (\Theta) =_{\mathsmaller{\prod_{i \in I}\mu_0(i)}} \Pi \Psi (\Theta{'}) & :\TOT
\forall_{i \in I}\big(\Psi_i (\Theta_i) =_{\mu_0(i)} \Psi_i (\Theta{'}_i)\big)\\
& \To \forall_{i \in I}\big(\Theta_i =_{\lambda_0(i)} \Theta{'}_i\big)\\
& :\TOT \Theta =_{\mathsmaller{\prod_{i \in I}\lambda_0(i)}} \Theta{'}.
\end{align*}
(iv) First we show that $\pi \Psi$ is well-defined i.e., $\pi \Psi (\Omega) \in \prod_{i \in I}\lambda_0(i)$. If 
$i =_I j$, and since $\Phi : M \To \Lambda$, by the commutativity of the above right diagram
%
%
%
%
%
\[ \big[\pi \Psi (\Omega)\big]_j := \Phi_j (\Omega_j) =_{\lambda_0(j)} \Phi_j \big(\mu_{ij}(\Omega_i)\big)
=_{\lambda_0(j)} \lambda_{ij}\big(\Phi_i (\Omega_i)\big) := 
\lambda_{ij}\big(\big[\pi \Psi (\Omega)\big]_i\big). \]  
It is immediate to show that $\pi \Psi$ is a function. Finally we show that $\Pi \Psi\big(\pi \Psi (\Omega)\big) 
= \Omega$,
for every $\Omega \in \prod_{i \in I}\mu_0(i)$. If $i \in I$, and since $\Phi_i$ is a modulus of 
surjectivity for $\Psi_i$,
we get  
\[ \big[\Pi \Psi\big(\pi \Psi (\Omega)\big)\big]_i := \Psi_i\big(\big[\pi \Psi (\Omega)\big]_i\big) := 
 \Psi_i\big(\Phi_i(\Omega_i)\big) = \Omega_i.\qedhere \]

\end{proof}

\begin{proposition}\label{prp: MAPdep}
If $\Lambda := (\lambda_0, \lambda_1)$, $M := (\mu_0, \mu_1) \in \Fam(I)$, then
\[ \prod_{i \in I}\big(\lambda_0(i) \times \mu_0(i)\big) =_{\D V_0} \bigg(\prod_{i \in I}\lambda_0(i)\bigg) \times 
\bigg(\prod_{i \in I}\mu_0(i)\bigg),   \]
\[\Map_I(\Lambda, M) =_{\D V_0} \prod_{i \in I}\D F\big(\lambda_0(i), \mu_0(i)\big). \]
\end{proposition}

\begin{proof}
Let the operation $f \colon \prod_{i \in I}\big(\lambda_0(i) \times \mu_0(i)\big) \sto \prod_{i \in I}\lambda_0(i) \times 
\prod_{i \in I}\mu_0(i)$ be defined by $f(\Phi) := \big(\pr_1(\Phi), \pr_2(\Phi)\big)$, for every 
$\Phi \in \prod_{i \in I}\big(\lambda_0(i) \times \mu_0(i)\big)$, where
$\pr_1(\Phi)_i := \pr_1(\Phi_i)$ and $\pr_2(\Phi)_i := \pr_2(\Phi_i)$, for every $i \in I$.
Using Definition~\ref{def: newfamsofsets1}(i),
\[  \pr_1(\Phi)_j := \pr_1(\Phi_j) = \pr_1\big((\lambda_1 \times \mu_1)_{ij}(\Phi_i)\big) :=
\pr_1\big(\lambda_{ij}(\pr_1(\Phi)_i), \mu_{ij}(\pr_2(\Phi)_i\big) := \lambda_{ij}(\pr_1(\Phi)_i),
\]
hence $\pr_1(\Phi) \in \prod_{i \in I}\lambda_0(i)$.
Similarly, $\pr_2(\Phi) \in \prod_{i \in I}\mu_0(i)$. It is immediate to show that the operation
$f$ is a function.
Let $g \colon \prod_{i \in I}\lambda_0(i) \times \prod_{i \in I}\mu_0(i) \sto 
 \prod_{i \in I}\big(\lambda_0(i) \times \mu_0(i)\big)$ be defined by 
 $g(\Psi, \Xi) := \Phi$, for every $\Psi \in \prod_{i \in I}\lambda_0(i)$ and $\Xi \in \prod_{i \in I}\mu_0(i)$,
 where $\Phi_i := (\Psi_i, \Xi_i)$, for every $i \in I$.
 We show that $g$ is well-defined i.e., $\Phi \in \prod_{i \in I}\big(\lambda_0(i) \times \mu_0(i)\big)$.
 If $i =_I j$, then
 \[ (\lambda_1 \times \mu_1)_{ij}(\Phi_i) := (\lambda_1 \times \mu_1)_{ij}(\Psi_i, \Xi_i) := \big(\lambda_{ij}(\Psi_i),
  \mu_{ij}(\Xi_i)\big) = (\Psi_j, \Xi_j) := \Phi_j.
 \]
Clearly, $f, g$ are inverse to each other. For the equality
$\Map_I(\Lambda, M) =_{\D V_0} \prod_{i \in I}\D F\big(\lambda_0(i), \mu_0(i)\big)$, we use Coroallry~\ref{cor: mapdependent}
and the corresponding identity maps are its witnesses.
\end{proof}

\section{Subfamilies of families of sets}
\label{sec: subfam}

\begin{definition}\label{def: subfamsets}
Let $\Lambda := (\lambda_0, \lambda_1) \in \Fam(I)$ and $h \colon J \to I$. The pair
$\Lambda \circ h := (\lambda_0 \circ h, \lambda_1 \circ h)$, defined in Definition~\ref{def: newfamsofsets1},
is called the $h$-subfamily of $\Lambda$\index{$h$-subfamily of a family of sets}, and we write
$(\Lambda \circ h)_J \leq \Lambda_I$\index{$(\Lambda \circ h)_J \leq \Lambda_I$}. If $J := \Nat$,
we call $\Lambda \circ h$ the $h$-subsequence of $\Lambda$\index{$h$-subsequence of $\Lambda$}.

\end{definition}

\begin{remark}\label{rem: sub1}
If $\Lambda \in \Set(I)$, then $\Lambda \circ h \in \Set(J)$ if and only $h$ is an embedding.
\end{remark}

\begin{proof}
Let $\Lambda \circ h \in \Set(J)$ and $h(j) =_I h(j{'})$, 
hence
$(\lambda_{h(j)h(j{'})}, \lambda_{h(j{'})h(j)}) \colon \lambda_0(h(j)) =_{\D V_0} \lambda_0(h(j{'}))$, 
and $j =_J j{'}$. If $h$ is an embedding and $(\lambda_0 \circ h)(j) =_{\D V_0} (\lambda \circ h)(j{'}) 
:\TOT \lambda_0(j(j)) =_{\D V_0} \lambda_0(h(j{'}))$, then $h(j) =_I h(j{'})$, since $\Lambda \in \Set(I)$, 
and hence $j =_J j{'}$.
\end{proof}

\begin{remark}\label{rem: sub2}
Let $\Lambda, M \in \Fam(I)$, $h \in \D F(J, I)$ and $g \in \D F(I, K)$. \\[1mm]
\normalfont (i)
\itshape $\Lambda \circ \id_I := \Lambda$.\\[1mm]
\normalfont (ii)
\itshape $(\Lambda \circ g) \circ h := \Lambda \circ (g \circ h)$.\\[1mm]
\normalfont (iii)
\itshape If $\Phi \colon \Lambda \To M$, then $\Phi \circ h \colon \Lambda \circ h \To M \circ h$, where 
\[ (\Phi \circ h)_j \colon \lambda_0(h(j)) \to \mu_0(h(j)), \ \ \ \ (\Phi \circ h)_j := \Phi_{h(j)}; \ \ \ \ j \in J. \]
\normalfont (iv)
\itshape If $\Phi \colon \Lambda \To M$, then $\Phi^h \colon \Lambda \circ h \stackrel{h} \To M$, where 
\[ \Phi_j^h \colon \lambda_0(h(j)) \to \mu_0(h(j)), \ \ \ \ \Phi_j^h := \Phi_{h(j)}; \ \ \ \ j \in J. \]
\normalfont (iv)
\itshape $(\Lambda \circ h) \times (M \circ h) := (\Lambda \times M) \circ h$.\\[1mm]
\normalfont (v)
\itshape $\D F\big((\Lambda \circ h), (M \circ h)\big) := \D F(\Lambda, M) \circ h$.
\end{remark}

\begin{proof}
All cases are straightforward to show.
\end{proof}

\begin{proposition}\label{prp: sub3}
Let $\Lambda \in \Fam(I)$, and $h \colon J \to I$. \\[1mm]
\normalfont (i)
\itshape The operation $\sum_h \colon \sum_{j \in J} \lambda_0(h(j)) \sto \sum_{i \in I} \lambda_0(i)$, defined by
\[   \sum_h(j, u) := (h(j), u); \ \ \ \ (j, u) \in \sum_{j \in J}\lambda_0(h(j)), \] 
is a function, and it is an embedding if $h$ is an embedding.\\[1mm]
\normalfont (ii)
\itshape The operation $\prod_h \colon \prod_{i \in I} \lambda_0(i) \sto \prod_{j \in J} \lambda_0(h(j))$,
defined by
\[   \Phi \mapsto \prod_h(\Phi), \ \ \ \ \bigg(\prod_h \Phi\bigg)_j := \Phi_{h(j)}; \ \ \ \ \Phi \in 
\prod_{i \in I}\lambda_0(i), \ j \in J, \] 
is a function, and if $h$ is an embedding, then $\prod_h$ is an embedding.
\end{proposition}

\begin{proof}
(i) By definition we have that 
\[ (j, u) =_{\mathsmaller{\sum_{j \in J}(\lambda_0 \circ h)(j)}} (j{'}, u{'}) :\TOT j =_J j{'} \ \& \ 
\lambda_{h(j)h(j{'})}(u) =_{\mathsmaller{\lambda_0(h(j{'})}} u{'}, \]
\[ (h(j), u) =_{\mathsmaller{\sum_{i \in I}\lambda_0 (i)}} (h(j{'}), u{'}) :\TOT h(j) =_I h(j{'}) \ \& \ 
\lambda_{h(j)h(j{'})}(u) =_{\mathsmaller{\lambda_0(h(j{'})}} u{'}. \]
Since $h$ is a function, the operation $\sum_h$ is a function. If $h$ is an embedding, it is immediate to
show that $\sum_h$ is an embedding. \\
(ii) First we show that $\prod_h$ is well-defined. If $j =_J j{'}$, then 
\[ \bigg(\prod_h \Phi\bigg)_{j'} := \Phi_{h(j{'})} =_{\mathsmaller{\lambda_0(h(j{'}))}} \lambda_{h(j)h(j{'})}
\big(\Phi_{h(j)}\big) := (\lambda_1 \circ h)_{jj{'}}\bigg(\bigg(\prod_h \Phi\bigg)_j\bigg). \]
It is immediate to show that $\prod_h$ is a function. Let $h$ be a surjection and let $\Phi, \Theta \in 
\prod_{i \in I}\lambda_0(i)$ such that $\prod_h(\Phi) =_{\mathsmaller{\prod_{j \in J}\lambda_0(h(j))}} \prod_h(\Theta)$.
If $i \in I$, let $j \in J$ with $h(j) =_I i$. As 
$\Phi_i =_{\lambda_0(i)} \lambda_{h(j)i}\Phi_{h(j)}$ and $\Theta_i =_{\lambda_0(i)} \lambda_{h(j)i}\Theta_{h(j)}$, 
and since $\Phi_{h(j)} =_{\lambda_0(h(j))} \Theta_{h(j)}$, we get $\Phi_i =_{\lambda_0(i)} \Theta_i$.
\end{proof}

\section{Families of sets over products}
\label{sec: famoverprod}

\begin{proposition}\label{prp: famsetsprod1}
Let $\Lambda := (\lambda_0, \lambda_1), K := (k_0, k_1) \in \Fam(I)$ and $M := (\mu_0, \mu_1), N := (\nu_0, \nu_1)
\in \Fam(J)$.\\[1mm]
\normalfont (i) 
\itshape $\Lambda \otimes M := (\lambda_0 \otimes \mu_0, \lambda_1 \otimes \mu_1) \in \Fam(I \times J)$, where 
$\lambda_0 \otimes \mu_0 \colon I \times J \sto \D V_0$ is defined by 
\[ (\lambda_0 \otimes \mu_0)(i,j) := \lambda_0(i) \times \mu_0(j); \ \ \ \ (i,j) \in I \times J, \]
\[  (\lambda_1 \otimes \mu_1)_{(i,j)(i{'}j{'})} \colon \lambda_0(i) \times \mu_0(j) \to \lambda_0(i{'})
\times \mu_0(j{'}),  \]
\[ (\lambda_1 \otimes \mu_1)_{(i,j)(i{'}j{'})}(u, w) := \big(\lambda_{ii{'}}(u), \mu_{jj{'}}(w)\big); \ \ \ \ 
(u,w) \in \lambda_0(i) \times \mu_0(j).
\]
\normalfont (ii) 
\itshape If $\Phi \colon \Lambda \To K$ and $\Psi \colon M \To N$, then $\Phi \otimes \Psi \colon \Lambda \otimes M \To 
K \otimes N$, where, for every $(i,j) \in I \times J$,
\[ (\Phi \otimes \Psi)_{(i,j)} \colon \lambda_0(i) \times \mu_0(j) \to k_0(i) \times \nu_0(j),  \]
\[  (\Phi \otimes \Psi)_{(i,j)}(u,w) := \big(\Phi_i(u), \Psi_j(w)\big); \ \ \ \ (u,w) \in \lambda_0(i) \times \mu_0(j).
\]
\normalfont (iii) 
\itshape The following equalities hold
\[ \sum_{(i,j) \in I \times J}\big(\lambda_0(i) \times \mu_0(j)\big) =_{\D V_0} \bigg(\sum_{i \in I}\lambda_0(i)\bigg)
 \times \bigg(\sum_{j \in J}\mu_0(j)\bigg), \]
\[ \prod_{(i,j) \in I \times J}\big(\lambda_0(i) \times \mu_0(j)\big) =_{\D V_0} \bigg(\prod_{i \in I}\lambda_0(i)\bigg)
 \times \bigg(\prod_{j \in J}\mu_0(j)\bigg). \]
\end{proposition}

\begin{proof}
(i) The proof is straightforward.\\
(ii) We show the required following commutativity by the following supposed ones by
\begin{center}
\begin{tikzpicture}

\node (E) at (0,0) {$k_0(i) \times \nu_0(j)$};
\node[right=of E] (K) {};
\node[right=of K] (L) {};
\node[right=of L] (F) {$k_0(i{'}) \times \nu_0(j{'})$};
\node[above=of F] (A) {$\lambda_0(i{'}) \times \mu_0(j{'})$};
\node [above=of E] (D) {$\lambda_0(i) \times \mu_0(j)$};

\draw[->] (E)--(F) node [midway,below]{$(k_1 \otimes \nu_1)_{(i,j)(i{'},j{'})}$};
\draw[->] (D)--(A) node [midway,above] {$(\lambda_1 \otimes \mu_1)_{(i,j)(i{'},j{'})}$};
\draw[->] (D)--(E) node [midway,left] {$(\Phi \otimes \Psi)_{(i,j)}$};
\draw[->] (A)--(F) node [midway,right] {$(\Phi \otimes \Psi)_{(i{'},j{'})}$};

\end{tikzpicture}
\end{center}
\begin{center}
\begin{tikzpicture}

\node (E) at (0,0) {$k_0(i)$};
\node[right=of E] (F) {$k_0(i{'})$};
\node[above=of F] (A) {$\lambda_0(i{'})$};
\node [above=of E] (D) {$\lambda_0(i)$};
\node [right=of F] (G) {$\mu_0(j)$};
\node [right=of G] (H) {$\mu_0(j{'})$,};
\node[above=of H] (K) {$\mu_0(j{'})$};
\node [above=of G] (L) {$\mu_0(j)$};

\draw[->] (E)--(F) node [midway,below]{$k_{ii{'}}$};
\draw[->] (D)--(A) node [midway,above] {$\lambda_{ii{'}}$};
\draw[->] (D)--(E) node [midway,left] {$\Phi_i$};
\draw[->] (A)--(F) node [midway,right] {$\Phi_{i'}$};
\draw[->] (G)--(H) node [midway,below] {$\nu_{jj{'}}$};
\draw[->] (L)--(K) node [midway,above] {$\mu_{jj{'}}$};
\draw[->] (L)--(G) node [midway,left] {$\Psi_j$};
\draw[->] (K)--(H) node [midway,right] {$\Psi_{j{'}}$};

\end{tikzpicture}
\end{center}
\begin{align*}
(\Phi \otimes \Psi)_{(i{'},j{'})}\big((\lambda_1 \otimes \mu_1)_{(i,j)(i{'},j{'})}\big) & :=
(\Phi \otimes \Psi)_{(i{'},j{'})}\big(\lambda_{ii{'}}(u), \mu_{jj{'}}(w)\big)\\
& := \big(\Phi_{i{'}}\big(\lambda_{ii{'}}(u)\big), \Psi_{j{'}}\big(\mu_{jj{'}}(w)\big)\big)\\
& = \big(k_{ii{'}}\big(\Phi_i(u)\big), \nu_{jj{'}}\big(\Psi_j(w)\big)\big)\\
& := (k_1 \otimes \nu_1)_{(i,j)(i{'},j{'})}\big(\Phi_i(u), \Psi_j(w)\big)\\
& := (k_1 \otimes \nu_1)_{(i,j)(i{'},j{'})}\big((\Phi \otimes \Psi)_{(i,j)}(u,w)\big).
\end{align*}
(iii) For the equality on $\sum_{(i,j) \in I \times J}\big(\lambda_0(i) \times \mu_0(j)\big)$ we have that 
\[ \big((i,j), (u,w)\big) =_{\mathsmaller{\sum_{(i,j) \in I \times J}(\lambda_0(i) \times \mu_0(j))}}
 \big((i{'},j{'}), (u{'},w{'})\big) :\TOT i =_I i{'} \ \& \ j =_J j{'} \ \& \ \]
\[ (\lambda_1 \otimes \mu_1)_{(i,j)(i{'},j{'})}(u,w) =_{\mathsmaller{\lambda_0(i{'}) \times \mu_0(j{'})}} (u{'}, w{'})
: \TOT \lambda_{ii{'}}(u) =_{\mathsmaller{\lambda_0(i{'})}} u{'} \  \& \ \mu_{jj{'}}(w) =_{\mathsmaller{\mu_0(j{'})}} w{'}. 
\]
For the equality on $\big(\sum_{i \in I}\lambda_0(i)\big) \times \big(\sum_{j \in J}\mu_0(j)\big)$ we have that 
\[ \big((i,u), (j,w)\big) =_{\mathsmaller{\big(\sum_{i \in I}\lambda_0(i)\big) \times \big(\sum_{j \in J}\mu_0(j)\big)}}
 \big((i{'},u{'}), (j{'},w{'})\big) :\TOT \]
 \[ (i, u) =_{\mathsmaller{\sum_{i \in I}\lambda_0(i)}} (i{'}, u{'}) \ \& \
 (j, w) =_{\mathsmaller{\sum_{j \in J}\mu_0(j)}} (j{'}, w{'}), \]
 i.e., if $i =_I i{'}$ and $\lambda_{ii{'}}(u) =_{\mathsmaller{\lambda_0(i{'})}} u{'}$, and $j =_J j{'}$
 and $\mu_{jj{'}}(w) =_{\mathsmaller{\mu_0(j{'})}} w{'}$. As the equality conditions for the two sets are equivalent, 
 the operation $\phi \colon \big(\sum_{i \in I}\lambda_0(i)\big) \times \big(\sum_{j \in J}\mu_0(j)\big) \sto 
 \sum_{(i,j) \in I \times J}\big(\lambda_0(i) \times \mu_0(j)\big)$, defined by the rule
 $\big((i,u), (j,w)\big) \mapsto \big((i,j), (u,w)\big)$, together with the operation $\theta \colon 
 \sum_{(i,j) \in I \times J}\big(\lambda_0(i) \times \mu_0(j)\big) \sto 
 \big(\sum_{i \in I}\lambda_0(i)\big) \times \big(\sum_{j \in J}\mu_0(j)\big)$, defined by the inverse rule
 $\big((i,j), (u,w)\big) \mapsto \big((i,u), (j,w)\big)$, are well-defined functions that witness the required equality
 of the two sets in $\D V_0$.\\
 (iv) We proceed similarly to the proof of Proposition~\ref{prp: MAPdep}. 
\end{proof}

Next we define new families of sets generated by a given family of sets indexed by the product
$X \times Y$ of $X$ and $Y$. These families will also be used in section~\ref{sec: bhkbst}.

\begin{definition}\label{def: ac1}
Let $X, Y$ be sets, and let $R := (\rho_0, \rho_1)$ be an $(X \times Y)$-family of sets.\\[1mm]
\normalfont (i) 
\itshape
If $x \in X$, the $x$-component of $R$\index{$x$-component of a family on $X \times Y$} is the pair\index{$R^x$} 
$R^x := (\rho_0^x, \rho_1^x)$, where the assignment routines $\rho_0^x \colon Y \sto \D V_0$ and
$\rho_1^x \colon \bigcurlywedge_{(y, y{'}) \in D(Y)}\D F\big(\rho_0^x(y), \rho_0^x(y{'})\big)$ are
defined by $\rho_0^x(y) := \rho_0(x, y)$, for every $y \in Y$, and 
$\rho_1^x (y, y{'}) := \rho_{yy{'}}^x := \rho_{(x,y)(x,y{'})}$, for every $(y, y{'}) \in D(Y)$.\\[1mm]
\normalfont (ii) 
\itshape If $y \in Y$, the $y$-component of $R$\index{$y$-component of a family on $X \times Y$}
is the pair\index{$R^y$} 
$R^y := (\rho_0^y, \rho_1^y)$, where the assignment routines $\rho_0^y \colon Y \sto \D V_0$ and
$\rho_1^y \colon \bigcurlywedge_{(x, x{'}) \in D(X)}\D F\big(\rho_0^y(x), \rho_0^y(x{'})\big)$ are
defined by $\rho_0^y(x) := \rho_0(x, y)$, for every $x \in X$, and
$\rho_1^y (x, x{'}) := \rho_{xx{'}}^y := \rho_{(x,y)(x{'},y)}$, for every $(x, x{'}) \in D(X)$.\\[1mm]
\normalfont (iii) 
\itshape
Let $\sum^1R := (\sum^1 \rho_0, \sum^1 \rho_1)$, where\index{$\sum^1R$} 
$\sum^1 \rho_0 : X \sto \D V_0$ and 
\[\sum^1 \rho_1 \colon \bigcurlywedge_{(x, x{'}) \in D(X)}\D F\bigg(\big(\sum^1 \rho_0\big)(x), 
\big(\sum^1 \rho_0\big)(x{'})\bigg) 
\ \ \ \ \mbox{are defined by } \]
\[ \bigg(\sum^1 \rho_0\bigg) (x) := \sum_{y \in Y}\rho_0^x(y) := \sum_{y \in Y}\rho_0(x, y); \ \ \ \ x \in X, \] 
\[ \bigg(\sum^1 \rho_1\bigg) (x, x{'}) := \bigg(\sum^1 \rho_1\bigg)_{xx{'}} \colon 
\sum_{y \in Y}\rho_0 (x, y) \to \sum_{y \in Y}\rho_0 (x{'}, y); \ \ \ \ 
(x, x{'}) \in D(X), \]
\[ \bigg(\sum^1 \rho_1\bigg)_{xx{'}}(y, u) := \big(y, \rho_{(x, y)(x{'},y)}(u)\big); \ \ \ \ (y, u) \in  
\sum_{y \in Y}\rho_0 (x, y). \] 
\normalfont (iv) 
\itshape
Let $\sum^2R := (\sum^2 \rho_0, \sum^2 \rho_1)$, where\index{$\sum^2R$} 
$\sum^2 \rho_0 : Y \sto \D V_0$ and 
\[\sum^2 \rho_1 \colon \bigcurlywedge_{(y, y{'}) \in D(X)}\D F\bigg(\big(\sum^2 \rho_0\big)(y), \big(\sum^2 
\rho_0\big)(y{'})\bigg)
\ \ \ \ \mbox{are defined by } \]
\[ \bigg(\sum^2 \rho_0\bigg) (y) := \sum_{x \in X}\rho_0^y(x) := \sum_{x \in X}\rho_0(x, y); \ \ \ \ y \in Y, \] 
\[ \bigg(\sum^2 \rho_1\bigg) (y, y{'}) := \bigg(\sum^2 \rho_1\bigg)_{yy{'}} \colon 
\sum_{x \in X}\rho_0 (x, y) \to \sum_{x \in X}\rho_0 (x, y{'}); \ \ \ \ 
(y, y{'}) \in D(Y), \]
\[ \bigg(\sum^2 \rho_1\bigg)_{yy{'}}(x, w) := \big(x, \rho_{(x, y)(x,y{'})}(w)\big); \ \ \ \ (x, w) \in  
\sum_{x \in X}\rho_0 (x, y). \] 
\normalfont (v) 
\itshape
Let $\prod^1R := (\prod^1 \rho_0, \prod^1 \rho_1)$, where\index{$\prod^1R$} 
$\prod^1 \rho_0 : X \sto \D V_0$ and 
\[\prod^1 \rho_1 \colon \bigcurlywedge_{(x, x{'}) \in D(X)}\D F\bigg(\big(\prod^1 \rho_0\big)(x), 
\big(\prod^1 \rho_0\big)(x{'})\bigg) 
\ \ \ \ \mbox{are defined by } \]
\[ \bigg(\prod^1 \rho_0\bigg) (x) := \prod_{y \in Y}\rho_0^x(y) := \prod_{y \in Y}\rho_0(x, y); \ \ \ \ x \in X, \] 
\[ \bigg(\prod^1 \rho_1\bigg) (x, x{'}) := \bigg(\prod^1 \rho_1\bigg)_{xx{'}} \colon 
\prod_{y \in Y}\rho_0 (x, y) \to \prod_{y \in Y}\rho_0 (x{'}, y); \ \ \ \ 
(x, x{'}) \in D(X), \]
\[ \bigg[\bigg(\prod^1 \rho_1\bigg)_{xx{'}}(\Theta)\bigg]_y := \rho_{(x, y)(x{'},y)}(\Theta_y)\big); \ \ \ \ \Theta \in  
\prod_{y \in Y}\rho_0 (x, y), \ y \in Y. \] 
\normalfont (vi) 
\itshape
Let $\prod^2R := (\prod^2 \rho_0, \prod^2 \rho_1)$, where\index{$\prod^2R$} 
$\prod^2 \rho_0 : Y \sto \D V_0$ and 
\[\prod^2 \rho_1 \colon \bigcurlywedge_{(y, y{'}) \in D(X)}\D F\bigg(\big(\prod^2 \rho_0\big)(y), 
\big(\prod^2 \rho_0\big)(y{'})\bigg)
\ \ \ \ \mbox{are defined by } \]
\[ \bigg(\prod^2 \rho_0\bigg) (y) := \prod_{x \in X}\rho_0^y(x) := \prod_{x \in X}\rho_0(x, y); \ \ \ \ y \in Y, \] 
\[ \bigg(\prod^2 \rho_1\bigg) (y, y{'}) := \bigg(\prod^2 \rho_1\bigg)_{yy{'}} \colon 
\prod_{x \in X}\rho_0 (x, y) \to \prod_{x \in X}\rho_0 (x, y{'}); \ \ \ \ 
(y, y{'}) \in D(Y), \]
\[ \bigg[\bigg(\prod^2 \rho_1\bigg)_{yy{'}}(\Phi)\bigg]_x := \rho_{(x, y)(x,y{'})}(\Phi_x)\big); \ \ \ \ \Phi \in  
\prod_{x \in X}\rho_0 (x, y), \ x \in X. \] 
\end{definition} 
 
It is easy to show
that $R^y, \sum^1R, \prod^1R  \in \Fam(X)$ and $R^x, \sum^2R, \prod^2R \in \Fam(Y)$.

\begin{proposition}\label{prp: newfamilymaps3}
Let $X, Y \in \D V_0$, $R := (\rho_0, \rho_1), S:= (\sigma_0, \sigma_1) \in \Fam(X \times Y)$, and
$\Phi \colon R \To S$.\\[1mm]
\normalfont (i) 
\itshape Let $\Phi^x \colon \bigcurlywedge_{y \in Y}\D F\big(\rho_0^x(y), \sigma_0^x(y)\big)$, where 
$\Phi^x_y := \Phi_{(x,y)} \colon \rho_0^x(y) \to \sigma^x_0(y)$\index{$\Phi^x$}.\\[1mm]
\normalfont (ii) 
\itshape Let $\Phi^y \colon \bigcurlywedge_{x \in X}\D F\big(R^y(x), S^y(x)\big)$, where 
$\Phi^y_x := \Phi_{(x,y)} \colon \rho_0^y(x) \to \sigma^y_0(x)$\index{$\Phi^y$}.\\[1mm]
\normalfont (iii) 
\itshape Let $\sum^1 \Phi \colon \bigcurlywedge_{x \in X}\D F\big(\big(\sum^1\rho_0\big)(x),
\big(\sum^1\sigma_0\big)(x)\big)$, where, for every $x \in X$, we define\index{$\sum^1 \Phi$} 
\[ \bigg(\sum^1 \Phi\bigg)_x \colon \sum_{y \in Y}\rho_0^y(x) \to \sum_{y \in Y}\sigma^x_0(y) \]
\[ \bigg(\sum^1 \Phi\bigg)_x(y, u) := \big(y, \Phi_{(x,y)}(u)\big); \ \ \ \ (y, u) \in \sum_{y \in Y}\rho_0(x, y). \]
\normalfont (iv) 
\itshape If $\sum^2 \Phi \colon \bigcurlywedge_{y \in Y}\D F\big(\big(\sum^2\rho_0\big)(y),
\big(\sum^2\sigma_0\big)(y)\big)$, where, for every $y \in Y$, we define\index{$\sum^2 \Phi$} 
\[ \bigg(\sum^2 \Phi\bigg)_y \colon \sum_{x \in X}\rho_0^x(y) \to \sum_{x \in X}\sigma^y_0(x) \]
\[ \bigg(\sum^2 \Phi\bigg)_y(x, w) := \big(x, \Phi_{(x,y)}(w)\big); \ \ \ \ (x, w) \in \sum_{x \in X}\rho(x, y). \]
\normalfont (v) 
\itshape Let $\prod^1 \Phi \colon \bigcurlywedge_{x \in X}\D F\big(\big(\prod^1\rho_0\big)(x),
\big(\prod^1\sigma_0\big)(x)\big)$, where, for every $x \in X$, we define\index{$\prod^1 \Phi$} 
\[ \bigg(\prod^1 \Phi\bigg)_x \colon \prod_{y \in Y}\rho_0^y(x) \to \prod_{y \in Y}\sigma^x_0(y) \]
\[ \bigg[\bigg(\prod^1 \Phi\bigg)_x(\Theta)\bigg]_y := \Phi_{(x,y)}\big(\Theta_y))\big); \ \ \ \ \Theta \in
\prod_{y \in Y}\rho_0(x, y). \]
\normalfont (vi) 
\itshape Let $\prod^2 \Phi \colon \bigcurlywedge_{y \in Y}\D F\big(\big(\prod^2\rho_0\big)(y),
\big(\prod^2\sigma_0\big)(y)\big)$, where, for every $y \in Y$, we define\index{$\prod^2 \Phi$} 
\[ \bigg(\prod^2 \Phi\bigg)_y \colon \prod_{x \in X}\rho_0^x(y) \to \prod_{x \in X}\sigma^y_0(x) \]
\[ \bigg[\bigg(\prod^2 \Phi\bigg)_y(\Theta)\bigg]_x := \Phi_{(x,y)}\big(\Theta_x)\big); \ \ \ \ \Theta \in 
\prod_{x \in X}\rho_0(x, y). \]
Then $\Phi^x \colon R^x \To S^x$, $\Phi^y \colon R^y \To S^y$, 
$\sum^1 \Phi \colon \big(\sum^1 R\big) \To \big(\sum^1 S\big)$,
$\sum^2 \Phi \colon \big(\sum^2 R\big) \To \big(\sum^2 S\big)$,
$\prod^1 \Phi \colon \big(\prod^1 R\big) \To \big(\prod^1 S\big)$, and
$\prod^2 \Phi \colon \big(\prod^2 R\big) \To \big(\prod^2 S\big)$.
\end{proposition}

\begin{proof}
The proofs of (ii), (iv) and (vi) are like the proofs of (i), (iii), and (v), respectively.\\
(i) It is immediate to show that the operation $\Phi^x_y \colon \rho_0^x(y) \sto \sigma^x_0(y)$ 
is a function. If $y =_Y y{'}$, the commutativity of the following left diagram from the hypothesis $\Phi \colon R \To S$
\begin{center}
\begin{tikzpicture}

\node (E) at (0,0) {$\sigma_0(x,y)$};
\node[right=of E] (X) {};
\node[right=of X] (F) {$\sigma_0(x,y{'})$};
\node[above=of F] (A) {$\rho_0(x,y{'})$};
\node [above=of E] (D) {$\rho_0(x,y)$};
\node [right=of F] (G) {$\sigma_0^x(y)$};
\node [right=of G] (H) {$\sigma_0^x(y{'})$};
\node[above=of H] (K) {$\rho_0^x(y{'})$};
\node [above=of G] (L) {$\sigma_0^x(y)$};

\draw[->] (E)--(F) node [midway,below]{$\sigma_{(x,y)(x,y{'})}$};
\draw[->] (D)--(A) node [midway,above] {$\rho_{(x,y)(x,y{'})}$};
\draw[->] (D)--(E) node [midway,left] {$\Phi_{(x,y)}$};
\draw[->] (A)--(F) node [midway,right] {$\Phi_{(x,y{'})}$};
\draw[->] (G)--(H) node [midway,below] {$\sigma^x_{yy{'}}$};
\draw[->] (L)--(K) node [midway,above] {$\rho^x_{yy{'}}$};
\draw[->] (L)--(G) node [midway,left] {$\Phi_y^x$};
\draw[->] (K)--(H) node [midway,right] {$\Phi^x_{y{'}}$};

\end{tikzpicture}
\end{center}
implies the required commutativity of the right above diagram, as these are the same diagrams.\\
(iii) First we explain why the operation $\big(\sum^1 \Phi\big)_x$ is a function. If
\[ (y, u) =_{\mathsmaller{\sum_{y \in Y}\rho_0^y(x)}} (y{'}, u{'}) :\TOT y =_Y y{'} \ \& \ \rho_{(x,y)(x,y{'})}(u) 
=_{\mathsmaller{\rho_0(x, y{'})}} u{'}, \]
\[  \big(y, \Phi_{(x,y)}(u)\big) =_{\mathsmaller{\sum_{y \in Y}\sigma_0^y(x)}}  \big(y, \Phi_{(x,y)}(u)\big) :\TOT 
 y =_Y y{'} \ \& \ \sigma_{(x,y)(x,y{'})}\big(\Phi_{(x,y)}(u)\big) =_{\mathsmaller{\sigma_0(x, y{'})}}
 \Phi_{(x,y{'})}(u{'}).
\]
From our hypothesis the second equality is equivalent to 
\[ \sigma_{(x,y)(x,y{'})}\big(\Phi_{(x,y)}(u)\big) =_{\mathsmaller{\sigma_0(x, y{'})}}
\Phi_{(x,y{'})}\big( \rho_{(x,y)(x,y{'})}(u)\big), \]
which is the commutativity of the above left diagram. If $x =_X x{'}$, and since
$\Phi_{(x{'},y)} \circ \rho_{(x,y)(x{'},y)} = \sigma_{(x,y)(x{'},y)} \circ \Phi_{(x,y)}$ 
%
%
%
we get the commutativity of the following left diagram by 
\begin{center}
\resizebox{14cm}{!}{%
\begin{tikzpicture}

\node (E) at (0,0) {$\mathsmaller{\sum_{y \in Y}\sigma_0(x,y)}$};
\node[right=of E] (X) {};
\node[right=of X] (F) {$\mathsmaller{\sum_{y \in Y}\sigma_0(x{'},y)}$};
\node[above=of F] (Y) {};
\node[above=of Y] (A) {$\mathsmaller{\sum_{y \in Y}\rho_0(x{'},y)}$};
\node [above=of E] (Z) {};
\node [above=of Z] (D) {$\mathsmaller{\sum_{y \in Y}\rho_0(x,y)}$};

\node[right=of F] (K) {$\mathsmaller{\prod_{y \in Y}\sigma_0(x,y)}$};
\node[right=of K] (Y) {};
\node[right=of Y] (L) {$\mathsmaller{\prod_{y \in Y}\sigma_0(x{'},y)}$};
\node[above=of L] (M) {};
\node[above=of M] (N) {$\mathsmaller{\prod_{y \in Y}\rho_0(x{'},y)}$};
\node [above=of K] (T) {};
\node [above=of T] (R) {$\mathsmaller{\prod_{y \in Y}\rho_0(x,y)}$};

\draw[->] (E)--(F) node [midway,below]{$\mathsmaller{\big(\sum^1 \sigma_1\big)_{xx{'}}}$};
\draw[->] (D)--(A) node [midway,above] {$\mathsmaller{\big(\sum^1 \rho_1\big)_{xx{'}}}$};
\draw[->] (D)--(E) node [midway,left] {$\mathsmaller{\big(\sum^1 \Phi\big)_x}$};
\draw[->] (A)--(F) node [midway,right] {$\mathsmaller{\big(\sum^1 \Phi\big)_{x{'}}}$};

\draw[->] (K)--(L) node [midway,below]{$\mathsmaller{\big(\prod^1 \sigma_1\big)_{xx{'}}}$};
\draw[->] (R)--(N) node [midway,above] {$\mathsmaller{\big(\prod^1 \rho_1\big)_{xx{'}}}$};
\draw[->] (R)--(K) node [midway,left] {$\mathsmaller{\big(\prod^1 \Phi\big)_x}$};
\draw[->] (N)--(L) node [midway,right] {$\mathsmaller{\big(\prod^1 \Phi\big)_{x{'}}}$};

\end{tikzpicture}
}
\end{center}
\begin{align*}
 \bigg(\sum^1 \Phi\bigg)_{x{'}}\bigg(\bigg(\sum^1 \rho_1\bigg)_{xx{'}}(y, u)\bigg) & := 
 \bigg(\sum^1 \Phi\bigg)_{x{'}}\big(y, \rho_{(x,y)(x{'},y)}(u)\big)\\
 & := \big(y, \Phi_{(x{'},y)}\big(\rho_{(x,y)(x{'},y)}(u)\big)\big)\\
 & =  \big(y, \sigma_{(x,y)(x{'},y)}\big(\Phi_{(x,y)}(u)\big)\big)\\
 & :=  \bigg(\sum^1 \sigma_1\bigg)_{xx{'}}\big(y, \Phi_{(x,y)}(u)\big)\big)\\
 & :=  \bigg(\sum^1 \sigma_1\bigg)_{xx{'}}\bigg(\bigg(\sum^1 \Phi\bigg)_x(y, u)\bigg).
\end{align*}
(v) First we explain why the operation $\big(\prod^1 \Phi\big)_x$ is well-defined. If $y =_Y y{'}$ and 
$\Theta \in \prod_{y \in Y}\rho_0^x(y)$, then by the 
commutativity of the above left diagram we have that
\begin{align*}
 \bigg[\bigg(\prod^1 \Phi\bigg)_{x}(\Theta)\bigg]_{y{'}} & := \Phi_{(x,y{'})}\big(\Theta_{y{'}}\big)\\
 & = \Phi_{(x,y{'})}\big(\rho_{(xy)(x,y{'})}(\Theta_y)\big)\\
 & = \sigma_{(x,y)(x,y{'})}\big(\Phi_{(x,y)}(\Theta_y)\big)\\
 & :=  \sigma^x_{yy{'}}\bigg(\bigg[\bigg(\prod^1 \Phi\bigg)_x(\Theta)\bigg]_y\bigg).
\end{align*}
Clearly, the operation $\big(\prod^1 \Phi\big)_x$ is a function. If $x =_X x{'}$, 
and by the commutativity of the first diagram in the proof of (iii) we get the commutativity of the
above right diagram 
%
%
%
\begin{align*}
 \bigg[\bigg(\prod^1 \Phi\bigg)_{x{'}}\bigg(\bigg(\prod^1 \rho_1\bigg)_{xx{'}}(\Theta)\bigg)\bigg]_y & := 
  \Phi_{(x{'},y)}\bigg[\bigg(\prod^1 \rho_1\bigg)_{xx{'}}(\Theta)\bigg]_y\bigg)\\
 & := \Phi_{(x{'},y)}\big(\rho_{(x,y)(x{'},y)}(\Theta_y)\big)\\
 & = \sigma_{(x,y)(x{'},y)}\big(\Phi_{(x,y)}(\Theta_y)\big)\\
 & := \sigma_{(x,y)(x{'},y)}\bigg(\bigg[\prod^1 \Phi\bigg)_x(\Theta)\bigg]_y\bigg)\\
 & := \bigg[\bigg(\prod^1 \sigma_1\bigg)_{xx{'}}\bigg(\prod^1 \Phi\bigg)_x(\Theta)\bigg)\bigg]_y.\qedhere
\end{align*}
\end{proof}

\begin{proposition}\label{prp: sigmasigma}
If $R := (\rho_0, \rho_1) \in \Fam(X \times Y)$, the following equalities hold.
\[ \sum_{x \in X}\sum_{y \in Y}\rho_0(x,y) =_{\D V_0} \sum_{y \in Y}\sum_{x \in X}\rho_0(x,y), \]
\[ \prod_{x \in X}\prod_{y \in Y}\rho_0(x,y) =_{\D V_0} \prod_{y \in Y}\prod_{x \in X}\rho_0(x,y). \]
\end{proposition}

\begin{proof}
The proof is straightforward. 
\end{proof}

\section{The distributivity of $\prod$ over $\sum$}
\label{sec: distributivity}

We prove the translation of the type-theoretic axiom of choice in $\BST$ (Theorem~\ref{thm: ac}). 
 
 \begin{lemma}\label{lem: ac3}
Let $R := (\rho_0, \rho_1)$, $R^x := (\rho_0^x, \rho_1^x)$ and $\sum^1R := (\sum^1 \rho_0, \sum^1 \rho_1)$ be
the families of sets of Definition~\ref{def: ac1}. 
If $\Phi \in \prod_{x \in X}\big(\sum^1 \rho_0\big)(x)$,
the operation $f_{\Phi} \colon X \sto Y$, defined by $x \mapsto \pr_1^{R^x}(\Phi_x)$, for every $x \in X$, 
is a function from $X$ to $Y$.
\end{lemma}

\begin{proof}
If $x =_X x{'}$, then $\Phi_{x{'}} = \big(\sum^1 \rho_1\big)_{xx{'}}(\Phi_x)$. 
Since $\Phi_x \in  \sum_{y \in Y}\rho_0(x, y)$, there are
$y \in Y$ and $u \in \rho_0(x, y)$ such that $\Phi_x := (y, u)$. Hence $f_{\Phi}(x) := y$ and
\[ f_{\Phi}(x{'}) := \pr_1^{R^x}(\Phi_{x{'}}) =_Y \pr_1^{R^x}\bigg(\big(\sum^1 \rho_1\big)_{xx{'}}(\Phi_x)\bigg)
:= \pr_1^{R^x}\big(y, \rho_{(x,y)(x{'},y)}(u)\big) := y.\qedhere \]  
\end{proof}

\begin{lemma}\label{lem: ac4}
Let $R := (\rho_0, \rho_1)$, $R^x := (\rho_0^x, \rho_1^x)$ and $\sum^1R := (\sum^1 \rho_0, \sum^1 \rho_1)$ be
as above. 
If $f \colon X \to Y$, the pair $N^f := \big(\nu_0^f, \nu_1^f\big)$ is an $X$-family of sets, where
the assignment routines $\nu_0^f : X \sto \D V_0$ and 
$\nu_1^f \colon  \bigcurlywedge_{(x, x{'}) \in D(X)}\D F\big(\nu_0^f (x), \nu_0^f (x{'})\big)$
are given by $\nu_0^f (x) := \rho_0 (x, f(x))$, for every $x \in X$, and 
$\nu_{xx{'}}^f := \rho_{(x, f(x))(x{'}, f(x{'}))}$, for every $(x, x{'}) \in D(X)$,
\end{lemma}

\begin{proof}
 The proof is straightforward (see also~\cite{Pe19c}, p.~12).
\end{proof}

\begin{lemma}\label{lem: ac5}
If $R := (\rho_0, \rho_1)$ and $N^f := \big(\nu_0^f, \nu_1^f\big)$ are families of sets as above,
then the pair $\Xi := (\xi_0, \xi_1)$ is an $\D F(X, Y)$-family of sets, where
the assignment routines 
$\xi_0 : \D F(X, Y) \sto \D V_0$ and $\xi_1 \colon \bigcurlywedge_{(f, f{'}) \in D(\D F(X, Y))}\D F\big(\xi_0(f),
\xi_0(f{'})\big)$ are defined by
\[ \xi_0(f) := \prod_{x \in X}\nu_0^f(x) := \prod_{x \in X}\rho_0(x, f(x)); \ \ \ \ f \in \D F(X, Y), \]
\[ \xi_{ff{'}} : \prod_{x \in X}\rho_0(x, f(x)) \to \prod_{x \in X}\rho_0(x, f{'}(x)); \ \ \ \ (f, f{'}) \in D(\D F(X, Y)), \]
\[ \big[\xi_{ff{'}}(H)\big]_x := \rho_{(x, f(x))(x, f{'}(x))}(H_x); \ \ \ \ H \in  \prod_{x \in X}\rho_0(x, f(x)), \  \
x \in X. \]
\end{lemma}
 
\begin{proof}
First we show that the operation $\xi_{ff{'}}$ is well-defined i.e., if 
\[ H \in \prod_{x \in X}\rho_0(x, f(x)) :\TOT \forall_{(x, x{'}) \in D(X)}\big(H_{x{'}} = \nu^f_{xx{'}}(H_x) := 
  \rho_{(x, f(x))(x{'}, f(x{'}))}(H_x)\big), \ \  \ \mbox{then}
\]
\[ \xi_{ff{'}}(H) \in \prod_{x \in X}\rho_0(x, f{'}(x)) :\TOT \forall_{(x, x{'}) 
\in D(X)}\big(\big[\xi_{ff{'}}(H)\big]_{x{'}} = \nu^f_{xx{'}}\big(\big[\xi_{ff{'}}(H)\big]_x\big) 
  \]
If $x =_X x{'}$, then $f(x) =_Y f(x{'}) =_Y f{'}(x{'}) =_Y f{'}(x)$, and
\begin{align*}
 \nu^f_{xx{'}}\big(\big[\xi_{ff{'}}(H)\big]_x\big) & :=  
 \rho_{(x, f{'}(x))(x{'}, f{'}(x{'}))}\big(\big[\xi_{ff{'}}(H)\big]_x\big)\\
 & := \rho_{(x, f{'}(x))(x{'}, f{'}(x{'}))}\big(\rho_{(x, f(x))(x, f{'}(x))}(H_x)\big)\\
 & = \rho_{(x, f(x))(x{'}, f{'}(x{'}))}(H_x)\\
 & =  \rho_{(x{'}, f(x{'}))(x{'}, f{'}(x{'}))}\big(\rho_{(x, f(x))(x{'}, f(x{'}))}(H_x)\big)\\
 & =  \rho_{(x{'}, f(x{'}))(x{'}, f{'}(x{'}))}(H_{x{'}})\\
 & := \big[\xi_{ff{'}}(H)\big]_{x{'}}.
\end{align*}
It is immediate to see to show that $\xi_{ff{'}}$ is a function. If $f \in \D F(X, Y)$, then 
\[ \big[\xi_{ff}(H)\big]_x := \rho_{(x, f(x))(x, f(x))}(H_x) := \id_{\rho_0(x, f(x))}(H_x) := H_x. \]
Moreover, if $f =_{\D F(X, Y)} f{'} =_{\D F(X, Y)} f{''}$, the following diagram is commutative:
\begin{center}
\begin{tikzpicture}

\node (E) at (0,0) {$\xi_0(f{'})$};
\node[right=of E] (F) {$\xi_0(f{''})$};
\node [above=of E] (D) {$\xi_0(f)$};

\draw[->] (E)--(F) node [midway,below] {$\xi_{f{'} f{''}}$};
\draw[->] (D)--(E) node [midway,left] {$\xi_{f f{'}}$};
\draw[->] (D)--(F) node [midway,right] {$\ \xi_{f f{''}}$};

\end{tikzpicture}
\end{center}
\begin{align*}
 \big[(\xi_{f{'} f{''}} \circ \xi_{ff{'}})(H)\big]_x & := \big[\xi_{f{'} f{''}}\big(\xi_{ff{'}}(H)\big)\big]_x\\
 & := \rho_{(x, f{'}(x))(x, f{''}(x))}\big(\big[\xi_{ff{'}}(H)\big]_x\big)\\
 & := \rho_{(x, f{'}(x))(x, f{''}(x))}\big(\rho_{(x, f(x))(x, f{'}(x))}(H_x)\big)\\
 & = \rho_{(x, f(x))(x, f{''}(x))}(H_x)\\
 & := \big[\xi_{ff{''}}(H)\big]_x.\qedhere
\end{align*}
\end{proof}

\begin{theorem}[Distributivity of $\prod$ over $\sum$]\label{thm: ac}
Let $X, Y$ be sets, $R := (\rho_0, \rho_1)$, $R^x := (\rho_0^x, \rho_1^x)$, and 
$\sum^1R := (\sum^1 \rho_0, \sum^1 \rho_1)$ as above. If
\[ \Phi \in \prod_{x \in X}\big(\sum^1 \rho_0\big)(x) := \prod_{x \in X}\sum_{y \in Y}\rho_0(x, y), \ \ \ \mbox{there is} \]
\[ \Theta_{\Phi} \in \prod_{x \in X}\nu_0^{f_{\Phi}}(x) := \prod_{x \in X}\rho_0(x, f_{\Phi}(x)), \]
where $f_{\Phi} \colon X \to Y$ is defined in Lemma~\ref{lem: ac3}.
The following operation is a function:
\[ \ac : \prod_{x \in X}\sum_{y \in Y}\rho_0(x, y) \sto 
\sum_{f \in \D F(X, Y)}\prod_{x \in X}\rho_0(x, f(x)) \]
\[ \Phi \mapsto \big(f_{\Phi}, \Theta_{\Phi}\big);  \ \ \ \ \Phi \in \prod_{x \in X}\sum_{y \in Y}\rho_0(x, y). \]
\end{theorem}
 
\begin{proof}
Since by Remark~\ref{rem: proj2}, 
\[ \pr_2^{R^x}  \in \prod_{w \in \sum_{y \in Y}\rho_0^x (y)}\rho_0^x (\pr_1^{R^x} (w))
:= \prod_{w \in \sum_{y \in Y}\rho_0 (x, y)}\rho_0 (x, \pr_1^{R^x} (w)), \]
\[ \pr_2^{R^x} (\Phi_x) \in \rho_0 (x, \pr_1^{R^x} (\Phi_x)) := \rho_0 (x, f_{\Phi}(x)). \]
Hence, the dependent operation
$\Theta_{\Phi} \in \bigcurlywedge_{x \in X}\nu_0^{f_{\Phi}}(x) := \bigcurlywedge_{x \in X}\rho_0(x, f_{\Phi}(x))$,
defined by
\[ \Theta_{\Phi}(x) := \pr_2^{R^x} (\Phi_x) := u; \ \ \ \ \Phi_x := (y, u), \ \ y := f_{\Phi}(x), \ \ x \in X, \]
is well-defined. To show that $\Theta_{\Phi} \in \prod_{x \in X}\nu_0^{f_{\Phi}}(x)$, let $x =_X x{'}$. Since,
\[ \Theta_{\Phi}(x{'}) := \pr_2^{R^{x{'}}} (\Phi_{x{'}}) := u{'}; \ \ \ \ \Phi_{x{'}} := (y{'}, u{'}), \ \ y{'}
:= f_{\Phi}(x{'}), \]
we need to show that $u{'} =_{\mathsmaller{\rho_0(x{'}, y{'})}} \nu_{xx{'}}^{f_{\Phi}}(u)$. Since 
$\Phi \in \prod_{x \in X}\big(\sum^1 \rho_0\big)(x)$, we have that
\[ (y{'}, u{'}) := (f_{\Phi}(x{'}), u{'}) := \Phi_{x{'}} 
 =_{\mathsmaller{\sum_{y \in Y}\rho_0 (x{'}, y)}} \big(\sum^1 \rho_1\big)_{xx{'}}(\Phi_x) =_{\mathsmaller{\sum_{y \in Y}\rho_0 (x{'}, y)}}
 \big(y, \rho_{(x,y)(x{'},y)}(u)\big). \]
 By the last equality we get $y =_Y y{'}$ and 
 \begin{align*}
  \rho_{yy{'}}^{x{'}}\big(\rho_{(x,y)(x{'},y)}(u)\big) =_{\mathsmaller{\rho_0^{x{'}}(y{'})}} u{'} & :\TOT
  \rho_{(x{'},y)(x{'},y{'})} \big(\rho_{(x,y)(x{'},y)}(u)\big) =_{\mathsmaller{\rho_0(x{'}, y{'})}} u{'}\\
  & \TOT \rho_{(x,y)(x{'},y{'})}(u) =_{\mathsmaller{\rho_0(x{'}, y{'})}} u{'}.
 \end{align*}
Hence, 
\[ \nu_{xx{'}}^{f_{\Phi}}(u) := \rho_{(x,f_{\Phi}(x))(x{'},f_{\Phi}(x{'}))}(u) := \rho_{(x,y)(x{'},y{'})}(u) 
=_{\mathsmaller{\rho_0(x{'}, y{'})}} u{'}. \]
To show that the operation $\ac$ is a function, we suppose that $\Phi  =_{\mathsmaller{\prod_{x \in X}\mu_0 (x)}}
\Phi{'}$, and we show that $\ac(\Phi) =_{\mathsmaller{\sum_{f \in \D F(X, Y)}\xi_0 (f)}} \ac(\Phi{'})$ i.e.,
\[ \big(f_{\Phi}, \Theta_{\Phi}\big) =_{\mathsmaller{\sum_{f \in \D F(X, Y)}\xi_0 (f)}} \big(f_{\Phi{'}}, 
\Theta_{\Phi{'}}\big) :\TOT f_{\Phi} =_{\D F(X, Y)} f_{\Phi{'}} \ \& \ \xi_{f_{\Phi},f_{\Phi{'}}}(\Theta_{\Phi}) =_{\xi_0(f_{\Phi{'}})}
\Theta_{\Phi{'}}. \]
By definition, $\Phi  =_{\mathsmaller{\prod_{x \in X}\big(\sum^1 \rho_0\big)(x)}} \Phi{'}$ if and only if 
$\Phi_x =_{\mathsmaller{\big(\sum^1 \rho_0\big)(x)}} \Phi_x{'}$,
for every $x \in X$. By Lemma~\ref{lem: ac3} 
\[ f_{\Phi}(x) := \pr_1^{R^x}(\Phi_{x}) := y; \ \ \ \ \Phi_x := (y, u), \]
\[ f_{\Phi{'}}(x) := \pr_1^{\Lambda^x}(\Phi{'}_{x}) := y{'}; \ \ \ \ \Phi{'}_x := (y{'}, u{'}). \]
Since $\Phi_x =_{\mathsmaller{\big(\sum^1 \rho_0\big)(x)}} \Phi_x{'}$, we get $y =_Y y{'}$,  and 
$\rho_{(x,y)(x,y{'})}(u) =_{\rho_0(x,y{'})} u{'}$. From the first equality we get and hence $f_{\Phi}(x) =_Y f_{\Phi{'}}(x)$,
and from the second we conclude that 
\[ \big[ \xi_{f_{\Phi},f_{\Phi{'}}}(\Theta_{\Phi}) \big]_x 
:= \rho_{(x,f_{\Phi}(x))(x,f_{\Phi{'}}(x))}\big([\Theta_f]_x\big)
:= \rho_{(x,y)(x,y{'})}(u)
=_{\rho_0(x,y{'})} u{'} := 
\big[\Theta_{\Phi{'}}\big]_x.\qedhere
\]
\end{proof}

\section{Sets of sets}
\label{sec: setofsets}

\begin{definition}\label{def: setofsets}
If $I$ is a set, a \textit{set of sets}\index{set of sets} indexed by $I$, or an $I$-\textit{set 
of sets}\index{$I$-set of sets}, is a pair $\Lambda := (\lambda_0, \lambda_1) \in \Fam(I)$ such that 
the following condition
is satisfied:
\[ Q(\Lambda) :\TOT \forall_{i, j \in I}\big(\lambda_0(i) =_{\D V_0} \lambda_0(j) \To i =_I j\big). \]
Let $\Set(I)$ be their totality, equipped with the canonical equality on $\Fam(I)$. 
\end{definition}

\begin{remark}\label{rem: Isets1}
If $\Lambda \in \Set(I)$ and $M \in \Fam(I)$ such that $\Lambda =_{\Fam(I)} M$, then $M \in \Set(I)$. 
\end{remark}

\begin{proof}
Let $i, j \in I$, $f \colon \mu_0(i) \to \mu_0(j)$ and $g \colon \mu_0(j) \to \mu_0(i)$, such that
$f \circ g = \id_{\mu_0(j)}$ and $g \circ f = \id_{\mu_0(i)}$. It suffices to show that 
$\lambda_0(i) =_{\D V_0} \lambda_0(j)$. Let $\Phi \in \Map_I(\Lambda, M)$ and $\Psi \in \Map_I(M, \Lambda)$ such that 
$\Phi \circ \Psi = \id_M$ and $\Psi \circ \Phi = \id_{\Lambda}$. We define
$f{'} \colon \lambda_0(i) \to \lambda_0(j)$ and $g{'} \colon \lambda_0(j) \to \lambda_0(i)$ by
\[ f{'} := \Psi_j \circ f \circ \Phi_i \ \ \ \& \ \ \ g{'} := \Psi_i \circ g \circ \Phi_j \]
\begin{center}
\begin{tikzpicture}

\node (E) at (0,0) {$\lambda_0(i)$};
\node[right=of E] (F) {$\lambda_0(j)$};
\node[above=of F] (A) {$\mu_0(j)$};
\node [above=of E] (D) {$\mu_0(i)$};
\node [right=of F] (G) {$\lambda_0(i)$};
\node [right=of G] (H) {$\lambda_0(j)$.};
\node[above=of H] (K) {$\mu_0(j)$};
\node [above=of G] (L) {$\mu_0(i)$};

\draw[->] (E)--(F) node [midway,below]{$f{'}$};
\draw[->] (D)--(A) node [midway,above] {$f$};
\draw[->] (E)--(D) node [midway,left] {$\Phi_i$};
\draw[->] (A)--(F) node [midway,right] {$\Psi_j$};
\draw[->] (H)--(G) node [midway,below] {$g{'}$};
\draw[->] (K)--(L) node [midway,above] {$g$};
\draw[->] (L)--(G) node [midway,left] {$\Psi_i$};
\draw[->] (H)--(K) node [midway,right] {$\Phi_j$};

\end{tikzpicture}
\end{center}
It is straightforward to show that $(f{'}, g{'}) \colon \lambda_0(i) =_{\D V_0} \lambda_0(j)$. 
\end{proof}

By the previous remark $Q(\Lambda)$ is an extensional property on $\Fam(I)$. Since $\Set(I)$ is defined by 
separation on $\Fam(I)$, which is impredicative, $\Set(I)$ is also 
an impredicative set. We can also see that by an argument similar to the one used for the 
impredicativity of $\Fam(I)$.

If $X, Y$ are not equal sets in $\D V_0$, then with Ex falsum we get that the $\D 2$-family
$\Lambda^{\D 2}$ of $X$ and $Y$ is a $\D 2$-set of sets.
Similarly, if $X_n$ and $X_m$ are not equal in $\D V_0$, for every $n \neq m$, then with Ex falsum we get that
the $\Nat$-family $\Lambda^{\Nat}$ of $(X_n)_{n \in \Nat}$ is an $\Nat$-set of sets.
If $I$ is a set with $(i, j) \in I \times I$ such that $\neg (i =_I j )$, 
then the constant $I$-family $A$, for some set $A$, is an $I$-family that is not an $I$-set of sets. 
We can easily turn an $I$-family of sets $\Lambda$ into an 
$I$-set of sets. 

\begin{definition}\label{def: lambdaI}
Let $\Lambda := (\lambda_0, \lambda_1) \in \Fam(I)$. The equality $=_I^{\Lambda}$\index{$=_I^{\Lambda}$} on $I$ induced 
by $\Lambda$\index{equality on the index-set induced by a family}
is given by $i =_I^{\Lambda} j : \TOT \lambda_0 (i) =_{\D V_0} \lambda_0 (j)$, for every $i, j \in I$.
The \textit{set $\lambda_0 I$ of sets 
generated by $\Lambda$}\index{set of sets generated by a family of sets}\index{$\lambda_0 I$} is the totality $I$ 
equipped with the equality $=_I^{\Lambda}$. For simplicity, we write $\lambda_0(i) \in \lambda_0I$, instead of 
$i \in I$, when $I$ is equipped with the equality $=_I^{\Lambda}$. The operation $\lambda_0^* : I \sto I$ from 
$(I, =_I)$ to $(I, =_I^{\Lambda})$, defined by $i \mapsto i$, for every $i \in I$, is denoted by
$\lambda_0^* : I \sto \lambda_0 I$, and its definition is rewritten as $\lambda_0^*(i) := \lambda_0(i)$,
for every $i \in I$. 
\end{definition}

Clearly, $\lambda_0^*$ is a function. 
In the next proof the hypothesis of a set of sets is crucial.

\begin{proposition}\label{prp: Famtoset1}
Let $\Lambda := (\lambda_0, \lambda_1)$ be an $I$-set of sets, and let $Y$ be a set.
If $f \colon I \to Y$, there is a unique function $\lambda_0 f \colon \lambda_0 I \to Y$ such that the 
following diagram commutes
\begin{center}
\begin{tikzpicture}

\node (E) at (0,0) {$\lambda_0 I$};
\node [above=of E] (D) {$I$};
\node[right=of D] (F) {$Y.$};

\draw[dashed, ->] (E)--(F) node [midway,right] {$\ \lambda_0 f$};
\draw [->]  (D)--(E) node [midway,left] {$\lambda_0$};
\draw[->] (D)--(F) node [midway,above] {$f$};

\end{tikzpicture}
\end{center}
Conversely, if $f \colon I \sto Y$ and $f^* : \lambda_0 I  \to Y$ such that 
the corresponding diagram commutes, then $f$ is a function and $f^*$ is equal to the function from 
$\lambda_0 I$ to $Y$ generated by $f$.
\end{proposition}

\begin{proof}
The operation $\lambda_0 f$ from $\lambda_0 I$ to $Y$ defined by
$\lambda_0 f (\lambda_0 (i)) := f(i)$, for every $\lambda_0(i) \in \lambda_0I$, is a function,
since, for every $i, j \in I$, we have that $\lambda_0(i) =_{\D V_0} \lambda_0(j) \To i =_I j$, hence 
$f(i) =_Y f(j) :\TOT \lambda_0 f (\lambda_0(i)) =_Y  \lambda_0 f (\lambda_0(j))$. The  commutativity 
of the diagram follows from the reflexivity of $=_Y$.
If $g \colon \lambda_0 I \to Y$ makes the above diagram commutative, then for every $\lambda_0(i)$ we have
that $g(\lambda_0(i)) =_Y f(i) =: \lambda_0 f (\lambda_0(i))$, hence $g =_{\D F(\lambda_0 I, Y)}
\lambda_0 f$. For the converse, if $i, j \in I$, then by the transitivity of $=_Y$ we have that
$ i =_I j \To \lambda_0 (i) =_{\D V_0} \lambda_0 (j)$, hence  
$\To f^* (\lambda_0(i)) =_Y  f^* (\lambda_0(j))$, and $f(i) =_Y f(j)$.
The proof of the fact that $f^*$ is the function from $\lambda_0 I$ to $Y$ generated by $f$ is immediate.
\end{proof}

\begin{proposition}\label{prp: Famtoset2}
Let $\Lambda := (\lambda_0, \lambda_1) \in \Fam(I)$, and let $Y$ be a set.
If $f^* : \lambda_0 I \to Y$, there is a unique function $f : I \to Y$ such that the following diagram commutes
\begin{center}
\begin{tikzpicture}

\node (E) at (0,0) {$\lambda_0 I$};
\node [above=of E] (D) {$I$};
\node[right=of D] (F) {$Y.$};

\draw[->] (E)--(F) node [midway,right] {$\ f^*$};
\draw [->]  (D)--(E) node [midway,left] {$\lambda_0$};
\draw[dashed, ->] (D)--(F) node [midway,above] {$f$};

\end{tikzpicture}
\end{center}
If $\Lambda \in \Set(I)$, then $f^*$ is equal to the function from $\lambda_0 I$ 
to $Y$ generated by $f$.
\end{proposition}

\begin{proof}
Let $f^* \colon I \sto Y$, defined by $f(i) := f^*(\lambda_0(i))$, for every $i \in I$.
Since $i =_I i{'}  \To \lambda_0(i) =_{\D V_0} \lambda_0(i{'})
\To f^*(\lambda_0(i)) =_{Y} f^*(\lambda_0(i{'}))
\TOT f(i) =_Y f(i{'}),
$
$f$ is the required function.
If $\Lambda \in \Set(I)$,
by Proposition~\ref{prp: Famtoset1} $f^*$ is generated by $f$. The uniqueness of $f$ follows immediately.
\end{proof}

\begin{remark}\label{rem: Composition1}
Let $f \in \D F(I, Y)$ and $g \in \D F(Y, Z)$. If $\Lambda \in \Set(I)$, then
$\lambda_0 (g \circ f) := g \circ \lambda_0 f$
\begin{center}
\begin{tikzpicture}

\node (E) at (0,0) {$\lambda_0 I$};
\node [above=of E] (D) {$I$};
\node[right=of D] (F) {$Y$};
\node [right=of F] (G) {$Z$.};

\draw[->] (E)--(F) node [midway,above] {$\lambda_0 f \ \ $};
\draw[->] (F)--(G) node [midway,above] {$g$};
\draw [->]  (D)--(E) node [midway,left] {$\lambda_0$};
\draw[->] (D)--(F) node [midway,above] {$f$};
\draw[->] (E)--(G) node [midway,right] {$\ \ \lambda_0 (g \circ f)$};

\end{tikzpicture}
\end{center}
\end{remark}

\begin{proof}
If $i \in I$, then by Proposition~\ref{prp: Famtoset1} we have that
$\lambda_0 (g \circ f)(\lambda_0 (i)) := (g \circ f)(i)
:= g(f(i)) 
=: g[\lambda_0f (\lambda_0 (i))]
=: [g \circ \lambda_0 f](\lambda_0 (i)).$
\end{proof}

\begin{proposition}\label{prp: FamtoFam1}
Let $\Lambda := (\lambda_0, \lambda_1) \in \Set(I)$ and 
$M := (\mu_0, \mu_1) \in \Set(J)$. If $f : I \to J$, there is a unique 
function $f^* : \lambda_0 I \to \mu_0 J$ such that the following diagram commutes
\begin{center}
\begin{tikzpicture}

\node (E) at (0,0) {$\lambda_0 I$};
\node[right=of E] (F) {$\mu_0 J.$};
\node[above=of F] (A) {$J$};
\node [above=of E] (D) {$I$};

\draw[dashed, ->] (E)--(F) node [midway,below] {$f^*$};
\draw[->] (D)--(A) node [midway,above] {$f $};
\draw[->] (D)--(E) node [midway,left] {$\lambda_0$};
\draw[->] (A)--(F) node [midway,right] {$\mu_0$};

\end{tikzpicture}
\end{center}
Conversely, if $f \colon I \sto J$, and 
$f^* : \lambda_0 I \to \mu_0 J$ such that the corresponding to the above diagram commutes, then $f \in \D F(I, J)$
and $f^*$ is equal to the function from $\lambda_0 I$ to $\mu_0 J$ generated by $f$.
\end{proposition}

\begin{proof}
Let $f^* \colon \lambda_0 I \sto \mu_0 J$ be defined by 
$f^* (\lambda_0(i)) := \mu_0(f(i))$, for every $\lambda_0(i) \in \lambda_0 I$. We show that 
$f^* \in \D F(\lambda_0 I, \mu_0 J)$.  If 
$i, j \in I$, such that $\lambda_0(i) =_{\D V_0} \lambda_0(j)$, then $i =_I j$, hence 
$f(i) =_J f(j)$, and consequently $\mu_0(f(i)) =_{\D V_0} \mu_0(f(j))$ i.e.,
$f^* (\lambda_0(i)) =_{\D V_0}  f^* (\lambda_0(j)).$
The uniqueness of $f^*$ is trivial. For the converse, by the 
transitivity of $=_{\D V_0}$, and since $M \in \Set(J)$, we have that
$ i =_I j \To \lambda_0(i) =_{\D V_0} \lambda_0(j) 
\To f^* (\lambda_0(i)) =_{\D V_0}  f^* (\lambda_0(j))$,
hence $\mu_0(f(i)) =_{\D V_0} \mu_0(f(j))$, which implies 
$f(i) =_J f(j)$.
Clearly, $f^*$ is equal to the function from $\lambda_0 I$ to $\mu_0 J$ generated by $f$.
\end{proof}

\begin{proposition}\label{prp: FamtoFam2}
Let $\Lambda := (\lambda_0, \lambda_1) \in \Fam(I)$ and 
$M := (\mu_0, \mu_1) \in \Set(J)$.
If $f^* \colon \lambda_0 I \to \mu_0 J$, there is a unique function $f$ from $I$ to $J$, such that
the following diagram commutes, and $f^*$ is equal to the function from $\lambda_0 I$ to $\mu_0 J$ generated by $f$
\begin{center}
\begin{tikzpicture}

\node (E) at (0,0) {$\lambda_0 I$};
\node[right=of E] (F) {$\mu_0 J.$};
\node[above=of F] (A) {$J$};
\node [above=of E] (D) {$I$};

\draw[->] (E)--(F) node [midway,below] {$f^*$};
\draw[dashed, ->] (D)--(A) node [midway,above] {$f $};
\draw[->]  (D)--(E) node [midway,left] {$\lambda_0$};
\draw[->]  (A)--(F) node [midway,right] {$\mu_0$};

\end{tikzpicture}
\end{center}
\end{proposition}

\begin{proof}
If $i \in I$, then $f^*(\lambda_0(i)) := \mu_0(j)$, for some $j \in J$.
We define the routine 
$f(i) := j$ i.e., the output of $f^*$ determines the output of $f$. Since
$i =_I i{'} \To \lambda_0(i) =_{\D V_0} \lambda_0(i{'}) 
\To f^* (\lambda_0(i)) =_{\D V_0}  f^* (\lambda_0(i{'}))$
we get $\mu_0(j) =_{\D V_0} \mu_0(j{'}) \To j =_J j{'} \TOT f(i) =_J f(i{'})$, hence
$f$ is a function. The required commutativity of the diagram follows immediately. If $g : I \to J$
such that the above diagram commutes, then 
$\mu_0(g(i)) =_{\D V_0} f^*(\lambda_0(i)) := \mu_0(j) =: \mu_0(f(i))$, hence $g(i) =_J f(i)$. 
\end{proof}


\section{Direct families of sets}
\label{sec: dirfamsets}


\begin{definition}\label{def: dfamilyofsets}
Let $(I, \lt_I)$ be a directed set, and 
$D^{\lt}(I) := \big\{(i, j) \in I \times I \mid i \lt_I j \big\}$
the \textit{diagonal} of $\lt_I$.
A \textit{direct family of sets}\index{direct family of sets indexed by} $(I, \lt_I)$, 
or an $(I, \lt_I)$-\textit{family of sets}\index{$(I, \lt_I)$-family of sets}, 
is a pair $\Lambda^{\lt} := (\lambda_0, \lambda_1^{\lt})$, where
$\lambda_0 : I \sto \D V_0$, and $\lambda_1^{\lt}$\index{$\lambda_1^{\lt}$},
a \textit{modulus of transport maps} for $\lambda_0$,\index{modulus of transport maps} is defined by 
\[ \lambda_1^{\lt} : \bigcurlywedge_{(i, j) \in D^{\lt}(I)}\D F\big(\lambda_0(i), \lambda_0(j)\big), \ \ \ 
\lambda_1^{\lt}(i, j) := \lambda_{ij}^{\lt}, \ \ \ (i, j) \in D^{\lt}(I), \]
such that the \textit{transport maps}\index{transport map of a direct family of sets} $\lambda_{ij}^{\prec}$
of $\Lambda^{\lt}$ satisfy the following conditions:\\[1mm]
\normalfont (a) 
\itshape
For every $i \in I$, we have that $\lambda_{ii}^{\lt} := \id_{\lambda_0(i)}$.\\[1mm]
\normalfont (b) 
\itshape If $i \lt_I j$ and $j \lt_I k$, the following diagram commutes
\begin{center}
\begin{tikzpicture}

\node (E) at (0,0) {$\lambda_0(j)$};
\node[right=of E] (F) {$\lambda_0(k).$};
\node [above=of E] (D) {$\lambda_0(i)$};

\draw[->] (E)--(F) node [midway,below] {$\lambda_{jk}^{\lt}$};
\draw[->] (D)--(E) node [midway,left] {$\lambda_{ij}^{\lt}$};
\draw[->] (D)--(F) node [midway,right] {$\ \lambda_{ik}^{\lt}$};

\end{tikzpicture}
\end{center}
If $X \in \D V_0$, the \textit{constant} $(I, \lt_I)$-\textit{family}\index{constant family over $(I, \lt_I)$} 
$X$ is the pair $C^{\lt,X} := (\lambda_0^X, \lambda_1^{\lt,X})$, where $\lambda_0^X (i) := X$,  
and $\lambda_1^{\lt,X} (i, j) := \id_X$, for every $i \in I$ and $(i, j) \in  D^{\lt}(I)$.

\end{definition}

Since in general $\lt_I$ is not symmetric, the transport map $\lambda_{ij}^{\lt}$
does not necessarily have an inverse. Hence $\lambda_1^{\lt}$ is only a modulus of transport for 
$\lambda_0$, in the sense that determines the transport maps of $\Lambda^{\lt}$, and not necessarily 
a modulus of function-likeness for $\lambda_0$.

\begin{definition}\label{def: dirsumpi}
If $\Lambda^{\lt} := (\lambda_0, \lambda_1^{\lt})$ and $M^{\lt} 
:= (\mu_0, \mu_1^{\lt})$ are $(I, \lt_I)$-families of sets, a
\textit{direct family-map}\index{direct family-map} $\Phi$ from $\Lambda^{\lt}$ to $M^{\lt}$, denoted by 
$\Phi \colon \Lambda^{\lt} \To M^{\lt}$\index{$\Phi \colon \Lambda^{\lt} \To M^{\lt}$},
their set $\Map_{(I, \lt_I)}(\Lambda^{\lt}, M^{\lt})$\index{$\Map_{(I, \lt_I)}(\Lambda^{\lt}, M^{\lt})$},
and the totality $\Fam(I, \lt_I)$\index{$\Fam(I, \lt_I)$} of $(I, \lt_I)$-families are defined as
in Definition~\ref{def: map}.
The \textit{direct sum}\index{direct sum} 
$\sum_{i \in I}^{\lt}\lambda_0(i)$\index{$\sum_{i \in I}^{\lt}\lambda_0(i)$} over 
$\Lambda^{\lt}$ is the totality $\sum_{i \in I}\lambda_0(i)$ equipped with the equality
\[ (i, x) =_{\sum_{i \in I}^{\lt} \lambda_0(i)} (j, y) : \TOT \exists_{k \in I}\big(i \lt_I k \ 
\& \ j \lt_I k \ \& \ \lambda_{ik}^{\lt}(x) =_{\lambda_0(k)} \lambda_{jk}^{\lt}(y)\big). \] 
The totality $\prod_{i \in I}^{\lt}\lambda_0(i)$\index{$\prod_{i \in I}^{\lt}\lambda_0(i)$} 
of \textit{dependent functions over}\index{dependent functions over $\Lambda^{\lt}$} $\Lambda^{\lt}$ is defined by
\[ \Phi \in \prod_{i \in I}^{\lt}\lambda_0(i) :\TOT \Phi \in \D A(I, \lambda_0) 
\ \& \ \forall_{(i,j) \in D^{\lt}(I)}\big(\Phi_j =_{\lambda_0(j)} \lambda_{ij}^{\lt}(\Phi_i)\big), \]
and it is equipped with the equality of $\D A(I, \lambda_0)$.

\end{definition}

Clearly, the property $P(\Phi) :\TOT \forall_{(i,j) \in D^{\lt}(I)}\big(\Phi_j 
=_{\lambda_0(j)} \lambda_{ij}^{\lt}(\Phi_i)\big)$
is extensional on $\D A(I, \lambda_0)$,  
the equality on $\prod_{i \in I}^{\lt}\lambda_0(i)$ is an equivalence relation.   
$\prod_{i \in I}^{\lt}\lambda_0(i)$ is considered to be a set.

\begin{proposition}\label{prp: direct1}
The relation $(i, x) =_{\sum_{i \in I}^{\lt} \lambda_0(i)} (j, y)$ 
is an equivalence relation.
\end{proposition}

\begin{proof}
If $i \in I$, and since $i \lt_I i$, there is $k \in I$
such that $i \lt_I k$, and by the reflexivity of the equality on $ \lambda_0(k)$ we get 
$\lambda_{ik}^{\lt}(x) =_{\lambda_0(k)} \lambda_{ik}^{\lt}(x)$. The symmetry of 
$=_{\sum_{i \in I}^{\lt} \lambda_0(i)}$ follows from the symmetry of the equalities
$=_{\lambda_0(k)}$. To prove transitivity, we suppose that 
\[ (i, x) =_{\sum_{i \in I}^{\lt} \lambda_0(i)} (j, y) : \TOT \exists_{k \in I}\big(i \lt_I k \ \& \ j 
\lt_I k \ \& \ \lambda_{ik}^{\lt}(x) =_{\lambda_0(k)} \lambda_{jk}^{\lt}(y)\big), \]
\[ (j, y) =_{\sum_{i \in I}^{\lt} \lambda_0(i)} (j{'}, z) : \TOT \exists_{k{'} \in I}\big(j \lt_I k{'} 
\ \& \ j {'} \lt_I k{'} \ \& \ \lambda_{jk{'}}^{\lt}(y) =_{\lambda_0(k{'})} 
\lambda_{j{'}k{'}}^{\lt}(z)\big), \]
and we show that
\[ (i, x) =_{\sum_{i \in I}^{\lt} \lambda_0(i)} (j{'}, z) : \TOT \exists_{k{''} \in I}\big(i \lt_I k{''} 
\ \& \ j{'} \lt_I k{''} \ \& \ \lambda_{ik{''}}^{\lt}(x) =_{\lambda_0(k{''})} 
\lambda_{j{'}k{''}}^{\lt}(z)\big). \]
By the definition of a directed set there is $k{''} \in I$ such that $k \lt_I k{''}$ and $k{'} \lt_I k{''}$
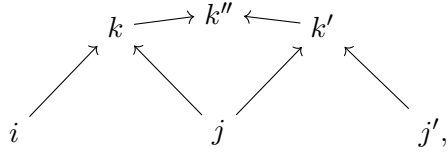
\begin{center}
\begin{tikzpicture}

\node (E) at (0,0) {$i$};
\node [right=of E] (F) {};
\node [above=of F] (A) {$k$};
\node [right=of F] (B) {$j$};
\node [right=of B] (C) {};
\node [right=of C] (D) {$j{'}$,};
\node [above=of C] (F) {$k{'}$};
\node [above=of B] (H) {$k{''}$};

\draw[->] (E)--(A) node [midway,below] {};
\draw[->] (B)--(A) node [midway,above] {};
\draw[->] (B)--(F) node [midway,above] {};
\draw[->] (D)--(F) node [midway,above] {};
\draw[->] (A)--(H) node [midway,above] {};
\draw[->] (F)--(H) node [midway,above] {};

\end{tikzpicture}
\end{center}
hence by transitivity $i \lt_I k{''}$ and $j{'} \lt_I k{''}$. Moreover,
\begin{align*}
\lambda_{ik{''}}^{\lt}(x) & \stackrel{i \lt_I k \lt_I k{''}} = \lambda_{kk{''}}^{\lt}
\big(\lambda_{ik}^{\lt}(x)\big) \\
& \ \ \  = \ \ \  \lambda_{kk{''}}^{\lt} \big(\lambda_{jk}^{\lt}(y)\big) \\
& \stackrel{j \lt_I k \lt_I k{''}} = \lambda_{jk{''}}^{\lt} (y) \\
& \stackrel{j \lt_I k{'} \lt_I k{''}} = \lambda_{k{'}k{''}}^{\lt} \big(\lambda_{jk{'}}^{\lt}(y)\big) \\
& \ \ \ = \ \ \ \lambda_{k{'}k{''}}^{\lt} \big(\lambda_{j{'}k{'}}^{\lt}(z)\big) \\
& \stackrel{j{'} \lt_I k{'} \lt_I k{''}} = \lambda_{j{'}k{''}}^{\lt}(z).\qedhere
\end{align*}
\end{proof}

Notice that the projection operation from $\sum_{i \in I}^{\lt}\lambda_0(i)$ to $I$ 
is not a function.

\begin{proposition}\label{prp: preorderfamilymap1}
If $(I, \lt_I)$ is a directed set, $\Lambda^{\lt} := (\lambda_0, \lambda_1^{\lt})$, $M^{\lt} 
:= (\mu_0, \mu_1^{\lt})$ are $(I, \lt_I)$-families of sets, and $\Psi^{\lt} : \Lambda^{\lt} \To M^{\lt}$, 
the following hold.\\[1mm]
\normalfont (i) 
\itshape
For every $i \in I$ the operation
$e_i^{\Lambda^{\lt}} : \lambda_0(i) \sto \sum_{i \in I}^{\lt}\lambda_0(i)$, defined by 
$x \mapsto (i, x)$, for every $x \in \lambda_0(i)$, 
is a function from $\lambda_0(i)$ to $\sum_{i \in I}^{\lt}\lambda_0(i)$.\\[1mm]
\normalfont (ii) 
\itshape The operation
$\Sigma^{\lt} \Psi : \sum_{i \in I}^{\lt}\lambda_0(i) \sto \sum_{i \in I}^{\lt}\mu_0(i)$,
defined by
$\big(\Sigma^{\lt} \Psi\big) (i, x) := (i, \Psi_i (x))$, for every $(i, x) \in \sum_{i \in I}^{\lt}\lambda_0(i)$,
is a function from $\sum_{i \in I}^{\lt}\lambda_0(i)$ to $\sum_{i \in I}^{\lt}\mu_0(i)$ such that, for every
$i \in I$, the following left diagram commutes
\begin{center}
\begin{tikzpicture}

\node (E) at (0,0) {$\sum_{i \in I}^{\lt}\lambda_0(i)$};
\node[right=of E] (F) {$\sum_{i \in I}^{\lt}\mu_0(i)$};
\node[above=of F] (A) {$\mu_0(i)$};
\node [above=of E] (D) {$\lambda_0(i)$};

\node [right=of F] (G) {$\prod_{i \in I}^{\lt}\lambda_0(i)$};
\node[right=of G] (K) {$\prod_{i \in I}^{\lt}\mu_0(i)$.};
\node[above=of K] (L) {$\mu_0(i)$};
\node [above=of G] (M) {$\lambda_0(i)$};

\draw[->] (E)--(F) node [midway,below]{$\Sigma^{\lt} \Psi$};
\draw[->] (D)--(A) node [midway,above] {$\Psi_{i}$};
\draw[->] (D)--(E) node [midway,left] {$e_i^{\Lambda^{\lt}}$};
\draw[->] (A)--(F) node [midway,right] {$e_i^{M^{\lt}}$};
\draw[->] (G)--(K) node [midway,below]{$\Pi^{\lt} \Psi$};
\draw[->] (M)--(L) node [midway,above] {$\Psi_{i}$};
\draw[->] (G)--(M) node [midway,left] {$\pi_i^{\Lambda^{\lt}}$};
\draw[->] (K)--(L) node [midway,right] {$\pi_i^{M^{\lt}}$};

\end{tikzpicture}
\end{center}
\normalfont (iii) 
\itshape If $\Psi_i$ is an embedding, for every $i \in I$, then $\Sigma^{\lt} \Psi$ is an embedding.\\[1mm]
\normalfont (iv) 
\itshape For every $i \in I$ the operation
$\pi_i^{\Lambda^{\lt}} : \prod_{i \in I}^{\lt}\lambda_0(i) \sto \lambda_0(i)$, defined by
$\Theta \mapsto \Theta_i$, for every $\Theta \in \prod_{i \in I}^{\lt}\lambda_0(i)$, 
is a function from $\prod_{i \in I}^{\lt}\lambda_0(i)$ to $\lambda_0(i)$.\\[1mm]
\normalfont (v) 
\itshape The operation
$\Pi^{\lt} \Psi : \prod_{i \in I}^{\lt}\lambda_0(i) \sto \prod_{i \in I}^{\lt}\mu_0(i)$, defined by
$[\Pi^{\lt} \Psi (\Theta)]_i := \Psi_i (\Theta_i)$, for every $i \in I$ and $\Theta \in \prod_{i \in I}^{\lt}\lambda_0(i)$,
is a function from $\prod_{i \in I}^{\lt}\lambda_0(i)$ to $\prod_{i \in I}^{\lt}\mu_0(i)$, such that, for every $i \in I$, 
the above right diagram commutes.\\[1mm]
%
%
%
%
%
\normalfont (vi) 
\itshape
If $\Psi_i$ is an embedding, for every $i \in I$, then $\Pi^{\lt} \Psi$ is an 
embedding.
\end{proposition}

\begin{proof}
(i) If $x, y \in \lambda_0(i)$ such that $x =_{\lambda_0 (i)} y$, then, since $\lt$ is reflexive, if we
take $k := i$, we get $\lambda_{ii}^{\lt}(x) := \id_{\lambda_0(i)}(x) := x =_{\lambda_0(i)} y :=
\id_{\lambda_0(i)}(y) := \lambda_{ii}^{\lt}(y)$, hence 
$(i, x) =_{\sum_{i \in I}^{\lt}\lambda_0(i)} (i, y)$.\\
(ii) If $(i, x) =_{\sum_{i \in I}^{\lt}\lambda_0(i)} (j, y)$, there is $k \in I$ such that $i \lt_I k$,
$j \lt_I k$ and $\lambda_{ik}^{\lt}(x) =_{\lambda_0(k)} \lambda_{jk}^{\lt}(y)$. We show the following equality: 
\begin{align*}
\big(\Sigma^{\lt} \Psi\big)(i, x) =_{\sum_{i \in I}^{\lt}\mu_0(i)} 
\big(\Sigma^{\lt} \Psi\big)(j, y) & : \TOT (i, \Psi_i (x)) 
=_{\mathsmaller{\sum_{i \in I}^{\lt}\mu_0(i)}} (j, \Psi_j (y))\\
& : \TOT \exists_{k{'} \in I}\big(i, j \lt_I k{'} \ \& \ 
\mu_{ik{'}}^{\lt}(\Psi_i (x)) =_{\mu_0(k{'})} \mu_{jk{'}}^{\lt}(\Psi_j(y))\big).
\end{align*} 
If we take $k{'} := k$, by the commutativity of the following diagrams, and  since $\Psi_k$ is a function,
\begin{center}
\begin{tikzpicture}

\node (E) at (0,0) {$\mu_0(i)$};
\node[right=of E] (F) {$\mu_0(k)$};
\node[above=of F] (A) {$\lambda_0(k)$};
\node [above=of E] (D) {$\lambda_0(i)$};
\node [right=of F] (G) {$\mu_0 (j)$};
\node [above=of G] (H) {$\lambda_0(j)$};
\node [right=of G] (K) {$\mu_0 (k)$};
\node [above=of K] (L) {$\lambda_0(k)$};

\draw[->] (E)--(F) node [midway,below]{$\mu_{ik}^{\lt}$};
\draw[->] (D)--(A) node [midway,above] {$\lambda_{ik}^{\lt}$};
\draw[->] (D)--(E) node [midway,left] {$\Psi_i$};
\draw[->] (A)--(F) node [midway,right] {$\Psi_k$};
\draw[->] (G)--(K) node [midway,below] {$\mu_{jk}^{\lt}$};
\draw[->] (H)--(L) node [midway,above] {$\lambda_{jk}^{\lt}$};
\draw[->] (H)--(G) node [midway,left] {$\Psi_j$};
\draw[->] (L)--(K) node [midway,right] {$\Psi_k$};

\end{tikzpicture}
\end{center}
\[ \mu_{ik}^{\lt}\big(\Psi_i(x)\big) =_{\mu_0 (k)} \ \Psi_k \big(\lambda_{ik}^{\lt}(x)\big)
=_{\mu_0 (k)} \ \Psi_k \big(\lambda_{jk}^{\lt}(y)\big)
=_{\mu_0 (k)} \ \mu_{jk}^{\lt}\big(\Psi_j (y)\big). \]
(iii) If we suppose $\big(\Sigma^{\lt} \Psi\big)(i, x) =_{\mathsmaller{\sum_{i \in I}^{\lt}\mu_0(i)}} 
\big(\Sigma^{\lt} \Psi\big)(j, y)$ i.e.,
$\mu_{ik}^{\lt}(\Psi_i (x)) =_{\mu_0(k)} \mu_{jk}^{\lt}(\Psi_j(y))\big)$,
for some $k \in I$ with $i, j \lt_I k$, by the proof of case (ii) we get  
$\Psi_k \big(\lambda_{ik}^{\lt}(x)\big) =_{\mu_0 (k)}  
\Psi_k \big(\lambda_{jk}^{\lt}(y)\big)$, and since
$\Psi_k$ is an embedding, we get $\lambda_{ik}^{\lt}(x) =_{\lambda_0 (k)} 
\lambda_{jk}^{\lt}(y)$ i.e., $(i, x) =_{\sum_{i \in I}^{\lt}\lambda_0(i)} (j, y)$.\\
(iv)-(vi) Their proof is omitted, since a proof of their contravariant version (see Note~\ref{not: preorder})
is given in the proof of Theorem~\ref{thm: inverselimitmap}. 
\end{proof}

Since the transport functions $\lambda_{ik}^{\lt}$ are not in general embeddings, 
we cannot show in general that $e_i^{\Lambda^{\lt}}$ is an embedding, as it is the case for the map 
$e_i^{\Lambda}$ in Proposition~\ref{prp: map1}(i). The study of direct families of sets can be extended following
the study of set-indexed families of sets.

\section{Set-relevant families of sets}
\label{sec: genfamsofsets}

In general, we may want to have more than one transport maps from $\lambda_0(i)$ to $\lambda_0(j)$, if $i =_I j$. 
In this case, to each $(i, j) \in D(I)$ we associate a set of transport maps.

\begin{definition}\label{def: genfamofsets}
If $I$ is a set, a \textit{set-relevant family of sets indexed by} $I$,\index{set-relevant family of sets}
is a triplet $\Lambda^* := \big(\lambda_0, \varepsilon_0^{\lambda}, \lambda_2)$,
where $\lambda_0 : I \sto \D V_0$, $\varepsilon_0^{\lambda} : D(I) \sto \D V_0$, and\index{$\Lambda^*$}\index{$\lambda_2$} 
\[ \lambda_2 : \bigcurlywedge_{(i, j) \in D(I)}\bigcurlywedge_{p \in \varepsilon_0^{\lambda}(i, j)}
\D F\big(\lambda_0(i), \lambda_0(j)\big), 
\ \ \ \lambda_2\big((i, j), p\big) := \lambda_{ij}^p, \ \ \ (i, j) \in D(I), \ \ \ p \in 
\varepsilon_0^{\lambda}(i, j), \]
such that the following conditions hold:\\[1mm]
\normalfont (i) 
\itshape For every $i \in I$ there is $p \in \varepsilon_0^{\lambda}(i, i)$ such that $\lambda_{ii}^p
=_{\D F(\lambda_0(i), \lambda_0(i))} \id_{\lambda_0(i)}$.\\[1mm]
\normalfont (ii) 
\itshape For every $(i, j) \in D(I)$ and every $p \in \varepsilon_0^{\lambda}(i, j)$ there is some
$q \in \varepsilon_0^{\lambda}(j, i)$ such that 
such that the following left diagram 
commutes
\begin{center}
\begin{tikzpicture}

\node (E) at (0,0) {$\lambda_0(j)$};
\node[right=of E] (F) {$\lambda_0(i)$};
\node [above=of E] (D) {$\lambda_0(i)$};
\node [right=of F] (G) {$\lambda_0(j)$};
\node [above=of G] (H) {$\lambda_0(i)$};
\node [right=of G] (K) {$\lambda_0(k)$.};

\draw[->] (E)--(F) node [midway,below] {$\lambda_{ji}^q$};
\draw[->] (D)--(E) node [midway,left] {$\lambda_{ij}^p$};
\draw[->] (D)--(F) node [midway,right] {$\ \id_{\lambda_0(i)}$};
\draw[->] (H)--(G) node [midway,left] {$\lambda_{ij}^p$};
\draw[->] (H)--(K) node [midway,right] {$\ \lambda_{ik}^r$};
\draw[->] (G)--(K) node [midway,below] {$\lambda_{jk}^q$};

\end{tikzpicture}
\end{center}
\normalfont (iii)
\itshape If $i =_I j =_I k$, then for every $p \in \varepsilon_0^{\lambda}(i, j)$ 
and every $q \in \varepsilon_0^{\lambda}(j, k)$ there is
$r \in \varepsilon_0^{\lambda}(i, k)$ such that the above right diagram commutes.\\[1mm]
We call $\Lambda^*$ \textit{function-like}\index{function-like set-relevant family of sets}, if 
$ \forall_{(i, j) \in D(I)}\forall_{p, p{'} \in \varepsilon_0^{\lambda}(i, j)}\big( p =_{\varepsilon_0^{\lambda} 
(i, j)} p{'} \To \lambda_{ij}^p =_{\D F(\lambda_0(i), \lambda_0(j))} \lambda_{ij}^{p{'}}\big)$.
\end{definition}

It is immediate to show that if $\Lambda := (\lambda_0, \lambda_1) \in \Fam(I)$, then $\Lambda$ 
generates a set-relevant family over $I$, where 
$\varepsilon_0^{\lambda}(i,j) := \D 1$, and $\lambda_2\big((i,j), p)\big) := \lambda_{ij}$, for every $(i,j) \in D(I)$.

\begin{definition}\label{def: setrelmap}
Let $\Lambda^* := (\lambda_0, \varepsilon_0^{\lambda}, \lambda_2)$ and $M := (\mu_0, 
\varepsilon_0^{\mu}, \mu_2)$ be set-relevamt families of sets over $I$.
A \textit{covariant set-relevant family-map} from $\Lambda^*$ to $M^*$\index{covariant set-relevant family-map},
in symbols $\Psi \colon \Lambda^* \To M^*$\index{$\Psi \colon \Lambda^* \To M^*$}, is a dependent operation
$\Psi \colon \bigcurlywedge_{i \in I}\D F \big(\lambda_0(i), \mu_0(i)\big)$ 
such that for every $(i, j) \in D(I)$ and for every $p \in \varepsilon_0^{\lambda}(i, j)$ there 
is $q \in \varepsilon_0^{\mu}(i,j)$ such that the following diagram commutes
\begin{center}
\begin{tikzpicture}

\node (E) at (0,0) {$\mu_0(i)$};
\node[right=of E] (F) {$\mu_0(j)$.};
\node[above=of F] (A) {$\lambda_0(j)$};
\node [above=of E] (D) {$\lambda_0(i)$};

\draw[->] (E)--(F) node [midway,below]{$\mu_{ij}^q$};
\draw[->] (D)--(A) node [midway,above] {$\lambda_{ij}^p$};
\draw[->] (D)--(E) node [midway,left] {$\Psi_i$};
\draw[->] (A)--(F) node [midway,right] {$\Psi_j$};

\end{tikzpicture}
\end{center}
A contravariant set-relevant family-map\index{contravariant set-relevant family-map} is defined by the property:
for every 
 $q \in \varepsilon_0^{\mu}(i,j)$, there is $p \in \varepsilon_0^{\lambda}(i, j)$ such that the above diagram commutes. 
Let $\Map_I(\Lambda^*, M^*)$\index{$\Map_I(\Lambda^*, M^*)$} be the totality of covariant set-relevant 
family-maps\index{totality of covariant set-relevant family-maps} 
from $\Lambda^*$ to $M^*$, which is equipped with the pointwise equality.
If $\Xi : M^* \To N^*$, the \textit{composition  set-relevant family-map}\index{composition of
covariant set-relevant family-maps} 
$\Xi \circ \Psi \colon \Lambda^* \To N^*$ is defined, for every $i \in I$, by 
$(\Xi \circ \Psi)_i := \Xi_i \circ \Psi_i$
\begin{center}
\begin{tikzpicture}

\node (E) at (0,0) {$\lambda_0(i)$};
\node[right=of E] (F) {$\lambda_0(j)$};
\node[below=of F] (A) {$\mu_0(j)$};
\node[below=of E] (B) {$\mu_0(i)$};
\node[below=of B] (K) {$\nu_0(i)$};
\node[below=of A] (L) {$\nu_0(j)$.};

\draw[->] (E)--(F) node [midway,above] {$\lambda_{ij}^p$};
\draw[->] (F)--(A) node [midway,left] {$\Psi_j$};
\draw[->] (B)--(A) node [midway,below] {$\mu_{ij}^q$};
\draw[->] (E) to node [midway,right] {$\Psi_i$} (B);
\draw[->] (B) to node [midway,right] {$\Xi_i$} (K);
\draw[->] (A)--(L) node [midway,left] {$\Xi_j$};
\draw[->] (K)--(L) node [midway,below] {$\nu_{ij}^r$};
\draw[->,bend right] (E) to node [midway,left] {$(\Xi \circ \Psi)_i  \ $} (K);
\draw[->,bend left] (F) to node [midway,right] {$\  (\Xi \circ \Psi)_j$} (L);

\end{tikzpicture}
\end{center} 
The composition of contravariant set-relevant family-maps is defined similarly
The \textit{identity set-relevant family-map}\index{identity set-relevant family-map}\index{$\Id_{\Lambda^*}$} 
 is defined by $\Id_{\Lambda^*}(i) := \id_{\lambda_0(i)}$, for every $i \in I$.
Let $\Fam^*(I)$\index{$\Fam^*(I)$} be the totality of set-relevant $I$-families\index{totality of
set-relevant $I$-families},
equipped with the obvious canonical equality.

\end{definition}

$\Id_{\Lambda^*}$ is both a covariant and a contravariant set-relevant family-map from $\Lambda^*$ to itself.

\begin{definition}\label{def: setrelsumset}
Let $\Lambda^* := \big(\lambda_0, \varepsilon_0^{\lambda}, \lambda_2\big) \in \Fam^*(I)$.
The exterior union $\sum_{i \in I}^*\lambda_0(i)$\index{$\sum_{i \in I}^*\lambda_0(i)$} of $\Lambda^*$ is
the  totality $\sum_{i \in I} \lambda_0 (i)$, equipped with the following equality
\[ (i, x) =_{\sum_{i \in I}^* \lambda_0 (i)} (j, y) : \TOT i =_I j \ \& \ \exists_{p \in 
\varepsilon_0^{\lambda} (i, j)}\big(
\lambda_{ij}^p (x) =_{\lambda_0 (j)} y\big). \]
The totality $\prod_{i \in I}^*\lambda_0(i)$\index{$\prod_{i \in I}^*\lambda_0(i)$} of dependent 
functions over $\Lambda^*$ is defined by 
\[ \Theta \in \prod_{i \in I}^*\lambda_0(i) :\TOT \Theta \in \D A(I, \lambda_0) \
\& \ \forall_{(i,j) \in D(I)}\forall_{p \in \varepsilon_0^{\lambda}(i,j)}\big(\Theta_j =_{\lambda_0(j)} \lambda_{ij}^p(\Theta_i)\big), \]
and it is equipped with the pointwise equality.
\end{definition}

A motivation for the definitions of $\sum^*_{i \in I}\lambda_0(i)$ and $\prod_{i \in I}^*\lambda_0(i)$ 
is provided in Note~\ref{not: onprsigma}.

\begin{remark}\label{rem: genexteriorunion1}
The equalities on $\sum_{i \in I}^* \lambda_0 (i)$ and $\prod_{i \in I}^*\lambda_0(i)$ satisfy 
the conditions of an equivalence relation.
\end{remark}

\begin{proof}
Let $(i, x), (j, y)$ and $(k, z) \in \sum_{i \in I}\lambda_0 (i)$. By definition there is
$p \in \varepsilon^{\lambda}_0 (i, i)$ such that $\lambda_{ii}^p = \id_{\lambda_0 (i)}$,
hence $(i, x) =_{\sum_{i \in I}^*  \lambda_0 (i)} (i, x)$.
If $(i, x) =_{\sum_{i \in I}^* \lambda_0 (i)} (j, y)$, then $j =_I i$ and there is $q \in 
\varepsilon^{\lambda}_0 (ji)$ such that 
$\lambda_{ji}^q (y) = \lambda_{ji}^q(\lambda_{ij}^p (x)) = \id_{\lambda_0 (i)}(x) := x$,
hence $(j, y) =_{\sum_{i \in I}^g*\lambda_0 (i)} (i, x)$. If $(i, x) =_{\sum_{i \in I}^* \lambda_0 (i)} (j, y)$
and $(j, y) =_{\sum_{i \in I}^* \lambda_0 (i)} (k, z)$, then from the hypotheses $i =_I j$ and $j =_I k$, we 
get $i =_I k$. From the hypotheses $\exists_{p \in \varepsilon_0^{\lambda} (i, j)}\big(
\lambda_{ij}^p (x) =_{\lambda_0 (j)} y\big)$ and $\exists_{q \in \varepsilon_0^{\lambda} (j, k)}\big(
\lambda_{jk}^q (y) =_{\lambda_0 (k)} z\big)$, let $r \in \varepsilon_0^{\lambda} (i, k)$ 
such that $\lambda_{ik}^r = \lambda_{jk}^q \circ \lambda_{ij}^p$. Hence 
$\lambda_{ik}^r (x) =_{\lambda_0(k)} \lambda_{jk}^q(\lambda_{ij}^p(x)) =_{\lambda_0 (k)}
\lambda_{jk}^q(y) = z$. The proof for the equality on $\prod_{i \in I}^*\lambda_0(i)$ is trivial.
\end{proof}

\begin{proposition}\label{prp: setrelmap1}
Let $\Lambda := (\lambda_0, \varepsilon_0^{\lambda}, \lambda_2)$, $M := (\mu_0, \varepsilon_0^{\mu},
\mu_2) \in \Fam^*(I)$, and $\Psi : \Lambda^* \To M^*$.\\[1mm]
\normalfont (i) 
\itshape For every $i \in I$ the operation
$e_i^{\Lambda^*} : \lambda_0(i) \sto \sum_{i \in I}^*\lambda_0(i)$, defined by 
$ e_i^{\Lambda^*}(x) := (i, x)$, for every $x \in \lambda_0(i)$,
is a function.\\[1mm]
\normalfont (ii) 
\itshape If $\Psi$ is covariant, the operation
$\Sigma^* \Psi : \sum_{i \in I}^*\lambda_0(i) \sto \sum_{i \in I}^*\mu_0(i)$, defined by
$ \Sigma^* \Psi (i, x) := (i, \Psi_i (x))$, for every $(i, x) \in \sum_{i \in I}^*\lambda_0(i)$, 
is a function,
such that for every
$i \in I$ the following left diagram commutes
\begin{center}
\begin{tikzpicture}

\node (E) at (0,0) {$\sum_{i \in I}^*\lambda_0(i)$};
\node[right=of E] (F) {$\sum_{i \in I}^*\mu_0(i)$.};
\node[above=of F] (A) {$\mu_0(i)$};
\node [above=of E] (D) {$\lambda_0(i)$};

\node[right=of F] (K) {$\prod_{i \in I}^*\lambda_0(i)$};
\node[right=of K] (L) {$\prod_{i \in I}^*\mu_0(i)$.};
\node[above=of L] (M) {$\mu_0(i)$};
\node [above=of K] (N) {$\lambda_0(i)$};

\draw[->] (E)--(F) node [midway,below]{$\Sigma^* \Psi$};
\draw[->] (D)--(A) node [midway,above] {$\Psi_{i}$};
\draw[->] (D)--(E) node [midway,left] {$e_i^{\Lambda^*}$};
\draw[->] (A)--(F) node [midway,right] {$e_i^{M^*}$};

\draw[->] (K)--(L) node [midway,below]{$\Pi^* \Psi$};
\draw[->] (N)--(M) node [midway,above] {$\Psi_{i}$};
\draw[->] (K)--(N) node [midway,left] {$\pi_i^{\Lambda^*}$};
\draw[->] (L)--(M) node [midway,right] {$\pi_i^{M^*}$};

\end{tikzpicture}
\end{center}
\normalfont (iii) 
\itshape If $i \in I$, the operation
$\pi_i^{\Lambda^*} : \prod_{i \in I}^*\lambda_0(i) \sto \lambda_0(i)$, defined by
$\Theta \mapsto \Theta_i,$ is a function.\\[1mm]
\normalfont (iv) 
\itshape If $\Psi$ is contravariant, the operation
$\Pi^* \Psi : \prod_{i \in I}^*\lambda_0(i) \sto \prod_{i \in I}^*\mu_0(i)$, defined by
$[\Pi^* \Psi (\Theta)]_i := \Psi_i (\Theta_i)$, for every $i \in I$, 
is a function, 
such that for every $i \in I$ the above right diagram commutes.\\[1mm]
\normalfont (v) 
\itshape If $\Psi_i$ is an embedding, for every $i \in I$, then $\Pi^* \Psi$ is an 
embedding.

\end{proposition}

\begin{proof}
We proceed as in the proofs of Propositions~\ref{prp: map1} and~\ref{prp: map2}.
\end{proof}

The definitions of operations on $I$-families of sets and their family-maps extend to operations
on set-relevant $I$-families and their family-maps. An important example of a set-relevant family 
of sets is that of a family of sets over a set with a proof-relevant equality 
(see Definition~\ref{def: prfamofsets}). For reasons that are going to be clear in the study of 
these families of sets, the first definitional clause of a set-relevant family of sets over $I$ involves 
the equality $=_{\D F(\lambda_0(i), \lambda_0(i))}$ instead of the definitional one. Next follows the
definition of the direct version of a set-relevant family of sets, the importance of which is explained 
in Note~\ref{not: setrelavantspectra}.

\begin{definition}\label{def: dirgenfamofsets}
If $(I, \lt_I)$ is a directed set, a \textit{set-relevant direct family of sets indexed by}
$I$,\index{set-relevant direct family of sets} is a quadruple $\Lambda^{*, \lt} := \big(\lambda_0, 
\varepsilon_0^{\lambda^{\lt}},  \Theta, \lambda_2^{\lt})$,
where $\lambda_0 : I \sto \D V_0$, $\varepsilon_0^{\lambda^{\lt}} : D^{\lt}(I) \sto \D V_0$, $\Theta 
\in \bigcurlywedge_{(i,j) \in D^{\lt}(I)}\varepsilon_0^{\lambda^{\lt}}(i,j)$ is a modulus of 
inhabitedness for $\varepsilon_0^{\lambda^{\lt}}$, 
and\index{$\Lambda^{*,\lt}$}\index{$\lambda_2^{\lt}$} 
\[ \lambda_2^{\lt} : \bigcurlywedge_{(i, j) \in D^{\lt}(I)}\bigcurlywedge_{p \in \varepsilon_0^{\lambda^{\lt}}(i, j)}
\D F\big(\lambda_0(i), \lambda_0(j)\big), 
\ \  \lambda_2^{\lt}\big((i, j), p\big) := \lambda_{ij}^{p, \lt}, \ \ (i, j) \in D(I), \ \  p \in
\varepsilon_0^{\lambda^{\lt}}(i, j), \]
such that the following conditions hold:\\[1mm]
\normalfont (i) 
\itshape For every $i \in I$ there is $p \in \varepsilon_0^{\lambda^{\lt}}(i, i)$ such that $\lambda_{ii}^{p, \lt}
=_{\D F(\lambda_0(i), \lambda_0(i))} \id_{\lambda_0(i)}$.\\[1mm]
\normalfont (ii) 
\itshape If $i \lt_I j \lt_I k$, then for every $p \in \varepsilon_0^{\lambda^{\lt}}(i, j)$ and every
$q \in \varepsilon_0^{\lambda^{\lt}}(j, k)$ there is
$r \in \varepsilon_0^{\lambda^{\lt}}(i, k)$ such that $\lambda_{jk}^{q, \lt} \circ \lambda_{ij}^{p, \lt} 
= \lambda_{ik}^{r, \lt}$.\\
\normalfont (iii)
\itshape For every $(i, j) \in D^{\lt}(I)$ and every $p, p{'} \in \varepsilon_0^{\lambda^{\lt}}(i, j)$ 
and every $x \in \lambda_0(i)$ there is $k \in I$ such that $j \lt_I k$ and there is
$q \in \varepsilon_0^{\lambda^{\lt}} (j, k)$ such that 
\[ \lambda_{jk}^{q, \lt}\big(\lambda_{ij}^{p, \lt}(x)\big) =_{\lambda_0 (k)}
\lambda_{jk}^{q, \lt}\big(\lambda_{ij}^{p{'}, \lt}(x)\big). \]

\end{definition}

The modulus of inhabitedness $\Theta$ for $\varepsilon_0^{\lambda^{\lt}}$ and the last condition in the
previous definition guarantee that the equality on the corresponding $\sum$-set of a set-relevant direct family of sets 
satisfies the conditions of an equivalence relation.

\begin{definition}\label{def: dirgenexteriorunion}
Let $\Lambda^{*, \lt} := \big(\lambda_0, \varepsilon_0^{\lambda^{\lt}}, \Theta, \lambda_2^{\lt}\big)$ a 
set-relevant family of sets over a directed set $(I, \lt)$. Its exterior union $\sum_{i \in I}^{*, \lt}
\lambda_0 (i)$ is the totality  
$\sum_{i \in I} \lambda_0 (i)$ equipped with the equality
\[ (i, x) =_{\mathsmaller{\sum_{i \in I}^{*, \lt}\lambda_0 (i)}} (j, y) : \TOT 
\exists_{k \in I}\bigg(i, j \lt_I k \ \& \ 
\exists_{p \in \varepsilon_0^{\lambda^{\lt}} (i, k)}\exists_{q \in \varepsilon_0^{\lambda^{\lt}} (j, k)}\big(
\lambda_{ik}^{p, \lt} (x) =_{\mathsmaller{\lambda_0 (k)}} \lambda_{jk}^{q, \lt} (y)\big)\bigg). \]
The set $\prod_{i \in I}^{*, \lt}\lambda_0(i)$ of dependent functions over $\Lambda^{*, \lt}$ is defined by 
\[ \Phi \in \prod_{i \in I}^{*,\lt}\lambda_0(i) :\TOT \Phi \in \D A(I, \lambda_0) \
\& \ \forall_{(i,j) \in D^{\lt}(I)}\forall_{p \in \varepsilon_0^{\lambda, \lt}(i,j)}\big(\Phi_j =_{\lambda_0(j)} \lambda_{ij}^{p,\lt}(\Phi_i)\big), \]
and it is equipped with the pointwise equality.
\end{definition}

\begin{proposition}\label{prp: genexteriorunion1}
The equality on $\sum_{i \in I}^{*, \lt} \lambda_0 (i)$ satisfies the conditions of an equivalence relation.
\end{proposition}

\begin{proof}
To show that $(i, x) =_{\mathsmaller{\sum_{i \in I}^{*, \lt}\lambda_0 (i)}} (i, x)$ we use the first 
definitional clause of a set-relevant directed family of sets.
The proof of the equality $(j, y) =_{\mathsmaller{\sum_{i \in I}^{*, \lt}\lambda_0 (i)}} (i, x)$ from the equality 
$(i, x) =_{\mathsmaller{\sum_{i \in I}^{*, \lt}\lambda_0 (i)}} (j, y)$ is trivial.
For transitivity we suppose that 
$(i, x) =_{\sum_{i \in I}^{*, \lt} \lambda_0(i)} (j, y)$ and 
$(j, y) =_{\mathsmaller{\sum_{i \in I}^{*, \lt} \lambda_0(i)}} (j{'}, z)$ i.e.,
\[ \exists_{k{'} \in I}\bigg(j \lt_I k{'} 
\ \& \ j {'} \lt_I k{'} \ \& \ \exists_{p{'} \in \varepsilon_0^{\lambda^{\lt}} (j, k{'})}\exists_{q{'} \in 
\varepsilon_0^{\lambda^{\lt}}(j{'}, k{'})}\big(\lambda_{jk{'}}^{p{'}, \lt}(y) =_{\mathsmaller{\lambda_0(k{'})}} 
\lambda_{j{'}k{'}}^{q{'}, \lt}(z)\big)\bigg). \]
There is $k{''} \in I$ such that $k \lt_I k{''}$ and $k{'} \lt_I k{''}$. 
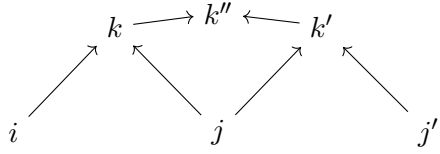
\begin{center}
\begin{tikzpicture}

\node (E) at (0,0) {$i$};
\node [right=of E] (F) {};
\node [above=of F] (A) {$k$};
\node [right=of F] (B) {$j$};
\node [right=of B] (C) {};
\node [right=of C] (D) {$j{'}$};
\node [above=of C] (F) {$k{'}$};
\node [above=of B] (H) {$k{''}$};

\draw[->] (E)--(A) node [midway,below] {};
\draw[->] (B)--(A) node [midway,above] {};
\draw[->] (B)--(F) node [midway,above] {};
\draw[->] (D)--(F) node [midway,above] {};
\draw[->] (A)--(H) node [midway,above] {};
\draw[->] (F)--(H) node [midway,above] {};
\%draw[->] (H)--(P) node [midway,above] {};

\end{tikzpicture}
\end{center}
Moreover, there are $r \in \varepsilon_0^{\lambda^{\lt}}(i, k{''})$ and 
$s \in \varepsilon_0^{\lambda^{\lt}}(j, k{''})$ such that
\[ \lambda_{ik{''}}^{r, \lt}(x)  \stackrel{\mathsmaller{i \lt_I k \lt_I k{''}}} = 
\lambda_{kk{''}}^{\Theta_{(k,k{''})},\lt} \big(\lambda_{ik}^{p, \lt}(x)\big)  = \ \ \  
\lambda_{kk{''}}^{\Theta_{(k,k{''})}, \lt} \big(\lambda_{jk}^{q, \lt}(y)\big) \ \ \ 
\stackrel{\mathsmaller{j \lt_I k \lt_I k{''}}} = \lambda_{jk{''}}^{s, \lt} (y). \]
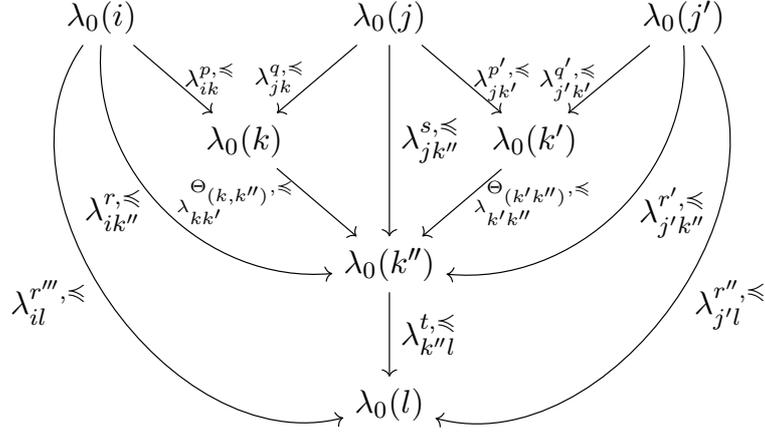
\begin{center}
\resizebox{11cm}{!}{%
\begin{tikzpicture}

\node (E) at (0,0) {$\lambda_0(i)$};
\node[right=of E] (B) {};
\node[right=of B] (F) {$\lambda_0(j)$};
\node[below=of F] (C) {};
\node[below=of C] (A) {$\lambda_0(k{''})$};
\node[right=of F] (D) {};
\node[right=of D] (I) {$\lambda_0(j{'})$};
\node[below=of B] (J) {$\lambda_0(k)$};
\node[below=of D] (K) {$\lambda_0(k{'})$};
\node[below=of A] (L) {$\lambda_0(l)$};

\draw[->] (E)--(J) node [midway,right] {$\mathsmaller{\lambda_{ik}^{p, \lt}}$};
\draw[->,bend right=50] (E) to node [midway,left] {$\lambda^{r, \lt}_{ik{''}} \ $} (A);
\draw[->,bend right=70] (E) to node [midway,left] {$\lambda_{il}^{r{'''}, \lt} \ $} (L);
\draw[->] (F)--(J) node [midway,left] {$\mathsmaller{\lambda_{jk}^{q, \lt}}$};
\draw[->] (F)--(K) node [midway,right] {$\mathsmaller{\lambda_{jk{'}}^{p{'}, \lt}}$};
\draw[->] (F)--(A) node [midway,right] {$\lambda_{jk{''}}^{s, \lt}$};
\draw[->] (I)--(K) node [midway,left] {$\mathsmaller{\lambda_{j{'}k{'}}^{q{'}, \lt}}$};
\draw[->,bend left=50] (I) to node [midway,right] {$ \ \lambda^{r{'}, \lt}_{j{'}k{''}}$} (A);
\draw[->,bend left=70] (I) to node [midway,right] {$ \ \lambda_{j{'}l}^{r{''}, \lt} $} (L);
\draw[->] (A)--(L) node [midway,right] {$\lambda_{k{''}l}^{t, \lt}$};
\draw[->] (J)--(A) node [midway,left] {$\mathsmaller{\mathsmaller{\lambda_{kk{'}}^{\Theta_{(k,k{''})}, \lt}}} \ $};
\draw[->] (K)--(A) node [midway,right] {$\mathsmaller{\mathsmaller{\lambda_{k{'}k{''}}^{\Theta_{(k{'}k{''})}, \lt}}}$};

\end{tikzpicture}
}
\end{center}
If we apply condition (iii) of Definition~\ref{def: dirgenexteriorunion} to $j \lt_I k{''}$, $y \in \lambda_0 (j)$ and
the transport maps $\lambda_{jk{''}}^{s, \lt} $ and $\lambda_{k{'}k{''}}^{\Theta_{(k{'},k{''})}, \lt} \circ 
\lambda_{jk{'}}^{p{'}, \lt}$ from $\lambda_0 (j)$ to $\lambda_0 (k{''})$, then there is 
$l \in I$, such that $k{''} \lt_I l$, and some $t \in \varepsilon_0^{\lambda^{\lt}} (k{''}, l)$ such that 
\begin{align*}
\lambda_{k{''}l}^{t, \lt}\big(\lambda_{jk{''}}^{s, \lt}(y)\big) & =
\lambda_{k{''}l}^{t, \lt}\bigg(\lambda_{k{'}k{''}}^{\Theta_{(k{'},k{''})}, 
\lt}\big(\lambda_{jk{'}}^{p{'}, \lt}(y)\big)\bigg)\\
& = \lambda_{k{''}l}^{t, \lt}\bigg(\lambda_{k{'}k{''}}^{\Theta_{(k{'},k{''})}, 
\lt}\big(\lambda_{j{'}k{'}}^{q{'}, \lt}(z)\big)\bigg)\\
& = \lambda_{k{''}l}^{t, \lt}\big(\lambda_{j{'}k{''}}^{r{'}, \lt}(z)\bigg)\\
& = \lambda_{j{'}l}^{r{''}, \lt}(z),
\end{align*}
for some $r{'} \in \varepsilon_0^{\lambda^{\lt}} (j{'}, k{''})$ and some $r{''} \in 
\varepsilon_0^{\lambda^{\lt}} (j{'}, l)$. From 
$\lambda_{ik{''}}^{r, \lt} (x) =_{\mathsmaller{\lambda_0 (k{''})}} \lambda_{jk{''}}^{s, \lt} (y)$ we get
\[ 
\lambda_{k{''}l}^{t, \lt}\big(\lambda_{jk{''}}^{s, \lt}(y)\big)  =
\lambda_{k{''}l}^{t, \lt}\big(\lambda_{ik{''}}^{r, \lt} (x)\big)
= \lambda_{il}^{r{'''}, \lt} (x),
\]
for some $r{'''} \in \varepsilon_0^{\lambda^{\lt}} (i, l)$. Hence,
$\lambda_{il}^{r{''i}, \lt} (x) =_{\mathsmaller{\lambda_0 (l)}} 
\lambda_{j{'}l}^{r{''}, \lt}(z)$
i.e., $(i, x) =_{\mathsmaller{\sum_{i \in I}^{*, \lt} \lambda_0(i)}} (j{'}, z)$.
\end{proof}

\section{Families of families of sets, an impredicative interlude}
\label{sec: famoffam}

We define the notion of a family of families of sets $(\Lambda^i)_{i \in I}$, where each family of sets 
$\Lambda^i$ is indexed by some set $\mu_0(i)$, and $i \in I$. As expected, the index-sets are given by 
some family $M \in \Fam(I)$, and $(\Lambda^i)_{i \in I}$ must must be a function-like object i.e., 
if $i =_I j$, the family $\Lambda^i$ of sets over the index-set $\mu_0(i)$ is ``equal'' to the 
family $\Lambda^j$ of sets over the index-set $\mu_0(j)$. This equality can be expressed through the 
notion of a family map from $\Lambda^i$ to $\Lambda^j$ over $\mu_{ij} $ 
(see Definition~\ref{def: hmap}). As in the case of the definition of a family of sets we provide the a 
priori given transport maps of $(\Lambda^i)_{i \in I}$ with certain properties that guarantee the existence 
of these family-maps.
As $\Fam(I)$ is an impredicative set, to define a family of families of sets, we need to introduce, in 
complete analogy to the introduction of $\D V_0$, the class\index{universe of sets and impredicative sets} $\D V_0^{\im}$\index{$\D V_0^{\im}$} \textit{of sets and
impredicative sets}. All notions of assignment routines defined in Chapter~\ref{chapter: BST} are defined 
in a similar way when the class $\D V_0^{\im}$ is used instead of $\D V_0$. We add the superscript
$^{\im}$ to a symbol in order to denote the version of the corresponding notion that requires the use of $\D V_0^{\im}$.

\begin{definition}\label{def: famoffamofsets}
Let $M := (\mu_0, \mu_1) \in \Fam(I)$ and, if $(i, j) \in D(I)$ let the 
set\index{$T_{ij}(M)$}\index{$(\Lambda^i)_{i \in I}$}
\[ T_{ij}(M) := \{ (m, n) \in \mu_0(i) \times \mu_0(j) \mid \mu_{ij}(m) =_{\mu_0(j)} n \}. \]
A \textit{family of families of sets over} $I$ and $M$\index{family
of families of sets}, or an\index{$(I,M)$-family of families of sets} $(I, M)$-\textit{family of families of sets},
is a pair
$(\Lambda^i)_{i \in I} := \big(\Lambda^{0,M}, \Lambda^{1, M}\big)$,
where 
\[ \Lambda^{0,M} \colon \bigcurlywedge^{\im}\Fam(\mu_0(i)), \ \ \ \ \Lambda^{0,M}_i := \big(\lambda_0^i, 
\lambda_1^1\big); \ \ \ \ i \in I, \]
\[ \Lambda^{1,M} \colon \bigcurlywedge_{(i,j) \in D(I)}\bigcurlywedge_{(m,n) \in T_{ij}(M)}\D F\big(\lambda_0^i(m),
\lambda_0^j(n)\big), \]
\[ \big(\Lambda^{1,M}_{(i,j)}\big)_{(m,n)} := \lambda_{mn}^{ij} \colon \lambda_0^i(m) \to \lambda_0^j(n); \ \ \ \ 
(i,j) \in D(I), \ (m,n) \in T_{ij}(M), \]
such that the \textit{transport maps} $\lambda^{ij}_{mn}$
\index{transport maps of a family of families of sets} of $(\Lambda^i)_{i \in I}$, satisfy the following 
conditions:\\[1mm]
\normalfont (i)
\itshape For every $i \in I$ and $(m,m{'}) \in T_{ij}(M)$, we have that $\lambda^{ii}_{mm{'}} := \lambda^i_{mm{'}}$.\\[1mm]
\normalfont (ii)
\itshape If $i =_I j =_I k$, for every $(m,n) \in T_{ij}(M)$ and $(n,l) \in T_{jk}(M)$, the following 
diagram commutes 
\begin{center}
\begin{tikzpicture}

\node (E) at (0,0) {$\lambda_0^j (n)$};
\node[right=of E] (F) {$\lambda_0^k (l)$.};
\node [above=of E] (D) {$\lambda_0^i (m)$};

\draw[->] (E)--(F) node [midway,below] {$\lambda^{jk}_{nl}$};
\draw[->] (D)--(E) node [midway,left] {$\lambda^{ij}_{mn}$};
\draw[->] (D)--(F) node [midway,right] {$\ \lambda^{ik}_{ml}$};

\end{tikzpicture}
\end{center}
Let $\Fam(I,M)$\index{$\Fam(I, M)$} be the totality of $(I,M)$-families of families of sets. 
\end{definition}

For condition (i) above, we have that $\mu_{ii}(m) = m{'} :\TOT m =_{\mu_(i)} m{'}$ and for condition (ii), 
from the hypotheses $(m,n) \in T_{ij}(M) :\TOT \mu_{ij}(m) =n$ and $(n,l) \in T_{jk}(M) :\TOT \mu_{jk}(n) = l$
we get $(m,l) \in T_{ik}(M) :\TOT \mu_{ik}(m) = l$, as $\mu_{ik}(m) = \mu_{jk}(\mu_{ij}(m)) = l$. The main
intuition behind this defintion is that if $i =_I j$, then $\mu_0(i) =_{\D V_0} \mu_0(j)$, hence, if 
$m \in \mu_0(i)$ and $n \in \mu_0(j)$, there is a transport map $\lambda^{ij}_{mn}$ from $\lambda_0^i(m)$
to $\lambda_0^j(n)$. It is easy to see that if $\Lambda \in \Fam(J)$, then by taking $I := \D 1$ and $M$ 
the constant family $J$ over $\D 1$, then $\Lambda$ can be viewed as an $(\D 1, M)$-family of families of sets.

\begin{lemma}\label{lem: famoffam1}
If $(\Lambda^i)_{i \in I} \in \Fam(I,M)$, for its transport maps $\lambda_{mn}^{ij}$ the following hold:\\[1mm]
\normalfont (i)
\itshape $\lambda^{ii}_{mm} := \id_{\lambda_0^i (m)}$.\\[1mm]
\normalfont (ii)
\itshape $\lambda^{ij}_{mn} \circ \lambda^{ji}_{nm} = \id_{\lambda_0^j (n)}$.\\[1mm]
\normalfont (iii)
\itshape If $\mu_{ij}(m) = n$, then $\lambda^{ij}_{m\mu_{ij}(m)} = \lambda^j_{n \mu_{ij}(m)} \circ \lambda^{ij}_{mn}$
\begin{center}
\begin{tikzpicture}

\node (E) at (0,0) {$\lambda^j_0(\mu_{ij}(m))  $};
\node[right=of E] (H) {};
\node[right=of H] (F) {$ \lambda^j_0(n)$.};
\node[above=of F] (A) {$\lambda^i_0(m)$};
\node [above=of E] (D) {$\lambda^i_0(m)$};

\draw[->] (F)--(E) node [midway,below] {$\lambda^j_{n \mu_{ij}(m)}$};
\draw[->] (D)--(A) node [midway,above] {$\lambda^i_{mm}$};
\draw[->] (D)--(E) node [midway,left] {$\lambda^{ij}_{m \mu_{ij}(m)}$};
\draw[->] (A)--(F) node [midway,right] {$\lambda^{ij}_{mn}$};

\end{tikzpicture}
\end{center}
\normalfont (iv)
\itshape If $m =_{\mathsmaller{\mu_0 (i)}} m{'}$, then $\lambda^j_{\mu_{ij}(m) \mu_{ij}(m{'})} \circ
\lambda^{ij}_{m \mu_{ij}(m)} =   \lambda^{ij}_{m{'}\mu_{ij}(m{'})} \circ \lambda^{i}_{mm{'}}$ 
\begin{center}
\begin{tikzpicture}

\node (E) at (0,0) {$\lambda^j_0(\mu_{ij}(m))  $};
\node[right=of E] (H) {};
\node[right=of H] (F) {$\lambda^j_0(\mu_{ij}(m{'}))$.};
\node[above=of F] (A) {$\lambda^i_0(m{'})$};
\node [above=of E] (D) {$\lambda^i_0(m)$};

\draw[->] (E)--(F) node [midway,below] {$\lambda^j_{\mu_{ij}(m) \mu_{ij}(m{'})}$};
\draw[->] (D)--(A) node [midway,above] {$\lambda^i_{mm{'}}$};
\draw[->] (D)--(E) node [midway,left] {$\lambda^{ij}_{m \mu_{ij}(m)}$};
\draw[->] (A)--(F) node [midway,right] {$\lambda^{ij}_{m{'}\mu_{ij}(m{'})}$};

\end{tikzpicture}
\end{center}
\normalfont (v)
\itshape If $i =_I j =_I k$, then 
$\lambda^k_{\mu_{ik}(m)\mu_{jk}(\mu_{ij}(m))} \circ \lambda^{ik}_{m \mu_{ik}(m)} = 
\lambda^{jk}_{\mu_{ij}(m)\mu_{jk}(\mu_{ij}(m))} \circ \lambda^{ij}_{m \mu_{ij}(m)}$
\begin{center}
\begin{tikzpicture}

\node (E) at (0,0) {$\lambda^j_0(\mu_{ij}(m)) $};
\node[right=of E] (H) {};
\node[right=of H] (L) {};
\node[right=of L] (F) {$\lambda^k_0(\mu_{jk}(\mu_{ij}(m))$.};
\node[above=of F] (A) {$\lambda^k_0(\mu_{ik}(m))$};
\node [above=of E] (D) {$\lambda^i_0(m)$};

\draw[->] (E)--(F) node [midway,below] {$\lambda^{jk}_{\mu_{ij}(m) \mu_{jk}(\mu_{ij}(m)}$};
\draw[->] (D)--(A) node [midway,above] {$\lambda^{ik}_{m \mu_{ik}(m)}$};
\draw[->] (D)--(E) node [midway,left] {$\lambda^{ij}_{m \mu_{ij}(m)}$};
\draw[->] (A)--(F) node [midway,right] {$\lambda^{ik}_{\mu_{ik}(m)\mu_{jk}(\mu_{ij}(m))}$};

\end{tikzpicture}
\end{center}
\end{lemma}

\begin{proof}
(i) By Definition~\ref{def:  famoffamofsets} $\lambda^{ii}_{mm} := \lambda^i_{mm} := \id_{\lambda_0^i (m)}$.\\
(ii) By the composition-rule and case (i) we get $\lambda^{ij}_{mn} \circ \lambda^{ji}_{nm} = \lambda^{jj}_{nn} = 
\id_{\lambda_0^j (n)}$.\\
(iii) By the composition-rule we get
$\lambda^{ij}_{m\mu_{ij}(m)} = \lambda^{jj}_{n \mu_{ij}(m)} \circ \lambda^{ij}_{mn} := \lambda^j_{n \mu_{ij}(m)}
\circ \lambda^{ij}_{mn}$.\\[1mm]
(iv) By Definition~\ref{def: famoffamofsets} we have that
\begin{align*}
\lambda^j_{\mu_{ij}(m) \mu_{ij}(m{'})} \circ \lambda^{ij}_{m \mu_{ij}(m)} & := \lambda^{jj}_{\mu_{ij}(m)
\mu_{ij}(m{'})} \circ \lambda^{ij}_{m \mu_{ij}(m)}\\
& = \lambda^{ij}_{m \mu_{ij}(m{'})}\\
& = \lambda^{ij}_{m{'} \mu_{ij}(m{'})} \circ \lambda^{ii}_{m m{'}}\\
& := \lambda^{ij}_{m{'}\mu_{ij}(m{'})} \circ \lambda^{i}_{mm{'}}.
\end{align*}
(v) If $l := \mu_{jk}(\mu_{ij}(m))$, then by case (iii) 
$\lambda^{ik}_{m \mu_{ik}(m)} = \lambda^k_{\mu_{jk}(\mu_{ij}(m)) \mu_{ik}(m)} \circ \lambda^{ik}_{m 
\mu_{jk}(\mu_{ij}(m))}$, hence
\begin{align*}
\lambda^k_{\mu_{ik}(m)\mu_{jk}(\mu_{ij}(m))} \circ \lambda^{ik}_{m \mu_{ik}(m)} & = 
\lambda^k_{\mu_{ik}(m)\mu_{jk}(\mu_{ij}(m))} \circ \big[\lambda^k_{\mu_{jk}(\mu_{ij}(m)) \mu_{ik}(m)} \circ \lambda^{ik}_{m \mu_{jk}(\mu_{ij}(m))} \big]\\
& := \big[\lambda^k_{\mu_{ik}(m)\mu_{jk}(\mu_{ij}(m))} \circ \lambda^k_{\mu_{jk}(\mu_{ij}(m)) \mu_{ik}(m)}\big] 
\circ \lambda^{ik}_{m \mu_{jk}(\mu_{ij}(m))}\\ 
& = \lambda^k_{\mu_{jk}(\mu_{ij}(m))\mu_{jk}(\mu_{ij}(m))} \circ \lambda^{ik}_{m \mu_{jk}(\mu_{ij}(m))}\\
& :=  \lambda^{ik}_{m \mu_{jk}(\mu_{ij}(m))}\\
& = \lambda^{jk}_{\mu_{ij}(m)\mu_{jk}(\mu_{ij}(m))} \circ \lambda^{ij}_{m \mu_{ij}(m)}.\qedhere
\end{align*}
\end{proof}

\begin{definition}\label{def: transportfammaps}
If $(\Lambda^i)_{i \in I} \in \Fam(I, M)$, and based on Definition~\ref{def: hmap}, its transport 
family-maps $\Phi_{ij}^{\Lambda}$\index{transport family-maps}\index{$\Phi_{ij}^{\Lambda}$} are the 
family-maps $\Phi_{ij}^{\Lambda} \colon \Lambda_i^{0,M} \stackrel{\mu_{ij}} \Longrightarrow \Lambda_j^{0,M}$, 
defined by the rule 
\[ \big[\Phi_{ij}^{\Lambda}\big]_m := \lambda_{m \mu_{ij}(m)}^{ij}; \ \ \ \ m \in \mu_0(i), (i,j) \in D(I).\]
\end{definition}

The fact that $\lambda_{m \mu_{ij}(m)}^{ij} \colon \Lambda_i^{0,M} \stackrel{\mu_{ij}} \Longrightarrow \Lambda_j^{0,M}$
is shown by the commutativity of the diagram in case (iv) of Lemma~\ref{lem: famoffam1}. In analogy to 
the transport maps $\lambda_{ij}$ of an $I$-family of sets $\Lambda$, the transport family-maps
$\Phi_{ij}^{\Lambda}$ witness the equality between
the $\mu_0(i)$-family of sets $\Lambda^{0,M}_i$ and the $\mu_0(j)$-family of sets $\Lambda_j^{0,M}$.

\begin{definition}\label{def: exteriorunionfamoffam}
If $(\Lambda^i)_{i \in I} \in \Fam(I, M)$,
its \textit{exterior union} $\sum_{i \in I}\sum_{m \in \mu_0 (i)}\lambda_0 ^i(m)$\index{exterior union of a
family of families of sets}\index{$\sum_{i \in I}\sum_{m \in \mu_0 (i)}\lambda_0 ^i(m)$}
is defined by
\[ w \in \sum_{i \in I}\sum_{m \in \mu_0 (i)}\lambda_0 ^i(m) : \TOT 
\exists_{i \in I}\exists_{m \in \mu_0 (i)}\exists_{x \in \lambda_0^i (m)}\big(w := (i, m, x)\big),\]
\[ (i, m, x) =_{\mathsmaller{\sum_{i \in I}\sum_{m \in \mu_0 (i)}\lambda_0 ^i(m)}} (j, n, y) 
:\TOT i =_I j \ \& \ \mu_{ij}(m) =_{\mathsmaller{\mu_0 (j)}} n \ \& \ \lambda^{ij}_{mn}(x)
=_{\mathsmaller{\lambda^j_0(n)}} y.\]
\end{definition}

\begin{remark}\label{rem: exteriorunion1}
The equality on $\sum_{i \in I}\sum_{m \in \mu_0 (i)}\lambda_0 ^i(m)$ satisfies the conditions of 
an equivalence relation.
\end{remark}

\begin{proof}
To show $(i, m, x) = (i, m, x) :\TOT i =_I i \ \& \ \mu_{ii}(m) = m \ \& \ \lambda^{ii}_{mm}(x) = x$,
we use Lemma~\ref{lem: famoffam1}(i). If $(i, m, x) = (j, n, y) : \TOT i =j \ \& \ \mu_{ij}(m) = n \ 
\& \ \lambda^{ij}_{mn}(x) = y,$ then $j = i$ and $\mu_{ji}(n) = m$, and, using  Lemma~\ref{lem: famoffam1}(ii),
$\lambda^{ji}_{nm}(y) = \lambda^{ji}_{nm}\big(\lambda^{ij}_{mn}(x)\big) = \lambda^{ii}_{mm}(x) = x$ i.e., 
$(j, n, y) = (i, m, x)$. If $(i, m, x) = (j, n, y)$ and $(j, n, xy = (k, l, z) : \TOT j = k \ \& \ \mu_{jk}(n) =
l \ \& \ \lambda^{jk}_{nl}(y) = z,$ then $i = k$ and $\mu_{ik}(m) = \mu_{jk}(\mu_{ij}(m)) = \mu_{jk}(n) = l$, 
and $\lambda^{ik}_{ml}(x) = \lambda^{jk}_{nl}\big(\lambda^{ij}_{mn}(x)\big) = \lambda^{jk}_{nl}(y) = z$ i.e.,
$(i, m, x) = (k, l, z)$.
\end{proof}

\begin{proposition}\label{prp: sumoffamoffam}
If $(\Lambda^i)_{i \in I} \in \Fam(I, M)$, then $\Sigma := (\sigma_0, \sigma_1) \in \Fam(I)$, where 
\[ \sigma_0 (i) := \sum_{m \in \mu_0 (i)}\lambda^i_0(m); \ \ \ \ i \in I,\]
\[ \sigma_1 (i, j) := \sigma_{ij} \colon \bigg(\sum_{m \in \mu_0 (i)}\lambda^i_0(m)\bigg) \to \sum_{n 
\in \mu_0 (j)}\lambda^j_0(n); \ \ \ \ (i,j) \in D(I), \]
\[ \sigma_{ij}(m, x) := \big(\mu_{ij}(m), \lambda^{ij}_{m\mu_{ij}(m)}(x)\big); \ \ \ \ m \in \mu_0(i),
\ x \in \lambda_0^i(m).\]
\end{proposition}

\begin{proof}
First we show that the operation $\sigma_{ij}$
 is a function. We suppose that $(m, x) =_{\mathsmaller{\sum_{m \in \mu_0 (i)}\lambda^i_0(m)}} (m{'}, x{'})
 : \TOT m =_{\mathsmaller{\mu_0 (i)}} m{'} \ \& \ \lambda^i_{mm{'}}(x) =_{\mathsmaller{\lambda^i_0(m{'})}} x$, 
 and we show that 
\[ \big(\mu_{ij}(m), \lambda^{ij}_{m\mu_{ij}(m)}(x)\big)  =_{\mathsmaller{\sum_{n \in \mu_0 (j)}\lambda^j_0(n)}}  
\big(\mu_{ij}(m{'}), \lambda^{ij}_{m\mu_{ij}(m{'})}(x{'})\big) \TOT \]
\[ \mu_{ij}(m) =_{\mathsmaller{\mu_0 (j)}} \mu_{ij}(m{'}) \ \& \ 
\lambda^j_{\mu_{ij}(m)\mu_{ij}(m{'})}\big(\lambda^{ij}_{m\mu_{ij}(m)}(x)\big) 
=_{\mathsmaller{\lambda^j_0(\mu_{ij}(m{'}))}} \lambda^{ij}_{m{'}\mu_{ij}(m{'})}(x{'}).\]
The first conjunct follows from $m =_{\mathsmaller{\mu_0 (i)}} m{'}$, and the second is
Lemma~\ref{lem: famoffam1}(iv). Since
\[ \sigma_{ii}(m, x) := \big(\mu_{ii}(m), \lambda^{ii}_{m \mu_{ii}(m)}(x)\big) := (m, \lambda^{ii}_{mm}(x)) 
:= (m, \id_{\lambda_0 ^i(m)}(x)) := (m, x), \]
we get $\sigma_{ii} := \id_{\sum_{m \in \mu_0 (i)}\lambda^i_0(m)}$. For the commutativity of the diagram
\begin{center}
\begin{tikzpicture}

\node (E) at (0,0) {$\sum_{n \in \mu_0 (j)}\lambda^j_0(n)$};
\node[right=of E] (F) {$\sum_{l \in \mu_0 (k)}\lambda^k_0(l)$};
\node [above=of E] (D) {$\sum_{m \in \mu_0 (i)}\lambda^i_0(m)$};

\draw[->] (E)--(F) node [midway,below] {$\sigma_{jk}$};
\draw[->] (D)--(E) node [midway,left] {$\sigma_{ij}$};
\draw[->] (D)--(F) node [midway,right] {$\ \sigma_{ik}$};

\end{tikzpicture}
\end{center}
we have that by definition
$\sigma_{ik}(m, x) := \big(\mu_{ik}(m), \lambda^{ik}_{m \mu_{ik}(m)}(x)\big),$ and
\begin{align*}
\sigma_{jk}\big(\sigma_{ij}(m, x)\big) & := \sigma_{jk}\big(\mu_{ij}(m), \lambda^{ij}_{m \mu_{ij}(m)}(x)\big)\\
& := \bigg(\mu_{jk}(\mu_{ij}(m)), \ 
\lambda^{jk}_{\mu_{ij}(m)\mu_{jk}(\mu_{ij}(m))}\big(\lambda^{ij}_{m \mu_{ij}(m)}(x)\big)\bigg).
\end{align*}
Hence, 
$\sigma_{ik}(m, x) =_{\mathsmaller{\sum_{l \in \mu_0 (k)}\lambda^k_0(l)}} \sigma_{jk}\big(\sigma_{ij}(m, x)\big) 
: \TOT \mu_{ik}(m)  =_{\mathsmaller{\mu_0 (k)}} \mu_{jk}(\mu_{ij}(m))$ and
\[ \lambda^k_{\mu_{ik}(m)\mu_{jk}(\mu_{ij}(m))}\big(\lambda^{ik}_{m \mu_{ik}(m)}(x)\big) = 
\lambda^{jk}_{\mu_{ij}(m)\mu_{jk}(\mu_{ij}(m))}\big(\lambda^{ij}_{m \mu_{ij}(m)}(x)\big). \]
The first conjunct is immediate to show, and the second is exactly Lemma~\ref{lem: famoffam1}(v).
\end{proof}

Clearly, for the exterior union of $(\Lambda^i)_{i \in I}$ we have that
\[ \sum_{i \in I}\sum_{m \in \mu_0 (i)}\lambda^i_0(m) =_{\D V_0} \sum_{i \in I}\sigma_0 (i) := 
\sum_{i \in I}\bigg(\sum_{m \in \mu_0 (i)}\lambda^i_0(m)\bigg).\]

If $(\Lambda^i)_{i \in I}$ and $(M^i)_{i \in I}$ are $(I, M)$-families of families, a map 
from $(\Lambda^i)_{i \in I}$ to $(M^i)_{i \in I}$ is an appropriate dependent function $(\Psi^i)_{i \in I}$ 
such that $\Psi^i$ is a family map from $\Lambda^i$ to $M^i$, for every $i \in I$. Before giving this definition 
we show a fact of independent interest.

\begin{proposition}\label{prp: Map1}
Let $(K^i)_{i \in I} := (K^{0,M}, K^{1,M}), (\Lambda^i)_{i \in I} := (\Lambda^{0,M}, \Lambda^{1,M}) \in \Fam(I, M)$, 
and let $i =_I j$.\\[1mm]
\normalfont (i)
\itshape The operation $\gamma_{ij} : \Map\big(K_i^{0,M}, \Lambda_i^{0,M}\big) \sto \Map\big(K_j^{0,M}, 
\Lambda_j^{0,M}\big)$, defined by the rule
$\Psi^i \mapsto \big[\gamma_{ij}(\Psi^i)\big]^j$, is a function, where, for every $n \in \mu_0(j)$, the 
map $\big[\gamma_{ij}(\Psi^i)\big]_n^j \colon \kappa_0^j(n) \to \lambda_0^j(n)$ is defined by
$\big[\gamma_{ij}(\Psi^i)\big]_n^j := \lambda^{ij}_{\mu_{ji}(n)n} \circ \Psi^i_{\mu_{ji}(n)} \circ
\kappa^{ji}_{n \mu_{ji}(n)}$
\begin{center}
\begin{tikzpicture}

\node (E) at (0,0) {$\kappa^i_0(\mu_{ji}(n))  $};
\node[right=of E] (L) {};
\node[right=of L] (F) {$\lambda^i_0(\mu_{ji}(n))$.};
\node[above=of F] (A) {$\lambda^j_0(n)$};
\node [above=of E] (D) {$\kappa^j_0(n)$};

\draw[->] (E)--(F) node [midway,below] {$\Psi^i_{\mu_{ji}(n)}$};
\draw[->] (D)--(A) node [midway,above] {$\big[\gamma_{ij}(\Psi^i)\big]_n^j$};
\draw[->] (D)--(E) node [midway,left] {$\kappa^{ji}_{n \mu_{ji}(n)}$};
\draw[->] (F)--(A) node [midway,right] {$\lambda^{ij}_{\mu_{ji}(n)n}$};

\end{tikzpicture}
\end{center}
\normalfont (ii)
\itshape The pair $\Gamma := (\gamma_0, \gamma_1) \in \Fam(I)$, where $\gamma_0 (i) := \Map\big(K_i^{0,M},
\Lambda_i^{0,M}\big)$, for every $i \in I$, and $\gamma_1 (i, j) := \gamma_{ij}$, for every 
$(i, j) \in D(I)$.
\end{proposition}

\begin{proof}
(i) First we show that $\gamma_{ij}$ is well-defined i.e., $\big[\gamma_{ij}(\Psi^i)\big]^j \in \Map(K^j, 
\Lambda^j)$. If $n, n{'} \in \mu_0 (j)$, we show that the following diagram commutes
\begin{center}
\begin{tikzpicture}

\node (E) at (0,0) {$\lambda_0^j(n)$};
\node[right=of E] (L) {};
\node[right=of L] (F) {$\lambda_0^j(n{'})$.};
\node[above=of F] (A) {$\kappa_0^j(n{'})$};
\node [above=of E] (D) {$\kappa_0^j(n)$};

\draw[->] (E)--(F) node [midway,below]{$\lambda^j_{nn{'}}$};
\draw[->] (D)--(A) node [midway,above] {$\kappa^j_{nn{'}}$};
\draw[->] (D)--(E) node [midway,left] {$\big[\gamma_{ij}(\Psi^i)\big]_n^j$};
\draw[->] (A)--(F) node [midway,right] {$\big[\gamma_{ij}(\Psi^i)\big]_{n{'}}^j$};

\end{tikzpicture}
\end{center}
By definition we have that
$ \big[\gamma_{ij}(\Psi^i)\big]_{n{'}}^j \circ \kappa^j_{nn{'}} := \big[\lambda^{ij}_{\mu_{ji}(n{'})n{'}} \circ
\Psi^i_{\mu_{ji}(n{'})} \circ \kappa^{ji}_{n{'} \mu_{ji}(n{'})}\big] \circ \kappa^j_{nn{'}}$.
Since $\Psi^i$ is in $\Map(K_i^{0,M}, \Lambda_i^{0,M})$, by the commutativity of the following diagram
\begin{center}
\begin{tikzpicture}

\node (E) at (0,0) {$\lambda_0^i(\mu_{ji}(n))$};
\node[right=of E] (L) {};
\node[right=of L] (F) {$\lambda_0^i(\mu_{ji}(n{'}))$,};
\node[above=of F] (A) {$\kappa_0^i(\mu_{ji}(n{'}))$};
\node [above=of E] (D) {$\kappa_0^i(\mu_{ji}(n))$};

\draw[->] (E)--(F) node [midway,below]{$\lambda^i_{\mu_{ji}(n)\mu_{ji}(n{'})}$};
\draw[->] (D)--(A) node [midway,above] {$\kappa^i_{\mu_{ji}(n)\mu_{ji}(n{'})}$};
\draw[->] (D)--(E) node [midway,left] {$\Psi^i_{\mu_{ji}(n)}$};
\draw[->] (A)--(F) node [midway,right] {$\Psi^i_{\mu_{ji}(n{'})}$};

\end{tikzpicture}
\end{center}
we get $\Psi^i_{\mu_{ji}(n{'})} \circ \kappa^i_{\mu_{ji}(n)\mu_{ji}(n{'})} =
\lambda^i_{\mu_{ji}(n)\mu_{ji}(n{'})} \circ \Psi^i_{\mu_{ji}(n)}$, and hence
\begin{align*}
\Psi^i_{\mu_{ji}(n{'})} \circ \kappa^{ji}_{n{'} \mu_{ji}(n{'})} & = \Psi^i_{\mu_{ji}(n{'})} 
\circ \big(\kappa^i_{\mu_{ji}(n)\mu_{ji}(n{'})} \circ \kappa^{ji}_{n{'}\mu_{ji}(n)}\big)\\
& := \big(\Psi^i_{\mu_{ji}(n{'})} \circ \kappa^i_{\mu_{ji}(n)\mu_{ji}(n{'})}\big) \circ \kappa^{ji}_{n{'}\mu_{ji}(n)}\\
& = \big(\lambda^i_{\mu_{ji}(n)\mu_{ji}(n{'})} \circ \Psi^i_{\mu_{ji}(n)}\big) \circ \kappa^{ji}_{n{'}\mu_{ji}(n)},
\end{align*}
\begin{align*}
\big[\gamma_{ij}(\Psi^i)\big]_{n{'}}  \circ \kappa^j_{nn{'}} & := \big[\lambda^{ij}_{\mu_{ji}(n{'})n{'}} 
\circ \big(\Psi^i_{\mu_{ji}(n{'})} \circ \kappa^{ji}_{n{'} \mu_{ji}(n{'})}\big)\big] \circ \kappa^j_{nn{'}}\\
& = \big[\lambda^{ij}_{\mu_{ji}(n{'})n{'}} \circ \big(\big(\lambda^i_{\mu_{ji}(n)\mu_{ji}(n{'})} \circ 
\Psi^i_{\mu_{ji}(n)}\big) \circ \kappa^{ji}_{n{'}\mu_{ji}(n)}\big)\big] \circ \kappa^j_{nn{'}}\\
& := \big(\lambda^{ij}_{\mu_{ji}(n{'})n{'}} \circ \lambda^i_{\mu_{ji}(n)\mu_{ji}(n{'})}\big) \circ
\Psi^i_{\mu_{ji}(n)} \circ \big(\kappa^{ji}_{n{'}\mu_{ji}(n)} \circ \kappa^j_{nn{'}}\big)\\
& = \lambda^{ij}_{\mu_{ji}(n)n{'}} \circ \Psi^i_{\mu_{ji}(n)}  \circ \kappa^{ji}_{n\mu_{ji}(n)}\\
& = \lambda^j_{nn{'}} \circ \big(\lambda^{ij}_{\mu_{ji}(n)n} \circ \Psi^i_{\mu_{ji}(n)}  \circ 
\kappa^{ji}_{n\mu_{ji}(n)}\big)\\
& := \lambda^j_{nn{'}} \circ \big[\gamma_{ij}(\Psi^i)\big]_{n}^j.
\end{align*}
If $\Psi^i = \Phi^i$, we show that $\big[\gamma_{ij}(\Psi^i)\big]^j = \big[\gamma_{ij}(\Phi^i)\big]^j$. 
As 
\[ [\gamma_{ij}(\Psi^i)]_n^j :=  \lambda^{ij}_{\mu_{ji}(n)n} \circ \Psi^i_{\mu_{ji}(n)}  \circ
\kappa^{ji}_{n \mu_{ji}(n)} \ \ \ \& \ \ \ 
[\gamma_{ij}(\Phi^i)]_n^j := \lambda^{ij}_{\mu_{ji}(n)n} \circ \Phi^i_{\mu_{ji}(n)} \circ \kappa^{ji}_{n \mu_{ji}(n)}, \]
and since $\Psi^i = \Phi^i$, we get $\Psi^i_{\mu_{ji}(n)} = \Phi^i_{\mu_{ji}(n)}$, and hence the 
following diagram commutes
\begin{center}
\begin{tikzpicture}

\node (E) at (0,0) {$\kappa_0^j(n)$};
\node[right=of E] (K) {};
\node[right=of K] (L) {};
\node[right=of L] (M) {};
\node[right=of M] (F) {$\lambda_0^j(n)$.};
\node[above=of F] (A) {$\lambda_0^j(n)$};
\node [above=of E] (D) {$\kappa_0^j(n)$};

\draw[->] (E)--(F) node [midway,below]{$\lambda^{ij}_{\mu_{ji}(n)n} \circ \Phi^i_{\mu_{ji}(n)}  
\circ \kappa^{ji}_{n \mu_{ji}(n)}$};
\draw[->] (D)--(A) node [midway,above] {$\lambda^{ij}_{\mu_{ji}(n)n} \circ \Psi^i_{\mu_{ji}(n)}  
\circ \kappa^{ji}_{n \mu_{ji}(n)}$};
\draw[->] (D)--(E) node [midway,left] {$\kappa^j_{nn}$};
\draw[->] (A)--(F) node [midway,right] {$\lambda^j_{nn}$};

\end{tikzpicture}
\end{center}
(ii) If $m \in \mu_0 (i)$, then $[\gamma_{ii}(\Psi^i)]_m^i := \lambda^{ii}_{\mu_{ii}(m)m} \circ 
\Psi^i_{\mu_{ii}(m)}  \circ \kappa^{ii}_{m \mu_{ii}(m)} := \lambda^{ii}_{mm} \circ \Psi^i_{m} 
\circ \kappa^{ii}_{mm}
:= \id_{\lambda_0^i (m)}  \circ \Psi^i_{m}  \circ \id_{\kappa_0^i (m)}
:= \Psi^i_{m},$
hence $[\gamma_{ii}(\Psi^i)]_m^i := \Psi^i_m$, and consequently $\big[\gamma_{ii}(\Psi^i)\big]^i 
:= \Psi^i$. For the commutativity of the diagram
\begin{center}
\begin{tikzpicture}

\node (E) at (0,0) {$\Map\big(K_j^{0,M}, \Lambda_j^{0,M}\big)$};
\node[right=of E] (L) {};
\node[right=of L] (F) {$\Map\big(K_k^{0,M}, \Lambda_k^{0,M}\big)$};
\node [above=of E] (D) {$\Map\big(K_i^{0,M}, \Lambda_i^{0,M}\big)$};

\draw[->] (E)--(F) node [midway,below] {$\gamma_{jk}$};
\draw[->] (D)--(E) node [midway,left] {$\gamma_{ij}$};
\draw[->] (D)--(F) node [midway,right] {$\ \gamma_{ik}$};

\end{tikzpicture}
\end{center}
we need to show the equality between the  maps
\[ \chi^k_l := \lambda^{jk}_{\mu_{kj}(l)l} \circ \big[\gamma_{ij}(\Psi^i)\big]^j_{\mu_{kj}(l)}
\circ \kappa^{kj}_{l \mu_{kj}(l)} \ \ \ \& \ \ \ 
\upsilon^k_l := \lambda^{ik}_{\mu_{ki}(l)l} \circ \Psi^j_{\mu_{ki}(l)}  \circ \kappa^{ki}_{l \mu_{ki}(l)}. \]
By the definition of $\big[\gamma_{ij}(\Psi^i)\big]^j_{\mu_{kj}(l)}$ we get
\begin{align*}
\chi^k_l & := \lambda^{jk}_{\mu_{kj}(l)l} \circ \big(\lambda^{ij}_{\mu_{ji}(\mu_{kj}(l))\mu_{kj}(l)}
\circ \Psi^i_{\mu_{ji}(\mu_{kj}(l))}  \circ \kappa^{ji}_{\mu_{kj}(l) \mu_{ji}(\mu_{kj}(l))}\big) 
\circ \kappa^{kj}_{l \mu_{kj}(l)}\\
& = \big(\lambda^{jk}_{\mu_{kj}(l)l} \circ \lambda^{ij}_{\mu_{ji}(\mu_{kj}(l))\mu_{kj}(l)}\big) 
\circ \Psi^i_{\mu_{ji}(\mu_{kj}(l))}  \circ \big(\kappa^{ji}_{\mu_{kj}(l) \mu_{ji}(\mu_{kj}(l))} 
\circ \kappa^{kj}_{l \mu_{kj}(l)}\big)\\
& = \lambda^{ik}_{\mu_{ji}(\mu_{kj}(l))l} \circ \Psi^i_{\mu_{ji}(\mu_{kj}(l))}  \circ
\kappa^{ki}_{l \mu_{ji}(\mu_{kj}(l))}.
\end{align*}
By the supposed commutativity of the following diagram 
\begin{center}
\begin{tikzpicture}

\node (E) at (0,0) {$\lambda_0^i (\mu_{ki}(l))$};
\node[right=of E] (K) {};
\node[right=of K] (L) {};
\node[right=of L] (F) {$\lambda_0^i (\mu_{ji}(\mu_{kj}(l)))$,};
\node[above=of F] (A) {$\kappa_0^i(\mu_{ji}(\mu_{kj}(l)))$};
\node [above=of E] (D) {$\kappa_0^i(\mu_{ki}(l))$};

\draw[->] (E)--(F) node [midway,below]{$\lambda^i_{\mu_{ki}(l)\mu_{ji}(\mu_{kj}(l))}$};
\draw[->] (D)--(A) node [midway,above] {$\kappa^i_{\mu_{ki}(l)\mu_{ji}(\mu_{kj}(l))}$};
\draw[->] (D)--(E) node [midway,left] {$\Psi^i_{\mu_{ki}(l)}$};
\draw[->] (A)--(F) node [midway,right] {$\Psi^i_{\mu_{ji}(\mu_{kj}(l))}$};

\end{tikzpicture}
\end{center}
\begin{align*}
\upsilon^k_l & := \lambda^{ik}_{\mu_{ki}(l)l} \circ \Psi^j_{\mu_{ki}(l)}  \circ 
\kappa^{ki}_{l \mu_{ki}(l)}\\
& = \lambda^{ik}_{\mu_{ji}(\mu_{kj}(l))l} \circ \big(\lambda^{i}_{\mu_{ki}(l)\mu_{ji}(\mu_{kj}(l))}
\circ \Psi^j_{\mu_{ki}(l)} \big) \circ \kappa^{ki}_{l \mu_{ki}(l)}\\
& = \lambda^{ik}_{\mu_{ji}(\mu_{kj}(l))l} \circ \big(\Psi^i_{\mu_{ji}(\mu_{kj}(l))}
\circ \kappa^i_{\mu_{ki}(l)\mu_{ji}(\mu_{kj}(l))} \big) \circ \kappa^{ki}_{l \mu_{ki}(l)}\\
& = \lambda^{ik}_{\mu_{ji}(\mu_{kj}(l))l} \circ \Psi^i_{\mu_{ji}(\mu_{kj}(l))} \circ 
\big(\kappa^i_{\mu_{ki}(l)\mu_{ji}(\mu_{kj}(l))} \circ \kappa^{ki}_{l \mu_{ki}(l)}\big)\\
& = \lambda^{ik}_{\mu_{ji}(\mu_{kj}(l))l} \circ \Psi^i_{\mu_{ji}(\mu_{kj}(l))} \circ
\kappa^{ki}_{l \mu_{ji}(\mu_{kj}(l))}\\
& = \chi^k_l.\qedhere
\end{align*}
\end{proof}

\begin{definition}\label{def: famfammap}
If $(K^i)_{i \in I}, (\Lambda^i)_{i \in I}  \in \Fam(I, M)$, \textit{a family of families-map}
from $(K^i)_{i \in I}$ to $(\Lambda^i)_{i \in I}$\index{family of families-map}, in symbols 
$\Psi \colon (K^i)_{i \in I} \To (\Lambda^i)_{i \in I}$,
is a dependent operation
$\Psi \colon \bigcurlywedge_{i \in I}\Map\big(K_i^{0,M}, \Lambda_i^{0,M}\big)$ 
such that for every $(i, j) \in D(I)$ the following diagram commutes
\begin{center}
\begin{tikzpicture}

\node (E) at (0,0) {$\Lambda_i^{0,M}$};
\node[right=of E] (L) {};
\node[right=of L] (F) {$\Lambda_j^{0,M}$,};
\node[above=of F] (A) {$K_j^{0,M}$};
\node [above=of E] (D) {$K_i^{0,M}$};

\draw[->,double equal sign distance] (E)--(F) node [midway,below]{$\Phi_{ij}^{\Lambda}$};
\draw[->,double equal sign distance] (D)--(A) node [midway,above] {$\Phi_{ij}^K$};
\draw[->,double equal sign distance] (D)--(E) node [midway,left] {$\Psi^i$};
\draw[->,double equal sign distance] (A)--(F) node [midway,right] {$\Psi^j$};

\end{tikzpicture}
\end{center}
where $\Phi_{ij}^K$ and $\Phi_{ij}^{\Lambda}$ are the transport family-maps of
$(K^i)_{i \in I}$ and $(\Lambda^i)_{i \in I}$, respectively, according to 
Definition~\ref{def: transportfammaps}. If $\Xi \colon(\Lambda^i)_{i \in I} \To (N^i)_{i \in I}$,
the composition $\Xi \circ \Psi \colon (K^i)_{i \in I} \To (N^i)_{i \in I}$ is defined, for every $i \in I$, by 
$(\Xi \circ \Psi)^i := \Xi^i \circ \Psi^i$
\begin{center}
\begin{tikzpicture}

\node (E) at (0,0) {$K_i^{0,M} \ $};
\node[right=of E] (X) {};
\node[right=of X] (F) {$K_j^{0,M}$};
\node[below=of F] (A) {$ \ \Lambda_j^{0,M}$};
\node[left=of A] (Y) {};
\node[left=of Y] (B) {$\Lambda_i^{0,M} \ $};
\node[below=of B] (K) {$N_i^{0,M}$};
\node[right=of K] (X) {};
\node[right=of X] (L) {$ N_j^{0,M} $.};

\draw[->,double equal sign distance] (E)--(F) node [midway,above] {$\Phi_{ij}^K$};
\draw[->,double equal sign distance] (F)--(A) node [midway,left] {$\Psi^j$};
\draw[->,double equal sign distance] (B)--(A) node [midway,below] {$\Phi_{ij}^{\Lambda}$};
\draw[->,double equal sign distance] (E) to node [midway,right] {$\Psi^i$} (B);
\draw[->,double equal sign distance] (B) to node [midway,right] {$\Xi^i$} (K);
\draw[->,double equal sign distance] (A)--(L) node [midway,left] {$\Xi^j$};
\draw[->,double equal sign distance] (K)--(L) node [midway,below] {$\Phi_{ij}^N$};
\draw[->,double equal sign distance,bend right=40] (E) to node [midway,left] {$(\Xi \circ \Psi)^i  \ $} (K);
\draw[->,double equal sign distance,bend left=40] (F) to node [midway,right] {$ \  (\Xi \circ \Psi)^j$} (L);

\end{tikzpicture}
\end{center} 
The \textit{identity family of families-map}\index{identity family of families-map}\index{$\Id_{(\Lambda^i)_{i \in I}}$} 
$\Id_{(\Lambda^i)_{i \in I}}$ is defined by the 
rule $\big[\Id_{(\Lambda^i)_{i \in I}}\big]^i := \Id_{\Lambda_i^{0,M}}$, for every $i \in I$. 
The totality of family of families-maps from $(K^i)_{i \in I}$ to $(\Lambda^i)_{i \in I}$, and 
the canonical equality on $\Fam(I,M)$ is defined in analogy to Definition~\ref{def: map}.
\end{definition}

If $\Psi \colon (K^i)_{i \in I} \To (\Lambda^i)_{i \in I}$, the commutativity of the 
diagram in Definition~\ref{def: famfammap} is unfolded as follows. If $i =_I j$ and $m \in \mu_0(i)$, then 
\begin{align*}
\big[\Psi^j \circ \Phi_{ij}^K\big]_m = \big[\Phi_{ij}^{\Lambda} \circ \Psi^i\big]_m & 
:\TOT \Psi_{\mu_{ij}(m)}^j \circ \big[\Phi_{ij}^K\big]_m = \big[\Phi_{ij}^{\Lambda}\big]_m \circ \Psi_m^i\\
& :\TOT \Psi_{\mu_{ij}(m)}^j \circ \kappa^{ij}_{m \mu_{ij}(m)} = \lambda^{ij}_{m \mu_{ij}(m)} \circ \Psi_m^i
\end{align*}
i.e., the following diagram commutes
\begin{center}
\begin{tikzpicture}

\node (E) at (0,0) {$\lambda_0^i(m)$};
\node[right=of E] (L) {};
\node[right=of L] (F) {$\lambda^j_0(\mu_{ij}(m))$.};
\node[above=of F] (A) {$\kappa^j_0(\mu_{ij}(m))$};
\node[left=of A] (M) {};
\node [left=of M] (D) {$\kappa^i_0(m) \ $};

\draw[->] (E)--(F) node [midway,below] {$\lambda^{ij}_{m \mu_{ij}(m)}$};
\draw[->] (D)--(A) node [midway,above] {$\kappa^{ij}_{m \mu_{ij}(m)}$};
\draw[->] (D)--(E) node [midway,left] {$\Psi^i_m$};
\draw[->] (A)--(F) node [midway,right] {$\Psi^j_{\mu_{ij}(m)}$};

\end{tikzpicture}
\end{center}

In analogy to Corollary~\ref{cor: mapdependent} we have the following.

\begin{corollary}\label{cor: Mapdependent}
If $(K^i)_{i \in I}, (\Lambda^i)_{i \in I}  \in \Fam(I, M)$ and 
$\Psi \colon \bigcurlywedge_{i \in I}\Map\big(K_i^{0M}, \Lambda_i^{0,M}\big)$, the following are equivalent:\\[1mm]
\normalfont (i)
\itshape $\Psi \colon (K^i)_{i \in I} \To (\Lambda^i)_{i \in I} $.\\[1mm]
\normalfont (ii) 
\itshape $\Psi \in \prod_{i \in I}\Map\big(K_i^{0M}, \Lambda_i^{0,M}\big)$.

\end{corollary}

\begin{proof}
If $i =_I j$, the commutativity of the diagram in the definition of a family of 
families-map $\Psi \colon (K^i)_{i \in I} \To (\Lambda^i)_{i \in I} $ is equivalent to
the membership condition $\Psi \in \prod_{i \in I}\Map\big(K_i^{0M}, \Lambda_i^{0,M}\big)$ 
using the above unfolding of the equality $\big[\Psi^j \circ \Phi_{ij}^K\big]_m = 
\big[\Phi_{ij}^{\Lambda} \circ \Psi^i\big]_m$.
\end{proof}

\begin{definition}\label{def: dependentfamoffam}
The totality $\prod_{i \in I}\prod_{m \in \mu_0 (i)}\lambda_0^i (m)$  of dependent functions over
a family of families of set $(\Lambda^i)_{i \in I} \in \Fam(I,M)$ is defined by 
\[ \Theta \in \prod_{i \in I}\prod_{m \in \mu_0 (i)}\lambda_0^i (m) :\TOT \Theta \colon
\bigcurlywedge_{i \in I}\bigcurlywedge_{m \in \mu_0(i)}\lambda_0^i(m) \  \& \ \forall_{(i,j)
\in D(I)}\forall_{(m,n) \in T_{ij}(M)}\big(\Theta^j_n =_{\lambda^j_0(n)} \lambda^{ij}_{mn}(\Theta^i_m)\big), \]
\[ \Theta =_{\mathsmaller{\prod_{i \in I}\prod_{m \in \mu_0 (i)}\lambda_0^i (m)}} \Phi 
:\TOT \forall_{i \in I}\forall_{m \in \mu_0(i)}\big(\Theta^i_m =_{\lambda_0^i(m)} \Phi^i_m\big).\]

\end{definition}

The theory of families of families of sets over $(I,M)$ within $\D V_0^{\im}$ can be developed 
further along the lines of the
theory of families of sets over $I$ within $\D V_0$.

\section{Notes}
\label{sec: notes3}

\begin{note}\label{not: deffamofsets}
\normalfont
The concept of a family of sets indexed by a (discrete) set was asked to be defined 
in~\cite{Bi67}, Exercise 2, p.~72, and the required definition, given by Richman, 
is included in~\cite{BB85}, Exercise 2, p.~78, where the discreteness hypothesis is omitted.
The definition has a strong type-theoretic flavour, although, Richman's motivation had categorical origin,
rather than type-theoretic. In a personal communication regarding this definition, Richman 
referred to the definition of a set-indexed family of objects of a category, given in~\cite{MRR88}, p.~18,
as the source of the definition attributed to him in~\cite{BB85}, p.~78. 
Given the categorical flavour of Bishop's notion of a subset, it might be that Bishop was also 
thinking in categorical terms, although Bishop, to our knowledge, neither used a purely categorical language to describe his concepts, nor he used general category theory as a foundational framework for $\BISH$.

Specifically, in~\cite{MRR88} Richman 
presented a set $I$ as a category with objects its elements and 
\[ \Hom_{=_I}(i, j) := \{x \in \{0\} \mid i =_I j\}, \]
for every $i, j \in I$. If we view $\D V_0$ 
as a category with objects its elements and 
\[ \Hom_{=_{\D V_0}}(X, Y) := \big \{(f, f{'}) : \D F(X, Y) \times \D F(Y, X) \mid (f, f{'}) : 
X =_{\D V_0} Y \big\}, \]
for every $X, Y \in \D V_0$, then
an $I$-family of sets is a functor from the category $I$ to the category $\D V_0$.
Notice that in the definitions of $\Hom_{=_I}(i, j)$
and of $\Hom_{=_{\D V_0}}(X, Y)$ the properties $P(x) := i =_I j$ and $Q(f, f{'}) := (f, f{'}) : 
X =_{\D V_0} Y$ are extensional.
In~\cite{Pe19c} we reformulated Richman's definition using the universe 
$\D V_0$ of sets and the universe $\D V_1$ of triplets $(A, B, f)$, where $A, B \in \D V_0$ and $f \colon A \to B$.
Definition~\ref{def: famofsets} rests on the notion of dependent operation, in order to be absolutely 
faithful to Bishop's account of sets and functions in~\cite{Bi67} and~\cite{BB85}.
For the definition of the concept of a family of sets in $\ZF$, or $\CZF$, see~\cite{Pa12a}, p.~35, and Note~\ref{not: onfams(C)ZF}.

The term ``transport map'' in Definition~\ref{def: famofsets} is drawn from MLTT. Actually, Definition~\ref{def: famofsets}
is a ``definitional form'' of the type-theoretic transport i.e., the existence of the transport map 
$p_* \colon P(x) \to P(y)$, where $p \colon x =_A y$ and $P \colon A \to \C U$ is a type-family over $A \colon \C U$ in
the universe of types $\C U$.
In $\MLTT$ the existence of $p_*$ follows from Martin-L\"of's $J$-rule, the induction principle that 
accommodates the indentity type-family $=_A \colon A \to A \to \C U$, for every type $A \colon \C U$.
In Definition~\ref{def: famofsets} we describe in a proof-irrelevant way i.e., using only the fact that $i =_I j$ and not referring to witnesses of this equality, a structure of transport maps. This structure in $\BST$ is defined, and not generated from the equality type family of $\MLTT$.

\end{note}

\begin{note}\label{not: defmap}
\normalfont
In the categorical setting of Richman (see Note~\ref{not: deffamofsets}), a family map $\Psi \in \Map_I(\Lambda, M)$ is
a natural transformation from the functor $\Lambda$ to the functor $M$. The fact that the most fundamental 
concepts of category theory, that of a functor and of a natural transformation, are formulated in a
natural way in $\BST$ through the notion of a dependent operation explains why category theory 
is so closely connected to $\BST$. For more on the connections between $\BST$, dependent type theory and category theory see section~\ref{sec: typescats}. 
\end{note}

\begin{note}\label{not: exteriorunion1}
\normalfont
The exterior union, is necessary to the definition of the infinite product of a sequence of sets.
In~\cite{BB85}, p.~125, the following is noted:
\begin{quote}
Within the main body of this text, we have only defined the product
of a family of subsets of a given set. However, with the aid of
Problem 2 of Chapter 3 we can define the product of an arbitrary
sequence of sets. Definition (1.7) then applies to such a product\footnote{This is the
definition of the countable
product of metric spaces.}.
\end{quote}
\end{note}

\begin{note}\label{not: defprod}
\normalfont
If $\Lambda^{\Nat}$ is the sequence of sets defined in Definition~\ref{def: basicfamilies}, 
the definitional clauses of the corresponding exterior union can be written as follows: 
\[ \sum_{n \in \Nat}X_n =: \{(n, x) \mid n \in \Nat \ \& \ x \in X_n\}, \]
\[ (n, x) =_{\sum_{n \in \Nat}X_n} (m, y) : \TOT n =_{\Nat} m \ \& \ x =_{X_n} y. \]
Traditionally, the countable product of this sequence of sets is defined by
\[ \prod_{n \in \Nat}X_n := \bigg\{\phi \colon \Nat \to \sum_{n \in \Nat}X_n  \mid
\forall_{n \in \Nat}\big(\phi(n) \in X_n\big) \bigg\}, \]
which is a rough writing of the following
\[ \prod_{n \in \Nat}X_n := \bigg\{\phi \colon \Nat \to \sum_{n \in \Nat}X_n  \mid 
\forall_{n \in \Nat}\big(\prb_1(\phi(n)) =_{\Nat} n\big) \bigg\}. \]
In the second writing
$\prb_1(\phi(n)) =_{\Nat} n$ implies that $\prb_1 (\phi(n) := n$, 
hence, if $\phi(n) := (m, y)$, then $m = n$ and $y \in X_n$. When the equality of $I$ though, is not like
that of $\Nat$, we cannot solve this problem in a satisfying way. Although Bishop did not consider
products other than countable ones, in more  abstract areas of mathematics, like e.g., the general topology
of Bishop spaces, arbitrary products are considered (see~\cite{Pe15}). One could have defined 
\[ \Phi \in  \prod_{i \in I}\lambda_0(i) : \TOT \Phi \in \D F \bigg(I, \sum_{i \in I}\lambda_0 (i)\bigg)
\ \& \ \forall_{i \in I}\big(\prb_1(\Phi(i)) := i\big). \]
This approach has the problem that
the property 
\[ Q(\Phi) : \TOT \forall_{i \in I}\big(\prb_1(\Phi(i)) := i\big) \]
is not necessarily extensional; let $\Phi =_{\mathsmaller{\D F (I, \sum_{i \in I}\lambda_0 (i))}} \Theta$ i.e.,
$\forall_{i \in I}\big(\Phi(i) =_{\mathsmaller{\sum_{i \in I}\lambda_0 (i)}} \Theta(i)\big)$, and suppose that 
$Q(\Phi)$. If we fix some $i \in I$, and $\Phi(i) := (i, x)$ and $\Theta(i) := (j, y)$, we only 
get that $j =_I i$. The use of dependent operations allows us to 
define the right analogue to the $\prod$-type of $\MLTT$ and being at the same time compatible 
with the use of dependent operations by Bishop in~\cite{Bi67}, p.~65.
\end{note}

\begin{note}\label{not: countableprodmetric}
\normalfont
A precise formulation of the definition in~\cite{BB85}, p.~85, 
of the countable product of a sequence $\big(X_n, \rho_n\big)_{n \in \Nat}$ of 
metric spaces, where $\rho_n$ is bounded by $1$, for every $n \in \Nat$, is the following. Let
$\Lambda^{\Nat} := (\lambda_0^{\Nat}, \lambda_1^{\Nat})$ be the $\Nat$-family of the sets $(X_n)_{n \in \Nat}$ (see
Definition~\ref{def: basicfamilies}). Notice that the dependent operation $\lambda_1^{\Nat}$ is compatible to 
the corresponding metric structures in the sense that each transport map $\lambda^{\Nat}_{nn} := \id_{X_n}$ is a 
morphism in any category of metric spaces considered. This is an example of a \textit{spectrum of metric spaces} 
over $\Lambda^{\Nat}$\index{spectrum of metric spaces over a family of sets} (see also the introduction to section~\ref{sec: spectra}). The 
\textit{countable product metric}\index{countable product metric} $\rho_{\infty}$\index{$\rho_{\infty}$}
on $\prod_{n \in \Nat}X_n$, for every $\Phi, \Theta \in \prod_{n \in \Nat}X_n$, is defined by
\[ \rho_{\infty}(\Phi, \Theta) := \sum_{n = 0}^{\infty}\frac{\rho_n\big(\Phi_n, \Theta_n\big)}{2^n}. \]
\end{note}

\begin{note}\label{not: dependentapplication}
\normalfont
The equality 
$ \Phi_j =_{\lambda_0 (j)} \lambda_{ij}(\Phi_i)$
in Definition~\ref{def: piset} is the proof-irrelevant version of dependent application of a dependent function in $\MLTT$ 
(see also Note~\ref{not: onprsigma}).
\end{note}

\begin{note}\label{not: AC}
\normalfont
As it is mentioned in~\cite{Pa13}, the axiom of choice is ``freely used in Bishop constructivism''. In
Theorem~\ref{thm: ac} we show only the formal version of the type-theoretic axiom choice within $\BST$ i.e., the 
the distributivity of $\prod$ over $\sum$. This term was suggested to us by M.~Maietti. In~\cite{Pe19c} a 
proof of this result is also given, where dependecy is formulated with the help of the universe $\D V_1$ 
of triplets $(A, B, f)$ (see Note~\ref{not: deffamofsets}).
As it was first noted to us by E.~Palmgren, this distributivity 
holds in every locally cartesian closed category. In~\cite{Wi13} it is mentioned that this fact
is attributed to Martin-L\"of and his work~\cite{ML84}. For a proof see~\cite{Aw95}.

\end{note}

\begin{note}\label{not: deflambdaI}
\normalfont
The notion of an $I$-set of sets is in accordance with Bishop's predicative spirit, and his need to avoid
the treatment of the universe $\D V_0$ as a set. This notion was not defined by Bishop, only its ``internal'' version, the notion
of an $I$-set of subsets, was defined similarly by him in~\cite{Bi67}, p.~65. The use of the term ``set of subsets'' was
a source of misreading of~\cite{Bi67} from the side of Myhill in~\cite{My75} (see also Notes~\ref{not: BCms} 
and~\ref{not: numerical}). 
The definition of the set $\lambda_0 I$ is in the spirit of the definition of the quotient group $G/H$ of the group 
$G$ by its normal subgroup $H$, given in~\cite{MRR88}, p.~38. 
If $I$ is equipped with the equality $=_I^{\Lambda}$, then $\Lambda$ does not become necessarily an 
$I$-set of sets. The reason for this is that the transport maps of $\Lambda$ are given beforehand, and if we equip $I$
with $=_I^{\Lambda}$ we need to add a transport map $\lambda_{ij}$ for every pair $(i,j)$ for which $\lambda_0(i) =_{\D V_0} 
\lambda_0(j)$ and $(i,j) \notin D(I)$, where $D(I)$ is understood here as the diagonal $D(I, =_I)$ with respect to the 
equality $=_I$. So, $\lambda_1$ has to be extended, and define a new family of sets over $(I, =_I^{\Lambda})$, 
which is going
to be an $(I, =_I^{\Lambda})$-set of sets.

\end{note}

\begin{note}\label{not: defdirfamofsets}
\normalfont
A direct family of sets is a useful variation of the notion of a set-indexed family of sets (see Chapter~\ref{chapter: bspaces}). 
A directed set $(I, \lt_I)$ can also be seen as
a category with objects
the elements of $I$, and $\Hom_{\lt_I}(i, j) := \{x \in \{0\} \mid i \lt_I j\}.$
If the universe $\D V_0$ is seen as a category with objects its elements and 
$\Hom_{\lt_{\D V_0}}(X, Y) :=  \D F(X, Y),$ an $(I, \lt_I)$-family of sets is a
functor from the category $(I, \lt_I)$ to this new category $\D V_0$.

\end{note}

\begin{note}\label{not: preorder}
\normalfont
A generalisation of the notion of a direct family of sets is that of a preorder family of sets. If $(I, \lt_I)$ is a 
preorder (see Definition~\ref{def: dirset}), a \textit{covariant preorder family of sets}
\index{covariant preorder family of sets} over $(I, \lt_I)$ is defined as a direct family of sets. One needs though the 
property of a directed set to define an interesting equality on the exterior union of the corresponding family.
A \textit{contravariant preorder family of sets} over $(I, \lt_I)$\index{contravariant preorder family of sets},
or an $(I, \mt_I)$-family of 
sets\index{$(I, \mt_I)$-family of sets}, is a pair 
$M^{\mt} := (\mu_0, \mu_1^{\mt})$, where if $(j, i) \in D^{\mt}(I)$, the transport maps
$\mu_1^{\mt}(j, i) \colon \mu_0 (j) \to \mu_0 (i)$ behave in a dual way i.e., for every
$i, j, k \in I$ with $k \mt_I j \mt_I i$, the following diagram commutes
\begin{center}
\begin{tikzpicture}

\node (E) at (0,0) {$\mu_0(j)$};
\node[right=of E] (F) {$\mu_0(k).$};
\node [above=of E] (D) {$\mu_0(i)$};

\draw[->] (F)--(E) node [midway,below] {$\mu_{jk}^{\mt}$};
\draw[->] (E)--(D) node [midway,left] {$\mu_{ij}^{\mt}$};
\draw[->] (F)--(D) node [midway,right] {$\ \mu_{ik}^{\mt}$};

\end{tikzpicture}
\end{center}
If $(I, \lt)$ is an inverse-directed set (see Definition~\ref{def: dirset}) and $M^{\mt}$ is an
$(I, \mt_I)$-contravariant direct family of sets, defined in the obvious way, the 
\textit{inverse-direct sum}\index{inverse-direct sum} $\sum_{i \in I}^{\mt} \mu_0(i)$\index{$\sum_{i \in I}^{\mt} \mu_0(i)$}
of $M^{\mt}$ is the totality $\sum_{i \in I}\mu_0(i)$, equipped with the equality
\[ (i, x) =_{\sum_{i \in I}^{\mt} \mu_0(i)} (j, y) : \TOT \exists_{k \in I}\big(i \mt_I k \ 
\& \ j \mt_I k \ \& \ \mu_{ik}^{\mt}(x) =_{\mu_0(k)} \mu_{jk}^{\mt}(y)\big). \]
The set $\prod_{i \in I}^{\mt}\mu_0(i)$ is defined in the expected way. 
Thinking classically, a topology $T$ of
open sets on a set $X$, equipped with the subset order $\subseteq$, is an inverse-directed set, and the notion of 
a \textit{presheaf of sets} on $(X, T)$ is an example of a $(T, \supseteq)$-contravariant direct family of sets. 
In the language of presheaves (see~\cite{Jo02}, p.~72) the transport maps $\mu_{ij}^{\mt}$ are called 
\textit{restriction maps}, and a family-map $\Phi \colon \Lambda^{\mt} \To M^{\mt}$ is called a morphism of presheaves.
It is natural to use also the term \textit{extension map} for the transport map $\lambda_{ij}^{\lt}$ of a 
covariant (direct) preorder family of sets.
The notion of a family of sets over a partial order is also used in the definition of a Kripke model 
for intuitionistic predicate logic. For that see~\cite{TD88I}, p.~85, where the 
transport maps $\lambda_{ij}^{\lt}$ are called there \textit{transition functions}.

\end{note}

\begin{note}\label{not: onsetrelfams}
\normalfont
If a set-relevant family-map $\Psi \colon \Lambda^* \To M^*$ was defined by the stronger condition: for every $(i,j) \in D(I)$, every $p \in \varepsilon_0^{\lambda}(i,j)$ and \textit{every} $q \in \varepsilon_0^{\mu}(i,j)$ the diagram in Definition~\ref{def: setrelmap} commutes, then the expected fact $\id_{\Lambda^*} \colon \Lambda^* \To \Lambda^*$ implies that $\lambda_{ij}^p = \lambda_{ij}^q$, for every $p, q \in \varepsilon_0^{\lambda}(i, j)$. This property is called proof-irrelevance in Definition~\ref{def: prfamofsets}.

\end{note}

\begin{note}\label{not: onfamsoffams}
\normalfont
The theory of families of families of sets over $(I,M)$ within $\D V_0^{\im}$ is the third rung of the
ladder of set-like objects in $\D V_0^{\im}$.  The first three rungs can be described as follows:
\[ X, Y \in \D V_0, \ \ \ \ \ \ \ \fXY,\]
\[ \Lambda, M \in \Fam(I), \ \ \ \ \Psi \colon \Lambda \To M \TOT \Psi \in \prod_{i \in I}\D F\big(\lambda_0(i),
\mu_0(i)\big),\]
\[ (K^i)_{i \in I}, (\Lambda^i)_{i \in I} \in \Fam(I,M), \ \ \ \ \Psi \colon (K^i)_{i \in I} 
\To (\Lambda^i)_{i \in I} \TOT \Psi \in \prod_{i \in I}\Map\big(K_i^{0,M}, \Lambda_i^{0,M}\big).\]
This hierarchy of universes and families can be extended further, if necessary.
\end{note}

\begin{note}[Small categories within $\BST$]\label{not: categories}
\normalfont
As it is mentioned in the introduction to Chapter 9 of~\cite{HoTT13}, where category theory is developed within $\HoTT$, categories do not fit well with set-based mathematics. Quit earlier, see e.g, in~\cite{HS98},  it is mentioned that ``type theory is adequate to represent faithfully categorical reasoning''. In~\cite{HS98} the objects are modelled as types and the Hom-sets as Hom-setoids of arrows, within the Calculus of Inductive Constructions. In~\cite{PW14} there are elements of such a development of category theory within type theory, where both the algebraic and the hom-definition are given. In~\cite{Pa17} are included interesting remarks on the formulation of category theory in~\cite{HoTT13}. For relations between category theory and Explicit Mathematics see~\cite{Ja19}. In this note we briefly explain why small categories fit well with $\BST$.

As we have already explained in Note~\ref{not: deffamofsets}, Richman used the notion of a functor to define the fundamental notion of a set-indexed family of sets, as a special case of a set-indexed family of objects in some category $\C C$.
Here we do the opposite. The notion of a set-indexed family of sets is fundamental and comes first. 
We use the basic theory of set-indexed families of sets to describe the basic notions of category theory within $\BST$.
In what follows we consider the objects of a category to be a set, although that could also be a class. The totality of arrows 
is always a set i.e., we could study locally small categories, but here we only present small categories.
A set is not  necessarily in the homotopy sense of the book-HoTT (see the corresponding notion 
of a strict category in~\cite{HoTT13}, section 9.6).
At this point we do not equip $\Ob_{\C C}$ with equality with evidence $(\EwE)$ that makes
possible the formulation of precategory and category in the sense of the book-HoTT (see section~\ref{sec: mlsets}).
For a general discussion on the relations between categories and sets in $\BST$ see section~\ref{sec: typescats}.

\begin{definition}\label{def: category}
A $($small$)$ \textit{category}\index{category} is a structure $\C C := \big(\Ob_{\C C}, \Morr^{\C C}_0, 
\Morr^{\C C}_1, \Comp^{\C C}, \Id^{\C C}\big)$, where $\Ob_{\C C}$ is a set, $\big(\Morr^{\C C}_0, 
\Morr^{\C C}_1\big) \in \Fam(\Ob_{\C C} \times \Ob_{\C C})$,
\[ \Comp^{\C C} : \bigcurlywedge_{x,y,z \in \Ob_{\C C}}\D F\big(\Morr^{\C C}_0(y,z) \times \Morr^{\C C}_0(x,y),
\Morr^{\C C}_0(x,z)\big),\] 
\[ \Comp^{\C C}_{xyz} := \Comp^{\C C}(x,y,z) : \Morr^{\C C}_0(y,z) \times \Morr^{\C C}_0(x,y) \to 
\Morr^{\C C}_0(x,z), \ \ \ \ 
\Comp^{\C C}_{xyz}(\phi, \psi) := \phi \circ \psi, \]
\[\Id^{\C C} : \bigcurlywedge_{x \in \Ob_{\C C}}\Morr^{\C C}_0(x,x), \ \ \ \ 
\Id^{\C C}_x := \Id^{\C C}(x); \ \ \ \  x \in \Ob_{\C C},\]
such that the following conditions are satisfied:\\[1mm]
$(\Cat_1)$ For every $x, y, z, w \in \Ob_{\C C}$, $\chi \in \Morr^{\C C}_0(z, w), \phi \in \Morr^{\C C}_0(y, z)$,
$\psi \in \Morr^{\C C}_0(x, y)$,
\begin{align*}
 \chi \circ (\phi \circ \psi) & := \Comp^{\C C}_{xzw}\big(\chi, \Comp^{\C C}_{xyz}(\phi, \psi)\big)\\
 & =_{\mathsmaller{\Morr^{\C C}_0(x, w)}} \Comp^{\C C}_{xyw}\big(\Comp^{\C C}_{yzw}(\chi, \phi), \psi\big)\\
 & := (\chi \circ \phi) \circ \psi.
\end{align*}
$(\Cat_2)$ For every $x, y \in \Ob_{\C C}$, and for every $\psi \in \Morr^{\C C}_0(x, y)$,
\[ \Id^{\C C}_y \circ \psi := \Comp^{\C C}_{xyy}(\Id^{\C C}_y, \psi) 
=_{\mathsmaller{\Mor^{\C C}_0(x, y)}} \psi \ \ \ \& \ \ \ 
\psi \circ \Id^{\C C}_x := \Comp^{\C C}_{xxy}(\psi, \Id^{\C C}_x, \psi) =_{\mathsmaller{\Morr^{\C C}_0(x, y)}} \psi. \]
$(\Cat_3)$ For every $x, y, z, x{'}, y{'}, z{'} \in \Ob_{\C C}$, with $x =_{\mathsmaller{\Ob_{\C C}}} x{'}$,
$y =_{\mathsmaller{\Ob_{\C C}}} y{'}$ and $z =_{\mathsmaller{\Ob_{\C C}}} z{'}$, for
every $\psi \in \Morr^{\C C}_0(x, y)$,
$\phi \in \Morr^{\C C}_0(y, z)$,
$\Morr^{\C C}_{(x,z)(x{'}, z{'})}(\phi \circ \psi) =_{\mathsmaller{\Morr^{\C C}_0(x{'}, z{'})}} 
\Morr^{\C C}_{(y,z)(y{'}, z{'})}(\phi) \circ \Morr^{\C C}_{(x,y)(x{'}, y{'})}(\psi)$

\begin{center}
\begin{tikzpicture}

\node (E) at (0,0) {$x{'}$};
\node[right=of E] (K) {};
\node[right=of K] (L) {};
\node[right=of L] (F) {$y{'}$};
\node [above=of E] (D) {$x$};
\node[right=of D] (X) {};
\node[right=of X] (W) {};
\node[right=of W] (A) {$ \ y$};
\node[right=of A] (Y) {};
\node[right=of Y] (T) {};
\node[right=of T] (P) {$ \ z$};
\node[right=of F] (M) {};
\node[right=of M] (N) {};
\node[right=of N] (O) {$z{'}$.};

\draw[->] (F)--(O) node [midway,above]{$\Morr_{(y,z)(y{'},z{'})}(\phi)$};
\draw[->] (E)--(F) node [midway,above]{$\Morr_{(x,y)(x{'},y{'})}(\psi)$};
\draw[->] (D)--(A) node [midway,above] {$\psi$};
\draw[->] (A)--(P) node [midway,above] {$\phi$};
\draw[-,double equal sign distance] (D)--(E) node [midway,left] {};
\draw[-,double equal sign distance] (A)--(F) node [midway,right] {};
\draw[-,double equal sign distance] (O)--(P) node [midway,right] {};
\draw[->,bend right=20] (E) to node [midway,below] {$\mathsmaller{\Morr^{\C C}_{(y,z)(y{'}, z{'})}(\phi) 
\ \circ  \ \Morr^{\C C}_{(x,y)(x{'}, y{'})}(\psi)}$} (O);
\draw[->,bend right=50] (E) to node [midway,below] {$\mathsmaller{\Morr^{\C C}_{(x,z)(x{'}, z{'})}(\phi
\circ \psi)}$} (O);

\end{tikzpicture}
\end{center}
$(\Cat_4)$ For every $x, x{'} \in \Ob_{\C C}$, with $x =_{\mathsmaller{\Ob_{\C C}}} x{'}$, 
$\Morr^{\C C}_{(x,x)(x{'}, x{'})}\big(\Id^{\C C}_x\big) =_{\mathsmaller{\Morr^{\C C}_0(x{'}, x{'})}} 
\Id^{\C C}_{x{'}}$
\begin{center}
\begin{tikzpicture}

\node (E) at (0,0) {$x{'}$};
\node[right=of E] (K) {};
\node[right=of K] (L) {};
\node[right=of L] (F) {$x{'}$.};
\node[above=of F] (A) {$x$};
\node [above=of E] (D) {$x$};

\draw[->] (E)--(F) node [midway,above]{$\Morr_{(x,x)(x{'},x{'})}(\Id_x)$};
\draw[->] (D)--(A) node [midway,above] {$\Id_x$};
\draw[-,double equal sign distance] (D)--(E) node [midway,left] {};
\draw[-,double equal sign distance] (A)--(F) node [midway,right] {};
\draw[->,bend right=40] (E) to node [midway,below] {$\Id_{x{'}}$} (F);

\end{tikzpicture}
\end{center}
\end{definition}

The last two conditions, which reflect a functorial behaviour of the transport maps of $\Morr^{\C C}_1$ and are not
found in the standard definition of a category, are necessary
compatibility conditions between these transport maps and the $\big(\Comp^{\C C}, \Id^{\C C}\big)$-structure
of the category
$\C C$. While in intensional $\MLTT$ these conditions follow from the transport, hence the $J$-rule, here we need to include them in our definition.

As a characteristic example of a category in the above sense, we consider the constructive analogue to the category of posets.
Classically, the category of posets has objects the collection of all posets and arrows the monotone functions. 
In order to formulate this constructively, we need to generalise Definition~\ref{def: category} to categories with
objects an abstract totality $\Ob_{\C C}$. In Definition~\ref{def: spectrumofposets} we define the category generated
by a spectrum of posets. We can define similarly the category generated by a spectrum of groups, rings, modules etc. 
(for the notion of an $S$-spectrum, where $S$ is a structure on a set $X$, see the introduction to section~\ref{sec: spectra}).

\begin{definition}\label{def: spectrumofposets}
A \textit{spectrum of posets}\index{spectrum of posets} over a set $I$ is an $I$-family of sets
$\Lambda := \big(\lambda_0, \lambda_1\big)$ such that 
$(\lambda_0(i), \leq_i)$ is a poset for every $i \in I$, and for every $(i, j) \in D(I)$ the transport
map $\lambda_{ij} :
\lambda_0(i) \to \lambda_0(j)$ is a monotone function. 
If $\D F^{\mon}\big(\lambda_0(i), \lambda_0(j)\big)$\index{$F^{\mon}(X, Y)$} is the set of monotone 
functions from $\lambda_0(i)$ to 
$\lambda_0(j)$, the category $\C C_{\Lambda}$ generated by the $I$-spectrum $\Lambda$ is the structure
$\C C_{\Lambda} := \big(\Ob^{\C C_{\Lambda}}, 
\Morr_0^{\C C_{\Lambda}}, \Morr_1^{\C C_{\Lambda}}, \Comp^{\C C_{\Lambda}}, \Id^{\C C_{\Lambda}}\big)$, where
$\Ob^{\C C_{\Lambda}} :=  \lambda_0 I$, and $\Morr_0^{\C C_{\Lambda}} := \D F^{\mon}\big(\lambda_0(i), \lambda_0(j)\big)$.
If $i =_I i{'}$ and $j =_I j{'}$, and since the composition of monotone functions is monotone,  let
$\Morr^{\C C_{\Lambda}}_{(i,j)(i{'}j{'})} : \D F^{\mon}\big(\lambda_0(i), \lambda_0(j)\big) \to 
\D F^{\mon}\big(\lambda_0(i{'}), \lambda_0(j{'})\big)$, defined by
\[ f \mapsto \Morr^{\C C_{\Lambda}}_{(i,j)(i{'},j{'})}(f), \ \ \ \ 
\Morr^{\C C_{\Lambda}}_{(i,j)(i{'},j{'})}(f) := \lambda_{jj{'}} \circ f \circ \lambda_{i{'}i}; \ \ \ \ f \in \D F^{\mon}\big(\lambda_0(i), \lambda_0(j)\big),\]
\begin{center}
\begin{tikzpicture}

\node (E) at (0,0) {$\lambda_0(i{'})$};
\node[right=of E] (H) {};
\node[right=of H] (F) {$\lambda_0(j{'})$.};
\node[above=of F] (A) {$\lambda_0(j)$};
\node [above=of E] (D) {$\lambda_0(i)$};

\draw[->] (E)--(F) node [midway,below]{$\Morr^{\C C_{\Lambda}}_{(i,j)(i{'}j{'})}(f)$};
\draw[->] (D)--(A) node [midway,above] {$f$};
\draw[->] (E)--(D) node [midway,left] {$\lambda_{i{'}i}$};
\draw[->] (A)--(F) node [midway,right] {$\lambda_{jj{'}}$};

\end{tikzpicture}
\end{center}
The dependent operations $\Id^{\C C_{\Lambda}}\big(\lambda_0(i)\big) := \Id_{\lambda_0(i)}$ and
$\Comp^{\C C_{\Lambda}}$ are defined as expected.
\end{definition}

Next we only show $(\Cat_3)$ and $(\Cat_4)$ for $\C C_{\Lambda}$. 
If $i =_I i{'}$, $f =_I j{'}$ and $k =_I k{'}$, we have that 
\begin{align*}
\Morr^{\C C_{\Lambda}}_{(k,i)(k{'}, i{'})}(\phi) \circ \Morr^{\C C_{\Lambda}}_{(j,k)(j{'}, k{'})}(\psi) & := 
\big(\lambda_{ii{'}} \circ \phi \circ \lambda_{k{'}k} \big) \circ \big(\lambda_{kk{'}} \circ \psi \circ
\lambda_{j{'}j}\big)\\
& = \lambda_{ii{'}} \circ \phi \circ \big( \lambda_{k{'}k} \circ \lambda_{kk{'}}\big) \circ \psi \circ
\lambda_{j{'}j}\\
& = \lambda_{ii{'}} \circ (\phi \circ \psi) \circ \lambda_{j{'}j}\\
& := \Morr^{\C C_{\Lambda}}_{(j,i)(j{'}, i{'})}(\phi \circ \psi),
\end{align*}
\[ \Morr^{\C C_{\Lambda}}_{(i,i)(i{'}, i{'})}\big(\Id^{\C C_{\Lambda}}\big(\lambda_0(i)\big)) 
:= \Morr^{\C C_{\Lambda}}_{(i,i)(i{'}, i{'})}\big(\Id_{\lambda_0(i)}\big) := \lambda_{ii{'}} \circ \Id_{\lambda_0(i)}
\circ \lambda_{i{'}i} = \Id_{\lambda_0(j)} := \Id^{\C C_{\Lambda}}\big(\lambda_0(j)\big).\]

\begin{definition}\label{def: functor}
A \textit{functor}\index{functor} $\B F : \C C \to \C D$ from $\C C := \big(\Ob_{\C C}, \Morr^{\C C}_0, 
\Morr^{\C C}_1, \Comp^{\C C}, \Id^{\C C}\big)$ to
$\C D := \big(\Ob_{\C D}, \Morr^{\C D}_0, \Morr^{\C D}_1, \Comp^{\C D}, \Id^{\C D}\big)$ is a pair
$\B F := (F_0, F_1)$, where 
$F_0 : \Ob_{\C C} \to \Ob_{\C D}$ and
\[ F_1 : \bigcurlywedge_{x, y \in \Ob_{\C C}} \D F\big(\Morr^{\C C}_0(x, y), \Morr^{\C D}_0(F_0(x), F_0(y))\big),\] 
\[ F_{xy} := F_1(x,y) : \Morr^{\C C}_0(x, y) \to \Morr^{\C D}_0(F_0(x), F_0(y)), \]
such that the following conditions are satisfied:\\[1mm]
$(\Funct_1)$ For every $x, y, z \in \Ob_{\C C}$, and for every $\psi \in \Morr^{\C C}_0(x, y)$,
$\phi \in \Morr^{\C C}_0(y, z)$ we have that
$$F_{xz}(\phi \circ \psi) =_{\mathsmaller{\Morr^{\C D}_0(F_0(x),F_0(z))}} F_{yz}(\phi) \circ F_{xy}(\psi).$$
$(\Funct_2)$ For every $x \in \Ob_{\C C}$ we have that
$F_{xx}\big(\Id^{\C C}_x\big) =_{\mathsmaller{\Morr^{\C D}_0(F_0(x),F_0(x))}} \Id^{\C D}_{F_0(x)}$.\\[1mm]
$(\Funct_3)$ For every $x, y, x{'}, y{'} \in \Ob_{\C C}$, such that $x =_{\mathsmaller{\Ob_{\C C}}} x{'}$ and
$y =_{\mathsmaller{\Ob_{\C C}}} y{'}$, hence $(x, y) = (x{'}, y{'})$ and $\big(F_0(x), F_0(y)\big) = \big(F_0(x{'}), 
F_0(y{'})\big)$, the following diagram commutes
\begin{center}
\begin{tikzpicture}

\node (E) at (0,0) {$\Morr^{\C D}_0(F_0(x), F_0(y))$};
\node[right=of E] (F) {};
\node[right=of F] (G) {};
\node[right=of G] (K) {};
\node[right=of K] (H) {$\Morr^{\C D}_0(F_0(x{'}), F_0(y{'}))$.};
\node[above=of E] (A) {$\Morr^{\C C}_0(x, y)$};
\node [above=of H] (B) {$\Morr^{\C C}_0(x{'}, y{'})$};

\draw[->] (E)--(H) node [midway,below]{$\Morr^{\C D}_{\big(F_0(x), F_0(y)\big)\big(F_0(x{'}), F_0(y{'})\big)}$};
\draw[->] (A)--(B) node [midway,above] {$\Morr^{\C C}_{(x, y)(x{'}, y{'})}$};
\draw[->] (A)--(E) node [midway,left] {$F_{xy}$};
\draw[->] (B)--(H) node [midway,right] {$F_{x{'}y{'}}$};

\end{tikzpicture}
\end{center}
\end{definition}

The last condition, which is not found in the standard definition of a functor, is a compatibility condition 
between the $F_1$-part of a functor $\B F : \C C \to \C D$, the transport maps $\Morr^{\C C}_1$ 
of $\C C$ and the transport 
maps $\Morr^{\C D}_1$ of $\C D$. As an example of a standard categorical construction in this framework, we formulate the
notion of slice category. If $\C C := \big(\Ob_{\C C}, \Morr^{\C C}_0, \Morr^{\C C}_1, 
\Comp^{\C C}, \Id^{\C C}\big)$ is a category and 
$x \in \Ob_{\C C}$, then $\Lambda^{x} := \big(\lambda_0^x, \lambda_1^x\big) \in \Fam(\Ob_{\C C})$, where 
\[ \lambda_0^x : \Ob_{\C C} \sto \D V_0, \ \ \ \ \lambda_0^x(y) :=  \Morr^{\C C}_0(y, x); \ \ \ \ y \in \Ob_{\C C},\]
\[\lambda_1^x : \bigcurlywedge_{(y, y{'}) \in D(\Ob_{\C C})}\D F\big(\Morr^{\C C}_0(y, x), \Morr^{\C C}_0(y{'}, x)\big),\]
\[\lambda^x_{yy{'}} := \lambda_1^x(y, y{'}) := \Morr^{\C C}_{(y,x)(y{'}, x)} : \Morr^{\C C}_0(y, x) \to \Morr^{\C C}_0(y{'}, 
x); \ \ \ \ (y,y{'}) \in \D (\Ob_{C C}).\]

Then we can prove the following fact.

\begin{proposition}\label{prp: slice2}
Let $\C C := \big(\Ob_{\C C}, \Morr^{\C C}_0, \Morr^{\C C}_1, \Comp^{\C C}, \Id^{\C C}\big)$ be a category and 
$x, z \in \Ob_{\C C}$. Let the structure $\C C/x := \big(\Ob_{\C C/x}, \Morr^{\C C/x}_0, \Morr^{\C C/x}_1, \Comp^{\C C/x}, 
\Id^{\C C/x}\big)$, where
\[ \Ob_{\C C/x} := \sum_{y \in \Ob_{\C C}}\lambda_0^x(y) := \sum_{y \in \Ob_{\C C}}\Morr^{\C C}_0(y, x),\]
\[\Morr^{\C C/x}_0\big((y, f), (z, g)\big) := \big\{h \in \Morr^{\C C}_0(y, z) \mid g \circ h = f\big\}.\]
If $(y, f) =_{\mathsmaller{\Ob_{\C C/x}}} (y{'}, f{'})$ and $(z, g) =_{\mathsmaller{\Ob_{\C C/x}}} (z{'}, g{'})$, the
function 
\[ \Morr^{\C C/x}_{\big((y, f), (z, g)\big), \big((y{'}, f{'}), (z{'}, g{'})\big)} \colon 
\Morr^{\C C/x}_0\big((y, f), (z, g)\big) 
\to \Morr^{\C C/x}_0\big((y{'}, f{'}), (z{'}, g{'})\big), \]
\[ h \mapsto \Morr^{\C C}_{(y, z)(y{'}, z{'})}(h); \ \ \ \ h \in \Morr^{\C C/x}_0\big((y, f), (z, g)\big), \]
is well-defined. If $\Comp^{\C C/x}$ 
is defined in the expected compositional way, and if $\Id^{\C C/x}\big((y, f)\big)
:= \Id^{\C C}(y)$, for every $(y, f) \in \Ob_{\C C/x}$, then $\C C/x$ is a category. Moreover, 
if $h \in \Morr_0^{\C C}(x,z)$, then $\B H := (H_0, H_1) \colon \C C/x \to \C C/z$, where
\[ H_0 \colon  \bigg(\sum_{y \in \Ob_{\C C}}\Morr^{\C C}_0(y, x)\bigg) \to \bigg(\sum_{y \in 
\Ob_{\C C}}\Morr^{\C C}_0(y, z)\bigg), \]
\[(y, f) \mapsto (y, h \circ f); \ \ \ \ (y, f) \in \sum_{y \in \Ob_{\C C}}\Morr^{\C C}_0(y, x), \]
\[ H_1\colon \bigcurlywedge_{(y,f), (y{'}, f{'}) \in \C C/x}\D F\bigg(\Morr_0^{\C C/x}\big((y,f)(y{'},f{'})\big), \Morr_0^{\C C/z}\big((y, h \circ f),(y{'}, h \circ f{'})\big)\bigg),\]
\[ H_{(y,f)(y{'},f{'})} \colon \Morr_0^{\C C/x}\big((y,f)(y{'},f{'})\big) \to  \Morr_0^{\C C/z}\big((y, h \circ f),(y{'}, h \circ f{'})\big), \]
\[  g \mapsto g; \  \ \ \ \ g \in \Morr_0^{\C C/x}((y,f)(y{'},f{'})).\]
\end{proposition}

\end{note}

\chapter{Families of subsets}
\label{chapter: familiesofsubsets}

We develop the basic theory of set-indexed families of subsets and of the corresponding
family-maps between them. In contrast to set-indexed families of sets, the properties of which 
are determined ``externally'' 
through their transport maps, the properties of a set-indexed family $\Lambda(X)$ of subsets of a given set $X$ 
are determined ``internally'' through the embeddings of the subsets of $\Lambda(X)$ to $X$. The interior union of 
$\Lambda(X)$ is the internal analogue to the $\sum$-set of a set-indexed family of sets $\Lambda$, 
and the intersection of $\Lambda(X)$ is the internal analogue to the $\prod$-set of $\Lambda$.
Families of sets over products, sets of subsets, and direct families of subsets are the internal analogue to the 
corresponding notions for families of sets. Set-indexed families of partial functions
and set-indexed families of complemented subsets, together with their corresponding family-maps, are studied.

\section{Set-indexed families of subsets}
\label{sec: famofsubsets}

Roughly speaking, a family of subsets of a set $X$ indexed by some set $I$ is an assignment routine 
$\lambda_0 : I \sto \C P(X)$ that behaves like a function i.e., if $i =_I j$, then $\lambda_0(i) =_{\C P(X)} 
\lambda_0 (j)$. The following definition is a formulation of this rough description that reveals the witnesses of the 
equality $\lambda_0(i) =_{\C P(X)} \lambda_0 (j)$. This is done ``internally'', through the embeddings of the subsets 
into $X$. The equality $\lambda_0(i) =_{\D V_0} \lambda_0 (j)$, which in the previous chapter
is defined ``externally'' through the transport maps, follows, and a family of subsets
is also a family of sets.

\begin{definition}\label{def: famofsubsets}
Let $X$ and $I$ be sets. A \textit{family of subsets}\index{family of subsets} of $X$ indexed by $I$,
or an $I$-\textit{family of subsets}\index{$I$-family of subsets} of $X$, is a triplet
$\Lambda(X) := (\lambda_0, \C E^X, \lambda_1)$, where
\index{$\Lambda(X)$}\index{$\C E^X$}
$\lambda_0 : I \sto \D V_0$,
\[ \C E^X : \bigcurlywedge_{i \in I}\D F\big(\lambda_0(i), X\big), \ \ \ \ \C E^X(i) := \C E_i^X; \ \ \ \ i \in I, \]
\[ \lambda_1 : \bigcurlywedge_{(i, j) \in D(I)}\D F\big(\lambda_0(i), \lambda_0(j)\big), \ \ \ \ 
\lambda_1(i, j) := \lambda_{ij}; \ \ \ \ (i, j) \in D(I), \]
such that the following conditions hold:\\[1mm]
\normalfont (a) 
\itshape For every $i \in I$, the function $\C E_i^X : \lambda_0(i) \to X$ is an embedding.\\[1mm]
\normalfont (b) 
\itshape For every $i \in I$, we have that $\lambda_{ii} := \id_{\lambda_0(i)}$.\\[1mm]
\normalfont (c) 
\itshape For every $(i, j) \in D(I)$ we have that 
$\C E_i^X = \C E_j^X \circ \lambda_{ij}$ and $\C E_j^X = \C E_i^X \circ \lambda_{ji}$
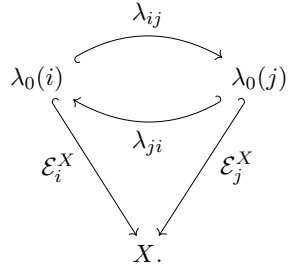
\begin{center}
\resizebox{4cm}{!}{%
\begin{tikzpicture}

\node (E) at (0,0) {$\lambda_0(i)$};
\node[right=of E] (B) {};
\node[right=of B] (F) {$\lambda_0(j)$};
\node[below=of B] (C) {};
\node[below=of C] (A) {$X$.};

\draw[left hook->,bend left] (E) to node [midway,above] {$\lambda_{ij}$} (F);
\draw[left hook->,bend left] (F) to node [midway,below] {$\lambda_{ji}$} (E);
\draw[right hook->] (E)--(A) node [midway,left] {$\C E_{i}^X \ $};
\draw[left hook->] (F)--(A) node [midway,right] {$ \ \C E_j^X$};

\end{tikzpicture}
}
\end{center} 
$\C E^X$ is a modulus of embeddings\index{modulus of embeddings} for $\lambda_0$, and $\lambda_1$ 
a modulus of transport maps for $\lambda_0$. Let 
$\Lambda := (\lambda_0, \lambda_1)$ be the $I$-family of sets 
that corresponds to $\Lambda(X)$\index{the $I$-family that corresponds to $\Lambda(X)$}.
If $(A, i_A) \in \C P(X)$, the \textit{constant} $I$-\textit{family of subsets}\index{constant family of subsets} 
$A$\index{$C^A_X$}
is the pair 
$C^{A}(X) := (\lambda_0^{A}, \C E^{X,A}, \lambda_1^A)$, where $\lambda_0 (i) := A$, $\C E_i^{X,A} := i_A$, and 
$\lambda_1 (i, j) := \id_A$, for every $i \in I$ and $(i, j) \in D(I)$ $($see the left diagram in
Definition~\ref{def: basicfamiliessub}$)$.
%
%
%
%
%
\end{definition}

\begin{proposition}\label{prp: famofsubsetsequiv}
Let $X$ and $I$ be sets, $\lambda_0 : I \sto \D V_0$, $\C E^X$ a modulus of embeddings for $\lambda_0$, and
$\lambda_1$ a modulus of transport maps for $\lambda_0$. The following are equivalent.\\[1mm]
\normalfont (i) 
\itshape $\Lambda(X) := (\lambda_0, \C E^X, \lambda_1)$ is an $I$-family of subsets of $X$.\\[1mm]
\normalfont (ii) 
\itshape $\Lambda := (\lambda_0, \lambda_1) \in \Fam(I)$ and $\C E^X \colon \Lambda \To C^X$, where
$C^X$ is the constant $I$-family $X$.
\end{proposition}

\begin{proof}
(i)$\To$(ii) First we show that $\Lambda \in \Fam(I)$. If $i =_I j =_I k$, 
then 
$\C E_k^X \circ (\lambda_{jk} \circ \lambda_{ij}) = (\C E_k^X \circ \lambda_{jk}) \circ \lambda_{ij}
= \C E_j^X \circ \lambda_{ij} = \C E_i^X$ and $\C E_k^X \circ \lambda_{ik} = \C E_i^X$ 
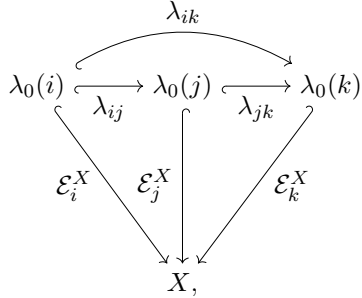
\begin{center}
\resizebox{5cm}{!}{%
\begin{tikzpicture}

\node (E) at (0,0) {$\lambda_0(i)$};
\node[right=of E] (B) {$\lambda_0(j)$};
\node[right=of B] (F) {$\lambda_0(k)$};
\node[below=of B] (C) {};
\node[below=of C] (A) {$X$,};

\draw[left hook->,bend left] (E) to node [midway,above] {$\lambda_{ik}$} (F);
\draw[left hook->] (E)--(B) node [midway,below] {$\lambda_{ij}$} (E);
\draw[left hook->] (B)--(F) node [midway,below] {$\lambda_{jk}$} (E);
\draw[right hook->] (E)--(A) node [midway,left] {$\C E_{i}^X \ $};
\draw[left hook->] (F)--(A) node [midway,right] {$ \ \C E_k^X$};
\draw[right hook->] (B)--(A) node [midway,left] {$\C E_{j}^X  $};

\end{tikzpicture}
}
\end{center} 
hence $\C E_k^X \circ (\lambda_{jk} \circ \lambda_{ij}) = \C E_k^X \circ \lambda_{ik}$, and since $\C E_k^X$ is 
an embedding, we get $\lambda_{jk} \circ \lambda_{ij} = \lambda_{ik}$. If $\C E^X \colon \Lambda \To C^X$, 
the following squares are commutative 
\begin{center}
\begin{tikzpicture}

\node (E) at (0,0) {$X$};
\node[right=of E] (F) {$X$};
\node[above=of F] (A) {$\lambda_0(j)$};
\node [above=of E] (D) {$\lambda_0(i)$};
\node [right=of F] (G) {$X$};
\node [above=of G] (H) {$\lambda_0(j)$};
\node [right=of G] (K) {$X$};
\node [above=of K] (L) {$\lambda_0(i)$};

\draw[->] (E)--(F) node [midway,below]{$\id_X$};
\draw[->] (D)--(A) node [midway,above] {$\lambda_{ij}$};
\draw[right hook->] (D)--(E) node [midway,left] {$\C E_i^X$};
\draw[left hook->] (A)--(F) node [midway,right] {$\C E_j^X$};
\draw[->] (G)--(K) node [midway,below] {$\id_X$};
\draw[->] (H)--(L) node [midway,above] {$\lambda_{ji}$};
\draw[right hook->] (H)--(G) node [midway,left] {$\C E_j^X$};
\draw[left hook->] (L)--(K) node [midway,right] {$\C E_i^X$};

\end{tikzpicture}
\end{center}
\begin{center}
\begin{tikzpicture}

\node (E) at (0,0) {};
\node[right=of E] (F) {$X$};
\node[above=of E] (A) {$\lambda_0(i)$};
\node [right=of F] (B) {};
\node [above=of B] (G) {$\lambda_0(j)$};
\node [right=of G] (K) {$\lambda_0(j)$};
\node [below=of K] (L) {};
\node [right=of L] (M) {$X$};
\node [above=of M] (N) {};
\node [right=of N] (P) {$\lambda_0(i)$};

\draw[right hook->] (A)--(F) node [midway,left] {$\C E_i^X \ \ $};
\draw[left hook->] (G)--(F) node [midway,right] {$\ \C E_j^X$};
\draw[right hook->] (A)--(G) node [midway,above] {$\lambda_{ij}$};
\draw[left hook->] (P)--(M) node [midway,right] {$ \ \C E_i^X$};
\draw[right hook->] (K)--(M) node [midway,left] {$\C E_j^X \ \ $};
\draw[right hook->] (K)--(P) node [midway,above] {$\lambda_{ji}$};

\end{tikzpicture}
\end{center}
if and only if the above triangles are commutative. The implication (ii)$\To$(i) follows immediately
from the equivalence between the commutativity of the above pairs of diagrams.
\end{proof}

\begin{definition}\label{def: basicfamiliessub}
Let $X$ be a set and $(A, i_A^X), (B, i_B^X) \subseteq X$.
The triplet \index{$\Lambda^{\D 2}(X)$}$\Lambda^{\D 2}(X) := (\lambda_0^{\D 2}, \C E^X, \lambda_1^{\D 2})$,
where $\Lambda^{\D 2} := \lambda_0^{\D 2}, \lambda_1^{\D 2}$ is the $\D 2$-family of $A, B$,
$\C E^X_0 := i_A^X$, and $\C E^X_1 := i_B^X$
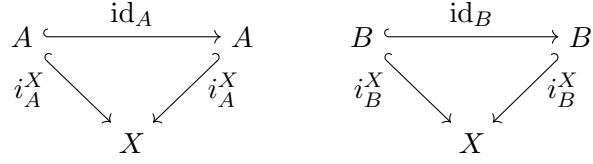
\begin{center}
\begin{tikzpicture}

\node (E) at (0,0) {};
\node[right=of E] (F) {$X$};
\node[above=of E] (A) {$A$};
\node [right=of F] (B) {};
\node [above=of B] (G) {$A$};
\node [right=of G] (K) {$B$};
\node [below=of K] (L) {};
\node [right=of L] (M) {$X$};
\node [above=of M] (N) {};
\node [right=of N] (P) {$B$};

\draw[right hook->] (A)--(F) node [midway,left] {$i_A^X \ \ $};
\draw[left hook->] (G)--(F) node [midway,right] {$\ i_A^X$};
\draw[right hook->] (A)--(G) node [midway,above] {$\id_A$};
\draw[left hook->] (P)--(M) node [midway,right] {$ \ i_B^X$};
\draw[right hook->] (K)--(M) node [midway,left] {$i_B^X \ \ $};
\draw[right hook->] (K)--(P) node [midway,above] {$\id_B$};

\end{tikzpicture}
\end{center}
is the $\D 2$-family of subsets $A$ and $B$\index{$\D 2$-family of subsets} of $X$.
The \index{$\Lambda^{\D n}(X)$}$\D n$-family $\Lambda^{\D n}(X)$\index{$\Lambda^{\D n}$} of the subsets 
$(A_1, i_1), \ldots, (A_n, i_n)$ of $X$\index{$\D n$-family of subsets}, and 
the $\Nat$-family of subsets $(A_n, i_n)_{n \in \Nat}$\index{$\Nat$-family of subsets} of $X$ 
are defined similarly.
\end{definition}

\begin{definition}\label{def: subfammap}
If $\Lambda(X) := (\lambda_0, \C E^X, \lambda_1), M(X) := (\mu_0, \C Z^X, \mu_1)$ and 
$N(X) := (\nu_0, \C H^X, \nu_1)$ are $I$-families of subsets of $X$, a 
\textit{family of subsets-map}\index{family of subsets-map}
$\Psi \colon \Lambda(X) \To M(X)$\index{$\Psi \colon \Lambda(X) \To M(X)$}
from $\Lambda(X)$ to $M(X)$ is a dependent operation
$\Psi : \bigcurlywedge_{i \in I}\D F\big(\lambda_0(i), \mu_0(i)\big)$, where $\Psi(i) := \Psi_i$,
for every $i \in I$, such that, for every $i \in I$, the following diagram commutes\footnote{Trivially, for 
every $i \in I$ the map $\Psi_i \colon \lambda_0(i) \to \mu_0(i)$ is an embedding.}
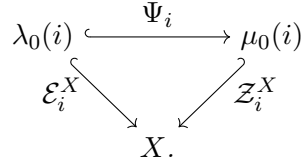
\begin{center}
\begin{tikzpicture}

\node (E) at (0,0) {};
\node[right=of E] (F) {$X$.};
\node[above=of E] (A) {$\lambda_0 (i)$};
\node [right=of F] (B) {};
\node [above=of B] (G) {$\mu_0 (i)$};

\draw[right hook->] (A)--(F) node [midway,left] {$\C E_{i}^X \ $};
\draw[left hook->] (G)--(F) node [midway,right] {$\ \C Z_i^X$};
\draw[right hook->] (A)--(G) node [midway,above] {$\Psi_{i}$};

\end{tikzpicture}
\end{center}
The totality $\Map_I(\Lambda(X), M(X))$\index{$\Map_I(\Lambda(X), M(X))$} of family of subsets-maps 
from $\Lambda(X)$ to $M(X)$ is
equipped with the pointwise equality. 
If $\Psi \colon \Lambda(X) \To M(X)$ and $\Xi \colon M(X) \To N(X)$, the \textit{composition family of subsets-map} 
\index{composition family of subsets-map} 
$\Xi \circ \Psi \colon \Lambda(X) \To N(X)$ is defined by 
$(\Xi \circ \Psi)(i) := \Xi_i \circ \Psi_i$,
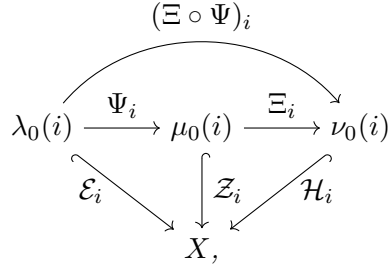
\begin{center}
\begin{tikzpicture}

\node (E) at (0,0) {$\lambda_0(i)$};
\node[right=of E] (B) {$\mu_0(i)$};
\node[right=of B] (F) {$\nu_0(i)$};
\node[below=of B] (C) {$X$,};

\draw[left hook->] (E)--(C) node [midway,left] {$\C E_i \ $};
\draw[right hook->] (B)--(C) node [midway,right] {$\C Z_i$};
\draw[right hook->] (F)--(C) node [midway,right] {$ \ \C  H_{i} $};
\draw[->] (E)--(B) node [midway,above] {$\Psi_i$};
\draw[->] (B)--(F) node [midway,above] {$\Xi_i$};
\draw[->,bend left=50] (E) to node [midway,above] {$(\Xi \circ \Psi)_{i}$} (F);

\end{tikzpicture}
\end{center} 
for every $i \in I$. The identity family of subsets-map $\Id_{\Lambda(X)} \colon \Lambda(X) \To \Lambda(X)$ 
and the equality on the totality $\Fam(I,X)$\index{$\Fam(I,X)$} 
of $I$-families of subsets of $X$ are defined as in Definition~\ref{def: map}.
\end{definition}

We see no obvious reason, like the one for $\Fam(I)$, not to consider 
$\Fam(I, X)$ to be a set. In the case of $\Fam(I)$ the constant $I$-family $\Fam(I)$ would be in $\Fam(I)$, while
the constant $I$-family $\Fam(I, X)$ is not clear how could be seen as a family of subsets of $X$. If 
$\nu_0(i) := \Fam(I, X)$, for every $i \in I$, we need to define a modulus of embeddings $\C N^X_i \colon \Fam(I, X) 
\eto X$, for every $i \in I$. From the given data one could define the assignment routine $\C N^X_i$ by the rule 
$\C N^X_i\big(\Lambda(X)\big) := \C E_i^X(u_i)$, if it is known that $u_i \in \lambda_0(i)$. Even in that case, the 
assignment routine $\C N^X_i$ cannot be shown to satisfy the expected properties. Clearly, if $\C N^X_i$ was defined 
by the rule $\C N^X_i\big(\Lambda(X)\big) := x_0 \in X$, then it cannot be an embedding.

\begin{definition}\label{def: prfeqfamofsubsets}
If $\Lambda(X) M(X) \in \Fam(I, X)$, let 
\[\Lambda(X) \leq M(X) :\TOT \exists_{\Phi \in \Map_I(\Lambda(X),
M(X))}\big(\Phi \colon \Lambda(X) \To M(X)\big), \] 
If $\Phi \in \Map_I(\Lambda(X), M(X)), \Psi \in \Map_I(M(X), \Lambda(X)),
\Phi{'} \in \Map_I(M(X), N(X))$ and $\Psi{'} \in \Map_I(N(X), M(X))$, let the following set and operations
\[ \Eq(\Lambda(X), M(X)) := \Map_I(\Lambda(X), M(X)) \times \Map_I(M(X), \Lambda(X)), \]
\[ \refl(\Lambda) := \big(\Id_{\Lambda_X}, \Id_{\Lambda_X}\big) \ \ \& \ \ (\Phi, \Psi)^{-1} := (\Psi, \Phi) \ \ \& \ \ 
(\Phi, \Psi) \ast (\Phi{'}, \Psi{'}) := (\Phi{'} \circ \Phi, \Psi \circ \Psi{'}). \]
\end{definition}

%

\begin{proposition}\label{prp: famsubmap1}
 Let $\Lambda(X) := (\lambda_0, \C E^X, \lambda_1), M(X) := (\mu_0, \C Z^X, \mu_1) \in \Fam(I, X)$.\\[1mm]
\normalfont (i) 
\itshape If $\Psi \colon \Lambda(X) \To M(X)$, then $\Psi \colon \Lambda \To M$.\\[1mm]
\normalfont (ii) 
\itshape If $\Psi \colon \Lambda(X) \To M(X)$ and $\Phi \colon \Lambda(X) \To M(X)$, then 
$\Phi =_{\Map_I(\Lambda(X), M(X))} \Psi$.
\end{proposition}

\begin{proof}
(i) By the commutativity of the following inner diagrams 
\begin{center}
\begin{tikzpicture}
  \matrix (m) [matrix of math nodes,row sep=3em,column sep=3em,minimum width=2em]
  {
     \lambda_0 (i) & {} & \mu_0 (i)\\
     {} & X & {}\\
     \lambda_0 (j) & {} & \mu_0 (j), \\};
\path[left hook->]
    (m-1-1) edge node [left] {$\lambda_{ij}$} (m-3-1)
    (m-1-1)   edge node [above] {$\Psi_i$} (m-1-3);
\path[right hook->]
    (m-1-1) edge node [right] {$ \ \C E_i^X$}  (m-2-2);
\path[left hook->]
    (m-1-3) edge node [left] {$\C Z_i^X \ $}  (m-2-2);
\path[right hook->]
    (m-3-1) edge node [below] {$\Psi_j$} (m-3-3);
\path[right hook->]
    (m-3-1) edge node [right] {$ \ \C E_j^X$} (m-2-2);
\path[left hook->]
    (m-3-3) edge node [left] {$\C Z_j^X \ \ $} (m-2-2);
\path[right hook->]
    (m-1-3) edge node [right] {$\mu_{ij}$} (m-3-3);
\end{tikzpicture}
\end{center}
we get the required commutativity of the above outer diagram. 
If $x \in \lambda_0(i)$, then
\[ (\C Z_j^X \circ \Psi_j)\big(\lambda_{ij}(x)\big)
= \C E_j^X \big(\lambda_{ij}(x)\big) = \C E_i^X (x) = (\C Z_i^X \circ \Psi_i)(x) = \C Z_j^X
\big(\mu_{ij}(\Psi_i (x))\big). \]
Since $\C Z_j^X \big(\Psi_j (\lambda_{ij}(x))\big) =  \C Z_j^X
\big(\mu_{ij}(\Psi_i (x))\big)$, we get 
$\Psi_j (\lambda_{ij}(x) = \mu_{ij}(\Psi_i (x))$.\\
(ii) If $i \in I$, then $\Psi_i \colon \lambda_0(i) \subseteq \mu_0(i)$, $\Phi_i \colon \lambda_0(i) \subseteq \mu_0(i)$
\begin{center}
\resizebox{4cm}{!}{%
\begin{tikzpicture}

\node (E) at (0,0) {$\lambda_0(i)$};
\node[right=of E] (B) {};
\node[right=of B] (F) {$\mu_0(i)$};
\node[below=of B] (C) {};
\node[below=of C] (A) {$X$,};

\draw[left hook->,bend left] (E) to node [midway,above] {$\Phi_i$} (F);
\draw[left hook->,bend right] (E) to node [midway,below] {$\Psi_i$} (F);
\draw[right hook->] (E)--(A) node [midway,left] {$\C E_i \ $};
\draw[left hook->] (F)--(A) node [midway,right] {$ \ E_i$};

\end{tikzpicture}
}
\end{center}
hence by Proposition~\ref{prp: subset1} we get $\Psi_i =_{\D F(\lambda_0(i), \mu_0(i))} \Phi_i$. 
\end{proof}

Because of Proposition~\ref{prp: famsubmap1}(ii) all the elements of 
$\Eq(\Lambda(X), M(X))$ are equal to each other, hence the groupoid- properties (i)-(iv) for $\Eq(\Lambda(X), M(X))$
hold trivially. Of course, $\Lambda(X) =_{\Fam(I, X)} M(X) :\TOT \Lambda(X) \leq M(X) \ \& \ M(X) \leq \Lambda(X)$.
The characterisation of a family of subsets given in Proposition~\ref{prp: famofsubsetsequiv} together with the operations
on family-maps help us define new families of subsets from given ones.

\begin{proposition}\label{prp: newfamsofsubsets1}
If $\Lambda(X) := (\lambda_0, \C E^X, \lambda_1) \in \Fam(I, X)$ and $M(Y) := (\mu_0, \C E^Y, \mu_1) \in \Fam(I, Y)$,
then 
\[ (\Lambda \times M)(X \times Y) := \Lambda(X) \times M(Y) := \big(\lambda_0 \times \mu_0, \C E^X \times \C E^Y, 
 \lambda_1, \mu_1\big) \in \Fam(I, X \times Y),
\]
where the $I$-family $\Lambda \times M := (\lambda_0 \times \mu_0, \lambda_1 \times \mu_1)$ is defined in 
Definition~\ref{def: newfamsofsets1} and the family-map $\C E^X \times \C E^Y 
\colon \Lambda \times M \To C^X \times C^Y$ is defined in Proposition~\ref{prp: newfamilymaps1}\normalfont(ii).
\end{proposition}

\begin{proof}
By Proposition~\ref{prp: newfamilymaps1}(ii) $\C E^X \times \C E^Y 
\colon \Lambda \times M \To C^X \times C^Y$, where $(\C E^X \times \C E^Y)_i : \lambda_0(i) \times \mu_0(i) 
\hookrightarrow X \times Y$ is defined by the rule $(u, w) \mapsto \big(\C E^X(u), \C E^Y(w)\big)$, for every
$(u, w) \in \lambda_0(i) \times \mu_0(i)$. By Proposition~\ref{prp: subset15} 
this is a well-defined subset of $X \times Y$. 
By Proposition~\ref{prp: constequal7}(i) 
$\C E^X \times \C E^Y \colon \Lambda \times M \To C^{X \times Y}$, and we use Proposition~\ref{prp: famofsubsetsequiv}.
\end{proof}

The operations on subsets induce operations on families of subsets.

\begin{proposition}\label{prp: newfamsofsubsets1.5}
 Let $\Lambda(X) := (\lambda_0, \C E^X, \lambda_1)$ and $M(X) := (\mu_0, \C Z^X, \mu_1) \in \Fam(I, X)$.\\[1mm]
\normalfont (i) 
\itshape $(\Lambda \cap M)(X) := (\lambda_0 \cap \mu_0, \C E^X \cap \C Z^X, \lambda_1 \cap \mu_1) 
\in \Fam(I, X)$,\index{$(\Lambda \cap M)(X)$}
where $\lambda_0 \cap \mu_0 \colon I \sto \D V_0$ is defined by $(\lambda_0 \cap \mu_0)(i) := \lambda_0(i) \cap \mu_0(i)$,
for every $i \in I$, and the dependent operations $\C E^X \cap \C Z^X \colon 
\bigcurlywedge_{i \in I}\D F\big(\lambda_0(i) \cap \mu_0(i), X\big)$,  
$\lambda_1 \cap \mu_1 \colon \bigcurlywedge_{(i,j) \in D(I)}\D F\big(\lambda_0(i) \cap \mu_0(i),
\lambda_0(j) \cap \mu_0(j)\big)$ are defined by
\[ \big(\C E^X \cap \C Z^X\big)_i := i_{\lambda_0(i) \cap \mu_0(i)}^X; \ \ \ \ i \in I, \]
\[ \big[\big(\lambda_1 \cap \mu_1\big)_1(i, j)\big](u,w) := \big(\lambda_1 \cap \mu_1\big)_{ij}(u,w) 
:= \big(\lambda_{ij}(u),
\mu_{ij}(w)\big); \ \ \ \ (u,w) \in \lambda_0(i) \cap \mu_0(i).
\]
\normalfont (ii) 
\itshape $(\Lambda \cup M)(X) := (\lambda_0 \cup \mu_0, \C E^X \cup \C Z^X, \lambda_1 \cup \mu_1) 
\in \Fam(I, X)$,\index{$(\Lambda \cup M)(X)$} 
where $\lambda_0 \cup \mu_0 \colon I \sto \D V_0$ is defined by $(\lambda_0 \cup \mu_0)(i) := \lambda_0(i) \cup \mu_0(i)$,
for every $i \in I$, and the dependent operations $\C E^X \cup \C Z^X \colon 
\bigcurlywedge_{i \in I}\D F\big(\lambda_0(i) \cup \mu_0(i), X\big)$,  
$\lambda_1 \cup \mu_1 \colon \bigcurlywedge_{(i,j) \in D(I)}\D F\big(\lambda_0(i) \cup \mu_0(i),
\lambda_0(j) \cup \mu_0(j)\big)$ are defined by
\[ \big(\C E^X \cup \C Z^X\big)_i(z) := \left\{ \begin{array}{ll}
                 \C E_i^X(z)    &\mbox{, $z \in \lambda_0(i)$}\\
                 \C Z_i^X(z)             &\mbox{, $z \in \mu_0(i)$,}
                 \end{array}
          \right.
          ; \ \ \ \ i \in I, \ z \in \lambda_0(i) \cup \mu_0(i) \]
\[  (\lambda_1 \cup \mu_1)_{ij}(z) := \left\{ \begin{array}{ll}
                 \lambda_{ij}(z)    &\mbox{, $z \in \lambda_0(i)$}\\
                 \mu_{ij}(z)             &\mbox{, $z \in \mu_0(i)$,}
                 \end{array}
          \right.
          ; \ \ \ \ (i, j) \in D(I), \ z \in \lambda_0(i) \cup \mu_0(i). \]
\end{proposition}

\begin{proof}
(i) By Definition~\ref{def: image} we have that
\[  \lambda_0(i) \cap \mu_0(i) := \{(u, w) \in \lambda_0(i) \times \mu_0(i) \mid \C E_i^X(u) =_X \C Z_i^X(w)\}, \ \ \ \ 
i_{\lambda_0(i) \cap \mu_0(i)}^X(u, w) := \C E_i^X(u), \]
\[ \lambda_0(j) \cap \mu_0(j) := \{(u{'}, w{'}) \in \lambda_0(j) \times \mu_0(j) \mid \C E_j^X(u{'}) =_X \C Z_j^X(w{'})\}, \ 
i_{\lambda_0(j) \cap \mu_0(j)}^X(u{'}, w{'}) := \C E_j^X(u{'}). \]
Since 
$\C E_j^X\big(\lambda_{ij}(u)\big) =_X \C E_i^X(u) =_X \C Z_i^X(w) =_X \C Z_j^X\big(\mu_{ij}(w)\big)$,
we get $\big(\lambda_1 \cap \mu_1\big)_{ij}(u,w) \in \lambda_0(j) \cap \mu_0(j)$.
Clearly, $\big(\lambda_1 \cap \mu_1\big)_{ij}$ is a function. The commutativity of the following left inner diagrams
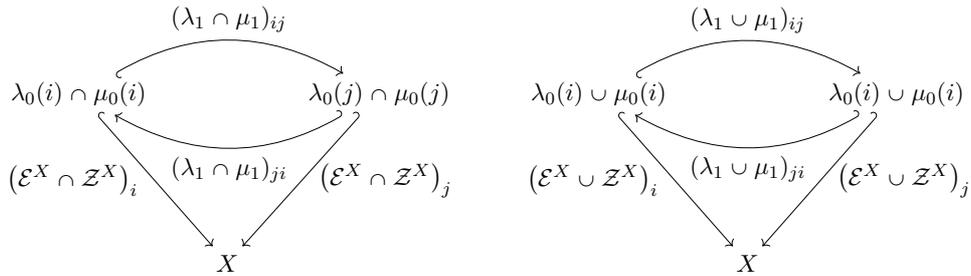
\begin{center}
\resizebox{13cm}{!}{%
\begin{tikzpicture}

\node (E) at (0,0) {$\lambda_0(i) \cap \mu_0(i)$};
\node[right=of E] (B) {};
\node[right=of B] (F) {$\lambda_0(j) \cap \mu_0(j)$};
\node[below=of B] (C) {};
\node[below=of C] (A) {$X$};

\node[right=of F] (G) {$\lambda_0(i) \cup \mu_0(i)$};
\node[right=of G] (H) {};
\node[right=of H] (J) {$\lambda_0(i) \cup \mu_0(i)$};
\node[below=of H] (K) {};
\node[below=of K] (L) {$X$};

\draw[left hook->,bend left] (E) to node [midway,above] {$(\lambda_1 \cap \mu_1)_{ij}$} (F);
\draw[left hook->,bend left] (F) to node [midway,below] {$(\lambda_1 \cap \mu_1)_{ji}$} (E);
\draw[right hook->] (E)--(A) node [midway,left] {$\big(\C E^X \cap \C Z^X\big)_i \ $};
\draw[left hook->] (F)--(A) node [midway,right] {$ \ \big(\C E^X \cap \C Z^X\big)_j$};
\draw[left hook->,bend left] (G) to node [midway,above] {$(\lambda_1 \cup \mu_1)_{ij}$} (J);
\draw[left hook->,bend left] (J) to node [midway,below] {$(\lambda_1 \cup \mu_1)_{ji}$} (G);
\draw[right hook->] (G)--(L) node [midway,left] {$\big(\C E^X \cup \C Z^X\big)_i \ $};
\draw[left hook->] (J)--(L) node [midway,right] {$ \ \big(\C E^X \cup \C Z^X\big)_j$};

\end{tikzpicture}
}
\end{center} 
follows by the equalities 
$\big(\C E^X \cap \C X^X\big)_j\big((\lambda_1 \cap \mu_1)_{ij}(u,w)\big) := 
 i_{\lambda_0(i) \cap \mu_0(i)}^X\big(\lambda_{ij}(u), \mu_{ij}(w)\big) := \C E_j^X\big(\lambda_{ij}(u)\big) 
 =_X \C E_i^X(u) := \big(\C E^X \cap \C X^X\big)_i(u,w)$.\\
(ii) First we show that $(\lambda_1 \cup \mu_1)_{ij}$ is a function. The more interesting case is $z \in \lambda_0(i)$,
$w \in \mu_0(i)$ and $\C E_i^X(z) =_X \C Z_i^X(w)$. Hence $\C E_j^X\big(\lambda_{ij}(u)\big) =_X =_X 
\C Z_j^X\big(\mu_{ij}(w)\big)$, and $\lambda_{ij}(z) =_{\mathsmaller{\lambda_0(i) \cup \mu_0(i)}} \mu_{ij}(w)$.
The commutativity of the above right inner diagrams is straightforward to show.
\end{proof}

\begin{proposition}\label{prp: newfamsofsubsets2}
Let $\Lambda(X) := (\lambda_0, \C E^X, \lambda_1) \in \Fam(I, X)$ and $M(Y) := (\mu_0, \C E^Y, \mu_1) \in \Fam(J, Y)$. 
If $\fXY$, let
$[f(\Lambda)](Y) := \big(f(\lambda_0), f(\C E^X)^Y, f(\lambda_1)\big)$\index{[$f(\Lambda)](Y)$}, where 
the non-dependent assignment routine $f(\lambda_0) \colon I \sto \D V_0$, and the 
dependent operations $f(\C E^X)^Y \colon \bigcurlywedge_{i \in I}
\D F(f(\lambda_0)(i), Y)$ and $f(\lambda_1) \colon \bigcurlywedge_{(i, i{'}) \in D(I)}\D F\big(f(\lambda_0)(i),
f(\lambda_0)(i{'})\big)$ are defined by
\begin{center}
\begin{tikzpicture}

\node (E) at (0,0) {$\lambda_0(i)$};
\node[right=of E] (F) {$X$};
\node[right=of F] (A) {$Y$};

\draw[right hook->] (E)--(F) node [midway,above] {$\C E_i^X$};
\draw[->] (F)--(A) node [midway,above] {$f$};
\draw[->,bend right] (E) to node [midway,below] {$f_i^Y$} (A);

\end{tikzpicture}
\end{center}
\[ [f(\lambda_0)](i) := f(\lambda_0(i)) := (\lambda_0(i), f_i), \ \ \ \ f_i^Y := f \circ \C E_i^X; \ \ \ \ i \in I, \]
\[ \big[f(\C E^X)^Y\big](i) := f_i^Y, \ \ \ \  f(\lambda_1)_{ii{'}} := \lambda_{ii{'}};  \ \ \ \ i \in I, \ (i, i{'}) \in D(I). \]
We call $[f(\Lambda)](Y)$ the image of $\Lambda$\index{image of a family of subsets} under $f$.
The pre-image of $M$\index{pre-image of a family of subsets} under $f$ is the triplet 
$[f^{-1}(M)](X) := \big(f^{-1}(\mu_0), f^{-1}(\C E^Y)^X, f^{-1}(\mu_1)\big)$\index{[$f^{-1}(M)](X)$}, where 
the non-dependent assignment routine $f^{-1}(\mu_0) \colon J \sto \D V_0$, and the 
dependent operations $f^{-1}(\C E^Y)^X \colon \bigcurlywedge_{j \in J}
\D F(f^{-1}(\mu_0)(j), X)$ and $f^{-1}(\mu_1) \colon \bigcurlywedge_{(j, j{'}) \in D(J)}\D F\big(f^{-1}(\mu_0)(j),
f^{-1}(\mu_0)(j{'})\big)$ are defined by
\[ [f^{-1}(\mu_0)](j) := f^{-1}(\mu_0(j)) := \big\{(x, y) \in X \times \mu_0(j) \mid f(x) =_Y \C E_j^Y(y)\big\};
\ \ \ \ j \in J \]
\[ e_j \colon f^{-1}(\mu_0(j)) \eto X \ \ \ \ e_j(x, y) := x; \ \ \ \ x \in X, \ y \in \mu_0(j), \ j \in 
J, \] 
\[ \big[f^{-1}(\C E^Y)^X\big](j) := e_j; \ \ \ \ j \in J, \]
\[ f^{-1}(\mu_1)_{jj{'}} \colon f^{-1}(\mu_0)(j) \to f^{-1}(\mu_0)(j{'}) \ \ \ \ f^{-1}(\mu_1)_{jj{'}}(x, y) :=
(x, \mu_{jj{'}}(y)); \ \ \ \ (j, j{'}) \in D(J). \]
Then $[f(\Lambda)](Y) \in \Fam(I, Y)$ and $[f^{-1}(M)](X) \in \Fam(J, X)$.
\end{proposition}

\begin{proof}
It suffices to show the commutativity of the following diagrams
\begin{center}
\resizebox{10cm}{!}{%
\begin{tikzpicture}

\node (E) at (0,0) {$\lambda_0(i)$};
\node[right=of E] (B) {};
\node[right=of B] (F) {$\lambda_0(i{'})$};
\node[below=of B] (C) {};
\node[below=of C] (A) {$X$};
\node[left=of A] (D) {$Y$};

\node[right=of F] (G) {$f^{-1}(\mu_0)(j)$};
\node[right=of G] (H) {};
\node[right=of H] (J) {$f^{-1}(\mu_0)(j{'})$};
\node[below=of H] (K) {};
\node[below=of K] (L) {$X$.};

\draw[left hook->,bend left] (E) to node [midway,above] {$\lambda_{ii{'}}$} (F);
\draw[left hook->,bend left] (F) to node [midway,below] {$\lambda_{i{'}i}$} (E);
\draw[right hook->] (E)--(A) node [midway,left] {$\C E_i^X \ $};
\draw[left hook->] (F)--(A) node [midway,right] {$ \ \C E^X_{i{'}}$};
\draw[->] (A)--(D) node [midway,above] {$f$};
\draw[left hook->,bend left] (G) to node [midway,above] {$f^{-1}(\mu_1)_{jj{'}}$} (J);
\draw[left hook->,bend left] (J) to node [midway,below] {$f^{-1}(\mu_1)_{j{'}j}$} (G);
\draw[right hook->] (G)--(L) node [midway,left] {$\mathsmaller{\big[f^{-1}(\C E^Y)^X\big](j)} \ $};
\draw[left hook->] (J)--(L) node [midway,right] {$ \ \mathsmaller{\big[f^{-1}(\C E^Y)^X\big](j{'})}$};

\end{tikzpicture}
}
\end{center} 
For the left, we use the supposed commutativity of the two diagrams without the arrow $\fXY$. For the above right
outer diagram we have that $\big[f^{-1}(\C E^Y)^X\big](j{'})\big(f^{-1}(\mu_1)_{jj{'}}(x, y)\big) := 
\big[f^{-1}(\C E^Y)^X\big](j{'})(x, \mu_{jj{'}}(y)) := e_{j{'}}((x, \mu_{jj{'}}(y)) := x := e_j(x, y) :=
\big[f^{-1}(\C E^Y)^X\big](j)(x, y)$. For the commutativity of the above right inner diagram we proceed similarly.  
\end{proof}

The operations on families of subsets generate operations on family of subsets-maps.

\begin{proposition}\label{prp: operationsfamsofsubsetsmaps}
 Let $\Lambda(X), K(X), M(X), N(X) \in \Fam(I, X)$, $P(Y), Q(Y) \in \Fam(J, Y)$, and $\fXY$.
 Let also $\Phi \colon \Lambda(X) \To K(X)$, $\Psi \colon M(X) \To N(X)$, and $\Xi \colon P(Y) \To Q(Y)$.\\[1mm]
 \normalfont (i) 
\itshape $\Phi \cap \Psi \colon (\Lambda \cap M)(X) \To (K \cap N)(X)$, where, for every $i \in I$ and 
$(u,w) \in \lambda_0(i) \cap \mu_0(i)$, 
\[ (\Phi \cap \Psi)_i \colon \lambda_0(i) \cap \mu_0(i) \to k_0(i) \cap \nu_0(i), \ \ \ \ 
 (\Phi \cap \Psi)_i(u,w) := \big(\Phi_i(u), \Psi_i(w)\big). \]
 \normalfont (ii) 
\itshape  $\Phi \cup \Psi \colon (\Lambda \cup M)(X) \To (K \cup N)(X)$, where, for every $i \in I$,
\[ (\Phi \cup \Psi)_i \colon \lambda_0(i) \cup \mu_0(i) \to k_0(i) \cup \nu_0(i), \]
\[  (\Phi \cup \Psi)_i(z) := \left\{ \begin{array}{ll}
                 \Phi_i(z)    &\mbox{, $z \in \lambda_0(i)$}\\
                 \Psi_i(z)             &\mbox{, $z \in \mu_0(i)$}
                 \end{array}
          \right.
\]
\normalfont (iii) 
\itshape $\Phi \times \Xi \colon (\Lambda \times P)(X \times Y) \To (K \times Q)(X \times Y)$, where, 
for every $i \in I$ and 
$(u,w) \in \lambda_0(i) \times p_0(i)$, 
\[ (\Phi \times \Xi)_i \colon \lambda_0(i) \times p_0(i) \to k_0(i) \times q_0(i), \ \ \ \ 
 (\Phi \times \Xi)_i(u,w) := \big(\Phi_i(u), \Xi_i(w)\big). \]
\normalfont (iv) 
\itshape $f(\Phi) \colon [f(\Lambda)](Y) \To [f(K)](Y)$, where, 
for every $i \in I$ and $u \in f(\lambda_0(i))$, 
\[ [f(\Phi)]_i \colon f(\lambda_0(i)) \to f(k_0(i)), \ \ \ \ 
  [f(\Phi)]_i(u) := \Phi_i(u). \]
\normalfont (v) 
\itshape $f^{-1}(\Xi) \colon [f^{-1}(P)](X) \To [f^{-1}(Q)](X)$, where, 
for every $j \in J$ and $(x,y) \in f^{-1}(p_0(j))$, 
\[ [f^{-1}(\Xi)]_i \colon f^{-1}(p_0(j)) \to f^{-1}(q_0(j)), \ \ \ \ 
  [f^{-1}(\Xi)]_j(x, y) := \big(x, \Xi_j(y)\big). \] 
\end{proposition}

\begin{proof}
It is straightforward to show that all family of subsets-maps above are well-defined.  
\end{proof}

\begin{definition}\label{def: subfamsubsets}
Let $\Lambda(X) := (\lambda_0, \C E^X, \lambda_1) \in \Fam(I, X)$ and $h \colon J \to I$. The triplet
$\Lambda(X) \circ h := (\lambda_0 \circ h, \C E^X \circ h, \lambda_1 \circ h)$, where $\Lambda \circ h := 
(\lambda_0 \circ h, \lambda_1 \circ h)$ is the $h$-subfamily of $\Lambda$, and the dependent operation
$\C E^X \circ h \colon \bigcurlywedge_{j \in J}\D F\big(\lambda_0(h(j)), X\big)$ is defined by 
$(\C E^X \circ h)_j := \C E^X_{h(j)}$, for every $j \in J$, is called the $h$-subfamily of $\Lambda(X)$\index{$h$-subfamily of a family of subsets}. If $J := \Nat$, we call $\Lambda(X) \circ h$ the $h$-subsequence of $\Lambda(X)$\index{$h$-subsequence of $\Lambda(X)$}.

\end{definition}

It is immediate to show that $\Lambda(X) \circ h \in \Fam(J, X)$, and if $\Lambda(X) \circ h \in \Set(J, X)$,
then $h$ is an embedding. All notions and results of section~\ref{sec: subfam} on
subfamilies of families of sets extend naturally to subfamilies of families of subsets.

\section{The interior union of a family of subsets}
\label{sec: union}

\begin{definition}\label{def: interiorunion} 
Let $\Lambda(X) := (\lambda_0, \C E^X, \lambda_1)$ be an $I$-family of subsets of $X$. 
The \textit{interior union}\index{interior union}, or simply the union\index{union of a family of subsets}
of $\Lambda(X)$ is the totality $\sum_{i \in I} \lambda_0(i)$, which we 
denote in this case by $\bigcup_{i \in I} \lambda_0(i)$. 
Let the non-dependent assignment routine
$e_{\mathsmaller{\bigcup}}^{\Lambda(X)} \colon \bigcup_{i \in I} \lambda_0(i) \sto X$ defined by
$(i, x) \mapsto \C E_i^X (x)$, for every $(i, x) \in \bigcup_{i \in I} \lambda_0(i)$, and let 
\[ (i, x) =_{\mathsmaller{\bigcup_{i \in I} \lambda_0(i)}} (j, y) :\Leftrightarrow 
e_{\mathsmaller{\bigcup}}^{\Lambda(X)}(i, x) =_X 
e_{\mathsmaller{\bigcup}}^{\Lambda(X)}(j, y) :\Leftrightarrow \C E_i^X (x) =_X \C E_j^X (y). \]
If $\neq_X$ is an inequality on $X$, let 
$(i, x) \neq_{\mathsmaller{\bigcup_{i \in I} \lambda_0(i)}} (j, y) :\Leftrightarrow 
\C E_i^X (x) \neq_X \C E_j^X (y)$.
The family $\Lambda(X)$ is called a \textit{covering}\index{covering} of $X$, or 
$\Lambda(X)$ covers $X$\index{$\Lambda(X)$ covers $X$},  if
\[ \bigcup_{i \in I} \lambda_0(i)  =_{\C P (X)} X. \]  
If $\neq_I$ is an inequality on $I$, and $\neq_X$ an inequality on $X$, we say that $\Lambda(X)$ is a 
family of disjoint subsets of $X$\index{family of disjoint subsets of $X$} $($with respect to $\neq_I)$,  if
\[ \forall_{i, j \in I}\big(i \neq_I j \To \lambda_0(i) \Disj \lambda_0(j)\big), \]
where by Definition~\ref{def: apartsubsets} 
$\lambda_0(i) \Disj \lambda_0(j) :\TOT \forall_{u \in \lambda_0(i)}\forall_{w \in 
\lambda_0(j)}\big(\C E_i^X(u) \neq_X \C E_j^X(w)\big)$.
$\Lambda(X)$ is called a \textit{partition}\index{partition} of $X$, if it covers $X$ and 
it is a family of disjoint subsets of $X$.
\end{definition}

Clearly, $=_{\mathsmaller{\bigcup_{i \in I} \lambda_0(i)}}$ is an equality on 
$\bigcup_{i \in I} \lambda_0(i)$, which is considered to be a set, and the
operation $e_{\mathsmaller{\bigcup}}^{\Lambda(X)}$
is an embedding of $\bigcup_{i \in I} \lambda_0(i)$ into $X$, hence 
$\big(\bigcup_{i \in I} \lambda_0(i), e_{\mathsmaller{\bigcup}}^{\Lambda(X)}\big) \subseteq X$. The inequality
$\neq_{\mathsmaller{\bigcup_{i \in I} \lambda_0(i)}}$ is the canonical inequality of the subset 
$\bigcup_{i \in I} \lambda_0(i)$ of $X$ (see Corollary~\ref{cor: corapartness2}). 
Hence, if $(X, =_X, \neq_X)$ is discrete, then $\big(\bigcup_{i \in I} \lambda_0(i), 
=_{\mathsmaller{\bigcup_{i \in I} \lambda_0(i)}}, \neq_{\mathsmaller{\bigcup_{i \in I} \lambda_0(i)}}\big)$ is discrete,
and if $\neq_X$ is tight, then $=_{\mathsmaller{\bigcup_{i \in I} \lambda_0(i)}}$ is tight.
As the following left diagram commutes, $\Lambda(X)$ covers $X$, if and only of the following right
diagram commutes i.e., if and only if $X \subseteq \bigcup_{i \in I} \lambda_0(i)$
\begin{center}
\resizebox{10cm}{!}{%
\begin{tikzpicture}

\node (E) at (0,0) {$\bigcup_{i \in I} \lambda_0(i)$};
\node[right=of E] (B) {};
\node[right=of B] (F) {$X$};
\node[below=of B] (C) {};
\node[below=of C] (A) {$X$};

\node[right=of F] (G) {$\bigcup_{i \in I} \lambda_0(i)$};
\node[right=of G] (H) {};
\node[right=of H] (J) {$X$};
\node[below=of H] (K) {};
\node[below=of K] (L) {$X$.};

\draw[left hook->,bend left] (E) to node [midway,above] {$e_{\mathsmaller{\bigcup}}^{\Lambda(X)}$} (F);
\draw[right hook->] (E)--(A) node [midway,left] {$e_{\mathsmaller{\bigcup}}^{\Lambda(X)} \ $};
\draw[left hook->] (F)--(A) node [midway,right] {$ \ \id_X$};
\draw[right hook->,bend right] (J) to node [midway,above] {$g$} (G);
\draw[right hook->] (G)--(L) node [midway,left] {$e_{\mathsmaller{\bigcup}}^{\Lambda(X)} \ $};
\draw[left hook->] (J)--(L) node [midway,right] {$ \ \id_X$};

\end{tikzpicture}
}
\end{center}

If $(i, x) =_{\mathsmaller{\bigcup_{i \in I} \lambda_0(i)}} (j, y)$, it is not
necessary that $i =_I j$, hence it is not necessary that 
$(i, x) =_{\mathsmaller{\sum_{i \in I} \lambda_0(i)}} (j, y)$ (as we show in the next proposition, the converse 
implication holds). Consequently, the first projection operation $\pr_1^{\Lambda(X)} 
:= \pr_1^{\Lambda}$\index{$\pr_1^{\Lambda(X)}$}, where $\Lambda$ is the $I$-family of sets induced by $\Lambda(X)$, 
is not necessarily a function! 
The second projection map on $\Lambda(X)$ is defined by 
$\pr_2^{\Lambda(X)} := \pr_2^{\Lambda}$\index{$\pr_2^{\Lambda(X)}$}.
Notice that $\neq_{\mathsmaller{\bigcup_{i \in I} \lambda_0(i)}}$ is an 
inequality on $\bigcup_{i \in I} \lambda_0(i)$, without supposing neither an inequality on $I$, nor an inequality 
on the sets $\lambda_0(i)$'s, as we did in Proposition~\ref{prp: sigmaset1}(ii). 
Moreover, $\neq_{\mathsmaller{\bigcup_{i \in I}}}$ is tight, if $\neq_X$ is tight. Cases (ii) and (iii) of the next
proposition are due to M.~Zeuner\index{Zeuner}.

\begin{proposition}\label{prp: interiorunion1}
Let $\Lambda(X) := (\lambda_0, \C E^X, \lambda_1) \in \Fam(I, X)$.\\[1mm]
\normalfont (i) 
\itshape If $(i, x) =_{\mathsmaller{\sum_{i \in I} \lambda_0(i)}} (j, y)$, then 
$(i, x) =_{\mathsmaller{\bigcup_{i \in I} \lambda_0(i)}} (j, y)$.\\[1mm]
\normalfont (ii) 
\itshape If $e_{\mathsmaller{\bigcup}}^{\Lambda(X)} \colon \sum_{i \in I} \lambda_0(i) \sto X$ is an embedding,
$\big(\sum_{i \in I} \lambda_0(i), e_{\mathsmaller{\bigcup}}^{\Lambda(X)}\big) =_{\C P(X)} \big(\bigcup_{i \in I} \lambda_0(i), e_{\mathsmaller{\bigcup}}^{\Lambda(X)}\big)$.\\[1mm]
\normalfont (iii) 
\itshape If $\neq_I$ is a tight inequality on $X$, and $\Lambda(X)$ is a family of disjoint subsets of $X$ 
with respect to $\neq_I$, then  $e_{\mathsmaller{\bigcup}}^{\Lambda(X)} \colon \sum_{i \in I} \lambda_0(i) \eto X$.

\end{proposition}

\begin{proof}
(i) If $i =_I j$, and since $\C E_j^X$ is a function, we get 
$\C E_i^X(x) = \C E_j^X(\lambda_{ij}(x)) = \C E_j^X(y)$.\\[1mm]
(ii) Let $\C E_i^X(x) =_X \C E_j^X(y) \TOT (i, x) =_{\mathsmaller{\sum_{i \in I}\lambda_0(i)}} (j, y)$. 
We define the operations $\id_1 \colon \sum_{i \in I}\lambda_0(i) \sto \bigcup_{i \in I}\lambda_0(i)$ and 
$\id_2 \colon \bigcup_{i \in I}\lambda_0(i) \sto \sum_{i \in I}\lambda_0(i)$, both 
defined by the identity map-rule. That $\id_1$ is a function, follows from (i).
That $\id_2$ is a function, follows from the hypothesis on  $e_{\mathsmaller{\bigcup}}^{\Lambda(X)}$.\\
(iii) We suppose that $e_{\mathsmaller{\bigcup}}^{\Lambda(X)}(i, x) := \C E_i^X(x) =_X \C E_j^X(y)
:= e_{\mathsmaller{\bigcup}}^{\Lambda(X)}(j, y)$ and we show that 
$(i,x) =_{\mathsmaller{\sum_{i \in I}\lambda_0(i)}} (j, y)$. The converse implication follows from (i).
If $\neg(i \neq_I j)$, then $\lambda_0(i) \Disj \lambda_0(j)$, hence 
$\C E_i^X(x) \neq_X \C E_j^X(y)$, which contradicts our hypothesis. By the tightness of $\neq_I$ we
get $i =_I j$, and it remains to show that $\lambda_{ij}(x) =_{\mathsmaller{\lambda_0(j)}} y$. By 
the equalities $\C E_j^X\big(\lambda_{ij}(x)\big) =_X \C E_i^X(y) =_X \C E_j^X(y)$, and as $\C E_j^X$
is an embedding, we get $\lambda_{ij}(x) =_{\mathsmaller{\lambda_0(j)}} y$.
\end{proof}

\begin{remark}\label{rem: unionA}
Let $i_0 \in I$, $(A, i_A^X) \subseteq X$, and $C^A(X) := (\lambda_0^A, \C E^{A, X}, \lambda_1^A)
\in \Fam(I, X)$ the constant family $A$ of subsets of $X$. Then 
\[ \bigcup_{i \in I}A := \bigcup_{i \in I}\lambda_0^A(i) =_{\C P(X)} A. \] 

\end{remark}

\begin{proof}
By definition $(i, a) =_{\mathsmaller{\bigcup_{i \in I}A}} (j, b) :\TOT \C E^{A,X}_i(a) =_X \C E^{A,X}_j(b) :\TOT
i_A^X(a) =_X i_A^X(b) \TOT a =_A b$. Let the operation $\phi \colon \bigcup_{i \in I}A \sto A$, defined by 
$\phi(i, a) := a$, for every $(i, a) \in \bigcup_{i \in I}A$, and let the operation 
$\theta \colon A \sto \bigcup_{i \in I}A$, defined by $\theta(a) := (i_0, a)$, for every 
$a \in A$. Clearly, $\phi$ and $\theta$ are functions. The required equality of these subsets 
follows from the following equalities: 
$i_A^X(a) =_X e^{C^A(X)}_{\mathsmaller{\bigcup}}(i, a) := \C E^{A, X}_i(a) := i_A^X(a)$, 
and $e^{C^A(X)}_{\mathsmaller{\bigcup}}(i_0, a) := \C E^{A, X}_{i_0}(a) := i_A^X(a)$.
\end{proof}

The interior union of a family of subsets generalises the union of two subsets.

\begin{proposition}\label{prp: interiorunion2}
If $\Lambda^{\D 2}(X)$ is the $\D 2$-family of subsets $A, B$ of $X$, 
$\bigcup_{i \in \D 2} \lambda_0^{\D 2}(i) =_{\C P(X)} A \cup B.$
\end{proposition}

\begin{proof}
The operation $g \colon \bigcup_{i \in \D 2}\lambda_0^{\D 2} (i) \sto A \cup B$, defined by $g(i, x) := \prb_2 (i, x)
:= x$, for every $(i, x) \in \bigcup_{i \in \D 2}\lambda_0^{\D 2} (i)$, is well-defined, and it is an embedding, 
since
$ g(i, x) =_{A \cup B} g(j, y) :\TOT x =_{A \cup B} y 
 :\TOT i^X_A(x) =_X i^X_B(y)
:\TOT (i, x) =_{\mathsmaller{\bigcup_{i \in \D 2}\lambda_0 (i)}} (j, y).
$
The operation $f : A \cup B \sto \bigcup_{i \in \D 2}\lambda_0^{\D 2} (i)$, defined by
$f(z) := (0,z)$, if $z \in A$, and $f(z) := (1, z)$, if $z \in B$, 
is easily seen to be a function. For the commutativity of the following inner diagrams
\begin{center}
\resizebox{5cm}{!}{%
\begin{tikzpicture}

\node (E) at (0,0) {$A \cup B$};
\node[right=of E] (B) {};
\node[right=of B] (F) {$\mathsmaller{\bigcup_{i \in \D 2} \lambda_0^{\D 2}(i)}$};
\node[below=of B] (C) {};
\node[below=of C] (A) {$X$};

\draw[left hook->,bend left] (E) to node [midway,above] {$f$} (F);
\draw[left hook->,bend left] (F) to node [midway,below] {$g$} (E);
\draw[right hook->] (E)--(A) node [midway,left] {$i_{A \cup B}^X \ $};
\draw[left hook->] (F)--(A) node [midway,right] {$ \ e_{\mathsmaller{\bigcup}}^{\Lambda^{\D 2}(X)}$};

\end{tikzpicture}
}
\end{center} 
%
%
%
%
%
%
we use the equalities $i_{A \cup B}^X(g(i, x)) := i_{A \cup B}^X(x) := i_A^X(x) 
:= e_{\mathsmaller{\bigcup}}^{\Lambda^{\D2}(X)}(i, x)$ and 
\[ e_{\mathsmaller{\bigcup}}^{\Lambda^{\D 2}(X)}f(z)) := \left\{ \begin{array}{lll}
                 e_{\mathsmaller{\bigcup}}^{\Lambda^{\D 2}(X)}(0, z)    &\mbox{, $z \in A$}\\
                 {}\\
                 e_{\mathsmaller{\bigcup}}^{\Lambda^{\D 2}(X)}(1, z)             &\mbox{, $z \in B$}
                 \end{array}
          \right.
          := \left\{ \begin{array}{lll}
                 i_A^X(z)    &\mbox{, $z \in A$}\\
                 {}\\
                 i_B^X(z)             &\mbox{, $z \in B$}
                 \end{array}
                 \right.
           := i_{A \cup B}^X(z).\qedhere \] 
\end{proof}

\begin{proposition}\label{prp: interiorunion3}
Let $\Lambda(X) := (\lambda_0, \C E^X, \lambda_1) \in \Fam(I, X)$ and $M(Y) := (\mu_0, \C E^Y, \mu_1) \in \Fam(I, Y)$.
If $\fXY$, the following hold:\\[1mm] 
\normalfont (i) 
\itshape $f \bigg(\bigcup_{i \in I} \lambda_0(i)\bigg) =_{\C P(Y)} \bigcup_{i \in I} f\big(\lambda_0(i)\big)$.\\[1mm]
\normalfont (ii) 
\itshape $f^{-1} \bigg(\bigcup_{i \in I} \mu_0(i)\bigg) =_{\C P(X)} \bigcup_{i \in I} f^{-1}\big(\mu_0(i)\big)$.
\end{proposition}

\begin{proof}
(i) By Definition~\ref{def: image} we have that $f\big(\bigcup_{i \in I} \lambda_0(i)\big) := 
\big(\bigcup_{i \in I} \lambda_0(i), f_{\mathsmaller{\bigcup}}^{\Lambda(X)}\big)$, 
where $f_{\mathsmaller{\bigcup}}^{\Lambda(X)} := f \circ e_{\mathsmaller{\bigcup}}^{\Lambda(X)}$, 
and by Proposition~\ref{prp: newfamsofsubsets2} we have that
\begin{align*}
 (i,x) =_{\mathsmaller{f\big(\bigcup_{i \in I} \lambda_0(i)\big)}} (j,y) 
 & :\TOT f_{\mathsmaller{\bigcup}}^{\Lambda(X)}(i,x) =_Y  f_{\mathsmaller{\bigcup}}^{\Lambda(X)}(j,y)\\
 & :\TOT \big(f \circ e_{\mathsmaller{\bigcup}}^{\Lambda(X)}\big)(i,x) =_Y 
 \big(f \circ e_{\mathsmaller{\bigcup}}^{\Lambda(X)}\big)(j,y)\\
 & :\TOT f\big(\C E_i^X(x)\big) =_Y f\big(\C E_j^X(y)\big)\\
 & :\TOT f_i^Y(x) =_Y f_j^Y(y)
\end{align*}
By Proposition~\ref{prp: newfamsofsubsets2} and Definition~\ref{def: interiorunion} for the subset
$\big(\bigcup_{i \in I} f\big(\lambda_0(i)\big), e_{\mathsmaller{\bigcup}}^{[f(\Lambda)](Y)}\big)$ of $Y$ we have that 
$f(\lambda_0(i)) := \big(\lambda_0(i), f_i^Y\big)$, where $f_i^Y := f \circ \C E_i^X$, for every $i \in I$, and 
$x =_{\mathsmaller{f(\lambda_0(i))}} x{'} :\TOT f_i^Y(x) =_Y f_o^Y(x{'})$. 
Moreover, $e_{\mathsmaller{\bigcup}}^{[f(\Lambda)](Y)}(i,x)
:= \big[f(\C E^X)^Y\big]_i(x) := f_i^Y(x)$. Let the operations $g \colon f\big(\bigcup_{i \in I} \lambda_0(i)\big)
\sto \bigcup_{i \in I} f\big(\lambda_0(i)\big)$ and $h \colon \bigcup_{i \in I} f\big(\lambda_0(i)\big) \sto
f\big(\bigcup_{i \in I} \lambda_0(i)\big)$, defined by the same rule $(i,x) \mapsto (i,x)$.
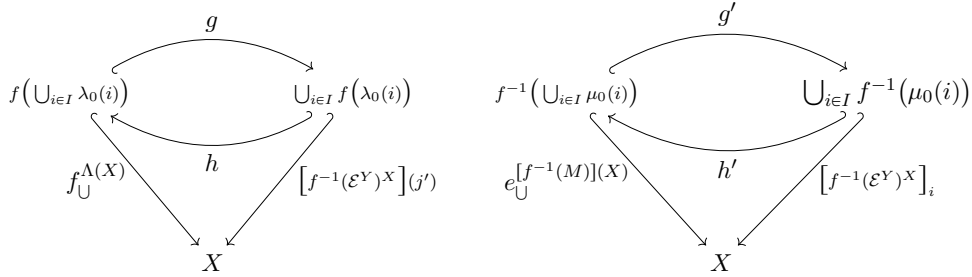
\begin{center}
\resizebox{13cm}{!}{%
\begin{tikzpicture}

\node (E) at (0,0) {$\mathsmaller{f\big(\bigcup_{i \in I} \lambda_0(i)\big)}$};
\node[right=of E] (B) {};
\node[right=of B] (F) {$\mathsmaller{\bigcup_{i \in I} f\big(\lambda_0(i)\big)}$};
\node[below=of B] (C) {};
\node[below=of C] (A) {$X$};

\node[right=of F] (G) {$\mathsmaller{f^{-1} \big(\bigcup_{i \in I} \mu_0(i)\big)}$};
\node[right=of G] (H) {};
\node[right=of H] (J) {$\bigcup_{i \in I} f^{-1}\big(\mu_0(i)\big)$};
\node[below=of H] (K) {};
\node[below=of K] (L) {$X$};

\draw[left hook->,bend left] (E) to node [midway,above] {$g$} (F);
\draw[left hook->,bend left] (F) to node [midway,below] {$h$} (E);
\draw[right hook->] (E)--(A) node [midway,left] {$f_{\mathsmaller{\bigcup}}^{\Lambda(X)} \ $};
\draw[left hook->] (F)--(A) node [midway,right] {$ \ \mathsmaller{\big[f^{-1}(\C E^Y)^X\big](j{'})}$};
\draw[left hook->,bend left] (G) to node [midway,above] {$g{'}$} (J);
\draw[left hook->,bend left] (J) to node [midway,below] {$h{'}$} (G);
\draw[right hook->] (G)--(L) node [midway,left] {$e_{\mathsmaller{\bigcup}}^{[f^{-1}(M)](X)} \ $};
\draw[left hook->] (J)--(L) node [midway,right] {$ \ \mathsmaller{\big[f^{-1}(\C E^Y)^X\big]_i}$};

\end{tikzpicture}
}
\end{center} 
%
%
%
It is immediate by the previous equalities that the above left diagrams commute.\\
(ii) By Definitions~\ref{def: interiorunion} and~\ref{def: image} for the subset 
$\big(\bigcup_{i \in I} \mu_0(i), e_{\mathsmaller{\bigcup}}^{M(Y)}\big)$ of $Y$ we have that the embedding
$e_{\mathsmaller{\bigcup}}^{M(Y)} \colon \bigcup_{i \in I} \mu_0(i) \eto Y$ is given
by the rule $(i,y) \mapsto \C E_i^Y(y)$, and 
\[ f^{-1} \bigg(\bigcup_{i \in I} \mu_0(i)\bigg) := \bigg\{\big(x, (i,y)\big) \in X \times \bigcup_{i \in I} \mu_0(i) 
\mid f(x) =_Y e_{\mathsmaller{\bigcup}}^{M(Y)}(i,y) \bigg\}, \]
with embedding into $X$ the mapping $e^X$, defined by the rule $e^X\big(x, (i,y)\big) := x$. Moreover,
\[ \big(x, (i,y)\big) =_{\mathsmaller{f^{-1}(\bigcup_{i \in I} \mu_0(i))}}
\big(x{'}, (i{'},y{'})\big) :\TOT x =_X x{'} \ \& \ \C E_i^Y(y) =_Y \C E_{i{'}}^Y(y{'}). \]
The subset $f^{-1}(\mu_0(i)) := \{(x,y) \in X \times \mu_0(i) \mid f(x) =_Y \C E_i^Y(y)\}$ of $X$ is equipped with
the embedding $e^X_{\mathsmaller{f^{-1}(\mu_0(i))}} \colon f^{-1}(\mu_0(i)) \eto X$, which is defined by 
$e^X_{\mathsmaller{f^{-1}(\mu_0(i))}}(x, y) := x$, for every $(x,y) \in f^{-1}(\mu_0(i))$. Moreover, we have that
\begin{align*}
\big(x, (i,y)\big) =_{\mathsmaller{\bigcup_{i \in I} f^{-1}\big(\mu_0(i)\big)}}
\big(x{'}, (i{'},y{'})\big) & :\TOT \big[f^{-1}(\C E^Y)^X\big]_i(x,y) =_X \big[f^{-1}(\C E^Y)^X\big]_{i{'}}(x{'},y{'})\\
& :\TOT e_{\mathsmaller{f^{-1}(\mu_0(i))}}^X(x,y) =_X f^{-1}(\mu_0(i{'}))(x{'},y{'})\\
& :\TOT x =_X x{'}.
\end{align*}
If the operation $g{'} \colon f^{-1} \big(\bigcup_{i \in I} \mu_0(i)\big) \sto \bigcup_{i \in I} f^{-1}\big(\mu_0(i)\big)$
is defined by the rule $\big(x, (i,y)\big) \mapsto \big(i, (x,y)\big)$ and the operation
$h{'} \colon  \bigcup_{i \in I} f^{-1}\big(\mu_0(i)\big) \sto f^{-1} \big(\bigcup_{i \in I} \mu_0(i)\big)$ 
is defined by the rule
inverse rule, then it is immediate to show that $g{'}$ is a function. To show that $h{'}$ is a function, 
we suppose that 
$x =_X x{'}$, hence $f(x) =_Y f(x{'})$, and by the definition of $f^{-1} \big(\bigcup_{i \in I} \mu_0(i)\big)$
we get
$e_{\mathsmaller{\bigcup}}^{M(Y)}(i,y) =_Y e_{\mathsmaller{\bigcup}}^{M(Y)}(i{'},y{'}) 
:\TOT \C E_i^Y(y) =_Y \C E_{i{'}}(y{'})$, hence 
$\big(x, (i,y)\big) =_{\mathsmaller{f^{-1} \big(\bigcup_{i \in I} \mu_0(i)\big)}}
\big(x{'}, (i{'},y{'})\big)$. It is immediate to show the commutativity of the 
above right diagrams.
\end{proof}

\index{Extension theorem for coverings}

\begin{theorem}[Extension theorem for coverings]\label{thm: extensioncoverings}
Let $X, Y$ be sets, and let $\Lambda(X) :=  (\lambda_0, \C E^X, \lambda_1) \in \Fam(I, X)$ be a covering of 
$X$. If $f_i \colon \lambda_0(i) \to Y$, for every $i \in I$, such that 
\[ {f_i}_{|\lambda_0(i) \cap \lambda_0(j)} =_{\mathsmaller{\D F(\lambda_0(i) \cap \lambda_0(j), Y)}}
{f_j}_{|\lambda_0(i) \cap \lambda_0(j)}, \]
for every $i, j \in I$, there is a unique 
$\fXY$ such that $f_{|\lambda_0(i)} 
=_{\mathsmaller{\D F(\lambda_0(i), Y)}} f_i$, for every $i \in I$.
\end{theorem}

\begin{proof}
Let $e \colon X \eto \bigcup_{i \in I} \lambda_0(i)$ such that the following diagram commutes
\begin{center}
\resizebox{4cm}{!}{%
\begin{tikzpicture}

\node (E) at (0,0) {$X$};
\node[right=of E] (B) {};
\node[right=of B] (F) {$\bigcup_{i \in I} \lambda_0(i)$};
\node[below=of B] (C) {};
\node[below=of C] (A) {$X$.};

\draw[left hook->,bend left] (E) to node [midway,above] {$e$} (F);
\draw[right hook->] (E)--(A) node [midway,left] {$\id_X \ $};
\draw[left hook->] (F)--(A) node [midway,right] {$ \ e_{\mathsmaller{\bigcup}}^{\Lambda(X)}$};

\end{tikzpicture}
}
\end{center} 
Let the operation 
$f \colon X \sto Y$ defined by  
\[ f(x) := f_{\pr_1^{\Lambda(X)}(e(x))}\big(\pr_2^{\Lambda(X)}(e(x)\big), \]
for every $x \in X$. Hence, if $x \in X$, and $e(x) := (i, u)$, for some $i \in I$ and $u \in \lambda_0(i)$, then 
$f(x) := f_i(u)$. We show that $f$ is a function. Recall that 
$ \lambda_0(i) \cap \lambda_0(j) := \big\{(u,w) \in \lambda_0(i) \times \lambda_0(j) \mid \C E_i^X(u) 
=_X \C E_j^X(w)\big\}$.
If $x, x{'} \in X$, let $e(x) := (i, u)$ and $e(x{'}) := (j, w)$. If $x =_X x{'}$, then
\[ e(x) =_{\mathsmaller{\bigcup_{i \in I} \lambda_0(i)}} e(x{'}) :\TOT
(i,u) =_{\mathsmaller{\bigcup_{i \in I} \lambda_0(i)}} (j, w) :\TOT \C E_i^X(u) =_X \C E_j^X(w) :\TOT
(u,w) \in \lambda_0(i) \cap \lambda_0(j). \]
By the definition of $f$ we have that $f(x) := f_i(u)$ and $f(x{'}) := f_j(w)$. We show that $f_i(u) =_Y f_j(w)$.
Since $\lambda_0(i) \cap \lambda_0(j) 
\subseteq \lambda_0(i)$ and $\lambda_0(i) \cap \lambda_0(j) \subseteq \lambda_0(j)$, and as we have explained right before 
Proposition~\ref{prp: subset3}, by Definition~\ref{def: image} we have that 
\[ {f_i}_{|\lambda_0(i) \cap \lambda_0(j)} := f_i \circ \pr_{\lambda_0(i)} \ \ \& \ \   
{f_j}_{|\lambda_0(i) \cap \lambda_0(j)} := f_j \circ \pr_{\lambda_0(j)}. \]
Since $(u,w) \in \lambda_0(i) \cap \lambda_0(j)$,
by the equality of the restrictions of  $f_i$ and $f_j$ to $\lambda_0(i) \cap \lambda_0(j)$
\[ f_i(u) := \big(f_i \circ \pr_{\lambda_0(i)}\big)(u,w) =_Y \big(f_j \circ \pr_{\lambda_0(j)}\big)(u,w) :=
f_j(w). \]
Next we show that, if $i \in I$, then $f_{|\lambda_0(i)} = f_i$. Since $\C E_i^X \colon (\lambda_0(i), \C E_i^X) \subseteq
(X, \id_X)$
\begin{center}
\resizebox{4cm}{!}{%
\begin{tikzpicture}

\node (E) at (0,0) {$\lambda_0(i)$};
\node[right=of E] (B) {};
\node[right=of B] (F) {$X$};
\node[below=of B] (C) {};
\node[below=of C] (A) {$X$,};

\draw[left hook->,bend left] (E) to node [midway,above] {$\C E_i^X$} (F);
\draw[right hook->] (E)--(A) node [midway,left] {$\C E_i^X \ $};
\draw[left hook->] (F)--(A) node [midway,right] {$ \ \id_X$};

\end{tikzpicture}
}
\end{center} 
by Definition~\ref{def: image} we have that $f_{|\lambda_0(i)} := f \circ \C E_i^X$. If $u \in \lambda_0(i)$, let 
$e\big(\C E_i^X(u)\big) := (j, w)$, for some $j \in I$ and $w \in \lambda_0(j)$. Hence, by the definition of $f$ we get
\[ f_{|\lambda_0(i)}(u) := f\big(\C E_i^X(u)\big) := f_j(w). \]
By the commutativity of the first diagram in this proof we get for $\C E_i^X(u) \in X$ 
\[ \C E_j^X(w) := e_{\mathsmaller{\bigcup}}^{\Lambda(X)}\big(e\big(\C E_i^X(u)\big)\big) =_X \id_X\big(\C E_i^X(u)\big) 
:= \C E_i^X(u) \]
i.e., $(u,w) \in \lambda_0(i) \cap \lambda_0(j)$. Hence, $f_{|\lambda_0(i)}(u) := f_j(w) =_Y f_i(u)$. Finally, let 
$f^* \colon X \to Y$ such that $f_{|\lambda_0(i)}^* := f^* \circ \C E_i^X  =_{\mathsmaller{\D F(\lambda_0(i), Y)}} f_i$, 
for every $i \in I$. If $x \in X$ let $e(x) := (i, u)$, for some $i \in I$ and $u \in \lambda_0(i)$. 
By the commutativity of the first diagram, and since $f^*$ is a function, we get
\begin{align*}
f^*(x) & =_Y f^*\big(e_{\mathsmaller{\bigcup}}^{\Lambda(X)}(e(x))\big)\\
& := f^*\big(e_{\mathsmaller{\bigcup}}^{\Lambda(X)}(i,u)\big)\\
& := f^*\big(\C E_i^X(u)\big)\\
& =_Y f_i(u)\\
& := f_{\pr_1^{\Lambda(X)}(e(x))}\big(\pr_2^{\Lambda(X)}(e(x)\big)\\
& := f(x).\qedhere
\end{align*}
\end{proof}

\begin{corollary}\label{cor: corextcoverings}
Let 
$\Lambda(X) :=  (\lambda_0, \C E^X, \lambda_1) \in \Fam(I, X)$ be a partition of 
$X$. If $f_i \colon \lambda_0(i) \to Y$, for every $i \in I$, there is a unique 
$\fXY$ with $f_{|\lambda_0(i)} 
=_{\mathsmaller{\D F(\lambda_0(i), Y)}} f_i$, for every $i \in I$. 
\end{corollary}

\begin{proof}
The condition ${f_i}_{|\lambda_0(i) \cap \lambda_0(j)} =_{\mathsmaller{\D F(\lambda_0(i) \cap \lambda_0(j), Y)}}
{f_j}_{|\lambda_0(i) \cap \lambda_0(j)}$ of Theorem~\ref{thm: extensioncoverings} is trivially satisfied using
the logical principle Ex falso quodlibet\index{Ex falso quodlibet}. If we suppose that $(u,w) \in 
\lambda_0(i) \cap \lambda_0(j)$, which is impossible as $\lambda_0(i) \Disj \lambda_0(j)$, the equality
$(f_i \circ \pr_{\lambda_0(i)})(u, w) =_Y (f_j \circ \pr_{\lambda_0(j)})(u, w)$, where $i, j \in I$, follows immediately.
\end{proof}

\begin{proposition}\label{prp: unionmap1}
Let $\Lambda(X) := (\lambda_0, \C E^X, \lambda_1)$, $M(X) := (\mu_0, \C Z^X, \mu_1) \in \Fam(I, X)$,
$\Psi : \Lambda(X) \To M(X)$, and $(B, i_B^X) \subseteq X$.\\[1mm]
\normalfont (i) 
\itshape For every $i \in I$ the operation
$e_i^{\Lambda(X)} : \lambda_0(i) \sto \bigcup_{i \in I}\lambda_0(i)$, defined by 
$x \mapsto (i, x)$, is an embedding, and $e_i^{\Lambda(X)} \colon \lambda_0(i) \subseteq 
\bigcup_{i \in I}\lambda_0(i)$.\\[1mm]
\normalfont (ii)
\itshape If $\lambda_0(i) \subseteq B$, for every $i \in I$, then $\bigcup_{i \in I}\lambda_0(i) \subseteq B$.\\[1mm]
\normalfont (iii) 
\itshape The operation
$\bigcup \Psi : \bigcup_{i \in I}\lambda_0(i) \sto \bigcup_{i \in I}\mu_0(i)$, defined by
$\bigcup \Psi (i, x) := (i, \Psi_i (x))$, 
is an embedding, 
such that for every $i \in I$ the following diagram commutes
\begin{center}
\begin{tikzpicture}

\node (E) at (0,0) {$\bigcup_{i \in I}\lambda_0(i)$};
\node[right=of E] (F) {$\bigcup_{i \in I}\mu_0(i)$.};
\node[above=of F] (A) {$\mu_0(i)$};
\node [above=of E] (D) {$\lambda_0(i)$};

\draw[right hook->] (E)--(F) node [midway,below]{$\bigcup \Psi$};
\draw[right hook->] (D)--(A) node [midway,above] {$\Psi_{i}$};
\draw[right hook->] (D)--(E) node [midway,left] {$e_i^{\Lambda(X)}$};
\draw[left hook->] (A)--(F) node [midway,right] {$e_i^{M(X)}$};

\end{tikzpicture}
\end{center}

\end{proposition}

\begin{proof}
(i) If $x, x{'} \in \lambda_0(i)$, and since $\C E_i^X$ is an embedding, we have that 
\[ e_i^{\Lambda(X)}(x) =_{\mathsmaller{\bigcup_{i \in I}\lambda_0(i)}} e_i^{\Lambda(X)}(x{'}) :\TOT 
(i, x) =_{\mathsmaller{\bigcup_{i \in I}\lambda_0(i)}} (i, x{'})
:\TOT \C E_i^X(x) = \C E_i^X(x{'})
\TOT x =_{\lambda_0(i)} x{'}.
\]
Moreover, $e_{\mathsmaller{\bigcup}}^X\big(e_i^{\Lambda(X)}(x)\big) := e_{\mathsmaller{\bigcup}}^X(i, x)
:= \C E_i^X(x)$, hence $e_i^{\Lambda(X)} \colon \lambda_0(i) \subseteq \bigcup_{i \in I}\lambda_0(i)$.\\
(ii)  If $i \in I$ and $e_i^B \colon \lambda_0(i) \subseteq B$, then $i_B^X(e_i^B(x)) =_X \C E^X_i(x)$, 
for every $x \in \lambda_0(i)$. Let the operation $e^B \colon \bigcup_{i \in I}\lambda_0(i) \sto B$, 
defined by $e^B(i,x) :=
e_i^B(x)$, for every $(i,x)  \in \bigcup_{i \in I}\lambda_0(i)$. The operation $e^B$ is a function:
\begin{align*}
(i, x) =_{\mathsmaller{\bigcup_{i \in I}\lambda_0(i)}} (j, y) & :\TOT \C E_i^X(x) =_X \C E_j^X(y)\\
& \To i_B^X(e_i^B(x)) =_X i_B^X(e_j^B(y))\\
& \To e_i^B(x) =_X e_j^B(y)\\
& :\TOT e^B(i,x) =_X e^B(j,y).
\end{align*}
Moroever, $i_B^X\big(e^B(i,x)\big) := i_B^X(e_i^B(x)) =_X \C E_i^X(x) := e_{\mathsmaller{\bigcup}}^{\Lambda(X)}(i,x)$,
hence $e^B \colon \bigcup_{i \in I}\lambda_0(i) \subseteq B$.\\[1mm]
(iii) The required commutativity of the diagram is immediate, and $\bigcup \Psi$ is an embedding, since 
\begin{center}
\begin{tikzpicture}

\node (E) at (0,0) {};
\node[right=of E] (F) {$X$};
\node[above=of E] (A) {$\lambda_0 (i)$};
\node [right=of F] (B) {};
\node [above=of B] (G) {$\mu_0 (i)$};
\node[right=of G] (C) {$\lambda_0 (j)$};
\node [right=of C] (D) {};
\node [right=of D] (E) {$\mu_0 (j)$};
\node [below=of D] (H) {$X$,};

\draw[right hook->] (A)--(F) node [midway,left] {$\C E_{i}^X \ $};
\draw[left hook->] (G)--(F) node [midway,right] {$\ \C Z_i^X$};
\draw[right hook->] (A)--(G) node [midway,above] {$\Psi_{i}$};
\draw[right hook->] (C)--(H) node [midway,left] {$\C E_{j}^X \ $};
\draw[left hook->] (E)--(H) node [midway,right] {$\ \C Z_j^X$};
\draw[right hook->] (C)--(E) node [midway,above] {$\Psi_{j}$};

\end{tikzpicture}
\end{center}
\begin{align*}
(i, x) =_{\mathsmaller{\bigcup_{i \in I} \lambda_0(i)}} (j, y) & :\Leftrightarrow \C E_i^X (x) =_X \C E_j^X (y)\\
& \TOT \C Z_i^X(\Psi_i(x)) =_X \C Z^X_j(\Psi_j(y))\\
& :\TOT (i, \Psi_i(x)) =_{\mathsmaller{\bigcup_{i \in I} \mu_0(i)}} (j, \Psi_j(y))\\
& :\TOT \bigcup \Psi (i, x) =_{\mathsmaller{\bigcup_{i \in I} \mu_0(i)}} \bigcup \Psi (j, y).\qedhere
\end{align*}
\end{proof}


\section{The intersection of a family of subsets}
\label{sec: intersection}

\begin{definition}\label{def: intfamilyofsubsets} 
Let $\Lambda(X) := (\lambda_0, \C E^X, \lambda_1) \in \Fam(I, X)$, and $i_0 \in I$. 
The intersection $\bigcap_{i \in I} \lambda_0 (i)$\index{intersection of a family of
subsets}\index{$\bigcap_{i \in I} \lambda_0 (i)$} of $\Lambda(X)$ is the 
totality defined by
\[ \Phi \in \bigcap_{i \in I} \lambda_0 (i) : \TOT \Phi \in \D A(I, \lambda_0)
 \ \& \ \forall_{i, j \in I}\big(\C E_i^X (\Phi_i) =_X \C E_{j}^X (\Phi_{j})\big). \]
Let $e_{\mathsmaller{\bigcap}}^{\Lambda(X)} : \bigcap_{i \in I} \lambda_0 (i) \sto X$ be defined by
$e_{\mathsmaller{\bigcap}}^{\Lambda(X)}(\Phi) := \C E_{i_0}^X\big(\Phi_{i_0}\big)$, for
every $\Phi \in \bigcap_{i \in I} \lambda_0 (i)$, and 
\[ \Phi =_{\mathsmaller{\bigcap_{i \in I} \lambda_0 (i)}} \Theta : \TOT e_{\mathsmaller{\bigcap}}^{\Lambda(X)}(\Phi) 
=_X e_{\mathsmaller{\bigcap}}^{\Lambda(X)}(\Theta) :\TOT 
\C E_{i_0}^X\big(\Phi_{i_0}\big) =_X \C E_{i_0}^X\big(\Theta_{i_0}\big),
\]
If $\neq_X$ is a given inequality on $X$, let 
$\Phi \neq_{\mathsmaller{\bigcap_{i \in I} \lambda_0 (i)}} \Theta :\TOT 
 \C E_{i_0}^X\big(\Phi_{i_0}\big) \neq_X \C E_{i_0}^X\big(\Theta_{i_0}\big)$.

\end{definition}

Clearly, $=_{\bigcap_{i \in I} \lambda_0(i)}$ is an equality on 
$\bigcap_{i \in I} \lambda_0(i)$, which is considered to be a set, and $e_{\mathsmaller{\bigcap}}^{\Lambda(X)}$
is an embedding, hence $\big(\bigcap_{i \in I} \lambda_0(i), 
e_{\mathsmaller{\bigcap}}^{\Lambda(X)}\big) \subseteq X$. Moreover, the inequality
$\neq_{\mathsmaller{\bigcap_{i \in I} \lambda_0(i)}}$ is the canonical inequality of the subset 
$\bigcap_{i \in I} \lambda_0(i)$ of $X$ (see Corollary~\ref{cor: corapartness2}).

\begin{proposition}\label{prp: intersection1}
Let $\Lambda(X) := (\lambda_0, \C E^X, \lambda_1) \in \Fam(I, X)$.\\[1mm]
\normalfont (i) 
\itshape $\Phi =_{\bigcap_{i \in I} \lambda_0 (i)} \Theta \TOT \Phi =_{\D A(I, \lambda_0)} \Theta$.\\[1mm]
\normalfont (ii) 
\itshape If $\Phi \in \bigcap_{i \in I} \lambda_0 (i)$, then $\Phi \in \prod_{i \in I}\lambda_0(i)$.\\[1mm]
\normalfont (iii) 
\itshape If $(X, =_X, \neq_X)$ is discrete, the set $\big(\bigcap_{i \in I} \lambda_0(i), 
=_{\mathsmaller{\bigcap_{i \in I} \lambda_0(i)}}, \neq_{\mathsmaller{\bigcap_{i \in I} \lambda_0(i)}}\big)$ is discrete.
\end{proposition}

\begin{proof}
(i) To show the implication $(\To)$, if $i \in I$, then 
$\C E_i^X(\Phi_i) =_X \C E_{i_0}^X (\Phi_{i_0}) =_X \C E_{i_0}^X (\Theta_{i_0}) =_X \C E_i^X(\Theta_i),$
and since $\C E_i^X$ is an embedding, $\Phi_i =_{\lambda_0(i)} \Theta_i$. For the converse implication,
the pointwise equality of $\Phi$ and $\Theta$ implies that $\Phi_{i_0} =_{\lambda_0(i_0)} \Theta_{i_0}$, hence
$\C E_{i_0}^X (\Phi_{i_0}) =_X \C E_{i_0}^X (\Theta_{i_0})$.\\
(ii) If $i =_I j$, then 
$\C E_j^X \big(\lambda_{ij}(\Phi_i)\big) =_X \C E_i^X(\Phi_i) =_X \C E_j(\Phi_j)$,
and as $\C E_j^X$ is an embedding, we get the required equality $\lambda_{ij}(\Phi_i) =_{\lambda_0(j)}
\Phi_j$. The proof of (iii) is immediate. 
\end{proof}

Since the equality of $\prod_{i \in I}\lambda_0(i)$ is the pointwise equality of $\D A(I, \lambda_0)$,
then, as we explained above, the equality of $\prod_{i \in I}\lambda_0(i)$ is the equality of
$\bigcap_{i \in I} \lambda_0 (i)$.

\begin{remark}\label{rem: intA}
Let $i_0 \in I$, $(A, i_A^X) \subseteq X$, and $C^A(X) := (\lambda_0^A, \C E^{A, X}, \lambda_1^A) \in \Fam(I, X)$
the constant family $A$ of subsets of $X$. Then 
\[ \bigcap_{i \in I}A := \bigcap_{i \in I}\lambda_0^A(i) =_{\C P(X)} A. \] 

\end{remark}

\begin{proof}
We proceed similarly to the proof of Remark~\ref{rem: unionA}.
\end{proof}


\begin{proposition}\label{prp: intersection2}
If $\Lambda^{\D 2}(X)$ is the $\D 2$-family of subsets $A, B$ of $X$, 
$\bigcap_{i \in \D 2} \lambda_0^{\D 2}(i) =_{\C P(X)} A \cap B.$
\end{proposition}

\begin{proof}
By definition $\Phi \in \bigcap_{i \in I} \lambda_0^{\D 2} (i) : \TOT \Phi \colon \bigcurlywedge_{i \in I}\lambda_0^{\D 2}
(i)$ and for every $i, j \in \D 2$ we have that $\C E_i^X (\Phi_i) =_X \C E_{j}^X (\Phi_{j})$, 
where $\C E_0^X := i_A^X$ and 
$\C E_1^X := i_B^X$. Moreover, $e_{\mathsmaller{\bigcap}}^{\Lambda^{\D 2}(X)} : \bigcap_{i \in I}
\lambda_0^{\D 2} (i) \sto X$ is 
given by $e_{\mathsmaller{\bigcap}}^{\Lambda^{\D 2}(X)}(\Phi) := \C E_{0}^X\big(\Phi_{0}\big)$, 
for every $\Phi \in \bigcap_{i \in I}
\lambda_0^{\D 2} (i)$, and $\Phi =_{\mathsmaller{\bigcap_{i \in I} \lambda_0^{\D 2} (i)}} \Theta :\TOT 
\C E_{0}^X\big(\Phi_{0}\big) =_X \C E_{0}^X\big(\Theta_{0}\big)$. Let $f \colon A \cap B 
\sto \bigcap_{i \in \D 2} \lambda_0^{\D 2}(i)$ be defined by $f(a,b) := \Phi_{(a,b)}$, for every $(a,b) \in A \cap B$,
where $\Phi_{(a,b)} \colon \bigcurlywedge_{i \in I}\lambda_0^{\D 2}(i)$, such that 
$\Phi_{(a,b)}(0) := a$ and $\Phi_{(a,b)}(1) := b$.
Since $\C E_0^X\big(\Phi_{(a,b)}(0)\big) =_X \C E_1^X\big(\Phi_{(a,b)}(1)\big) \TOT \C E_0^X(a) =_X 
\C E_1(b) :\TOT i_A^X(a) =_X i_B^X(b)$, where the last equality holds by the definition of $A \cap B$ 
(see Definition~\ref{def: intofsubsets}), the operation $f$ is well-defined. It is straightforward 
to show that $f$ is a 
function. Let the operation $g \colon \bigcap_{i \in \D 2} \lambda_0^{\D 2}(i) \sto A \cap B$, defined by 
$g(\Phi) := (\Phi_0, \Phi_1)$, for every $\Phi \in \bigcap_{i \in \D 2} \lambda_0^{\D 2}(i)$.
Since $i_Z^X(\Phi_0) := \C E_0^X(\Phi_0) =_X \C E_1^X(\Phi_1) := i_B^X(\Phi_1)$, we have that 
$g$ is well-defined. It is easy to show that $g$ is a function, 
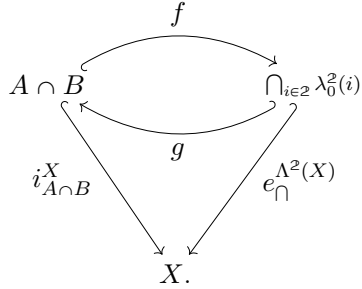
\begin{center}
\resizebox{5cm}{!}{%
\begin{tikzpicture}

\node (E) at (0,0) {$A \cap B$};
\node[right=of E] (B) {};
\node[right=of B] (F) {$\mathsmaller{\bigcap_{i \in \D 2} \lambda_0^{\D 2}(i)}$};
\node[below=of B] (C) {};
\node[below=of C] (A) {$X$.};

\draw[left hook->,bend left] (E) to node [midway,above] {$f$} (F);
\draw[left hook->,bend left] (F) to node [midway,below] {$g$} (E);
\draw[right hook->] (E)--(A) node [midway,left] {$i_{A \cap B}^X \ $};
\draw[left hook->] (F)--(A) node [midway,right] {$ \ e_{\mathsmaller{\bigcap}}^{\Lambda^{\D 2}(X)}$};

\end{tikzpicture}
}
\end{center}
and the above inner diagrams commute.
\end{proof}

\begin{proposition}\label{prp: intersection3}
Let $\Lambda(X) := (\lambda_0, \C E^X, \lambda_1) \in \Fam(I, X)$ and $M(Y) := (\mu_0, \C E^Y, \mu_1) \in \Fam(I, Y)$.
If $\fXY$, the following hold:\\[1mm] 
\normalfont (i) 
\itshape $f \bigg(\bigcap_{i \in I} \lambda_0(i)\bigg) \subseteq \bigcap_{i \in I} f\big(\lambda_0(i)\big)$.\\[1mm]
\normalfont (ii) 
\itshape $f^{-1} \bigg(\bigcap_{i \in I} \mu_0(i)\bigg) =_{\C P(X)} \bigcap_{i \in I} f^{-1}\big(\mu_0(i)\big)$.
\end{proposition}

\begin{proof}
 We proceed similarly to the proof of Proposition~\ref{prp: interiorunion3}.
\end{proof}

\begin{proposition}\label{prp: internalmap1}
Let $\Lambda(X) := (\lambda_0, \C E^X, \lambda_1), M_X := (\mu_0, \C Z^X, \mu_1) \in \Fam(I,X)$, let $i_0 \in I$,
$\Psi : \Lambda(X) \To M(X)$, and $(B, i_B^X) \subseteq X$.\\[1mm]
\normalfont (i) 
\itshape The operation
$\pi_i^{\Lambda(X)} : \bigcap_{i \in I}\lambda_0(i) \sto \lambda_0(i)$, defined by
$\Theta \mapsto \Theta_i,$ is a function, and $\pi_i^{\Lambda(X)} \colon \bigcap_{i \in I}\lambda_0(i) 
\subseteq \lambda_0(i)$, for every $i \in I$.\\[1mm]
\normalfont (ii) 
\itshape If $B \subseteq \lambda_0(i)$, for every $i \in I$, then $B \subseteq \bigcap_{i \in I}\lambda_0(i)$.\\[1mm]
\normalfont (iii)
\itshape The operation
$\bigcap \Psi : \bigcap_{i \in I}\lambda_0(i) \sto \bigcap_{i \in I}\mu_0(i)$, defined by
$[\bigcap \Psi (\Theta)]_i := \Psi_i (\Theta_i)$, for every $i \in I$, is an embedding, 
such that for every $i \in I$ 
the following diagram commutes
\begin{center}
\begin{tikzpicture}

\node (E) at (0,0) {$\bigcap_{i \in I}\lambda_0(i)$};
\node[right=of E] (F) {$\bigcap_{i \in I}\mu_0(i)$.};
\node[above=of F] (A) {$\mu_0(i)$};
\node [above=of E] (D) {$\lambda_0(i)$};

\draw[right hook->] (E)--(F) node [midway,below]{$\bigcap \Psi$};
\draw[right hook->] (D)--(A) node [midway,above] {$\Psi_{i}$};
\draw[->] (E)--(D) node [midway,left] {$\pi_i^{\Lambda(X)}$};
\draw[->] (F)--(A) node [midway,right] {$\pi_i^{M(X)}$};

\end{tikzpicture}
\end{center}
\end{proposition}

\begin{proof}
(i) Since $\Phi =_{\bigcap_{i \in I} \lambda_0 (i)} \Theta \To \Phi =_{\D A(I, \lambda_0)} \Theta$, we get
$\Phi_i = \Theta_i$, for every $i \in I$. Since $\C E_i^X\big(\pi_i^{\Lambda(X)}(\Theta)\big) := \C E_i^X(\Theta_i) =_X 
\C E_{i_0}^X(\Theta_{i_0}) := e^X_{\mathsmaller{\bigcap}}(\Theta)$, we get 
$\pi_i^{\Lambda(X)} \colon \bigcap_{i \in I}\lambda_0(i) \subseteq \lambda_0(i)$.\\[1mm]
(ii) If $i \in I$, let $e_B^i \colon B \subseteq \lambda_0(i)$, hence $\C E_i^X(e_B^i(b)) =_X i_B^X(b)$, 
for every $b \in B$.
Let the operation $e_B \colon B \sto \bigcap_{i \in I}\lambda_0(i)$, defined by the rule $b \mapsto e_B(b)$, 
where
$[e_B(b)]_i := e_B^i(b)$, for every $b \in B$ and $i \in I$. First we show that $e_B$ is well defined.
If $i, j \in I$, then 
\[\C E_i^X\big([e_B(b)]_i\big) =_X \C E_j^X\big([e_B(b)]_j\big) :\TOT \C E_i^X(e_B^i(b)) =_X \C E_j^X(e_B^j(b)) \TOT 
i_B^X(b) =_X i_B^X(b). \] 
Clearly, $e_B$ is a function. Moreover, $e_B \colon B \subseteq \bigcap_{i \in I}\lambda_0(i)$, since, 
for every $b \in B$,
\[ e_{\mathsmaller{\bigcap}}(e_B(b)) := \C E_{i_0}^X\big([e_B(b)]_{i_0}\big) := 
\C E_{i_0}^X(e_B^{i_0}(b)) =_X i_B^X(b). \]
(iii) It suffices to show that $\bigcap \Psi$ is an embedding.
If $\Phi, \Theta \in \bigcap_{i \in I} \lambda_0 (i)$, then 
\begin{align*}
\Phi =_{\mathsmaller{\bigcap_{i \in I} \lambda_0 (i)}} \Theta & :\TOT 
\C E_{i_0}^X\big(\Phi_{i_0}\big) =_X \C E_{i_0}^X\big(\Theta_{i_0}\big) \\
& \TOT \C Z_{i_0}^X(\Psi_{i_0}(\Phi_{i_0})) =_X \C Z^X_{i_0}(\Psi_{i_0}(\Theta_{i_0}))\\
& :\TOT \C Z_{i_0}^X\bigg(\bigg[\bigcap \Psi(\Phi)\bigg]_{i_0}\bigg) =_X 
\C Z_{i_0}^X\bigg(\bigg[\bigcap \Psi(\Theta)\bigg]_{i_0}\bigg)\\
& :\TOT \bigg(\bigcap \Psi\bigg) (\Phi) =_{\mathsmaller{\bigcap_{i \in I} \mu_0 (i)}} 
\bigg(\bigcap \Psi\bigg)(\Theta).\qedhere
\end{align*}
\end{proof}

The above notions and results can be generalised as follows.

\begin{definition}\label{def: genmapfromltom}
Let $X$ and $Y$ be sets, and $h : X \to Y$. 
If $\Lambda(X) := (\lambda_0, \C E^X, \lambda_1) \in \Fam(I,X)$, and 
$M(Y) := (\mu_0, \C Z^Y, \mu_1) \in \Fam(I,Y)$, a 
\textit{family of subsets-map} from $\Lambda(X)$ to $M(Y)$ is a dependent operation
$\Psi : \bigcurlywedge_{i \in I}\D F\big(\lambda_0(i), \mu_0(i)\big)$, where if
$\Psi(i) := \Psi_i$, for every $i \in I$, then, for every $i \in I$, the following diagram commutes
\begin{center}
\begin{tikzpicture}

\node (E) at (0,0) {$\lambda_0 (i)$};
\node[right=of E] (F) {$\mu_0 (i)$};
\node[below=of E] (A) {$X$};
\node[below=of F] (B) {$Y$.};

\draw[right hook->] (E)--(A) node [midway,left] {$\C E_{i}^X \ $};
\draw[left hook->] (F)--(B) node [midway,right] {$\ \C Z^Y_i$};
\draw[->] (E)--(F) node [midway,above] {$\Psi_{i}$};
\draw[->] (A)--(B) node [midway,below] {$h$};

\end{tikzpicture}
\end{center}
The totality $\Map_{I,h}(\Lambda(X), M(Y))$\index{$\Map_{I,h}(\Lambda(X), M(Y))$} of family of
subsets-maps from $\Lambda(X)$ to $M(Y)$ is
equipped with the pointwise equality, 
and we write $\Psi \colon \Lambda(X) \stackrel{h \ } \To M(Y)$\index{$\Psi \colon \Lambda(X) \stackrel{h \ } \To M(Y)$},
if $\Psi \in \Map(\Lambda(X), M(X))$.
If
$\Xi \colon M(Y) \stackrel{g \ } \To N(Z)$, where $g 
\colon Y \to Z$, the \textit{composition family of subsets-map} 
$\Xi \circ \Psi : \Lambda(X) \stackrel{g \circ h} \Longrightarrow N(Z)$ is defined by 
$(\Xi \circ \Psi)(i) := \Xi_i \circ \Psi_i$, for every $i \in I$
\begin{center}
\begin{tikzpicture}

\node (E) at (0,0) {$\lambda_0(i)$};
\node[right=of E] (F) {$\mu_0(i)$};
\node[right=of F] (G) {$\nu_0(i)$};
\node[below=of E] (A) {$X$};
\node[below=of F] (B) {$Y$};
\node[below=of G] (C) {$Z$.};

\draw[right hook->] (E)--(A) node [midway,left] {$\C E_{i}^X \ $};
\draw[left hook->] (F)--(B) node [midway,right] {$\ \C Z_i^Y$};
\draw[left hook->] (G)--(C) node [midway,right] {$\C H_i^Z$};
\draw[->] (E)--(F) node [midway,above] {$\Psi_{i}$};
\draw[->] (A)--(B) node [midway,below] {$h$};
\draw[->] (F)--(G) node [midway,above] {$\Xi_i$};
\draw[->] (B)--(C) node [midway,below] {$g$};

\end{tikzpicture}
\end{center} 

\end{definition}

If $Y := X$, and $h := \id_X$, and if $\Psi \colon \Lambda(X) \stackrel{\id_X} \Longrightarrow
M(X)$, then $\Psi \colon \Lambda(X) \To M(X)$. In the general case, if 
$\Psi \colon \Lambda(X) \stackrel{h \ } \To M(Y)$, then $\Psi_i$ is an embedding, if $h$ is an embedding.  

\begin{proposition}\label{prp: geninternalmap1}
Let $X$ and $Y$ be sets, and $h : X \to Y$. Let also $\Lambda(X) := (\lambda_0, \C E^X, \lambda_1) \in \Fam(I,X)$, 
$M(Y) := (\mu_0, \C Z^Y, \mu_1) \in \Fam(I,Y)$, and $\Psi \colon \Lambda(X) \stackrel{h \ } \To M(Y)$.\\[1mm]
\normalfont (i) 
\itshape The operation
$\bigcup_h \Psi : \bigcup_{i \in I}\lambda_0(i) \sto \bigcup_{i \in I}\mu_0(i)$, defined by
$\big(\bigcup_h \Psi\big) (i, u) := (i, \Psi_i (u))$, for every $(i, u) \in \bigcup_{i \in I}\lambda_0(i)$, is  
a  function,
and for every $i \in I$ the following left diagram commutes
\begin{center}
\begin{tikzpicture}

\node (E) at (0,0) {$\bigcup_{i \in I}\lambda_0(i)$};
\node[right=of E] (F) {$\bigcup_{i \in I}\mu_0(i)$};
\node[above=of F] (A) {$\mu_0(i)$};
\node [above=of E] (D) {$\lambda_0(i)$};
\node[right=of F] (G) {$\bigcap_{i \in I}\lambda_0(i)$};
\node[right=of G] (K) {$\bigcap_{i \in I}\mu_0(i)$.};
\node[above=of K] (L) {$\mu_0(i)$};
\node [above=of G] (M) {$\lambda_0(i)$};

\draw[->] (E)--(F) node [midway,below]{$\bigcup_h \Psi$};
\draw[->] (D)--(A) node [midway,above] {$\Psi_{i}$};
\draw[right hook->] (D)--(E) node [midway,left] {$e_i^{\Lambda(X)}$};
\draw[left hook->] (A)--(F) node [midway,right] {$e_i^{M(Y)}$};
\draw[->] (G)--(K) node [midway,below]{$\bigcap_h \Psi$};
\draw[->] (M)--(L) node [midway,above] {$\Psi_{i}$};
\draw[->] (G)--(M) node [midway,left] {$\pi_i^{\Lambda(X)}$};
\draw[->] (K)--(L) node [midway,right] {$\pi_i^{M(Y)}$};

\end{tikzpicture}
\end{center}
\normalfont (ii) 
\itshape If $i_0 \in I$, the operation
$\bigcap_h \Psi : \bigcap_{i \in I}\lambda_0(i) \sto \bigcap_{i \in I}\mu_0(i)$, defined by
$[\bigcap_h \Psi (\Theta)]_i := \Psi_i (\Theta_i)$, for every $i \in I$, is a function,
such that for every $i \in I$ the above right diagram commutes.
%
%
%
%
%
\end{proposition}

\begin{proof}
(i) The commutativity of the diagram is trivial, and we show that $\bigcup_h \Psi$ is a function:
\begin{align*}
(i, u) =_{\mathsmaller{ \bigcup_{i \in I}\lambda_0(i)}} (j, w) & :\TOT \C E_i(u) =_X \C E_j(w)\\
& \To h(\C E_i(u)) =_Y h(\C E_j(w))\\
& \TOT E_i(\Psi_i(u)) =_Y E_j(\Psi_j(w))\\
& :\TOT (i, \Psi_i(u)) =_{\mathsmaller{\bigcup_{i \in I}\mu_0(i)}} (j, \Psi_j(w))\\
& :\TOT \bigg(\bigcup_h \Psi\bigg) (i, u) =_{\mathsmaller{\bigcup_{i \in I}\mu_0(i)}} \bigg(\bigcup_h \Psi\bigg)
(j, w).
\end{align*}
(ii) The commutativity of the diagram is trivial, and we show that $\bigcap_h \Psi$ is a function:
\begin{align*}
\Phi =_{\mathsmaller{\bigcap_{i \in I}\lambda_0(i)}} \Theta & \TOT \C E_{i_0}(\Phi_{i_0}) =_X 
\C E_{i_0}(\Theta_{i_0})\\
& \To h\big(\C E_{i_0}(\Phi_{i_0})\big) =_Y h \big(\C E_{i_0}(\Theta_{i_0})\big)\\
& \TOT E_{i_0}\big(\Psi_{i_0}(\Phi_{i_0})\big) =_Y E_{i_0}\big(\Psi_{i_0}(\Theta_{i_0})\big)\\
& :\TOT E_{i_0}\bigg(\bigg[\bigcap_h \Psi (\Phi)\bigg]_{i_0}\bigg) =_Y
E_{i_0}\bigg(\bigg[\bigcap_h \Psi (\Theta)\bigg]_{i_0}\bigg)\\
& :\TOT \bigg(\bigcap_h \Psi\bigg) (\Phi) =_{\mathsmaller{\bigcap_{i \in I}\mu_0(i)}}
\bigg(\bigcap_h \Psi\bigg) (\Theta).\qedhere
\end{align*}
\end{proof}


\section{Families of subsets over products}
\label{sec: famsofsubsetsprod}

\begin{proposition}\label{prp: famsubsetsprod1}
Let $\Lambda(X) := (\lambda_0, \C E^X, \lambda_1), K(X) := (k_0, \C H^X, k_1) \in \Fam(I, X)$, 
$M(Y) := (\mu_0, \C E^Y, \mu_1)$, and $N(Y) := (\nu_0, \C H^Y, \nu_1) \in \Fam(J, Y)$.\\[1mm]
\normalfont (i) 
\itshape $(\Lambda \otimes M)(X \times Y) := (\lambda_0 \otimes \mu_0, \C E^X \otimes \C E^Y, \lambda_1 \otimes \mu_1)
\in \Fam(I \times J, X 
\times Y)$, where 
\[ (\lambda_0 \otimes \mu_0)(i,j) := \lambda_0(i) \times \mu_0(j); \ \ \ \ (i,j) \in I \times J, \]
\[ \big(\C E^X \otimes \C E^Y\big)_{(i,j)} \colon \lambda_0(i) \times \mu_0(j) \eto X \times Y, \]
\[  \big(\C E^X \otimes \C E^Y\big)_{(i,j)}(u, w) := \big(\C E^X_i(u), \C E^Y_j(w)\big); \ \ \ \ 
(u,w) \in \lambda_0(i) \times \mu_0(j), \ (i,j) \in I \times J,  
\]
\[  (\lambda_1 \otimes \mu_1)_{(i,j)(i{'}j{'})} \colon \lambda_0(i) \times \mu_0(j) \to \lambda_0(i{'})
\times \mu_0(j{'}),  \]
\[ (\lambda_1 \otimes \mu_1)_{(i,j)(i{'}j{'})}(u, w) := \big(\lambda_{ii{'}}(u), \mu_{jj{'}}(w)\big); \ \ \ \ 
(u,w) \in \lambda_0(i) \times \mu_0(j).
\]
\normalfont (ii) 
\itshape If $\Phi \colon \Lambda(X) \To K(X)$ and $\Psi \colon M(X) \To N(X)$, then $\Phi \otimes \Psi \colon 
(\Lambda \otimes M)(X \times Y) \To (K \otimes N)(X \times Y)$, where, for every $(i,j) \in I \times J$,
\[ (\Phi \otimes \Psi)_{(i,j)} \colon \lambda_0(i) \times \mu_0(j) \to k_0(i) \times \nu_0(j),  \]
\[  (\Phi \otimes \Psi)_{(i,j)}(u,w) := \big(\Phi_i(u), \Psi_j(w)\big); \ \ \ \ (u,w) \in \lambda_0(i) \times \mu_0(j).
\]
\normalfont (iii) 
\itshape The following equality holds
\[ \bigcup_{(i,j) \in I \times J}\big(\lambda_0(i) \times \mu_0(j)\big) 
=_{\C P(X \times Y)} \bigg(\bigcup_{i \in I}\lambda_0(i)\bigg)
 \times \bigg(\bigcup_{j \in J}\mu_0(j)\bigg). \]
\normalfont (iv) 
\itshape  If $i_0 \in I$ and $j_0 \in J$, the following equality holds
\[ \bigcap_{(i,j) \in I \times J}\big(\lambda_0(i) \times \mu_0(j)\big) =_{\C P(X \times Y)}
\bigg(\bigcap_{i \in I}\lambda_0(i)\bigg)
 \times \bigg(\bigcap_{j \in J}\mu_0(j)\bigg). \]
 \normalfont (v) 
\itshape If $\Lambda(X)$ covers $X$ and $M(Y)$ covers $Y$, then $(\Lambda \otimes M)(X \times Y)$ 
covers $X \times Y$.\\[1mm]
\normalfont (vi) 
\itshape Let the inequalities $\neq_I, \neq_J, \neq_X$ and $\neq_Y$ on $I, J, X$ and $Y$, respectively. 
If $\Lambda(X)$ is a partition of $X$ and $M(Y)$ is a partition of $Y$, then $(\Lambda \otimes M)(X \times Y)$ 
is a partition of $X \times Y$.
\end{proposition}

\begin{proof}
The proofs of (i)-(iv) are the internal analogue to the proofs of Proposition~\ref{prp: famsetsprod1}(i)-(iii).\\
(v) Since $X =_{\C P(X)} \bigcup_{i \in I} \lambda_0(i)$ and $Y =_{\C P(Y)} \bigcup_{j \in J} \mu_0(j)$, by case 
(iii) and by Proposition~\ref{prp: subset15}(iv) we have that
\[ X \times Y =_{\C P(X \times Y)} \bigg(\bigcup_{i \in I}\lambda_0(i)\bigg)
 \times \bigg(\bigcup_{j \in J}\mu_0(j)\bigg) =_{\C P(X \times Y)} 
 \bigcup_{(i,j) \in I \times J}\big(\lambda_0(i) \times \mu_0(j)\big).
 \]
(vi) By Definition~\ref{def: interiorunion} we have that
\[ i \neq_I i{'} \To \C E_i^X(u) \neq_X \C E_{i{'}}^X (u{'}); \ \ \ \ i, i{'} \in I, \ u \in \lambda_0(i), \ u{'} \in
 \lambda_0(i{'}), 
\]
\[ j \neq_J j{'} \To \C E_i^Y(w) \neq_Y \C E_{j{'}}^Y (w{'}); \ \ \ \ j, j{'} \in J, \ w \in \mu_0(j), \ w{'} \in
 \mu_0(j{'}). 
\]
Let $(i,j) \neq_{I \times J} (i{'}, j{'}) :\TOT i \neq_I i{'} \vee j \neq_J j{'}$. If $i \neq_I i{'}$ is the case, then
$\C E_i^X(u) \neq_X \C E_{i{'}}^X (u{'})$, hence $ \big(\C E^X_i(u), \C E^Y_j(w)\big) \neq_{X \times Y} 
 \big(\C E^X_i(u{'}), \C E^Y_j(w{'})\big)$. If $j \neq_J j{'}$, we proceed similarly.
\end{proof}

\begin{proposition}\label{prp: famsubsetsprod2}
Let $\Lambda(X) := (\lambda_0, \C E^X, \lambda_1), K(X) := (k_0, \C H^X, k_1) \in \Fam(I, X)$,  
$M(X) := (\mu_0, \C Z^X, \mu_1)$, and $N(X) := (\nu_0, \C F^X, \nu_1) \in \Fam(J, X)$.\\[1mm]
\normalfont (i) 
\itshape $(\Lambda \wedge M)(X) := (\lambda_0 \wedge \mu_0, \C E^X \wedge \C Z^X, \lambda_1 \wedge \mu_1)
\in \Fam(I \times J, X)$, where 
\[ (\lambda_0 \wedge \mu_0)(i,j) := \lambda_0(i) \cap \mu_0(j); \ \ \ \ (i,j) \in I \times J, \]
\[ \big(\C E^X \wedge \C Z^X\big)_{(i,j)} \colon \lambda_0(i) \cap \mu_0(j) \eto X, \]
\[  \big(\C E^X \wedge \C Z^X\big)_{(i,j)}(u, w) := \C E^X_i(u); \ \ \ \ 
(u,w) \in \lambda_0(i) \cap \mu_0(j), \ (i,j) \in I \times J,  
\]
\[  (\lambda_1 \wedge \mu_1)_{(i,j)(i{'}j{'})} \colon \lambda_0(i) \cap \mu_0(j) \to \lambda_0(i{'})
\cap \mu_0(j{'}),  \]
\[ (\lambda_1 \otimes \mu_1)_{(i,j)(i{'}j{'})}(u, w) := \big(\lambda_{ii{'}}(u), \mu_{jj{'}}(w)\big); \ \ \ \ 
(u,w) \in \lambda_0(i) \cap \mu_0(j).
\]
\normalfont (ii) 
\itshape $(\Lambda \vee M)(X) := (\lambda_0 \vee \mu_0, \C E^X \vee \C Z^X, \lambda_1 \vee \mu_1)
\in \Fam(I \times J, X)$, where 
\[ (\lambda_0 \vee \mu_0)(i,j) := \lambda_0(i) \cup \mu_0(j); \ \ \ \ (i,j) \in I \times J, \]
\[ \big(\C E^X \vee \C Z^X\big)_{(i,j)} \colon \lambda_0(i) \cup \mu_0(j) \eto X, \]
\[ \big(\C E^X \vee \C Z^X\big)_{(i,j)}(z) := \left\{ \begin{array}{ll}
                 \C E_i^X(z)    &\mbox{, $z \in \lambda_0(i)$}\\
                 \C Z_j^Y(z)             &\mbox{, $z \in \mu_0(j)$}
                 \end{array}
          \right.
          ; \ \ \ \ i \in I, \ z \in \lambda_0(i) \cup \mu_0(i) \]
\[  (\lambda_1 \vee \mu_1)_{(i,j)(i{'}j{'})}(z) := \left\{ \begin{array}{ll}
                 \lambda_{ii{'}}(z)    &\mbox{, $z \in \lambda_0(i)$}\\
                 \mu_{jj{'}}(z)             &\mbox{, $z \in \mu_0(j)$}
                 \end{array}
          \right.
          ; \ \ \ \ \big((i,j),(i{'}j{'})\big) \in D(I \times J).
\]
\normalfont (iii) 
\itshape If $\Phi \colon \Lambda(X) \To K(X)$ and $\Psi \colon M(X) \To N(X)$, then $\Phi \wedge \Psi \colon 
(\Lambda \wedge M)(X) \To (K \wedge N)(X)$, where, for every $(i,j) \in I \times J$,
\[ (\Phi \wedge \Psi)_{(i,j)} \colon \lambda_0(i) \cap \mu_0(j) \to k_0(i) \cap \nu_0(j),  \]
\[  (\Phi \wedge \Psi)_{(i,j)}(u,w) := \big(\Phi_i(u), \Psi_j(w)\big); \ \ \ \ (u,w) \in \lambda_0(i) \cap \mu_0(j).
\]
\normalfont (iv) 
\itshape If $\Phi \colon \Lambda(X) \To K(X)$ and $\Psi \colon M(X) \To N(X)$, then $\Phi \vee \Psi \colon 
(\Lambda \vee M)(X) \To (K \vee N)(X)$, where, for every $(i,j) \in I \times J$,
\[ (\Phi \wedge \Psi)_{(i,j)} \colon \lambda_0(i) \cup \mu_0(j) \to k_0(i) \cup \nu_0(j),  \]
\[ (\Phi \vee \Psi)_{(i,j)}(z) := \left\{ \begin{array}{ll}
                 \Phi_i(z)    &\mbox{, $z \in \lambda_0(i)$}\\
                 \Psi_j(z)             &\mbox{, $z \in \mu_0(j)$.}
                 \end{array}
          \right.
          \]
\normalfont (v) 
\itshape The following equality holds
\[ \bigcup_{(i,j) \in I \times J}\big(\lambda_0(i) \cap \mu_0(j)\big) =_{\C P(X)} \bigg(\bigcup_{i \in I}\lambda_0(i)\bigg)
 \cap \bigg(\bigcup_{j \in J}\mu_0(j)\bigg). \]
\normalfont (vi) 
\itshape  If $(i_0, j_0) \in I \times J$, the following equality holds
\[ \bigcap_{(i,j) \in I \times J}\big(\lambda_0(i) \cup \mu_0(j)\big) =_{\C P(X)} \bigg(\bigcap_{i \in I}\lambda_0(i)\bigg)
 \cup \bigg(\bigcap_{j \in J}\mu_0(j)\bigg). \]
 \normalfont (vii) 
\itshape If $\Lambda(X)$ covers $X$ and $M(Y)$ covers $Y$, then $(\Lambda \wedge M)(X)$ 
covers $X$.\\[1mm]
\normalfont (viii) 
\itshape Let the inequalities $\neq_I, \neq_J, \neq_X$ and $\neq_Y$ on $I, J, X$ and $Y$, respectively.
If $\Lambda(X)$ is a partition of $X$ and $M(Y)$ is a partition of $Y$, then $(\Lambda \wedge M)(X)$ 
is a partition of $X$.
\end{proposition}

\begin{proof}
We proceed as in the proof of Proposition~\ref{prp: famsubsetsprod1}. 
\end{proof}

Let $M(Y) := (\mu_0, \C Z^Y, \mu_1) \in \Fam(J, Y)$, $(A, i_A^X) \subseteq X$, $(B, i_B^Y) \subseteq Y$, and let
$\Lambda^A(X) := (\lambda_0^A, \C E^{A,X}, \lambda_1^X) \in \Fam(\D 1, X)$ the constant family $A$ of 
subsets of $X$, and $\Lambda^B := (\lambda_0^B, \C E^{B, Y}, \lambda_1^B)$ the constant family $B$ of subsets of $Y$. By 
Propositions~\ref{prp: famsubsetsprod1} and~\ref{prp: famsubsetsprod2} we have that
\begin{align*}
\bigcup_{j \in J}(A \times \mu_0(j)) & := \bigcup_{(i,j) \in \D 1 \times J}A \times \mu_0(j)\\
& := \bigcup_{(i,j) \in \D 1 \times J}\big(\lambda_0^A(i) \times \mu_0(j)\big)\\
& =_{\C P(X \times Y)} \bigg(\bigcup_{i \in \D 1}\lambda_0^A(i)\bigg) \times \bigg(\bigcup_{j \in J}\mu_0(j)\bigg)\\
& =_{\C P(X \times Y)} A \times \bigg(\bigcup_{j \in J}\mu_0(j)\bigg),
\end{align*}
\begin{align*}
\bigcap_{j \in J}(A \times \mu_0(j)) & := \bigcap_{(i,j) \in \D 1 \times J}A \times \mu_0(j)\\
& := \bigcap_{(i,j) \in \D 1 \times J}\big(\lambda_0^A(i) \times \mu_0(j)\big)\\
& =_{\C P(X \times Y)} \bigg(\bigcap_{i \in \D 1}\lambda_0^A(i)\bigg) \times \bigg(\bigcap_{j \in J}\mu_0(j)\bigg)\\
& =_{\C P(X \times Y)} A \times \bigg(\bigcap_{j \in J}\mu_0(j)\bigg),
\end{align*}
\begin{align*}
\bigcup_{j \in J}(B \cap \mu_0(j)) & := \bigcup_{(i,j) \in \D 1 \times J}B \cap \mu_0(j)\\
& := \bigcup_{(i,j) \in \D 1 \times J}\big(\lambda_0^B(i) \cap \mu_0(j)\big)\\
& =_{\C P(Y)} \bigg(\bigcup_{i \in \D 1}\lambda_0^B(i)\bigg) \cap \bigg(\bigcup_{j \in J}\mu_0(j)\bigg)\\
& =_{\C P(Y)} B \cap \bigg(\bigcup_{j \in J}\mu_0(j)\bigg),
\end{align*}
\begin{align*}
\bigcap_{j \in J}(B \cup \mu_0(j)) & := \bigcap_{(i,j) \in \D 1 \times J}B \cup \mu_0(j)\\
& := \bigcap_{(i,j) \in \D 1 \times J}\big(\lambda_0^B(i) \cup \mu_0(j)\big)\\
& =_{\C P(Y)} \bigg(\bigcap_{i \in \D 1}\lambda_0^B(i)\bigg) \cup \bigg(\bigcap_{j \in J}\mu_0(j)\bigg)\\
& =_{\C P(Y)} B \cup \bigg(\bigcap_{j \in J}\mu_0(j)\bigg).
\end{align*}

\begin{definition}\label{def: subprod1}
Let $X, Y, Z \in \D V_0$, $x_0 \in X, y_0 \in Y$, and $R(Z) := (\rho_0, \C E^Z \rho_1) \in \Fam(X \times Y, Z)$.\\[1mm]
\normalfont (i) 
\itshape
If $x \in X$, the $x$-component of $R$\index{$x$-component of a family of subsets on $X \times Y$} is 
the triplet\index{$R^x(Z)$} 
$R^x(Z) := (\rho_0^x, \C E^{x,Z}, \rho_1^x)$, where the assignment routines $\rho_0^x, \rho_1^x$ are 
as in Definition~\ref{def: ac1}, and the dependent operation $\C E^{x,Z} \colon 
\bigcurlywedge_{y \in Y}\D F\big(\rho_0^x(y), Z\big)$ is defined by 
$\C E^{x,Z}_y := \C E_{(x,y)}^Z$, for every $y \in Y$.\\[1mm]
\normalfont (ii) 
\itshape If $y \in Y$, the $y$-component of $R$\index{$y$-component of a family of subsets on $X \times Y$}
is the triplet\index{$R^y(Z)$} 
$R^y(Z) := (\rho_0^y, \C E^{y,Z}, \rho_1^y)$, where the assignment routines $\rho_0^y, \rho_1^y$ are 
as in Definition~\ref{def: ac1}, and the dependent operation $\C E^{y,Z} \colon 
\bigcurlywedge_{x \in X}\D F\big(\rho_0^y(x), Z\big)$ is defined by 
$\C E^{y,Z}_x := \C E_{(x,y)}^Z$, for every $x \in Y$.\\[1mm]
\normalfont (iii) 
\itshape
Let $\bigcup^1R := (\bigcup^1 \rho_0, \big(\bigcup^1 \C E\big)^Z, \bigcup^1 \rho_1)$, where\index{$\bigcup^1R$} 
$\bigcup^1 \rho_0 : X \sto \D V_0$, 
\[\bigcup^1 \rho_1 \colon \bigcurlywedge_{(x, x{'}) \in D(X)}\D F\bigg(\big(\bigcup^1 \rho_0\big)(x), 
\big(\bigcup^1 \rho_0\big)(x{'})\bigg), 
 \big(\bigcup^1 \C E\big)^Z \colon \bigcurlywedge_{x \in X}\D F\big(\big(\bigcup^1\rho_0\big)(x), Z\big)
\ \mbox{are defined by} \]
\[ \bigg(\bigcup^1 \rho_0\bigg) (x) := \bigcup_{y \in Y}\rho_0^x(y) := \bigcup_{y \in Y}\rho_0(x, y); 
\ \ \ \ x \in X, \] 
\[ \bigg(\bigcup^1 \rho_1\bigg) (x, x{'}) := \bigg(\bigcup^1 \rho_1\bigg)_{xx{'}} \colon 
\bigcup_{y \in Y}\rho_0 (x, y) \to \bigcup_{y \in Y}\rho_0 (x{'}, y); \ \ \ \ 
(x, x{'}) \in D(X), \]
\[ \bigg(\bigcup^1 \rho_1\bigg)_{xx{'}}(y, u) := \big(y, \rho_{(x, y)(x{'},y)}(u)\big); \ \ \ \ (y, u) \in  
\bigcup_{y \in Y}\rho_0 (x, y), \]
\[ \big(\bigcup^1 \C E\big)^Z_x (y, u) := \C E_{(x,y)}^Z(u); \ \ \ \ x \in X, \ (y, u) \in \bigcup_{y \in Y}\rho_0(x,y). \]
\normalfont (iv) 
Let $\bigcup^2R := (\bigcup^2 \rho_0, \big(\bigcup^2 \C E\big)^Z, \bigcup^2 \rho_1)$, where\index{$\bigcup^2R$} 
$\bigcup^2 \rho_0 : Y \sto \D V_0$, 
\[\bigcup^2 \rho_1 \colon \bigcurlywedge_{(y, y{'}) \in D(Y)}\D F\bigg(\big(\bigcup^2 \rho_0\big)(x),
\big(\bigcup^2 \rho_0\big)(x{'})\bigg), 
 \big(\bigcup^2 \C E\big)^Z \colon \bigcurlywedge_{y \in Y}\D F\big(\big(\bigcup^2\rho_0\big)(y), Z\big)
\ \mbox{are defined by} \]
\[ \bigg(\bigcup^2 \rho_0\bigg) (y) := \bigcup_{x \in X}\rho_0^y(x) := \bigcup_{x \in X}\rho_0(x, y); \ \ \ \ y \in Y, \] 
\[ \bigg(\bigcup^2 \rho_1\bigg) (y, y{'}) := \bigg(\bigcup^2 \rho_1\bigg)_{yy{'}} \colon 
\bigcup_{x \in X}\rho_0 (x, y) \to \bigcup_{x \in X}\rho_0 (x{'}, y); \ \ \ \ 
(y, y{'}) \in D(Y), \]
\[ \bigg(\bigcup^2 \rho_1\bigg)_{yy{'}}(x, w) := \big(x, \rho_{(x, y)(x{'},y)}(w)\big); \ \ \ \ (x, w) \in  
\bigcup_{x \in X}\rho_0 (x, y), \]
\[ \big(\bigcup^2 \C E\big)^Z_y (x, w) := \C E_{(x,y)}^Z(w); \ \ \ \ y \in Z, \ (x, w) \in \bigcup_{x \in X}\rho_0(x,y). \]
\normalfont (v) 
\itshape
Let $\bigcap^1R := (\bigcap^1 \rho_0, \big(\bigcap^1 \C E\big)^Z, \bigcap^1 \rho_1)$, where\index{$\bigcap^1R$} 
$\bigcap^1 \rho_0 : X \sto \D V_0$, 
\[\bigcap^1 \rho_1 \colon \bigcurlywedge_{(x, x{'}) \in D(X)}\D F\bigg(\big(\bigcap^1 \rho_0\big)(x), 
\big(\bigcap^1 \rho_0\big)(x{'})\bigg), 
 \big(\bigcap^1 \C E\big)^Z \colon \bigcurlywedge_{x \in X}\D F\big(\big(\bigcap^1\rho_0\big)(x), Z\big)
\ \mbox{are defined by} \]
\[ \bigg(\bigcap^1 \rho_0\bigg) (x) := \bigcap_{y \in Y}\rho_0^x(y) := \bigcap_{y \in Y}\rho_0(x, y); \ \ \ \ x \in X, \] 
\[ \bigg(\bigcap^1 \rho_1\bigg) (x, x{'}) := \bigg(\bigcap^1 \rho_1\bigg)_{xx{'}} \colon 
\bigcap_{y \in Y}\rho_0 (x, y) \to \bigcap_{y \in Y}\rho_0 (x{'}, y); \ \ \ \ 
(x, x{'}) \in D(X), \]
\[ \bigg[\bigg(\bigcap^1 \rho_1\bigg)_{xx{'}}(\Phi)\bigg]_y := \rho_{(x, y)(x{'},y)}(\Phi_y); \ \ \ \ \Phi \in  
\bigcap_{y \in Y}\rho_0 (x, y), \]
\[ \big(\bigcap^1 \C E\big)^Z_x (\Phi) := \C E_{(x,y_0)}^Z(\Phi_{y_0}); \ \ \ \Phi \in \bigcap_{y \in Y}\rho_0(x,y). \]
\normalfont (vi) 
\itshape
Let $\bigcap^2R := (\bigcap^2 \rho_0, \big(\bigcap^2 \C E\big)^Z, \bigcap^2 \rho_1)$, where\index{$\bigcap^2R$} 
$\bigcap^2 \rho_0 : X \sto \D V_0$, 
\[\bigcap^2 \rho_1 \colon \bigcurlywedge_{(y, y{'}) \in D(Y)}\D F\bigg(\big(\bigcap^2 \rho_0\big)(y), 
\big(\bigcap^2 \rho_0\big)(y{'})\bigg), 
 \big(\bigcap^2 \C E\big)^Z \colon \bigcurlywedge_{y \in Z}\D F\big(\big(\bigcap^2\rho_0\big)(y), Z\big)
\ \mbox{are defined by} \]
\[ \bigg(\bigcap^2 \rho_0\bigg) (y) := \bigcap_{x \in X}\rho_0^y(x) := \bigcap_{x \in X}\rho_0(x, y); \ \ \ \ y \in Y, \] 
\[ \bigg(\bigcap^2 \rho_1\bigg) (y, y{'}) := \bigg(\bigcap^2 \rho_1\bigg)_{yy{'}} \colon 
\bigcap_{x \in X}\rho_0 (x, y) \to \bigcap_{x \in X}\rho_0 (x, y{'}); \ \ \ \ 
(y, y{'}) \in D(Y), \]
\[ \bigg[\bigg(\bigcap^2 \rho_1\bigg)_{yy{'}}(\Phi)\bigg]_x := \rho_{(x, y)(x,y{'})}(\Phi_x); \ \ \ \ \Phi \in  
\bigcap_{x \in X}\rho_0 (x, y), \]
\[ \big(\bigcap^2 \C E\big)^Z_y (\Phi) := \C E_{(x_0,y)}^Z(\Phi_{x_0}); \ \ \ \Phi \in \bigcap_{x \in X}\rho_0(x,y). \] 
\end{definition} 
 
 Clearly, 
$R^y(Z), \bigcup^1R(Z), \bigcap^1R(Z)  \in \Fam(X, Z)$ and  
$R^x(Z), \bigcup^2R(Z), \bigcap^2R(Z) \in \Fam(Y, Z)$.

\begin{proposition}\label{prp: newfamilysubmaps3}
Let $X, Y, Z \in \D V_0$, $R(Z) := (\rho_0, \C E^Z \rho_1), S(Z):= (\sigma_0, \C A^Z, \sigma_1)
\in \Fam(X \times Y, Z)$, and
$\Phi \colon R(Z) \To S(Z)$.\\[1mm]
\normalfont (i) 
\itshape Let $\Phi^x \colon \bigcurlywedge_{y \in Y}\D F\big(\rho_0^x(y), \sigma_0^x(y)\big)$, where 
$\Phi^x_y := \Phi_{(x,y)} \colon \rho_0^x(y) \to \sigma^x_0(y)$.\\[1mm]
\normalfont (ii) 
\itshape Let $\Phi^y \colon \bigcurlywedge_{x \in X}\D F\big(R^y(x), S^y(x)\big)$, where 
$\Phi^y_x := \Phi_{(x,y)} \colon \rho_0^y(x) \to \sigma^y_0(x)$.\\[1mm]
\normalfont (iii) 
\itshape Let $\bigcup^1 \Phi \colon \bigcurlywedge_{x \in X}\D F\big(\big(\bigcup^1\rho_0\big)(x),
\big(\bigcup^1\sigma_0\big)(x)\big)$, where, for every $x \in X$, we define\index{$\bigcup^1 \Phi$} 
\[ \bigg(\bigcup^1 \Phi\bigg)_x \colon \sum_{y \in Y}\rho_0^y(x) \to \bigcup_{y \in Y}\sigma^x_0(y) \]
\[ \bigg(\bigcup^1 \Phi\bigg)_x(y, u) := \big(y, \Phi_{(x,y)}(u)\big); \ \ \ \ (y, u) \in \bigcup_{y \in Y}\rho_0(x, y). \]
\normalfont (iv) 
\itshape Let $\bigcup^2 \Phi \colon \bigcurlywedge_{y \in Y}\D F\big(\big(\bigcup^2\rho_0\big)(y),
\big(\bigcup^2\sigma_0\big)(y)\big)$, where, for every $y \in Y$, we define\index{$\bigcup^2 \Phi$} 
\[ \bigg(\bigcup^2 \Phi\bigg)_y \colon \bigcup_{x \in X}\rho_0^x(y) \to \bigcup_{x \in X}\sigma^y_0(x) \]
\[ \bigg(\bigcup^2 \Phi\bigg)_y(x, w) := \big(x, \Phi_{(x,y)}(w)\big); \ \ \ \ (x, w) \in \bigcup_{x \in X}\rho(x, y). \]
\normalfont (v) 
\itshape Let $\bigcap^1 \Phi \colon \bigcurlywedge_{x \in X}\D F\big(\big(\bigcap^1\rho_0\big)(x),
\big(\bigcap^1\sigma_0\big)(x)\big)$, where, for every $x \in X$, we define\index{$\bigcap^1 \Phi$} 
\[ \bigg(\bigcap^1 \Phi\bigg)_x \colon \bigcap_{y \in Y}\rho_0^y(x) \to \bigcap_{y \in Y}\sigma^x_0(y) \]
\[ \bigg[\bigg(\bigcap^1 \Phi\bigg)_x(\Theta)\bigg]_y := \Phi_{(x,y)}\big(\Theta_y))\big); \ \ \ \ \Theta \in
\bigcap_{y \in Y}\rho_0(x, y). \]
\normalfont (vi) 
\itshape Let $\bigcap^2 \Phi \colon \bigcurlywedge_{y \in Y}\D F\big(\big(\bigcap^2\rho_0\big)(y),
\big(\bigcap^2\sigma_0\big)(y)\big)$, where, for every $y \in Y$, we define\index{$\bigcap^2 \Phi$} 
\[ \bigg(\bigcap^2 \Phi\bigg)_y \colon \bigcap_{x \in X}\rho_0^x(y) \to \bigcap_{x \in X}\sigma^y_0(x) \]
\[ \bigg[\bigg(\bigcap^2 \Phi\bigg)_y(\Theta)\bigg]_x := \Phi_{(x,y)}\big(\Theta_x)\big); \ \ \ \ \Theta \in 
\bigcap_{x \in X}\rho_0(x, y). \]
Then we have that $\Phi^x \colon R^x(Z) \To S^x(Z)$ and $\Phi^y \colon R^y(Z) \To S^y(Z)$ and
$\bigcup^1 \Phi \colon \big(\bigcup^1 R\big)(Z) \To \big(\bigcup^1 S\big)(Z)$ and
$\bigcup^2 \Phi \colon \big(\bigcup^2 R\big)(Z) \To \big(\bigcup^2 S\big)(Z)$ and
$\bigcap^1 \Phi \colon \big(\bigcap^1 R\big)(Z) \To \big(\bigcap^1 S\big)(Z)$ and
$\bigcap^2 \Phi \colon \big(\bigcap^2 R\big)(Z) \To \big(\bigcap^2 S\big)(Z)$.
\end{proposition}

\begin{proof}
We proceed similarly to the proof of Proposition~\ref{prp: newfamilymaps3}. 
\end{proof}

\begin{proposition}\label{prp: unionunion}
If $R := (\rho_0, \rho_1) \in \Fam(X \times Y, Z)$, the following equalities hold.
\[ \bigcup_{x \in X}\bigcup_{y \in Y}\rho_0(x,y) =_{\C P(Z)} \bigcup_{y \in Y}\bigcup_{x \in X}\rho_0(x,y), \]
\[ \bigcap_{x \in X}\bigcap_{y \in Y}\rho_0(x,y) =_{\C P(Z)} \bigcap_{y \in Y}\bigcap_{x \in X}\rho_0(x,y). \]
\end{proposition}

\begin{proof}
The proof is straightforward. 
\end{proof}

\section{The semi-distributivity of $\bigcap$ over $\bigcup$}
\label{sec: bigcapcup}

Section~\ref{sec: famsofsubsetsprod} is the ``internal'' analogue to section~\ref{sec: famoverprod}, 
as the presentation of the families of subsets over products follows the presentation of the families 
of sets over products.
The distributivity of $\bigcap$ over $\bigcup$ though, cannot be approached as the distributivity of 
$\Sigma$ over $\prod$, as the crucial Lemma~\ref{lem: ac3} depends on the fact that the operation $\pr_1^{R^x}$ is a
function, something which is not the case, as we have already explained in section~\ref{sec: union},
when the totality of the exterior union is equipped with the equality of the 
interior union. 

\begin{definition}\label{def: subcomposition}
If $\Lambda(X) := (\lambda_0, \C E^X, \lambda_1) \in \Fam(I,X)$ and $h \colon J \to I$, the composition family 
\index{$\Lambda(X) \circ h$} of $\Lambda(X)$\index{composition family of $\Lambda(X)$ with $h$}
with $h$ is the triplet $\Lambda(X) \circ h := (\lambda_0 \circ h, \C E^X \circ h, \lambda_1 \circ h)$, where 
$\lambda_0 \circ h \colon J \sto \D V_0$ and $\lambda_1 \circ h \colon \bigcurlywedge_{(j,j{'}) \in D(J)}
\D F\big(\lambda_0(h(j)), \lambda_0(h(j{'}))\big)$ are given in Definition~\ref{def: newfamsofsets1}$($iii$)$, and 
the dependent operation $\C E^X \circ h \colon  \bigcurlywedge_{j \in J}\D F\big(\lambda_0(h(j)), X\big)$ is defined by
$(\C E^X \circ h)_j := \C E^X_{h(j)}$, for every $j \in J$.
\end{definition}

Clearly, $\Lambda(X) \circ h \in \Fam(J, X)$. To formulate the distributivity of 
$\bigcap$ over $\bigcup$ in the language of $\BST$ we need to introduce a family of subsets $P(I)$ 
of the index-set $I$ of a given family of subsets of a set $X$. 
Throughout this section let the following data:\\[2mm]
\textbf{(a)} $\Lambda(X) := (\lambda_0, \C E^X, \lambda_1) \in \Fam(I,X)$.\\[1mm]
\textbf{(b)} $(K, =_K, \neq_K)$ is a set, and $k_0 \in K$.\\[1mm]
\textbf{(c)} $P(I) := (p_0, \C Z^I, p_1) \in \Fam(K, I)$.\\[1mm]
\textbf{(d)} $\Lambda(X) \circ \C Z^I_k := (\lambda_0 \circ \C Z^I_k, \C E^X \circ \C Z^I_k, 
\lambda_1 \circ \C Z^I_k) \in 
\Fam(p_0(k), X)$, for every $k \in K$.\\[1mm] 
\textbf{(e)} $T := \bigcap_{k \in K}p_0(k)$.

\begin{proposition}\label{prp: subdistr1}
$N(K) := (\nu_0, \C N^X, \nu_1) \in \Fam(K, X)$, where $\nu_0 \colon K \sto \D V_0$ is defined by 
\[\nu_0(k) := \bigcup_{j \in p_0(k)}(\lambda_0 \circ \C Z^I_k)(j) := \bigcup_{j \in p_0(k)}\lambda_0\big(\C Z^I_k(j)\big);
 \ \ \  \ k \in K,
\]
and $\C N^X \colon \bigcurlywedge_{k \in K}\D F(\nu_0(k), X)$,  
$\nu_1 \colon \bigcurlywedge_{(k,k{'}) \in D(K)}\D F\big(\nu_0(k), \nu_0(k{'})\big)$ are defined by
\[ \C N^X_k \colon \bigg(\bigcup_{j \in p_0(k)}\lambda_0\big(Z^I_k(j)\big)\bigg) \eto X, \ \ 
\C N^X_k(j, u) := \C E^X_{\C Z^I_k(j)}(u); \ \ \ \ j \in p_0(k), \ u \in \lambda_0\big(Z^I_k(j)\big), \]
\[ \nu_1(k, k{'}) := \nu_{kk{'}} \colon \bigcup_{j \in p_0(k)}\lambda_0\big(Z^I_k(j)\big) \to 
\bigcup_{j \in p_0(k{'})}\lambda_0\big(Z^I_{k{'}}(j)\big), \]
\[ \nu_{kk{'}}(j, u) := \big(p_{kk{'}}(j), \lambda_{Z^I_k(j) Z^I_{k{'}}(p_{kk{'}}(j))}(u)\big); \ \ \ \ 
j \in p_0(k), \ u \in \lambda_0\big(Z^I_k(j)\big). 
\]
\end{proposition}

\begin{proof}
The operation $\C N^X_k$ is an embedding, since by Definition~\ref{def: subcomposition} 
\[ (j, u) =_{\mathsmaller{\bigcup_{j \in p_0(k)}\lambda_0(Z^I_k(j))}} (j{'}, u{'}) :\TOT
\C E^X_{\C Z^I_k(j)}(u) =_X \C E^X_{\C Z^I_k(j{'})}(u{'}) :\TOT \C N^X_k(j, u) =_X \C N^X_k(j{'}, u{'}). \] 
Let $k =_K k{'}$, $j \in p_0(k)$ and $u \in \lambda_0\big(Z^I_k(j)\big)$. By the commutativity of the left inner diagrams
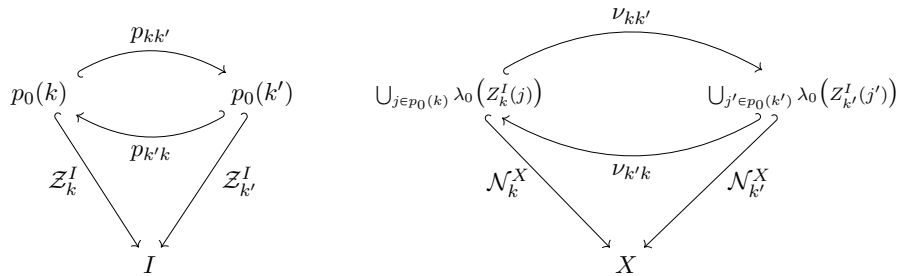
\begin{center}
\resizebox{12cm}{!}{%
\begin{tikzpicture}

\node (E) at (0,0) {$p_0(k)$};
\node[right=of E] (B) {};
\node[right=of B] (F) {$p_0(k{'})$};
\node[below=of B] (C) {};
\node[below=of C] (A) {$I$};

\node[right=of F] (K) {$\mathsmaller{\bigcup_{j \in p_0(k)}\lambda_0\big(Z^I_k(j)\big)}$};
\node[right=of K] (L) {};
\node[right=of L] (M) {$\mathsmaller{\bigcup_{j{'} \in p_0(k{'})}\lambda_0\big(Z^I_{k{'}}(j{'})\big)}$};
\node[below=of L] (N) {};
\node[below=of N] (X) {$X$};

\draw[left hook->,bend left] (E) to node [midway,above] {$p_{kk{'}}$} (F);
\draw[left hook->,bend left] (F) to node [midway,below] {$p_{k{'}k}$} (E);
\draw[right hook->] (E)--(A) node [midway,left] {$\C Z^I_k \ $};
\draw[left hook->] (F)--(A) node [midway,right] {$ \ \C Z^I_{k{'}}$};

\draw[left hook->,bend left] (K) to node [midway,above] {$\nu_{kk{'}}$} (M);
\draw[left hook->,bend left] (M) to node [midway,below] {$\nu_{k{'}k}$} (K);
\draw[right hook->] (K)--(X) node [midway,left] {$\C N^X_k \ $};
\draw[left hook->] (M)--(X) node [midway,right] {$ \ \C N^X_{k{'}}$};

\end{tikzpicture}
}
\end{center}
we have that $\C Z^I_{k{'}}\big(p_{kk{'}}(j)\big) =_I \C Z^I_k(j)$. Hence $\lambda_{Z^I_k(j) Z^I_{k{'}}(p_{kk{'}}(j))}
\colon \lambda_0\big(\C Z^I_k(j)\big) \to \lambda_0\big(\C Z^I_{k{'}}\big(p_{kk{'}}(j)\big)\big)$ and
$\nu_{kk{'}}(j, u)$ is well defined. Next we show that the above right inner diagrams commute. If 
\[ j{'} := p_{kk{'}}(j) \ \ \& \ \ i :=  Z^I_k(j) \ \ \& \ \ i{'} := \C Z^I_{k{'}}(j{'}), \ \ \ \ \mbox{then} \]
\[ \C N^X_{k{'}}\big(\nu_{kk{'}}(j, u)\big) := \C N^X_{k{'}}\big(j{'}, \lambda_{ii{'}}(u)\big) := 
 \C E^X_{i{'}}\big( \lambda_{ii{'}}(u)\big) =_X \C E^X_i(u) := \C N^X_{k}(j,u), \]
using the commutativity of the following diagram
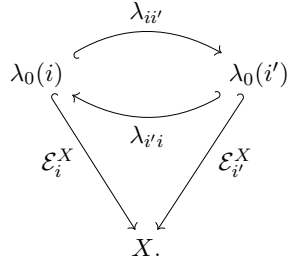
\begin{center}
\resizebox{4cm}{!}{%
\begin{tikzpicture}

\node (E) at (0,0) {$\lambda_0(i)$};
\node[right=of E] (B) {};
\node[right=of B] (F) {$\lambda_0(i{'})$};
\node[below=of B] (C) {};
\node[below=of C] (A) {$X$.};

\draw[left hook->,bend left] (E) to node [midway,above] {$\lambda_{ii{'}}$} (F);
\draw[left hook->,bend left] (F) to node [midway,below] {$\lambda_{i{'}i}$} (E);
\draw[right hook->] (E)--(A) node [midway,left] {$\C E^X_i \ $};
\draw[left hook->] (F)--(A) node [midway,right] {$ \ \C E^X_{i{'}}$};
 
\end{tikzpicture}
}
\end{center} 
For the other above right inner diagram we proceed similarly. Clearly, $\nu_{kk}(j, u) := (j, u)$.
\end{proof}

\begin{proposition}\label{prp: subdistr2}
If $\tau \in T$, then 
$\tau(X) := (\tau_0, \C T^X, \tau_1) \in \Fam(K, X)$, where $\tau_0 \colon K \sto \D V_0$ is defined by 
$\tau_0(k) := \lambda_0\big(\C Z^I_k(\tau_k)\big)$, for every $k \in K$, and the dependent operations 
$\C T^X \colon \bigcurlywedge_{k \in K}\D F(\tau_0(k), X)$,  
$\tau_1 \colon \bigcurlywedge_{(k,k{'}) \in D(K)}\D F\big(\tau_0(k), \tau_0(k{'})\big)$ are defined by
\[ \C T^X_k \colon\lambda_0\big( \C Z^I_k(\tau_k)\big) \eto X, \ \ 
\C T^ X_k := \C E^X_{\C Z^I_k(\tau_k)},  \]
\[ \tau_1(k, k{'}) := \tau_{kk{'}} \colon \lambda_0\big(\C Z^I_k(\tau_k)\big) \to 
\lambda_0\big(\C Z^I_{k{'}}(\tau_{k{'}})\big), \ \ \ \ \tau_{kk{'}} := \lambda_{Z^I_k(\tau_k) Z^I_{k{'}}(\tau_{k{'}})}.
\]
\end{proposition}

\begin{proof}
What we want follows in a straightforward way from the fact that $\Lambda(X) \in \Fam(I, X)$.
\end{proof}

\begin{proposition}\label{prp: subdistr3}
$\Xi (X) := (\xi_0, \C H^X, \xi_1) \in \Fam(T, X)$, where $\xi_0 \colon T \sto \D V_0$ is defined by 
\[ \xi_0(\tau) := \bigcap_{k \in K}\tau_0(k) := \bigcap_{k \in K}\lambda_0\big(\C Z_k^I(\tau_k)\big); \ \ \ \ \tau \in T,\]
and the dependent operations 
$\C H^X \colon \bigcurlywedge_{\tau \in T}\D F(\xi_0(\tau), X)$,  
$\xi_1 \colon \bigcurlywedge_{(\tau,\tau{'}) \in D(T)}\D F\big(\xi_0(\tau), \xi_0(\tau{'})\big)$ are defined, respectively, by
$\C H^X_{\tau} \colon \bigg(\bigcap_{k \in K}\tau_0(k)\bigg) \eto X$, where $\C H_{\tau}^X := e_{\mathsmaller{\bigcap}}^{\tau(X)}$, for every $\tau \in T$,
\[ \xi_1(\tau, \tau{'}) := \xi_{\tau \tau{'}} \colon \bigcap_{k \in K}\lambda_0\big(\C Z_k^I(\tau_k)\big) \to 
\bigcap_{k \in K}\lambda_0\big(\C Z_k^I(\tau{'}_k)\big),\]
\[ \Phi \mapsto \xi_{\tau \tau{'}}(\Phi); \ \ \ \ \big[\xi_{\tau \tau{'}}(\Phi)\big]_k := \lambda_{\C Z_k^I(\tau_k)\C Z_k^I(\tau{'}_k)}(\Phi_k); \ \ \ \ \Phi \colon \bigcap_{k \in K}\lambda_0\big(\C Z_k^I(\tau_k)\big), \   k \in K.\]  
\end{proposition}

\begin{proof}
If $\tau \in T$, then by the definition of the embedding $e_{\mathsmaller{\bigcap}}^{\tau(X)}$ we get
\[ \C H_{\tau}^X(\Phi) := \C T_{k_0}^X(\Phi_{k_0}) := \C E_{\C Z_{k_0}^I(\tau_{k_0})}^X(\Phi_{k_0}); \ \ \ \
\Phi \colon \bigcap_{k \in K}\lambda_0\big(\C Z_k^I(\tau_k)\big).\]
$\C H^X_{\tau}$ is an embedding. Next we show that $\xi_{\tau \tau{'}}(\Phi) \in 
\bigcap_{k \in K}\lambda_0\big(\C Z_k^I(\tau{'}_k)\big)$. As $\Phi \colon 
\bigcap_{k \in K}\lambda_0\big(\C Z_k^I(\tau_k)\big)$, 
\begin{align*}
\C T_k^X\big(\big[\xi_{\tau \tau{'}}(\Phi)\big]_k\big) & := \C E^X_{\C Z^I_k(\tau{'}_k)}\bigg(\big[\xi_{\tau \tau{'}}(\Phi)\big]_k\bigg)\\
& := \C E^X_{\C Z^I_k(\tau{'}_k)}\bigg(\lambda_{\C Z_k^I(\tau_k)\C Z_k^I(\tau{'}_k)}(\Phi_k)\bigg)\\
& =_X \C E^X_{\C Z^I_k(\tau_k)}(\Phi_k)\\
& =_X \C E^X_{\C Z^I_l(\tau_l)}(\Phi_l)\\
& =_X \C E^X_{\C Z^I_l(\tau{'}_l)}\bigg(\lambda_{\C Z_k^I(\tau_l)\C Z_k^I(\tau{'}_l)}(\Phi_l)\bigg)\\
& := \C T_l^X\big(\big[\xi_{\tau \tau{'}}(\Phi)\big]_l\big),
\end{align*}
for every $k, l \in K$. Similarly we show that $\xi_{\tau \tau{'}}$ is a function. If $\tau =_T \tau{'}$, then
\begin{center}
\begin{tikzpicture}

\node (E) at (0,0) {$\bigcap_{k \in K}\lambda_0(\C Z_k^I(\tau_k)))$};
\node[right=of E] (B) {};
\node[right=of B] (F) {$\bigcap_{k \in K}\lambda_0(\C Z^I_k(\tau{'}_k))$};
\node[below=of B] (C) {};
\node[below=of C] (A) {$X$};

\draw[left hook->,bend left] (E) to node [midway,above] {$\xi_{\tau \tau{'}}$} (F);
\draw[left hook->,bend left] (F) to node [midway,below] {$\xi_{\tau{'}\tau}$} (E);
\draw[right hook->] (E)--(A) node [midway,left] {$\C H^X_{\tau} \ $};
\draw[left hook->] (F)--(A) node [midway,right] {$ \ \C H^X_{\tau{'}}$};
 
\end{tikzpicture}
\end{center} 
\begin{align*}
\C H^X_{\tau{'}}\big(\xi_{\tau \tau{'}}(\Phi)\big) & := \C E_{\C Z_{k_0}^I(\tau{'}_{k_0})}^X\bigg(\big[\xi_{\tau \tau{'}}(\Phi)\big]_{k_0}\bigg)\\
& := \C E_{\C Z_{k_0}^I(\tau{'}_{k_0})}^X\bigg( \lambda_{\C Z_{k_0}^I(\tau_{k_0})\C Z_{k_0}^I(\tau{'}_{k_0})}(\Phi_{k_0})  \bigg)\\
& =_X \C E_{\C Z_{k_0}^I(\tau_{k_0})}^X\big(\Phi_{k_0}\big)\\
& := \C H_{\tau}^X(\Phi).\qedhere
\end{align*}
\end{proof}

The set 
\[ W:= \bigcap_{k \in K}\nu_0(k) := \bigcap_{k \in K}\bigg[\bigcup_{j \in p_0(k)}\lambda_0\big(\C Z^I_k(j)\big)\bigg]\]
is embedded into $X$ through the map $e_{\mathsmaller{\bigcap}}^{N(X)}$, where $e_{\mathsmaller{\bigcap}}^{N(X)}(A) := \C N^X_{k_0}(A_{k_0})$, for every $A \in \bigcap_{k \in K}\nu_0(k)$. By definition, if $A \colon \bigcurlywedge_{k \in K}\nu_0(k)$, then
\[ A \in \bigcap_{k \in K}\nu_0(k) \TOT \forall_{k,l \in K}\big(\C N^X_k(A_k) =_X \C N^X_l(A_l)\big), \]
\[ A_k \in \bigcup_{j \in p_0(k)}\lambda_0\big(\C Z^I_k(j)\big), \ \ \ \ \mbox{i.e.,} \ \ \ \ A_k := (j, u), \ j \in p_0(k), \ u \in  \lambda_0\big(\C Z^I_k(j)\big),\]
\[ \C N^X_k(A_k) := \C N^X_k(j, u) := \C E^X_{\C Z^I_k(j)}(u),\]
\[ A =_{\mathsmaller{\bigcap_{k \in K}\nu_0(k)}} B :\TOT \C N^X_{k_0}(A_{k_0}) =_X \C N^X_{k_0}(B_{k_0}).\]
The set 
\[ V := \bigcup_{\tau \in T}\xi_0(\tau) := \bigcup_{\tau \in T}\bigg[\bigcap_{k \in K}\lambda_0\big(\C Z^I_k(\tau_k)\big)\bigg]\]
is embedded into $X$ through the map $e_{\mathsmaller{\bigcup}}^{\Xi(X)}$, where 
\[e_{\mathsmaller{\bigcup}}^{\Xi(X)}(\tau, \Phi) := \C H^X_{\tau}(\Phi) := \C E^X_{\C Z^I_{k_0}(\tau_{k_0})}(\Phi_{k_0}); \ \ \ \ (\tau, \Phi)  \in \bigcup_{\tau \in T}\xi_0(\tau),\]
\[(\tau, \Phi) =_{\mathsmaller{\bigcup_{\tau \in T}\xi_0(\tau)}} (\tau{'}, \Phi{'}) :\TOT \C H^X_{\tau}(\Phi) =_X \C H^X_{\tau{'}}(\Phi{'}).\]

\begin{proposition}[Semi-distributivity of $\bigcap$ over $\bigcup$]\label{prp: subdistr4}
$\big(V, e_{\mathsmaller{\bigcup}}^{\Xi(X)}\big) \subseteq \big(W, e_{\mathsmaller{\bigcap}}^{N(X)}\big)$.
\end{proposition}

\begin{proof}
Let the operation $\theta \colon V \sto W$, defined by
\[ \theta(\tau, \Phi) := A^{(\tau, \phi)}, \ \ \ \ (\tau, \Phi) \in V,\]
\[ A^{(\tau, \Phi)}_k := (\tau_k, \Phi_k); \ \ \ \ k \in K.\]
By definition $\tau_k \in p_0(k)$ and $\Phi_k \in \lambda_0\big(\C Z^I_k(\tau_k)\big)$. We show that $\theta$ is well-defined i.e., $A^{(\tau, \Phi)} \in \bigcap_{k \in K}\nu_0(k)$. If $k, l \in K$, by the above unfolding of $A \in W$ we need to show that 
\[ \C E^X_{\C Z^I_k(\tau_k)}(\Phi_k) =_X \C E^X_{\C Z^I_k(\tau_l)}(\Phi_l),\]
which follows immediately from the unfolding of the membership $\Phi \in \bigcap_{k \in K}\lambda_0\big(\C Z^I_k(\tau_k)\big)$. If 
\[(\tau, \Phi) =_{\mathsmaller{\bigcup_{\tau \in T}\xi_0(\tau)}} (\tau{'}, \Phi{'}) :\TOT \C H^X_{\tau}(\Phi) =_X \C H^X_{\tau{'}}(\Phi{'}) :\TOT \C E_{\C Z_{k_0}^I(\tau_{k_0})}^X(\Phi_{k_0}) =_X \C E_{\C Z_{k_0}^I(\tau_{k_0})}^X(\Phi{'}_{k_0}),\]
\[A^{(\tau, \Phi)} =_{\mathsmaller{\bigcap_{k \in K}\nu_0(k)}} A^{(\tau{'}, \Phi{'})} :\TOT \C N^X_{k_0}(A_{k_0}^{(\tau, \Phi)}) =_X \C N^X_{k_0}(A_{k_0}^{(\tau{'}, \Phi{'})}) \]
\[ :\TOT \C N^X_{k_0}(A_{k_0}^{(\tau_k, \Phi_k)}) =_X \C N^X_{k_0}(A_{k_0}^{(\tau{'}_k, \Phi{'}_k)}) :\TOT \C E_{\C Z_{k_0}^I(\tau_{k_0})}^X(\Phi_{k_0}) =_X \C E_{\C Z_{k_0}^I(\tau_{k_0})}^X(\Phi{'}_{k_0}),\]
hence $\theta$ is a function. The commutativity of the following diagram is shown by the equalities
\begin{center}
\begin{tikzpicture}

\node (E) at (0,0) {$\bigcup_{\tau \in T}\bigg[\bigcap_{k \in K}\lambda_0\big(\C Z^I_k(\tau_k)\big)\bigg]$};
\node[right=of E] (B) {};
\node[right=of B] (F) {$\bigcap_{k \in K}\bigg[\bigcup_{j \in p_0(k)}\lambda_0\big(\C Z^I_k(j)\big)\bigg]$};
\node[below=of B] (C) {};
\node[below=of C] (A) {$X$};

\draw[left hook->,bend left] (E) to node [midway,above] {$\theta$} (F);
\draw[right hook->] (E)--(A) node [midway,left] {$e_{\mathsmaller{\bigcup}}^{\Xi(X)} \  \ $};
\draw[left hook->] (F)--(A) node [midway,right] {$ \  e_{\mathsmaller{\bigcap}}^{N(X)}$};
 
\end{tikzpicture}
\end{center} 
\[ e_{\mathsmaller{\bigcap}}^{N(X)}\big(\theta(\tau, \Phi)\big) := \C N^X_{k_0}\big(\tau_{k_0}, \Phi_{k_0}) := \C E_{\C Z_{k_0}^I(\tau_{k_0})}^X(\Phi_{k_0}) := e_{\mathsmaller{\bigcup}}^{\Xi(X)}(\tau, \Phi). \qedhere \]
\end{proof}

For the converse inclusion see Note~\ref{not: converseinclusion}.

\section{Sets of subsets}
\label{sec: setofsubsets}

\begin{definition}\label{def: setofsubsets}
If $I, X \in \D V_0$, a \textit{set of subsets}\index{set of subsets} of $X$ indexed by $I$, or an $I$-\textit{set 
of subsets}\index{$I$-set of subsets} of $X$, is triplet $\Lambda(X) := (\lambda_0, \C E^X, \lambda_1) \in \Fam(I, X)$
such that the following condition is satisfied:
\[ Q(\Lambda(X)) :\TOT \forall_{i, j \in I}\big(\lambda_0(i) =_{\C P(X)} \lambda_0(j) \To i =_I j\big). \]
Let $\Set(I, X)$ be their totality, equipped with the canonical equality on $\Fam(I, X)$. 
\end{definition}

\begin{remark}\label{rem: Isetssubsets1}
If $\Lambda(X) \in \Set(I, X)$ and $M(X) \in \Fam(I, X)$ such that $\Lambda(X) =_{\Fam(I, X)} M(X)$, then 
$M(X) \in \Set(I, X)$. 
\end{remark}

\begin{proof}
Let $\Phi \colon \Lambda(X) \To M(X)$ and $\Psi \colon M(X) \To \Lambda(X)$ such that 
$(\Phi, \Psi) \colon \Lambda(X) =_{\Fam(I, X)} M(X)$. Let also $(f, g) \colon \mu_0(i) =_{\C P(X)} \mu_0(j)$. It 
suffices to show that $\lambda_0(i) =_{\C P(X)} \lambda_0(j)$.  
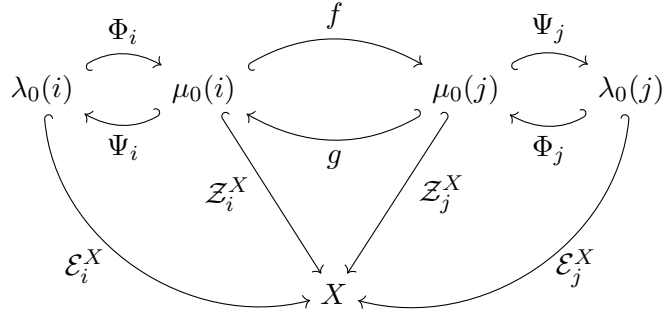
 \begin{center}
\resizebox{9cm}{!}{%
\begin{tikzpicture}

\node (E) at (0,0) {$\mu_0(i)$};
\node[right=of E] (B) {};
\node[right=of B] (F) {$\mu_0(j)$};
\node[below=of B] (C) {};
\node[below=of C] (A) {$X$};
\node[right=of F] (K) {$\lambda_0(j)$};
\node[left=of E] (L) {$\lambda_0(i)$};

\draw[left hook->,bend left] (E) to node [midway,above] {$f$} (F);
\draw[left hook->,bend left] (F) to node [midway,below] {$g$} (E);
\draw[right hook->] (E)--(A) node [midway,left] {$\C Z^X_i \ $};
\draw[left hook->] (F)--(A) node [midway,right] {$ \ \C Z^X_j$};

\draw[left hook->,bend left] (L) to node [midway,above] {$\Phi_i$} (E);
\draw[left hook->,bend left] (E) to node [midway,below] {$\Psi_i$} (L);
\draw[right hook->,bend right=50] (L) to node [midway,left] {$\C E^X_i \ \ $} (A);

\draw[left hook->,bend left] (F) to node [midway,above] {$\Psi_j$} (K);
\draw[left hook->,bend left] (K) to node [midway,below] {$\Phi_j$} (F);
\draw[left hook->,bend left=50] (K) to node [midway,right] {$\ \C E^X_j $} (A);

\end{tikzpicture}
}
\end{center}
If we define $f{'} := \Psi_j \circ f \circ \Phi_i$ and $g{'} := \Psi_i \circ g \circ \Phi_j$, it is straightforward to 
show that $(f{'}, g{'}) \colon \lambda_0(i) =_{\C P(X)} \lambda_0(j)$, hence $i =_I j$.
\end{proof}

By the previous remark $Q(\Lambda(X))$ is an extensional property on $\Fam(I,X)$. Since $\Set(I,X)$ is defined by 
separation on $\Fam(I,X)$, and since we see no objection to consider $\Fam(I,X)$ to be a set, we also see no objection  
to consider $\Set(I,X)$ to be a set.

\begin{definition}\label{def: sublambdaI}
Let $\Lambda(X) := (\lambda_0, \C E^X \lambda_1) \in \Fam(I,X)$. Let the equality 
$=_I^{\Lambda(X)}$\index{$=_I^{\Lambda(X)}$} on $I$
\index{equality on the index-set induced by a family of subsets}
given by $ i =_I^{\Lambda(X)} j : \TOT \lambda_0 (i) =_{\C P(X)} \lambda_0 (j)$, for every $i, j \in I$.
The \textit{set $\lambda_0 I(X)$ of subsets of $X$
generated by $\Lambda(X)$}\index{set of subsets generated by a family of subsets}\index{$\lambda_0 I(X)$} is 
the totality $I$ equipped with the equality $=_I^{\Lambda(X)}$. We write $\lambda_0(i) \in \lambda_0I(X)$,
instead of $i \in I$, when $I$ is equipped with the equality $=_I^{\Lambda(X)}$. 
The operation $\lambda_0^* : I \sto I$ from $(I, =_I)$ to $(I, =_I^{\Lambda(X)})$ is 
defined as in Definition~\ref{def: lambdaI}.
\end{definition}

Clearly, $\lambda_0^*$ is a function. 
All results in section~\ref{sec: setofsets} are shown similarly for sets of subsets, and for convenience 
we include them here without proof.

\begin{proposition}\label{prp: subFamtoset1}
Let $\Lambda(X) := (\lambda_0, \C E^X \lambda_1) \in \Set(I, X)$, and let $Y$ be a set.
If $f \colon I \to Y$, there is a unique function $\lambda_0 f \colon \lambda_0 I(X) \to Y$ such that the 
following diagram commutes
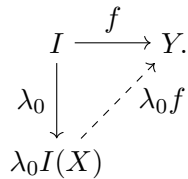
\begin{center}
\begin{tikzpicture}

\node (E) at (0,0) {$\lambda_0 I(X)$};
\node [above=of E] (D) {$I$};
\node[right=of D] (F) {$Y.$};

\draw[dashed, ->] (E)--(F) node [midway,right] {$\ \lambda_0 f$};
\draw [->]  (D)--(E) node [midway,left] {$\lambda_0$};
\draw[->] (D)--(F) node [midway,above] {$f$};

\end{tikzpicture}
\end{center}
Conversely, if $f \colon I \sto Y$ and $f^* \colon \lambda_0 I(X)  \to Y$ such that 
the corresponding diagram commutes, then $f$ is a function and $f^*$ is equal to the function from 
$\lambda_0 I(X)$ to $Y$ generated by $f$.
\end{proposition}

\begin{proposition}\label{prp: subFamtoset2}
Let $\Lambda(X) := (\lambda_0, \C E^X, \lambda_1) \in \Fam(I, X)$, and let $Y$ be a set.
If $f^* : \lambda_0 I(X) \to Y$, there is a unique function $f \colon I \to Y$ such that the following diagram commutes
\begin{center}
\begin{tikzpicture}

\node (E) at (0,0) {$\lambda_0 I(X)$};
\node [above=of E] (D) {$I$};
\node[right=of D] (F) {$Y.$};

\draw[->] (E)--(F) node [midway,right] {$\ f^*$};
\draw [->]  (D)--(E) node [midway,left] {$\lambda_0$};
\draw[dashed, ->] (D)--(F) node [midway,above] {$f$};

\end{tikzpicture}
\end{center}
If $\Lambda \in \Set(I,X)$, then $f^*$ is equal to the function from $\lambda_0 I(X)$ 
to $Y$ generated by $f$.
\end{proposition}

\begin{definition}\label{def: comp}
Let  $\Lambda(X) := (\lambda_0, \C E^X \lambda_1) \in \Set(I, X)$, and let $Y$ be a set.
If $f^* \colon \lambda_0 I(X) \to Y$, we denote the unique function $f \colon I \to Y$ generated by $f^*$ by 
$f^* \circ \lambda_0.$
\end{definition}

\begin{corollary}\label{cor: subfamtoset3}
Let  $\Lambda(X) := (\lambda_0, \C E^X \lambda_1) \in \Fam(I, X)$, and let $Y$ be a set.\\[1mm]
\normalfont (i) 
\itshape The operation $\Phi \colon \D F(\lambda_0 I, Y) \sto \D F(I, Y)$, defined by
$\Phi(f^*) := f^* \circ \lambda_0$, for every $f^* \in \D F(\lambda_0 I, Y)$, is an embedding.\\[1mm]
\normalfont (ii) 
\itshape If $\Lambda(X) \in \Set(I, X)$, then $\Phi$ is a surjection, the operation 
$\Theta \colon \D F(I, Y) \sto \D F(\lambda_0 I, Y)$, defined by
$\Theta(f) := \lambda_0 f$, for every $f \in \D F(I, Y)$, is an embedding, and
$(\Theta, \Phi) \colon (\D F(I, Y) =_{\D V_0} \D F(\lambda_0 I, Y)$.
\end{corollary}

\begin{proof}
(i) By definition of the corresponding equalities we have that
\begin{align*}
f^* =_{\D F(\lambda_0 I(X), Y)} g^* & \TOT \forall_{i \in I}\big(f^*(\lambda_0 (i)) =_Y 
g^*(\lambda_0 (i))\big)\\
& \TOT \forall_{i \in I}\big((f^* \circ \lambda_0)(i)) =_Y (g^* \circ \lambda_0)(i))\big)\\
& \TOT f^* \circ \lambda_0 =_{\D F(I, Y)} g^* \circ \lambda_0.
\end{align*}
(ii) If $f \in \D F(I, Y)$, then by Proposition~\ref{prp: subFamtoset1} there is unique $\lambda_0 f
\in \D F(\lambda_0 I, Y)$ such that $\Phi(\lambda_0 f) := \lambda_0 f \circ \lambda_0 =_{\D F(I, Y)} f$.
By definition of the corresponding equalities we have that
\begin{align*}
f =_{\D F(I, Y)} g & \TOT \forall_{i \in I}\big(f(i)) =_Y g(i)\big)\\
& \TOT \forall_{i \in I}\big(\lambda_0 f (\lambda_0(i)) =_Y \lambda_0 g(\lambda_0(i)\big)\\
& \TOT \lambda_0 f =_{\D F(\lambda_0 I, Y)} \lambda_0 g.
\end{align*}
Moreover, we have that 
$(\Theta \circ \Phi)(f^*) := \Theta(f^* \circ \lambda_0) := \lambda_0 (f^* \circ \lambda_0) 
=_{\D F(\lambda_0 I(X), Y)} f^*$, and 
$(\Phi \circ \Theta)(f) := \Phi(\lambda_0 f) := (\lambda_0 f) \circ  \lambda_0 =_{\D F(I, Y)} f$.
\end{proof}

\begin{proposition}\label{prp: subFamtoFam1}
Let $\Lambda(X) := (\lambda_0, \C E^X, \lambda_1) \in \Set(I,X)$ and 
$M(X) := (\mu_0, \C Z^X \mu_1) \in \Set(J,Y)$. If $f \colon I \to J$, there is a unique 
function $f^* \colon \lambda_0 I(X) \to \mu_0 J(Y)$ such that the following diagram commutes
\begin{center}
\begin{tikzpicture}

\node (E) at (0,0) {$\lambda_0 I(X)$};
\node[right=of E] (F) {$\mu_0 J(Y).$};
\node[above=of F] (A) {$J$};
\node [above=of E] (D) {$I$};

\draw[dashed, ->] (E)--(F) node [midway,below] {$f^*$};
\draw[->] (D)--(A) node [midway,above] {$f $};
\draw[->] (D)--(E) node [midway,left] {$\lambda_0$};
\draw[->] (A)--(F) node [midway,right] {$\mu_0$};

\end{tikzpicture}
\end{center}
If $f \colon I \sto J$, and 
$f^* : \lambda_0 I(X) \to \mu_0 J(Y)$ such that the corresponding to the above diagram commutes, then $f \in \D F(I, J)$
and $f^*$ is equal to the map in $\D F\big(\lambda_0 I(X), \mu_0 J(Y)\big)$ generated by $f$.
\end{proposition}

\begin{remark}\label{rem: detachablefam}
Let the set $\big(X, =_X, \neq_X^{\mathsmaller{\D F(X, \D 2)}}\big)$, and\index{$\Delta^1(X)$}
$\Delta^1(X) := \big(\delta_0^1, \C E^{1,X}, \delta_1^1\big)$, where the non-dependent assignment routine
$\delta_0^1 \colon \D F(X, \D 2) \sto \D V_0$ is defined by the rule $f \mapsto \delta_0^1(f)$, for every 
$f \in \D F(X, \D 2)$ $($see Definition~\ref{def: detachable}$)$, and the dependent operations 
$\C E^{1,X} \colon \bigcurlywedge_{f \in D F(X, \D 2)}\D F(\delta_0^1(f), X)$ and 
$\delta_1^1 \colon \bigcurlywedge_{(f,g) \in D(\D F(X, \D 2))}\D F(\delta_0^1(f), \delta_0^1(g))$ are defined, respectively, by
\[ \C E^{1,X}_f \colon \delta_0^1(f) \eto X \ \ \ \ x \mapsto x; \ \ \ \ x \in \delta_0^1(f), \]
\[ \delta_1^1(f,g) := \delta_{fg}^1 \colon \delta_0^1(f) \to \delta_0^1(g) \ \ \ \ x \mapsto x; \ \ \ \ 
x \in \delta_0^1(f). \]
If $\Delta^0(X) := \big(\delta_0^0, \C E^{0,X}, \delta_1^0\big)$, where\index{$\Delta^0(X)$} 
$\delta_0^0 \colon \D F(X, \D 2) \sto \D V_0$ is defined by the rule $f \mapsto \delta_0^0(f)$, 
for every $f \in \D F(X, \D 2)$, and the dependent operations $\C E^{0,X}, \delta_1^0$ are defined similarly,
then $\Delta^1(X), \Delta^0(X) \in \Set(\D F(X, \D 2), X)$, and they are called
the $\D F(X, \D 2)$-sets of detachable subsets of $X$\index{sets of detachable subsets}.
\end{remark}

\begin{proof}
We give the proof only for $\Delta^1(X)$. It is easy to show that $\Delta^1(X) \in \Fam(\D F(X, \D 2), X)$.
Let $f, g \in \D F(X, \D 2)$ such that $\delta_0^1(f) =_{\C P(X)} \delta_0^1(g)$ i.e., there are 
$e \in \D F\big(\delta_0^1(f), \delta_0^1(g)\big)$ and $k \in \D F\big(\delta_0^1(g), \delta_0^1(f)\big)$ such that
$(e, k) \colon \delta_0^1(f) =_{\C P(X)} \delta_0^1(g)$
 \begin{center}
 \resizebox{4cm}{!}{%
\begin{tikzpicture}

\node (E) at (0,0) {$\delta_0^1(f)$};
\node[right=of E] (B) {};
\node[right=of B] (F) {$\delta_0^1(g)$};
\node[below=of B] (C) {};
\node[below=of C] (A) {$X$.};

\draw[left hook->,bend left] (E) to node [midway,above] {$e$} (F);
\draw[left hook->,bend left] (F) to node [midway,below] {$k$} (E);
\draw[right hook->] (E)--(A) node [midway,left] {$\C E^{1,X}_f \ $};
\draw[left hook->] (F)--(A) node [midway,right] {$ \ \C E^{1,X}_g$};

\end{tikzpicture}
}
\end{center} 
Let $x \in X$. By the commutativity of the above diagram $x := \C E^{1,X}_g(x) =_X \C E^{1,X}_f\big(k(x)\big) := k(x)$.
Hence, if $f(x) =_{\D 2} 1$, then $f(k(x)) =_{\D 2} 1$. Since $k(x) \in \delta_0^1(g)$, we get
$g(k(x)) =_{\D 2} 1$, and since $x =_X k(x)$, we get $g(x) =_{\D 2} 1$. If $f(x) =_{\D 2} 0$, we use proceed similarly.
\end{proof}

Clearly, $\delta_0^1(\overline{1}) = X$, 
$\delta_0^1(f) \cap \delta_0^1(g) = \delta_0^1(f \cdot g)$, and $\delta_0^1(f) \cup \delta_0^1(g) = 
\delta_0^1(f + g - f \cdot g)$. 

\begin{proposition}\label{prp: famdetachable1}
Let the family $\Delta^1(X) := \big(\delta_0^1, \C E^{1,X}, \delta_1^1\big)$ of detachable subsets of $X$.\\[1mm]
If $\compl \colon \D F(X, \D 2) \to \D F(X, \D 2)$ is defined by $f \mapsto 1-f,$
for every $f \in \D F(X, \D 2)$, then the operation $\Compl \colon [\delta \D F(X, \D 2)](X) \sto
[\delta \D F(X, \D 2)](X)$, defined by
\[ \Compl(\delta_0^1(f)) := \delta(\compl(f)) =: \delta_0^1(1 - f) = \delta_0^0(f); \ \ \ \ \delta_0^1(f) \in 
[\delta \D F(X, \D 2)](X), \]
is a function such that the following conditions hold:\\[1mm]
\normalfont (a) 
\itshape  $\Compl(\Compl(\delta_0^1(f)) = \delta_0^1(f)$.\\[1mm]
\normalfont (b) 
\itshape  $\Compl(\delta_0^1(f) \cap \delta_0^1(g)) = \Compl(\delta_0^1(f)) \cup \Compl_0^1(\delta(g))$.\\[1mm]
\normalfont (c) 
\itshape  $\Compl(\delta_0^1(f) \cup \delta_0^1(g)) = \Compl(\delta_0^1(f)) \cap \Compl(\delta_0^1(g))$.
\end{proposition}

\begin{proof}
(i) By Proposition~\ref{prp: subFamtoset1} the operation $\Compl$ is the unique function 
from $[\delta_0^1 \D F(X, \D 2)](X)$
to $[\delta_0^1 \D F(X, \D 2)](X)$ that makes the following diagram commutative
\begin{center}
\begin{tikzpicture}

\node (E) at (0,0) {$[\delta_0^1 \D F(X, \D 2)](X)$};
\node[right=of E] (F) {$[\delta_0^1 \D F(X, \D 2)](X).$};
\node[above=of F] (A) {$\D F(X, \D 2)$};
\node [above=of E] (D) {$\D F(X, \D 2)$};

\draw[dashed, ->] (E)--(F) node [midway,below] {$\Compl$};
\draw[->] (D)--(A) node [midway,above] {$\compl $};
\draw[->] (D)--(E) node [midway,left] {$\delta_0^1$};
\draw[->] (A)--(F) node [midway,right] {$\delta_0^1$};

\end{tikzpicture}
\end{center}
The proofs of conditions (a)-(c) are easy to show.
\end{proof}

\begin{proposition}\label{prp: famdetachable2}
Let $X, Y$ be sets, and let the sets of detachable subsets 
$\Delta^1(X) := \big(\delta_0^{1,X}, \C E^{1,X}, \delta_1^{X}\big)$, 
$\Delta^1(Y) := \big(\delta_0^{1,Y}, \C E^{1,Y}, \delta_1^{1,Y}\big)$  
of $X$ and $Y$, respectively. 
If $h : Y \to X$, then the operation $\tilde{h} : \D F(X, \D 2) \sto \D F(Y, \D 2)$, defined by $f \mapsto f \circ h$,
for every $f \in \D F(X, \D 2)$, is
a function,
and there is a unique function 
$\delta_0^1 \tilde{h} : [\delta_0^1 \D F(X, \D 2)](X) \to [\delta_0^1 D F(Y, \D 2)](Y)$ 
such that the following diagram commutes
\begin{center}
\begin{tikzpicture}

\node (E) at (0,0) {$[\delta_0^1 \D F(X, \D 2)](X)$};
\node[right=of E] (F) {$[\delta_0^1 \D F(Y, \D 2)](Y).$};
\node[above=of F] (A) {$\D F(X, \D 2)$};
\node [above=of E] (D) {$\D F(Y, \D 2)$};

\draw[dashed, ->] (E)--(F) node [midway,below] {$\tilde{h}$};
\draw[->] (D)--(A) node [midway,above] {$\delta_0^1 \tilde{h} $};
\draw[->] (D)--(E) node [midway,left] {$\delta_{0,Y}^1$};
\draw[->] (A)--(F) node [midway,right] {$\delta_{0,X}^1$};

\end{tikzpicture}
\end{center}
\end{proposition}

\begin{proof}
It follows immediately from Proposition~\ref{prp: subFamtoFam1}.
\end{proof}

\begin{proposition}\label{prp: subfamtofam2}
Let $\Lambda(X) := (\lambda_0, \C E^X, \lambda_1) \in \Fam(I,X)$ and 
$M(Y) := (\mu_0, \C Z^Y \mu_1) \in \Set(J,Y)$.
If $f^* \colon \lambda_0 I(X) \to \mu_0 J(Y)$, there is a unique $f \colon I \to J$, such that
the following diagram commutes
\begin{center}
\begin{tikzpicture}

\node (E) at (0,0) {$\lambda_0 I$};
\node[right=of E] (F) {$\mu_0 J.$};
\node[above=of F] (A) {$J$};
\node [above=of E] (D) {$I$};

\draw[->] (E)--(F) node [midway,below] {$f^*$};
\draw[dashed, ->] (D)--(A) node [midway,above] {$f $};
\draw[->]  (D)--(E) node [midway,left] {$\lambda_0$};
\draw[->]  (A)--(F) node [midway,right] {$\mu_0$};

\end{tikzpicture}
\end{center}
Moreover, $f^*$ is equal to the function from $\lambda_0 I(X)$ to $\mu_0 J(Y)$ generated by $f$. 
\end{proposition}

\begin{corollary}\label{cor: subfamtofam3}
Let $\Lambda(X) := (\lambda_0, \C E^X, \lambda_1) \in \Set(I,X)$ and 
$M(Y) := (\mu_0, \C Z^Y \mu_1) \in \Fam(J,Y)$. The operation $\Theta \colon \D F(I, J)
\sto \D F(\lambda_0 I, \mu_0 J)$, defined by $f \mapsto f^*$, for every $f \in \D F(I, J)$,
is a function. If $M(Y) \in \Set(J,Y)$, then $\Theta$ is an embedding, and a surjection.
\end{corollary}

\begin{proof}
By definition of the corresponding equalities we have that
\begin{align*}
f =_{\D F(I, J)} g & \TOT \forall_{i \in I}\big(f(i) =_J g(i)\big)\\
& \To \forall_{i \in I}\big(\mu_0 (f(i)) =_{\C P(Y)} \mu_0 (g(i))\big)\\
& \TOT \forall_{i \in I}\big(f^*(\lambda_0(i)) =_{\C P(Y)} g^*(\lambda_0(i))\big)\\
& \TOT f^* =_{\D F(\lambda_0 I(X), \mu_0 J(Y))} g^*.
\end{align*}
If $M(Y) \in \Set(J,Y)$, the above implication is also an equivalence, 
hence $\Theta$ is an embedding. By Proposition~\ref{prp: subFamtoset2} we have that $\Theta$ is a surjection. 
\end{proof}

The notions of fiber and cofiber of a function were introduced in Definition~\ref{def: surjective}.

\begin{proposition}\label{prp: fibcofib}
Let the sets $(X, =_X, \neq_X)$ and $(Y, =_Y, \neq_Y)$, and let $\fXY$.\\[1mm]
\normalfont (i) 
\itshape Let $\fib^f(X) := \big(\fib^f_0, \C E^{\fib,X}, \fib_1^f\big)$, where\index{$\fib^f(X)$}
$\fib^f_0 \colon Y \sto \D V_0$ is
defined by the rule $\fib^f_0(y) := \fib^f(y)$, for every $y \in Y$, and the dependent operations 
$\C E^{\fib,X} \colon \bigcurlywedge_{y \in Y}\D F(\fib^f_0(y), X)$ and 
$\fib_1^f \colon \bigcurlywedge_{(y,y{'}) \in D(Y)}\D F(\fib_0^f(y), \fib_0^f(y{'}))$ are defined, respectively, by 
\[ \C E^{\fib,X} \colon \fib^f_0(y) \eto X \ \ \ \ x \mapsto x; \ \ \ \ x \in \fib^f_0(y), \]
\[ \fib_1^1(y,y{'}) := \fib_{yy{'}}^1 \colon \fib^f(y) \to \fib^f(y{'}) \ \ \ \ x \mapsto x; \ \ \ \ x \in 
\fib^f_0(y). \] 
Then $\fib^f(X) \in \Fam(Y, X)$ and if $f$ is a surjection, then $\fib^f(X) \in \Set(Y, X)$.\\[1mm]
\normalfont (ii) 
\itshape $f$ is strongly extensional if and only if $\cofib^f_0(y) \Disj \fib^f_0(y)$, for every $y \in Y$.\\[1mm]
\normalfont (iii) 
\itshape Let $\cofib^f(X) := \big(\cofib^f_0, \C E^{\cofib,X}, \cofib_1^f\big)$, where\index{$\cofib^f(X)$}
$\cofib^f_0 \colon Y \sto \D V_0$ is
defined by the rule $\cofib^f_0(y) := \cofib^f(y)$, for every $y \in Y$, and 
$\C E^{\cofib,X} \colon \bigcurlywedge_{y \in Y}\D F(\cofib^f_0(y), X)$,  
$\cofib_1^f \colon \bigcurlywedge_{(y,y{'}) \in D(Y)}\D F(\cofib_0^f(y), \cofib_0^f(y{'}))$ are defined, respectively, by 
\[ \C E^{\cofib,X}_y \colon \cofib^f_0(y) \eto X \ \ \ \ x \mapsto x; \ \ \ \ x \in \cofib^f_0(y), \]
\[ \cofib_1^1(y,y{'}) := \cofib_{yy{'}}^1 \colon \cofib^f(y) \to \cofib^f(y{'}) \ \ \ \ x \mapsto x; \ \ \ \ x \in 
\cofib^f_0(y). \] 
Then $\cofib^f(X) \in \Fam(Y, X)$, and if $f$ is a surjection, then $\cofib^f(X) \in \Set(Y, X)$ if and only if
the inequality $\neq_Y$ is tight.
\end{proposition}

\begin{proof}
(i)  If $y =_Y y{'}$ and $x \in \fib^f(y)$, then $x \in \fib^f(y{'})$. Since the functions $\C E^{\fib,X}_y, 
\C E^{\fib,X}_{y{'}}$, and $\fib_{yy{'}}^1$ are
defined through the identity map-rule, we get $\fib^f(X) \in \Fam(Y, X)$. Let $y, y{'} \in Y$ and 
functions $g \in \D F\big(\fib^f(y), \fib^f(y{'})\big)$
and $h \in \D F\big(\fib^f(y{'}), \fib^f(y)\big)$, such that $(g, h) \colon \fib^f(y) =_{\C P(X)} \fib^f(y{'})$.
Let $x \in X$ such that $f(x) =_Y y$ i.e., $x \in \fib^f(y)$. By the commutativity of one of the following left inner
diagrams we have that $g(x) =_X x$, and, of course, $g(x) \in \fib^f(y{'})$ i.e., $f(g(x)) =_Y y{'}$.
Hence, $y{'} =_Y f(g(x)) =_Y f(x) =_Y y$.
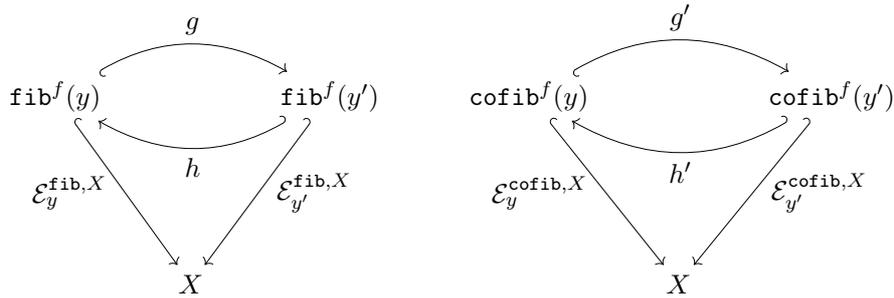
\begin{center}
\resizebox{12cm}{!}{%
\begin{tikzpicture}

\node (E) at (0,0) {$\fib^f(y)$};
\node[right=of E] (B) {};
\node[right=of B] (F) {$\fib^f(y{'})$};
\node[below=of B] (C) {};
\node[below=of C] (A) {$X$};

\node[right=of F] (G) {$\cofib^f(y)$};
\node[right=of G] (H) {};
\node[right=of H] (J) {$\cofib^f(y{'})$};
\node[below=of H] (K) {};
\node[below=of K] (L) {$X$};

\draw[left hook->,bend left] (E) to node [midway,above] {$g$} (F);
\draw[left hook->,bend left] (F) to node [midway,below] {$h$} (E);
\draw[right hook->] (E)--(A) node [midway,left] {$\C E^{\fib,X}_y \ $};
\draw[left hook->] (F)--(A) node [midway,right] {$ \ \C E^{\fib,X}_{y{'}}$};
\draw[left hook->,bend left] (G) to node [midway,above] {$g{'}$} (J);
\draw[left hook->,bend left] (J) to node [midway,below] {$h{'}$} (G);
\draw[right hook->] (G)--(L) node [midway,left] {$\C E^{\cofib,X}_y \ $};
\draw[left hook->] (J)--(L) node [midway,right] {$ \ \C E^{\cofib,X}_{y{'}}$};

\end{tikzpicture}
}
\end{center} 
(ii) Suppose that $f$ is strongly extensional and let $x \in \cofib^f(y)$ and $z \in \fib^f(y)$ i.e., $f(x) \neq_Y y$
and $f(x) =_Y y$. By the extensionality of $\neq_Y$ (Remark~\ref{rem: apartness1}) we get $f(x) \neq_Y f(z)$, and
as $f$ is strongly extensional, we conclude that $x \neq_X z$. Suppose next that 
$\cofib^f_0(y) \Disj \fib^f_0(y)$, for every $y \in Y$, and let $x,z \in X$ with $f(x) \neq_Y f(z)$. In this case,
we get $x \in \fib^f(f(x))$ and $z \in \cofib^f(f(x))$. Since $\cofib^f_0(f(x)) \Disj \fib^f_0(f(x))$ and the corresponding
embeddings into $X$ are given by the identity map-rule, we get $x \neq_X z$.\\
(iii)  If $y =_Y y{'}$ and $x \in \cofib^f(y)$, then $f(x) \neq_Y y$, and by the extensionality of $\neq_Y$, we get
$f(x) \neq_Y y{'}$ i.e., $x \in \fib^f(y{'})$. Since the functions $\C E^{\cofib,X}_y, 
\C E^{\cofib,X}_{y{'}}$, and $\cofib_{yy{'}}^1$ are
defined through the identity map-rule, we get $\cofib^f(X) \in \Fam(Y, X)$. 
Let $f$ be a surjection. We suppose first that $\cofib^f(X) \in \Set(Y, X)$. If $\neg(y \neq_Y y{'})$, we show that
$y =_Y y{'}$, by showing that $\cofib^f(y) =_{\C P(X)} \cofib^f(y{'})$. If $x \in \cofib^f(y)$, then $f(x) \neq_Y y$.
By condition $(\Ap_3)$ either $y{'} \neq_Y f(x)$ or $y{'} \neq_Y y$. Since the latter contradicts our hypothesis 
$\neg(y \neq_Y y{'})$, we conclude that $y{'} \neq_Y f(x)$ i.e., $x \in \cofib^f(y{'})$. Similarly we show that
if $x \in \cofib^f(y{'})$, then $x \in \cofib^f(y)$. Hence, the functions between $\cofib^f(y)$ and $\cofib^f(y{'})$
that are given by the identity map-rule witness the equality $\cofib^f(y) =_{\C P(X)} \cofib^f(y{'})$.
Suppose next that the inequality $\neq_Y$ is tight. Let $y, y{'} \in Y$ and let functions 
$g{'} \in \D F\big(\cofib^f(y), \cofib^f(y{'})\big)$ and $h{'} \in \D F\big(\cofib^f(y{'}), \cofib^f(y)\big)$, 
such that $(g{'}, h{'}) \colon \cofib^f(y) =_{\C P(X)} \cofib^f(y{'})$. We show that $y =_Y y{'}$ by showing 
$\neg(y \neq_Y y{'})$. For that suppose $y \neq_Y y{'}$, and let $x, x{'} \in X$ such that $f(x) =_Y y$ and
$f(x{'}) =_Y y{'}$. By the extensionality of $\neq_Y$ we get $f(x) \neq_Y y{'}$ i.e., $x \in \cofib^f(y{'})$.
Since $h{'}(x) \in \cofib^f(y)$, and since by the commutativity of one of the above right inner diagrams $h{'}(x) =_X x$,
we get $x \in \cofib^f(y)$. Since $f(x) \neq_Y y$ and $y =_Y f(x)$, by the extensionality of $\neq_Y$ we get
$f(x) \neq_Y f(x)$, which leads to the required contradiction.
\end{proof}

If $f$ is not a surjection, it is possible that $\fib^f(y), \fib^f(y{'})$ are not inhabited, and $y \neq_Y y{'}$. 
If $f$ is not a surjection, like the function $f \colon X \to \{0, 1, 2\}$, defined by $f(x) := 0$, for every $x \in X$,
then $\cofib^f(1) = X = \cofib^f(2)$ and $1 \neq 2$. Notice that it is not necessary that a family of subsets 
is a family of fibers or a family of cofibers, as the moduli of embeddings of the latter
are given through the identity map-rule. 

\begin{definition}\label{def: famofcontrsets}
An $I$-family of sets $\Lambda := (\lambda_0, \lambda_1)$ is a family of 
contractible sets\index{family of contractible sets}, if $\lambda_0(i)$ is contractible, for every $i \in I$.
A modulus of centres of contraction for $\Lambda$\index{modulus of centres of contraction} is a 
dependent operation $\centr^{\Lambda} \colon \bigcurlywedge_{i \in I}\lambda_0(i)$, with 
$\centr^{\Lambda}_i$ a centre of contraction for $\lambda_0(i)$, for every $i \in I$.
\end{definition}

In Proposition~\ref{prp: fiber1} we saw that if $(f, g) \colon X =_{\D V_0} Y$,
the set $\fib^f(y)$ is contractible with $\centr^f_y := g(y)$, for every $y \in Y$ i.e., the dependent operation
$\centr^f$ is a modulus of centres of contractions for the family $\fib^f(X)$. Next follows a kind of inverse to
Proposition~\ref{prp: fiber1}. 

\begin{proposition}\label{prp: centerfiber}
Let the sets $(X, =_X, \neq_X)$, $(Y, =_Y, \neq_Y)$, and $\fXY$. If  
$\fib^f(X) := \big(\fib^f_0, \C E^{\fib,X}, \fib_1^f\big)$ is a family of contractible subsets of $X$ with 
$\centr^{\fib^f(X)} \colon \bigcurlywedge_{y \in Y}\fib^f(y)$ a modulus of centres of contraction for $\fib^f(X)$,
there is $g \in \D F(Y,X)$ with $(f, g) \colon X =_{\D V_0} Y$.
\end{proposition}

\begin{proof}
Let the operation $g \colon Y \sto X$, defined by $g(y) := \centr^{\fib^f(X)}(y)$, for every $y \in Y$. 
Since $g(y) \in \fib^f(y)$, we have that $f(g(y)) =_Y y$. Since $g(y)$ is a centre of contraction for $\fib^f(y)$, we
have that $\forall_{x \in X}\big(f(x) =_Y y \To x =_X g(y)\big)$. First we show that the operation $g$ is a function.
For that, let $y =_Y y{'}$, and we show that $g(y) =_X g(y{'})$. Since the map $\fib_{yy{'}}^1$ in 
Proposition~\ref{prp: fibcofib} is given by the identity map-rule, and since $g(y{'}) \in \fib^f(y{'})$, we get
$g(y{'}) \in \fib^f(y)$. Since $g(y)$ is a centre of contraction for $\fib^f(y)$, we get $g(y{'}) =_X g(y)$.
It remains to show that if $x \in X$, then $g(f(x)) =_X x$. By the definition of $g$ we have that $g(f(x)) :=
\centr^{\fib^f(X)}(f(x))$. As $x \in \fib^f(f(x))$, we get $x =_X g(f(x))$.
\end{proof}

\section{Families of equivalence classes}
\label{sec: fameqclass}

In this section we extend results on sets of subsets to families of equivalence classes.
Although a family of equivalence classes is not, in general, a set of subsets, 
we can define functions on them, if we use appropriate functions on their index-set.

\begin{definition}\label{def: equivstructure}
If $X$ is a set and $R_X(x, x{'})$ is an extensional property on $X \times X$ that 
satisfies the conditions of an equivalence relation, we call the pair $(X,R_X)$  
an \textit{equivalence structure}\index{equivalence structure}. If $(Y, S_Y)$ is an equivalence structure,
a function $f \colon X \to Y$ is an \textit{equivalence preserving function}\index{equivalence preserving
function}, or an $(R_X, S_Y)$-function\index{$(R_X, S_Y)$-function}, if
\[ \forall_{x, x{'} \in X}\big(R(x, x{'}) \To S(f(x), f(x{'}))\big). \]
If, for every $x,x{'} \in X$, the converse implication holds, we say that  
$f$ is an $(R_X, S_Y)$-embedding\index{$(R_X, S_Y)$-embedding}.
Let $\D F(R_X, S_Y)$\index{$\D F(R_X, S_Y)$} be the set of $(R_X, S_Y)$-functions\footnote{By the extensionality of $S_Y$
the property of being an $(R_X, S_Y)$-function is extensional on $\D F(X, Y)$.}. 
\end{definition}

\begin{proposition}\label{prp: equivstr1}
If $(X, R_X)$ is an  equivalence structure, let $R(X) := \big(\rho_0, \C R^X, \rho_1\big)$, where\index{$R(X)$}
$\rho_0 \colon X \sto \D V_0$ is defined by 
$\rho_0(x) := \{y \in X \mid R_X(y, x)\}$, for every $x \in X$,
and the dependent operations 
$\C R^X \colon \bigcurlywedge_{x \in X}\D F\big(\rho_0(x), X\big)$, 
$\rho_1 \colon \bigcurlywedge_{(x,x{'}) \in D(X)}\D F\big(\rho_0(x), \rho_0(x{'})\big)$ are defined by 
\[ \C R^X_x \colon \rho_0(x) \eto X \ \ \ \ y \mapsto y; \ \ \ \ y \in \rho_0(x), \]
\[ \rho_1(x,x{'}) := \rho_{xx{'}} \colon \rho_0(x) \to \rho_0(x{'}) \ \ \ \ y \mapsto y; \ \ \ \ y \in 
\rho_0(x). \] 
Then $R(X) \in \Fam(X, X)$, such that $\forall_{x x{'} \in X}\big(\rho_0(x) =_{\C P(X)} \rho_0(x{'}) \To 
R(x, x{'})\big)$.

\end{proposition}

\begin{proof}
By the extensionality of $R_X$ the set $\rho_0(x)$ is a well-defined extensional subset of $X$. If $x =_X x{'}$ and 
$R_X(y, x)$, then by the extensionality of $R_X$ we get $R_X(y, x{'})$, hence $\rho_{xx{'}}$ is well-defined. 
Let $(f, g) \colon \rho_0(x) =_{\C P(X)} \rho_0(x{'})$
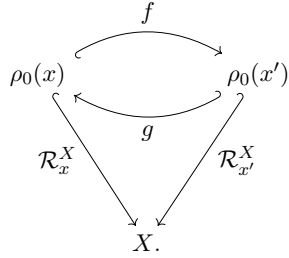
\begin{center}
\resizebox{4cm}{!}{%
\begin{tikzpicture}

\node (E) at (0,0) {$\rho_0(x)$};
\node[right=of E] (B) {};
\node[right=of B] (F) {$\rho_0(x{'})$};
\node[below=of B] (C) {};
\node[below=of C] (A) {$X$.};

\draw[left hook->,bend left] (E) to node [midway,above] {$f$} (F);
\draw[left hook->,bend left] (F) to node [midway,below] {$g$} (E);
\draw[right hook->] (E)--(A) node [midway,left] {$\C R^X_x \ $};
\draw[left hook->] (F)--(A) node [midway,right] {$ \ \C R^X_{x{'}}$};

\end{tikzpicture}
}
\end{center} 
If $y \in \rho_0(x) :\TOT R_X(y, x)$, then $f(y) \in \rho_0(x{'}) :\TOT R_X(f(y), x{'})$, and by the commutativity of 
the corresponding above diagram we get $f(y) =_X y$. Hence by the extensionality of $R_X$ we get $R_X(y, x{'})$.
Since $R_X(y, x)$ implies $R_X(x, y)$, by transitivity we get $R_X(x, x{'})$.
\end{proof}

\begin{corollary}\label{cor: corequivstr1}
 Let $\Eql(X) := \big(\eql_0^X, \C E^X, \eql_1^X\big)$\index{$\Eql(X)$} be the $X$-family of subsets of $X$ 
 induced by the equivalence relation $=_X$\index{$\eql_0^X(x)$} i.e., $\eql_0^X(x) := \{y \in X \mid y =_X x\}$. 
 Then $\Eql(X) \in \Set(X,X)$.
\end{corollary}

\begin{proof}
It follows immediately from Proposition~\ref{prp: equivstr1}. 
\end{proof}


\begin{proposition}\label{prp: quotient1}
If $(X, R_X)$ is an equivalence structure, and $f \colon X \to Y$ is an $(R_X, =_Y)$-function
there is a unique $\rho_0 f : \rho_0 X(X) \to Y$ such that the following diagram commutes
\begin{center}
\begin{tikzpicture}

\node (E) at (0,0) {$\rho_0 X(X)$};
\node [above=of E] (D) {$X$};
\node[right=of D] (F) {$Y.$};

\draw[dashed, ->] (E)--(F) node [midway,right] {$\rho_0 f$};
\draw [->] (D)--(E) node [midway,left] {$\rho^*_0$};
\draw[->] (D)--(F) node [midway,above] {$f$};

\end{tikzpicture}
\end{center}
Conversely, if $f \colon X \to Y$ and $f^* : \rho_0 X(X)  \to Y$ such that the above diagram commutes, 
then $f$ is an $(R_X, =_Y)$-function and $f^*$ is equal to the function from $\rho_0 X(X)$ to $Y$ 
generated by $f$.
\end{proposition}

\begin{proposition}\label{prp: quotient2}
If $(X, R_X)$ is an equivalence structure, and $f^* : \rho_0 X(X) \to Y$, there is a unique 
$f \colon X \to Y$, which is an $(R_X, =_Y)$-function, such that the following diagram commutes
\begin{center}
\begin{tikzpicture}

\node (E) at (0,0) {$\rho X(X)$};
\node [above=of E] (D) {$X$};
\node[right=of D] (F) {$Y.$};

\draw[->] (E)--(F) node [midway,right] {$f^*$};
\draw [->] (D)--(E) node [midway,left] {$\rho^*_0$};
\draw[dashed, ->] (D)--(F) node [midway,above] {$f$};

\end{tikzpicture}
\end{center}
Moreover, $f^*$ is equal to the function from $\rho_0 X(X)$ to $Y$ generated by $f$.
\end{proposition}


\begin{proposition}\label{prp: quotient4}
Let $(X, R_X)$ and $(Y, S_Y)$ be equivalence structures and $f \colon X \to Y$ an $(R_X, S_Y)$-function. 
If $R(X)$ and $S(Y)$ are the corresponding families of equivalence classes, there is a unique
function $f^* \colon \rho_0 X (X) \to \sigma_0 Y(Y)$ such that the following diagram commutes
\begin{center}
\begin{tikzpicture}

\node (E) at (0,0) {$\rho_0 X(X)$};
\node[right=of E] (F) {$\sigma_0 Y(Y).$};
\node[above=of F] (A) {$Y$};
\node [above=of E] (D) {$X$};

\draw[dashed, ->] (E)--(F) node [midway,below] {$f^*$};
\draw[->] (D)--(A) node [midway,above] {$f $};
\draw[->] (D)--(E) node [midway,left] {$\rho_0^*$};
\draw[->] (A)--(F) node [midway,right] {$\sigma_0^*$};

\end{tikzpicture}
\end{center}
If $f \colon X \to Y$ and $f^* : \rho_0 X(X) \to \sigma_0 Y(Y)$ such that the above diagram commutes,
then $f$ is an $(R_X, S_Y)$-function and $f^*$ is equal to the function from $\rho_0 X(X)$ to $\sigma_0 Y(Y)$
generated by $f$.
\end{proposition}

\begin{proof}
The assignment routine$f^*$ from $\rho_0 X(X)$ to $\sigma_0 Y(Y)$ defined by 
$f^* (\rho_0(x)) := \sigma_0 (f(x))$, for every $\rho_0(x) \in \rho_0 X(X)$
is extensional, since for every $x, x{'} \in X$ we have that
$\rho_0(x) =_{\C P(X)} \rho_0(x{'}) \To R_X(x, x{'}) \To S_Y(f(x), f(x{'}))$, hence 
$\sigma_0(f(x)) =_{\C P(Y)} \sigma_0(f(x{'})) : \TOT f^* (\rho_0(x)) =_{\C P(Y)}  f^* (\rho_0(x{'}))$.
The uniqueness of $f^*$ is immediate. For the converse, if $x, x{'} \in X$, then by the transitivity 
of $=_{\C P(Y)}$ we have that
$ R_X(x, x{'}) \To \rho_0 (x) =_{\C P(X)} \rho_0 (x{'}) \To f^* (\rho_0(x)) =_{\C P(Y)}  f^* (\rho_0(x{'}))
\To \sigma_0(f(x)) =_{\C P(Y)} \sigma_0(f(x{'}))$, hence $S_Y(f(x), f(x{'}))$.
The proof that $f^*$ is equal to the function from $\rho_0 X(X)$ to $\sigma_0 Y(Y)$ generated by $f$ is immediate.
\end{proof}

The previous is the constructive analogue to a standard classical fact (see~\cite{Du66}, p.~17).
A function $f^* \colon \rho_0 X(X) \to \sigma_0 Y(Y)$ 
does not generate a function from $X$ to $Y$.

\begin{proposition}\label{prp: quotient5}
Let $(X, R_X)$ and $(Y, S_Y)$ be equivalence structures and $R(X), S(Y)$ the families of their
equivalence classes. If $f^* : \rho_0 X(X) \to \sigma_0 Y(Y)$, there is 
$f \colon X \sto Y$, 
which is $(R_X, S_Y)$-preserving and $(=_X, S_Y)$-preserving, such that the following diagram commutes
\begin{center}
\begin{tikzpicture}

\node (E) at (0,0) {$\rho_0 X(X)$};
\node[right=of E] (F) {$\sigma_0 Y(Y).$};
\node[above=of F] (A) {$Y$};
\node [above=of E] (D) {$X$};

\draw[->] (E)--(F) node [midway,below] {$f^*$};
\draw[->] (D)--(A) node [midway,above] {$f $};
\draw[->] (D)--(E) node [midway,left] {$\rho^*$};
\draw[->]  (A)--(F) node [midway,right] {$\sigma^*$};

\end{tikzpicture}
\end{center}
\end{proposition}

\begin{proof}
If $x \in X$, then $f^*(\rho_0(x)) := \sigma_0(y)$, for some $y \in Y$. We define the routine $f(x) := y$ i.e., 
the output of $f^*$ determines the output of $f$. Since
$R(x, x{'})  \To \rho_0 (x) =_{\C P(X)} \rho_0 (x{'}) \To f^* (\rho_0(x)) =_{\C P(Y)}  f^* (\rho_0(x{'}))$,
hence $\sigma_0(y) =_{\C P(Y)} \sigma_0(y{'})$ and $S_Y(y, y{'})$, we get $S_Y(f(x), f(x{'}))$, and
the operation $f$ is $(R_X, S_Y)$-preserving. Although we cannot show that $f$ is a function, we can show
that it is $(=_X, S)$-preserving, since $x =_X x{'} \To R_X(x, x{'})$, and we work as above.
\end{proof}


\section{Families of partial functions}
\label{sec: fampartial}

\begin{definition}\label{def: famofpartial}
Let $X, Y$ and $I$ be sets. A \textit{family of partial functions}\index{family of partial functions}
from $X$ to $Y$ indexed by $I$,
or an $I$-\textit{family of partial functions from $X$ to $Y$}\index{$I$-family of partial functions}, is a triplet
$\Lambda(X,Y) := (\lambda_0, \C E^X, \lambda_1, \C P^Y)$, where\index{$\Lambda(X,Y)$}\index{$\C P^Y$}
$\Lambda(X) := (\lambda_0, \C E^X, \lambda_1) \in \Fam(I, X)$ and  
$\C P^Y : \bigcurlywedge_{i \in I}\D F\big(\lambda_0(i), Y\big)$ with  
$\C P^Y(i) := \C P_i^Y$, for every $i \in I$, 
such that, for every $(i, j) \in D(I)$, the following inner diagrams commute 
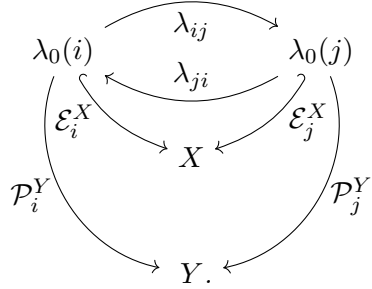
\begin{center}
\begin{tikzpicture}

\node (E) at (0,0) {$\lambda_0(i)$};
\node[right=of E] (B) {};
\node[right=of B] (F) {$\lambda_0(j)$};
\node[below=of B] (A) {$X$};
\node[below=of A] (C) {$ \ Y$.};

\draw[->,bend left] (E) to node [midway,below] {$\lambda_{ij}$} (F);
\draw[->,bend left] (F) to node [midway,above] {$\lambda_{ji}$} (E);
\draw[->,bend right=50] (E) to node [midway,left] {$\C P_i^Y$} (C);
\draw[->,bend left=50] (F) to node [midway,right] {$\C P_j^Y$} (C);
\draw[right hook->,bend right=20] (E) to node [midway,left] {$ \ \C E_i^X \ $} (A);
\draw[left hook->,bend left=20] (F) to node [midway,right] {$\ \C E_j^X \ $} (A);

\end{tikzpicture}
\end{center} 
We call $\C P^Y$\index{$\C P^Y$} a modulus of partial functions\index{modulus of partial functions} for $\lambda_0$, 
and $\Lambda(X)$ the $I$-family of domains of $\Lambda(X,Y)$\index{family of domains of $\Lambda(X,Y)$}. 
If $M(X,Y) := (\mu_0, \C Z^X, \mu_1, \C Q^Y)$ and 
$N(X,Y) := (\nu_0, \C H^X, \nu_1, \C R^Y)$ are $I$-families of partial functions from $X$ to $Y$, a 
\textit{family of partial functions-map}\index{family of partial functions-map}
$\Psi \colon \Lambda(X,Y) \To M(X,Y)$\index{$\Psi \colon \Lambda(X,Y) \To M(X,Y)$}
from $\Lambda(X,Y)$ to $M(X,Y)$ is a dependent operation
$\Psi \colon \bigcurlywedge_{i \in I}\D F\big(\lambda_0(i), \mu_0(i)\big)$, where $\Psi(i) := \Psi_i$,
for every $i \in I$, such that, for every $i \in I$, the following inner diagrams commute
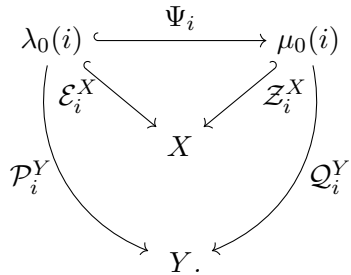
\begin{center}
\begin{tikzpicture}

\node (E) at (0,0) {$\lambda_0(i)$};
\node[right=of E] (B) {};
\node[right=of B] (F) {$\mu_0(i)$};
\node[below=of B] (A) {$X$};
\node[below=of A] (C) {$ \ Y$.};

\draw[right hook->] (E)--(F) node [midway,above] {$\Psi_i$};
\draw[->,bend right=40] (E) to node [midway,left] {$\C P_i^Y$} (C);
\draw[->,bend left=40] (F) to node [midway,right] {$\C Q_i^Y$} (C);
\draw[right hook->] (E)--(A) node [midway,left] {$\C E_i^X \ $};
\draw[left hook->] (F)--(A) node [midway,right] {$\ \C Z_i^X  $};

\end{tikzpicture}
\end{center} 
The totality $\Map_I(\Lambda(X,Y), M(X,Y))$\index{$\Map_I(\Lambda(X,Y), M(X,Y))$} of the family of partial functions-maps 
from $\Lambda(X,Y)$ to $M(X,Y)$ is equipped with the pointwise equality. 
If $\Psi \colon \Lambda(X,Y) \To M(X,Y)$ and if $\Xi \colon M(X,Y) \To N(X,Y)$, the
\textit{composition family of partial functions-map} 
\index{composition family of partial functions-map} 
$\Xi \circ \Psi \colon \Lambda(X,Y) \To N(X,Y)$ is defined by 
$(\Xi \circ \Psi)(i) := \Xi_i \circ \Psi_i$,
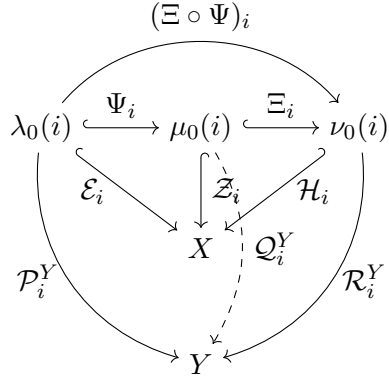
\begin{center}
\begin{tikzpicture}

\node (E) at (0,0) {$\lambda_0(i)$};
\node[right=of E] (B) {$\mu_0(i)$};
\node[right=of B] (F) {$\nu_0(i)$};
\node[below=of B] (C) {$X$};
\node[below=of C] (D) {$Y$};

\draw[right hook->] (E)--(C) node [midway,left] {$\C E_i \ $};
\draw[right hook->] (B)--(C) node [midway,right] {$\C Z_i$};
\draw[left hook->] (F)--(C) node [midway,right] {$ \ \C  H_{i} $};
\draw[right hook->] (E)--(B) node [midway,above] {$\Psi_i$};
\draw[right hook->] (B)--(F) node [midway,above] {$\Xi_i$};
\draw[->,bend left=50] (E) to node [midway,above] {$(\Xi \circ \Psi)_{i}$} (F);
\draw[->,bend right=40] (E) to node [midway,left] {$\C P_i^Y$} (D);
\draw[->,bend left=40] (F) to node [midway,right] {$\C R_i^Y$} (D);
\draw[dashed, ->,bend left] (B) to node [midway,right] {$\C Q_i^Y$} (D);

\end{tikzpicture}
\end{center} 
for every $i \in I$. The identity family of partial functions-map $\Id_{\Lambda(X,Y)} \colon \Lambda(X,Y) \To \Lambda(X,Y)$ 
and the equality on the totality $\Fam(I,X,Y)$\index{$\Cam(I,X,Y)$} of $I$-families of partial functions from 
$X$ to $Y$ are defined as in Definition~\ref{def: map}.
\end{definition}

Clearly, if $\Lambda(X,Y) \in \Fam(I,X,Y)$ and $(i,j) \in D(I)$, then $(\lambda_{ij}, \lambda_{ji}) \colon 
\C P_i^Y =_{\mathsmaller{\C F(X,Y)}} \C P_j^Y$.

\begin{proposition}\label{prp: compfampartial}
Let $\Lambda(X,Y) := (\lambda_0, \C E^X, \lambda_1, \C P^Y) \in \Fam(I,X,Y)$ and let
$M(Y,Z) := (\mu_0, \C Z^Y, \mu_1, \C Q^Z) \in \Fam(I,Y,Z)$. Their composition\index{composition of families of
partial functions} $(M \pcirc \Lambda)(X, Z)$\index{$(M \pcirc \Lambda)(X, Z)$} is defined by
\[ (M \mathsmaller{\pcirc} \Lambda)(X, Z) := \big(\mu_0 \mathsmaller{\pcirc} \lambda_0, \C Z^Y 
\mathsmaller{\pcirc} \C E^X, \mu_1 \mathsmaller{\pcirc} \lambda_1, 
(\C Q \mathsmaller{\pcirc} \C P)^Z\big) \]
\[ (\mu_0 \mathsmaller{\pcirc} \lambda_0)(i) := \big(\C P_i^Y\big)^{-1}(\mu_0(i)); \ \ \ \ i \in I, \]
\[ \big(\C Z^Y \mathsmaller{\pcirc} \C E^X\big)_i := \C E_i^X \circ e_{\mathsmaller{(\C P_i^Y)^{-1}(\mu_0(i))}}^{\lambda_0(i)}
  \colon (\mu_0 \mathsmaller{\pcirc} \lambda_0)(i) \eto X, 
\]
\[ (\mu_1 \mathsmaller{\pcirc} \lambda_1)_{ij} \colon \big(\C P_i^Y\big)^{-1}(\mu_0(i)) \to \big(\C P_j^Y\big)^{-1}(\mu_0(j)), \]
\[ (\mu_1 \mathsmaller{\pcirc} \lambda_1)_{ij}(u,w) := \big(\lambda_{ij}(u), \mu_{ij}(w)\big); \ \ \ \ (u,w) \in 
\big(\C P_i^Y\big)^{-1}(\mu_0(i)),  
\]
\[ (\C Q \mathsmaller{\pcirc} \C P)^Z_i := \C Q_i^Z \mathsmaller{\pcirc} \C P_i^Y, \ \ \ \ i \in I. \]
Then $M(Y,Z) \in \Fam(I, X, Z)$.
\end{proposition}

\begin{proof}
By Definition~\ref{def: image} we have that
\[ \big(\C P_i^Y\big)^{-1}(\mu_0(i)) := \big\{(u,w) \in \lambda_0(i) \times \mu_0(i) \mid \C P_i^Y(u) =_Y \C Z_i^Y(w) \big\}, \]
\[ \big(\C P_j^Y\big)^{-1}(\mu_0(j)) := \big\{(u{'},w{'}) \in \lambda_0(j) \times \mu_0(j) \mid \C P_j^Y(u{'}) =_Y
\C Z_j^Y(w{``}) \big\}. \]
If $(i,j) \in D(I)$, then
$ \C P_j^Y\big(\lambda_{ij}(u)\big) =_Y \C P^Y(u) =_Y \C Z_i^Y(w) =_Y \C Z_j^Y\big(\mu_{ij}(w)\big)$,
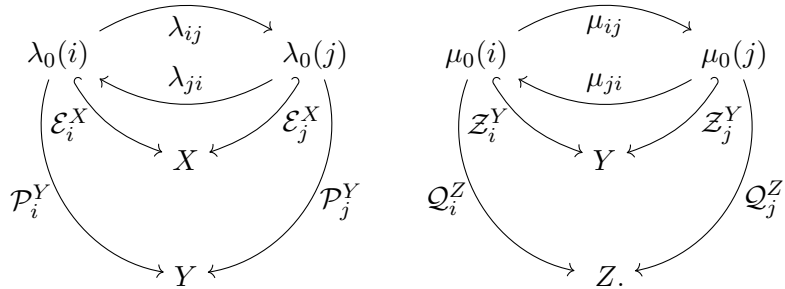
\begin{center}
\begin{tikzpicture}

\node (E) at (0,0) {$\lambda_0(i)$};
\node[right=of E] (B) {};
\node[right=of B] (F) {$\lambda_0(j)$};
\node[below=of B] (A) {$X$};
\node[below=of A] (C) {$Y$};

\node[right=of F] (K) {$\mu_0(i)$};
\node[right=of K] (L) {};
\node[right=of L] (M) {$\mu_0(j)$};
\node[below=of L] (N) {$Y$};
\node[below=of N] (J) {$ \ Z$.};

\draw[->,bend left] (E) to node [midway,below] {$\lambda_{ij}$} (F);
\draw[->,bend left] (F) to node [midway,above] {$\lambda_{ji}$} (E);
\draw[->,bend right=50] (E) to node [midway,left] {$\C P_i^Y$} (C);
\draw[->,bend left=50] (F) to node [midway,right] {$\C P_j^Y$} (C);
\draw[right hook->,bend right=20] (E) to node [midway,left] {$ \ \C E_i^X \ $} (A);
\draw[left hook->,bend left=20] (F) to node [midway,right] {$\ \C E_j^X \ $} (A);

\draw[->,bend left] (K) to node [midway,below] {$\mu_{ij}$} (M);
\draw[->,bend left] (M) to node [midway,above] {$\mu_{ji}$} (K);
\draw[->,bend right=50] (K) to node [midway,left] {$\C Q_i^Z$} (J);
\draw[->,bend left=50] (M) to node [midway,right] {$\C Q_j^Z$} (J);
\draw[right hook->,bend right=20] (K) to node [midway,left] {$ \ \C Z_i^Y \ $} (N);
\draw[left hook->,bend left=20] (M) to node [midway,right] {$\ \C Z_j^Y \ $} (N);

\end{tikzpicture}
\end{center} 
hence the operation $(\mu_1 \mathsmaller{\pcirc} \lambda_1)_{ij}$ is well-defined, and it is immediate to show that it 
is a function. For the commutativity of the following inner diagrams we have that  
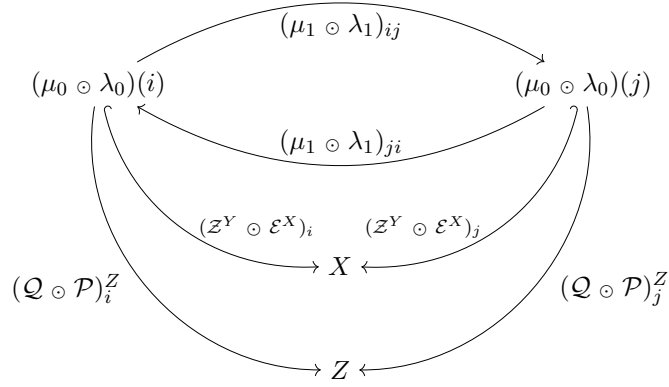
\begin{center}
\resizebox{9cm}{!}{%
\begin{tikzpicture}

\node (E) at (0,0) {$(\mu_0 \mathsmaller{\pcirc} \lambda_0)(i)$};
\node[right=of E] (L) {};
\node[right=of L] (M) {};
\node[right=of M] (B) {};
\node[right=of B] (F) {$(\mu_0 \mathsmaller{\pcirc} \lambda_0)(j)$};
\node[below=of M] (N) {};
\node[below=of N] (A) {$X$};
\node[below=of A] (C) {$Z$};

\draw[->,bend left] (E) to node [midway,below] {$(\mu_1 \mathsmaller{\pcirc} \lambda_1)_{ij}$} (F);
\draw[->,bend left] (F) to node [midway,above] {$(\mu_1 \mathsmaller{\pcirc} \lambda_1)_{ji}$} (E);
\draw[->,bend right=50] (E) to node [midway,left] {$(\C Q \mathsmaller{\pcirc} \C P)^Z_i \ \ $} (C);
\draw[->,bend left=50] (F) to node [midway,right] {$ \ \ (\C Q \mathsmaller{\pcirc} \C P)^Z_j$} (C);
\draw[right hook->,bend right=40] (E) to node [midway,right] {$ \ \mathsmaller{(\C Z^Y \mathsmaller{\pcirc} \C E^X)_i} $} (A);
\draw[left hook->,bend left=40] (F) to node [midway,left] {$\mathsmaller{(\C Z^Y \mathsmaller{\pcirc} \C E^X)_j} \ $} (A);

\end{tikzpicture}
}
\end{center}  
\begin{align*}
 \big(\C Z^Y \mathsmaller{\pcirc} \C E^X\big)_j\big((\mu_1 \mathsmaller{\pcirc} \lambda_1)_{ij}(u,w)\big) & :=
 \big(\C Z^Y \mathsmaller{\pcirc} \C E^X\big)_j\big(\lambda_{ij}(u), \mu_{ij}(w)\big)\\
 & := \bigg[\C E_j^X \circ e_{\mathsmaller{(\C P_j^Y)^{-1}(\mu_0(j))}}^{\lambda_0(j)}\bigg]
 \big(\lambda_{ij}(u), \mu_{ij}(w)\big)\\
 & := \C E_j^X\big(\lambda_{ij}(u)\big)\\
 & =_X \C E_i^X(u)\\
 & := \bigg[\C E_i^X \circ e_{\mathsmaller{(\C P_i^Y)^{-1}(\mu_0(i))}}^{\lambda_0(i)}\bigg](u,w)\\
 & := \big(\C Z^Y \mathsmaller{\pcirc} \C E^X\big)_i(u,w),
\end{align*}
\begin{align*}
 \big(\C Q \mathsmaller{\pcirc} \C P\big)^Z_j\big(\lambda_{ij}(u), \mu_{ij}(w)\big) & := 
 \big(\C Q_j^Z \mathsmaller{\pcirc} \C P_j^Y\big)\big(\lambda_{ij}(u), \mu_{ij}(w)\big)\\
 & := \C Q_j^Z\big(\mu_{ij}(w)\big)\\
 & =_Z \C Q_i^Z(w)\\
 & := \big(\C Q_i^Z \mathsmaller{\pcirc} \C P_i^Y\big)(u,w)\\
 & := \big(\C Q \mathsmaller{\pcirc} \C P\big)^Z_i(u,w).
\end{align*}
For the other two inner diagrams we proceed similarly.
\end{proof}

The basic properties of the composition of partial functions extend to equalities for the corresponding
families of partial functions. E.g., we get
\[ N(Z,W) \mathsmaller{\pcirc} \big[M(Y,Z) \mathsmaller{\pcirc} \Lambda(X,Y)\big] =_{\Fam(I,X,W)} \big[N(Z,W) 
\mathsmaller{\pcirc} M(Y,Z)\big] \mathsmaller{\pcirc} \Lambda(X,Y). \]
Suppose that  $\Lambda(X,Y) := (\lambda_0, \C E^X, \lambda_1, \C P^Y) \in \Fam(I,X,Y)$
and $M(X,Y) := (\mu_0, \C Z^X, \mu_1, \C Q^Y) \in \Fam(I,X,Y)$. We can define in the expected way the following
families of partial functions:
\[ (\Lambda \cap_l M)(X,Y) := \big(\lambda_0 \cap_l \mu_0, \C E^X \cap_l \C Z^X, \lambda_1 \cap_l \mu_1, 
 (\C P \cap_l \C Q)^Y\big), 
\]
\[ (\Lambda \cap_r M)(X,Y) := \big(\lambda_0 \cap_r \mu_0, \C E^X \cap_r \C Z^X, \lambda_1 \cap_r \mu_1, 
 (\C P \cap_r \C Q)^Y\big), 
\]
\[ (\Lambda \cup M)(X,Y) := \big(\lambda_0 \cup \mu_0, \C E^X \cup \C Z^X, \lambda_1 \cup \mu_1, 
 (\C P \cup \C Q)^Y\big).
\]
The basic properties of the intersections and union of partial functions extend to equalities for the corresponding
families of partial functions. E.g., we get
\[ (\Lambda \cup M)(X,Y) =_{\Fam(I,X,Y)} (M \cup \Lambda)(X,Y). \]

Various notions and results on families of subsets extend to families of partial functions.

\section{Families of complemented subsets}
\label{sec: famcomplsubsets}

\begin{definition}\label{def: famofcompl}
Let the sets $(X, =_X, \neq_X)$ and $(I, =_I)$. 
A \textit{family of complemented subsets}\index{family of complemented subsets}
of $X$ indexed by $I$, or an $I$-\textit{family of complemented subsets of $X$}\index{$I$-family of complemented subsets},
is a structure\index{$\B \Lambda(X)$} 
$\B \Lambda(X) := \big(\lambda_0^1, \C E^{1,X}, \lambda_1^1, \lambda_0^0, \C E^{0,X}, \lambda_1^0)$, such that 
$\Lambda^1(X) := \big(\lambda_0^1, \C E^{1,X}, \lambda_1^1\big) \in \Fam(I, X)$ and 
$\Lambda^0(X) := \big(\lambda_0^0, \C E^{0,X}, \lambda_1^0\big) \in \Fam(I, X)$
i.e., for every $(i, j) \in D(I)$, the following inner diagrams commute
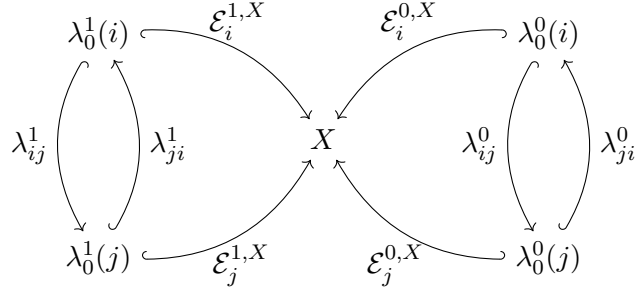
\begin{center}
\begin{tikzpicture}

\node (E) at (0,0) {$X$};
\node[above=of E] (F) {};
\node[right=of F] (P) {};
\node[right=of P] (G) {$\lambda_0^0(i)$};
\node[below=of G] (S) {};
\node[below=of S] (H) {$\lambda_0^0(j)$};
\node[left=of F] (X) {};
\node[left=of X] (J) {$\lambda_0^1(i)$};
\node[below=of J] (T) {};
\node[below=of T] (K) {$\lambda_0^1(j)$};

\draw[right hook->,bend right] (K) to node [midway,below] {$\C E_j^{1,X}$} (E);
\draw[left hook->,bend left] (J) to node [midway,above] {$\C E_i^{1,X}$} (E);
\draw[right hook->,bend right] (G) to node [midway,above] {$\C E_i^{0,X}$} (E);
\draw[left hook->,bend left] (H) to node [midway,below] {$\C E_j^{0,X} \ \ $} (E);
\draw[right hook->,bend right] (J) to node [midway,left] {$\lambda_{ij}^1$} (K);
\draw[right hook->,bend right] (K) to node [midway,right] {$\lambda_{ji}^1$} (J);
\draw[right hook->,bend right] (G) to node [midway,left] {$\lambda_{ij}^0$} (H);
\draw[right hook->,bend right] (H) to node [midway,right] {$\lambda_{ji}^0$} (G);

\end{tikzpicture}
\end{center}
such that 
\[ \forall_{i \in I}\big(\B \lambda_0(i) := \big(\lambda_0^1(i), \lambda_0^0(i)\big) \in \C P^{\Disj}(X)\big) \]
If $\B M(X) := (\mu_0^1, \C Z^{1,X}, \mu_1^1, \mu_0^0, \C Z^{0,X}, \mu_1^0)$, 
$\B N(X) := (\nu_0^1, \C H^{1,X}, \nu_1^1, \nu_0^0, \C H^{0,X}, \nu_1^0)$ are $I$-families of complemented 
subsets of $X$, a \textit{family of complemented subsets-map}\index{family of complemented subsets-map}
$\B \Psi \colon \B \Lambda(X) \To \B M(X)$\index{$\B \Psi \colon \B \Lambda(X) \To \B M(X)$}
from $\B \Lambda(X)$ to $\B M(X)$ is a pair $\B \Psi := (\Psi^1, \Psi^0)$, where $\Psi^1 \colon 
\Lambda^1(X) \To M^1(X)$ and $\Psi^0 \colon \Lambda^0(X) \To M^0(X)$ i.e., 
for every $i \in I$, the following inner diagrams commute
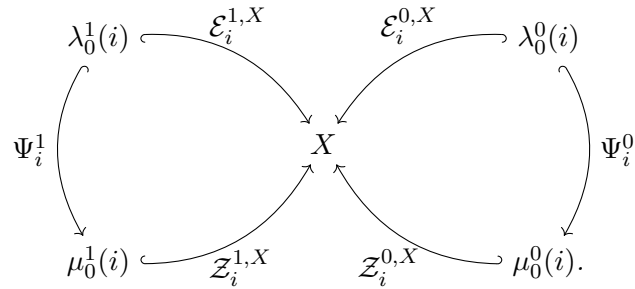
\begin{center}
\begin{tikzpicture}

\node (E) at (0,0) {$X$};
\node[above=of E] (F) {};
\node[right=of F] (P) {};
\node[right=of P] (G) {$\lambda_0^0(i)$};
\node[below=of G] (S) {};
\node[below=of S] (H) {$\mu_0^0(i)$.};
\node[left=of F] (X) {};
\node[left=of X] (J) {$\lambda_0^1(i)$};
\node[below=of J] (T) {};
\node[below=of T] (K) {$\mu_0^1(i)$};

\draw[right hook->,bend right] (K) to node [midway,below] {$\C Z_i^{1,X}$} (E);
\draw[left hook->,bend left] (J) to node [midway,above] {$\C E_i^{1,X}$} (E);
\draw[right hook->,bend right] (G) to node [midway,above] {$\C E_i^{0,X}$} (E);
\draw[left hook->,bend left] (H) to node [midway,below] {$\C Z_i^{0,X} \ \ \ $} (E);
\draw[right hook->,bend right] (J) to node [midway,left] {$\Psi_i^1$} (K);
\draw[left hook->,bend left] (G) to node [midway,right] {$\Psi_i^0$} (H);

\end{tikzpicture}
\end{center}
The totality $\Map_I(\B \Lambda(X), \B M(X))$\index{$\Map_I(\B \Lambda(X), \B M(X))$} of the
family of complemented subsets-maps 
from $\B \Lambda(X)$ to $\B M(X)$ is equipped with the pointwise equality. 
If $\B \Psi \colon \B \Lambda(X) \To \B M(X)$ and if $\B \Xi \colon \B M(X) \To \B N(X)$, the
\textit{composition family of complemented subsets-map} 
\index{composition family of complemented subsets-map} 
$\B \Xi \circ \B \Psi \colon \B \Lambda(X) \To \B N(X)$ is defined by $\B \Xi := \big((\Xi \circ \Psi)^1, 
(\Xi \circ \Psi)^0\big)$, where $(\Xi \circ \Psi)^1 := \Xi^1 \circ \Psi^1$ and 
$(\Xi \circ \Psi)^0 := \Xi^0 \circ \Psi^0$.
Moreover, $\Id_{\B \Lambda(X)} := (\Id_{\Lambda^1(X)}, \Id_{\Lambda^0(X)})$, and the totality $\Fam(I, \B X)$ of families of complemented subsets of $X$ over $I$\index{$\Fam(I, \B X)$} is equipped with the equality
$\B \Lambda(X) =_{\mathsmaller{\Fam(I, \B X)}} \B M(X) $ if and only if 
\[ \exists_{\B \Psi \in \Map_I(\B \Lambda(X), \B M(X))}
\exists_{\B \Xi \in \Map_I(\B M(X), \B \Lambda(X))}\big( \B \Psi \circ \B \Xi = \Id_{\B M(X)} \ \& \ \B \Xi \circ \B \Psi = 
\Id_{\B \Lambda(X)}\big). \] 
\end{definition}

As in the case of $\Fam(I, X)$, we see no reason not to consider $\Fam(I, \B X)$ a set.\index{$\Set(I, \B X)$}
Clearly, the obviously defined set $\Eq(\B \Lambda(X), \B M(X))$ is a subsingleton. A family $\B \Lambda(X) \in \Fam(I, \B X)$ is in $\Set(I, \B X)$, if $\B \lambda_0(i) =_{\C P^{\Disj}(X)} \B \Lambda_0(j) \To i =_I j$, for every $i, j \in I$. Trivially, if $\Lambda^1(X) \in \Set(I, X)$, or if $\Lambda^0(X) \in \Set(I, X)$, then $\B \Lambda(X) \in \Set(I, \B X)$. Clearly, if $\B \Lambda(X) \in \Set(I, \B X)$ and $\B M(X) \in \Fam(I, \B X)$ such that $\B M(X) =_{\mathsmaller{\Fam(I, \B X)}} \B \Lambda(X)$, then $\B M(X) \in \Set(I, \B X)$.
The operations between complemented subsets induce new families of complemented subsets and family-maps between them. If
$\B \Lambda(X) := \big(\lambda_0^1, \C E^{1,X}, \lambda_1^1, \lambda_0^0, \C E^{0,X}, \lambda_1^0\big)$ and
$\B M(X) := \big(\mu_0^1, \C Z^{1,X}, \mu_1^1, \mu_0^0, \C Z^{0,X}, \mu_1^0\big) \in \Fam(I, \B X)$, let the following new elements of $\Fam(I, \B X)$:
\[(- \B \Lambda)(X) := \big(\lambda_0^0, \C E^{0,X}, \lambda_1^0, \lambda_0^1, \C E^{1,X}, \lambda_1^1\big), \]
\[ (\B \Lambda \cap \B M)(X) := \big(\lambda_0^1 \cap \mu_0^1, \C E^{1,X} \cap \C Z^{1,X}, \lambda_1^1 \cap \mu_1^1, \lambda_0^0 \cup \mu_0^0, \C E^{0,X} \cup \C Z^{0,X} , \lambda_1^0 \cup \mu_1^0\big), \]
\[ (\B \Lambda \cup \B M)(X) := \big(\lambda_0^1 \cup \mu_0^1, \C E^{1,X} \cup \C Z^{1,X}, \lambda_1^1 \cup \mu_1^1, \lambda_0^0 \cap \mu_0^0, \C E^{0,X} \cap \C Z^{0,X} , \lambda_1^0 \cap \mu_1^0\big), \]
\[ (\B \Lambda - \B M)(X) := [\B \Lambda \cap (- \B M)](X). \]
If $\B N(Y) := (\nu_0^1, \C H^{1,Y}, \nu_1^1, \nu_0^0, \C H^{0,Y}, \nu_1^0) \in \Fam(I, \B Y)$ and $\fXY$, then using  
Proposition~\ref{prp: newfamsofsubsets2} we define
\[ f^{-1}(\B N)(X) := \big(f^{-1}(\nu_0^1), f^{-1}\big(\C H^{1,Y}\big)^X, f^{-1}(\nu_1^1), f^{-1}(\nu_0^0), f^{-1}\big(\C H^{0,Y}\big)^X, f^{-1}(\nu_1^0)\big) \in \Fam(I, \B X), \]
\[ f(\B \Lambda)(Y) := \big(f^{-1}(\lambda_0^1), f\big(\C E^{1,X}\big)^Y, f(\lambda_1^1), f(\lambda_0^0), f\big(\C E^{0,X}\big)^Y, f(\lambda_1^0)\big) \in \Fam(I, \B Y). \]
Properties between complemented subsets induce equalities between their families e.g.,
\[ \big[f^{-1}\big[(\B N \cup \B K)(Y)\big]\big](X) =_{\Fam(I, \B X)} \big[f^{-1}(\B N)(Y) \cup f^{-1}(\B K)(Y)\big](X). \]
Using definitions from section~\ref{sec: famofsubsets}, if $\B \Phi \colon \B \Lambda(X) \To \B M(X)$, let 
$- \B \Phi \colon (- \B \Lambda)(X) \To (- \B M)(X)$ 
\[ - \B \Phi := (\Phi^0, \Phi^1); \ \ \ \ \B \Phi := (\Phi^1, \Phi^0). \]
If $\B \Psi \colon \B P(X) \To \B Q(X)$, then $\B \Phi \cap \B \Psi \colon (\B \Lambda \cap \B P)(X) \To (\B M \cap \B R)(X)$, where 
\[ \B \Phi \cap \B \Psi := (\Phi^1 \cap \Psi^1, \Phi^0 \cup^0), \]
and $\B \Phi \cup \B \Psi \colon (\B \Lambda \cup \B P)(X) \To (\B M \cup \B R)(X)$, where 
\[ \B \Phi \cup \B \Psi := (\Phi^1 \cup \Psi^1, \Phi^0 \cap^0), \]
and $\B \Phi - \B \Psi \colon (\B \Lambda - \B P)(X) \To (\B M - \B R)(X)$, where 
$\B \Phi - \B \Psi := \B \Phi \cap (- \B \Psi)$. If $\B S(Y) \in \Fam(J, Y)$, 
\[ (\B \Lambda \times \B S)(X \times Y) := \big(\lambda_0^1 \times s_0^1, \C E^{1,X} \times \C S^{1,Y}, \lambda_1^1 \times s_1^1, \lambda_0^0 \times s_0^0, \C E^{0,X} \times \C S^{0,Y}, \lambda_0^0 \times s_0^0\big) \in \Fam(I \times J, \B X \times \B Y), \]
\[ (\B \lambda_0 \times \B s_0)(i, j) := \B \lambda_0(i) \times \B s_0(i). \]
If $\B \Xi \colon \B S(Y) \To \B T(Y)$, then $\B \Phi \times \B \Xi \colon (\B \Lambda \times \B S)(X \times Y) \To (\B M \times \B T)(X \times Y)$, where 
\[ \B \Phi \times \B \Xi := \big(\Phi^1 \times \Xi^1, \Phi^0 \times \Xi^0\big). \]
Due to the above families of complemented subsets the following proposition is well-formulated.

\begin{proposition}\label{prp: unioncompl}
Let $\B \Lambda(X) := \big(\lambda_0^1, \C E^{1,X}, \lambda_1^1, \lambda_0^0, \C E^{0,X}, \lambda_1^0\big) \in \Fam(I, \B X)$, $i_0 \in I$, and let
\[ \bigcup_{i \in I} \B \lambda_0 (i) := \bigg(\bigcup_{i \in I}\lambda_0^1(i), \bigcap_{i \in I}\lambda_0^0(i)\bigg) 
\ \ \ \ \& \ \ \ \ \bigcap_{i \in I} \B \lambda_0 (i) := \bigg(\bigcap_{i \in I}\lambda_0^1(i), 
\bigcup_{i \in I}\lambda_0^0(i)\bigg). \]
\normalfont (i) 
\itshape $\bigcup_{i \in I} \B \lambda_0 (i), \bigcap_{i \in I} \B \lambda_0 (i) \in \C P^{\Disj}(X)$.\\[1mm]
\normalfont (ii) 
\itshape $- \bigcup_{i \in I} \B \lambda_0 (i) =_{\mathsmaller{\C P^{\Disj}(X)}} \bigcap_{i \in I} \big(- \B \lambda_0 (i)\big)$.\\[1mm]
\normalfont (iii) 
\itshape $- \bigcap_{i \in I} \B \lambda_0 (i) =_{\mathsmaller{\C P^{\Disj}(X)}} \bigcup_{i \in I} 
\big(- \B \lambda_0 (i)\big)$.\\[1mm]
\normalfont (iv)
\itshape If $i \in I$, then $\B \lambda_0(i) \subseteq \bigcup_{i \in I}\B \lambda_0(i)$.\\[1mm]
\normalfont (v)
\itshape If 
$\B A \subseteq \B \lambda_0(i)$, for some $i \in I$, then 
$\B A \subseteq \bigcup_{i \in I}\B \lambda_0(i)$.\\[1mm]
\normalfont (vi)
\itshape If $\B \lambda_0(i) \subseteq \B A$, for every $i \in I$, then $\bigcup_{i \in I} \B \lambda_0(i) \subseteq \B A$.\\[1mm]
\normalfont (vii)
\itshape If $\B \lambda_0(i) \supseteq \B A$, for every $i \in I$, then $\bigcap_{i \in I} \B \lambda_0(i) \supseteq \B A$.\\[1mm]
\normalfont (viii) 
\itshape If $\B M(X) := \big(\mu_0^1, \C Z^{1,X}, \mu_1^1, \mu_0^0, \C Z^{0,X}, \mu_1^0) \in \Fam(I, \B Y)$ and $\fXY$, 
then
\[ f^{-1}\bigg(\bigcup_{i \in I} \B \mu_0 (i)\bigg) =_{\mathsmaller{\C P^{\Disj}(X)}} \bigcup_{i \in I}f^{-1}\big(\B \mu_0 (i)\big) \ \ \ \ \& \ \ \ \ 
f^{-1}\bigg(\bigcap_{i \in I} \B \mu_0 (i)\bigg) =_{\mathsmaller{\C P^{\Disj}(X)}} \bigcap_{i \in I}f^{-1}\big(\B \mu_0 (i)\big). \]
\end{proposition}

\begin{proof}
(i) We show the first membership only. If $(i, x) \in \bigcup_{i \in I}\lambda_0^1(i)$
and $\Phi \in \bigcap_{i \in I}\lambda_0^0(i)$, then $e_{\mathsmaller{\bigcup}}^{\Lambda(X)}(i, x)
:= \C E_i^X(x)$ and $e_{\mathsmaller{\bigcap}}^{\Lambda(X)}(\Phi) := \C E_{i_0}^X(\Phi_{i_0})$. 
Since $\C E_i^X(\Phi_i) =_X \C E_{i_0}^X(\Phi_{i_0})$ and $\lambda_0^1(i) \Disj \lambda_0^0(i)$,
we have that $\C E_i^X(x) \neq_X \C E_i^X(\Phi_i)$, and by the extensionality of $\neq_X$ we get
$\C E_i^X(x) \neq_X \C E_{i_0}^X(\Phi_{i_0})$.\\
(ii) and (iii) are straightforward to show. For (iv) we need to show that $\lambda_0^1(i) \subseteq 
\bigcup_{i \in I}\lambda_0^1(i)$ and $\bigcap_{i \in I}\lambda_0^0(i) \subseteq \lambda_0^0(i)$, 
which follow from Propositions~\ref{prp: unionmap1}(ii) and~\ref{prp: internalmap1}(ii), respectively.
Case (v) follows from (iv) and the transitivity of $\B A \subseteq \B B$.\\[1mm]
(vi) If $\lambda_0^1(i) \subseteq A^1$, for every $i \in I$, then $\bigcup_{i \in I}\lambda_0^1(i) \subseteq
A^1$, and if $A^0 \subseteq \lambda_0^0(i)$, for every $i \in I$, then $A^0 \subseteq \bigcap_{i \in I}\lambda_0^0(i)$.
Case (vii) is shown similarly.\\
(viii) We show the first equality only. By Propositions~\ref{prp: interiorunion3} and~\ref{prp: intersection3} 
we have that 
\begin{align*}
f^{-1}\bigg(\bigcup_{i \in I} \B \mu_0 (i)\bigg) & := \bigg(f^{-1}\bigg(\bigcup_{i \in I}\mu_0^1(i)\bigg),
f^{-1}\bigg(\bigcap_{i \in I}\mu_0^0(i)\bigg)\bigg)\\
& =_{\mathsmaller{\C P^{\Disj}(X)}} \bigg(\bigcup_{i \in I}f^{-1}\big(\mu_0^1(i)\big), 
\bigcap_{i \in I}f^{-1}\big(\mu_0^0(i)\big)\bigg)  \\
& := \bigcup_{i \in I}\bigg(f^{-1}\big(\mu_0^1(i)\big), f^{-1}\big(\mu_0^0(i)\big)\bigg)\\
& := \bigcup_{i \in I}f^{-1}\big(\B \mu_0 (i)\big).\qedhere
\end{align*}
\end{proof}

Let $\B \Lambda(X), \B M(X)$ and $\Psi \colon \B \Lambda(X) \To \B M(X)$. Since $\Psi^1 \colon 
\Lambda^1(X) \To M^1(X)$ and $\Psi^0 \colon \Lambda^0(X) \To M^0(X)$, the following maps between 
complemented subsets (see Definition~\ref{def: complementedsubset}) are defined
\[ \bigcup \B \Psi := \big( \bigcup \Psi^1, \bigcap \Psi^0\big) \colon \bigcup_{i \in I} \B \lambda_0(i) \to 
\bigcup_{i \in I}\B \mu_0(i),  \]
\[ \bigcap \B \Psi := \big( \bigcap \Psi^0, \bigcup \Psi^1\big) \colon \bigcap_{i \in I} \B \lambda_0(i) \to 
\bigcap_{i \in I}\B \mu_0(i),  \ \ \ \ \mbox{where} \]
\[ \bigcup \Psi^1 \colon \bigcup_{i \in I} \lambda_0^1(i) \to \bigcup_{i \in I}\mu_0^1(i) \ \ \ \& \ \ \ 
\bigcap \Psi^0 \colon \bigcap_{i \in I} \lambda_0^0(i) \to \bigcap_{i \in I}\mu_0^0(i) \]
are defined according to Proposition~\ref{prp: unionmap1}(ii) and~\ref{prp: internalmap1}(ii).

\begin{proposition}\label{prp: prodfam}
Let $\B A \in \C P^{\Disj}(X), \B B \in \C P^{\Disj}(Y)$, $\B \Lambda(X) \in \Fam(I, \B X)$ and $\B M(Y) \in \Fam(J, \B Y)$. The following properties hold:
\[ \B A \times \bigcup_{j \in J}\B \mu_0(j) =_{\C P^{\Disj}(X \times Y)} \bigcup_{j \in J} (\B A \times \B \mu_0(j)), \]
\[ \B A \times \bigcap_{j \in J}\B \mu_0(j) =_{\C P^{\Disj}(X \times Y)} \bigcap_{j \in J} (\B A \times \B \mu_0(j)), \]
\[ \bigg(\bigcup_{i \in I}\B \lambda_0(i)\bigg) \times \B B=_{\C P^{\Disj}(X \times Y)} \bigcup_{i \in I} (\B \lambda_0(i) \times \B B), \]
\[ \bigg(\bigcap_{i \in I}\B \lambda_0(i)\bigg) \times \B B=_{\C P^{\Disj}(X \times Y)} \bigcap_{i \in I} (\B \lambda_0(i) \times \B B). \]
\end{proposition}

\begin{proof}
We show the first equality, and for the rest we proceed similarly. By the equalities shown after 
Propositions~\ref{prp: famsubsetsprod2} 
 we have that   
\begin{align*}
\B A \times \bigcup_{j \in J}\B \mu_0(j) & := \bigg(A^1 \times \bigcup_{j \in J}\mu_0^1(j), (A^0 \times Y) \cup \bigg[X \times \bigg(\bigcap_{j \in J}\mu_0^0(j)\bigg)\bigg]\bigg)\\
& =_{\C P^{\Disj}(X \times Y)} \bigg(\bigcup_{j \in J}(A^1 \times \mu_0^1(j)), (A^0 \times Y) \cup \bigg(\bigcap_{j \in J}(X \times \mu_0^0(j))\bigg)\bigg)\\
& =_{\C P^{\Disj}(X \times Y)} \bigg(\bigcup_{j \in J}(A^1 \times \mu_0^1(j)), \bigcap_{j \in J}\big[ (A^0 \times Y) \cup (X \times \mu_0^0(j))\big]\bigg)\\
& :=  \bigcup_{j \in J} (\B A \times \B \mu_0(j)).\qedhere
\end{align*}
\end{proof}


\section{Direct families of subsets}
\label{sec: dirfamsubsets}

 \begin{definition}\label{def: dirfamsubsets}
 Let $(I, \lt)$ be a directed set, and $X \in \D V_0$. A $($covariant$)$ \textit{direct family of subsets}
 \index{direct family of subsets} of $X$ indexed by $I$, or 
 an $(I, \lt)$-\textit{family of subsets}\index{$(I, \lt)$-family of subsets} of $X$, is a triplet 
 $\Lambda^{\lt}(X) := (\lambda_0, \C E^X, \lambda_1^{\lt})$, where
 $\lambda_0 : I \sto \D V_0$, $\C E^X$ is a modulus of embeddings for $\lambda_0$ $($see 
 Definition~\ref{def: famofsubsets}$)$
 \[ \lambda_1^{\preccurlyeq} : \bigcurlywedge_{(i, j) \in \ltI}\D F\big(\lambda_0(i), \lambda_0(j)\big), \ \ \ 
 \lambda_1^{\lt}(i, j) := \lambda_{ij}^{\lt}, \ \ \ (i, j) \in  \ \ltI, \]
 a modulus of covariant transport maps for $\lambda_0$, such that 
 \index{modulus of covariant transport maps}
 $\lambda_{ii} := \id_{\lambda_0(i)}$, for every $i \in I$, and, for every $(i, j) \in \ \ltI$, 
 the following left diagram 
 commutes
 \begin{center}
 \begin{tikzpicture}
 
 \node (E) at (0,0) {$\lambda_0(i)$};
 \node[right=of E] (B) {};
 \node[right=of B] (F) {$\lambda_0(j)$};
 \node[right=of F] (K) {$\lambda_0(i)$};
 \node[right=of K] (L) {};
 \node[right=of L] (M) {$\lambda_0(j)$};
 \node[below=of L] (N) {$X$.};
\node[below=of B] (A) {$X$};
 
 \draw[right hook->] (E) to node [midway,above] {$\lambda_{ij}^{\lt}$} (F);
 \draw[right hook->] (E)--(A) node [midway,left] {$\C E_{i} \ $};
 \draw[left hook->] (F)--(A) node [midway,right] {$ \ \C E_j$};
 \draw[left hook->] (M) to node [midway,above] {$\lambda_{ji}^{\mt}$} (K);
 \draw[right hook->] (K)--(N) node [midway,left] {$\C E_{i} \ $};
 \draw[left hook->] (M)--(N) node [midway,right] {$ \ \C E_j$};

 \end{tikzpicture}
 \end{center} 
 A \textit{contravariant} $(I, \mt)$-family of subsets of $X$ is defined dually i.e., 
 \index{contravariant direct family of subsets}
 \[ \lambda_1^{\mt} : \bigcurlywedge_{(i, j) \in \lt(I)}\D F\big(\lambda_0(j), \lambda_0(i)\big), \ \ \ 
 \lambda_1^{\mt}(i, j) := \lambda_{ji}^{\mt}, \ \ \ (i, j) \in \
  \ltI, \]
is a modulus of contravariant transport maps for $\lambda_0$,
\index{modulus of contravariant transport maps}
such that for every $(i, j) \in \ltI$, the above right diagram commutes. 
\end{definition}
 
%
%
%

\begin{proposition}\label{prp: dirfamofsubsetsequiv}
Let $X \in \D V_0$, $(I, \lt_I)$ a directed set, $\lambda_0 : I \sto \D V_0$, $\C E^X$ a modulus 
of embeddings for $\lambda_0$, and $\lambda_1$ a modulus of transport maps for $\lambda_0$. The 
following are equivalent.\\[1mm]
\normalfont (i) 
\itshape $\Lambda^{\mt}(X) := (\lambda_0, \C E^X, \lambda_1^{\lt})$ is an $(I, \lt_I)$-family of subsets of $X$.\\[1mm]
\normalfont (ii) 
\itshape $\Lambda^{\lt} := (\lambda_0, \lambda_1) \in \Fam(I, \lt_I)$ and $\C E^X \colon \Lambda^{\lt} \To C^{\lt,X}$, 
where $C^{\lt,X}$ is the constant $(I, \lt_I)$-family $X$ $($see Definition~\ref{def: dfamilyofsets}$)$.
\end{proposition}

\begin{proof}
We proceed exactly as in the proof of Proposition~\ref{prp: famofsubsetsequiv}.
\end{proof}

If $\Lambda^{\lt} := (\lambda_0, \C E^X, \lambda_1^{\lt})$ is an $(I, \lt_I)$-family of subsets of $X$, and 
if $i \lt_I j$, then $\lambda^{\lt}_{ij} \colon \lambda_0(i) \subseteq \lambda_0(j)$ i.e., $\lambda_1^{\lt}$
is a \textit{modulus of subset-witnesses} for $\lambda_0$.

\begin{definition}\label{def: dirsubfammap}
If $\Lambda^{\lt}(X) := (\lambda_0, \C E^X, \lambda_1^{\lt}) and  M^{\lt}(X) := (\mu_0, \C Z^X, \mu_1^{\lt})$
 are $(I, \lt_I)$-families of subsets of $X$, a 
\textit{direct family of subsets-map}\index{direct family of subsets-map}
$\Psi \colon \Lambda^{\lt}(X) \To M^{\lt}(X)$\index{$\Psi \colon \Lambda^{\lt}(X) \To M^{\lt}(X)$}
from $\Lambda^{\lt}(X)$ to $M^{\lt}(X)$ is a family of subsets-map $\Phi \colon \Lambda(X) \To M(X)$.
Their set $\Map_{(I, \lt_I)}(\Lambda^{\lt}(X), M^{\lt}(X))$\index{$\Map_{(I, \lt_I)}(\Lambda^{\lt}(X), 
M^{\lt}(X))$} is the set 
$\Map_I(\Lambda(X), M(X))$. The composition of direct family of subsets-maps, and the totality 
$\Fam(I, \lt_I, X)$ of $(I, \lt_I)$-families of subsets of $X$ are defined as the composition of
family of subsets-maps, and as the totality $\Fam(I, X)$, respectively. The
totality $\Fam(I, \mt, X)$\index{$\Fam(I, \mt, X)$} of contravariant direct families of subsets
of $X$ over $(I, lt_I)$ and the corresponding family-maps are defined similarly.
\end{definition}

 \begin{proposition}\label{prp: dirfamsubmap1}
 Let $\Lambda^{\lt}(X) := (\lambda_0, \C E^X, \lambda_1^{\lt}), M(X) := (\mu_0, \C Z^X, \mu_1^{\lt}) \in \Fam(I, \lt_I, X)$.\\[1mm]
\normalfont (i) 
\itshape If $\Psi \colon \Lambda^{\lt}(X) \To M^{\lt}(X)$, then $\Psi \colon \Lambda^{\lt} \To M^{\lt}$.\\[1mm]
\normalfont (ii) 
\itshape If $\Psi \colon \Lambda^{\lt}(X) \To M^{\lt}(X)$ and $\Phi \colon \Lambda^{\lt}(X) \To M^{\lt}(X)$, then 
$\Phi =_{\Map_{(I, \lt_I)}(\Lambda^{\lt}(X), M^{\lt}(X))} \Psi$.
\end{proposition}

\begin{proof}
We proceed exactly as in the proof of Proposition~\ref{prp: famsubmap1} 
\end{proof}

 The interior union and intersection of $\Lambda^{\lt}(X) (\Lambda^{\mt}(M))$, are defined as for  
 an $I$-family of subsets $\Lambda(X)$. As in the case of $\sum_{i \in I}\lambda_0(i)$ and $\bigcup_{i \in I}\lambda_0(i)$, the equality of $\bigcup_{i \in I}\lambda_0(i)$ does not imply the externally defined equality of 
 $\sum_{i \in I}^{\lt}\lambda_0(i)$, only the converse is true i.e.,
 \[ (i, x) =_{\mathsmaller{\sum_{i \in I}^{\lt}\lambda_0(i)}} (j, y) \To 
 (i, x) =_{\mathsmaller{\bigcup_{i \in I}\lambda_0(i)}} (j, y), \]
 as, if there is some $k \in I$ such that $i \lt_I k, j \lt_I k$, and $\lambda_{ik}^{\lt}(x) =_{\lambda_0(k)}
 \lambda_{jk}^{\lt}(y)$, then by the equalities $\C E_i = \C E_k \circ \lambda_{ik}^{\lt}$ and 
 $\C E_j = \C E_k \circ \lambda_{jk}^{\lt}$ we get
 $\C E_i(x) = \C E_k\big(\lambda_{ik}^{\lt}(x)\big) = \C E_k\big(\lambda_{jk}^{\lt}(y)\big) =
 \C E_j (y)$.


\section{Notes}
\label{sec: notes4}

\begin{note}\label{not: deffamofsubsets}
\normalfont
The definition of a family of subsets given by Bishop in~\cite{Bi67}, p.~65, was the rough description we gave at 
the beginning of this chapter. Our definition~\ref{def: famofsubsets} highlights the witnessing data of the 
rough description, and it is in complete analogy to Richman's definition of a set-indexed family of sets,
included later by Bishop and Bridges in~\cite{BB85}, p.~78.
In~\cite{BB85}, p.~80, and in~\cite{Bi67}, p.~65, an alternative definition of a family
of subsets of $X$ indexed by $I$ is given, as a subset $\Lambda$ of $X \times I$. The fact 
that  $(x, i) \in \Lambda$ can be interpreted as $x \in \lambda_0(i)$
This definition though, which was never used by Bishop, does not reveal the witnessing data for the 
equality $\lambda_0(i) =_{\C P(X)} \lambda_0(j)$, if $i =_I j$, and it is not possible to connect with
the notion of a family of sets. 
The definition of a set of subsets is given by Bishop in~\cite{Bi67}, p.~65, and 
it is repeated in~\cite{BB85}, p.~69. The example of the set of detachable subsets of a set is given in~\cite{Bi67}, 
p.~65, where the term \textit{free} subsets\index{free subset} is used instead, 
and it is repeated in~\cite{BB85}, p.~70. 
\end{note}

\begin{note}\label{not: exfamofsubsets}
\normalfont
There are many examples of families of subsets in the literature of Bishop-style
constructive mathematics. In topology a neighborhood space (in~\cite{BB85}, p.~75, the reference
to the indices is  omitted for simplicity) is a pair $(X, N)$, where
 $X$ is a set and $N$ is a family $\nu$ of subsets of $X$ indexed by some set $I$ such that
 \[ \forall_{i, j \in I}\forall_{x \in X}\big(x \in \nu(i) \cap \nu(j) \To \exists_{k \in I}\big(x \in 
 \nu(k) \subseteq \nu(i) \cap \nu(j)\big)\big). \]
 The covering property is not mentioned there. 
 If $(X, F)$ is a Bishop space (see~\cite{BB85}, chapter 3, and~\cite{Pe15}),
 the neighborhood structure $N_F$ on $X$ generated by the Bishop topology $F$ on $X$ is the 
 family $U$ of subsets of $X$ indexed by $F$ that assigns to every element $f \in F$ the set
 \[ U(f) := \{x \in X \mid f(x) > 0\}. \]
 If $f = g$, then $U(f) = U(g)$, while the converse is not true (take e.g., $f = \id_{\Real}$
 and $g = 2 \id_{\Real}$, where $X = \Real$ and $F = \BR$). In real analysis 
sequences of bounded intervals of $\Real$ are considered in~\cite{BB85} Problem 1, p.~292.
In the theory of normed linear spaces a sequence of 
 bounded, located, open, convex sets is constructed in the proof of the separation theorem (see~\cite{BB85},
 pp.~336--340). A family $N_t^*$, for every $t > 0$,  of subsets of the unit sphere of the dual space 
 $X^*$ of a separable normed space $X$ occurs in the proof of Theorem (6.8) in~\cite{BB85}, p.~354.
In constructive algebra families of ideals and families of submodules of an $R$-module are studied 
(see~\cite{MRR88}, p.~44, and p.~53, respectively). 

\end{note}

\begin{note}\label{not: interior}
\normalfont
In~\cite{BB85}, p.~69, the interior union $\bigcup_{i \in I} \lambda_0(i)$ is defined as the totality
\[ \bigcup_{i \in I} \lambda_0(i) := \big\{x \in X \mid \exists_{i \in I}\big(x \in \lambda_0(i)\big)\big\}. \]
Using our notation though, in~\cite{BB85}, pp.~69--70 it is written that
\begin{quote} 
$\ldots$ to construct an element $u$ of
$\bigcup_{i \in I} \lambda_0(i)$ we first 
construct an element $i$ of $I$, and then construct an element $x$ of $\lambda_0 (i)$.
\end{quote}
Clearly, what is meant by the totality $\bigcup_{i \in I} \lambda_0(i)$ 
is what is written in Definition~\ref{def: interiorunion}. The intersection of an $
I$-family $\lambda$ of subsets of $X$ is roughly 
defined in~\cite{Bi67}, p.~70, as 
\[ \bigcap_{i \in I} \lambda_0(i) := \big\{x \in X \mid \forall_{i \in I}\big(x \in \lambda_0(i)\big)\big\}, \]
while the more precise definition that follows this simplified notation is different,
and it is based on the undefined in~\cite{Bi67} and~\cite{BB85} 
notion of a dependent operation over $\lambda$, hence it is not that precise. Moreover, the definition of 
$\prod_{i \in I}^{\mathsmaller{X}}\lambda_0 (i)$, given in~\cite{BB85}, p.~70, as the set
\[ \bigg\{f : I \to \bigcup_{i \in I}\lambda_0 (i) \mid \forall_{i \in I}\big(f(i) \in \lambda_0 (i)\big)\bigg\} \]
is not compatible with the precise definition of $\bigcup_{i \in I}\lambda_0 (i)$, and 
it is not included in~\cite{Bi67}.

\end{note}

\begin{note}\label{not: onfamofdisjointsubsets}
\normalfont
One could have defined an $I$-family of disjoint subsets of $X$ with respect to given
inequalities $\neq_I$ and $\neq_X$ (Definition~\ref{def: interiorunion}) by
\[ \forall_{i, j \in I}\big(i \neq_I j \To \neg \big(\lambda_0(i) \between \lambda_0(j)\big)\big), \ \ \ \ \mbox{or} \]
\[ \forall_{i, j \in I}\big(\lambda_0(i) \between \lambda_0(j) \To i =_I j\big). \]
The first definition is negativistic, while the second, which avoids $\neq_I$, is too strong. 
\end{note}

\begin{note}\label{not: onextensioncoverings}
\normalfont
The classical proof of the extension theorem of coverings (Theorem~\ref{thm: extensioncoverings}) 
is based on the definition of the interior union as the set 
$\bigcup_{i \in I} \lambda_0(i) := \big\{x \in X \mid \exists_{i \in I}\big(x \in \lambda_0(i)\big)\big\}$.
As a result, the required function $\fXY$ is defined as follows: If $x \in \bigcup_{i \in I} \lambda_0(i)$, 
there is $i \in I$ such that 
$x \in \lambda_0(i)$. Then, one defines $f(x) = f_i(x)$, and shows that the value $f(x)$ does not depend on the choice 
of $i$ (see~\cite{Du66}, p.~13). The use of choice is avoided in our proof, because of the embedding $e \colon X \eto 
\bigcup_{i \in I} \lambda_0(i)$.
Theorem~\ref{thm: extensioncoverings} is related to the notion of a \textit{sheaf of sets}\index{sheaf of sets}.
The sheaf-property added to the notion of a presheaf\index{presheaf} is exactly the main condition of 
Theorem~\ref{thm: extensioncoverings}, where the covering of $X$ is an open covering i.e., a covering of 
open subsets (see~\cite{Go73}). 

\end{note}

\begin{note}\label{not: converseinclusion}
\normalfont
If $P(I)$ is a partition of $I$, such that $p_0(k) \neq \emptyset$, for every $k \in K$, and if 
\[ T := \prod_{k \in K}p_0(k),\]
then the converse inclusion to the semi-distributivity of $\bigcap$ over $\bigcup$ (Proposition~\ref{prp: subdistr4}) holds classically, and the distributivity of $\bigcap$ over $\bigcup$ holds classically. The converse inclusion to the semi-distributivity of $\bigcap$ over $\bigcup$ is equivalent to the axiom of choice (see~\cite{Du66}, p.~25). It is expected that this converse inclusion is constructively provable only if non-trivial data are added to the hypotheses.

\end{note}

\begin{note}\label{not: onmodulusofcenters}
\normalfont
In the hypothesis of Proposition~\ref{prp: centerfiber} we need to suppose the existence of a modulus of centres of 
contraction to avoid choice in the definition of function $g$. Proposition~\ref{prp: centerfiber} is our translation
of Theorem 4.4.3 of book-HoTT into $\BST$. In the formulation of Theorem 4.4.3 of~\cite{HoTT13} no modulus of centres of
contraction is mentioned, as the type-theoretic axiom of choice is provable in $\MLTT$.
\end{note}

\begin{note}\label{not: ineqequiv}
 \normalfont
 As an equivalence structure $(X, R_X)$ is the analogue to the set $(X, =_X)$, one can equip $(X, R_X)$ 
 with an extensional relation $I_X$ on $X \times Y$ satisfying the properties of an inequality. 
 In this way the structure $(X, R_X, I_X)$ becomes the equivalence relation-analogue to the set $(X, =_X, \neq_X)$.
  \end{note}

\begin{note}\label{not: examplesfamspartial}
\normalfont
Examples of families of partial functions are found in the predicative reconstruction of the Bishop-Cheng measure 
theory in~\cite{Ze19} and~\cite{PZ20}.
\end{note}

\begin{note}\label{not: exfamofcomplsubsets}
\normalfont
There are many examples of families of complemented subsets in the literature of Bishop-style
constructive mathematics. In the theory of normed linear spaces, sequences of complemented subsets
occur in the formulation
 of the constructive version of Lebesgue's decomposition of measures (see~\cite{BB85}, pp.~329--331), and in the 
 formulation of the constructive Radon-Nikodym theorem (see~\cite{BB85}, pp.~333--334).
 In the integration theory of~\cite{BB85}, the sequences of integrable sets in an 
integrable space $X$ (see~\cite{BB85}, pp.~234--235) are families of subsets
of $X$ indexed by $\Nat$. Sequences of measurable sets are considered in~\cite{BB85}, pp.~269--271.
Moreover, a measure space (see~\cite{BB85}, p.~282) is defined as a triplet $(X, M, \mu)$, where 
$M$ is a set of complemented sets in an inhabited set $X$. In the definition of complete measure
space in~\cite{BB85}, pp.~288--289, the notion of a sequence of elements of $M$ is also used.

\end{note}

\begin{note}\label{not: Myhill1}
\normalfont
In the measure theory developed in~\cite{Bi67} certain families $\C F$ (and subfamilies $\C M$ of $\C F$)
of complemented subsets of some set $X$
are considered in the definition of a measure space (see~\cite{Bi67}, p. 183). 
For the definition of a measure space found in~\cite{Bi67}, p.~183, Myhill writes
in~\cite{My75}, p. 351, the following:
\begin{quote}
 The only one of the classical set-existence axioms (not counting choice) which is 
 missing\footnote{He means from his system $\CST$.} is power set. Certainly there is no hint of this axiom
 in Bishop's book (except for $\C F$ on p.~183, \textit{surely a slip\footnote{Our emphasis.}}), 
 or for that matter anywhere in Brouwer's writings prior to 1974.
\end{quote}
In our view, Myhill is wrong to believe first, that the use of family of $\C F$ requires the powerset axiom,
and, second, that its use from Bishop is surely a slip. 
The notion of family of subsets does not imply the use of the powerset as a set, since a family of subsets
is a certain assignment routine from $I$ to $\D V_0$ that behaves like a function, without being one.
Moreover, it is not a slip, as it is repeatedly 
used by Bishop in the new measure theory, also found in~\cite{BB85}, and by practicioners of Bishop-style constructive
mathematics, like Bridges and Richman. 
It is not a coincidence that the notion of family of subsets is not a fundamental function-like object in Myhill's system $\CST$. 
\end{note}

\begin{note}\label{not: propertiescompl}
\normalfont
In~\cite{Bi67}, p.~68, the following properties of complemented subsets are mentioned 
\[(\B A \cup - \B A) \cap \bigg(\B A \cup \bigcap_{i \in I} \B \lambda_0 (i)\bigg) =_{\mathsmaller{\C P^{\Disj}(X)}}
 (\B A \cup - \B A) \cap \bigg[\bigcap_{i \in I} \big(\B A \cup \B \lambda_0(i)\big)\bigg], \]
 \[ (\B A \cap - \B A) \cup \bigg(\B A \cap \bigcup_{i \in I} \B \lambda_0 (i)\bigg) =_{\mathsmaller{\C P^{\Disj}(X)}}
 (\B A \cap - \B A) \cup \bigg[\bigcup_{i \in I} \big(\B A \cap \B \lambda_0 (i)\big)\bigg]. \]
These equalities are the constructive analogue of the classical properties
\[ \B A \cup \bigcap_{i \in I} \B \lambda_0 (i) =_{\mathsmaller{\C P^{\Disj}(X)}} \bigcap_{i \in I} \big(\B A \cup \B \lambda_0 (i)\big), \]
\[ \B A \cap \bigcup_{i \in I} \B \lambda_0 (i) =_{\mathsmaller{\C P^{\Disj}(X)}} \bigcup_{i \in I} \big(\B A \cap \B \lambda_0 (i)\big). \]
\end{note}

\begin{note}\label{not: BCfamoperations}
\normalfont
In~\cite{BC72}, pp.~16--17, and in~\cite{BB85}, p.~73, 
the join and meet of a countable family of complemented subsets are defined by
\[ \bigvee_{n = 1}^{\infty}\B \lambda_0(n) := \bigg(\bigg[\bigcap_{n = 1}^{\infty}\big(\lambda_0^1(n) 
\cup \lambda_0^0(n)\big)\bigg] \cap \bigg[\bigcup_{n = 1}^{\infty}\lambda_0^1(n)\bigg], \
\ \bigcap_{n = 1}^{\infty}\lambda_0^0(n)\bigg),\]
\[ \bigwedge_{n = 1}^{\infty}\B \lambda_0(n) := \bigg(\bigcap_{n = 1}^{\infty}\lambda_0^1(n) , 
\ \ \bigg[\bigcap_{n = 1}^{\infty}\big(\lambda_0^1(n) \cup \lambda_0^0(n)\big)\bigg] 
\cap \bigg[\bigcup_{n = 1}^{\infty}\lambda_0^0(n)\bigg]\bigg).\]
These definitions can be generalised to arbitrary families of complemented subsets and properties
similar to the ones shown for $\bigcup_{i \in I}\B \lambda_0(i)$ and $\bigcap_{i \in I}\B \lambda_0(i)$ hold.
\end{note}

\begin{note}\label{not: setrelfamsub}
\normalfont
Set-relevant families of subsets over some set $I$, and set-relevant direct families of subsets over 
some directed set $(I, \lt_I)$ can be studied in a way similar to set-relevant families of sets over $I$
and set-relevant direct families of sets over $(I, \lt_I)$ in section~\ref{sec: genfamsofsets}. As a
consequence, a theory of generalised direct spectra of subspaces can be developed. Families of families
of subsets of $X$ can also be studied, in analogy to families of families of sets (see Section~\ref{sec: famoffam}).
As $\Fam(I, X)$ is in $\D V_0$, the families of families of subsets of $X$ are defined in $\D V_0$.

\end{note}

\chapter{Families of sets and spectra of Bishop spaces}
\label{chapter: bspaces}

We connect various notions and results from the theory of families of sets and subsets to the 
theory of Bishop spaces, a function-theoretic approach to constructive topology. Associating in
an appropriate way to each set $\lambda_0(i)$ of an $I$-family of sets $\Lambda$ a Bishop 
topology $F_i$ a spectrum $S(\Lambda)$ of Bishop spaces is generated. The $\sum$-set and the $\prod$-set 
of a spectrum $S(\Lambda)$ are equipped with canonical Bishop topologies. 
A direct spectrum of Bishop spaces is a family of Bishop spaces associated to a direct family of sets.
The direct and inverse limits of direct spectra of Bishop spaces are studied. Direct spectra of Bishop 
subspaces are also examined.
For all notions and facts on Bishop spaces mentioned in this chapter we refer to section~\ref{sec: bishop}
of the Appendix. Many Bishop topologies are defined inductively within the extension $\BISH^*$ of $\BISH$ 
with inductive definitions with rules of countably many premises. For all notions and facts on directed
sets mentioned in this chapter we refer to section~\ref{sec: directed} of the Appendix.

\section{Spectra of Bishop spaces}
\label{sec: spectra}

Roughly speaking, if $S$ is a structure on some set, an $S$-spectrum\index{$S$-spectrum} is 
an $I$-family of sets $\Lambda$ such that each set $\lambda_0 (i)$ is equipped with a structure $S_i$, which is
compatible with the transport maps $\lambda_{ij}$ of $\Lambda$. 
Accordingly, a spectrum of Bishop spaces is an $I$-family of sets $\Lambda$ such that each set $\lambda_0 (i)$
is equipped with a Bishop topology, which is compatible with the transport maps of $\Lambda$. As expected, 
in the case of a spectrum of Bishop spaces this compatibility condition is that the transport
maps $\lambda_{ij}$ are Bishop morphisms i.e. $\lambda_{ij} \in \Mor(\C F_i, \C F_j)$. It is natural to associate to 
$\Lambda$ an $I$-family of sets $\Phi := (\phi_0^{\Lambda}, \phi_1^{\Lambda})$ such that $\C F_i := 
\big(\lambda_0 (i), \phi_0^{\Lambda} (i)\big)$ is the Bishop space corresponding to $i \in I$.
If $i =_I j$, and if we put no restriction to the definition of $\phi_{ij}^{\Lambda} : F_i \to F_j$,
we need to add extra data in the definition of a map between spectra of
Bishop spaces. Since the map $\lambda_{ji}^* : F_i \to F_j$, where $\lambda_{ji}^*$ is the element 
of $\D F(F_i, F_j)$ induced by the Bishop morphism $\Lambda_{ji} \in \Mor(\C F_j, \C F_i)$, is 
generated by the data of $\Lambda$, it is natural to define $\phi_{ij} := \lambda_{ji}^*$. 
In this way proofs of properties of maps between spectra of Bishop spaces become easier. If $X$ 
is a set, we use the notation $\D F(X) := \D F(X, \Real)$, and every subset of $\D F(X)$ considered 
in this chapter is an extensional subset of it.

\begin{definition}\label{def: spectrum}
Let $\Lambda := (\lambda_0, \lambda_1)$, $M := (\mu_0, \mu_1) \in \Fam(I)$. 
A \textit{family of Bishop topologies associated to} $\Lambda$\index{family of Bishop topologies associated
to a family of sets} is a pair $\Phi^{\Lambda} := \big(\phi_0^{\Lambda}, \phi_1^{\Lambda}\big)$, where
$\phi_0^{\Lambda} \colon I \sto \D V_0$ and $\phi^{\Lambda} \colon \bigcurlywedge_{(i,j) \in D(I)}\D F\big(\phi_0^{\lambda}(i), \phi_0^{\Lambda}(j)\big)$, such that the following conditions hold:\\[1mm]
\normalfont (i)
\itshape  $\phi_0^{\Lambda}(i) := F_i \subseteq \D F(\lambda_0(i))$, and $\C F_i := (\lambda_0(i), F_i)$
is a Bishop space, for every $i \in I$.\\[1mm]
\normalfont (ii) 
\itshape $\lambda_{ij} \in \Mor(\C F_i, \C F_j)$, for every $(i, j) \in D(I)$.\\[1mm]
\normalfont (iii) $\phi_1^{\Lambda}(i, j) := \lambda_{ji}^*$, for every $(i, j) \in D(I)$,
where, if $f \in F_i$, the induced map $\lambda_{ji}^* \colon F_i \to F_j$ from $\lambda_{ji}$ is defined by 
 $\lambda_{ji}^* (f) := f \circ \lambda_{ji}$, for every $f \in F_i$.\\[1mm]
We call the structure $S(\Lambda) := (\lambda_0, \lambda_1, \phi_0^{\Lambda}, \phi_1^{\Lambda})$
is called a \textit{spectrum of Bishop spaces} over $I$\index{spectrum over a set}, or an 
$I$-\textit{spectrum}\index{$I$-spectrum} with Bishop spaces $(\C F_i)_{i \in I}$ and Bishop isomorphisms
$(\lambda_{ij})_{(i,j) \in D(I)}$. \index{$S(\Lambda)$}
If 
$S(M) := (\mu_0, \mu_1, \phi_0^M, \phi_1^M)$ is an $I$-spectrum with Bishop spaces $(\C G_i)_{i \in I}$ and
Bishop isomorphisms $(\mu_{ij})_{(i,j) \in D(I)}$, 
a \textit{spectrum-map} $\Psi$\index{spectrum-map} from $S(\Lambda)$ to $S(M)$\index{spectrum-map}, in
symbols $\Psi \colon S(\Lambda) \To S(M)$\index{$\Phi \colon S(\Lambda) \To S(M)$}, is a family-map
$\Psi \colon \Lambda \To M$. The totality of spectrum-maps from $S(\Lambda)$ to $S(M)$ is denoted by
$\Map_I(S(\Lambda), S(M))$ and it is equipped with the equality of $\Map_I(\Lambda, M)$.
A spectrum-map $\Phi \colon S(\Lambda) \To S(M)$ is called \textit{continuous}\index{continuous spectrum-map}, if 
$\Psi_i \in\Mor(\C F_i, \C G_i)$, for every $i \in I$, and we denote by $\Cnt_I(S(\Lambda), S(M))$ their totality, which is equipped with the equality of $\Map_I(\Lambda, M)$. The totality $\Spec(I)$ of $I$-spectra of Bishop spaces 
is equipped with the equality $S(\Lambda) =_{\mathsmaller{\Spec(I)}} S(M)$ if and only if there exist
continuous spectrum-maps $\Phi \colon S(\Lambda) \To S(M)$ and $\Psi \colon S(M) \To S(\Lambda)$ such
that $\Phi \circ \Psi =_{\mathsmaller{\Map_I(M, M)}} \Id_M$ and  $\Psi \circ 
\Phi =_{\mathsmaller{\Map_I(\Lambda, \Lambda)}} \Id_{\Lambda}$.
\end{definition}

As the identity map $\id_X \in \Mor(\C F, \C F)$, where $\C F := (X, F)$ is a Bishop space, the identity family-map 
$\id_{\Lambda} \colon \Lambda \To \Lambda$ is a continuous spectrum-map from $S(\Lambda)$ to $S(\Lambda)$.
As the composition of Bishop morphism is a Bishop morphism, if $\Phi \colon S(\Lambda) \To S(M)$ and $\Xi \colon S(M) \To S(N)$ are continuous spectrum-maps, then $\Xi \circ \Phi \colon S(\Lambda) \To S(N)$ is a continuous spectrum-map.

\begin{definition}\label{def: 2spectrum}
The structure $S(\D 2) := \big(\lambda_0^{\D 2}, \lambda_1^{\D 2}, 
\phi_0^{\Lambda^{\D 2}}, \phi_1^{\Lambda^{\D 2}}\big)$, where $\Lambda^{\D 2} 
:= (\lambda_0^{\D 2}, \lambda_1^{\D 2})$ is the $\D 2$-family of 
$X$ and $Y$, and $\Phi^{\Lambda^{\D 2}} := \big(\phi_0^{\Lambda^{\D 2}}, \phi_1^{\Lambda^{\D 2}}\big)$ is the 
$\D 2$-family of the sets $F$ and $G$, $\phi_0^{\Lambda^{\D 2}}(0) := F$ is a topology on $X$, and
$\phi_0^{\Lambda^{\D 2}}(1) := G$ is a topology on $Y$, is the $\D 2$-\textit{spectrum of} $\C F$ \textit{and}
$\C G$\index{$\D 2$-spectrum of $\C F$ and $\C G$}. 
\end{definition}

Since $\id_X \in \Mor(\C F, \C F)$, $\id_Y \in \Mor(\C G, \C G)$, $\phi_1^{\Lambda^{\D 2}}(0, 0)
:= \id_X^*$ with $\id_X^* := \id_F$, and similarly, $\phi_1^{\Lambda^{\D 2}}(1, 1)
:= \id_Y^*$ with $\id_Y^* := \id_G$, we conclude that  $S(\D 2)$ is a $\D 2$-spectrum with Bishop 
spaces $\C F, \C G$ and Bishop isomorphisms $\id_X, \id_Y$.

\begin{remark}\label{rem: spectrum1}
 Let $S(\Lambda) := (\lambda_0, \lambda_1, \phi_0^{\Lambda}, \phi_1^{\Lambda})$ be an $I$-spectrum 
 with Bishop spaces $\C F_i$ and Bishop isomorphisms $\lambda_{ij}$, 
 $S(M) := (\mu_0, \mu_1, \phi_0^M, \phi_1^M)$ an $I$-spectrum with Bishop spaces $\C G_i$ and Bishop
 isomorphisms $\mu_{ij}$,  and $\Psi \colon S(\Lambda) \To S(M)$. Then $\Phi^{\Lambda} := \big(\phi_0^{\Lambda}, 
 \phi_1^{\Lambda}\big) \in \Fam(I)$, and if $\Psi$ is continuous, then, for every $(i, j) \in D(I)$, the 
 following diagram commutes
\begin{center}
\begin{tikzpicture}

\node (E) at (0,0) {$F_i$};
\node[right=of E] (F) {$F_j$.};
\node[above=of F] (A) {$G_j$};
\node [above=of E] (D) {$\mathlarger{\mathlarger{G_i}}$};

\draw[->] (E)--(F) node [midway,below]{$\lambda_{ji}^*$};
\draw[->] (D)--(A) node [midway,above] {$\mu_{ji}^*$};
\draw[->] (D)--(E) node [midway,left] {$\Psi_i^*$};
\draw[->] (A)--(F) node [midway,right] {$\Psi_j^*$};

\end{tikzpicture}
\end{center}
\end{remark}

\begin{proof}
If $i \in I$, then $\phi_{ii}^{\Lambda}(f) := f \circ \lambda_{ii} := f \circ \id_{\lambda_0(i)} := f$.
If $i =_I j =_I k$ and $f \in F_i$, then
\begin{center}
\begin{tikzpicture}

\node (E) at (0,0) {$F_j$};
\node[right=of E] (F) {$F_k$};
\node [above=of E] (D) {$F_i$};

\draw[->] (E)--(F) node [midway,below] {$\lambda_{kj}^*$};
\draw[->] (D)--(E) node [midway,left] {$\lambda_{ji}^*$};
\draw[->] (D)--(F) node [midway,right] {$\ \lambda_{ki}^*$};

\end{tikzpicture}
\end{center}
$ \lambda_{kj}^*\big(\lambda_{ji}^*(f)\big) := \lambda_{kj}^*(f \circ \lambda_{ji})
:= (f \circ \lambda_{ji}) \circ \lambda_{kj} = f \circ (\lambda_{ji} \circ \lambda_{kj})
 = f \circ \lambda_{ki}
 := \lambda_{ki}^* (f)$.  By the definition of a continuous spectrum-map we have that if $(j, i) \in D(I)$, then 
\[ \Psi_j^*(\mu_{ji}^*(g))  := \Psi_j^*(g \circ \mu_{ji})
 := (g \circ \mu_{ji}) \circ \Psi_j 
 = g \circ (\mu_{ji} \circ \Psi_j) 
 = g \circ (\Psi_i \circ \lambda_{ji}) \]
 \[
 = (g \circ \Psi_i) \circ \lambda_{ji} 
  := \lambda_{ji}^* (g \circ \Psi_i)
 := \lambda_{ji}^* (\Psi_i^*(g)). \ \ \ \ \ \ \ \qedhere \]
\end{proof}

\section{The topology on the $\sum$- and the $\prod$-set of a spectrum}
\label{sec: sumtop}

\begin{remark}\label{rem: beforesumtop}
Let $S(\Lambda) := (\lambda_0, \lambda_1, \phi_0^{\Lambda}, \phi_1^{\Lambda}) \in \Spec(I)$ 
with Bishop spaces $(\C F_i)_{i \in I}$ and Bishop isomorphisms $(\lambda_{ij})_{(i,j) \in D(I)}$.
If $\Theta \in \prod_{i \in I}F_i$,
the following operation is a function
\[ f_{\Theta} : \bigg(\sum_{i \in I}\lambda_0 (i) \bigg) \sto \Real, \ \ \ \ 
f_{\Theta}(i, x) := \Theta_i(x); \ \ \ \ (i, x) \in \sum_{i \in I}\lambda_0 (i). \]

\end{remark}

\begin{proof}
If $(i, x) =_{\sum_{i \in I}\lambda_0 (i)} (j, y) : \TOT i =_I j \ \& \ \lambda_{ij}(x) =_{\lambda_0 (j)} y$,
by the definition of $\prod_{i \in I}F_i$ we have that $\Theta_i =_{F_i} \phi_{ji}^{\Lambda}(\Theta_j) := \lambda_{ij}^*(\Theta_j) := \Theta_j \circ \lambda_{ij}$, hence 
$f_{\Theta}(i, x) := \Theta_i(x) =_{\Real} [\Theta_j \circ \lambda_{ij}](x) =_{\Real} \Theta_j(y) := f_{\Theta}(j, y).$
\end{proof}

\begin{definition}\label{def: toponsigma}
Let $S(\Lambda) := (\lambda_0, \lambda_1, \phi_0^{\Lambda}, \phi_1^{\Lambda}) \in \Spec(I)$ 
with Bishop spaces $(\C F_i)_{i \in I}$ and Bishop isomorphisms $(\lambda_{ij})_{(i,j) \in D(I)}$. 
The \textit{sum Bishop space}\index{sum Bishop space} of $S(\Lambda)$
is the pair \index{$\sum_{i \in I}\C F_i$} \index{$\int_{i \in I}F_i$}
\[ \sum_{i \in I}\C F_i := \bigg(\sum_{i \in I}\lambda_0 (i), \int_{i \in I}F_i \bigg), \ \ \ \ \mbox{where} \ \
\int_{i \in I}F_i := \bigvee_{\Theta \in \prod_{i \in I}F_i}f_{\Theta}, \]
and the \textit{dependent product Bishop space}\index{dependent 
product Bishop space} of $S(\Lambda)$ is the pair\index{$\prod_{i \in I}\C F_i$} \index{$\oint_{i \in I}F_i$}
\[\prod_{i \in I}\C F_i := \bigg(\prod_{i \in I}\lambda_0(i), \oint_{i \in I}F_i\bigg), \ \ \ \ \mbox{where} \ \
\oint_{i \in I}F_i := \bigvee_{i \in I}^{f \in F_i}\big(f \circ \pi_i^{\Lambda}\big),\] 
and $\pi_i^{\Lambda}$ is the projection function defined in
Proposition~\ref{prp: map2}$($i$)$.
\end{definition}

\begin{proposition}\label{prp: spectrum1}
Let $S(\Lambda) := (\lambda_0, \lambda_1, \phi_0^{\Lambda}, \phi_1^{\Lambda}) \in \Spec(I)$ with Bishop spaces 
$(\C F_i)_{i \in I}$ and Bishop isomorphisms $(\lambda_{ij})_{(i,j) \in D(I)}$,  
$S(M) := (\mu_0, \mu_1, \phi_0^{M}, \phi_1^{M}) \in Spec(I)$
with Bishop spaces $(\C G_i)_{i \in I}$ and Bishop isomorphisms $(\mu_{ij})_{(i,j) \in D(I)}$, and $\Psi \colon 
S(\Lambda) \To S(M)$.\\[1mm]
\normalfont (i)
\itshape  If $i \in I$, then $e_i^{\Lambda} \in \Mor\big(\C F_i, \sum_{i \in I}\C F_i\big)$.\\[1mm]
\normalfont (ii)
\itshape If $\Psi$ is continuous, then 
$\Sigma \Psi \in \Mor\big(\sum_{i \in I}\C F_i, \sum_{i \in I}\C G_i \big)$.\\[1mm]
\normalfont (iii)
\itshape If $\Psi$ is continuous, then
$\Pi \Psi \in \Mor\big(\prod_{i \in I}\C F_i, \prod_{i \in I}\C G_i \big)$. 
\end{proposition}

\begin{proof}
(i) By the $\bigvee$-lifting of morphisms it suffices to show that 
$\forall_{\Theta \in \prod_{i \in I}F_i}\big(f_{\Theta} \circ e_i^{\Lambda} \in F_i\big)$. If 
$x \in \lambda_0(i)$, then
$\big(f_{\Theta} \circ e_i^{\Lambda}\big)(x) := f_{\Theta}(i, x) := \Theta_i(x),$ and
$f_{\Theta} \circ e_i^{\Lambda} := \Theta_i \in F_i$.\\
(ii) By the $\bigvee$-lifting of morphisms it suffices to show that 
$\forall_{\Theta{'} \in \prod_{i \in I}G_i}\big(f_{\Theta{'}} \circ \Sigma \Psi \in \int_{i \in I}F_i\big).$
If $i \in I$ and $x \in \lambda_0(i)$, we have that 
$\big(f_{\Theta{'}} \circ \Sigma \Psi\big)(i, x) := f_{\Theta{'}}\big(i, \Psi_i(x)\big) 
:= \Theta_i {'}\big(\Psi_i(x)\big) := f_{\Theta}(i, x),$
where $\Theta : \bigcurlywedge_{i \in I}F_i$ is defined by $\Theta_i := \Theta_i{'} \circ \Psi_i$, 
for every $i \in I$. By the continuity of $\Psi$ we get $\Theta_i \in F_i$. We show that $\Theta
\in \prod_{i \in I}F_i$. If $i =_I j$, by the commutativity of the diagram of Remark~\ref{rem: spectrum1} we get
$\phi_{ij}^{\Lambda}(\Theta_i) := \lambda_{ji}^* (\Theta_i) := \Theta_i \circ \lambda_{ji} := 
(\Theta_i{'} \circ \Psi_i) \circ \lambda_{ji} :=  (\Theta_i{'} \circ \mu_{ji}) \circ \Psi_j = \Theta_j{'} \circ \Psi_j 
:= \Theta_j$.\\
(iii) By the $\bigvee$-lifting of morphisms it suffices to show that 
$\forall_{i \in I}\forall_{g \in G_i}\big((g \circ \pi_i^M) \circ \prod \Psi \in \oint_{i \in I} F_i\big)$.
If $\Theta \in \prod_{i \in I}\lambda_0(i)$, then
$\big[(g \circ \pi_i^M) \circ \Pi \Psi \big](\Theta) := g(\Psi_i (\Theta_i)) := 
\big[(g \circ \Psi_i) \circ \pi_i^{\Lambda}\big](\Theta),$
hence $(g \circ \pi_i^M) \circ \Pi \Psi = (g \circ \Psi_i) \circ \pi_i^{\Lambda}$. By the 
continuity of $\Psi$ we have that $\Psi_i \in \Mor(\C F_i, \C G_i)$, hence $g \circ \Psi_i \in F_i$, and
$(g \circ \Psi_i) \circ \pi_i^{\Lambda} \in \oint_{i \in I}F_i$.
\end{proof}

If $S(\Lambda^{\D 2})$ is the spectrum of the Bishop spaces $\C F$ and $\C G$,
its sum Bishop space 
\[ \C F + \C G  :=  \bigg(\sum_{i \in \D 2}\lambda_0^{\D 2}(i),  
\int_{i \in \D 2}\phi_0^{\Lambda^{\D2}}(i)\bigg) := (X + Y, F + G) \]
is called the \textit{coproduct} of $\C F$ and $\C G$\index{coproduct of Bishop spaces}. 
By definition of the sum Bishop topology
\[ F + G := \bigvee_{f \in F, g \in G}f \oplus g , \]
\[ (f \oplus g)(w) := \left\{ \begin{array}{ll}
                 f(x)   &\mbox{, $\exists_{x \in X}\big(w := (0, x)\big)$}\\
                 g(y)             &\mbox{, $\exists_{y \in Y}\big(w := (1, y)\big)$}
                 \end{array}
          \right. ;
          \ \ \ \ f \in F, \ g \in G.
\]

The coproduct Bishop space is the coproduct in the category of Bishop spaces.

\begin{proposition}\label{prp: coproduct1}
Let $\C F := (X, F)$ and $\C G := (Y, G)$ be Bishop spaces.\\[1mm]
\normalfont (i)
\itshape  The function $i_X : X \to X + Y$, defined by $x \mapsto (0, x)$, for every $x \in X$,
is in $\Mor(\C F, \C F + \C G)$.\\[1mm]
\normalfont (ii)
\itshape  The function $i_Y : Y \to X + Y$, defined by $y \mapsto (1, y)$, for every $y \in Y$, 
is in $\Mor(\C G, \C F + \C G)$.\\[1mm]
\normalfont (iii)
\itshape If $\C H := (Z, H)$ is a Bishop space, $\phi_X \in \Mor(\C F, \C H)$ and $\phi_Y \in
\Mor(\C G, \C H)$, there is a unique $\phi \in \Mor(\C F + \C G, \C H)$ such that the following
inner diagrams commute
\begin{center}
\begin{tikzpicture}

\node (E) at (0,0) {$X$};
\node[right=of E] (A) {$X + Y$};
\node [above=of A] (D) {$Z$};
\node[right=of A] (F) {$Y.$};

\draw[dashed, ->] (A)--(D) node [midway,right] {$\phi$};
\draw[->] (E)--(D) node [midway,left] {$\phi_X \ $};
\draw [->] (E)--(A) node [midway,below] {$i_X$};
\draw[->] (F)--(A) node [midway,below] {$i_Y$};
\draw[->] (F)--(D) node [midway,right] {$ \ \phi_Y$};

\end{tikzpicture}
\end{center}

\end{proposition}

\begin{proof}
(i) By definition $i_X \in \Mor(\C F, \C F + \C G)$ if and only if 
$\forall_{f \in F}\forall_{g \in G}\big((f \oplus g) \circ i_X \in F\big)$. It is
immediate to see that $(f \oplus g) \circ i_X := f \in F$. Case (ii) is shown similarly.\\
(iii) We define $\phi : X + Y \to \Real$ by
$$\phi(w) := \left\{ \begin{array}{ll}
                 \phi_X(x)   &\mbox{, $\exists_{x \in X}\big(w := (0, x)\big)$}\\
                 \phi_Y(y)             &\mbox{, $\exists_{y \in y}\big(w := (1, y)\big)$,}
                 \end{array}
          \right.$$ 
  and since $\phi \circ i_X := \phi_X$ and $\phi \circ i_Y := \phi_Y$, the diagrams commute.
  If $h \in H$, then  
$$(h \circ \phi)(w) := \left\{ \begin{array}{ll}
                 h(\phi_X(x))   &\mbox{, $\exists_{x \in X}\big(w := (0, x)\big)$}\\
                 h(\phi_Y(y))             &\mbox{, $\exists_{y \in Y}\big(w := (1, y)\big)$,}
                 \end{array}
          \right.$$
and since $h \circ \phi_X \in F$ and $h \circ \phi_Y \in G$, we get   
$h \circ \phi := (h \circ \phi_X) \oplus (h \circ \phi_Y) \in F + G$. 
The uniqueness of $\phi$ is immediate to show.
\end{proof}

\begin{proposition}\label{prp: coproduct2}
If $F$ is a topology on $X$, $G$ is a topology on $Y$, $F_0 \subseteq \D F(X, \Real)$,
and $G_0 \subseteq \D F(Y, \Real)$ are inhabited, then
\[ \bigg(\bigvee F_0\bigg) + G = \bigvee_{f_0 \in F_0, g \in G}f_0 \oplus g := F_0 + G, \]
\[ F + \bigg(\bigvee G_0\bigg) =  \bigvee_{f \in F, g_0 \in G_0}f \oplus g_0 := F + G_0. \]
\end{proposition}

\begin{proof}
We prove only the first equality, and the proof of the second is similar. Clearly, 
$F_0 + G \subseteq \big(\bigvee F_0\big) + G.$
Since $\big(\bigvee F_0\big) + G := \bigvee_{f \in \bigvee F_0, g \in G}f \oplus g,$
and since $F_0 + G$ is a topology, for the converse inclusion 
it suffices to show inductively
that $\forall_{f \in \bigvee F_0}P(f)$, where
$P(f) : \TOT \big(\forall_{g \in G}\big(f \oplus g \in F_0 + G\big)\big).$
If $f_0 \in F_0$, then $P(f_0)$ follows immediately. If $a \in \Real$ and $g \in G$, we show 
that $\overline{a}^X \oplus g \in F_0 + G$. 
Since $\frac{g}{3} -\overline{a}^Y \in G$, by the inductive hypothesis on $f_0 \in F_0$, we
get $f_0 \oplus \big(\frac{g}{3} -\overline{a}^Y\big) \in F_0 + G$. Since
$$(*) \ \ \ \ (f_1 \oplus g_1) + (f_2 \oplus g_2) = (f_1 + f_2) \oplus (g_1 + g_2),$$
and since $\overline{a}^X \oplus \overline{a}^Y = \overline{a}^{X+Y} \in F_0 + G$, by $(*)$ we get 
$\big(f_0 + \overline{a}^X\big) \oplus \frac{g}{3} = \big(f_0 \oplus 
\big(\frac{g}{3} -\overline{a}^Y\big)\big) + 
\big(\overline{a}^X \oplus \overline{a}^Y\big) \in F_0 + G$.
Since by the inductive hypothesis $\big(f_0 \oplus - \frac{2g}{3}\big) \in
F_0 + G$, and since $-(f \oplus g) = (-f) \oplus (-g)$, we also  get $\big(- f_0 \oplus \frac{2g}{3}\big)
\in F_0 + G$, hence by $(*)$  
$$\overline{a} \oplus g = \big[(f_0 + \overline{a}^X) \oplus \frac{g}{3}\big] + \big[\big(-f_0
\oplus \frac{2g}{3}\big)\big]
\in F_0 + G.$$
Let $f_1, f_2 \in \bigvee F_0$ such that $P(f_1)$ and $P(f_2)$. If $g \in G$, by these hypotheses 
we get $f_1 \oplus \frac{g}{2} \in F_0 + G$ and $f_2 \oplus \frac{g}{2} 
\in F_0 + G$. Hence by $(*)$
$$(f_1 + f_2) \oplus g = \big(f_1 \oplus \frac{g}{2}\big) + \big(f_2 \oplus \frac{g}{2}\big) \in 
F_0 + G.$$
If $\phi \in \BR$ and $f \in \bigvee F_0$ such that $P(f)$, we show $P(\phi \circ f)$. If $g \in G$, then
$$(**) \ \ \ \ \phi \circ (f \oplus g) := (\phi \circ f) \oplus (\phi \circ g).$$
By $P(f)$ we get $f \oplus \overline{0}^Y \in F_0 + G$,
and since $\phi \circ \overline{0}^Y = \overline{\phi(0)}^Y \in F_0 + G$,
by $(**)$  
$$(\phi \circ f) \oplus \overline{\phi(0)}^Y = \big(\phi \circ f\big) \oplus \big(\phi \circ
\overline{0}^Y\big) = \phi \circ \big(f \oplus \overline{0}^Y\big) \in F_0 + G.$$
By the case of constant functions $\overline{0}^X \oplus \big(g - \overline{\phi(0)}^Y\big)
\in F_0 + G$, hence by $(*)$
$$(\phi \circ f) \oplus g =  \big[(\phi \circ f) \oplus \overline{\phi(0)}^Y\big] + \big[\overline{0}^X
\oplus \big(g - \overline{\phi(0)}^Y\big)\big] \in F_0 + G.$$
If $f \in \bigvee F_0$ such that for every $n \geq 1$ there is some $f_n \in \bigvee F_0$ such 
that $P(f_n)$ and $U\big(X; f, f_n, \frac{1}{n}\big)$, then, for every $g \in G$, we get
$U\big(X + Y ; f \oplus g, f_n \oplus g, \frac{1}{n}\big),$
and since $F_0 + G$ is a Bishop topology, by $\BS_4$ we get $f \oplus g \in F_0 + G$, hence $P(f)$. 
\end{proof}

\section{Direct spectra of Bishop spaces}
\label{sec: dirspectra}

As in the case of a family of Bishop spaces associated to an $I$-family of sets,
the family of Bishop spaces associated to an $(I, \lt)$-family of sets is defined in a minimal way from 
the data of $\Lambda^{\lt}$. According to these
data, the corresponding functions $\phi_{ij}^{\lt}$ behave necessarily in a contravariant manner i.e., 
$\phi_{ij}^{\Lambda^{\lt}} \colon F_j \to F_i$. Moreover, the transport maps $\lambda_{ij}^{\lt}$ are
Bishop morphisms, and not necessarily Bishop isomorphisms.

\begin{definition}\label{def: preorderedspectrum}
Let $(I, \lt)$ be a directed set, and let $\Lambda^{\lt} := (\lambda_0, \lambda_1^{\lt}),
M^{\lt} := (\mu_0, \mu_1^{\lt}) \in \Fam(I, \lt_I)$. 
A \textit{family of Bishop topologies associated to} $\Lambda^{\lt}$ is a pair $\Phi^{\Lambda^{\lt}} := 
\big(\phi_0^{\Lambda^{\lt}}, \phi_1^{\Lambda^{\lt}}\big)$, where $\phi_0^{\Lambda^{\lt}} : I \sto \D V_0$
and $\phi_1^{\Lambda^{\lt}} : \bigcurlywedge_{(i,j) \in \lt(I)} \D F \big(\phi_0^{\Lambda^{\lt}}(j), 
\phi_0^{\Lambda^{\lt}}(i)\big)$, 
such that the following conditions hold:\\[1mm]
\normalfont (i)
\itshape  $\phi_0^{\Lambda^{\lt}}(i) := F_i \subseteq \D F(\lambda_0(i))$,
and $\C F_i := (\lambda_0(i), F_i)$ is a Bishop space, for every $i \in I$.\\[1mm]
\normalfont (ii)
\itshape $\lambda_{ij}^{\lt} \in \Mor(\C F_i, \C F_j)$, for every $(i, j) \in \ltI$.\\[1mm]
\normalfont (iii)
\itshape $\phi_1^{\Lambda^{\lt}}(i, j) := \big(\lambda_{ij}^{\lt}\big)^*$, 
for every $(i, j) \in \ltI$, where, if $f \in F_j$, 
$\big(\lambda_{ij}^{\lt}\big)^* (f) := f \circ \lambda_{ij}^{\lt}$.\\[1mm]
The structure $S(\Lambda^{\lt}) := (\lambda_0, \lambda_1^{\lt}, \phi_0^{\Lambda^{\lt}}, \phi_1^{\Lambda^{\lt}})$
is called a \textit{direct spectrum}\index{direct spectrum}\index{$S(\Lambda^{\lt})$} over $(I, \lt)$, or an 
$(I, \lt)$-\textit{spectrum} with Bishop spaces $(\C F_i)_{i \in I}$ and Bishop morphisms
 $(\lambda_{ij}^{\lt})_{(i,j) \in \ltI}$. If 
$S(M^{\lt}) := (\mu_0, \mu_1, \phi_{0}^{M^{\lt}}, \phi_{1}^{M^{\lt}})$ is an $(I, \lt)$-spectrum
with Bishop spaces $(\C G_i)_{i \in I}$ and Bishop morphisms $(\mu_{ij}^{\lt})_{(i,j) \in \ltI}$, 
a \textit{direct spectrum-map}\index{direct spectrum-map} $\Psi$ from $S(\Lambda^{\lt})$ to $S(M^{\lt})$,
 in symbols $\Psi \colon S(\Lambda^{\lt}) \To S(M^{\lt})$,\index{$\Psi \colon S(\Lambda^{\lt}) \To S(M^{\lt})$} 
is a direct family-map $\Psi : \Lambda^{\lt} \To M^{\lt}$. The totality of direct spectrum-maps 
from $S(\Lambda^{\lt})$ to $S(M^{\lt})$ is denoted by $\Map_{(I, \lt_I)}(S(\Lambda^{\lt}), 
S(M^{\lt}))$\index{$\Map_{(I, \lt_I)}(S(\Lambda^{\lt}), S(M^{\lt}))$} and it is equipped with the 
equality of $\Map_I(\Lambda^{\lt}, M^{\lt})$.
A direct spectrum-map $\Psi : S(\Lambda^{\lt}) \To S(M^{\lt})$
is called \textit{continuous}, if\index{continuous direct spectrum-map} 
$\forall_{i \in I}\big(\Psi_i \in\Mor(\C F_i, \C G_i)\big)$, and let $\Cnt_{(I, \lt_I)}(S(\Lambda^{\lt}), S(M^{\lt}))$
be their totality, equipped with the equality of $\Map_I(\Lambda^{\lt}, M^{\lt})$. The 
totality $\Spec(I, \lt_I)$\index{$\Spec(I, \lt_I)$} of direct spectra over $(I, \lt_I)$ is equipped with an 
equality defined similarly to the equality on $\Spec(I)$.
A \textit{contravariant} direct spectrum\index{contravariant direct spectrum}\index{$S(\Lambda^{\mt})$} 
 $S(\Lambda^{\mt}) := (\lambda_0, \lambda_1^{\mt} ; \phi_0^{\Lambda^{\mt}}, \phi_1^{\Lambda^{\mt}})$
over $(I, \lt)$, a \textit{contravariant} direct spectrum-map $\Psi : S(\Lambda^{\mt}) \To S(M^{\mt})$, and 
their totalities 
 $\Map_{(I, \lt_I)}(S(\Lambda^{\mt}), S(M^{\mt}))$\index{$\Map_{(I, \lt_I)}(S(\Lambda^{\mt}), S(M^{\mt}))$}, 
 $\Spec(I, \mt_I)$\index{$\Spec(I, \mt_I)$} are defined similarly.
\end{definition}

\begin{remark}\label{rem: preorderedspectrum1}
 Let $(I, \lt)$ be a directed set,  
 $S(\Lambda^{\lt}) := (\lambda_0, \lambda_1 ; \phi_0^{\Lambda^{\lt}}, \phi_1^{\Lambda^{\lt}}) \in \Spec(I, \lt_I)$  
 with Bishop spaces $(\C F_i)_{i \in I}$ and Bishop morphisms 
 $(\lambda_{ij}^{\lt})_{(i,j) \in \ltI}$, 
 $S(M^{\lt}) := (\mu_0, \mu_1, \phi_0^{M^{\lt}}, \phi_1^{M^{\lt}}) \in \Spec(I, \lt_I)$ 
 with Bishop spaces $(\C G_i)_{i \in I}$ and Bishop morphisms $(\mu_{ij}^{\lt})_{(i,j) \in \ltI}$,
 and $\Psi \colon S(\Lambda^{\lt}) \To S(M^{\lt})$. Then $\Phi^{\Lambda^{\lt}} := 
 \big(\phi_0^{\Lambda^{\lt}}, \phi_1^{\Lambda^{\lt}}\big)$ is an $(I, \mt)$-family of sets, 
 and if $\Psi$ is continuous, then, for every $(i, j) \in \ltI$,
 the following diagram commutes
\begin{center}
\begin{tikzpicture}

\node (E) at (0,0) {$F_j$};
\node[right=of E] (K) {};
\node[right=of K] (F) {$F_i$.};
\node[above=of F] (A) {$\mathlarger{\mathlarger{G_i}}$};
\node [above=of E] (D) {$G_j$};

\draw[->] (E)--(F) node [midway,below]{$\big(\lambda_{ij}^{\lt}\big)^*$};
\draw[->] (D)--(A) node [midway,above] {$\big(\mu_{ij}^{\lt}\big)^*$};
\draw[->] (D)--(E) node [midway,left] {$\big(\Psi_j \big)^*$};
\draw[->] (A)--(F) node [midway,right] {$\big(\Psi_i \big)^*$};

\end{tikzpicture}
\end{center}
\end{remark}

\begin{proof}
Since $(\lambda_{ii}^{\lt})^*(f) := f \circ \lambda_{ii}^{\lt} := f \circ \id_{\lambda_0(i)} := f$, for every
$f \in F_i$, we get $(\lambda_{ii}^{\lt})^* := \id_{F_i}$.  If $i \lt j \lt k$ and $f \in F_k$, 
the required commutativity of the following diagram is shown:
\begin{center}
\begin{tikzpicture}

\node (E) at (0,0) {$F_j$};
\node[right=of E] (F) {$F_k$};
\node [above=of E] (D) {$F_i$};

\draw[->] (F)--(E) node [midway,below] {$\big(\lambda_{jk}^{\lt}\big)^*$};
\draw[->] (E)--(D) node [midway,left] {$\big(\lambda_{ij}^{\lt}\big)^*$};
\draw[->] (F)--(D) node [midway,right] {$\ \big(\lambda_{ik}^{\lt}\big)^*$};

\end{tikzpicture}
\end{center}
\[ \big(\lambda_{ij}^{\lt}\big)^* \big( \big(\lambda_{jk}^{\lt}\big)^* (f)\big) := 
 \big(\lambda_{ij}^{\lt}\big)^* (f \circ \lambda_{jk}^{\lt})
:= (f \circ \lambda_{jk}^{\lt}) \circ \lambda_{ij}^{\lt}
:= f \circ (\lambda_{jk}^{\lt} \circ \lambda_{ij}^{\lt})
= f \circ \lambda_{ik}^{\lt}
:= \big(\lambda_{ik}\big)^* (f). \]
To show the required commutativity, if $g \in G_j$, then
\[ \big(\lambda_{ij}^{\lt}\big)^*\big(\big(\Psi_j^{\lt}\big)^*(g)\big) := 
 \big(\Psi_j\big)^*(g) \circ \lambda_{ij}^{\lt}
 := \big(g \circ \Psi_{j}\big) \circ \lambda_{ij}^{\prec} 
 =_{F_i} g \circ \big(\Psi_{j} \circ \lambda_{ij}^{\lt}\big) 
 =_{F_i} g \circ \big(\mu_{ij}^{\lt} \circ \Psi_i\big) \]
\[=_{F_i} \big(g \circ \mu_{ij}^{\lt}\big) \circ \Psi_i 
  := \big(\Psi_i\big)^*\big(g \circ \mu_{ij}^{\lt}\big)
 := \big(\Psi_i\big)^*\big(\big(\mu_{ij}^{\lt}\big)^*(g)\big).\qedhere \]
\end{proof}

\section{The topology on the $\sum$-set of a direct spectrum}
\label{sec: dirsumtop}

\begin{remark}\label{rem: beforetopondirectedsigma}
Let $(I, \lt)$ be a directed set and
$S(\Lambda^{\lt}) := (\lambda_0, \lambda_1, \phi_0^{\Lambda^{\lt}}, \phi_1^{\Lambda^{\lt}}) \in \Spec(I, \ltI)$  
with Bishop spaces $(\C F_i)_{i \in I}$ and Bishop morphisms 
$(\lambda_{ij}^{\lt})_{(i,j) \in \ \ltI}$. If $\Theta \in \prod_{i \in I}^{\mt}F_i$, 
the following operation is a function
\[ f_{\Theta} : \bigg(\sum_{i \in I}^{\lt}\lambda_0 (i) \bigg) \sto \Real, \ \ \ \ f_{\Theta}(i, x) := 
\Theta_i(x), \ \ \ \ 
(i,x) \in \sum_{i \in I}^{\lt}\lambda_0 (i). \]
 
\end{remark}

\begin{proof}
Let $(i, x) =_{\mathsmaller{\sum_{i \in I}^{\lt}\lambda_0 (i)}} (j, y) :\TOT
\exists_{k \mt i, j}\big(\lambda_{ik}^{\lt}(x) =_{\mathsmaller{\lambda_0 (k)}} 
\lambda_{jk}^{\lt}(y)\big)$. Since 
$\Theta_i = \phi_{ki}^{\mt}(\Theta_k) := (\lambda_{ik}^{\lt})^*(\Theta_k) := 
\Theta_k \circ \lambda_{ik}^{\lt},$
and similarly $\Theta_j = \Theta_k \circ \lambda_{jk}^{\lt}$, we have that
\[ \Theta_i(x) = \big[\Theta_k \circ \lambda_{ik}^{\lt}\big](x) := \Theta_k\big(\lambda_{ik}^{\lt}(x)\big)
= \Theta_k\big(\lambda_{jk}^{\lt}(y)\big) := \big[\Theta_k \circ \lambda_{jk}^{\lt}\big](y) = 
\Theta_j(y).\qedhere  \]
\end{proof}


\begin{definition}\label{def: topondirectedsigma}
Let $(I, \lt)$ be a directed set and
$S(\Lambda^{\lt}) := (\lambda_0, \lambda_1, \phi_0^{\Lambda^{\lt}}, \phi_1^{\Lambda^{\lt}}) \in \Spec(I, \lt_I)$  
with Bishop spaces $(\C F_i)_{i \in I}$ and Bishop morphisms 
$(\lambda_{ij}^{\lt})_{(i,j) \in\ltI}$.
The Bishop space 
\[ \sum_{i \in I}^{\lt}\C F_i := \bigg(\sum_{i \in I}^{\lt}\lambda_0 (i), \int_{i \in I}^{\lt}F_i \bigg) 
\ \ \ \ \mbox{where} \ \ 
\int_{i \in I}^{\lt}F_i := \bigvee_{\mathsmaller{\Theta \in \prod_{i \in I}^{\mt}F_i}}f_{\Theta}, \]
is the \textit{sum Bishop space} of $S(\Lambda^{\lt})$\index{sum Bishop space of $S^{\lt}$}. If $S^{\mt}$ is 
a contravariant direct spectrum over $(I, \lt)$, the sum Bishop space of $S(\Lambda^{\mt})$ is defined 
dually.
\end{definition}

\begin{lemma}\label{lem: preorderdependent}
Let $S(\Lambda^{\lt}) := (\lambda_0, \lambda_1^{\lt}, \phi_0^{\Lambda^{\lt}}, \phi_1^{\Lambda^{\lt}}),   
S(M^{\lt}) := (\mu_0, \mu_1^{\lt}, \phi_0^{M^{\lt}}, \phi_1^{M^{\lt}}) \in \Spec(I, \lt_I)$,
and let $\Psi : S(\Lambda^{\lt}) \to S(M^{\lt})$ be continuous.  
If $H \in \prod_{i \in I}^{\mt}G_i$, the dependent operation
$H^* : \bigcurlywedge_{i \in I}F_i$, defined by $H^*_i := \Psi^*_i(H_i) := H_i \circ \Psi_i$, for every $i \in I$, 
is in $\prod_{i \in I}^{\mt}F_i$.
\end{lemma}

\begin{proof}
If $i \lt j$, we need to show that $H_i^* = (\lambda_{ij}^{\lt})^*(H_j^*) = H_j^* \circ \lambda_{ij}^{\lt}$. Since
$H \in \prod_{i \in I}^{\mt}G_i$, we have that $H_i = H_j \circ \mu_{ij}^{\lt}$, and by the continuity of $\Psi$
and the commutativity of the diagram
\begin{center}
\begin{tikzpicture}

\node (E) at (0,0) {$\mu_0(i)$};
\node[right=of E] (F) {$\mu_0(j)$,};
\node[above=of F] (A) {$\lambda_0(j)$};
\node [above=of E] (D) {$\lambda_0(i)$};

\draw[->] (E)--(F) node [midway,below]{$\mu_{ij}^{\lt}$};
\draw[->] (D)--(A) node [midway,above] {$\lambda_{ij}^{\lt}$};
\draw[->] (D)--(E) node [midway,left] {$\Psi_i$};
\draw[->] (A)--(F) node [midway,right] {$\Psi_j$};

\end{tikzpicture}
\end{center}
\[ H_j^* \circ \lambda_{ij}^{\lt} := \Psi_j^*(H_j) \circ \lambda_{ij}^{\lt}
:= \big(H_j \circ \Psi_j \big) \circ \lambda_{ij}^{\lt}
:= H_j(0) \circ \big(\Psi_j \circ \lambda_{ij}^{\lt}\big) \]
\[ = H_j \circ \big(\mu_{ij}^{\lt} \circ \Psi_i \big)
= \big(H_j \circ \big(\mu_{ij}^{\lt}\big) \circ \Psi_i
= H_i \circ \Psi_i
:= \Psi^*_i(H_i)
:= H^*_i.\qedhere \]
\end{proof}

\begin{proposition}\label{prp: Spectrum1}
Let $S(\Lambda^{\lt}) := (\lambda_0, \lambda_1^{\lt}, \phi_0^{\Lambda^{\lt}}, \phi_1^{\Lambda^{\lt}})$ and  
$S(M^{\lt}) := (\mu_0, \mu_1^{\lt}, \phi_0^{M^{\lt}}, \phi_1^{M^{\lt}})$ be spectra over $(I, \lt_I)$, 
and let $\Psi : S(\Lambda^{\lt}) \To S(M^{\lt})$. \\[1mm]
\normalfont (i)
\itshape If $i \in I$, then $e_i^{\Lambda^{\lt}} \in \Mor\big(\C F_i, \sum_{i \in I}^{\lt}\C F_i\big)$.\\[1mm]
\normalfont (ii)
\itshape If $\Psi$ is continuous, then 
$\Sigma^{\lt} \Psi \in \Mor\big(\sum_{i \in I}^{\lt}\C F_i, \sum_{i \in I}^{\lt}\C G_i \big)$.\\[1mm]
\end{proposition}

\begin{proof}
(i) By the $\bigvee$-lifting of morphisms it suffices to show that 
$\forall_{\Theta \in \prod_{i \in I}^{\mt}F_i}\big(f_{\Theta} \circ e_i^{\Lambda^{\lt}}
\in F_i\big)$.
If $x \in \lambda_0(i)$, then
$\big(f_{\Theta} \circ e_i^{\Lambda^{\lt}}\big)(x) := f_{\Theta}(i, x) 
:= \Theta_i(x),$ hence
$f_{\Theta} \circ e_i^{\Lambda^{\lt}} := \Theta_i \in F_i$.\\
(ii) By the $\bigvee$-lifting of morphisms it suffices to show that 
\[ \forall_{H \in \prod_{i \in I}^{\mt}G_i}\bigg(g_{H} \circ \Sigma^{\lt} \Psi \in \int_{i \in I}^{\lt}F_i\bigg). \]
If $i \in I$ and $x \in \lambda_0(i)$, and if $H^* \in \prod_{i \in I}^{\mt}F_i$, defined 
in Lemma~\ref{lem: preorderdependent}, then
$ \big(g_{H} \circ \Sigma^{\lt} \Psi\big)(i, x) := g_{H}(i, \Psi_i(x))
:= H_i (\Psi_i(x))
:= (H_i \circ \Psi_i)(x)
:= f_{H^*}(i, x),
$
and $g_{H} \circ \Sigma^{\lt} \Psi := f_{H^*} \in \int_{i \in I}^{\lt}F_i$.
\end{proof}

\section{Direct limit of a covariant spectrum of Bishop spaces}
\label{sec: directlimit}

If $X$ is a set, by Corollary~\ref{cor: corequivstr1} the family $\Eql(X) := \big(\eql_0^X, \C E^X, 
\eql_1^X\big) \in \Set(X, X)$,
where $\eql_0^X(x) := \{y \in X \mid y =_X x\}$. Consequently, if $f \colon X \to Y$ , there is unique 
$\eql_0 f \colon \eql_0 X(X) \to Y$ such that the following diagram commutes
\begin{center}
\begin{tikzpicture}

\node (E) at (0,0) {$\eql_0 X(X)$};
\node [above=of E] (D) {$X$};
\node[right=of D] (F) {$Y,$};

\draw[dashed, ->] (E)--(F) node [midway,right] {$ \ \eql_0 f$};
\draw [->] (D)--(E) node [midway,left] {$\eql^*_0$};
\draw[->] (D)--(F) node [midway,above] {$f$};

\end{tikzpicture}
\end{center}
where $\eql_0 X(X)$ is the totality $X$ with the equality $x =_{\eql_0 X (X)} x{'} 
:\TOT \eql_0^X(x) =_{\C P(X)} \eql_0^X(x{'})$. As $\Eql(X) \in \Set(X, X)$, we get $\eql_0^X(x) 
=_{\C P(X)} \eql_0^X(x{'}) \TOT x =_X x{'}$.
The map $\eql_0^* \colon X \to \eql_0 X(X)$ is defined by the identity map-rule, written in the 
form $x \mapsto \eql_0^X(x)$, for every $x \in X$. We use the set $\eql_0 X(X)$ to define the direct
limit of a direct spectrum of Bishop spaces. In what follows we avoid including the superscript $X$ in our notation.

\begin{definition}\label{def: Lim}
Let $S(\Lambda^{\lt}) := (\lambda_0, \lambda_1^{\lt}, \phi_0^{\Lambda^{\lt}}, 
\phi_1^{\Lambda^{\lt}}) \in \Spec(I, \lt_I)$ and 
$ \eql_0 \colon \sum_{i \in I}^{\lt} \lambda_0(i) \sto  \D V_0$, defined by 
\[ \eql_0 (i, x) := \bigg\{(j, y) \in \sum_{i \in I}^{\lt} \lambda_0(i) \mid (j, y)
=_{\mathsmaller{\sum_{i \in I}^{\lt} \lambda_0(i)}} (i, x)\bigg\}; \ \ \ \ (i, x) \in 
\sum_{i \in I}^{\lt} \lambda_0(i),.
\]
The \textit{direct limit} $\underset{\to} \Lim \lambda_0 (i)$\index{$\underset{\to} \Lim \lambda_0 (i)$} 
of\index{direct limit of a spectrum} $S(\Lambda^{\lt})$ is the set
\[  \underset{\to} \Lim \lambda_0 (i) := \eql_0 \sum_{i \in I}^{\lt} \lambda_0(i)
\bigg(\sum_{i \in I}^{\lt} \lambda_0(i)\bigg), \]
\[ \eql_0 (i, x) =_{\underset{\to} \Lim \lambda_0 (i)} \eql_0 (j, y) : \TOT 
 \eql_0 (i, x) =_{\C P \big(\sum_{i \in I}^{\lt} \lambda_0(i)\big)} \eql_0 (j, y)
\TOT (i, x) =_{\sum_{i \in I}^{\lt} \lambda_0(i)} (j, y). \]
We write $\eql_0^{\Lambda^{\lt}}$\index{$\eql_0^{\Lambda^{\prec}}$} when we need to express 
the dependence of $\eql_0$ from $\Lambda^{\lt}$.
\end{definition}

\begin{remark}\label{rem: directLim1}
If $S(\Lambda^{\lt}) := (\lambda_0, \lambda_1^{\lt}, \phi_0^{\Lambda^{\lt}},
\phi_1^{\Lambda^{\lt}}) \in \Spec(I, \lt_I)$ and $i \in I$, the operation
$\eql_i \colon \lambda_0 (i) \sto 
\underset{\to} \Lim \lambda_0 (i)$, defined by $\eql_i (x) := \eql_0 (i, x)$, for every $x \in \lambda_0 (i)$,
is a function.
\end{remark}

\begin{proof}
 If $x, x{'} \in \lambda_0 (i)$ such that $x =_{\lambda_0 (i)} x{'}$, then
 \begin{align*}
  \eql_i (x) =_{\underset{\to} \Lim \lambda_0 (i)} \eql_i (x{'}) & : \TOT 
  \eql_0 (i, x) =_{\underset{\to} \Lim \lambda_0 (i)} \eql_0 (i, x{'})\\
  & \TOT (i, x) =_{\sum_{i \in I}^{\lt} \lambda_0(i)} (i, x{'})\\
  & : \TOT \exists_{k \in I}\big(i \lt k \ \& \ \lambda_{ik}^{\lt}(x) =_{\lambda_0(k)}
  \lambda_{ik}^{\lt}(x{'})\big),
 \end{align*}
which holds, since $\lambda_{ik}^{\lt}$ is a function, and hence 
if $x =_{\lambda_0 (i)} x{'}$, then
$\lambda_{ik}^{\lt}(x) =_{\lambda_0(k)} \lambda_{ik}^{\lt}(x{'})$, for every $k \in I$ such that
$i \lt k$. Such a $k \in I$ always exists e.g., one can take $k := i$.
\end{proof}

\begin{definition}\label{def: Limtop}
Let $S(\Lambda^{\lt}) := (\lambda_0, \lambda_1^{\lt}, \phi_0^{\Lambda^{\lt}}, 
\phi_1^{\Lambda^{\lt}}) \in \Spec(I, \lt_I)$ with Bishop spaces $(F_i)_{i \in I}$ and Bishop morphisms 
$(\lambda_{ij}^{\lt})_{(i,j) \in \ltI}$.
The \textit{direct limit of} $S(\Lambda^{\lt})$\index{direct limit of a direct spectrum} 
is the Bishop space\index{$\underset{\to} \Lim \C F_i$}\index{$\underset{\to} \Lim F_i$}
\[ \underset{\to} \Lim \C F_i := \big(\underset{\to} \Lim \lambda_0 (i), \underset{\to} \Lim F_i\big), 
\ \ \ \ \mbox{where} \] 
\[ \underset{\to} \Lim F_i := 
\bigvee_{\mathsmaller{\Theta \in \prod_{i \in I}^{\mt}F_i}}\eql_0  f_{\Theta}, \]
\[ \eql_0 f_{\Theta} \big(\eql_0 (i, x)\big) := f_{\Theta}(i, x) := \Theta_i(x); \ \ \ \ \eql_0 (i, x)
\in \underset{\to} \Lim \lambda_0 (i) \]
\begin{center}
\begin{tikzpicture}

\node (E) at (0,0) {$\underset{\to} \Lim \lambda_0 (i)$};
\node [above=of E] (D) {$\sum_{i \in I}^{\lt} \lambda_0(i)$};
\node[right=of D] (F) {$\Real$.};

\draw[dashed, ->] (E)--(F) node [midway,right] {$ \ \eql_0 f_{\Theta}$};
\draw [->] (D)--(E) node [midway,left] {$\eql^*_0$};
\draw[->] (D)--(F) node [midway,above] {$f_{\Theta}$};

\end{tikzpicture}
\end{center}\end{definition}

\begin{remark}\label{rem: Limtop1}
 If $(I, \lt)$ is a directed set, $\C G := (Y, G)$ is a Bishop space, and 
 $S(\Lambda^{\lt, Y})$
 is the constant direct spectrum over $(I, \lt_I)$ with Bishop space $\C G$ and Bishop morphism $\id_Y$, the direct limit 
 $\underset{\to} \Lim \C G$ of $S(\Lambda^{\lt, Y})$
 is Bishop-isomorphic to $\C G$. Moreover, every Bishop 
 space is Bishop-isomorphic to the direct limit of a direct spectrum over any given directed set.
\end{remark}

\begin{proof}
 The proof is straightforward.
\end{proof}

\begin{proposition}[Universal property of the direct limit]\label{prp: universaldirect}
 If $S(\Lambda^{\lt}) := (\lambda_0, \lambda_1^{\lt}, \phi_0^{\Lambda^{\lt}}, \phi_1^{\Lambda^{\lt}}) \in \Spec(I, \lt_I)$ 
with Bishop spaces $(\C F_i)_{i \in I}$ and Bishop morphisms $(\lambda_{ij}^{\lt})_{(i,j) \in \lt(I)}$, 
its direct limit $\underset{\to} \Lim \C F_i$ satisfies the universal property of
direct limits\index{universal property of direct limits} i.e., \\[1mm]
\normalfont (i)
\itshape For every $i \in I$, we have that $\eql_i \in \Mor(\C F_i, \underset{\to} \Lim \C F_i)$.\\[1mm]
\normalfont (ii)
\itshape If $i \lt_I j$, the following left diagram commutes
\begin{center}
\begin{tikzpicture}

\node (E) at (0,0) {$\underset{\to} \Lim \lambda_0(i)$};
\node[below=of E] (F) {};
\node [right=of F] (B) {$\lambda_0(j)$};
\node [left=of F] (C) {$\lambda_0(i)$};
\node [right=of B] (D) {$\lambda_0(i)$};
\node[right=of D] (K) {};
\node[right=of K] (L) {$\lambda_0(j)$.};
\node[above=of K] (M) {$Y$};

\draw[->] (B)--(E) node [midway,right] {$ \ \eql_j$};
\draw[->] (C)--(E) node [midway,left] {$\eql_i \ $};
\draw[->] (C)--(B) node [midway,below] {$\lambda_{ij}^{\lt}$};

\draw[->] (L)--(M) node [midway,right] {$ \ \varepsilon_j$};
\draw[->] (D)--(M) node [midway,left] {$\varepsilon_i \ $};
\draw[->] (D)--(L) node [midway,below] {$\lambda_{ij}^{\lt}$};

\end{tikzpicture}
\end{center}
\normalfont (iii)
\itshape If $\C G := (Y, G)$ is a Bishop space and 
$\varepsilon_i : \lambda_0(i) \to Y \in \Mor(\C F_i, \C G)$, for every $i \in I$, such that if 
$i \lt j$, the
above right diagram commutes,
there is a unique function
$h \colon \underset{\to} \Lim \lambda_0(i) \to Y \in \Mor(\underset{\to} \Lim \C F_i, \C G)$ such that 
the following diagrams commute
\begin{center}
\begin{tikzpicture}

\node (E) at (0,0) {$Y$};
\node[below=of E] (F) {};
\node [right=of F] (B) {$\lambda_0(j)$,};
\node [left=of F] (C) {$\lambda_0(i)$};
\node[below=of F] (G) {$\underset{\to} \Lim \lambda_0(i)$.};

\draw[->] (B)--(E) node [midway,right] {$ \ \  \varepsilon_j$};
\draw[->] (C)--(E) node [midway,left] {$\varepsilon_i \ $};
\draw[->] (C)--(B) node [midway,below] {$\lambda_{ij}^{\lt} \ \ \ \ \ \  $};
\draw[->] (B)--(G) node [midway,right] {$ \ \eql_j$};
\draw[->] (C)--(G) node [midway,left] {$\eql_i \ $};
\draw[dashed,->] (G)--(E) node [midway,near end] {$\ \ \ h$};

\end{tikzpicture}
\end{center}
\end{proposition}

\begin{proof}
For the proof of (i), we use the $\bigvee$-lifting of morphisms. We have that
\[ \eql_i \in \Mor(\C F_i, \underset{\to} \Lim \C F_i) \TOT \forall_{\Theta \in 
\prod_{i \in I}^{\lt}F_i}\big(\eql_0 f_{\Theta} \circ \eql_i \in F_i\big). \]
If $x \in \lambda_0 (i)$, then 
$\big(\eql_0 f_{\Theta} \circ \eql_i\big)(x) := \eql_0 f_{\Theta}\big(\eql_0 (i, x)\big)
:= f_{\Theta} (i, x) := \Theta_i(x)$
hence $\eql_0 f_{\Theta} \circ \eql_i := \Theta_i \in F_i$. For the proof of (ii), if $x \in \lambda_0 (i)$, then 
\begin{align*}
 \eql_j (\lambda_{ij}^{\lt}(x)) =_{\underset{\to} \Lim \lambda_0 (i)} \eql_i (x) & : \TOT
 \eql_0 \big(j, \lambda_{ij}^{\lt}(x)\big) =_{\underset{\to} \Lim \lambda_0 (i)} \eql_0 (i, x)\\
 & \TOT \big(j, \lambda_{ij}^{\lt}(x)\big) =_{\sum_{i \in I}^{\lt}\lambda_0 (i)} (i, x)\\
 & : \TOT \exists_{k \in I} \big(i \lt k \ \& \ j \lt k \ \& \ \lambda_{ik}^{\lt}(x) =_{\lambda_0 (k)} 
 \lambda_{jk}^{\lt}(\lambda_{ij}^{\lt}(x))\big),
\end{align*}
which holds, since if $k \in I$ with $j \lt k$, the equality  
$ \lambda_{ik}^{\lt}(x) =_{\lambda_0 (k)} \lambda_{jk}^{\lt}(\lambda_{ij}^{\lt}(x))$ holds by
the definition of a direct family of sets, and by the definition of a directed set such a $k$ always exists.  To prove  
(iii) let the operation $h \colon \underset{\to} \Lim \lambda_0(i) \sto Y$, defined by 
$h\big(\eql_0(i, x)\big) := \varepsilon_i(x)$, for every $\omega(i, x) \in \underset{\to} \Lim \lambda_0(i)$. 
First we show that $h$ is a function. Let 
$$ \eql_0 (i, x) =_{\underset{\to} \Lim \lambda_0 (i)} \eql_0 (j, y) \TOT 
\exists_{k \in I}\big(i, j \lt k \ \& \ \lambda_{ik}^{\lt}(x) =_{\lambda_0(k)}
  \lambda_{jk}^{\lt}(y)\big).$$
By the supposed commutativity of the following diagrams
\begin{center}
\begin{tikzpicture}

\node (E) at (0,0) {$Y$};
\node[below=of E] (F) {};
\node [right=of F] (B) {$\lambda_0(k)$};
\node [left=of F] (C) {$\lambda_0(i)$};
\node [right=of B] (D) {$\lambda_0(j)$};
\node[right=of D] (K) {};
\node[right=of K] (L) {$\lambda_0(k)$};
\node[above=of K] (M) {$Y$};

\draw[->] (B)--(E) node [midway,right] {$ \ \varepsilon_k$};
\draw[->] (C)--(E) node [midway,left] {$\varepsilon_i \ $};
\draw[->] (C)--(B) node [midway,below] {$\lambda_{ik}^{\lt}$};

\draw[->] (L)--(M) node [midway,right] {$ \ \varepsilon_k$};
\draw[->] (D)--(M) node [midway,left] {$\varepsilon_j \ $};
\draw[->] (D)--(L) node [midway,below] {$\lambda_{jk}^{\lt}$};

\end{tikzpicture}
\end{center}
we get
$ h\big(\omega(i, x)\big) := \varepsilon_i(x) = \varepsilon_k\big(\lambda_{ik}^{\lt}(x)\big) =
\varepsilon_k\big(\lambda_{jk}^{\lt}(y)\big) = \varepsilon_j(y) := h\big(\omega(j, y)\big)$.
Next we show that $h$ is a Bishop morphism. By the $\bigvee$-lifting of morphisms we have that
$h \in \Mor(\underset{\to} \Lim \C F_i, \C G) \TOT \forall_{g \in G}\big(g \circ h \in 
\underset{\to} \Lim F_i\big)$.
If $g \in G$, we show that the dependent operation 
$\Theta_g : \bigcurlywedge_{i \in I}F_i $, defined by $\Theta_g (i) := g \circ \varepsilon_i$, for every $i \in I$,
is well-defined, since $\varepsilon_i \in \Mor(\C F_i, \C G)$, and $\Theta_g \in \prod_{i \in I}^{\mt}F_i$. 
To prove the latter, if $i \lt k$, we show that $\Theta_g (i) = \Theta_g (k) 
\circ \lambda_{ik}^{\lt}$. By the commutativity of the above left diagram we have that
$\Theta_g (k) \circ \lambda_{ik}^{\lt} := (g \circ \varepsilon_k) \circ \lambda_{ik}^{\lt} :=
g \circ (\varepsilon_k \circ \lambda_{ik}^{\lt}) = g \circ \varepsilon_i := \Theta_g (i)$, 
Hence $f_{\Theta_g} \in \underset{\to} \Lim F_i$. Since
$(g \circ h)\big(\eql_0(i, x)\big) := g(\varepsilon_i(x)) := (g \circ \varepsilon_i)(x) :=
[\Theta_g (i)](x) := f_{\Theta_g}\big((\eql_0(i, x)\big)$, 
we get $g \circ h := f_{\Theta_g} \in \underset{\to} \Lim F_i$. The uniqueness of $h$, and the commutativity
of the diagram in property (iii) follow immediately. 
\end{proof}

The uniqueness of $\underset{\to} \Lim \lambda_0 (i)$, up to Bishop isomorphism, is shown easily from 
its universal property.
Note that if $i, j \in I$, $x \in \lambda_0 (i)$ and $y \in \lambda_0 (j)$, we have that
\begin{align*}
 \eql_i (x) =_{\underset{\to} \Lim \lambda_0 (i)} \eql_j (y) & : \TOT 
 \eql_0 (i, x) =_{\underset{\to} \Lim \lambda_0 (i)} \eql_0 (j, y)\\
 & \TOT (i, x) =_{\sum_{i \in I}^{\lt} \lambda_0(i)} (j, y)\\
 & : \TOT \exists_{k \in I}\big(i \lt k \ \& \ j \lt k \ \& \ \lambda_{ik}^{\lt}(x) =_{\lambda_0 (k)} 
\lambda_{jk}^{\lt}(y)\big).
\end{align*}

\begin{definition}\label{def: representative}
Let $S^{\lt} := (\lambda_0, \lambda_1^{\lt} ; \phi_0^{\Lambda^{\lt}}, \phi_1^{\Lambda^{\lt}})$ be 
a direct spectrum over $(I, \lt)$. 
If $i \in I$, an element $x$ of $\lambda_0 (i)$ is a
\textit{representative} of $\omega(z) \in \underset{\to} \Lim \lambda_0 (i)$,
if $\omega_i (x) =_{\mathsmaller{\underset{\to} \Lim \lambda_0 (i)}} \omega(z)$.
\end{definition}

Although an element $\eql_0(z) \in \underset{\to} \Lim \lambda_0 (i)$ may not have a representative
in every $\lambda_0 (i)$, it surely has one at some $\lambda_0 (i)$. Actually, the following holds. 

\begin{proposition}\label{prp: representative1}
For every $n \geq 1$ and every $\eql_0(z_1), \ldots, \eql_0(z_n) \in \underset{\to}
\Lim \lambda_0 (i)$
there are $i \in I$ and $x_1, \ldots, x_n \in \lambda_0 (i)$ such that $x_l$ represents $\eql_0 (z_l)$, for
every $l \in \{1, \ldots, n\}$.
\end{proposition}

\begin{proof}
The proof is by induction on $\Nat^+$. We present only the case $n := 2$. Let $z := (j, y), z{'} := 
(j{'}, y{'}) \in \sum_{i \in I}^{\lt} \lambda_0(i)$, and $k \in I$ with $j \lt k$ and $j{'} 
\lt k$. By definition we have that $\lambda_{jk}^{\lt}(y) \in \lambda_0 (k)$ and 
$\lambda_{j{'}k}^{\lt}(y{'}) \in \lambda_0 (k)$. We show that $\lambda_{jk}^{\lt}(y)$ represents $\eql_0(z)$
and $\lambda_{j{'}k}^{\lt}(y{'})$ represents $\eql_0(z{'})$. By our remark right before 
Definition~\ref{def: representative} for the first representation 
we need to show that
\[ \omega_k \big(\lambda_{jk}^{\lt}(y)\big) =_{\underset{\to} \Lim \lambda_0 (i)} \omega_j (y) 
\TOT \exists_{k{'} \in I}\big(k \lt k{'} \ \& \ j \lt k{'} \ \& \ 
\lambda_{kk{'}}^{\lt}(\lambda_{jk}^{\lt}(y)) =_{\lambda_0 (k{'})} \lambda_{jk{'}}^{\lt}(y)\big). \]
By the composition of the transport maps it suffices to take any $k{'} \in I$ with $k \lt k{'} \ \& \ j \lt k{'}$,
and for the second representation it suffices to take any $k{''} \in I$ with
$k \lt k{''} \ \& \ j{'} \lt k{''}$.
\end{proof}

\begin{theorem}\label{thm: directspectrummap1}
Let $S(\Lambda^{\lt}) := (\lambda_0, \lambda_1^{\lt}, \phi_0^{\Lambda^{\lt}}, \phi_1^{\Lambda^{\lt}})
\in \Spec(I, \lt_I)$ with Bishop spaces $(F_i)_{i \in I}$ and Bishop morphisms $(\lambda_{ij}^{\lt})_{(i,j) \in \ltI}$, 
$S(M^{\lt}) := (\mu_0, \mu_1^{\lt}, \phi_0^{M^{\lt}}, \phi_1^{M^{\lt}}) \in \Spec(I, \lt_I)$ 
with Bishop spaces $(G_i)_{i \in I}$
and Bishop morphisms $(\mu_{ij}^{\lt})_{(i,j) \in \ltI}$, and $\Psi \colon S(\Lambda^{\lt}) \To S(M^{\lt})$. \\[1mm]
\normalfont (i)
\itshape There is a unique function $\Psi_{\to} \colon \underset{\to} \Lim \lambda_0 (i) \to 
\underset{\to} \Lim \mu_0 (i)$ such that, for every $i \in I$, the following diagram commutes
\begin{center}
\begin{tikzpicture}

\node (E) at (0,0) {$\underset{\to} \Lim \lambda_0 (i)$};
\node[right=of E] (F) {$\underset{\to} \Lim \mu_0 (i)$.};
\node[above=of F] (A) {$\mu_0(i)$};
\node [above=of E] (D) {$\lambda_0(i)$};

\draw[dashed,->] (E)--(F) node [midway,below]{$\Psi_{\to}$};
\draw[->] (D)--(A) node [midway,above] {$\Psi_i$};
\draw[->] (D)--(E) node [midway,left] {$\eql_i^{\Lambda^{\lt}}$};
\draw[->] (A)--(F) node [midway,right] {$\eql_i^{M^{\lt}}$};

\end{tikzpicture}
\end{center}
\normalfont (ii)
\itshape If $\Psi$ is continuous, then 
$\Psi_{\to} \in \Mor(\underset{\to} \Lim \C F_i, \underset{\to} \Lim \C G_i)$.\\[1mm]
\normalfont (iii)
\itshape If $\Psi_i$ is an embedding, for every $i \in I$,
then $\Psi_{\to}$ is an embedding.
\end{theorem}

\begin{proof}
(i) The following well-defined operation $\Psi_{\to} \colon \underset{\to} \Lim \lambda_0 (i) \sto 
\underset{\to} \Lim \mu_0 (i)$, given by 
\[ \Psi_{\to} \big(\eql_0^{\Lambda^{\lt}}(i, x)\big) := \eql_0^{M^{\lt}}(i, \Psi_i (x));  \ \ \ 
\ \eql_0^{\Lambda^{\lt}}(i, x) \in 
\underset{\to} \Lim \lambda_0 (i) \]
is a function,  
since, if  
$
 \eql_0^{\Lambda^{\lt}} (i, x) =_{\mathsmaller{\underset{\to} \Lim \lambda_0 (i)}} \eql_0^{\Lambda^{\lt}}(j, y)   \TOT 
 (i, x)=_{\mathsmaller{\sum_{i \in I}^{\prec} \lambda_0(i)}} (j, y)$, which is equivalent to
 $\exists_{k \in I}\big(i \lt k \ \& \ j \lt k \ \& \ \lambda_{ik}^{\lt}(x) 
 =_{\mathsmaller{\lambda_0 (k)}} \lambda_{jk}^{\lt}(y)\big)$, 
we show that
\begin{align*}
 \Psi_{\to} \big(\eql_0^{\Lambda^{\lt}}(i, x)\big) =_{\mathsmaller{\underset{\to} \Lim \mu_0 (i)}} 
 \Psi_{\to} \big(\eql_0^{\Lambda^{\lt}}(j, y)\big) & : \TOT \eql_0^{M^{\lt}}\big(i, \Psi_i^{\lt} (x)\big)
 =_{\mathsmaller{\underset{\to} \Lim \mu_0 (i)}}  \eql_0^{M^{\lt}}\big(j, \Psi_j (y)\big)\\
 & \TOT (i, \Psi_i (x)) =_{\mathsmaller{\sum_{i \in I}^{\lt}\mu_0 (i)}} (j, \Psi_j (y))\\
 & : \TOT \exists_{i \in I}\big(i, j \lt k \ \& \ \mu_{ik}^{\lt}(\Psi_i (x)) 
 =_{\mathsmaller{\mu_0 (k)}} \mu_{jk}^{\lt}(\Psi_j (y)\big).
\end{align*}
By the commutativity of the following diagrams,  and since $\Psi_k$ is a function, 
\begin{center}
\begin{tikzpicture}

\node (E) at (0,0) {$\mu_0(i)$};
\node[right=of E] (F) {$\mu_0(k)$};
\node[above=of F] (A) {$\lambda_0(k)$};
\node [above=of E] (D) {$\lambda_0(i)$};
\node [right=of F] (G) {$\mu_0 (j)$};
\node [above=of G] (H) {$\lambda_0(j)$};
\node [right=of G] (K) {$\mu_0 (k)$,};
\node [above=of K] (L) {$\lambda_0(k)$};

\draw[->] (E)--(F) node [midway,below]{$\mu_{ik}^{\lt}$};
\draw[->] (D)--(A) node [midway,above] {$\lambda_{ik}^{\lt}$};
\draw[->] (D)--(E) node [midway,left] {$\Psi_i$};
\draw[->] (A)--(F) node [midway,right] {$\Psi_k$};
\draw[->] (G)--(K) node [midway,below] {$\mu_{jk}^{\lt}$};
\draw[->] (H)--(L) node [midway,above] {$\lambda_{jk}^{\lt}$};
\draw[->] (H)--(G) node [midway,left] {$\Psi_j$};
\draw[->] (L)--(K) node [midway,right] {$\Psi_k$};

\end{tikzpicture}
\end{center}
we get
$\mu_{ik}^{\lt}\big(\Psi_i(x)\big) =_{\mathsmaller{\mu_0 (k)}} \ \Psi_k \big(\lambda_{ik}^{\lt}(x)\big)
 =_{\mathsmaller{\mu_0 (k)}} \ \Psi_k \big(\lambda_{jk}^{\lt}(y)\big)
 =_{\mathsmaller{\mu_0 (k)}} \ \mu_{jk}^{\lt}\big(\Psi_j (y)\big)
$.\\
(ii) By the $\bigvee$-lifting of morphisms it suffices to show that 
$\forall_{H \in \prod_{i \in I}^{\mt}G_i}\big((\eql_0^{M^{\lt}}g_H)
\circ \Psi_{\to} \in \underset{\to} \Lim F_i \big).$
By Definition~\ref{def: Limtop} we have that
\[ \big((\eql_0^{M^{\lt}}g_H) \circ  \Psi_{\to}\big)\big(\eql_0^{\Lambda^{\lt}}(i, x)\big) 
:= \big(\eql_0^{\Lambda^{\lt}}g_H\big)\big(\eql_0^{M^{\lt}}(i, \Psi_i (x)) \big)
:= g_H(i, \Psi_i(x)) \]
\[ := H_i (\Psi_i(x))
= (H_i \circ \Psi_i)(x)
 := H_i^* (x)
:= f_{H^*}(i, x)
:= \big(\eql_0^{\Lambda^{\lt}}f_{H^*}\big)\big(\eql_0^{\Lambda^{\lt}}(i, x)\big), \]
where $H^* \in \prod_{i \in I}^{\mt}F_i$ is defined in Lemma~\ref{lem: preorderdependent}, and
$(\eql_0^{M^{\lt}}g_{H^*}) \circ \Psi_{\to} := \eql_0^{\Lambda^{\lt}}f_{H^*} \in \underset{\to} \Lim F_i$. \\
(iii) If $\Psi_{\to} \big(\eql_0^{\Lambda^{\lt}}(i, x) =_{\mathsmaller{\underset{\to} \Lim \mu_0 (i)}} 
 \Psi_{\to} \big(\eql_0^{\Lambda^{\lt}}(j, y)$ i.e., 
$\mu_{ik}^{\lt}(\Psi_i (x)) =_{\mathsmaller{\mu_0(k)}} \mu_{jk}^{\lt}(\Psi_j(y))\big)$,
for some $k \in I$ with $i, j \lt k$, by the proof of case (ii) we get  
$\Psi_k \big(\lambda_{ik}^{\lt}(x)\big) =_{\mathsmaller{\mu_0 (k)}}  
\Psi_k \big(\lambda_{jk}^{\lt}(y)\big)$, and since
$\Psi_k$ is an embedding, we conclude that $\lambda_{ik}^{\lt}(x) =_{\mathsmaller{\lambda_0 (k)}} 
\lambda_{jk}^{\lt}(y)$ i.e., $(i, x) =_{\mathsmaller{\sum_{i \in I}^{\prec}\lambda_0(i)}} (j, y)$.
\end{proof}

\begin{proposition}\label{prp: transitivity}
Let $S(\Lambda^{\lt}) := (\lambda_0, \lambda_1^{\lt}, \phi_0^{\Lambda^{\lt}}, \phi_1^{\Lambda^{\lt}}) \in 
\Spec(I, \lt_I)$ with Bishop spaces $(F_i)_{i \in I}$ and Bishop morphisms $(\lambda_{ij}^{\lt})_{(i,j) \in \ltI}$, 
$S(M^{\lt}) := (\mu_0, \mu_1^{\lt}, \phi_0^{M^{\lt}}, \phi_1^{M^{\lt}}) \in \Spec(I, \lt_I)$ with 
Bishop spaces $(G_i)_{i \in I}$
and Bishop morphisms $(\mu_{ij}^{\lt})_{(i,j) \in \ltI}$, and 
$S(N^{\lt}) := (\nu_0, \nu_1^{\lt}, \phi_0^{N^{\lt}}, \phi_1^{N^{\lt}}) \in \Spec(I, \lt_I)$ with 
Bishop spaces $(H_i)_{i \in I}$
and Bishop morphisms $(\nu_{ij}^{\lt})_{(i,j) \in \ltI}$. If 
$\Psi \colon S(\Lambda^{\lt}) \To S(M^{\lt})$ and $\Xi \colon S(M^{\lt}) \To S(N^{\lt})$, then 
$(\Xi \circ \Psi)_{\to} :=  \Xi_{\to} \circ \Psi_{\to}$
\begin{center}
\begin{tikzpicture}

\node (E) at (0,0) {$\lambda_0(i)$};
\node[right=of E] (H) {};
\node[right=of H] (F) {$\mu_0(i)$};
\node[below=of E] (A) {$\underset{\to} \Lim \lambda_0 (i)$};
\node[below=of F] (B) {$\underset{\to} \Lim \mu_0 (i)$};
\node[right=of F] (K) {};
\node[right=of K] (G) {$\nu_0(i)$};
\node [below=of G] (C) {$\underset{\to} \Lim \nu_0 (i)$.};

\draw[->] (E)--(F) node [midway,above] {$ \Psi_i$};
\draw[->] (F)--(G) node [midway,above] {$\Xi_i$};
\draw[->] (B)--(C) node [midway,below] {$\Xi_{\to} \ \ $};
\draw[->] (A)--(B) node [midway,below] {$ \ \ \Psi_{\to}$};
\draw[->] (E)--(A) node [midway,left] {$\eql_i^{\Lambda^\lt}$};
\draw[->] (F) to node [midway,left] {$\eql_i^{M^\lt}$} (B);
\draw[->] (G)--(C) node [midway,right] {$\eql_i^{N^\lt}$};
\draw[->,bend right] (A) to node [midway,below] {$(\Xi \circ \Psi)_{\to}$} (C) ;
\draw[->,bend left] (E) to node [midway,above] {$(\Xi \circ \Psi)_i$} (G) ;

\end{tikzpicture}
\end{center}

\end{proposition}

\begin{proof}
If $\eql_0^{\Lambda^{\lt}}(i, x) \in \underset{\to} \Lim \lambda_0 (i)$, then
\begin{align*}
(\Xi \circ \Psi)_{\to}[\eql_0^{\Lambda^{\lt}}(i, x)] & := \eql_0^{N^{\lt}}(i, (\Xi \circ \Psi)_i (x))\\
& := \eql_0^{N^{\lt}}(i, (\Xi_i \circ \Psi_i) (x))\\
& := \eql_0^{N^{\lt}}(i, (\Xi_i (\Psi_i (x))))\\
& := \Xi_{\to}\big(\eql_0^{M^{\lt}}(i, \Psi_i(x))\big)\\
& := \Xi_{\to}\big(\Psi_{\to}\big(\eql_0^{\Lambda^{\lt}}(i, x)\big)\big)\\
& := (\Xi_{\to} \circ \Psi_{\to})\big(\eql_0^{\Lambda^{\lt}}(i, x)\big).\qedhere
\end{align*}
\end{proof}

\begin{definition}\label{def: directspectrummap2}
Let  $S(\Lambda^{\lt}) := (\lambda_0, \lambda_1^{\lt}, \phi_0^{\Lambda^{\lt}}, \phi_1^{\Lambda^{\lt}}) 
\in \Spec(I, \lt_I)$
and $(J, e, \cof_J) \subseteq^{\cof} I$, a cofinal subset of $I$ with modulus of cofinality $e \colon J \eto I$.
The relative spectrum of 
$S(\Lambda^{\lt})$ to $J$\index{relative direct spectrum} is the $e$-subfamily 
$S(\Lambda^{\lt}) \circ e := \big(\lambda_0 \circ e, \lambda_1 \circ e, \phi_0^{\Lambda^{\lt}} \circ e, 
\phi_1^{\Lambda^{\lt}} \circ e \big)$ of $S(\Lambda^{\lt})$, where $\Phi^{\Lambda^{\lt}} \circ e 
:= \big(\phi_0^{\Lambda^{\lt}} \circ e, \phi_1^{\Lambda^{\lt}} \circ e \big)$ is the $e$-subfamily 
of $\Phi^{\Lambda^{\lt}}$.
\end{definition}

\begin{lemma}\label{lem: cofinallemma}
Let $S(\Lambda^{\lt}) := (\lambda_0, \lambda_1^{\lt}, \phi_0^{\Lambda^{\lt}}, \phi_1^{\Lambda^{\lt}}) 
\in \Spec(I, \lt_I)$,
$(J, e, \cof_J) \subseteq^{\cof} I$, and $S(\Lambda^{\lt}) \circ e := \big(\lambda_0 \circ e,
\lambda_1 \circ e, \phi_0^{\Lambda^{\lt}} \circ e, \phi_1^{\Lambda^{\lt}} \circ e \big)$ the relative spectrum of 
$S(\Lambda^{\lt})$ to $J$.\\[1mm]
\normalfont (i)
\itshape  If $\Theta \in \prod_{i \in I}^{\mt} F_i$, then $\Theta^J \in \prod_{j \in J}^{\mt} F_j$, where
for every $j \in J$ we define $\Theta^J_j := \Theta_{e(j)}.$\\[1mm]
\normalfont (ii)
\itshape If $H^J \in \prod_{j \in J}^{\mt} F_j$, then $H \in \prod_{i \in I}^{\mt} F_i$, 
where, for every $i \in I$, let $H_i := H^J_{\cof_J(i)} \circ \lambda_{i e(\cof_J(i))}^{\lt}$
\begin{center}
\begin{tikzpicture}

\node (E) at (0,0) {$\Real$.};
\node[above=of E] (F) {$\lambda_0 (e(\cof_J(i)))$};
\node[left=of F] (B) {};
\node[left=of B] (A) {$\lambda_0 (i)$};

\draw[->] (F)--(E) node [midway,right] {$H^J_{\cof_J(i)}$};
\draw[->] (A)--(E) node [midway,left] {$H_i \ \ \ $};
\draw[->] (A)--(F) node [midway,above] {$\lambda_{i e(\cof_J(i))}^{\lt}$};

\end{tikzpicture}
\end{center}
\end{lemma}

\begin{proof}
(i) It suffices to show that if $j \lt j{'} :\TOT e(j) \lt e(j{'})$, then $\Theta^J_j = \Theta^J_{j{'}} \circ 
\lambda_{j j{'}}^{\lt}$. Since $\Theta \in \prod_{i \in I}^{\mt} F_i$ we have that
$\Theta^J_j := \Theta_{e(j)} = \Theta_{e(j{'})} \circ \lambda_{e(j) e(j{'})}^{\lt} := 
\Theta^J_{j{'}} \circ \lambda_{j j{'}}^{\lt}$.\\[1mm]
(ii) By definition $H^J_{\cof_J(i)} \in F_{\cof_J(i)} := F_{e(\cof_J(i))}$, and since $i \lt e(\cof_J(i))$, 
we get $H_i \in \Mor(\C F_i, \C R) = \C F_i$ i.e., $H : \bigcurlywedge_{i \in I}F_i$. Next we show that if
$i \lt i{'}$, then $H_i = H_{i{'}} \circ \lambda_{i i{'}}^{\lt}$. By $(\Cf_3)$ and $(\Cf_2)$ we have that
\begin{equation}\label{eq: eq1}
i \lt i{'} \lt e(\cof_J(i{'})),
\end{equation}
and $i \lt i{'} \To \cof_J(i) \lt \cof_J(i{'}) :\TOT e(\cof_J(i)) \lt e(\cof_J(i{'})),$
hence we also get 
\begin{equation}\label{eq: eq2}
i \lt e(\cof_J(i)) \lt e(\cof_J(i{'})).
\end{equation}
Since  $H^J \in \prod_{j \in J}^{\mt} F_j$, we have that
\begin{align*}
H_{i{'}} \circ \lambda_{i i{'}}^{\lt} & := \big[H^J_{\cof_J(i{'})} \circ \lambda_{i{'} e(\cof_J(i{'}))}^{\lt}\big] \circ 
\lambda_{ii{'}}^{\lt}\\
& := H^J_{\cof_J(i{'})} \circ \big[\lambda_{i{'} e(\cof_J(i{'}))}^{\lt} \circ \lambda_{ii{'}}^{\lt}\big]\\
& \stackrel{(\ref{eq: eq1})} =  H^J_{\cof_J(i{'})} \circ \lambda_{i e(\cof_J(i{'}))}^{\lt}\\
& \stackrel{(\ref{eq: eq2})} = H^J_{\cof_J(i{'})} \circ \big[\lambda_{e(\cof_J(i)) e(\cof_J(i{'}))}^{\lt} \circ
\lambda_{i e(\cof_J(i))}^{\lt} \big]\\
& := \big[H^J_{\cof_J(i{'})} \circ \lambda_{e(\cof_J(i)) e(\cof_J(i{'}))}^{\lt}\big] \circ
\lambda_{i e(\cof_J(i))}^{\lt}\\
& := \big[H^J_{\cof_J(i{'})} \circ \lambda_{\cof_J(i) \cof_J(i{'})}^{\lt}\big] \circ
\lambda_{i e(\cof_J(i))}^{\lt}\\
& := H^J_{\cof_J(i)} \circ \lambda_{i e(\cof_J(i))}^{\lt}\\
& := H_i.\qedhere
\end{align*}
\end{proof}

\begin{theorem}\label{thm: cofinal2}
Let $S(\Lambda^{\lt}) := (\lambda_0, \lambda_1^{\lt}, \phi_0^{\Lambda^{\lt}}, \phi_1^{\Lambda^{\lt}})
\in \Spec(I, \lt_I)$,
$(J, e, \cof_J) \subseteq^{\cof} I$, and $S(\Lambda^{\lt}) \circ e := \big(\lambda_0 \circ e, \lambda_1 
\circ e, \phi_0^{\Lambda^{\lt}} \circ e, \phi_1^{\Lambda^{\lt}} \circ e \big)$ the relative spectrum of 
$S(\Lambda^{\lt})$ to $J$.
Then 
\[ \underset{\to} \Lim \C F_j \simeq \underset{\to} \Lim \C F_i. \]
\end{theorem}

\begin{proof}
We define the operation $\phi : \underset{\to} \Lim \lambda_0 (j) \sto \underset{\to} 
\Lim \lambda_0 (i)$ by
$\phi \big(\eql_0^{\Lambda^{\lt} \circ e}(j, y)\big) := \eql_0^{\Lambda^{\lt}}(e(j), y)$
\begin{center}
\begin{tikzpicture}

\node (E) at (0,0) {$\underset{\to} \Lim \lambda_0 (j)$};
\node[right=of E] (A) {};
\node[right=of A] (F) {$\underset{\to} \Lim \lambda_0 (i)$,};
\node [above=of A] (D) {$\lambda_0 (j)$};

\draw[->] (E)--(F) node [midway,below] {$\phi$};
\draw[->] (D)--(E) node [midway,left,near start] {$\eql_j^{\Lambda^{\lt} \circ e} \ $};
\draw[->] (D)--(F) node [midway,right,near start] {$ \ \eql_{e(j)}^{\Lambda^{\lt}}$};

\end{tikzpicture}
\end{center}
for every $\eql_0^{\Lambda^{\lt} \circ e}(j, y) \in \underset{\to} \Lim \lambda_0 (j)$, 
where, if $j \in J$ and $y \in \lambda_0 (j)$, we have that
\[ \eql_0^{\Lambda^{\lt} \circ e}(j, y) := \bigg\{(j{'}, y{'}) \in \sum_{j \in J}^{\lt}\lambda_0 (j) 
\mid (j{'}, y{'})
=_{\mathsmaller{\sum_{j \in J}^{\lt}\lambda_0 (j)}} (j, y) \bigg\}, \]
\[ \eql_0^{\Lambda^{\lt}}(e(j), y) := \bigg\{(i, x) \in \sum_{i \in I}^{\lt}\lambda_0 (i) \mid (i, x)
=_{\mathsmaller{\sum_{i \in I}^{\lt}\lambda_0 (i)}} (e(j), y) \bigg\}. \]
First we show that $\phi$ is a function.
 By definition we have that
\begin{align*}
\eql_0^{\Lambda^{\lt} \circ e}(j, y) =_{\mathsmaller{\underset{\to} \Lim \lambda_0 (j)}} 
\eql_0^{\Lambda^{\lt} \circ e}(j{'}, y{'}) & \TOT 
(j, y) =_{\mathsmaller{\sum_{j \in J}^{\lt}\lambda_0 (j)}} (j{'}, y{'}) \\
& \TOT \exists_{j{''} \in J}\big(j, j{'} \lt j{''} \ \& \ \lambda_{jj{''}}^{\lt} (y) 
=_{\mathsmaller{\lambda_0 (j{''})}} \lambda_{j{'}j{''}}^{\lt} (y{'})\big)  \ \ \ (1)
\end{align*}
\begin{align*}
\eql_0^{\Lambda^{\lt}}(e(j), y) =_{\mathsmaller{\underset{\to} \Lim \lambda_0 (i)}} 
\eql_0^{\Lambda^{\lt}}(e(j{'}), y{'}) & \TOT (e(j), y) 
=_{\mathsmaller{\sum_{i \in I}^{\lt}\lambda_0 (i)}} (e(j{'}), y{'}) \\
& \TOT \exists_{k \in I}\big(e(j), e(j{'}) \lt k \ \& \ \lambda_{e(j)k}^{\lt} (y) =_{\mathsmaller{\lambda_0 (k)}}
\lambda_{e(j{'})k}^{\lt} (y{'})\big). \  (2)
\end{align*}
If $k := e(j{''})$, then (1) implies (2), and hence $\phi$ is a function. To show that $\phi$ is an embedding,
we show that (2) implies (1). Since $e(j), e(j{'}) \lt k \lt e(\cof_J(k))$, 
we get $j, j{'} \lt \cof_J(k) := j{''}$.
By the commutativity of the following diagrams 
\begin{center}
\begin{tikzpicture}

\node (E) at (0,0) {$\lambda_0(e(j))$};
\node[below=of E] (F) {$\lambda_0(k)$};
\node[right=of F] (H) {};
\node[right=of H] (G) {$\lambda_0(e(\cof_J(k)))$};
\node[right=of G] (A) {$\lambda_0(k)$};
\node[above=of A] (B) {$\lambda_0(e(j{'}))$};
\node[right=of A] (D) {};
\node[right=of D] (C) {$\lambda_0(e(\cof_J(k)))$};

\draw[->] (E)--(F) node [midway,left] {$\lambda_{e(j)k}^{\lt}$};
\draw[->] (E)--(G) node [midway,right, near start] {$\ \ \ \lambda_{e(j)e(\cof_J(k))}^{\lt}$};
\draw[->] (F)--(G) node [midway,below] {$\lambda_{ke(\cof_J(k))}^{\lt}$};
\draw[->] (B)--(A) node [midway,left] {$\lambda_{e(j{'})k}^{\lt}$};
\draw[->] (A)--(C) node [midway,below] {$\lambda_{ke(\cof_J(k))}^{\lt}$};
\draw[->] (B)--(C) node [midway,right, near start] {$\ \ \ \lambda_{e(j{'})e(\cof_J(k))}^{\lt}$};

\end{tikzpicture}
\end{center}
\begin{align*}
\lambda_{jj{''}}^{\lt} (y) & := \lambda_{e(j)e(\cof_J(k))}^{\lt}(y)\\
& = \big[\lambda_{ke(\cof_J(k))} \circ \lambda_{e(j)k}^{\lt}\big] (y) \\  
& = \lambda_{ke(\cof_J(k))} \big(\lambda_{e(j)k}^{\lt}(y)\big)\\ 
& = \lambda_{ke(\cof_J(k))} \big(\lambda_{e(j{'})k}^{\lt}(y{'})\big)\\
& := \big[\lambda_{ke(\cof_J(k))} \circ \lambda_{e(j{'})k}^{\lt}\big](y{'})\\
& = \lambda_{e(j{'})e(\cof_J(k))}^{\lt}(y{'})\\
& := \lambda_{j{'} j{''}}(y{'}).
\end{align*}
By the $\bigvee$-lifting of morphisms we have that 
\[ \phi \in \Mor(\underset{\to} \Lim \C F_j, \underset{\to} \Lim \C F_i) : \TOT
\forall_{\mathsmaller{\Theta \in \prod_{i \in I}^{\mt}F_i}}\big(\eql_0 f_{\Theta}
\circ \phi \in \underset{\to} \Lim F_j\big). \]
If $\Theta \in \prod_{i \in I}^{\mt}F_i$, we have that 
\[  (\eql_0 f_{\Theta} \circ \phi) \big(\eql_0^{\Lambda^{\lt} \circ e}(j, y)\big) := 
 (\eql_0 f_{\Theta})\big(\eql_0^{\Lambda^{\lt}}(e(j), y)\big) \]
\[ := \Theta_{e(j)}(y) := \Theta^J_j(y) := (\eql_0 f_{\Theta^{J}})\big(\eql_0^{\Lambda^{\lt} \circ e}(j, y)\big), \]
where $\Theta^{J} \in \prod_{j \in J}^{\lt} F_j$ is defined in Lemma~\ref{lem: cofinallemma}(i). 
Hence, $\eql_0 f_{\Theta} \circ \phi = \eql_0 f_{\Theta^{J}} \in \underset{\to} \Lim F_j$.
Next we show that $\phi$ is a surjection. If $\eql_0^{\Lambda^{\lt}}(i, x) \in 
\underset{\to} \Lim \lambda_0 (i)$, we find $\eql_0^{\Lambda^{\lt} \circ e}(j, y) \in 
\underset{\to} \Lim \lambda_0 (j)$ such that 
$\phi \big(\eql_0^{\Lambda^{\lt} \circ e}(j, y)\big) := \eql_0^{\Lambda^{\lt}}(e(j), y) 
=_{\mathsmaller{\underset{\to} \Lim \lambda_0 (i)}} \eql_0^{\Lambda^{\lt}}(i, x)$ i.e., 
we find $k \in I$ such that
$i, e(j) \lt k$ and $\lambda_{ik}^{\lt}(x) =_{\mathsmaller{\lambda_0 (k)}} \lambda_{e(j)k}^{\lt}(y)$.
If $j := \cof_J(i)$, by $(\Cf_3)$ we have that $i \lt e(\cof_J(i))$,
and by the reflexivity of $\lt$ we have that $e(\cof_J(i)) \lt e(\cof_J(i)) := k$. If 
$y := \lambda_{i e(\cof_J(i))}{^\lt}(x) \in \lambda_0(e(\cof_J(i))) := (\lambda_0 \circ e)(\cof_J(i))$, then
\[ \lambda_{e(\cof_J(i))e(\cof_J(i))}^{\lt}\big(\lambda_{ie(\cof_J(i))}^{\lt}(x)\big) 
=_{\mathsmaller{\lambda_0 (k)}} \lambda_{ie(\cof_J(i))}^{\lt}(x). \]
We can use the $\bigvee$-lifting of openness to show that $\phi$ is an open morphism, and hence a 
Bishop isomorphism, but it is better to define directly its inverse Bishop morphism using the previous
proof of the surjectivity
 of $\phi$. Let the operation $\theta \colon \underset{\to} \Lim \lambda_0(i) \sto \underset{\to} \Lim \lambda_0(j)$,
 defined by 
\[ \theta \big(\eql_0^{\Lambda^{\lt}}(i, x)\big) := \eql_0^{\Lambda^{\lt} \circ e}\big(\cof_J(i), 
\lambda_{ie(\cof_J(i))}(x)\big);  \ \ \ \
\eql_0^{\Lambda^{\lt}}(i, x) \in \underset{\to} \Lim \lambda_0(i). \]
First we show that $\theta$ is a function. We have that 
\[ \eql_0^{\Lambda^{\lt}}(i, x)  =_{\mathsmaller{\underset{\to} \Lim \lambda_0(i)}} 
\eql_0^{\Lambda^{\lt}}(i{'}, x{'}) \TOT
\exists_{k \in I}\big(i \lt k \ \& \ i{'} \lt k \ \& \ \lambda_{ik}^{\lt}(x) =_{\mathsmaller{\lambda_0(k)}}
\lambda_{i{'}k}^{\lt}(x{'})\big), \]
\[ \eql_0^{\Lambda^{\lt} \circ e}\big(\cof_J(i), \lambda_{ie(\cof_J(i))}(x)\big) =_{\mathsmaller{\underset{\to}
\Lim \lambda_0(j)}} \eql_0^{\Lambda^{\lt} \circ e}\big(\cof_J(i{'}), \lambda_{i{'}e(\cof_J(i{'}))}(x{'})\big) \TOT \]
\[ \exists_{j{'} \in J}\bigg(\cof_J(i) \lt j{'} \ \& \ \cof_J(i{'}) \lt j{'} \ \&  \] 
\[ \lambda^{\lt}_{e(\cof_J(i))e(j{'})}\big(\lambda^{\lt}_{ie(\cof_J(i))}(x)\big) =_{\mathsmaller{\lambda_0(e(j{'})}}
\lambda^{\lt}_{e(\cof_J(i{'}))e(j{'})}\big(\lambda^{\lt}_{i{'}e(\cof_J(i{'}))}(x{'})\bigg). \]
If $j{'} := \cof_J(k)$, then by $(\Cf_2)$ we get $\cof_J(i) \lt j{'}$ and $\cof_J(i{'}) \lt j{'}$. 
Next we show that
\[ \lambda^{\lt}_{e(\cof_J(i))e(\cof_J(k))}\big(\lambda^{\lt}_{ie(\cof_J(i))}(x)\big) 
=_{\mathsmaller{\lambda_0(e(\cof_J(k))}}
\lambda^{\lt}_{e(\cof_J(i{'}))e(\cof_J(k))}\big(\lambda^{\lt}_{i{'}e(\cof_J(i{'}))}(x{'}). \]
By the following order relations, the two terms of the required equality are written as  
\begin{center}
\begin{tikzpicture}

\node (E) at (0,0) {$e(\cof_J(k))$};
\node[above=of E] (F) {$k$};
\node[right=of F] (A) {$e(\cof_J(i{'}))$};
\node[left=of F] (B) {$e(\cof_J(i))$};
\node[above=of B] (C) {$i$};
\node[above=of A] (D) {$i{'}$};
\

\draw[->] (F)--(E) node [midway,below] {};
\draw[->] (B)--(E) node [midway,left] {};
\draw[->] (A)--(E) node [midway,below] {};
\draw[->] (C)--(B) node [midway,right] {};
\draw[->] (D)--(A) node [midway,right] {};
\draw[->] (D)--(F) node [midway,right] {};
\draw[->] (C)--(F) node [midway,right] {};

\end{tikzpicture}
\end{center}
$\lambda^{\lt}_{ie(\cof_J(k))}(x) = \lambda^{\lt}_{ke(\cof_J(k))}\big(\lambda^{\lt}_{ik}(x)\big)$, 
and $\lambda^{\lt}_{i{'}e(\cof_J(k))}(x{'}) = \lambda^{\lt}_{ke(\cof_J(k))}\big(\lambda^{\lt}_{i{'}k}(x{'})\big)$.
By the equality 
$\lambda_{ik}^{\lt}(x) =_{\mathsmaller{\lambda_0(k)}} \lambda_{i{'}k}^{\lt}(x{'})$ we get the required equality.
Next we show that
\[ \theta \in \Mor\big(\underset{\to} \Lim \C F_i, \underset{\to} \Lim \C F_j\big) \TOT \forall_{H^J 
\in \prod_{j \in J}^{\mt}F_j}\bigg(\eql_0 f_{H^J} \circ \theta \in \bigvee_{\Theta \in 
\prod_{i \in I}^{\mt}F_i}\eql_0 f_{\Theta}\bigg). \]
If we fix $H^J \in \prod_{j \in J}^{\mt} F_j$, and if $H \in \prod_{i \in I}^{\mt}F_i$, defined 
in Lemma~\ref{lem: cofinallemma}(ii), then 
\begin{align*}
\big(\eql_0 f_{H^J} \circ \theta \big)\big(\eql_0^{\Lambda^{\lt}}(i,x)\big) & 
:= \eql_0 f_{H^J}\bigg(\eql_0^{\Lambda^{\lt} \circ e}\big(\cof_J(i), \lambda_{ie(\cof_J(i))}(x)\big)\bigg)\\
& := f_{H^J}\big(\cof_J(i), \lambda_{ie(\cof_J(i))}(x)\big)\\
& := H^J_{\cof_J(i)}\big(\lambda_{ie(\cof_J(i))}(x)\big)\\
& := \big[H^J_{\cof_J(i)} \circ \lambda_{ie(\cof_J(i))}\big](x)\\
& := H_i(x)\\
& := f_H(i,x)\\
& := \eql_0 f_H \big(\eql_0^{\Lambda^{\lt}}(i,x)\big),
\end{align*}
hence $\eql_0 f_{H^J} \circ \theta := \eql_0 f_H \in \underset{\to} \Lim F_i$. Next we show that
$\phi$ and $\theta$ are inverse to each other.
\begin{align*}
\phi\big(\theta\big(\eql_0^{\Lambda^{\lt}}(i,x)\big)\big) & := \phi\big(\eql_0^{\Lambda^{\lt} 
\circ e}\big(\cof_J(i), \lambda_{ie(\cof_J(i))}(x)\big) \\
& := \eql_0^{\Lambda^{\lt}}\big(e(\cof_J(i)), \lambda_{ie(\cof_J(i))}(x)\big),
\end{align*}
which is equal to $\eql_0^{\Lambda^{\lt}}(i,x)$ if and only if there is $k \in I$ with $i \lt k$ and 
$e(\cof_J(i)) \lt k$ and 
\[ \lambda_{ik}^{\lt}(x) =_{\mathsmaller{\lambda_0(k)}} 
\lambda^{\lt}_{e(\cof_J(i))k}\big(\lambda_{ie(\cof_J(i))}(x)\big), \]
which holds for every such $k \in I$. As by $(\Cf_3)$ we have that $i \lt e(\cof_J(i))$, the existence 
of such a $k \in I$ follows trivially. Similarly,
\begin{align*}
\theta\big(\phi\big(\eql_0^{\Lambda^{\lt} \circ e}(j,y)\big)\big) & := 
\theta\big(\eql_0^{\Lambda^{\lt}}\big(e(j), y)\big)\\
& := \eql_0^{\Lambda^{\lt} \circ e}\big(\cof_J(e(j)), \lambda_{e(j)e(\cof_J(e(j)))}(y)\big),
\end{align*}
which is equal to $\eql_0^{\Lambda^{\lt} \circ e}(j,y)$ if and only if there is $j{'} \in J$ 
with $j \lt j{'}$, $(\cof_J(e(j))) \lt j{'}$ and 
\[ \lambda_{e(j)e(j{'})}^{\lt}(y) =_{\mathsmaller{\lambda_0(e(j{'}))}} 
\lambda^{\lt}_{e(\cof_J(e(j)))e(j{'})}\big(\lambda_{e(j)e(\cof_J(e(j)))}(y)\big), \]
which holds for every such $j{'} \in J$. As by $(\Cf_1)$ we have that $j =_J \cof_J(e(j))$, the existence of such a 
$j{'} \in J$ follows trivially. 
\end{proof}

For simplicity we use next the same symbol for different orderings.

\begin{proposition}\label{prp: prodspec}
If $(I, \lt), (J, \lt)$ are directed sets, $i \in I$ and $j \in J$, let
\[ (i, j) \lt (i{'}, j{'}) : \TOT i \lt i{'} \ \& \ j \lt j{'}. \]
If $(K, i_K, \cof_K) \subseteq^{\cof} I$ and $(L, i_L, \cof_L) \subseteq^{\cof} J$, let 
$i_{K \times L} : K \times L \hookrightarrow I \times J$ and $\cof_{K \times L} : I \times J \to K \times L$, 
defined, for every $k \in K$ and $l \in L$, by
\[ i_{K \times L}(k, l) := \big(i_K(k), i_L(l)\big) \ \ \ \& \ \ \ \cof_{K \times L}(i, j) := \big(\cof_K(i),
\cof_L(j)\big). \]
Let $\Lambda^{\lt} := (\lambda_0, \lambda_1^{\lt}) \in \Fam(I, \lt)$  
and $M^{\lt} := (\mu_0, \mu_1^{\lt}) \in \Fam(J, \lt)$ an $(J, \lt)$. Let also 
$S(\Lambda^{\lt}) := (\lambda_0, \lambda_1^{\lt}, \phi_0^{\Lambda^{\lt}},  \phi_1^{\Lambda^{\lt}}) \in \Spec(I, \lt)$ 
with Bishop spaces $(\C F_i)_{i \in I}$ and Bishop morphisms $(\lambda_{ii{'}})_{(i, i{'})
\in D^{\lt}(I \times J)}$,
and $S(M^{\lt}) := (\mu_0, \mu_1^{\lt}, \phi_0^{M^{\lt}},  \phi_1^{M^{\lt}}) \in \Spec(J, \lt)$ 
with Bishop spaces $(\C G_j)_{j \in J}$ and Bishop morphisms $(\mu_{jj{'}})_{(j, j{'}) \in \lt(J)}^{\lt}$.\\[1mm]
\normalfont (i)
\itshape $(I \times J, \lt)$ is a directed set, and $(K \times L, i_{K \times L}, \cof_{K \times L}) 
\subseteq^{\cof} I \times J$.\\[1mm]
\normalfont (ii)
\itshape The pair $\Lambda^{\lt} \times M^{\lt} := (\lambda_0 \times \mu_0 , (\lambda_1 \times \mu_1)^{\lt}) \in 
\Fam(I \times J, \lt)$, where  
\[ (\lambda_0 \times \mu_0)\big((i, j)\big) := \lambda_0 (i) \times \mu_0(j), \]
\[ (\lambda_1 \times \mu_1)^{\lt}\big((i, j), (i{'}, j{'})\big) := 
(\lambda_1 \times \mu_1)^{\lt}_{(i, j), (i{'}, j{'})}, \]
\[ (\lambda_1 \times \mu_1)^{\lt}_{(i, j), (i{'}, j{'})}\big((x, y)\big) := 
\big(\lambda_{ii{'}}^{\lt}(x), \mu_{jj{'}}^{\lt}(y)\big). \]
\normalfont (iii)
\itshape The structure 
$S(\Lambda^{\lt} \times M^{\lt}) := \big(\lambda_0 \times \mu_0, \lambda_1^{\lt} \times \mu_1^{\lt} ; 
\phi_0^{\Lambda^{\lt} \times M^{\lt}},  \phi_1^{\Lambda^{\lt} \times M^{\lt}}\big) \in \Spec(I \times J, \lt)$
with Bishop spaces $(F_i \times G_j)_{(i, j) \in I \times J}$ and Bishop morphisms 
$(\lambda_1 \times \mu_1)_{\big((i,j)(i{'},j{'})\big) \in D^{\lt}(I \times J)}$, where
\[ \phi_0^{\Lambda^{\lt} \times M^{\lt}}(i, j) := F_i \times G_j, \]
\[ \phi_1^{\Lambda^{\lt} \times M^{\lt}}\big((i, j), (i{'}, j{'})\big) := 
[(\lambda_1 \times \mu_1)_{(i,j)(i{'},j{'})}^{\lt}]^* : F_{i{'}} \times G_{j{'}} \to F_i \times G_j.\]
\end{proposition}

\begin{proof}
(i) is immediate to show. For the proof of case (ii) we have that 
$(\lambda_1 \times \mu_1)^{\lt}_{(i, j), (i, j)}\big((x, y)\big) := 
\big(\lambda_{ii}^{\lt}(x), \mu_{jj}^{\lt}(y)\big) := (x, y),$
and if $(i, j) \lt (i{'}, j{'}) \lt (i{''}, j{''})$, then the commutativity of the 
\begin{center}
\begin{tikzpicture}

\node (E) at (0,0) {$\lambda_0(i) \times \mu_0(j)$};
\node[below=of E] (F) {$\lambda_0(i{'}) \times \mu_0(j{'})$};
\node[right=of F] (H) {};
\node[right=of H] (G) {$\lambda_0(i{''}) \times \mu_0(j{''})$};

\draw[->] (E)--(F) node [midway,left] {$\mathsmaller{(\lambda_1 \times \mu_1)^{\lt}_{(i, j), (i{'}, j{'})}}$};
\draw[->] (E)--(G) node [midway,right, near start] {$\ \ \ \mathsmaller{(\lambda_1 \times \mu_1)^{\lt}_{(i, j), (i{''}, j{''})}}$};
\draw[->] (F)--(G) node [midway,below] {$\mathsmaller{(\lambda_1 \times \mu_1)^{\lt}_{(i{'}, j{'}), (i{''}, j{''})}}$};

\end{tikzpicture}
\end{center}
above diagram follows from the equalities $\lambda_{ii{''}}^{\lt} = \lambda_{i{'}i{''}}^{\lt} \circ \lambda_{ii{'}}^{\lt}$
and $\mu_{jj{''}}^{\lt} = \mu_{j{'}j{''}}^{\lt} \circ \mu_{jj{'}}^{\lt}$.\\
(iii) We show that $(\lambda_1 \times \mu_1)_{(i,j)(i{'},j{'})}^{\lt} \in \Mor(\C F_i \times \C G_j, 
\C F_{i{'}} \times \C G_{j{'}})$. By the $\bigvee$-lifting of morphisms it suffices to show that
$\forall_{f \in F_{i{'}}}\big((f \circ \pi_1) \circ  (\lambda_1 \times
\mu_1)_{(i,j)(i{'},j{'})}^{\lt} \in F_i \times G_j\big)$ and $\forall_{g \in G_{j{'}}}\big((g \circ \pi_2) \circ 
(\lambda_1 \times \mu_1)_{(i,j)(i{'},j{'})}^{\lt} \in F_i \times G_j \big).$
If $f \in F_{i{'}}$, then $(f \circ \pi_1) \circ  (\lambda_1 \times \mu_1)_{(i,j)(i{'},j{'})}^{\lt} := (f \circ 
\lambda_{ii{'}}^{\lt}) \circ \pi_1 \in F_i \times G_j$, as $f \circ \lambda_{ii{'}}^{\lt} \in F_i$ and
$[(f \circ \pi_1) \circ  (\lambda_1 \times \mu_1)_{(i,j)(i{'},j{'})}^{\lt}](x, y) := 
(f \circ \pi_1)\big(\lambda_{ii{'}}^{\lt}(x), \mu_{jj{'}}^{\lt}(y)\big) := f \big(\lambda_{ii{'}}^{\lt}(x)\big)
:= [(f \circ \lambda_{ii{'}}^{\lt}) \circ \pi_1](x, y)$.
If $g \in G_{j{'}}$, we get $(g\circ \pi_2) \circ  (\lambda_1 \times \mu_1)_{(i,j)(i{'},j{'})}^{\lt} := 
(g \circ \lambda_{jj{'}}^{\lt}) \circ \pi_2 \in F_i \times G_j$.
\end{proof}

\begin{lemma}\label{lem: prodlemma}
Let $S(\Lambda^{\lt}) := (\lambda_0, \lambda_1^{\lt},
\phi_0^{\Lambda^{\lt}},  \phi_1^{\Lambda^{\lt}}) \in \Spec(I, \lt)$ 
with Bishop spaces $(\C F_i)_{i \in I}$ and Bishop morphisms
$(\lambda_{ii{'}})_{(i, i{'})\in \ltI}^{\lt}$,
$S(M^{\lt}) := (\mu_0, \mu_1^{\lt}, \phi_0^{M^{\lt}},  \phi_1^{M^{\lt}}) \in \Spec(J, \lt)$ 
with Bishop spaces $(\C G_j)_{j \in J}$ and Bishop morphisms $(\mu_{jj{'}})_{(j, j{'}) \in D^{\lt}(J)}$,
$\Theta \in \prod_{i \in I}^{\mt} F_i$ and $\Phi \in \prod_{j \in J}^{\mt} G_j$,. Then 
\[ \Theta_1 \in \prod_{(i, j) \in I \times J}^{\mt} F_i \times G_j \ \ \ \& \ \ \ \ \Phi_2 \in
\prod_{(i, j) \in I \times J}^{\mt} F_i \times G_j, \]
\[ \Theta_1 (i, j) := \Theta_i \circ \pi_1 \in F_i \times G_j \ \ \ \& \ \ \
\Phi_2 (i, j) := \Phi_j \circ \pi_2 \in F_i \times G_j; \ \ \ \ (i, j) \in I \times J. \]
\end{lemma}

\begin{proof}
We prove that $\Theta_1 \in \prod_{(i, j) \in I \times J}^{\mt} F_i \times G_j$, and for 
$\Phi_2$ we proceed similarly. If $(i, j) \lt (i{'}, j{'})$, we need to show that 
$\Theta_1 (i, j) = \Theta_1 (i{'}, j{'}) \circ (\lambda_1 \times \mu_1)^{\lt}_{(i, j), (i{'}, j{'})}$.
Since $\Theta \in \prod_{i \in I}^{\mt} F_i$, we have that 
$\Theta_i = \Theta_{i{'}} \circ \lambda_{i i{'}}^{\lt}$.
If $x \in \lambda_0(i)$ and $y \in \mu_0(j)$, we have that
\begin{align*}
\big[\Theta_1 (i{'}, j{'}) \circ (\lambda_1 \times \mu_1)^{\lt}_{(i, j), (i{'}, j{'})}](x, y) & := 
\big[\Theta_{i{'}} \circ \pi_1\big]\big(\lambda_{ii{'}}^{\lt}(x), \mu_{jj{'}}^{\lt}(y)\big)\\ 
& := \Theta_{i{'}}\big(\lambda_{ii{'}}^{\lt}(x)\big)\\
& := \big[\big(\Theta_{i{'}} \circ \lambda_{ii{'}}^{\lt}\big) \circ \pi_1\big](x, y)\\
& := \big(\Theta_i \circ \pi_1\big)(x, y)\\
& := \big[\Theta_1(i, j)\big](x, y).\qedhere
\end{align*}
\end{proof}

\begin{proposition}\label{prp: proddirect}
If $S(\Lambda^{\lt}) := (\lambda_0, \lambda_1^{\lt}, \phi_0^{\Lambda^{\lt}},  \phi_1^{\Lambda^{\lt}}) \in \Spec(I, \lt)$ 
with Bishop spaces $(\C F_i)_{i \in I}$ and Bishop morphisms $(\lambda_{ii{'}})_{(i, i{'})\in \ltI}^{\lt}$,
and $S(M^{\lt}) := (\mu_0, \mu_1^{\lt}, \phi_0^{M^{\lt}},  \phi_1^{M^{\lt}}) \in \Spec(J, \lt)$ 
with Bishop spaces $(\C G_j)_{j \in J}$ and Bishop morphisms $(\mu_{jj{'}})_{(j, j{'}) \in D^{\lt}(J)}$,
there is a bijection
\[\theta : \underset{\to} \Lim \big(\lambda_0(i) \times \mu_0 (j)\big) \to \underset{\to} \Lim 
\lambda_0(i) \times 
\underset{\to} \Lim \mu_0(j) \in \Mor \big(\underset{\to} \Lim (\C F_i \times \C G_j), \underset{\to} 
\Lim \C F_i \times \underset{\to} \Lim \C G_j\big). \]
\end{proposition}

\begin{proof}
Let the operation  $\theta : \underset{\to} \Lim \big(\lambda_0(i) \times \mu_0 (j)\big) \sto 
\underset{\to} \Lim \lambda_0(i) \times \underset{\to} \Lim \mu_0(j)$, defined by 
\[ \theta\big(\eql_0^{\Lambda^{\lt} \times M^{\lt}}\big((i, j), (x, y)\big)\big) := \big(\eql_0^{\Lambda^{\lt}}(i, x), 
\eql_0^{M^{\lt}}(j, y)\big). \]
First we show that $\theta$ is an embedding as follows:
\begin{align*}
& \  \eql_0^{\Lambda^{\lt} \times M^{\lt}}\big((i, j), (x, y)\big) = \eql_0^{\Lambda^{\lt} 
\times M^{\lt}}\big((i{'}, j{'}), (x{'}, y{'})\big) :\TOT \\
& : \TOT \exists_{(k, l) \in I \times J}\big((i, j), (i{'}, j{'}) \lt (k, l) \ \& \ 
(\lambda_1 \times \mu_1)_{(i,j)(k,l)}^{\lt}(x, y) = 
(\lambda_1 \times \mu_1)_{(i{'},j{'})(k, l)}^{\lt}(x{'}, y{'})\big)\\
& : \TOT \exists_{(k, l) \in I \times J}\big((i, j), (i{'}, j{'}) \lt (k, l) \ \& \ 
\big(\lambda_{ik}^{\lt}(x), \mu_{jl}^{\lt}(y)\big) = \big(\lambda_{i{'}k}^{\lt}(x{'}), 
\mu_{j{'}k}^{\lt}(y{'})\big)\\
& \TOT \exists_{k \in I}\big(i, i{'} \lt k \ \& \ \lambda_{ik}^{\lt}(x) = \lambda_{i{'}k}^{\lt}(x{'})\big)
\ \& \
\exists_{l \in J}\big(j, j{'} \lt l \ \& \ \lambda_{jk}^{\lt}(y) = \lambda_{j{'}k}^{\lt}(y{'})\big)\\
& : \TOT \eql_0^{\Lambda^{\lt}}(i, x) = \eql_0^{\Lambda^{\lt}}(i{'}, x{'}) \ \ \& \ \ \eql_0^{M^{\lt}}(j, y) = 
\eql_0^{M^{\lt}}(j{'}, y{'})\\
& : \TOT \big(\eql_0^{\Lambda^{\lt}}(i, x), \eql_0^{M^{\lt}}(j, y)\big) = \big(\eql_0^{\Lambda^{\lt}}(i{'}, x{'}), 
\eql_0^{M^{\lt}}(j{'}, y{'})\big)\\
& : \TOT \theta\big(\eql_0^{\Lambda^{\lt} \times M^{\lt}}\big((i, j), (x, y)\big)\big) = 
\theta\big(\eql_0^{\Lambda^{\lt} \times M^{\lt}}\big((i{'}, j{'}), (x{'}, y{'})\big)\big).
\end{align*}
The fact that $\theta$ is a surjection 
is immediate to show. 
By definition of the direct limit and by the $\vee$-lifting of the product Bishop topology we have that
\[ \underset{\to} \Lim (\C F_i \times \C G_j) := \bigg(\underset{\to} \Lim \big(\lambda_0(i) \times 
\mu_0 (j)\big), \bigvee_{\Xi \in \prod_{(i, j) \in I \times J}^{\mt}F_i \times 
G_j}\eql_0 f_{\Xi}\bigg),\]
\[ \underset{\to} \Lim \C F_i \times \underset{\to} \Lim \C G_j := \bigg(\underset{\to} \Lim \lambda_0(i) 
\times \underset{\to} \Lim \mu_0(j), \bigvee_{\Theta \in \prod_{i \in I}^{\mt}F_i}^{H 
\in \prod_{j \in J}^{\mt}G_j}\eql_0 f_{\Theta} \circ \pi_1, \eql_0 f_{H} \circ \pi_2 \bigg). \]
To show that $\theta \in \Mor \big(\underset{\to} \Lim (\C F_i \times \C G_j), \underset{\to} \Lim \C F_i 
\times \underset{\to} \Lim \C G_j\big)$ it suffices to show that
\[ \forall_{\Theta \in \prod_{i \in I}^{\mt}F_i}\forall_{H \in \prod_{j 
\in J}^{\mt}G_j}\big(\eql_0 f_{\Theta} \circ \pi_1) \circ \theta \in \underset{\to} \Lim 
(F_i \times G_j) \ \& \ (\eql_0 f_{H} \circ \pi_2) \circ \theta \in \underset{\to} \Lim (F_i 
\times G_j)\big). \]
If $\Theta \in \prod_{i \in I}^{\mt}F_i$, we show that $(\eql_0 f_{\Theta} \circ \pi_1) 
\circ \theta \in \underset{\to} \Lim (F_i \times G_j)$  
From the equalities
\begin{align*}
[(\eql_0 f_{\Theta^{\mt}} \circ \pi_1) \circ \theta]\big(\eql_0^{\Lambda^{\lt} \times M^{\lt}}\big((i, j),
(x, y)\big)\big) & := (\eql_0 f_{\Theta^{\mt}} \circ \pi_1)\big(\eql_0^{\Lambda^{\lt}}(i, x), 
\eql_0^{M^{\lt}}(j, y)\big)\\
& := \eql_0 f_{\Theta^{\mt}}\big(\eql_0^{\Lambda^{\lt}}(i, x)\big)\\
& := \Theta_i(x)\\
& := \big(\Theta_i \circ \pi_1\big)(x, y)\\
& := \big[\Theta_1(i, j)\big](x, y)\\
& := \eql_0 f_{\Theta_1}\big(\eql_0^{\Lambda^{\lt} \times M^{\lt}}\big((i, j), (x, y)\big),
\end{align*}
where $\Theta_1 \in \prod_{i \in I}^{\mt} F_i \times G_j$ is defined in Lemma~\ref{lem: prodlemma},
we conclude that $(\eql_0 f_{\Theta^{\mt}} \circ \pi_1) \circ \theta := \eql_0 f_{\Theta_1} \in
\underset{\to} \Lim (F_i \times G_j)$. For the second case we work similarly.
\end{proof}

\section{Inverse limit of a contravariant spectrum of Bishop spaces}
\label{sec: inverselimit}

\begin{definition}\label{def: inverselimit}
If $S(\Lambda^{\mt}) := (\lambda_0, \lambda_1^{\mt}, \phi_0^{\Lambda^{\mt}}, \phi_1^{\Lambda^{\mt}})$ 
is a contravariant $(I, \lt)$-spectrum with Bishop spaces $(\C F_i)_{i \in I}$ and Bishop morphisms 
$(\lambda_{ji}^{\mt})_{(i,j) \in \ltI}$, the 
\textit{inverse limit} of $S\Lambda(^{\mt})$\index{inverse limit of a contravariant direct spectrum of 
Bishop spaces} is the Bishop space
\[ \underset{\ot} \Lim \C F_i := \big(\underset{\ot} \Lim \lambda_0 (i), \ \underset{\ot} \Lim F_i\big), \]
\[ \underset{\ot} \Lim \lambda_0 (i) := \prod_{i \in I}^{\mt}\lambda_0 (i) \ \ \ \& \ \ \ \
\underset{\ot} \Lim F_i := \bigvee_{i \in I}^{f \in F_i}f \circ \pi_i^{\Lambda^{\mt}}. \]
\end{definition}

For simplicity we write $\pi_i$ instead of $\pi_i^{\Lambda^{\mt}}$ for the 
function $\pi_i^{\Lambda^{\mt}} : \prod_{i \in I}^{\mt}\lambda_0 (i) \to \lambda_0(i)$, which is 
defined, as its dual $\pi_i^{\Lambda^{\lt}}$ in the Proposition~\ref{prp: preorderfamilymap1}(iv), by
the rule $\Phi \mapsto \Phi_i$, for every $i \in I$..

\begin{proposition}[Universal property of the inverse limit]\label{prp: universalinverse}
 If  $S(\Lambda^{\mt}) := (\lambda_0, \lambda_1^{\mt}, \phi_0^{\Lambda^{\mt}}, \phi_1^{\Lambda^{\mt}})$ 
is a contravariant direct spectrum over $(I, \lt)$ with Bishop spaces $(\C F_i)_{i \in I}$ and 
Bishop morphisms $(\lambda_{ji}^{\mt})_{(i,j) \in \lt(I)}$, its inverse limit $\underset{\ot} \Lim \C F_i$
satisfies the universal property of inverse limits i.e.,\\[1mm]
\normalfont (i)
\itshape  For every $i \in I$, we have that $\pi_i \in \Mor(\underset{\ot} \Lim \C F_i, \C F_i)$.\\[1mm]
\normalfont (ii)
\itshape If $i \lt j$, the following left diagram commutes 
\begin{center}
\begin{tikzpicture}

\node (E) at (0,0) {$\prod_{i \in I}^{\mt}\lambda_0(i)$};
\node[below=of E] (F) {};
\node [right=of F] (B) {$\lambda_0(j)$};
\node [left=of F] (C) {$\lambda_0(i)$};
\node [right=of B] (D) {$\lambda_0(i)$};
\node[right=of D] (K) {};
\node[right=of K] (L) {$\lambda_0(j)$.};
\node[above=of K] (M) {$Y$};

\draw[->] (E)--(B) node [midway,right] {$ \ \  \pi_j$};
\draw[->] (E)--(C) node [midway,left] {$\pi_i \ $};
\draw[->] (B)--(C) node [midway,below] {$\lambda_{ji}^{\mt}$};

\draw[->] (M)--(L) node [midway,right] {$ \ \varpi_j$};
\draw[->] (M)--(D) node [midway,left] {$\varpi_i \ $};
\draw[->] (L)--(D) node [midway,below] {$\lambda_{ji}^{\lt}$};

\end{tikzpicture}
\end{center}
\normalfont (iii)
\itshape If $\C G := (Y, G)$ is a Bishop space and 
$\varpi_i : Y \to \lambda_0(i) \in \Mor(\C G, \C F_i)$, for every $i \in I$, such that if 
$i \lt j$,
the above right diagram commutes, there is a 
unique function $h : Y \to \prod_{i \in I}\lambda_0(i) \in \Mor(\C G, \underset{\ot} \Lim \C F_i)$
such that the following diagrams commute
\begin{center}
\begin{tikzpicture}

\node (E) at (0,0) {$Y$};
\node[below=of E] (F) {};
\node [right=of F] (B) {$\lambda_0(j)$,};
\node [left=of F] (C) {$\lambda_0(i)$};
\node[below=of F] (G) {$\prod_{i \in I}\lambda_0(i)$};

\draw[->] (E)--(B) node [midway,right] {$ \ \  \varpi_j$};
\draw[->] (E)--(C) node [midway,left] {$\varpi_i \ $};
\draw[left hook->] (B)--(C) node [midway,below] {$\lambda_{ji}^{\mt} \ \ \ \ \ \  $};
\draw[->] (G)--(B) node [midway,right] {$ \ \  \pi_j$};
\draw[->] (G)--(C) node [midway,left] {$\pi_i \ $};
\draw[dashed,->] (E)--(G) node [midway,near start] {$\ \ \ h$};

\end{tikzpicture}
\end{center}
\end{proposition}

\begin{proof}
The condition
$\pi_i \in \Mor(\underset{\ot} \Lim \C F_i, \C F_i) :\TOT \forall_{f \in F_i}\bigg(f \circ \pi_i \in
\bigvee_{i \in I}^{f \in F_i}f \circ \pi_i\bigg)$ is trivially satisfied, and (i) follows.
For (ii), the required equality $\lambda_{ji}^{\mt}\big(\pi_j(\Phi)\big) 
=_{\mathsmaller{\lambda_0(i)}} \pi_i(\Phi) :\TOT \lambda_{ji}^{\mt}(\Phi_j) 
= _{\mathsmaller{\lambda_0(i)}} \Phi_i$ holds by 
the definition of $\prod_{i \in I}^{\mt}\lambda_0(i)$.
To show (iii), let the operation $h : Y \sto \prod_{i \in I}^{\mt}\lambda_0(i)$, defined by 
$h(y) := \Phi_y$, where $\Phi_y(i) := \varpi_i(y)$, for every $y \in Y$ and $i \in I$.
First we show that $h$ is well-defined i.e., $h(y) \in \prod_{i \in I}^{\mt}\lambda_0(i)$. If $i \lt j$, by the 
supposed commutativity of the above right diagram we have that
$\lambda_{ji}^{\mt}\big(\Phi_y(j)\big) :=  \lambda_{ji}^{\mt}\big(\varpi_j(y)\big) = \varpi_i(y) 
:= \Phi_y(i)$.
Next we show that $h$ is a function.
If $y =_Y y{'}$, 
the last formula in the following equivalences
\[ \Phi_y =_{\mathsmaller{\prod_{i \in I}^{\mt}\lambda_0(i)}} \Phi_{y{'}} :\TOT \forall_{i \in I}\big(\Phi_y(i)
=_{\mathsmaller{\lambda_0(i)}} \Phi_{y{'}}(i)\big) : \TOT \forall_{i \in I}\big(\varpi_i(y)
=_{\mathsmaller{\lambda_0(i)}} \varpi_i(y{'})\big) \]
holds by the fact that $\varpi_i$ is a function, for every $i \in I$. By the 
$\bigvee$-lifting of morphisms we have that 
$h \in \Mor(\C G, \underset{\ot} \Lim \C F_i) \TOT \forall_{i \in I}\forall_{f \in F_i}\big((f \circ \pi_i)
\circ h \in G\big)$. 
If $i \in I$, $f \in F_i$, and $y \in Y$, then
\[ \big[(f \circ \pi_i) \circ h \big](y) := (f \circ \pi_i)(\Phi_y) := f\big(\Phi_y(i)\big) := 
f\big(\varpi_i(y)\big) := (f \circ \varpi_i)(y), \]
hence $(f \circ \pi_i) \circ h := f \circ \varpi_i \in G$, since $\varpi_i \in \Mor(\C G, \C F_i)$. The required
commutativity of the last diagram above, and the uniqueness of $h$ follow immediately.
\end{proof}

The uniqueness of $\underset{\ot} \Lim \lambda_0 (i)$, up to Bishop isomorphism, follows easily from 
its universal property. Next follows the inverse analogue to the Theorem~\ref{thm: directspectrummap1}.

\begin{theorem}\label{thm: inverselimitmap}
Let $S(\Lambda^{\mt}) := (\lambda_0, \lambda_1^{\mt}, \phi_0^{\Lambda^{\mt}}, \phi_1^{\Lambda^{\mt}})$ be
a contravariant $(I, \lt)$-spectrum with Bishop spaces $(\C F_i)_{i \in I}$ and Bishop morphisms 
$(\lambda_{ji}^{\mt})_{(i,j) \in \ltI}$, $S(M^{\mt}) := (\mu_0, \mu_1, \phi_{0}^{M^{\mt}}, \phi_{1}^{M^{\mt}})$
a contravariant $(I, \lt)$-spectrum with Bishop spaces $(\C G_i)_{i \in I}$ and Bishop morphisms 
$(\mu_{ji}^{\mt})_{(i,j) \in \ltI}$, and $\Psi \colon S(\Lambda^{\mt}) \To S(M^{\mt})$.\\[1mm]
\normalfont (i)
\itshape There is a unique function $\Psi_{\ot} \colon \underset{\ot} \Lim \lambda_0 (i) \to 
\underset{\ot} \Lim \mu_0 (i)$ such that, for every $i \in I$, the following diagram commutes
\begin{center}
\begin{tikzpicture}

\node (E) at (0,0) {$\underset{\ot} \Lim \lambda_0 (i)$};
\node[right=of E] (F) {$\underset{\ot} \Lim \mu_0 (i)$.};
\node[above=of F] (A) {$\mu_0(i)$};
\node [above=of E] (D) {$\lambda_0(i)$};

\draw[dashed,->] (E)--(F) node [midway,below]{$\Psi_{\ot}$};
\draw[->] (D)--(A) node [midway,above] {$\Psi_i$};
\draw[->] (E)--(D) node [midway,left] {$\pi^{\Lambda^{\mt}}_i$};
\draw[->] (F)--(A) node [midway,right] {$\pi^{M^{\mt}}_i$};

\end{tikzpicture}
\end{center}
\normalfont (ii)
\itshape If $\Psi$ is continuous, then 
$\Psi_{\ot} \in \Mor(\underset{\ot} \Lim \C F_i, \underset{\ot} \Lim \C G_i)$.\\[1mm]
\normalfont (iii)
\itshape If $\Psi_i$ is an embedding, for every $i \in I$,
then $\Psi_{\ot}$ is an embedding.
\end{theorem}

\begin{proof}
(i) Let the assignment routine $\Psi_{\ot} \colon \underset{\ot} \Lim \lambda_0 (i) \sto 
\underset{\ot} \Lim \mu_0 (i)$, defined by 
\[ \Theta \mapsto \Psi_{\ot}(\Theta), \ \ \ \  \big[ \Psi_{\ot}(\Theta) \big]_i := \Psi_i(\Theta_i);  \ \ \ \ 
\Theta \in \underset{\ot} \Lim \lambda_0(i), \ i \in I. \]
First we show that $\Psi_{\ot}$ is well-defined i.e., $\Psi_{\ot}(\Theta) \in \prod_{i \in I}^{\mt}\mu_0 (i).$
If $i \lt j$, since $\Theta \in \prod_{i \in I}^{\mt}\lambda_0 (i)$, we have that 
$\Theta_i = \lambda_{ji}^{\mt}(\Theta_j)$, and since $\Psi \colon S(\Lambda^{\mt}) \To S(M^{\mt})$
\begin{center}
\begin{tikzpicture}

\node (E) at (0,0) {$\mu_0(i)$};
\node[right=of E] (F) {$\mu_0(j)$,};
\node[above=of E] (A) {$\lambda_0(i)$};
\node [right=of A] (D) {$\lambda_0(j)$};

\draw[->] (F)--(E) node [midway,below]{$\mu_{ji}^{\mt}$};
\draw[->] (D)--(A) node [midway,above] {$\lambda_{ji}^{\mt}$};
\draw[->] (A)--(E) node [midway,left] {$\Psi_i$};
\draw[->] (D)--(F) node [midway,right] {$\Psi_j$};

\end{tikzpicture}
\end{center}
$\big[\Psi_{\ot}(\Theta)\big]_i := \Psi_i(\Theta_i)
= \Psi_i\big(\lambda_{ji}^{\mt}(\Theta_j)\big)
:= \big(\Psi_i \circ \lambda_{ji}^{\mt}\big)(\Theta_j) 
= \big(\mu_{ji}^{\mt} \circ \Psi_j\big)(\Theta_j) 
 := \mu_{ji}^{\mt}\big(\Psi_j(\Theta_j)\big)
:= \mu_{ji}^{\mt}\big(\big[\Psi_{\ot}(\Theta)\big]_j\big).
$
Next we show that $\Psi_{\ot}$ is a function: 
$\Theta =_{\mathsmaller{\underset{\ot} \Lim \lambda_0 (i)}} \Phi : \TOT
\forall_{i \in I}\big(\Theta_i =_{\mathsmaller{\lambda_0 (i)}} \Phi_i\big)
 \To \forall_{i \in I}\big(\Psi_i(\Theta_i) =_{\mathsmaller{\mu_0 (i)}} \Psi_i(\Phi_i)\big) 
: \TOT \forall_{i \in I}\big(\big[\Psi_{\ot}(\Theta)\big]_i =_{\mathsmaller{\mu_0 (i)}}
 \big[\Psi_{\ot}(\Phi)\big]_i\big)
 : \TOT \Psi_{\ot}(\Theta) =_{\mathsmaller{\underset{\ot} \Lim \mu_0 (i)}}  \Psi_{\ot}(\Phi)$.
The commutativity of the diagram and the uniqueness of $\Psi_{\ot}$ are immediate to show.\\
(ii) By the $\bigvee$-lifting of morphisms we have that
$ \Psi_{\ot} \in \Mor(\underset{\ot} \Lim \C F_i, \underset{\ot} \Lim \C G_i) \TOT 
\forall_{i \in I}\forall_{g \in G_i}\big((g \circ \pi_i^{M^{\mt}}) \circ \Psi_{\ot}
\in \underset{\ot} \Lim F_i\big). $
If $i \in I$ and $g \in G_i$, then
$[(g \circ \pi_i^{M^{\mt}}) \circ \Psi_{\ot}](\Theta) :=
g \big([\Psi_{\ot}(\Theta)\big]_i\big)
:= g \big(\Psi_{i}(\Theta_i)\big)
:= \big(g \circ \Psi_{i}\big)(\Theta_i)
:= \big[\big(g \circ \Psi_{i}\big) \circ \pi_i^{\Lambda^{\mt}}\big](\Theta),
$
and $g \circ \Psi_{i} \in F_i$, by the continuity of $\Psi$, hence
$(g \circ \pi_i^{M^{\mt}}) \circ \Psi_{\ot} := 
\big(g \circ \Psi_{i}\big) \circ \pi_i^{\Lambda^{\mt}} \in \underset{\ot} \Lim F_i$.\\
(iii) By definition we have that
$
\Psi_{\ot}(\Theta) =_{\mathsmaller{\underset{\ot} \Lim \mu_0 (i)}}  \Psi_{\ot}(\Phi)
: \TOT \forall_{i \in I}\big(\Psi_i(\Theta_i) =_{\mathsmaller{\mu_0 (i)}} \Psi_i(\Phi_i)\big)
\To $
$\forall_{i \in I}\big(\Theta_i =_{\mathsmaller{\lambda_0 (i)}} \Phi_i\big)
: \TOT \Theta =_{\mathsmaller{\underset{\ot} \Lim \lambda_0 (i)}} \Phi.$
\end{proof}

\begin{proposition}\label{prp: transitivity2}
If $S(\Lambda^{\mt}) := (\lambda_0, \lambda_1^{\mt}, \phi_0^{\Lambda^{\mt}}, \phi_1^{\Lambda^{\mt}})$,  
$S(M^{\mt}) := (\mu_0, \mu_1^{\mt}, \phi_0^{M^{\mt}}, \phi_1^{M^{\mt}})$ and
$S(N^{\mt} ):= (\nu_0, \nu_1^{\mt}, \phi_0^{N^{\mt}}, \phi_1^{N^{\mt}})$
are contravariant direct spectra over $(I, \lt)$, and if $\Psi \colon S(\Lambda^{\mt}) \To S(M^{\mt})$ and 
$\Xi \colon S(M^{\mt}) \To S(N^{\mt})$, then 
$(\Xi \circ \Psi)_{\ot} :=  \Xi_{\ot} \circ \Psi_{\ot}$
\begin{center}
\begin{tikzpicture}

\node (E) at (0,0) {$\lambda_0(i)$};
\node[right=of E] (H) {};
\node[right=of H] (F) {$\mu_0(i)$};
\node[below=of E] (A) {$\underset{\ot} \Lim \lambda_0 (i)$};
\node[below=of F] (B) {$\underset{\ot} \Lim \mu_0 (i)$};
\node[right=of F] (K) {};
\node[right=of K] (G) {$\nu_0(i)$};
\node [below=of G] (C) {$\underset{\ot} \Lim \nu_0 (i)$.};

\draw[->] (E)--(F) node [midway,above] {$ \Psi_i$};
\draw[->] (F)--(G) node [midway,above] {$\Xi_i$};
\draw[->] (B)--(C) node [midway,below] {$\Xi_{\ot} \ \ $};
\draw[->] (A)--(B) node [midway,below] {$ \ \ \Psi_{\ot}$};
\draw[->] (A)--(E) node [midway,left] {$\pi^{\Lambda^{\lt}}_i$};
\draw[->] (B) to node [midway,left] {$\pi^{M^{\lt}}_i$} (F);
\draw[->] (C)--(G) node [midway,right] {$\pi^{N^{\lt}}_i$};
\draw[->,bend right] (A) to node [midway,below] {$(\Xi \circ \Psi)_{\ot}$} (C) ;
\draw[->,bend left] (E) to node [midway,above] {$(\Xi \circ \Psi)_i$} (G) ;

\end{tikzpicture}
\end{center}
\end{proposition}

\begin{proof}
The required equality is reduced to
$\forall_{i \in I}\big(\big[(\Xi \circ \Psi)_{\ot}(\Theta)\big]_i =_{\mathsmaller{\nu_0 (i)}}
\big[\Xi_{\ot}(\Psi_{\ot}(\Theta))\big]_i\big)$. 
If $i \in I$, then 
$\big[(\Xi \circ \Psi)_{\ot}(\Theta)\big]_i := (\Xi \circ \Psi)_i(\Theta_i) := \Xi_i(\Psi_i(\Theta_i))
:= \Xi_i\big(\big[\Psi_{\ot}(\Theta)\big]_i\big) := \big[\Xi_{\ot}(\Psi_{\ot}(\Theta))\big]_i$.
\end{proof}

\begin{theorem}\label{thm: cofinal3}
Let  $S(\Lambda^{\mt}) := (\lambda_0, \lambda_1^{\mt}, \phi_0^{\Lambda^{\mt}}, \phi_1^{\Lambda^{\mt}})$
be a contravariant direct spectrum over $(I, \lt)$, $(J, e, \cof_J)$ a cofinal subset of $I$, and $S(\Lambda^{\mt} 
\circ e) :=
\big((\lambda_0 \circ e, \lambda_1 \circ e, \phi_0^{\Lambda^{\mt} \circ e}, \phi_1^{\Lambda^{\mt} \circ e}\big)$ 
the relative spectrum of $S(\Lambda^{\mt})$ to $J$. Then  
\[ \underset{\ot} \Lim \C F_j \simeq \underset{\ot} \Lim \C F_i.\]
\end{theorem}

\begin{proof}
If $\Theta \in \prod_{j \in J}^{\mt}\lambda_0(j)$, then, if $j \lt j{'}$, we have that
$\Theta_j = \lambda_{j{'}j}^{\mt}\big(\Theta_{j{'}}\big) := \lambda_{e(j{'})e(j)}^{\mt}\big(\Theta_{j{'}}\big)$.
If $i \in I$, then $\cof_J(i) \in J$ and $\Theta_{\cof_J(i)} \in \lambda_0(e(\cof_J(i)))$. Since
$i \lt e(\cof_J(i))$, we define the operation  
$\phi : \underset{\ot} \Lim \lambda_0 (j) \sto \underset{\ot} 
\Lim \lambda_0 (i)$, by the rule $\Theta \mapsto \phi(\Theta)$, for every $\Theta \in \underset{\ot} \Lim \lambda_0 (j)$,
where
\[ [\phi(\Theta)]_i := \lambda_{e(\cof_J(i)) i}^{\mt}\big(\Theta_{\cof_J(i)}\big) \in \lambda_0(i); \ \ \ \
i \in I. \]
First we show that $\phi$ is well-defined i.e., $\phi(\Theta) \in \prod_{i \in I}^{\mt}\lambda_0(i)$ i.e., 
for every $i, i{'} \in I$, 
$i \lt i{'} \To [\phi(\Theta)]_i = \lambda_{i{'}i}^{\mt}\big([\phi(\Theta)]_{i{'}}\big)$.
Working as in the proof of Lemma~\ref{lem: cofinallemma}(ii), we get
\begin{align*}
\lambda_{i{'}i}^{\mt}\big([\phi(\Theta)]_{i{'}}\big) & := 
\lambda_{i{'}i}^{\mt}\big(\lambda_{e(\cof_J(i{'})) i{'}}^{\mt}\big(\Theta_{\cof_J(i{'})}\big)\big)\\
& := \big[\lambda_{i{'}i}^{\mt} \circ \lambda_{e(\cof_J(i{'})) i{'}}^{\mt}\big]\big(\Theta_{\cof_J(i{'})}\big)\\
& \stackrel{(\ref{eq: eq1})} = \lambda_{e(\cof_J(i{'})) i}^{\mt} \big(\Theta_{\cof_J(i{'})}\big)\\
& \stackrel{(\ref{eq: eq2})} = \big[\lambda_{e(\cof_J(i))i}^{\mt} \circ 
\lambda_{e(\cof_J(i{'})) e(\cof_J(i))}^{\mt}\big] \big(\Theta_{\cof_J(i{'})}\big)\\
& := \lambda_{e(\cof_J(i))i}^{\mt}\big(\lambda_{e(\cof_J(i{'})) 
e(\cof_J(i))}^{\mt}\big(\Theta_{\cof_J(i{'})}\big)\big)\\
& = \lambda_{e(\cof_J(i))i}^{\mt}\big(\Theta_{\cof_J(i)}\big)\\
& := [\phi(\Theta)]_i.
\end{align*}
To show that $\phi$ is a function we consider the following equivalences:
\begin{align*}
\phi(\Theta) =_{\mathsmaller{\underset{\ot} \Lim \lambda_0 (i)}} \phi(H) & :\TOT
\forall_{i \in I}\big(\big[\phi(\Theta)\big]_i =_{\mathsmaller{\lambda_0 (i)}} \big[\phi(H)\big]_i\big)\\
& : \TOT \forall_{i \in I}\big(\lambda_{e(\cof_J(i)) i}^{\mt}\big(\Theta_{\cof_J(i)}\big)
=_{\mathsmaller{\lambda_0 (i)}} \lambda_{e(\cof_J(i)) i}^{\mt}\big(H_{\cof_J(i)}\big)\big), \ \ \ (1)
\end{align*}
\[ \Theta =_{\mathsmaller{\underset{\ot} \Lim \lambda_0 (j)}} H : 
\TOT \forall_{j \in J}\big(\Theta_j =_{\mathsmaller{\lambda_0 (e(j))}} H_j\big) \ \ \ \ \ (2). \]
To show that (1) $\To$ (2) we use the fact that $e(\cof_J(j)) = j$, and since $j \lt j$, by the 
extensionality of $\lt$ we get $j \lt e(\cof_J(j))$. Since  
$\Theta_j = \lambda_{e(\cof_J(i)) i}^{\mt}\big(\Theta_{\cof_J(i)}\big)$, and 
$H_j = \lambda_{e(\cof_J(i)) i}^{\mt}\big(H_{\cof_J(i)}\big)$, we get (2).
By the $\bigvee$-lifting of morphisms
$\phi \in \Mor(\underset{\ot} \Lim \C F_j, \underset{\ot} \Lim \C F_i) \TOT \forall_{i \in I}\forall_{f \in F_i}
\big((f \circ \pi_i^{S^{\mt}}) \circ \phi \in \underset{\ot} \Lim F_j\big)$. 
If $\Theta \in \prod_{j \in J}^{\mt}\lambda_0(j)$, we have that
\begin{align*}
[(f \circ \pi_i^{S^(\Lambda{\mt})}) \circ \phi](\Theta) & := f\big(\pi_i^{S(\Lambda^{\mt})}(\phi(\Theta))\big)\\
& := f\big(\big[\phi(\Theta)\big]_i\big)\\
& := f\big(\lambda_{e(\cof_J(i)) i}^{\mt}\big(\Theta_{\cof_J(i)}\big)\big)\\
& := \big(f \circ \lambda_{e(\cof_J(i)) i}^{\mt}\big)\big(\Theta_{\cof_J(i)}\big)\\
& := \big[\big(f \circ \lambda_{e(\cof_J(i)) i}^{\mt}\big) \circ \pi_{\cof_J(i)}^{S(\Lambda^{\mt}) \circ e}\big](\Theta),
\end{align*}
hence 
$\big(f \circ \pi_i^{S(\Lambda^{\mt})}\big) \circ \phi :=
\big(f \circ \lambda_{e(\cof_J(i)) i}^{\mt}\big) \circ \pi_{\cof_J(i)}^{S(\Lambda^{\mt}) \circ e}
\in \underset{\ot} \Lim F_j$, 
as by definition $\lambda_{e(\cof_J(i)) i}^{\mt} \in \Mor(\C F_{e(\cof_J(i))}, \C F_i)$, and hence
\begin{center}
\begin{tikzpicture}

\node (E) at (0,0) {$\lambda_0(e(\cof_J(i)))$};
\node[right=of E] (H) {};
\node[right=of H] (F) {$\lambda_0(i)$};
\node[below=of F] (A) {$\Real$};

\draw[->] (E)--(F) node [midway,above] {$\lambda_{e(\cof_J(i)) i}^{\mt} $};
\draw[->] (F)--(A) node [midway,right] {$f$};
\draw[->] (E)--(A) node [midway,left] {$f \circ \lambda_{e(\cof_J(i)) i}^{\mt} \ \  $};

\end{tikzpicture}
\end{center}
$f \circ \lambda_{e(\cof_J(i)) i}^{\mt} \in F_{e(\cof_J(i))} := F_{\cof_J(i)}$. 
Let the operation $\theta \colon \underset{\ot} \Lim \lambda_0 (i) \sto \underset{\ot} \Lim \lambda_0(j)$,
defined by the rule $H \mapsto \theta(H) := H^J$, for every  
$H \in \prod_{i \in I}^{\mt}\lambda_0(i) $, where
$H_j^J := H_{e(j)} \in \lambda_0(e(j))$, for every $j \in J$. We show that $H^J \in \prod_{j \in J}\lambda_0(j)$. 
If $j \lt j{'}$, then 
\[ H_j^J = \lambda_{e(j{'})e(j)}^{\mt}\big(H_{j{'}}^J\big) :\TOT   
H_{e(j)} = \lambda_{e(j{'})e(j)}^{\mt}\big(H_{e(j{'})}\big), \]
which holds by the hypothesis $H \in \prod_{i \in I}^{\mt}\lambda_0(i)$. Moreover,
we have that $\phi(H^J) = H :\TOT \forall_{i \in I}\big(\big[\phi(H^J)\big]_i =_{\mathsmaller{\lambda_0(i)}} H_i\big)$.
If $i \in I$, and since $i \lt e(\cof_J(i))$, we have that 
\[ \big[\phi(H^J)\big]_i := \lambda_{e(\cof_J(i)) i}^{\mt}\big(H^J_{\cof_J(i)}\big)
:= \lambda_{e(\cof_J(i)) i}^{\mt}\big(H_{e(\cof_J(i))}\big)
= H_i. \]
It is immediate to show that $\theta$ is a function. Moreover, $\theta(\phi(\Theta) = \Theta$, as if $j \in J$, then
\[ \phi(\Theta)_j^J := \phi(\Theta)_{e(j)}
:= \lambda^{\mt}_{e(\cof_J(e(j))e(j)}\big(\Theta_{\cof_J(e(j))}\big)
= \Theta_j,
\]
as by hypothesis $\Theta_j = \lambda_{e(j{'})e(j)}^{\mt}(\Theta_{j{'}})$, with $j \lt j{'}$, and by $(\Cf_1)$ we 
have that
$j =_J (\cof_J(e(j))$, hence by the extensionality of $\lt$ we get $j \lt (\cof_J(e(j))$. Finally,
$\theta \in \Mor \big(\underset{\ot} \Lim \C F_i, \underset{\ot} \Lim \C F_j \big) \TOT \forall_{j
\in J}\forall_{f \in F_j}\big(\big(f \circ \pi_j^{S(\Lambda^{\mt})}\big) \circ \theta \in \underset{\ot} 
\Lim F_i\big)$, which follows from the equalities
\begin{align*}
\big[\big(f \circ \pi_j^{S(\Lambda^{\mt})}\big) \circ \theta \big](H) & := \big(f \circ
\pi_j^{S(\Lambda^{\mt})}\big)(H^J)\\
& := f\big(H^J_j\big)\\
& := f \big(H_{e(j)}\big)\\
& := \big(f \circ \pi_{e(j)}^{S(\Lambda^{\mt})}\big)(H).\qedhere
\end{align*}
\end{proof}

\begin{proposition}\label{prp: prodinverse}
If $(I, \lt), (J, \lt)$ are directed sets, $S(\Lambda^{\mt}) := (\lambda_0, \lambda_1^{\mt},
\phi_0^{\Lambda^{\mt}},  \phi_1^{\Lambda^{\mt}})$ is a contravariant direct spectrum over $(I, \lt)$
with Bishop spaces $(\C F_i)_{i \in I}$ and Bishop morphisms $(\lambda_{i{'}i}^{\mt})_{(i, i{'}) \in \ltI}$,
and $S(M^{\mt}) := (\mu_0, \mu_1^{\mt}, \phi_0^{M^{\mt}},  \phi_1^{M^{\mt}})$ is a contravariant direct
spectrum over $(J, \lt)$ with Bishop spaces $(\C G_j)_{j \in J}$ and Bishop morphisms 
$(\mu_{j{'}j}^{\mt})_{(j, j{'}) \in D^{\lt}(J)}$, there is a function
\[ \times : \prod_{i \in I}^{\mt}\lambda_0 (i) \times \prod_{j \in J}^{\mt}\mu_0 (j) \to 
\prod_{(i, j) \in I \times J}^{\mt}\lambda_0 (i) \times \mu_0 (j) \in \Mor\big(\underset{\ot} 
\Lim \C F_i \times \underset{\ot} \Lim \C G_j, \underset{\ot} \Lim (\C F_i \times \C G_j)\big). \]
\end{proposition}

\begin{proof}
We proceed as in the proof of Proposition~\ref{prp: proddirect}. 
\end{proof}

\section{Duality between direct and inverse limits of spectra}
\label{sec: duality}

\begin{proposition}\label{prp: plus}
Let $\C F := (X, F), \C G := (Y, G)$ and $\C H := (Z, H)$ be Bishop spaces, and let
$\lambda \in \Mor(\C G, \C H),
\mu \in \Mor(\C H, \C G)$. We define the mappings
\[ \lambda^{+} \colon \Mor(\C H, \C F) \to \Mor(\C G, \C F),  \ \ \ \  \lambda^{+} (\phi) := \phi 
\circ \lambda; \ \ \ \ \phi \in \Mor(\C H, \C F), \]
\[ \mu^{-} \colon \Mor(\C F, \C H) \to \Mor(\C F, \C G), \ \ \ \ \mu^{-} (\theta) := \mu \circ \theta; \ \ \ \ 
\theta \in \Mor(\C F, \C H), \]
\begin{center}
\begin{tikzpicture}

\node (E) at (0,0) {$Y$};
\node[above=of E] (F) {$Z$};
\node[right=of F] (A) {$X$};
\node[right=of A] (B) {$Z$};
\node[below=of B] (C) {$Y$.};

\draw[->] (E)--(F) node [midway,left] {$\lambda$};
\draw[->] (F)--(A) node [midway,above] {$\phi$};
\draw[->] (E)--(A) node [midway,right] {$\mathsmaller{\phi \circ \lambda}$};
\draw[->] (A)--(C) node [midway,left] {$\mathsmaller{\mu \circ \theta} \ $};
\draw[->] (A)--(B) node [midway,above] {$\theta$};
\draw[->] (B)--(C) node [midway,right] {$\mu$};

\end{tikzpicture}
\end{center}
\[ ^{+} \colon \Mor(\C G, \C H) \to \Mor(\C H \to \C F, \C G \to \C F) , \ \ \ \ \lambda \mapsto \lambda^+; \ \ \ \ 
\lambda \in \Mor(\C G, \C H), \]
\[ ^{-} \colon \Mor(\C H, \C G) \to \Mor(\C F \to \C H, \C F \to \C G), \ \ \ \ \mu \mapsto \mu^-; \ \ \ \ 
\mu \in \Mor(\C H, \C G). \]
Then $^{+} \in \Mor\big(\C G \to \C H, (\C H \to \C F) \to (\C G \to \C F)\big)$ and 
$^{-} \in \Mor\big(\C H \to\C G, (\C F \to \C H) \to (\C F \to \C G)\big)$.
\end{proposition}

\begin{proof}
By definition and the $\bigvee$-lifting of the exponential topology we have that
\[ \C G \to \C H := \bigg(\Mor(\C G, \C H), \bigvee_{y \in Y}^{h \in H}\phi_{y, h}\bigg), \ \ \ \
\C H \to \C F := \bigg(\Mor(\C H, \C F), \bigvee_{z \in Z}^{f\in F}\phi_{z, f}\bigg), \]
\[ \C G \to \C F := \bigg(\Mor(\C G, \C F), \bigvee_{y \in Y}^{f \in F}\phi_{y, f}\bigg), \]
\[ (\C H \to \C F) \to (\C G \to \C F) := \bigg(\Mor((\C H, \C F) \to \C G \to \C F)), 
\bigvee_{\varphi \in \Mor(\C H, \C F)}^{e \in G \to F}\phi_{\varphi, e}\bigg), \]
\[ \bigvee_{\varphi \in \Mor(\C H, \C F)}^{e \in G \to F}\phi_{\varphi, e} = 
\bigvee_{\varphi \in \Mor(\C H, \C F)}^{y \in Y, f \in F}\phi_{\varphi, \phi_{y, f}}. \]
By the  $\bigvee$-lifting of morphisms we have that 
\[ ^{+} \in \Mor\big(\C G \to \C H, (\C H \to \C F) \to (\C G \to \C F)\big) \TOT 
\forall_{\varphi \in \Mor(\C H, \C F)}\forall_{y \in Y}\forall_{f \in F}\big(\phi_{\varphi, 
\phi_{y, f}} \circ \ ^{+} \in G \to H\big). \]
If $\lambda \in \Mor(\C G, \C H)$, we have that
$ \phi_{\varphi, \phi_{y, f}} \circ \ ^{+}](\lambda) := \phi_{\varphi, \phi_{y, f}}(\lambda^{+})
 := (\phi_{y, f} \circ \lambda^+)(\varphi) := (\phi_{y, f}(\varphi \circ \lambda)
 := [f \circ (\varphi \circ \lambda)](y) := [(f \circ \varphi) \circ \lambda](y)
:= \phi_{y, f \circ \varphi}(\lambda)$
i.e., $\phi_{\varphi, \phi_{y, f}} \circ \ ^{+} :=  \phi_{y, f \circ \varphi} \in G \to H$,
since $\varphi\in \Mor(\C H, \C F)$ and hence $f \circ \varphi \in H$. For the mapping $^{-}$ we work similarly.
\end{proof}

Next we see how with the use of the exponential Bishop topology we can get a contravatiant spectrum
from a covariant one, and vice versa.

\begin{proposition}\label{prp: dirtoproj}
\normalfont (A)
\itshape Let $S(\Lambda^{\lt}) := (\lambda_0, \lambda_1^{\lt}, \phi_0^{\Lambda^{\lt}}, \phi_1^{\Lambda^{\lt}}) \in \Spec(I, \lt)$ and $\C F := (X, F)$ a Bishop space. \\[1mm]
\normalfont (i)
\itshape If $S(\Lambda^{\lt}) \to \C F := (\mu_0, \mu_1^{\mt}, \phi_0^{M^{\mt}}, \phi_1^{M^{\mt}})$, where 
$M^{\mt} := (\mu_0, \mu_1^{\mt})$ 
is a contravariant direct family of sets over $(I, \lt)$ with $\mu_0 (i) := \Mor(\C F_i, \C F)$ and 
\[ \mu_1^{\mt}(i, j) := \big(\Mor(\C F_j, \C F), \Mor(\C F_i, \C F), (\lambda_{ij}^{\lt})^{+} \big), \]
and if $\phi_0^{M^{\mt}}(i) := F_i \to F$ and $\phi_1^{M^{\mt}}(i, j) := \big(F_i \to F, F_j \to F, 
[(\lambda_{ij}^{\lt})^{+}]^*\big)$, then $S^{\lt} \to \C F$ is a contravariant $(I, \lt)$-spectrum
with Bishop spaces $(\Mor(\C F_i, \C F))_{i \in I}$ and Bishop morphisms 
$\big((\lambda_{ij}^{\lt})^{+}\big)_{(i, j) \in \ltI}$.\\[1mm]
\normalfont (ii)
\itshape If $\C F \to S^{\lt} := (\nu_0, \nu_1^{\lt}, \phi_0^{N^{\lt}}, \phi_1^{N^{\lt}})$, where 
$N^{\lt} := (\nu_0, \nu_1^{\lt})$ is a direct family of sets over $(I, \lt)$ with $\nu_0 (i) := \Mor(\C F, \C F_i)$ and 
\[ \nu_1^{\lt}(i, j) := \big(\Mor(\C F, \C F_i), \Mor(\C F, \C F_j), (\lambda_{ij}^{\lt})^{-} \big), \]
and if $\phi_0^{N^{\lt}}(i) := F \to F_i$ and $\phi_1^{N^{\lt}}(i, j) := \big(F \to F_j, F \to F_i,
[(\lambda_{ij}^{\lt})^{-}]^*\big)$, then $\C F \to S^{\lt}$ is a covariant $(I, \lt)$-spectrum with
Bishop spaces $(\Mor(\C F, \C F_i))_{i \in I}$ and Bishop morphisms 
$\big((\lambda_{ij}^{\lt})^{-}\big)_{(i, j) \in \ltI}$.\\[1mm]
\normalfont (B)
\itshape Let $S(\Lambda^{\mt}) := (\lambda_0, \lambda_1^{\mt}, \phi_0^{\Lambda^{\mt}}, \phi_1^{\Lambda^{\mt}})$
be a contravariant $(I, \lt)$-spectrum, and $\C F := (X, F)$ a Bishop space. \\[1mm]
\normalfont (i)
\itshape If $S(\Lambda^{\mt}) \to \C F := (\mu_0, \mu_1^{\lt}, \phi_0^{M^{\lt}}, \phi_1^{M^{\lt}})$, where 
$M^{\lt} := (\mu_0, \mu_1^{\lt})$ 
is a direct family of sets over $(I, \lt)$ with $\mu_0 (i) := \Mor(\C F_i, \C F)$ and 
\[ \mu_1^{\lt}(i, j) := \big(\Mor(\C F_i, \C F), \Mor(\C F_j, \C F), (\lambda_{ji}^{\mt})^{+} \big), \]
and if $\phi_0^{M^{\lt}}(i) := F_i \to F$ and $\phi_1^{M^{\lt}}(i, j) := \big(F_j \to F, F_i \to F, 
[(\lambda_{ji}^{\mt})^{+}]^*\big)$, then $S^{\mt} \to \C F$ is an $(I, \lt)$-spectrum with
Bishop spaces $(\Mor(\C F_i, \C F))_{i \in I}$ and Bishop morphisms 
$\big((\lambda_{ji}^{\mt})^{+}\big)_{(i, j) \in \ltI}$.\\[1mm]
\normalfont (ii)
\itshape If $\C F \to S(N^{\mt}) := (\nu_0, \nu_1^{\mt}, \phi_0^{N^{\mt}}, \phi_1^{N^{\mt}})$, where 
$N^{\mt} := (\nu_0, \nu_1^{\mt})$ 
is a contravariant direct family of sets over $(I, \lt)$ with $\nu_0 (i) := \Mor(\C F, \C F_i)$ and 
\[ \nu_1^{\mt}(i, j) := \big(\Mor(\C F, \C F_j), \Mor(\C F, \C F_i), (\lambda_{ji}^{\mt})^{-} \big), \]
and if $\phi_0^{N^{\mt}}(i) := F \to F_i$ and $\phi_1^{N^{\mt}}(i, j) := \big(F \to F_i, F \to F_j, 
[(\lambda_{ji}^{\mt})^{-}]^*\big)$, then $\C F \to S^{\lt}$ is a contravariant $(I, \lt)$-spectrum with
Bishop spaces $(\Mor(\C F, \C F_i)_{i \in I}$ and Bishop morphisms 
$\big((\lambda_{ij}^{\lt})^{-}\big)_{(i, j) \in \ltI}$.
\end{proposition}

\begin{proof}
We prove only the case (A)(i) and for the other cases we work similarly. 
It suffices to show that if $i \lt j \lt k$, then the following diagram commutes
\begin{center}
\begin{tikzpicture}

\node (E) at (0,0) {$\Mor(F_j, \C F)$};
\node[right=of E] (F) {$\Mor(\C F_k, \C F).$};
\node [above=of E] (D) {$\Mor(\C F_i, \C F)$};

\draw[->] (F)--(E) node [midway,below] {$\big(\lambda_{jk}^{\lt}\big)^+$};
\draw[->] (E)--(D) node [midway,left] {$\big(\lambda_{ij}^{\lt}\big)^+$};
\draw[->] (F)--(D) node [midway,right] {$\ \ \ \big(\lambda_{ik}^{\lt}\big)^+$};

\end{tikzpicture}
\end{center}
If $\phi \in \Mor(\C F_k, \C F)$, then
$\big(\lambda_{ij}^{\lt}\big)^+\big[\big(\lambda_{jk}^{\lt}\big)^+(\phi)\big]  := 
\big(\lambda_{ij}^{\lt}\big)^+[\phi \circ \lambda_{jk}^{\lt}]
:= (\phi \circ \lambda_{jk}^{\lt}) \circ \lambda_{ij}^{\lt}
:= \phi \circ (\lambda_{jk}^{\lt} \circ \lambda_{ij}^{\lt})
= \phi \circ \lambda_{ik}^{\lt}
:= \big(\lambda_{ik}^{\lt}\big)^+ (\phi)$.
\end{proof}

Similarly to the $\bigvee$-lifting of the product topology, if 
$S(\Lambda^{\mt}) := (\lambda_0, \lambda_1^{\mt} ,\phi_0^{\Lambda^{\mt}},
\phi_1^{\Lambda^{\mt}})$ a contravariant direct spectrum over $(I, \lt)$ with Bishop spaces
$\big(F_i = \bigvee F_{0i}\big)_{i \in I}$, then 
\[ \prod_{ i \in I}^{\mt} F_i = \bigvee_{i \in I}^{f \in F_{0i}}\big(f \circ \pi_i^{\Lambda^{\mt}}\big). \]

\begin{theorem}[Duality principle]\label{thm: duality1}
Let $S(\Lambda^{\lt}) := (\lambda_0, \lambda_1^{\lt}, \phi_0^{\Lambda^{\lt}},
\phi_1^{\Lambda^{\lt}}) \in \Spec(I, \lt)$ with Bishop spaces $(\C F_i)_{i \in I}$
and Bishop morphisms $(\lambda_{ij}^{\lt})_{(i, j) \in \ltI}$. If $\C F := (X, F)$ is a Bishop space and 
$S(\Lambda^{\lt}) \to \C F := (\mu_0, \mu_1^{\mt}, \phi_0^{M^{\mt}}, \phi_1^{M^{\mt}})$ is the contravariant 
direct spectrum over $(I, \lt)$ defined in Proposition~\ref{prp: dirtoproj} \normalfont (A)(i),
\itshape then 
\[ \underset{\ot} \Lim (\C F_i \to \C F) \simeq [(\underset{\to} \Lim \C F_i) \to \C F]. \]
\end{theorem}

\begin{proof}
First we determine the topologies involved in the required  Bishop isomorphism.
By definition and by the above remark on the $\bigvee$-lifting of the $\prod^{\mt}$-topology we have that
\[ \underset{\ot} \Lim (\C F_i \to \C F) := \bigg(\prod_{i \in I}^{\mt}\mu_0 (i), 
\bigvee_{i \in I}^{g \in F_i \to F}g \circ \pi_i^{S(\Lambda^{\lt}) \to \C F}\bigg), \]
\[ F_i \to F := \bigvee_{x \in \lambda_0 (i)}^{f \in F}\phi_{x, f}, \]
\[ \bigvee_{i \in I}^{g \in F_i \to F}g \circ \pi_i^{S^{\lt} \to \C F} = 
\bigvee_{i \in I}^{x \in \lambda_0 (i), f \in F}\phi_{x, f} \circ \pi_i^{S(\Lambda^{\lt}) \to \C F}, \]
\[ \underset{\to} \Lim \C F_i := \bigg(\underset{\to} \Lim \lambda_0 (i), \bigvee_{\mathsmaller{\Theta
\in \prod_{i \in I}^{\lt}F_i}}\eql_0 f_{\Theta}\bigg), \]
\[ \big(\underset{\to} \Lim \C F_i\big) \to \C F := \bigg(\Mor(\underset{\to} 
\Lim \C F_i, \C F), \bigvee_{\mathsmaller{\eql_0^{\Lambda^{\lt}}(i, x) \in 
\underset{\to} \Lim \lambda_0 (i)}}^{f \in F}
\phi_{\mathsmaller{\eql_0^{\Lambda^{\lt}}(i, x)}, f}\bigg), \]
\[ \phi_{\mathsmaller{\eql_0^{\Lambda^{\lt}}(i, x)}, f}(h) := (f \circ h)\big(\eql_0^{\Lambda^{\lt}}(i, x)\big) \]
\begin{center}
\begin{tikzpicture}

\node (E) at (0,0) {$\Real$.};
\node[above=of E] (F) {$X$};
\node[left=of F] (B) {};
\node[left=of B] (A) {$\underset{\to} \Lim \lambda_0 (i)$};

\draw[->] (F)--(E) node [midway,right] {$f \in F$};
\draw[->] (A)--(E) node [midway,left] {$f \circ h \ \ \ $};
\draw[->] (A)--(F) node [midway,above] {$\mathsmaller{h \in \Mor(\underset{\to} \Lim \C F_i, \C F)}$};

\end{tikzpicture}
\end{center}
If $H \in \prod_{i \in I}^{\mt}\Mor(\C F_i, \C F)$, let the operation $\theta(H) :
\underset{\to} \Lim \lambda_0(i) \sto X$, defined by
\[ \theta(H)\big(\eql_0^{\Lambda^{\lt}}(i, x)\big) := H_i(x); \ \ \ \ \eql_0^{\Lambda^{\lt}}(i, x) \in 
\underset{\to} \Lim \lambda_0(i). \]
We show that $\theta(H)$ is a function. If
\[ \eql_0^{\Lambda^{\lt}}(i, x) =_{\mathsmaller{\underset{\to} \Lim \lambda_0 (i)}} \eql_0^{\Lambda^{\lt}}(j, y) \TOT
\exists_{k \in I}\big(i, j \lt k \ \& \ \lambda_{ik}^{\lt}(x) =_{\mathsmaller{\lambda_0 (k)}} 
\lambda_{jk}^{\lt}(y)\big),\]
we show that
$ \theta(H)\big(\eql_0^{\Lambda^{\lt}}(i, x)\big) :=  H_i(x) =_X H_j(y) =:
\theta(H)\big(\eql_0^{\Lambda^{\lt}}(j, y)\big)$. 
By the equalities
$H_i = \big(\lambda_{ik}^{\lt}\big)^+(H_k) = H_k \circ \lambda_{ik}^{\lt}$ and  
$H_j = \big(\lambda_{jk}^{\lt}\big)^+(H_k) = H_k \circ \lambda_{jk}^{\lt}$ we get
\[ H_i(x) = \big(H_k \circ \lambda_{ik}^{\lt}\big)(x) := H_k\big(\lambda_{ik}^{\lt}(x)\big) =_X 
H_k\big(\lambda_{jk}^{\lt}(y)\big) := \big(H_k \circ \lambda_{jk}^{\lt}\big)(y) := H_j(y). \]
Next we show that 
$\theta(H) \in \Mor(\underset{\to} \Lim \C F_i, \C F) :\TOT \forall_{f \in F}\big(f 
\circ \theta(H) \in \underset{\to} \Lim  F_i\big)$. 
If $f \in F$, then the dependent assignment routine
$\Theta :\bigcurlywedge_{i \in I}F_i$, defined by $\Theta_i := f \circ H_i$, for every $i \in I$ 
\begin{center}
\begin{tikzpicture}

\node (E) at (0,0) {$\Real$};
\node[above=of E] (F) {$X$};
\node[left=of F] (B) {};
\node[left=of B] (A) {$\lambda_0 (i)$};

\draw[->] (F)--(E) node [midway,right] {$f \in F$};
\draw[->] (A)--(E) node [midway,left] {$f \circ H_i \ \ \ $};
\draw[->] (A)--(F) node [midway,above] {$\mathsmaller{H_i \in \Mor(\C F_i, \C F)}$};

\end{tikzpicture}
\end{center}
is in $\prod_{i \in I}^{\lt}F_i$ i.e., if $i \lt j$, then $\Theta_i = \big(\lambda_{ij}^{\lt}\big)^*(\Theta_j)
= \Theta_j \circ \lambda_{ij}^{\lt}$, since 
$\Theta_i := f \circ H_i = f \circ \big(H_j \circ \lambda_{ij}^{\lt}\big) = (f \circ H_j) \circ 
\lambda_{ij}^{\lt} := \Theta_j \circ \lambda_{ij}^{\lt}$.
Hence $f \circ \theta(H) := \eql_0 f_{\Theta} \in \underset{\to} \Lim  F_i$, since 
\[
[f \circ \theta(H)]\big(\eql_0^{\Lambda^{\lt}}(i, x)\big) := f\big(H_i(x)\big)
:= (f \circ H_i)(x)
:= f_{\Theta}(i, x)
:= \eql_0 f_{\Theta}\big(\eql_0^{\Lambda^{\lt}}(i, x)\big).
\]
Consequently, the operation $\theta \colon \prod_{i \in I}^{\mt}\Mor(\C F_i, \C F)
\sto \Mor(\underset{\to} \Lim \C F_i, \C F)$, defined by the rule $H \mapsto \theta(H)$,
is well-defined. Next we show that $\theta$ is an embedding. 
\begin{align*}
\theta(H) = \theta(K) & : \TOT \forall_{\mathsmaller{\eql_0^{\Lambda^{\lt}}(i, x) \in 
\underset{\to} \Lim \lambda_0 (i)}}\big(\theta(H)(\eql_0^{\Lambda^{\lt}}(i, x)) =  
\theta(K)(\omega_{S^{\lt}}(i, x))\big)\\
& : \TOT \forall_{i \in I}\big(H_i(x) =_X K_i(x)\big)\\
& : \TOT H = K.
\end{align*}
Next we show that $\theta \in \Mor\big(\underset{\ot} \Lim (\C F_i \to \C F), 
(\underset{\to} \Lim \C F_i) \to \C F\big)$ i.e.,
\[ \forall_{\mathsmaller{\eql_0^{\Lambda^{\lt}}(i, x) \in \underset{\to} \Lim 
\lambda_0 (i)}}\forall_{f \in F}\bigg(\phi_{\mathsmaller{\eql_0^{\Lambda^{\lt}}(i, x)}, f} \circ 
\theta \in  \bigvee_{i \in I, x \in \lambda_0 (i)}^{f \in F}\phi_{x, f} \circ \pi_i^{S(\Lambda^{\lt})
\to \C F} \bigg). \]
By the equalities
 \[ [\phi_{\mathsmaller{\eql_0^{\Lambda^{\lt}(i, x)}, f}} \circ \theta](H) := 
\phi_{\mathsmaller{\eql_0^{\Lambda^{\lt}}(i, x), f}}(\theta(H))
:= [f \circ \theta(H)]\big(\eql_0^{\Lambda^{\lt}}(i, x)\big) 
:= f \big(H_i(x)\big), \]
\[ [\phi_{x, f} \circ \pi_i^{S^{\lt} \to \C F}](H) := \phi_{x, f}\big(H_i)\big)
:= f\big(H_i(x)\big),  \]
we get $\phi_{\mathsmaller{\eql_0^{\Lambda^{\lt}}(i, x)}, f} \circ \theta = \phi_{x, f} \circ
\pi_i^{S(\Lambda^{\lt}) \to \C F}$. Let $\phi \colon \Mor(\underset{\to} \Lim \C F_i, \C F) \sto
\prod_{i \in I}^{\mt}\Mor(\C F_i, \C F)$ be defined by $h \mapsto \phi(h) := H^h$, where 
$ H^{h} : \bigcurlywedge_{i \in I}\Mor(\C F_i, \C F)$ is defined by  
$H^{h}_i := h \circ \eql_i$, for every $i \in I$
\begin{center}
\begin{tikzpicture}

\node (E) at (0,0) {$X$.};
\node[above=of E] (F) {$\underset{\to} \Lim \lambda_0 (i)$};
\node[left=of F] (B) {};
\node[left=of B] (A) {$\lambda_0 (i)$};

\draw[->] (F)--(E) node [midway,right] {$h$};
\draw[->] (A)--(E) node [midway,left] {$H_i^h \ \ \ $};
\draw[->] (A)--(F) node [midway,above] {$\eql_{i}$};

\end{tikzpicture}
\end{center}
By Proposition~\ref{prp: universaldirect}(i) $H_i \in \Mor(\C F_i, \C F) $,
as a composition of Bishop morphisms. To show that $H^{h} \in 
\prod_{i \in I}\Mor(\C F_i, \C F)$, let $i \lt j$, and by Proposition~\ref{prp: universaldirect}(ii)
we get
$H^{h}_i := h \circ \eql_i = h \circ \big(\eql_j \circ \lambda_{ij}^{\lt}\big) :=
(h \circ \eql_j) \circ \lambda_{ij}^{\lt} := H_j \circ \lambda_{ij}^{\lt}$.
Moreover, $\theta(H_{h}) := h$, since
$ \theta(H_{h})\big(\eql_0^{\Lambda^{\lt}}(i, x)\big) := H_{i}(x) 
:= (h \circ \eql_{i}(x) := h \big(\eql_0^{\Lambda^{\lt}}(i, x)\big)$.
Clearly, $\phi$ is a function. Moreover $H^{\theta(H)} := H$, as, for every $i \in I$ we have that
$ \big(H_i^{\theta(H)}\big)(x) := (\theta(H) \circ \eql_i)(x) 
:= \theta(H)\big(\eql_0^{\Lambda^{\lt}}(i,x)\big) := H_i(x)$.
Finally we show that $\phi \in \Mor\big((\underset{]to} \Lim \C F_i) \to \C F, \underset{\ot}
\Lim (\C F_i \to \C F)\big)$
if and only if 
\[\forall_{i \in I}\forall_{x \in \lambda_0(i)}\forall_{f \in F}\big(\phi_{x,f} \circ \phi \in
\bigvee_{\eql_0^{\Lambda^{\lt}}(i,x) \in \underset{\to}  \Lim \lambda_0(i)}^{f \in F}
\phi_{\eql_0^{\Lambda^{\lt}}(i,x), f}\big). \]
If $h \in \Mor(\underset{\to} \Lim \C F_i, \C F)$, then
\[ \big[\phi_{x,f} \circ \pi_i^{S(\Lambda^{\lt}) \to F}) \circ \phi\big](h := (\phi_{x,f} \circ \pi_i^{S(\Lambda^{\lt}) 
\to F})(H^h)
 := \phi_{x,f}\big(H^h_i\big) \]
\[ := \phi_{x, f}(h \circ \eql_i)
:= f\big[(h \circ \eql_i)(x)\big]
 := (f \circ h)\big(\eql_0^{\Lambda^{\lt}}(i,x)\big)
  := \phi_{\eql_0^{\Lambda^{\lt}}(i,x),f}(h).\qedhere
\]
\end{proof}

With respect to the possible dual to the previous theorem i.e., the isomorphism 
$\underset{\to} \Lim (\C F_i \to \C F) \simeq [(\underset{\ot} \Lim \C F_i) \to \C F]$, what we can 
show is the following proposition.

\begin{proposition}\label{prp: conversedual}
Let
$S(\Lambda^{\mt}) := (\lambda_0, \lambda_1^{\mt}, \phi_0^{\Lambda^{\mt}},
\phi_1^{\Lambda^{\mt}})$ be a contravariant direct spectrum over $(I, \lt)$ with Bishop spaces 
$(\C F_i)_{i \in I}$ and Bishop morphisms
$(\lambda_{ji}^{\succ})_{(i,j) \in \ltI}$. If $\C F := (X, F)$ is a Bishop space and 
$S(\Lambda^{\mt}) \to \C F := (\mu_0, \mu_1^{\lt}, \phi_0^{M^{\lt}}, \phi_1^{M^{\lt}})$ is the $(I, \lt)$-directed
spectrum defined in Proposition~\ref{prp: dirtoproj} \normalfont (B)(i), 
\itshape
there is  a function 
$\ \widehat{} \ \colon \underset{\to} \Lim [\Mor(\C F_i, \C F)] \to \Mor(\underset{\ot} \Lim \C F_i, \C F)$
such that the following hold:\\[1mm]
\normalfont (i)
\itshape $\ \widehat{} \ \in \Mor\big(\underset{\to} \Lim (\C F_i \to \C F), 
(\underset{\ot} \Lim \C F_i) \to \C F\big)$.\\[1mm]
\normalfont (ii)
\itshape  If for every $j \in J$ and every $y \in \lambda_0 (j)$ there is 
$\Theta_y \in \prod_{i \in I}^{\mt}\lambda_0 (i)$ 
such that $\Theta_y(j) =_{\mathsmaller{\lambda_0 (j)}} y$, 
then $\ \widehat{} \ $ is an embedding of $\underset{\to} \Lim [\Mor(\C F_i, \C F)]$ into 
$\Mor(\underset{\ot} \Lim \C F_i, \C F)$.
\end{proposition}

\begin{proof}
We proceed similarly to the proof of Theorem~\ref{thm: duality1}. 
\end{proof}

\begin{theorem}\label{thm: duality2}
Let $S(\Lambda^{\mt}) := (\lambda_0, \lambda_1^{\mt}, \phi_0^{\Lambda^{\mt}},
\phi_1^{\Lambda^{\mt}})$ be a contravariant  direct spectrum over $(I, \lt)$ with Bishop spaces 
$(\C F_i)_{i \in I}$ and Bishop morphisms $(\lambda_{ji}^{\prec})_{(i,j) \in \ltI}$.
If $\C F := (X, F)$ is a Bishop space and 
$\C F \to S(\Lambda^{\mt}) := (\nu_0, \nu_1^{\mt}, \phi_0^{^{\mt}}, \phi_1^{N^{\mt}})$ is the contravariant
direct spectrum over $(I, \lt)$, defined in Proposition~\ref{prp: dirtoproj} \normalfont (B)(ii), 
\itshape then 
\[ \underset{\ot} \Lim (\C F \to \C F_i) \simeq [\C F \to \underset{\ot} \Lim \C F_i]. \]
\end{theorem}

\begin{proof}

First we determine the topologies involved in the required Bishop isomorphism:
\[ \underset{\ot} \Lim (\C F \to \C F_i) := \bigg(\prod_{i \in I}^{\mt}\Mor(\C F, \C F_i), 
\bigvee_{i \in I}^{g \in F \to F_i}g \circ \pi_i^{\C F \to S(\Lambda^{\mt})}\bigg), \]
\[ \bigvee_{i \in I}^{g \in F \to F_i}g \circ \pi_i^{\C F \to S^{\mt}} = \bigvee_{i \in I, x \in 
\lambda_0 (i)}^{f \in F_i}\phi_{x, f} \circ \pi_i^{\C F \to S(\Lambda^{\mt})}, \]
\[ \underset{\ot} \Lim \C F_i := \bigg(\prod_{i \in I}^{\mt}\lambda_0 (i), 
\bigvee_{i \in I}^{f \in F_i}f \circ \pi_i^{S(\Lambda^{\mt})}\bigg), \]
\[ \C F \to \underset{\ot} \Lim \C F_i := \bigg(\Mor(\C F, \underset{\ot} \Lim \C F_i), 
\bigvee_{x \in X}^{g \in \underset{\ot} \Lim \C F_i}\phi_{x, g}\bigg), \]
 \[ \bigvee_{x \in X}^{g \in \underset{\ot} \Lim \C F_i}\phi_{x, g} = 
\bigvee_{x \in X, i \in I}^{f \in F_i}\phi_{x, f \circ \pi_i^{S(\Lambda^{\mt})}}. \]
If $H \in \prod_{i \in I}^{\mt}\Mor(\C F, \C F_i)$, and if $i \lt j$, then
$H_i = \nu_{ji}^{\mt}(H_j) = \big(\lambda_{ji}^{\mt}\big)^-(H_j) = \lambda_{ji}^{\mt} \circ H_j $
\begin{center}
\begin{tikzpicture}

\node (E) at (0,0) {$\lambda_0(i)$.};
\node[above=of E] (F) {$\lambda_0 (j)$};
\node[left=of F] (A) {$X$};

\draw[->] (F)--(E) node [midway,right] {$\lambda_{ji}^{\mt}$};
\draw[->] (A)--(E) node [midway,left] {$H_i \  $};
\draw[->] (A)--(F) node [midway,above] {$H_j$};

\end{tikzpicture}
\end{center}
Let the operation $e(H) : X \sto \prod_{i \in I}^{\mt} \lambda_0 (i)$, defined by
$x \mapsto [e(H)](x)$, where $\big[[e(H)](x)\big]_i := H_i(x)$, for every $i \in I$.
First we show that $[e(H)](x) \in \prod_{i \in I}^{\mt} \lambda_0 (i)$. If $i \lt j$, then 
$\big[[e(H)](x)\big]_i := H_i(x) = \big(\lambda_{ji}^{\mt} \circ H_j\big)(x) := 
\lambda_{ji}^{\mt}\big(H_j(x)\big) := \lambda_{ji}^{\mt}\big(\big[[e(H)](x)\big]_j\big)$.
Next we show that $e(H)$ is a function. If $x =_X x{'}$, then
$\forall_{i \in I}\big(H_i(x) =_{\mathsmaller{\lambda_0 (i)}} H_i(x{'})\big) : \TOT
\forall_{i \in I}\big(\big[[e(H)](x)\big]_i  =_{\mathsmaller{\lambda_0 (i)}} \big[[e(H)](x{'})\big]_i\big) 
 : \TOT [e(H)](x) =_{\mathsmaller{\prod_{i \in I}^{\mt}\lambda_0 (i)}} [e(H)](x{'})$.
By the $\bigvee$-lifting of morphisms  
$e(H) \in \Mor(\C F, \underset{\ot} \Lim \C F_i) \TOT 
\forall_{i \in I}\forall_{f \in F_i}\big(\big(f \circ \pi_i^{S(\Lambda^{\mt})}\big) \circ e(H) \in F\big).$ 
Since $ [(f \circ \pi_i^{S(\Lambda^{\mt})}) \circ e(H)](x) := 
\big(f \circ \pi_i^{S(\Lambda^{\mt})}\big)\big([e(H)](x)\big) 
:= f \big(H_i(x)\big) := (f \circ H_i)(x)$, 
we get $\big(f \circ \pi_i^{S(\Lambda^{\mt})}\big) \circ e(H) :=  f \circ H_i \in F$, since $f \in F_i$ 
and $H_i \in \Mor(\C F, \C F_i)$. Hence, the 
operation $e : \prod_{i \in I}^{\mt}\Mor(\C F, \C F_i) \sto \Mor(\C F, 
\underset{\ot} \Lim \C F_i)$, defined by the rule $H \mapsto e(H)$, 
is well-defined. Next we show that $e$ is an embedding. If $H, K \in \prod_{i \in I}^{\mt}\Mor(\C F, \C F_i)$,
then
\begin{align*}
e(H) = e(K) & : \TOT \forall_{x \in X}\big([e(H)](x) =_{\mathsmaller{\prod_{i \in I}^{\mt}\lambda_0 (i)}}
[e(K)](x)\big)\\
& : \TOT \forall_{x \in X}\forall_{i \in I}\big(H_i(x)  =_{\mathsmaller{\lambda_0 (i)}} K_i(x)\big)\\
& : \TOT \forall_{i \in I}\forall_{x \in X}\big(H_i(x)  =_{\mathsmaller{\lambda_0 (i)}} K_i(x)\big)\\
& :\TOT \forall_{i \in I}\big(H_i =_{\mathsmaller{\Mor(\C F, \C F_i)}} K_i\big)\\
& :\TOT H =_{\mathsmaller{\prod_{i \in I}^{\mt}\Mor(\C F, \C F_i)}} K.
\end{align*}
By the $\bigvee$-lifting of morphisms we show that 
\[ e \in \Mor(\underset{\ot} \Lim (\C F \to \C F_i), \C F \to \underset{\ot} \Lim \C F_i) 
\TOT \forall_{i \in I}\forall_{f \in F_i}\big(\phi_{x, f \circ \pi_i^{S(\Lambda^{\mt})}} \circ e 
\in \underset{\ot} \Lim (F \to F_i)\big) \]
\begin{align*}
\big(\phi_{x, f \circ \pi_i^{S(\Lambda^{\mt})}} \circ e \big)(H) & := 
\phi_{x, f \circ \pi_i^{S(\Lambda^{\mt})}}\big(e(H)\big)\\
& := \big[(f \circ \pi_i^{S(\Lambda^{\mt})}) \circ e(H)\big](x)\\
& := (f \circ \pi_i^{S(\Lambda^{\mt})})\big([e(H)](x)\big)\\
& := f \big([e(H)(x)]_i\big)\\
& := f \big(H_i(x)\big)\\
& := (f \circ  H_i)(x)\\
& := \phi_{x, f}\big(H_i\big)\\
& := \big[\phi_{x, f} \circ \pi_i^{\C F \to S(\Lambda^{\mt})}\big](H)
\end{align*}
we get $\phi_{x, f \circ \pi_i^{S(\Lambda^{\mt})}} \circ e := \phi_{x, f} \circ \pi_i^{\C F \to S(\Lambda^{\mt})} 
\in \underset{\ot} \Lim (F \to F_i)$. Let 
$\phi \colon \Mor(\C F, \underset{\ot} \Lim \C F_i) \sto \prod_{i \in I}^{\mt} \Mor(\C F, \C F_i)$, 
defined by the rule $\mu \mapsto H^{\mu}$, where 
for every $\mu : X \to \prod_{i \in I}^{\mt}\lambda_0 (i) \in \Mor(\C F, \underset{\ot} \Lim \C F_i)$ i.e.,
$\forall_{i \in I}\forall_{f \in F_i}\big(\big(f \circ \pi_i^{S(\Lambda^{\mt})}\big) \circ \mu \in F\big)$, let
\[ H^{\mu} : \bigcurlywedge_{i \in I}\Mor(\C F, \C F_i), \ \ \ [H_{\mu}]_i \colon X \to \lambda_0(i),  \ \ \ \ 
 H^{\mu}_i(x) := [\mu(x)]_i;  \ \ \ \ x \in X, \ i \in I. \]
First we show that $H^{\mu}_i \in \Mor(\C F, \C F_i) :\TOT
\forall_{f \in F_i}\big(f \circ H^{\mu}_i \in F\big)$. If $f \in F_i$, and $x \in X$, then 
$[f \circ H^{\mu}_i(x) := f\big(H^{\mu}_i\big) := f\big([\mu(x)]_i\big)
:= \big[\big(f \circ \pi_i^{S^(\Lambda{\mt})}\big) \circ \mu\big](x)$
i.e., $f \circ H^{\mu}_i := \big(f \circ \pi_i^{S(\Lambda^{\mt})}\big) \circ \mu \in F$,
as $\mu \in \Mor(\C F, \underset{\ot} \Lim \C F_i)$.
Since $\mu(x) \in \prod_{i \in I}^{\mt}\lambda_0 (i)$, 
$[\mu(x)]_i = \lambda_{ji}^{\mt}\big([\mu(x)]_j\big),$
for every $i, j \in I$ such that $i \lt j$. To show that $H_{\mu} \in \prod_{i \in I}^{\mt}\Mor(\C F, \C F_i)$,
let $i \lt j$. Then
\begin{align*}
H^{\mu}_i = \lambda_{ji}^{\mt} \circ H^{\mu}_j & \TOT \forall_{x \in X}\big(H^{\mu}_i(x) 
=_{\mathsmaller{\lambda_0(i)}} \big[\lambda_{ji}^{\mt} \circ H^{\mu}_j\big](x)\big)\\
& : \TOT \forall_{x \in X}\big([\mu(x)]_i 
=_{\mathsmaller{\lambda_0(i)}} \big[\lambda_{ji}^{\mt}\big([\mu(x)]_j\big)\big),
\end{align*}
which holds by the previous remark on $\mu(x)$. It is immediate to show that $\phi$ is a function. 
To show that 
$\phi \in \Mor \big([ \C F \to \underset{\ot} \Lim \C F_i ], \underset{\ot} \Lim (\C F \to \C F_i)\big)$, 
we show that
\[ \forall_{i \in I}\forall_{f \in F_i}\forall_{x \in \lambda_0(i)}\bigg(\big[\phi_{x,f} \circ 
\pi_i^{\C F \to S(\Lambda^{\mt})}\big] \circ \phi \in \bigvee_{x \in X, i \in I}^{f \in F_i}\phi_{x, f 
\circ \pi_i^{S(\Lambda^{\mt})}}\bigg), \]
\begin{align*}
\big[\big[\phi_{x,f} \circ \pi_i^{\C F \to S(\Lambda^{\mt})}\big] \circ \phi \big](\mu) & := 
\big[\phi_{x,f} \circ \pi_i^{\C F \to S(\Lambda^{\mt})}\big](H^{\mu})\\
& := \phi_{x,f}\big(H^{\mu}_i\big)\\
& := (f \circ H_i^{\mu})(x) \\
& := f \big(\mu(x)_i\big)\\
& := (f \circ \pi_i^{S(\Lambda^{\mt})} \circ \mu)(x)\\
& := \big[\phi_{x, f \circ \pi_i^{S(\Lambda^{\mt})}}\big](\mu).
\end{align*}
Moreover, $\phi(e(H)) := H$, as $H^{e(H)}_i(x) := [e(H)(x)]_i := H_i(x)$, and $e(\phi(\mu)) = \mu$, as
\begin{align*}
e(H^{\mu}) = \mu & :\TOT \forall_{x \in X}\big([e(H^{\mu})](x) 
=_{\mathsmaller{\prod_{i \in I}^{\mt}\lambda_0 (i)}}
\mu(x)\big)\\
& :\TOT \forall_{x \in X}\forall_{i \in I}\big(H^{\mu}_i(x) =_{\mathsmaller{\lambda_0 (i)}}
[\mu(x)]_i\big)\\
& :\TOT \forall_{x \in X}\forall_{i \in I}\big([\mu(x)]_i =_{\mathsmaller{\lambda_0 (i)}}
[\mu(x)]_i\big).\qedhere
\end{align*}
\end{proof}

With respect to the possible dual to the previous theorem i.e., the isomorphism
$\underset{\to} \Lim (\C F \to \C F_i) \simeq (\C F \to \underset{\to} \Lim \C F_i)$, 
what we can show is the following proposition.

\begin{proposition}\label{prp: conversedual2}
Let $S(\Lambda^{\lt}) := (\lambda_0, \lambda_1^{\lt}, \phi_0^{\Lambda^{\lt}},
\phi_1^{\Lambda^{\lt}}) \in \Spec(I)$ with Bishop spaces $(\C F_i)_{i \in I}$
and Bishop morphisms $(\lambda_{ij}^{\prec})_{(i,j) \in \ltI}$. If $\C F := (X, F)$ is a Bishop space and 
$\C F \to S(\Lambda^{\lt}) := (\nu_0, \nu_1^{\lt}, \phi_0^{N^{\lt}}, \phi_1^{N^{\lt}})$ is the 
$(I, \lt)$-direct spectrum defined in Proposition~\ref{prp: dirtoproj} \normalfont (A)(ii),
\itshape 
there is  a map $\ \widehat{} \ : \underset{\to} \Lim [\Mor(\C F, \C F_i)] \to 
\Mor(\C F, \underset{\to} \Lim \C F_i)$
with $\ \widehat{} \ \in \Mor\big((\underset{\to} \Lim (\C F \to \C F_i), \C F \to \underset{\to} 
\Lim \C F_i\big)$.
\end{proposition}

\begin{proof}
We proceed similarly to the proof of Theorem~\ref{thm: duality2}. 
\end{proof}

\section{Spectra of Bishop subspaces}
\label{sec: spectrasub}

\begin{definition}\label{def: intspectrum}
If $\Lambda(X) := (\lambda_0, \C E^X, \lambda_1) \in \Fam(I, X)$, 
a family of Bishop subspaces of the Bishop space $\C F := (X, F)$ associated to $\Lambda(X)$\index{family
of Bishop subspaces associated to a family of subsets} is a pair $\Phi^{\Lambda(X)} := \big(\phi_0^{\Lambda(X)}, 
\phi_1^{\Lambda(X)}\big)$, where 
$\phi_0^{\Lambda(X)} \colon I \sto \D V_0$ and $\phi_1^{\Lambda(X)} : \bigcurlywedge_{(i,j) \in D(I)}
\D F\big(\phi_0^{\Lambda(X)}(i), \phi_0^{\Lambda(X)}(j)\big)$ such that the following conditions hold:\\[1mm]
\normalfont (i) 
\itshape $\phi_0^{\Lambda(X)}(i) := F_i := F_{|\lambda_0(i)} := \bigvee _{f \in F} f \circ \C E_i^X$,
for every $i \in I$.\\[1mm]
\normalfont (ii)
\itshape $\phi_1^{\Lambda(X)} (i,j) := \lambda_{ji}^*$, for every $(i, j) \in D(I)$. \\[1mm]
We call the structure
$S_F(\Lambda(X)) := (\lambda_0, \C E^X, \lambda_1, F, \phi_0^{\Lambda(X)}, \phi_1^{\Lambda(X)})$ 
a \textit{spectrum of subspaces of $\C F$} over $I$,\index{spectrum of subspaces} or an $I$-spectrum 
of subspaces of $\C F$ with Bishop subspaces $(\C F_i)_{i \in I}$ and Bishop morphisms\index{$S_F(\Lambda(X))$}
$(\C E_i^X)_{i \in I}$. If $S_F(M(X)) := (\mu_0, \C Z^X, \mu_1, F, \phi_0^{M(X)}, \phi_1^{M(X)})$
is an $I$-spectrum of subspaces of $\C F$ with Bishop subspaces $(\C G_i)_{i \in I}$ and Bishop morphisms
$(\C Z^X_i)_{i \in I}$, a \textit{subspaces spectrum-map}\index{subspaces spectrum-map} $\Psi$ from
$S_F(\Lambda(X))$ to $S_F(M(X))$, in symbols 
$\Psi \colon S_F(\Lambda(X)) \To S_F(M(X))$, is a family of subsets-map $\Psi \colon \Lambda(X) \To M(X)$. 
If $F$ is clear from the context, we may omit the symbol $F$ as a subscript in the above notations.
\end{definition}

The topology $F_i$ on $\lambda_0(i)$ is the relative Bishop topology of $F$ to $\lambda_0(i)$, and 
it is the least topology that makes the embedding $\C E_i^X$ a Bishop morphism from $\C F_i$ to $\C F$. 
In contrast to the external framework of a spectrum of Bishop spaces, we can prove that  the transport
maps $\lambda_{ij}$ of a spectrum of subspaces $S_F(\Lambda(X))$ are always Bishop morphisms. 
The extensionality of a Bishop topology $F$ on a set $X$ as a subset of $\D F(X)$ is crucial to the next proof.

\begin{remark}\label{rem: intspectrum1}
Let $S_F(\Lambda(X)) := (\lambda_0, \C E^X, \lambda_1, F, \phi_0^{\Lambda(X)}, \phi_1^{\Lambda(X)})$ be an
$I$-spectrum of subspaces of $\C F := (X, F)$ with Bishop subspaces $(\C F_i)_{i \in I}$ and Bishop morphisms
$(\C E_i^X)_{i \in I}$, and $S_F(M(X)) := (\mu_0, \C Z^X, \mu_1, F, \phi_0^{M(X)}, \phi_1^{M(X)})$
an $I$-spectrum of subspaces of $\C F$ with Bishop subspaces $(\C G_i)_{i \in I}$ and Bishop morphisms
$(\C Z^X_i)_{i \in I}$.\\[1mm]
\normalfont (i)
\itshape $S(\Lambda) := (\lambda_0, \lambda_1, \phi_0^{\Lambda(X)}, \phi_1^{\Lambda(X)})$ is 
an $I$-spectrum with Bishop spaces $(\C F_i)_{i \in I}$ and Bishop isomorphisms $(\lambda_{ij})_{(i,j) \in D(I)}$.\\[1mm]
\normalfont (ii)
\itshape If $\Psi \colon S(\Lambda(X)) \To S(M(X))$, then $\Psi$ is continuous i.e., $\Psi_i \in \Mor(\C F_i, \C G_i)$,
for every $i \in I$.
\end{remark}

\begin{proof}
(i) It suffices to show that $\lambda_{ij} \in \Mor(\C F_i, \C F_j)$, for every $(i, j) \in D(I)$. By the 
$\bigvee$-lifting of morphisms we have that
$
\lambda_{ij} \in \Mor(\C F_i, \C F_j) \TOT \forall_{f \in F}\big((f \circ \C E_j^X) \circ 
\lambda_{ij} \in F_i\big)
 :\TOT \forall_{f \in F}\big(f \circ (\C E_j^X \circ 
\lambda_{ij}) \in F_i\big)$. 
If we fix some $f \in F$, and as $\C E_j^X \circ \lambda_{ij} =_{\D F(\lambda_0(i), X)} \C E_i^X$, we get 
$f \circ (\C E_j^X \circ \lambda_{ij}) =_{\D F(\lambda_0(i))} f \circ \C E_i^X$. Since $f \circ \C E_i^X \in F_i$
by the extensionality of $F_i$ we get $f \circ (\C E_j^X \circ \lambda_{ij}) \in F_i$.\\
(ii) By the $\bigvee$-lifting of morphisms we have that
$\Psi_{i} \in \Mor(\C F_i, \C G_i) \TOT \forall_{f \in F}\big((f \circ \C Z^X_i) \circ 
\Psi_{i} \in F_i\big)
 :\TOT \forall_{f \in F}\big(f \circ (\C Z^X_i \circ \Psi_{i}) \in F_i\big)$. 
Since $\Psi \colon \Lambda(X) \To M(X)$, we get $\C Z_i^X \circ \Psi_i =_{\D F(\lambda_0(i), X)} \C E_i^X$, and 
hence $f \circ (\C Z_i^X \circ \Psi_i) =_{\D F(\lambda_0(i))} f \circ \C E_i^X$, for every $i \in I$ and $f \in F$.
By the definition of $F_i$ we have that $f \circ \C E_i^X \in F_i$, and hence by the extensionality of $F_i$ 
we conclude that $f \circ (\C Z_i^X \circ \Psi_i) \in F_i$.
\end{proof}

\begin{definition}\label{def: specsub}
Let $\Spec_F(I, X)$ be the totality of spectra of subspaces of the Bishop space $\C F := (X, F)$ over $I$,
equipped with the equality of $\Spec(I, X)$.
\end{definition}

\begin{definition}\label{def: canonicaltopsinternal}
Let $S_F(\Lambda(X)) := \big(\lambda_0, \C E^X, \lambda_1, F, \phi_0^{\Lambda(X)}, \phi_1^{\Lambda(X)}\big) 
\in \Spec_F(I, X)$ with Bishop subspaces $(\C F_i)_{i \in I}$ and Bishop morphisms $(\C E_i^X)_{i \in I}$.
The \textit{canonical Bishop topology on the interior union} $\bigcup_{i \in I}\lambda_0(i)$\index{canonical 
Bishop topology on the interior union of a spectrum of subsets} is the relative topology of $F$ to it i.e., 
\[ \bigcup_{i \in I}\C F_i := \bigg(\bigcup_{i \in I}\lambda_0(i), \bigcup_{i \in I}F_i\bigg), \]
\[ \bigcup_{i \in I}F_i := \bigvee_{f \in F} f \circ e_{\mathsmaller{\bigcup}}^{\Lambda(X)}, \]
\[ \big(f \circ e_{\mathsmaller{\bigcup}}^{\Lambda(X)}\big) (i, x) := f(\C E_i^X(x)); \ \ \ \ (i, x) 
\in \bigcup_{i \in I}\lambda_0(i). \]
The \textit{canonical Bishop topology on}\index{canonical Bishop topology on the intersection of a spectrum of subsets}
$\bigcap_{i \in I}\lambda_0(i)$ is the relative topology of $\C F$ to it i.e., 
\[ \bigcap_{i \in I}\C F_i := \bigg(\bigcap_{i \in I}\lambda_0(i), \bigcap_{i \in I}F_i\bigg), \]
\[ \bigcap_{i \in I}F_i := \bigvee_{f \in F} f \circ e_{\mathsmaller{\bigcap}}^{\Lambda(X)}, \] 
\[ \big(f \circ e_{\mathsmaller{\bigcap}}^{\Lambda(X)} \big) (\Phi) := f(\C E_{i_0}^X(\Phi_{i_0})) ; 
\ \ \ \ \Phi \in \bigcap_{i \in I}\lambda_0(i). \]
\end{definition}

Next follows the continuous-analogue to Proposition~\ref{prp: internalmap1}, using 
repeatedly the $\bigvee$-lifting of morphisms and the extensionality of a Bishop topology.

\begin{proposition}\label{prp: continternalmap1}
Let $S(\Lambda(X)) := (\lambda_0, \C E^X, \lambda_1, F, \phi_0^{\Lambda(X)}, \phi_1^{\Lambda(X)}) \in \Spec_F(I, X)$ 
with Bishop subspaces $(\C F_i)_{i \in I}$ and Bishop morphisms $(\C E_i^X)_{i \in I}$, 
$S(M(X)) := (\mu_0, \C Z^X, \mu_1, F, \phi_0^{M(X)}, \phi_1^{M(X)}) \in \Spec_F(I, X)$
with Bishop subspaces $(\C G_i)_{i \in I}$ and Bishop morphisms $(\C Z_i^X)_{i \in I}$,
and $\Psi \colon S(\Lambda(X)) \To S(M(X))$.\\[1mm]
\normalfont (i)
\itshape  $e_i^{\Lambda(X)} \in \Mor(\C F_i, \bigcup_{i \in I}\C F_i)$, for every $i \in I$.\\[1mm]
\normalfont (ii)
\itshape $\bigcup \Psi \in \Mor(\bigcup_{i \in I}\C F_i, \bigcup_{i \in I}\C G_i)$.\\[1mm]
\normalfont (iii)
\itshape $\pi_i^{\Lambda(X)} \in \Mor(\bigcap_{i \in I}\C F_i, \C F_i)$, for every $i \in I$.\\[1mm]
\normalfont (iv)
\itshape  $\bigcap \Psi \in \Mor(\bigcap_{i \in I}\C F_i, \bigcap_{i \in I}\C G_i)$.
\end{proposition}

\begin{proof}
(i) $e_i^{\Lambda(X)} \in \Mor(\C F_i, \bigcup_{i \in I}\C F_i) \TOT \forall_{f \in F}\big(\big(f \circ 
e_{\mathsmaller{\bigcup}}^{\Lambda(X)}\big) \circ  e_i^{\Lambda(X)} \in F_i\big)$. If $f \in F$, then  
$\big(f \circ e_{\mathsmaller{\bigcup}}^{\Lambda(X)}\big) \circ  e_i^{\Lambda(X)} := f \circ \C E_i^X \in F_i$,
since, for every $x \in \lambda_0(i)$, we have that
$\big[ \big(f \circ e_{\mathsmaller{\bigcup}}^{\Lambda(X)}\big) \circ  e_i^{\Lambda(X)}\big](x) := 
\big( f \circ e_{\mathsmaller{\bigcup}}^{\Lambda(X)}\big) (i,x) := \big(f \circ \C E_i^X)(x)$.\\
(ii) $\bigcup \Psi \in \Mor(\bigcup_{i \in I}\C F_i, \bigcup_{i \in I}\C G_i) \TOT 
\forall_{f \in F}\big(\big(f \circ e_{\mathsmaller{\bigcup}}^{M(X)}\big) \circ \bigcup \Psi \in 
\bigcup_{i \in I}F_i\big)$, and
$\big(f \circ e_{\mathsmaller{\bigcup}}^{M(X)}\big) \circ \bigcup \Psi = f \circ
e_{\mathsmaller{\bigcup}}^{\Lambda(X)}  \in \bigcup_{i \in I}F_i$, as
\[  \big[\big(f \circ e_{\mathsmaller{\bigcup}}^{M(X)}\big) \circ \bigcup \Psi\big](i, x) := 
\big( f \circ e_{\mathsmaller{\bigcup}}^{M(X)}\big) \big((i, \Psi_i(x))\big)
:= f \big(\C Z^X_i(\Psi_i(x))\big) \]
\[ = \big((f \circ \C E^X_i) \circ \Psi_i\big)(x)
 = (f \circ \C E_i^X)(x)\\
 := \big( f \circ e_{\mathsmaller{\bigcup}}^{\Lambda(X)}\big)(i, x). \]
(iii) $\pi_i^{\Lambda(X)} \in \Mor(\bigcap_{i \in I}\C F_i, \C F_i) \TOT 
\forall_{f \in F}\big((f \circ \C E_i^X) \circ \pi_i^{\Lambda(X}) \in \bigcap_{i \in I}F_i \big)$, and
$(f \circ \C E_i^X) \circ \pi_i^{\Lambda(X)} = f \circ e_{\mathsmaller{\bigcap}}^{\Lambda(X)} \in \bigcap_{i \in I}F_i$,
as
$\big[(f \circ \C E_i^X) \circ \pi_i^{\Lambda(X)}\big](\Phi) := f(\C E_i^X(\Phi_i)) = f(\C E^X_{i_0}(\Phi_{i_0})) :=
\big( f \circ e_{\mathsmaller{\bigcap}}^{\Lambda(X)} \big) (\Phi)$.\\
(iv) $\bigcap \Psi \in \Mor(\bigcap_{i \in I}\C F_i, \bigcap_{i \in I}\C G_i) \TOT 
\forall_{f \in F}\big(\big( f \circ e_{\mathsmaller{\bigcap}}^{M(X)} \big) \circ \bigcap \Psi \in 
\bigcap_{i \in I}F_i\big)$, and
$\big( f \circ e_{\mathsmaller{\bigcap}}^{M(X)} \big) \circ \bigcap \Psi = 
f \circ e_{\mathsmaller{\bigcup}}^{\Lambda(X)} \in \bigcup_{i \in I}F_i$, as
$
 \big[\big(f \circ e_{\mathsmaller{\bigcap}}^{M(X)}\big) \circ \bigcap \Psi\big](\Phi) := 
 \big( f \circ e_{\mathsmaller{\bigcap}}^{M(X)} \big)\big(\C Z^X_{i_0}(\Psi_{i_0}(\Phi_{i_0}))\big)
 = \big( f \circ e_{\mathsmaller{\bigcap}}^{M(X)} \big)\big(\C E^X_{i_0}(\Phi_{i_0})\big)
 := f \big(\C E^X_{i_0}(\Phi_{i_0})\big)
 := \big( \big( f \circ e_{\mathsmaller{\bigcap}}^{M(X)} \big)(\Phi).
$
\end{proof}

The notions mentioned in the next proposition were defined in Proposition~\ref{prp: geninternalmap1}.

\begin{proposition}\label{prp: gencontinternalmap1}
Let $\C F := (X, F)$, $(\C G := (Y, G)$ be Bishop spaces, $h : X \to Y \in \Mor(\C F, \C G)$,  
$S(\Lambda(X)) := (\lambda_0, \C E^X, \lambda_1, F ; \phi_0^{\Lambda(X)}, \phi_1^{\Lambda(X)}) \in \Spec_F(I, X)$ 
with Bishop subspaces $(\C F_i)_{i \in I}$ and Bishop morphisms $(\C E^X_i)_{i \in I}$,
$S(M(Y)) := (\mu_0, \C Z^Y, \mu_1, G, \phi_0^{M(Y)}, \phi_1^{M(Y)}) \in \Spec_G(I, Y)$
with Bishop subspaces $(\C G_i)_{i \in I}$ and Bishop morphisms $(\C Z^Y_i)_{i \in I}$,
and $\Psi : \Lambda(X) \stackrel{ h \ } \To M(Y)$.\\[1mm]
\normalfont (i)
\itshape  $\Psi$ is continuous i.e., $\Psi_i \in \Mor(\C F_i, \C G_i)$, for every $i \in I$.\\[1mm]
\normalfont (ii)
\itshape $\bigcup_h \Psi \in \Mor(\bigcup_{i \in I}\C F_i, \bigcup_{i \in I}\C G_i)$.\\[1mm]
\normalfont (iii)
\itshape  $\bigcap_h \Psi \in \Mor(\bigcap_{i \in I}\C F_i, \bigcap_{i \in I}\C G_i)$.
\end{proposition}

\begin{proof}
(i) $\Psi_i \in \Mor(\C F_i, \C G_i) \TOT \forall_{g \in G}\big((g \circ E_i) \circ \Psi_i \in F_i$,
and if $g \in G$, then
$(g \circ \C Z^Y_i) \circ \Psi_i := g \circ (\C Z^Y_i \circ \Psi_i) = g \circ (h \circ \C E_i^X) := 
(g \circ h) \circ \C E_i^X
\in F_i$, as $h \in \Mor(\C F, \C G)$, and hence $g \circ h \in F$.\\
(ii) and (iii) Working as in the proof of the Proposition~\ref{prp: continternalmap1}(ii) and (iv), 
we get
$ \big( g \circ e_{\mathsmaller{\bigcup}}^{M(X)}\big) \circ \bigcup_h \Psi = 
(g \circ h) \circ e_{\mathsmaller{\bigcup}}^{\Lambda(X)}$ and 
$\big(f \circ e_{\mathsmaller{\bigcap}}^{M(X)}\big) 
 \circ \bigcap_h \Psi = 
 (g \circ h) \circ e_{\mathsmaller{\bigcap}}^{\Lambda(X)}$,
for every $g \in G$.
\end{proof}

\section{Direct spectra of Bishop subspaces}
\label{sec: dirspectrasub}

\begin{definition}\label{def: dirintspectrum}
If $\Lambda^{\lt}(X) := (\lambda_0, \C E^X, \lambda_1^{\lt}) \in \Fam(I, \lt, X)$, 
a family of Bishop subspaces of the Bishop space $\C F := (X, F)$
associated to $\Lambda^{\lt}(X)$\index{family of Bishop subspaces associated to a direct family of subsets} 
is a pair $\Phi^{\Lambda^{\lt}(X)} := \big(\phi_0^{\Lambda^{\lt}(X)}, \phi_1^{\Lambda^{\lt}(X)}\big)$, where 
$\phi_0^{\Lambda^{\lt}(X)} \colon I \sto \D V_0$ and $\phi_1^{\Lambda^{\lt}(X)} : \bigcurlywedge_{(i,j) \in \ltI}
\D F\big(\phi_0^{\Lambda^{\lt}(X)}(j), \phi_0^{\Lambda^{\lt}(X)}(i)\big)$ such that the following conditions hold:\\[1mm]
\normalfont (i) 
\itshape $\phi_0^{\Lambda^{\lt}(X)}(i) := F_i := F_{|\lambda_0(i)} := \bigvee _{f \in F} f \circ \C E_i^X$, 
for every $i \in I$.\\[1mm]
\normalfont (ii)
\itshape $\phi_1^{\Lambda^{\lt}(X)} (i,j) := \big(\lambda^{\lt}_{ij}\big)^*$, for every $(i, j) \in D^{\lt}(I)$. \\[1mm]
We call the structure
$S_F(\Lambda^{\lt}(X)) := (\lambda_0, \C E^X, \lambda_1^{\lt}, F, 
\phi_0^{\Lambda^{\lt}(X)}, \phi_1^{\Lambda^{\lt}(X)})$ 
a $($covariant$)$\textit{direct spectrum of subspaces of $\C F$} over $I$,\index{direct spectrum of subspaces} or 
an $(I, \lt_I)$-spectrum of subspaces of $\C F$ with Bishop subspaces $(\C F_i)_{i \in I}$ and 
Bishop morphisms\index{$S_F(\Lambda^{\lt}(X))$}
$(\C E_i^X)_{i \in I}$. If $S_F(M^{\lt}(X)) := (\mu_0, \C Z^X, \mu_1^{\lt}, F, \phi_0^{M^{\lt}(X)}, \phi_1^{M^{\lt}(X)})$
is an $(I, \lt_I)$-spectrum of subspaces of $\C F$ with Bishop subspaces $(\C G_i)_{i \in I}$ and 
Bishop morphisms $(\C Z^X_i)_{i \in I}$, a \textit{subspaces direct spectrum-map}\index{subspaces direct spectrum-map} 
$\Psi$ from $S_F(\Lambda^{\lt}(X))$ to $S_F(M^{\lt}(X))$, in symbols 
$\Psi \colon S_F(\Lambda^{\lt}(X)) \To S_F(M^{\lt}(X))$, is a direct family of
subsets-map $\Psi \colon \Lambda^{\lt}(X) \To M^{\lt}(X)$ $($see Definition~\ref{def: dirsubfammap}$)$. 
If $F$ is clear from the context, we may omit he symbol $F$ as a subscript in the above notations.
A contravariant  direct spectrum $S_F(\Lambda^{\mt}(X)) := (\lambda_0, \C E^X, \lambda_1^{\mt}, 
F, \phi_0^{\Lambda^{\mt}(X)}, \phi_1^{\Lambda^{\mt}(X)})$ of subspaces of $\C F$ over $(I, \lt_I)$
and a subspaces contravariant direct 
spectrum-map\index{contravariant direct spectrum of subspaces}\index{a subspaces contravariant direct spectrum-map} 
are defined similarly.
\end{definition}

\begin{remark}\label{rem: dirintspectrum1}
Let $S_F(\Lambda^{\lt}(X)) := (\lambda_0, \C E^X, \lambda_1^{\lt},
F, \phi_0^{\Lambda^{\lt}(X)}, \phi_1^{\Lambda^{\lt}(X)})$ be an
$(I, \lt_I)$-spectrum of subspaces of $\C F := (X, F)$ with Bishop subspaces $(\C F_i)_{i \in I}$ 
and Bishop morphisms $(\C E_i^X)_{i \in I}$, and $S_F(M^{\lt}(X)) := (\mu_0, \C Z^X, \mu^{\lt}_1, F, 
\phi_0^{M^{\lt}(X)}, \phi_1^{M^{\lt}(X)})$
an $(I, \lt_I)$-spectrum of subspaces of $\C F$ with Bishop subspaces $(\C G_i)_{i \in I}$ and Bishop 
morphisms $(\C Z^X_i)_{i \in I}$.\\[1mm]
\normalfont (i)
\itshape $S(\Lambda^{\lt}) := (\lambda_0, \lambda_1^{\lt}, \phi_0^{\Lambda^{\lt}(X)}, \phi_1^{\Lambda^{\lt}(X)})$ is 
an $(I, \lt_I)$-spectrum with Bishop spaces $(\C F_i)_{i \in I}$ and Bishop morphisms 
$(\lambda^{\lt}_{ij})_{(i,j) \in D^{\lt}(I)}$.\\[1mm]
\normalfont (ii)
\itshape If $\Psi \colon S(\Lambda^{\lt}(X)) \To S(M^{\lt}(X))$, then $\Psi$ is continuous.
\end{remark}

\begin{proof}
We proceed as in the proof of Remark~\ref{rem: intspectrum1}.
\end{proof}

\begin{definition}\label{def: dirspecsub}
Let $\Spec_F(I, \lt, X)$ be the totality of covariant direct spectra of subspaces of the 
Bishop space $\C F := (X, F)$ and let
 $\Spec_F(I, \mt, X)$ be the totality of contravariant direct spectra of subspaces of $\C F$ over $(I, \lt_I$,
 equipped with the equality of $\Fam(I, \lt, X)$ and $\Fam(I, \mt, X)$, respectively.
\end{definition}


\index{inverse limit of a contravariant direct spectrum of subspaces} \index{}
\begin{definition}\label{def: invlimdirint} If 
$S(\Lambda^{\mt}(X)) := (\lambda_0, \C E, \lambda_1^{\mt}, F, \phi_0^{\Lambda^{\mt}(X)}, 
\phi_1^{\Lambda^{\mt}(X)})$ is contravariant direct spectrum of subspaces of the Bishop space 
$\C F := (X, F)$ over $(I, \lt)$ with Bishop subspaces $(\C F_i)_{i \in I}$ and Bishop morphisms 
$(\C E_i)_{i \in I}$, its \textit{inverse limit} is the following
Bishop space
\[ \underset{\ot} \Lim S(\Lambda^{\mt}(X)):= \underset{\ot_X} \Lim \C F_i := \bigg(\bigcap_{i \in I}\lambda_0(i), 
\bigvee_{f \in F} f \circ e_{\mathsmaller{\bigcap}}^{\Lambda^{\mt}(X)}\bigg). \]
\end{definition}

Next we show the universal property of the inverse limit for $\underset{\ot_X} \Lim \C F_i $

\begin{proposition}\label{prp: universalinverse}
 If 
$S(\Lambda^{\mt}(X)) := \big(\lambda_0, \C E, \lambda_1^{\mt}, F, \phi_0^{\Lambda^{\mt}(X)}, 
\phi_1^{\Lambda^{\mt}(X)}\big) \in \Spec(I, \mt, X)$  
with Bishop subspaces $(\C F_i)_{i \in I}$ and Bishop morphisms $(\C E_i^X)_{i \in I}$, 
its inverse limit $\underset{\ot} \Lim \C F_i$ satisfies the universal property of inverse limits
i.e., if $i \lt_ j$, the following left diagram commutes
\begin{center}
\begin{tikzpicture}

\node (E) at (0,0) {$\bigcap_{i \in I}\lambda_0(i)$};
\node[below=of E] (F) {};
\node [right=of F] (B) {$\lambda_0(j)$};
\node [left=of F] (C) {$\lambda_0(i)$};
\node [right=of B] (K) {$\lambda_0(i)$};
\node[right=of K] (L) {};
\node [right=of L] (M) {$\lambda_0(j),$};
\node [above=of L] (N) {$Y$};

\draw[->] (E)--(B) node [midway,right] {$ \ \  \pi_j^{\Lambda^{\mt}(X)}$};
\draw[->] (E)--(C) node [midway,left] {$\pi_i^{\Lambda^{\mt}(X)} \ $};
\draw[left hook->] (B)--(C) node [midway,below] {$\lambda_{ji}^{\mt}$};

\draw[->] (N)--(M) node [midway,right] {$ \  \varpi_j$};
\draw[->] (N)--(K) node [midway,left] {$\varpi_i \ $};
\draw[left hook->] (M)--(K) node [midway,below] {$\lambda_{ji}^{\mt}$};

\end{tikzpicture}
\end{center}
and for every Bishop space $\C G := (Y, G)$ and a family $(\varpi_i)_{i \in I}$, where
$\varpi_i \in \Mor(\C G, \C F_i)$, for every $i \in I$, such that the above right diagram commutes, there is a 
unique Bishop morphism $h : Y \to \bigcap_{i \in I}\lambda_0(i)$ such that the following diagrams commute
\begin{center}
\begin{tikzpicture}

\node (E) at (0,0) {$Y$};
\node[below=of E] (F) {};
\node [right=of F] (B) {$\lambda_0(j)$,};
\node [left=of F] (C) {$\lambda_0(i)$};
\node[below=of F] (G) {$\bigcap_{i \in I}\lambda_0(i)$};

\draw[->] (E)--(B) node [midway,right] {$ \ \  \varpi_j$};
\draw[->] (E)--(C) node [midway,left] {$\varpi_i \ $};
\draw[left hook->] (B)--(C) node [midway,below] {$\lambda_{ji}^{\mt} \ \ \ \ \ \  $};
\draw[->] (G)--(B) node [midway,right] {$ \ \  \pi_j^{\Lambda^{\mt}(X)}$};
\draw[->] (G)--(C) node [midway,left] {$\pi_i^{\Lambda^{\mt}(X)} \ $};
\draw[dashed,->] (E)--(G) node [midway,near start] {$\ \ \ h$};

\end{tikzpicture}
\end{center}
\end{proposition}

\begin{proof}
For the commutativity of the first diagram, we have that if $\Phi \in \bigcap_{i \in I}\lambda_0(i)$, then
$\pi_i^{\Lambda^{\mt}(X)}(\Phi) := \Phi_i$, and $\lambda_{ji}^{\mt}\big(\pi_j^{\Lambda^{\mt}}(\Phi)\big) :=
\lambda_{ji}^{\mt}(\Phi_j)$, and since $\C E_j^X = \C E_i^X \circ \lambda_{ji}^{\mt}$, we have that  
$\C E_j^X (\Phi_j) = \C E_i^X \big(\lambda_{ji}^{\mt}(\Phi_j)\big)$, hence by the definition of
$\bigcap_{i \in I}\lambda_0(i)$ we get $\C E_i^X (\Phi_i) = \C E_j^X(\Phi_j) = 
\C E_i^X \big(\lambda_{ji}^{\mt}(\Phi_j)\big)$, and since $\C E_i^X$ is an embedding we get $\Phi_i = 
\lambda_{ji}^{\mt}(\Phi_j) =: \lambda_{ji}^{\mt}\big(\pi_j^{\Lambda^{\mt}(X)}(\Phi)\big)$.
Let a Bishop space $\C G := (Y, G)$ and a family of Bishop morphisms $(\varpi_i)_{i \in I}$, where
$\varpi_i : Y \to \lambda_0(i)$, for every $i \in I$, such that the above right diagram commutes. Let also the 
operation $h : Y \sto \bigcap_{i \in I}\lambda_0(i)$, defined by the rule $y \mapsto h(y)$, where
\[ h(y) : \bigcurlywedge_{i \in I}\lambda_0(i), \ \ \ \ h(y)_i := \varpi_i(y); \ \ \ \ i \in I. \]
To show that $h(y) \in \bigcap_{i \in I}\lambda_0(i)$ we need to show that 
\[ \C E_i^X\big(h(y)_i\big) =_X \C E_{i{'}}^X\big(h(y)_{i{'}}\big) \TOT \C E_i^X(\varpi_i(y)) =_X 
\C E_{i{'}}^X(\varpi_{i{'}}(y)), \]
for every $i, i{'} \in I$. Since $(I, \lt_I)$ is directed, there is $k \in I$ such that $i \lt_I k$ and $i{'} \lt_I k$,
hence
\[ \C E_i^X(\varpi_i(y)) =_X \C E_i^X \big(\lambda_{ki}^{\mt}(\varpi_k(y))\big) =_X \C E_k^X \big(\varpi_k(y)\big)
=_X \C E_{i{'}}^X \big(\lambda_{ki{'}}^{\mt}(\varpi_k(y))\big) =_X \C E_{i{'}}^X(\varpi_{i{'}}(y)). \]
It is immediate to show that $h$ is a function. Finally, we show that
$ h \in \Mor(\C G,  \underset{\ot_X} \Lim \C F_i) \TOT \forall_{f \in F}\big(\big(f \circ
e_{\mathsmaller{\bigcap}}^{\Lambda^{\mt}(X)}\big) \circ h \in G\big). $
If $y \in Y$, then 
\[ \big[\big(f \circ e_{\mathsmaller{\bigcap}}^{\Lambda^{\mt}(X)}\big) \circ h\big](y) :=
f\big(\C E_{i_0}^X(\varpi_{i_0}(y))\big) := [(f \circ \C E_{i_0}^X) \circ \varpi_{i_0}](y), \]
hence 
$\big(f \circ e_{\mathsmaller{\bigcap}}^{\Lambda^{\mt}(X)}\big) \circ h := (f \circ \C E_{i_0}^X) 
\circ \varpi_{i_0} \in G$, as by our hypothesis $\varpi_{i_0} \in \Mor(\C G, \C F_{i_0})$. 
\end{proof}

\section{Notes}
\label{sec: notes5}

\begin{note}\label{not: historybs}
\normalfont 
The theory of Bishop spaces, that was only sketched by Bishop in~\cite{Bi67}, and revived 
by Bridges in~\cite{Br12}, and Ishihara in~\cite{Is13}, was developed by the author in~\cite{Pe15}-\cite{Pe19d}
and~\cite{Pe20b}-\cite{Pe20d}.
Since inductive definitions with rules of countably many premises are used, for the study of Bishop spaces
we work within $\BST^*$, which is $\BST$ extended 
with such inductive definitions. A formal system for $\BISH$ extended with such definitions  
is Myhill's formal system $\CST^*$
with dependent choice, where $\CST^*$ is Myhill's extension of his formal system of constructive set theory 
$\CST$ with inductive definitions (see~\cite{My75}). A variation of $\CST^*$ is Aczel's system CZF together
with a very weak version of Aczel's regular extension axiom\index{regular extension axiom} (REA), 
to accommodate these inductive definitions (see~\cite{AR10}).

\end{note}

\begin{note}\label{not: onbs}
\normalfont 
In contrast to topological spaces, in the theory of Bishop spaces continuity of functions is  
an a priori notion, while the concept of an open set comes a posteriori, through the neighbourhood 
space induced by a Bishop topology. The theory of Bishop spaces can be seen as an abstract and
constructive approach to the theory of the ring $C(X)$ of continuous functions of a topological 
space $(X, \C T)$ (see~\cite{GJ60} for a classical treatment of this subject). 
\end{note}

\begin{note}\label{not: onlimitstheory}
\normalfont 
The results on the direct and inverse limits of direct spectra of Bishop spaces are the constructive 
analogue to the classical theory of direct and inverse limits of (spectra of) topological spaces, as 
this is developed e.g., in the Appendix of~\cite{Du66}. As in the case of the classic textbook of Dugundji,
we avoid here possible, purely categorical arguments in our proofs. One of the advantages of working with a proof-relevant 
definition of a cofinal subset is that the proof of the cofinality theorem~\ref{thm: cofinal2} is choice-free.

\end{note}

\begin{note}\label{not: setrelavantspectra}
\normalfont 
The notion of a spectrum of Bishop spaces can be generalised by considering a family of Bishop spaces
associated to a set-relevant family of sets over some set $I$. In this case, all transport maps $\lambda_{ij}^m$
are taken to be Bishop morphisms. The direct versions of set-relevant spectra of Bishop spaces can be defined, 
and their theory can be developed in complete analogy to the theory of direct spectra of Bishop spaces, 
as in the case of generalised direct spectra of topological spaces (see~\cite{Du66}, p.~426).

\end{note}

\begin{note}\label{not: onuniversal}
\normalfont 
The formulation of the universal properties of the various limits of spectra of Bishop spaces included
here is impredicative, as it requires quantification over the class of Bishop spaces.  A predicative 
formulation of a universal property can be given, if one is restricted to a given set-indexed family of Bishop spaces.
\end{note}

\begin{note}\label{not: ondirectlimitofsubspaces}
\normalfont 
The study of the direct limit of a spectrum of Bishop subspaces is postponed for future work. The natural candidate $\bigcup_{i \in I}\lambda_0(i)$, equipped with the relative topology, ``almost'' satisfies the universal property of the direct limit.

\end{note}

\chapter{Families of subsets in measure theory}
\label{chapter: measure}

We study the Borel and Baire sets within Bishop spaces as a constructive counterpart to the study 
of Borel and Baire algebras within topological spaces. As we use the inductively defined least Bishop 
topology, and as the Borel and Baire sets over a family of $F$-complemented subsets are defined inductively, 
we work within the extension $\BISH^*$ of $\BISH$ with inductive definitions with rules of countably many premises. 
In contrast to the classical theory, we show that the Borel and the Baire sets of a Bishop space coincide. 
Our reformulation within $\BST$ of the Bishop-Cheng definition of a measure space and of an integration space,
based on the notions of families of complemented subsets and of families of partial functions, facilitates 
a predicative reconstruction of the originally impredicative Bishop-Cheng measure theory.

\section{The Borel sets of a Bishop space}
\label{sec: borel}

The Borel sets of a topological space $(X, \C T)$ is the least set of subsets of $X$ that includes the 
open (or, equivalently the closed) sets in $X$ and it is closed under countable unions, countable 
intersections and relative complements. The Borel sets of a Bishop space $(X, F)$ is the least set
of complemented subsets of $X$ that includes the \textit{basic} $F$-complemented  subsets of $X$ that
are generated by $F$, and it is closed under countable unions and countable intersections. As 
the Borel sets of $(X, F)$ are complemented subsets, it is not a coincidence that their closure under 
complements is provable. 
\textit{In the next two sections} $F$ denotes a Bishop topology on a set $X$ and $G$ a Bishop topology
on a set $Y$. For simplicity, we denote the constant function on $X$ with value $a \in \Real$ also by $a$, 
and we may write equalities between elements of $\C P^{\Disj_F}(X)$ and equalities between elements of $F$
without denoting the corresponding subscripts.

\begin{definition}\label{def: Fcomplemented}
If $a, b \in \Real$, let\index{$a \neq_{\Real} b$} $a \neq_{\Real} b : \TOT |a - b| > 0 \TOT a > b \vee b < a$.
For simplicity we may write $a \neq b$, instead of $a \neq_{\Real} b$. The inequality $x \neq_X^F y$ on $X$
generated by $F$ is defined by\index{$f \colon x \neq_X^F y$}
\[ x \neq^F_X y : \TOT \exists_{f \in F}\big(f \colon x \neq_X^F y\big), \ \ \ \ \mbox{where} \ \ \ \  
f \colon x \neq_X^F y :\TOT f(x) \neq_{\Real} f(y). \]
A complemented subset $\B A$ of $X$ with respect to $\neq_X^F$ is called an $F$-complemented subset of
$X$\index{$F$-complemented subset}, and their totality is denoted by $\C P^{\Disj_F} (X)$\index{$\C P^{\Disj_F}(X)$}. An 
$F$-complemented subset $\B A$ of $X$ is uniformly $F$-complemented\index{uniformly $F$-complemented subset}, if 
\[ \exists_{f \in F}\big(f \colon A^1 \Disj_F A^0 \big), \ \ \ \ \mbox{where} \ \ \ \ f 
\colon A^1 \Disj_F A^0 :\TOT \forall_{x \in A^1}\forall_{y \in A^0}\big( f \colon x \neq_X^F y\big), \]
and $\B A$\index{$f \colon A^1 \Disj_F A^0$} is strongly $F$-complemented\index{strongly $F$-complemented},
if there is $f \in F$ such that $f \colon A^1 \Disj_F A^0$, $f(x) = 1$, for every $x \in A^1$, and $f(y) = 0$,
for every $y \in A^0$. 
\end{definition}

\begin{remark}\label{rem: compl2}
If $h \in \Mor(\C F, \C G)$ and $\B A \in \C P^{\Disj_G} (Y)$, then 
$h^{-1}(\B A) 
\in \C P^{\Disj_F} (X)$.
\end{remark}

\begin{proof}
Let $x \in h^{-1}(A^1)$ and $y \in h^{-1}(A^0)$ i.e., $h(x) \in A^1$ and $h(y) \in A^0$. Let $g \in G$ 
such that $g(h(x)) \neq g(h(y))$. Hence, $g \circ h \in F$ and $(g \circ h)(x) \neq (g \circ h)(y)$.
\end{proof}

\begin{definition}\label{def: ofam}
We denote by $\Fam(I, F, \B X)$\index{$\Fam(I, F, \B X)$} and $\Set(I, F, \B X)$ the sets of families 
and sets of $F$-complemented subsets of $X$, respectively. Let $\B O_F(X) := \big(o_0^{1, F}, \C O^{1,X},
o_1^{1, F}, o_0^{0,F}, \C O^{0,X}, o_1^{0,F}\big) \in \Fam(F, F, \B X)$\index{$\B O_F(X)$} be the family 
of basic open $F$-complemented subsets\index{family of basic open $F$-complemented subsets} of $X$, where 
\[ \B o_F (f) := \big(o_0^{1,F}(f), o_0^{0,F}(f)\big) := \big([f > 0], [f \leq 0]\big), \]
\[ [f > 0] := \{x \in X \mid f(x) > 0\}, \ \ \ \ [f \leq 0] := \{x \in X \mid f(x) \leq 0\},\]
and, as $[f > 0$, $[f \leq 0]$ are extensional subsets of $X$, the dependent operations $\C O^{1,X}, 
\C O^{0,X}, o_1^{1, F}$, and $o_1^{0, F}$ are defined by the identity map-rule. If $F$ is clear from 
the context, we may write $\B O_F(X) := \big(o_0^{1}, \C O^{1,X}, o_1^{1}, o_0^{0}, \C O^{0,X}, o_1^{0}\big)$.
\end{definition}

Clearly, $f \colon o_0^{1}(f) \Disj_F o_0^{0}(f)$, for every $f \in F$. Recall that a sequence of 
$F$-complemented subsets of $X$ is a structure $\B B(X) := \big(\beta_0^1, \C B^{1,X}, \beta_1^1, 
\beta_0^0, \C B^{0,X}, \beta_1^0\big) \in \Fam(\Nat^+, F, \B X)$, where $\B \beta(n) := \big(\beta_0^1(n), 
\beta_0^0(n)\big) \in \C P^{\Disj_F}(X)$, and $\beta^1_{nn} \colon \beta_0^1(n) \to \beta_0^1(n)$ and 
$\beta^0_{nn} \colon \beta_0^0(n) \to \beta_0^0(n)$ are given by $\beta^1_{nn} := \id_{\beta_0^1(n)}$ and 
$\beta^0_{nn} := \id_{\beta_0^0(n)}$, respectively, for every $n \in \Nat^+$. We also write 
$\bigcup_{n = 1}^{\infty} \B \beta_0(n)$ and $\bigcap_{n = 1}^{\infty} \B \beta_0(n)$, instead of 
$\bigcup_{n \in \Nat^+} \B \beta_0(n)$ and $\bigcap_{n \in \Nat^+}\B \beta_0(n)$, respectively. A family
$\B A(X) := \big(\alpha_0^1, \C A^{1,X}, \alpha_1^1, \alpha_0^0, \C A^{0,X}, \alpha_1^0\big) \in \Fam(\D 1, F, \B X)$ 
is defined similarly.

\begin{definition}\label{def: borel}
If $\B \Lambda(X) := \big(\lambda_0^1, \C E^{1,X}, \lambda_1^1, \lambda_0^0, \C E^{0,X}, \lambda_1^0\big) 
\in \Fam(I, F, \B X)$,  the set $\Borel(\B \Lambda(X))$ of Borel sets generated by
$\B \Lambda(X)$\index{Borel sets of $\B \Lambda(X)$}\index{$\Borel(\B \Lambda(X))$}
 is defined inductively by the following rules:
\[ (\Borel_1) \ \  \ \ \ \ \ \ \ \ \ \ \ \ \ \ \ \ \ \ \ \ \ \ \ \ \ \ \ \ \ \ \ \ \  \ \frac{i \in I}{\B \lambda_0 (i) \in \Borel(\B \Lambda(X))}, \ \ \ \ \ \ \ \ \ \ \ \ \ \ \ \ \ \ \ \ \ \ \ \ \ \ \ \ \ \ \ \ \ \ \ \ \ \ \ \ \ \ \  \]
\[ (\Borel_2) \ \ \ \ \ \ \ \ \  \frac{ \B \beta_0(1) \in \Borel(\B \Lambda(X)),  \ \B \beta_0(2) \in \Borel(\B \Lambda(X)), \ldots}{\bigcup_{n = 1}^{\infty} \B \beta_0(n) \in \Borel(\B \Lambda(X)) \ \ \& \ \ \bigcap_{n = 1}^{\infty} \B \beta_0(n) \in \Borel(\B \Lambda(X))}\mathsmaller{\mathsmaller{\B B(X) \in \Fam(\Nat^+, F, \B X)}}, \ \ \ \ \ \ \ \ \ \ \ \ \]
\[ (\Borel_3) \ \ \ \ \ \ \ \ \ \ \ \ \ \ \ \ \ \  \frac{\B B \in \Borel(\B \Lambda(X)), \ \B \alpha_0(0) =_{\C P^{\Disj_F}(X)} \B B}{ \B \alpha_0(0) \in \Borel(\B \Lambda(X))}\mathsmaller{\mathsmaller{\B A(X) \in \Fam(\D 1, F, \B X)}}. \ \ \ \ \ \ \ \ \ \ \ \ \ \ \ \ \ \ \  \]
The corresponding induction principle $\Ind_{\Borel(\B \Lambda(X))}$ is the formula 
\[ \forall_{i \in I}\big(P(\B \lambda_0 (i))\big) \  \& \  \forall_{\B B(X) \in \Fam(\Nat^+, F, \B X)}\bigg[\forall_{n \in \Nat^+}\big(\B \beta_0 (n) \in \Borel(\B \Lambda(X)) \ \&  \ P(\B \beta_0 (n))\big) \To \]
\[ P\bigg(\bigcup_{n = 1}^{\infty} \B \beta_0(n)\bigg) \ \& \ P\bigg(\bigcap_{n = 1}^{\infty} \B \beta_0(n)\bigg)\bigg] \ \& \ \]
\[ \forall_{\B A(X) \in \Fam(\D 1, F, \B X)}\forall_{\B B \in \Borel(\B \Lambda(X))}\bigg([P(\B B) \ \& \ \B \alpha_0(0) =_{\C P^{\Disj_F}(X)} \B B] \To P(\B \alpha_0(0))\bigg)  \]
\[ \To \forall_{\B B \in  \Borel(\B \Lambda(X))}\big(P(\B B)\big), \]
where $P$ is any bounded formula. 
Let
\[ \Borel(\C F) := \Borel(\B O_F(X)), \]
and we call its elements the Borel sets of $\C F$\index{Borel set of $\C F$}.
 
\end{definition}

In $\Ind_{\Borel(\B \Lambda(X))}$ we quantify over the sets $\Fam(\Nat^+, F, \B X)$ and $\Fam(\D 1, F, \B X)$,
avoiding quantification over $\C P^{\Disj_F}(X)$ in condition $(\Borel_3)$ and treating
$\C P^{\Disj_F}(X)$ as a set in $(\Baire_2)$.

\begin{proposition}\label{prp: borel1}
\normalfont (i)
\itshape $\B O_F(X)$ is not in $\Set(F, F, \B X)$.\\[1mm]
\normalfont (ii)
\itshape If $f, g \in F$, then $\B o (f) \cup \B o (g) = \B o (f \vee g)$.\\[1mm]
\normalfont (iii)
\itshape If $\B B \in \Borel(\C F)$, then $- \B B \in \Borel(\C F)$.\\[1mm]
\normalfont (iv)
\itshape There are Bishop space $(X, F)$ and $f \in F$ such that $\neg[- \B o (f) = \B o (g)]$, for every
$g \in F$. \\[1mm]
\normalfont (v)
\itshape $\B o(f) = \B o([f \vee 0] \wedge 1)$.
\end{proposition}

\begin{proof}
(i) If $f \in F$, then $\B o (f) = \B o (2f)$, but $\neg (f = 2 f)$.\\
(iii) This equality is implied from the following properties for reals
$a \vee b > 0 \TOT a > 0 \vee b > 0 $ and $ a \vee b \leq 0 \TOT a \leq 0 \wedge b \leq 0$.\\
(iv) If $a \in \Real$, then $a \leq 0 \TOT \forall_{n \geq 1}\big(a < \frac{1}{n}\big)$ and $a > 0 
\TOT \exists_{n \geq 1}\big(a \geq \frac{1}{n}\big)$, hence
\begin{align*}
- \B o (f) & := \big([f \leq 0], [f > 0]\big) \\
& = \bigg(\bigcap_{n = 1}^{\infty}\big[\big(\frac{1}{n} - f\big) > 0\big], 
\bigcup_{n = 1}^{\infty}\big[\big(\frac{1}{n} - f\big) \leq 0\big]\bigg)\\
& := \bigcap_{n = 1}^{\infty} \B o \big(\frac{1}{n} - f\big) \in \Borel(F).
\end{align*}
If $P(\B B) :\TOT  - \B B \in \Borel(\C F)$, the above equality proves the first step of the
corresponding induction on $\Borel(\C F)$. The rest of the inductive proof is straightforward.\\
(v) Let the Bishop space $(\Real, \BR)$. If we take $\B o (\id_{\Real}) := \big([x > 0], [x \leq 0]\big)$, and
if we suppose that $- \B o (\id_{\Real}) := \big([x \leq 0], [x > 0]\big) = \big([\phi > 0], [\phi \leq 0]\big)
=: \B o (\phi)$, for some $\phi \in \BR$, then $\phi(0) > 0$ and $\phi$ is not continuous at $0$, which 
contradicts the fact that $\phi$ is uniformly continuous, hence pointwise continuous, on $[-1, 1]$. \\
(vi) The proof is based on basic properties of $\Real$, like $a \wedge 1 = 0 \To a = 0$. 
\end{proof}

Since $\Borel(\C F)$ is closed under intersections and complements, if $\B A, \B B \in \Borel(\C F)$, 
then $\B A - \B B \in \Borel(\C F)$. Constructively, we cannot show, in general,  that $\B o(f) \cap 
\B o(g) = \B o(f \wedge g)$. If $f := \id_{\Real} \in \Bic(\Real)$ and $g := - \id_{\Real} \in \BR$, 
then $\B o(\id_{\Real}) \cap \B o(- \id_{\Real}) = \big([x > 0] \cap [x < 0], [x \leq 0] \cup [-x \leq 0]\big) 
= \big(\emptyset , [x \leq 0] \cup [x \geq 0]\big)$ Since $x \wedge (-x) = - |x|$, we get $\B o(\id_{\Real} \wedge (- \id_{\Real})) = \B o (- |x|) = \big(\emptyset, [|x| \geq 0]\big)$. The supposed equality implies that $|x| \geq 0 \TOT x \leq 0 \vee x \geq 0$. Since $|x| \geq0$ is always the case, we get $\forall_{x \in \Real}\big(x \leq 0 \vee x \geq 0\big),$ which implies LLPO (see~\cite{BV06}, p.~20). If one add the condition $|f| + |g| > 0$, then $\B o(f) \cap \B o(g) = \B o(f \wedge g)$ follows constructively. The condition $(\BS_4)$ in the definition of a Bishop space is crucial to the next proof.

\begin{proposition}\label{prp: inf}
If $(f_n)_{n = 1}^{\infty} \subseteq F$, then $f := \sum_{n = 1}^{\infty}(f_n \vee 0) \wedge 2^{-n} \in F$ and
$$\B o(f) = \bigcup_{n = 1}^{\infty}\B  o(f_n) = \bigg(\bigcup_{n = 1}^{\infty}[f_n > 0], 
\bigcup_{n = 1}^{\infty}[f_n \leq 0]\bigg).$$
\end{proposition}

\begin{proof}
The function $f$ is well-defined by the comparison test (see~\cite{BB85}, p.~32). If $g_n := (f_n \vee 0)
\wedge 2^{-n}$, for every $n \geq 1$, then 
$$\bigg|\sum_{n = 1}^{\infty}g_n - \sum_{n = 1}^N g_n\bigg| = \bigg|\sum_{n = N+1}^{\infty}g_n \bigg|
\leq \sum_{n = N+1}^{\infty}|g_n| \leq \sum_{n = N+1}^{\infty}\frac{1}{2^n} \stackrel{N} \longrightarrow 0,$$ 
the sequence of the partial sums $\sum_{n = 1}^N g_n \in F$ converges uniformly to $f$, hence by $\BS_4$ 
we get $f \in F$. Next we show that $[f > 0] \subseteq \bigcup_{n = 1}^{\infty}[f_n > 0]$. If 
$x \in X$ such that $f(x) > 0$, there is $N \geq 1$ such that $\sum_{n = 1}^N g_n(x) > 0$. By 
Proposition (2.16) in~\cite{BB85}, p.~26, there is $n \geq 1$ and $n \leq N$ with $g_n(x) > 0$,
hence $(f_n (x) \vee 0) \geq g_n(x) > 0$, which implies $f_n (x) > 0$. For the converse inclusion,
if $f_n (x) > 0$, for some $n \geq 1$, then $g_n (x) > 0$, hence $f(x) > 0$. To show $[f \leq 0]
\subseteq  \bigcup_{n = 1}^{\infty}[f_n \leq 0]$, let $x \in X$ such that $f(x) \leq 0$, and suppose
that $f_n (x) > 0$, for some $n \geq 1$. By the previous argument we get $f(x) > 0$, which
contradicts our hypothesis $f(x) \leq 0$. For the converse inclusion, let $f_n (x) \leq 0$, for
every $n \geq 1$, hence $f_n (x) \vee 0 = 0$ and $g_n (x) = 0$, for every $n \geq 1$. Consequently, $f(x) = 0$.
\end{proof}

\begin{proposition}\label{prp: borel2}
If $h \in \Mor(\C F, \C G)$ and $\B B \in \Borel(\C G)$, then $h^{-1}(\B B) \in \Borel(\C F)$.
\end{proposition}

\begin{proof}
By the definition of $h^{-1}(\B B)$, if $g \in G$, then 
\begin{align*}
h^{-1}(\B o_G(g)) & := h^{-1}\big([g > 0], [g \leq 0]\big)\\
& := \big(h^{-1}[g > 0], h^{-1}[g \leq 0]\big)\\
& = \big([(g \circ h) > 0], [(g \circ h) \leq 0]\big)\\
& := \B o_F(g \circ h) \in \Borel(\C F).
\end{align*}
If $P(\B B) := h^{-1}(\B B) \in \Borel(\C F)$, the above equality is the first step of the
corresponding inductive proof on $\Borel(\C G)$. The rest of the proof follows from the 
properties $h^{-1}\big(\bigcup_{n = 1}^{\infty}\B B_n\big) = \bigcup_{n = 1}^{\infty}h^{-1}(\B B_n)$ 
and $h^{-1}\big(\bigcap_{n = 1}^{\infty}\B B_n\big) = \bigcap_{n = 1}^{\infty}h^{-1}(\B B_n)$ of complemented subsets.
\end{proof}

\begin{definition}\label{def: borelB}
If $\Phi$ is an extensional subset of $F$ and if $\id_{\Phi}^F \colon \Phi \eto F$ is defined by the
identity map-rule, let
$\B O_{\Phi}(X) := \B O_F(X) \circ \id_{\Phi}^{F}$ be the $\id_{\Phi}^F$-subfamily of $\B O_F(X)$. We write  
$\B o_{\Phi} (f) := \B o_F (f)$, for every $f \in \Phi$, and let 
\[ \Borel(\Phi) := \Borel(\B O_{\Phi}(X)). \] 
\end{definition}

If $F_0$ is a subbase of $F$, then, $\Borel(F_0) \subseteq \Borel(F)$. More can be said on the relation
between $\Borel(\Phi)$ and $\Borel(\C F)$, when $\Phi$ is a base of $F$. 

\begin{proposition}\label{prp: base1}
Let $\Phi$ be a base of $F$.\\[1mm]
\normalfont (i)
\itshape If for every $f \in F$, $\B o_F(f) \in \Borel (\Phi)$, then $\Borel(\C F) = \Borel(\Phi)$.\\[1mm]
\normalfont (ii)
\itshape If for every $g \in \Phi$ and $f \in F$, $f \wedge g \in \Phi$, then 
$\Borel(\C F) = \Borel(\Phi)$.\\[1mm]
\normalfont (iii)
\itshape If for every $g \in B$ and every $n \geq 1$, $g - \frac{1}{n} \in \Phi$, then $\Borel(\C F) = \Borel(\Phi)$. 
\end{proposition}

\begin{proof}
(i) It follows by a straightforward induction on $\Borel(F)$.\\
(ii) and (iii) Let $f \in F$ and $(g_n)_{n = 1}^{\infty} \subseteq \Phi$ 
such that $\forall_{n \geq 1}\big(U(f, g_n, \frac{1}{n})\big)$. Then we have tha
\[ \B o_F (f) \subseteq \bigcup_{n = 1}^{\infty}\B o_{\Phi}(g_n) := \bigg(\bigcup_{n = 1}^{\infty}[g_n > 0],
\bigcap_{n = 1}^{\infty}[g_n \leq 0]\bigg) \]
i.e., $[f > 0] \subseteq \bigcup_{n = 1}^{\infty}[g_n > 0]$ and $\bigcap_{n = 1}^{\infty}[g_n \leq 0] 
\subseteq [f \leq 0]$; if $x \in X$ with $f(x) > 0$ there is $n \geq 1$ with $g_n (x) > 0$, and if
$\forall_{n \geq 1}\big(g_n (x) \leq 0\big)$, then for the same reason $\neg[ f(x) > 0$, hence $f(x) \leq 0$.\\[1mm]
Because of (i), for (ii), it suffices to show that $\B o_F(f) \in \Borel (\Phi)$. We show that 
\[ \B o_F (f) = \bigcup_{n = 1}^{\infty}\B o_{\Phi}(f \wedge g_n) := 
\bigg(\bigcup_{n = 1}^{\infty}[(f \wedge g_n) > 0], \bigcap_{n = 1}^{\infty}[(f \wedge g_n) \leq 0]\bigg)
\in \Borel(\Phi). \]
If $f(x) > 0$, then we can find $n \geq 1$ such that $g_n (x) > 0$, hence $f(x) \wedge g_n (x) > 0$. Hence
we showed that $[f > 0] \subseteq \bigcup_{n = 1}^{\infty}[(f \wedge g_n) > 0]$. For the converse inclusion,
let $x \in X$ and $n \geq 1$ such that $(f \wedge g_n)(x) > 0$. Then $f(x) > 0$ and $x \in [f > 0]$. 
If $f(x) \leq 0$, then $\forall_{n \geq 1}\big(f(x) \wedge g_n (x) \leq 0\big)$. Suppose next that
$\forall_{n \geq 1}\big(f(x) \wedge g_n (x) \leq 0\big)$. If $f(x) > 0$, there is $n \geq 1$ with $g_n(x) > 0$,
hence $f(x) \wedge g_n (x) > 0$, which contradict the hypothesis $f(x) \wedge g_n (x) \leq 0$. Hence $f(x) \leq 0$.\\
Because of (i), for (iii), it suffices to show that $\B o_F(f) \in \Borel (\Phi)$. We show that  
\[ \B o_F (f) = \bigcup_{n = 1}^{\infty}\B o_{\Phi}\big(g_n - \frac{1}{n}\big) 
:= \bigg(\bigcup_{n = 1}^{\infty}\big[\big(g_n - \frac{1}{n}\big) > 0 \big],
\bigcap_{n = 1}^{\infty}\big[\big(g_n - \frac{1}{n}\big) \leq 0\big]\bigg) \in \Borel(\Phi). \]
First we show that $[f > 0] \subseteq \bigcup_{n = 1}^{\infty}\big[\big(g_n - \frac{1}{n}\big) > 0 \big]$. 
If $f(x) > 0$, there is $n \geq 1$ with $f(x) > \frac{1}{n}$, hence, since $- \frac{1}{2n} \leq g_{2n}(x) - f(x)
\leq \frac{1}{2n}$, we get 
\[ g_{2n}(x) - \frac{1}{2n} \geq \bigg(f(x) - \frac{1}{2n}\bigg) - \frac{1}{2n} = f(x) - \frac{1}{n} > 0 \]
i.e., $x \in \big[\big(g_{2n} - \frac{1}{2n}\big) > 0 \big]$. For the converse inclusion, let $x \in X$ and
$n \geq 1$ such that $g_n (x) - \frac{1}{n} > 0$. Since $0 < g_n (x) - \frac{1}{n} \leq f(x)$, we get $x \in [f > 0]$.
Next we show that $[f \leq 0] \subseteq \bigcap_{n = 1}^{\infty}[\big(g_n - \frac{1}{n}\big) \leq 0]$. 
Let $x \in X$ with $f(x) \leq 0$, and suppose that $n \geq 1$ with $g_n (x) - \frac{1}{n} > 0$. Then 
$0 \geq f(x) >0$. By this contradiction we get $g_n (x) - \frac{1}{n} \leq 0$. For the converse inclusion 
let $x \in X$ such that $g_n (x) - \frac{1}{n} \leq 0$, for every $n \geq 1$, and suppose that $f(x) > 0$.
Since we have already shown that 
$[f > 0] \subseteq \bigcup_{n = 1}^{\infty}\big[\big(g_n - \frac{1}{n}\big) > 0 \big]$, 
there is  some $n \geq 1$ with $g_n (x) - \frac{1}{n} > 0$, which contradicts our hypothesis,
hence $f(x) \leq 0$.
\end{proof}


\section{The Baire sets of a Bishop space}
\label{sec: baire}

One of the definitions\footnote{A different definition is given in~\cite{Ha74}. See~\cite{RS65}
for the relations between these two definitions.} of the set of Baire sets in a topological space 
$(X, \C T)$, which was given by Hewitt in~\cite{He50}, is that it is the least $\sigma$-algebra of
subsets of $X$ that includes the zero sets of $X$ i.e., the sets of the form $f^{-1}(\{0\})$, where 
$f \in C(X)$. Clearly, a Baire set in $(X, \C T)$ is a Borel set in $(X, \C T)$, and for many topological spaces,
like the metrisable ones, the two classes coincide. In this section we adopt Hewitt's notion in Bishop spaces 
and the framework of $F$-complemented subsets.

\begin{definition}\label{def: baire}
Let $\B Z_F(X) := \big(\zeta_0^{1, F}, \C Z^{1,X}, \zeta_1^{1, F}, \zeta_0^{0,F}, \C Z^{0,X}, \zeta_1^{0,F}\big)
\in \Fam(F, F, \B X)$\index{$\B Z_F(X)$} be the family of zero $F$-complemented subsets\index{family of 
zero $F$-complemented subsets} of $X$, where 
\[ \B \zeta_F (f) := \big(\zeta_0^{1,F}(f), \zeta_0^{0,F}(f)\big) := \big([f = 0], [f \neq 0]\big), \]
\[ [f = 0] := \{x \in X \mid f(x) = 0\}, \ \ \ \ [f \leq 0] := \{x \in X \mid f(x) \neq 0\},\]
and, as $[f = 0$, $[f \neq 0]$ are extensional subsets of $X$, the dependent operations $\C Z^{1,X}, \C Z^{0,X}, 
\zeta_1^{1, F}$, and $\zeta_1^{0, F}$ are defined by the identity map-rule. If $F$ is clear from the context,
we may write $\B Z_F(X) := \big(\zeta_0^{1}, \C Z^{1,X}, \zeta_1^{1}, \zeta_0^{0}, \C Z^{0,X}, \zeta_1^{0}\big)$.
Let 
\[ \Baire(\C F) := \Borel(\B Z_F(X)), \]
and we call its elements the Baire sets of $\C F$.

\end{definition}

Since $a \neq 0 : \TOT |a| > 0 \TOT a < 0 \vee a > 0$, for every $a \in \Real$, we get
\[ \B \zeta (f) = \big([f = 0], [|f| > 0]\big) = \big([f = 0], [f > 0] \cup [f < 0]\big). \]

\begin{proposition}\label{prp: baire1}
\normalfont (i)
\itshape $\Baire(\C F)$ is not in $\Set(F, F, \B X)$.\\[1mm]
\normalfont (ii)
\itshape If $f, g \in F$, then $\B \zeta (f) \cap \B \zeta (g) = \B \zeta (|f| \vee |g|)$.\\[1mm]
\normalfont (iii)
\itshape If $\B B \in \Baire(\C F)$, then $- \B B \in \Baire(\C F)$.\\[1mm]
\normalfont (iv)
\itshape There are Bishop space $(X, F)$ and $f \in F$ such that $\neg[ - \B \zeta (f) = \B \zeta (g)$, 
 for every $g \in F$. \\[1mm]
 \normalfont (vi)
 \itshape $\B \zeta(f) = \B \zeta (|f| \wedge 1)$.
\end{proposition}

\begin{proof}
(i) If $f \in F$, then $\B \zeta (f) = \B \zeta (2f)$, but $\neg (f = 2 f)$.\\
(ii) This equality is implied from the following property for reals
$|a| \vee |b| = 0 \TOT |a| = 0 \wedge |b| = 0$ and $|a| \vee |b| \neq 0 \TOT |a| > 0 \vee |b| > 0$.\\
(iii) If $f \in F$, then $- \B \zeta _0(f) := \big([f \neq 0], [f = 0]\big)$. If, for every $n \geq 1$, 
\[ g_n := \bigg(|f| \wedge \frac{1}{n}\bigg) - \frac{1}{n} \in F, \]
\[ \bigcup_{n = 1}^{\infty}\B \zeta (g_n) := \bigg(\bigcup_{n = 1}^{\infty}[g_n = 0], 
\bigcap_{n = 1}^{\infty}[g_n \neq 0]\bigg) = - \B \zeta (f) \in \Baire(\C F). \]
First we show that $[f \neq 0] = \bigcup_{n = 1}^{\infty}[g_n = 0]$. If $|f(x)| > 0$, there 
is $n \geq 1$ such that $|f(x)| > \frac{1}{n}$, hence $|f(x)| \wedge \frac{1}{n} = \frac{1}{n}$,
and $g_n(x) = 0$. For the converse inclusion, let $x \in X$ and $n \geq 1$ such that $g_n (x) = 0 \TOT
|f(x)| \wedge \frac{1}{n} = \frac{1}{n}$, hence $|f(x)| \geq \frac{1}{n} > 0$. Next we show that $[f = 0] 
= \bigcap_{n = 1}^{\infty}[g_n \neq 0]$. If $x \in X$ such that $f(x) = 0$, and $n \geq 1$, then $g_n (x)
= - \frac{1}{n} < 0$. For the converse inclusion, let $x \in X$ such that for all $n \geq 1$ we have that 
$g_n (x) \neq 0$. If $|f(x)| > 0$, there is $n \geq 1$ such that $|f(x)| > \frac{1}{n}$, hence $g_n (x) = 0$, 
which contradicts our hypothesis. Hence, $|f(x)| \leq 0$, which implies that $|f(x)| = 0 \TOT f(x) = 0$.
If $P(\B B) := - \B B \in \Baire(\C F)$, the above equality proves the first step of the corresponding 
induction on $\Baire(F)$. The rest of the inductive proof is straightforward\footnote{Hence, if we define the 
set of Baire sets over an arbitrary family $\Theta$ of functions from $X$ to $\Real$, a sufficient condition
so that $\Baire(\Theta)$ is closed under complements is that $\Theta$ is closed under $|.|$, under wedge with 
$\frac{1}{n}$ and under subtraction with $\frac{1}{n}$, for every $n \geq 1$. If $\Theta := \D F(X, 2)$, then
$- \B o_{\D F(X, 2)}(f) = \B o_{\D F(X, 2)}(1-f) = \B \zeta_{\D F(X, 2)}(f)$, hence by 
Proposition~\ref{prp: borel2}(ii) we get $\Borel(\D F(X, 2)) = \Baire(\D F(X, 2))$.}.\\ 
(v) Let the Bishop space $(\Real, \BR)$. If we take $\B \zeta (\id_{\Real}) := \big([x = 0], [x \neq 0]\big)$,
and if we suppose that $- \B \zeta (\id_{\Real}) := \big([x \neq 0], [x = 0]\big) = \big([\phi = 0], 
[\phi \neq 0]\big) =: \B \zeta (\phi)$, for some $\phi \in \BR$, then $\phi(0) > 0 \vee \phi(0) < 0$ and 
$\phi(x) = 0$, if $x < 0$ or $x > 0$. Hence $\phi$ is not continuous at $0$, which contradicts the fact that
$\phi$ is uniformly continuous on $[-1, 1]$. \\
(v) Using basic properties of $\Real$, this proof is straightforward.
\end{proof}

As in the case of $\Borel(\C F)$, we cannot show constructively that $\B \zeta (f) \cup \B \zeta(g) 
= \B \zeta (|f| \wedge |g|)$. If we add the condition $|f| + |g| > 0$ though, this equality is constructively provable.

\begin{proposition}\label{prp: inf2}
If $(f_n)_{n = 1}^{\infty} \subseteq F$, then $f := \sum_{n = 1}^{\infty}|f_n| \wedge 2^{-n} \in F$ and
$$\B \zeta(f) = \bigcap_{n = 1}^{\infty}\B \zeta (f_n) = \bigg(\bigcap_{n = 1}^{\infty}[f_n = 0], 
\bigcup_{n = 1}^{\infty}[f_n \neq 0]\bigg).$$
\end{proposition}

\begin{proof}
Proceeding as in the proof of Proposition~\ref{prp: inf}, $f$ is well-defined, and if $g_n := |f_n| \wedge 2^{-n}$,
for every $n \geq 1$, 
the sequence of the partial sums $\sum_{n = 1}^N g_n \in F$ converges uniformly to $f$, and by $(\BS_4)$ 
we get $f \in F$. Since $f(x) = 0 \TOT \forall_{n \geq 1}(g_n (x) = 0) \TOT  \forall_{n \geq 1}(f_n (x) = 0)$,
we get $[f = 0] = \bigcap_{n = 1}^{\infty}[f_n = 0]$. To show $[f \neq 0] \subseteq \bigcup_{n = 1}^{\infty}[f_n \neq 0]$,  
if $|f(x)| > 0$, there is $N \geq 1$ such that $\sum_{n = 1}^N g_n(x) > 0$. By Proposition (2.16) 
in~\cite{BB85}, p.~26, there is some $n \geq 1$ and $n \geq N$ such that $g_n(x) > 0$, hence
$|f_n (x)| \geq g_n(x) > 0$. The converse inclusion follows trivially. 
\end{proof}

Let $\C F^* := (X, F^*)$ be the Bishop space generated by the bounded functions $F^*$ in $F$.

\begin{theorem}\label{thm: baire2}
\normalfont (i)
\itshape If $\B B \in \Baire(\C F)$, then $\B B \in \Borel(\C F)$.\\[1mm]
\normalfont (ii)
\itshape If $\B o(f) \in \Baire(\C F)$, for every $f \in F$, then $\Baire(\C F) = \Borel(\C F)$.\\[1mm]
\normalfont (iii)
\itshape If $f \in F$, then $\B o (f) = - \B \zeta \big((-f) \wedge 0\big)$.\\[1mm]
\normalfont (iv)
\itshape $\Baire(\C F^*) = \Baire(\C F) = \Borel(\C F) = \Borel(\C F^*)$.
\end{theorem}

\begin{proof}
(i) By Proposition~\ref{prp: borel1}(iv) $- \B o (f) = \big([f \leq 0], [f > 0]\big) \in \Borel(\C F)$, 
for every $f \in F$, hence $- \B o (-f) = \big([f \geq 0], [f < 0] \in \Borel(\C F)$ too. Consequently
\[ - \B o (f) \cap - \B o (-f) = \big([f \leq 0] \cap [f \geq 0], [f > 0] \cup [f < 0]\big) =
\B \zeta(f)\in \Borel(\C F). \]
If $P(\B B) := \B B \in \Borel(\C F)$, the above equality is the first step of the corresponding 
inductive proof on $\Baire(\C F)$. The rest of the inductive proof is straightforward.\\
(ii) The hypothesis is the first step of the obvious inductive proof on $\Borel(\C F)$, which shows 
that $\Borel(\C F) \subseteq \Baire(\C F)$. By (i) we get $\Baire(\C F) \subseteq \Borel(\C F)$.\\
(iii) We show that 
\[ \big([f > 0], [f \leq 0]\big) = \big([(-f) \wedge 0 \neq 0], [(-f) \wedge 0 = 0]\big). \]
First we show that $[f > 0] \subseteq [(-f) \wedge 0 \neq 0]$; if $f(x) > 0$, then $-f(x) \wedge 0 = -f(x) < 0$.
For the converse inclusion, let $-f(x) \wedge 0 \neq 0 \TOT -f(x) \wedge 0 > 0$ or $-f(x) \wedge 0 < 0$. 
Since $0 \geq -f(x) \wedge 0$, the first option is impossible. If $-f(x) \wedge 0 < 0$, then $-f(x) < 0$ or
$0 < 0$, hence $f(x) > 0$.
Next we show that $[f \leq 0] = [(-f) \wedge 0 = 0]$; since $f(x) \leq 0 \TOT - f(x) \geq 0 \TOT - f(x) 
\wedge 0 = 0$ (see~\cite{BV06}, p.~52), the equality follows.\\
(iv) Clearly, $\Baire (\C F^*) \subseteq \Baire(\C F)$. By Proposition~\ref{prp: baire1}(vi) $\B \zeta(f)
= \B \zeta(|f| \wedge 1)$, where $|f| \wedge 1 \in F^*$. Continuing with the obvious induction we
get $\Baire(\C F) \subseteq \Baire(\C F^*)$. By case (iii) and Proposition~\ref{prp: baire1}(iv) we 
get $\B o (f) \in \Baire(\C F)$, hence by case (ii) we conclude that $\Baire(\C F) = \Borel(\C F)$. 
Clearly, $\Borel (\C F^*) \subseteq \Borel(\C F)$. By Proposition~\ref{prp: borel1}(vi) $\B o(f) =
\B o((f \vee 0) \wedge 1)$, where $(f \vee 0) \wedge 1 \in F^*$. Continuing with the obvious induction
we get $\Borel(\C F) \subseteq \Borel(\C F^*)$. 
\end{proof}

Either by definition, as in the proof of Proposition~\ref{prp: borel2}, or by Theorem~\ref{thm: baire2}(iii) 
and Proposition~\ref{prp: borel2}, if $h \in \Mor(\C F, \C G)$ and $\B B \in \Baire(\C G)$, then
$h^{-1}(\B B) \in \Baire(\C F)$.
Suppose next that $\B A$ is strongly $F$-complemented i.e., there is $f \in F$ such that 
$f \colon A^1 \Disj_F A^0$ and $f(x) = 1$, for every $x \in A^1$, and $f(y) = 0$, for every $y \in A^0$.
If $g := (f \vee 0) \wedge 1 \in F$, then $0 \leq g \leq 1$ and $\forall_{x \in A^1}\forall_{y \in A^0}\big(g(x)
= 1 \ \& \ g(y) = 0\big)$. In~\cite{BC72}, p. 55, the following relation between complemented subsets is defined:
\[ \B A \leq \B B : \TOT A^1 \subseteq B^1 \ \ \& \ \ A^0 \subseteq B^0. \]  
If $\B A$ is strongly $F$-complemented, then $\B A \leq \B o (f)$. According to the classical
Urysohn lemma for $C(X)$-zero sets, the disjoint zero sets of a topological
space $X$ are separated by some $f \in C(X)$ (see~\cite{GJ60}, p.~17). Next we show a constructive version
of this result, where disjointness is replaced by a stronger, but positively defined form of it.

\begin{theorem}[Urysohn lemma for zero complemented subsets]\label{thm: urysohn}
If $\B A := (A^1, A^0) \in \C P^{][_F}(X)$, then $\B A$ is strongly $F$-complemented if and only if 
\[ \exists_{f, g \in F}\exists_{c > 0}\big(\B A \leq \B \zeta(f) \ \& \ - \B A \leq \B \zeta(g) \ \& \ 
|f| + |g| \geq c\big). \]
\end{theorem}

\begin{proof}
$(\To)$ Let $h \in F$ such that $0 \leq h \leq 1$, $A^1 \subseteq [h = 1]$ and $A^0 \subseteq [h = 0]$. 
We take $f := 1 - h \in F, g := h$ and $c := 1$. First we show that $\B A \leq \B \zeta(f)$. If 
$x \in A^1$, then $h(x) = 1$, and $f(x) = 0$. If $y \in A^0$, then $h(y) = 0$, hence $f(y) = 1$ and 
$y \in [f \neq 0]$. Next we show that $- \B A \leq \B \zeta(g)$. If $y \in A^0$, then $h(y) = 0 = g(y)$. 
If $x \in A^1$, then $h(x) = 1 = g(y)$ i.e., $x \in [g \neq 0]$. If $x \in X$, then $1 = |1 - h(x) + h(x)| 
\leq |1 - h(x)| + |h(x)|$.\\
$(\oT)$ Let $h := 1 - \big(\frac{1}{c}|f| \wedge 1\big) \in F$. If $x \in A^1$, then $f(x) = 0$, and
hence $h(x) = 1$. If $y \in A^0$, then $g(y) = 0$, hence $|f(y)| \geq c$, and consequently $h(y) = 0$.
\end{proof}

The condition $(\BS_3)$ of a Bishop space
is crucial to the next proof.

\begin{corollary}\label{cor: corur}
ILet$\B A := (A^1, A^0) \in \C P^{][_F}(X)$ and $f \in F$. If $f(\B A) := \big(f(A^1), f(A^0)\big)$ is 
strongly $\BR$-complemented, then $\B A$ is strongly $F$-complemented.
\end{corollary}

\begin{proof}
By the Urysohn lemma for zero complemented subsets there are $\phi, \theta \in \BR$ and $c > 0$ with 
$f(\B A) \leq \B \zeta(\phi), - f(\B A) \leq \B \zeta(\theta)$ and $|\phi| + |\theta| \geq c$.
Consequently, $\B A \leq \B \zeta(\phi \circ f), - \B A \leq \B \zeta(\theta \circ f)$ and $|\phi \circ f| 
+ |\theta \circ f| \geq c$. Since by $(\BS_3)$ we have that $\phi \circ f \in F$ and $\theta \circ f \in F$,
by the other implication of the Urysohn lemma for zero complemented subsets we conclude that $\B A$ is
strongly $F$-complemented.
\end{proof}

\section{Measure and pre-measure spaces}
\label{sec: premeasure}

There are two, quite different, notions of measure space in traditional Bishop-style constructive mathematics. The 
first, which was introduced in~\cite{Bi67} as part of \textit{Bishop's measure theory} $(\BMT)$ 
(see Note~\ref{not: Bms}), is an abstraction of the measure function $A \mapsto \mu(A)$, where $A$ is a
member of a family of complemented subsets of a locally compact metric space $X$. The use of complemented 
subsets in order to overcome the difficulties 
generated in measure theory by the use of negation and negatively defined concepts is one of Bishop's great
conceptual achievements, while the use of the concept of a family of complemented subsets 
is crucial to the predicative character of this notion of measure space\footnote{Myhill's impredicative 
interpretation in~\cite{My75} of Bishop's first definition is discussed in Note~\ref{not: onMy75}.}. The 
indexing required behind this first notion of measure space is evident in~\cite{Bi67}, and sufficiently 
stressed in~\cite{Bi70} (see Note~\ref{not: numerical}).
The second notion of measure space, introduced in~\cite{BC72} and repeated in~\cite{BB85} as part of 
the far more general \textit{Bishop-Cheng measure theory} $(\BCMT)$, is highly impredicative, as the 
necessary indexing for its predicative reformulation is missing. A lack of predicative concern is evident 
also in the integration theory of $\BCMT$. 
Next we define a predicative variation of the Bishop-Cheng notion of measure space using the predicative
conceptual ingredients of the initial Bishop notion of measure space. We also keep the operations of 
complemented subsets introduced in~\cite{Bi67}, and not the operations used in~\cite{BC72} and~\cite{BB85}. 
Following Bishop's views in~\cite{Bi70}, we introduce the notion of pre-measure space, which is understood though, 
in a way different from the classical term.

As in Definition~\ref{def: sublambdaI}, if $\B \Lambda(X) := (\lambda_0^1, \C E^{1,X}, \lambda_1^1, 
\lambda_0^0, \C E^{0,X}, \lambda_1^0) \in \Fam(I, \B X)$, the set  $\lambda_0 I(\B X)$ of complemented 
subsets of $X$ is the totality $I$, equipped with the equality 
$i =_I^{\Lambda(X)} j :\TOT \B \lambda_0(i) =_{\C P^{\Disj}(X)} \B \lambda_0(j)$, for every 
$i, j \in I$. For simplicity we write $\B \lambda_0(i) := (\lambda_0^1(i), \lambda_0^0(i))$ instead of $i$
for an element of $\lambda_0I(\B X)$.

\begin{definition}[Measure space within $\BST$]\label{def: measurespace}
Let $(X, =_X, \neq_X)$ be an inhabited set, $\B \Lambda(X) := (\lambda_0^1, \C E^{1,X}, \lambda_1^1 ; \lambda_0^0, \C E^{0,X}, \lambda_1^0) \in \Fam(I, \B X)$, and let $\mu \colon \lambda_0 I(\B X) \to [0, + \infty)$ such that the following conditions hold:
\vspace{-2mm}
\[ (\MS_1) \  \ \ \  \ \ \ \ \ \ \ \ \ \  \forall_{i, j \in I}\exists_{k, l \in I}\bigg(\B \lambda_0 (i) \cup \B \lambda_0 (j) = \B \lambda_0 (k) \ \& \ \B \lambda_0 (i) \cap \B \lambda_0 (j) = \B \lambda_0 (l) \ \& \ \ \ \ \ \ \ \ \ \ \ \ \ \ \ \ \ \ \ \   \]
\[ \ \ \ \ \ \ \ \ \ \ \ \ \ \ \  \mu \big(\B \lambda_0 (i)\big) + \mu \big(\B \lambda_0 (j)\big) = \mu \big(\B \lambda_0 (k)\big) + \mu \big(\B \lambda_0 (l)\big)\bigg). \]
\[ (\MS_2)  \ \ \ \  \ \ \ \ \ \ \ \ \ \  \forall_{i \in I}\forall_{\B A(X) \in \Fam(\D 1, \B X)}\bigg[\exists_{k \in I}\bigg(\B \lambda_0 (i) \cap \B \alpha_0 (0) = \B \lambda_0 (k)\bigg) \To  \ \ \ \ \ \ \ \ \ \ \ \ \ \ \ \ \ \ \ \ \ \ \ \ \ \ \ \ \ \ \]
\[ \ \ \ \ \ \ \ \ \  \big(\exists_{l \in I}\big(\B \lambda_0 (i) - \B \alpha_0 (0) = \B \lambda_0 (l)\big) \ \& \ \mu \big(\B \lambda_0 (i)\big) = \mu \big(\B \lambda_0 (k)\big) + \mu \big(\B \lambda_0 (l)\big)\bigg]. \]
\[ (\MS_3) \ \ \ \ \ \ \ \ \ \ \ \ \ \ \ \ \ \ \  \ \ \  \ \ \  \ \ \ \ \ \ \ \ \  \ \ \ \    \exists_{i \in I}\big(\mu \big(\B \lambda_0 (i)\big) > 0. 
\ \ \ \ \ \ \ \ \ \ \ \ \ \ \ \   \ \ \ \ \ \ \ \ \ \ \ \ \ \ \ \ \ \ \ \ \ \ \  \ \ \ \ \ \ \ \ \ \ \    \]
\[ (\MS_4) \  \ \ \ \ \ \ \   \ \ \ \ \ \ \ \ \ \ \forall_{\alpha \in \D F(\Nat, I)}\bigg\{\forall_{\beta \in \D F(\Nat, I)}\bigg[\forall_{m \in \Nat}\bigg(\bigcap_{n = 1}^m \B \lambda_0 (\alpha(n)) = \B \lambda_0 (\beta(m))\bigg) \ \& \ \ \ \ \ \ \ \ \ \ \ \ \ \ \ \ \ \ \  \ \ \ \ \ \ \  \ \ \ \ \  \]
\[ \ \ \ \ \ \ \ \exists  \lim_{\mathsmaller{m \to + \infty}} \mu \big(\B \lambda_0 (\beta(m))\big) \  \& \ \lim_{\mathsmaller{m \to + \infty}} \mu \big(\B \lambda_0 (\beta(m))\big) > 0 \To \exists_{x \in X}\bigg(x \in \bigcap_{n \in \Nat}\lambda_0^1 (\alpha (n))\bigg)\bigg]\bigg\}.\]
The triplet $\C M := (X, \lambda_0 I(\B X), \mu)$ is called a \textit{measure space}\index{measure space} with $ \lambda_0 I(\B X)$ its set of \textit{integrable}, or \textit{measurable}\index{integrable set} sets\index{measurable set}, and $\mu$ its \textit{measure}\index{measure of a measure space}.
\end{definition}

With respect to condition $(\MS_1)$, we do not say that the set $\lambda_0I (\B X)$ is closed
under the union or intersection of complemented subsets (as Bishop-Cheng do in their definition).
This amounts to the rather strong condition
$\forall_{i, j \in I}\exists_{k, l \in I}\big(\B \lambda_0 (i) \cup \B \lambda_0 (j) :=
\B \lambda_0 (k) \ \& \ \B \lambda_0 (i) \cap \B \lambda_0 (j) := \B \lambda_0 (l)\big)$. 
The weaker condition $(\MS_1)$ states that the complemented subsets $\B \lambda_0 (i) \cup 
\B \lambda_0 (j)$ and $\B \lambda_0 (i) \cap \B \lambda_0 (j)$ ``pseudo-belong'' to $\lambda_0I (\B X)$ i.e., 
there are elements of it, which are equal to them in $\C P^{\Disj}(X)$.
In contrast to the formulation of condition $(\MS_2)$ by Bishop and Cheng, we avoid quantification over
the class $\C P^{\Disj}(X)$, by quantifying over the set $\Fam(\D 1, \B X)$. In our formulation of
$(\MS_4)$ we quantify over $\D F(\Nat, I)$, in order to avoid the use of some choice principle.
If we had written   
\[ \forall_{m \in \Nat}\exists_{k \in \Nat}\bigg(\bigcap_{n = 1}^m \B \lambda_0 (\alpha(n)) 
=_{\mathsmaller{\C P^{\Disj_{\neq}}}} \B \lambda_0 (k)\bigg) \]
instead, we would need countable choice to express the limit to infinity of the terms
$\mu\big(\B \lambda_0 (k)\big)$. 
Next we define the notion of a pre-measure space, giving an explicit formulation of Bishop's idea,
expressed in~\cite{Bi70}, p.~67, and quoted in Note~\ref{not: numerical}, to formalise his
first definition of measure space, applied though, to Definition~\ref{def: measurespace}.
The main idea is to define operations on $I$ that correspond to the operations on complemented subsets,
and reformulate accordingly the clauses for the measure $\mu$. The fact that $\mu$ is defined on 
the index-set is already expressed in the definition of the set $\lambda_0I (\B X)$. The notion of a 
pre-measure space provides us a method to generate measure spaces.

\begin{definition}[Pre-measure space within $\BST$]\label{def: premeasurespace}
Let $(X, =_X, \neq_X)$ be an inhabited set, and let $(I, =_I)$ be equipped with operations 
$\vee \colon I \times I \sto I, \wedge \colon I \times I \sto I$, and $\sim \colon I \sto I$.
If $i, j \in I$ and $i_1, \ldots, i_m \in I$, where $m \geq 1$, let\footnote{The operations $\bigvee_{n = 1}^m i_n$
and $\bigvee_{n = 1}^m i_n$ are actually recursively defined.}
\[ i \sim j := i \wedge (\sim j) \ \ \&  \ \ i \leq j : \TOT i \wedge j = i,\]
\[ \bigvee_{n = 1}^m i_n := i_1 \vee \ldots \vee i_m \ \ \& \ \ \bigwedge_{n = 1}^m i_n := i_1 \wedge 
\ldots \wedge i_m.\]
Let $\B \Lambda(X) := (\lambda_0^1, \C E^{1,X}, \lambda_1^1 ; \lambda_0^0, \C E^{0,X}, \lambda_1^0)
\in \Set(I, \B X)$, and $\mu \colon I \to [0, + \infty)$ such that the following conditions hold:
\vspace{-2mm}
\[ (\PMS_1) \  \  \forall_{i, j \in I}\bigg(\B \lambda_0 (i) \cup \B \lambda_0 (j) = \B \lambda_0 (i \vee j)
\ \& \ \B \lambda_0 (i) \cap \B \lambda_0 (j) = \B \lambda_0 (i \wedge j) \ \& \  - \B \lambda_0 (i) = 
\B \lambda_0 (\sim i) \ \& \]
\[ \mu(i) + \mu(j) = \mu(i \vee j) + \mu(i \wedge j)\bigg).\]
\[ (\PMS_2) \ \ \   \ \ \ \ \ \ \ \ \ \ \ \ \ \ \ \ \   \ \ \ \ \ \ 
\forall_{i \in I}\forall_{\B A(X) \in \Fam(\D 1, \B X)}\bigg[\exists_{k \in I}\bigg(\B \lambda_0 (i) \ \cap \ \B \alpha_0 (0) = \B \lambda_0 (k)\bigg) \ \To \ \ \ \ \ \ \ \ \ \ \ \ \ \ \ \ \ \ \ \  \ \ \ \ \ \ \ \ \ \ \ \  \]
\[ \ \ \ \ \ \ \ \ \ \ \ \ \ \B \lambda_0 (i) - \B \alpha_0(0) = \B \lambda_0 (i \sim k) \ \& \ \mu(i) = \mu(k) + \mu(i \sim k)\bigg]. \ \ \ \ \ \ \ \ \ \ \]
\[ (\PMS_3) \ \ \ \ \ \ \ \ \ \ \  \ \ \ \ \ \  \ \ \  \ \ \ \ \ \ \ \ \ \ \ \ \ \  \ \ \ \ \ \ \    \exists_{i \in I}\big(\mu (i)\big) > 0. \ \ \ \ \ \ \ \  \ \ \ \ \ \ \ \ \ \ \  \ \ \ \ \ \ \ \  \ \ \ \ \ \ \ \ \ \ \ \ \ \ \ \ \ \ \ \  \]
\[ (\PMS_4) \ \ \ \ \ \ \ \ \  \ \ \ \ \ \ \  \forall_{\alpha \in \D F(\Nat, I)}\bigg[\exists  \lim_{\mathsmaller{m \to + \infty}} \mu \bigg(\bigwedge_{n = 1}^m \alpha(n)\bigg) \  \& \ \lim_{\mathsmaller{m \to + \infty}} \mu \bigg(\bigwedge_{n = 1}^m \alpha(n)\bigg) > 0 \To \ \ \ \ \ \ \ \ \  \ \ \ \ \ \ \ \ \ \]
\[ \ \ \ \ \ \  \ \ \ \ \  \To \exists_{x \in X}\bigg(x \in \bigcap_{n \in \Nat}\lambda_0^1 (\alpha (n))\bigg)\bigg].\]
The triplet $\C M (\B \Lambda(X)) := (X, I, \mu)$ is called 
a \textit{pre-measure space}\index{pre-measure space}, the function $\mu$ a \textit{pre-measure}\index{pre-measure}, 
and the index-set $I$ a set of integrable, or measurable \index{set of integrable indices}indices.
\end{definition}

\begin{corollary}\label{cor: premeasure1}
Let $\C M (\B \Lambda) := (X, I, \mu)$ be a pre-measure space and $i, j \in I$.\\[1mm]
\normalfont (i)
\itshape The operations $\vee$, $\wedge$ and $\sim$ are functions.\\[1mm]
\normalfont (ii)
\itshape The triplet $(I, \vee, \wedge)$ is a distributive lattice.\\[1mm]
\normalfont (iii)
\itshape $\sim (\sim i) =_I i$.\\[1mm]
\normalfont (iv)
\itshape $\sim (i \wedge j) =_I (\sim i) \vee (\sim j)$.\\[1mm]
\normalfont (v)
\itshape $i \leq j \TOT \ \sim j \leq \ \sim i$.\\[1mm]
\normalfont (vi) 
\itshape $i \leq j \TOT \B \lambda_0 (i) \subseteq \B \lambda_0 (j)$.\\[1mm]
\normalfont (vii)
\itshape  $\B \lambda_0 (i) - \B \lambda_0 (j) = \B \lambda_0 (i \sim j)$.
\end{corollary}

\begin{proof}
We show that $\vee$ is a function, and for $\wedge$ and $\sim$ we proceed similarly.
\begin{align*}
i = i{'} \ \& \ j = j{'} & \To \B \lambda_0 (i) = \B \lambda_0 (i{'}) \ \& \ \B \lambda_0 (j) = \B \lambda_0 (j{'})\\
& \To \B \lambda_0 (i) \cup \B \lambda_0 (j) = \B \lambda_0 (i{'}) \cup \B \lambda_0 (j{'})\\
& \To \B \lambda_0 (i \vee j) = \B \lambda_0 (i{'} \vee j{'})\\
& \To i \vee j = i{'} \vee j{'}.
\end{align*}
(ii) The defining clauses of a distributive lattice
follow from the corresponding properties of complemented subsets for $\B A \cap \B B$ and $\B A \cup \B B$,
from $(\PMS_1)$, and from the fact that $\B \Lambda(X) \in \Set(I, \B X)$. E.g., to show $i \vee j = j \vee i$, 
we use the equalities $\B \lambda_0 (i \vee j) = \B \lambda_0 (i) \cup \B \lambda_0 (j) = \B \lambda_0 (j) 
\cup \B \lambda_0 (i) = \B \lambda_0 (j \vee i)$. For the rest of the proof we proceed similarly.
\end{proof}

In the next example of a pre-measure space the index-set $I$ is a Boolean algebra. 

\begin{proposition}\label{prp: predetachable1}
Let $\big(X, =_X, \neq_X^{\mathsmaller{\D F(X, \D 2)}}\big)$ be a set, and
$\B \Delta(X) := \big(\delta_0^1, \C E^{1,X}, \delta_1^1, \delta_0^0, \C E^{0,X}, \delta_1^0\big) 
\in \Set(\D F(X, \D 2), \B X)$ the family of complemented detachable subsets of $X$, where 
by Remark~\ref{rem: detachablefam}
\[ \B \delta_0(f) := \big(\delta_0^1(f), \delta_0^0(f)\big) := \big([f = 0], [f = 1]\big).\]
If $x_0 \in X$ and $\mu_{x_0} \colon \D F (X, \D 2) \sto [0, + \infty)$ 
is defined by the rule
\[ \mu_{x_0} (f) := f(x_0); \ \ \ \ f \in \D F(X, \D 2), \]
then the triplet $\C M(\B \Delta(X)) := (X, \D F(X, \D2), \mu_{x_0})$ is a pre-measure space.
\end{proposition}

\begin{proof}
We define the maps $\vee, \wedge : \D F(X, \D2) \times \D F(X, \D2) \to \D F(X, \D2)$ and $\sim \colon \D F(X, \D2)
\to \D F(X, \D2)$ by
\[ f \vee g := f + g - fg, \ \ \ f \wedge g := fg, \ \ \ \sim f := 1 - f; \ \ \ \ f, g \in \D F(X, \D 2),\]
where $1$ also denotes the constant function on $X$ with value $1$.
By definition of the union and intersection of complemented subsets we have that
\begin{align*}
\B \delta_0(f) \cup \B \delta_0(g) & := \big(\delta_0^1(f) \cup \delta_0^1(g), \delta_0^0(f) 
\cap \delta_0^0(g)\big)\\
& = \big(\delta_0^1(f) \cup \delta_0^1(g), \delta_0^1(1-f) \cap \delta_0^1(1-g)\big)\\
& = \big(\delta_0^1(f+g-fg), \delta_0^1((1-f)(1-g))\big)\\
& = \big(\delta_0^1(f+g-fg), \delta_0^1(1 -(f + g - fg))\big)\\
& = \big(\delta_0^1(f+g-fg), \delta_0^0(f + g - fg)\big)\\
& := \B \delta_0 (f \vee g).
\end{align*}
\begin{align*}
\B \delta_0(f) \cap \B \delta_0(g) & := \big(\delta_0^1(f) \cap \delta_0^1(g), \delta_0^0(f) \cup \delta_0^0(g)\big)\\
& = \big(\delta_0^1(f) \cap \delta_0^1(g), \delta_0^1(1-f) \cup \delta_0^1(1-g)\big)\\
& = \big(\delta_0^1(fg), \delta_0^1((1-f) + (1-g) - (1-f)(1-g))\big)\\
& = \big(\delta_0^1(fg), \delta_0^1(1 -fg)\big)\\
& = \big(\delta_0^1(fg), \delta_0^0(fg)\big)\\
& := \B \delta_0 (f \wedge g).
\end{align*}
Clearly, $ \B \delta_0 (\sim f) :=  \B \delta_0 (1 - f) = -  \B \delta_0(f)$.
Clearly, the operation $\mu_{x_0}$ is a function. As
\[ \mu_{x_0}(f) + \mu_{x_0}(g) = \mu_{x_0}(f + g - fg) + \mu_{x_0}(fg) \TOT \]
\[ f(x_0) + g(x_0) = f(x_0) + g(x_0) - f(x_0)g(x_0) + f(x_0)g(x_0), \]
which is trivially the case, $(\PMS_1)$ follows. 
Let $f \in \D F (X, \D 2)$ and $\B B := (B^1, B^0)$ a given complemented subset of $X$ with
$\B \alpha_0 (0) := \B B$. If $g \in \D F(X, \D 2)$ such that 
\[\B \delta_0(f) \cap \B B := \big(\delta_0^1(f) \cap B^1, \delta_0^1(1-f) \cup B^0\big) = 
\big(\delta_0^1(g), \delta_0^1(1-g)\big) \TOT\]
\[ \delta_0^1(f) \cap B^1 = \delta_0^1(g) \ \ \& \ \  \delta_0^1(1-f) \cup B^0 = \delta_0^1(1-g),\]
\begin{align*}
\B \delta_0 (f (1 - g)) &:= \big(\delta_0^1 (f (1-g)), \ \delta_0^0 (f (1-g))\big)\\
& = \big(\delta_0^1 (f) \cap \delta_0^1 (1-g), \ \delta_0^0 (f) \cup \delta_0^1 (1-g)\big)\\
& = \big(\delta_0^1 (f) \cap [\delta_0^1(1-f) \cup B^0], \ \delta_0^0 (f) \cup [\delta_0^1 (f)
\cap B^1]\big)\\
& = \big([\delta_0^1 (f) \cap \delta_0^1(1-f)] \cup [\delta_0^1 (f) \cap B^0], \ [\delta_0^0 (f) 
\cup \delta_0^1 (f) \cap [\delta_0^0 (f) \cup B^1]\big)\\
& = \big(\emptyset \cup [\delta_0^1 (f) \cap B^0], \ X \cap [\delta_0^0 (f) \cup B^1]\big)\\
& = \big(\delta_0^1 (f) \cap B^0, \ \delta_0^0 (f) \cup B^1\big)\\
& := \B \delta_0(f) - \B B.
\end{align*}
To complete the proof of $(\PMS_2)$, we need to show
\[ \mu_{x_0}(f) = \mu_{x_0}(g) + \mu_{x_0}(f(1-g)) \TOT f(x_0) = g(x_0) + f(x_0)(1 - g(x_0)).\]
If $g(x_0) = 0$, the equality holds trivially. If $g(x_0) = 1$, and since $\delta_0^1 (g) = \delta_0^1 (f) \cap B^1$,
we also have that $f(x_0) = 1$, and the required equality holds. 
As $\mu_{x_0}(1) = 1 > 0$, $(\PMS_3)$ follows. For the proof of $(\PMS_4)$ we fix $\alpha : \Nat \to \D F(X, \D 2)$,
 and we suppose that 
\begin{align*}
& \exists  \lim_{\mathsmaller{m \to + \infty}} \mu_{x_0}\bigg(\bigwedge_{n = 0}^m \alpha_n\bigg) \  \& \ 
\lim_{\mathsmaller{m \to + \infty}} \mu_{x_0}\bigg(\bigwedge_{n = 0}^m \alpha_n\bigg) > 0 \ \TOT \\
& \exists  \lim_{\mathsmaller{m \to + \infty}} \bigg(\bigwedge_{n = 0}^m \alpha_n\bigg)(x_0) \ \& \ 
\lim_{\mathsmaller{m \to + \infty}} \bigg(\bigwedge_{n = 0}^m \alpha_n\bigg)(x_0) > 0 \ \TOT \\
& \exists  \lim_{\mathsmaller{m \to + \infty}} \prod_{n = 0}^m \alpha_n (x_0)  \ \& \ 
\lim_{\mathsmaller{m \to + \infty}} \prod_{n = 0}^m \alpha_n(x_0) > 0.
\end{align*}

Finally, we have that
\begin{align*}
\lim_{\mathsmaller{m \to + \infty}} \prod_{n = 0}^m \alpha_n(x_0) > 0 & \To 
\lim_{\mathsmaller{m \to + \infty}} \prod_{n = 0}^m \alpha_n(x_0) =1\\
& \TOT \exists_{m_0 \in \Nat}\forall_{m \geq m_0}\bigg(\prod_{n = 0}^m \alpha_n(x_0) = 1\bigg)\\
& \To \forall_{n \in \Nat}\big(\alpha_n (x_0) = 1\big)\\
& \TOT x_0  \in \bigcap_{n \in \Nat}\delta_0^1 (\alpha_n).\qedhere
\end{align*}
\end{proof}

\begin{proposition}\label{prp: lifting}
Let $\C M(\B \Lambda(X)) := (X, I, \mu)$ be a pre-measure space. If $\mu^* \colon \B \lambda_0 I \sto [0, + \infty)$,
where $\mu^* (\B \lambda_0 (i)) := \mu (i),$ for every $\B \lambda_0 (i) \in \B \lambda_0 I$, then 
$\C M := (X, \lambda_0 I(\B X), \mu^*)$ is a measure space.
\end{proposition}

\begin{proof}
By Proposition~\ref{prp: subFamtoset1} $\mu^*$ is a function.
For the proof of $(\MS_1)$ we fix $i, j \in I$ and we take $k := i \vee j$ and $l := i \wedge j$. From 
$(\PMS_1)$  and $(\PMS_4)$ we get
\begin{align*}
\mu^* (\B \lambda_0 (i)) + \mu^* (\B \lambda_0 (j)) & := \mu(i) + \mu(j)\\
& = \mu(i \vee j) + \mu(i \wedge j)\\
& := \mu^* (\B \lambda_0 (i \vee j)) + \mu^* (\B \lambda_0 (i \wedge j)).
\end{align*} 
For the proof of $(\MS_2)$ we fix $i \in I$ and $\B A \in \Fam (\D 1, \B X)$ with $\B \alpha_0 (0) := \B B$. 
If $\B \lambda_0 (i) \ \cap \ \B B = \B \lambda_0 (k)$, for some $k \in I$, we take $l := i \sim k \in I$ and
by $(\PMS_2)$
$ \mu^* (\B \lambda_0 (i)) := \mu_0 (i) = \mu(k) + \mu(i \sim k) :=  \mu^* (\B \lambda_0 (k)) + \mu^* 
(\B \lambda_0 (i \sim k))$. Condition $(\MS_3)$ follows immediately from $(\PMS_3)$. 
For the proof of $(\MS_4)$, we fix $\alpha, \beta \in \D F(\Nat, I)$, and we suppose that
\[ \forall_{m \in \Nat}\bigg(\bigcap_{n = 1}^m \B \lambda_0 (\alpha(n)) = \B \lambda_0 (\beta(m))\bigg) \ \& 
\ \exists  \lim_{\mathsmaller{m \to + \infty}} \mu^* \big(\B \lambda_0 (\beta(m))\big) \  \& \ 
\lim_{\mathsmaller{m \to + \infty}} \mu^* \big(\B \lambda_0 (\beta(m))\big) > 0 \TOT \]
\[ \forall_{m \in \Nat}\bigg(\bigcap_{n = 1}^m \B \lambda_0 (\alpha(n)) = \B \lambda_0 (\beta(m))\bigg) \ \&  
\ \exists  \lim_{\mathsmaller{m \to + \infty}} \mu \big(\beta(m)\big) \  \& \ \lim_{\mathsmaller{m \to + \infty}}
\mu \big(\beta(m)\big) > 0. \] 
If $m \geq 1$, by $(\PMS_1)$ we have that
\[ \B \lambda_0 (\beta(m)) = \bigcap_{n = 1}^m \B \lambda_0 (\alpha(n)) = \B \lambda_0 
\bigg(\bigwedge_{n = 1}^m \alpha(n)\bigg), \]
hence, since $\B \lambda_0$ is a set of complemented subsets, $\beta(m) = \bigwedge_{n = 1}^m \alpha(n)$,
and consequently
$\mu\big(\beta(m)\big) = \mu \big(\bigwedge_{n = 1}^m \alpha(n)\big)$. Hence
\[ \exists  \lim_{\mathsmaller{m \to + \infty}} \mu \bigg(\bigwedge_{n = 1}^m \alpha(n)\bigg) \  \& \ 
\lim_{\mathsmaller{m \to + \infty}} \mu \bigg(\bigwedge_{n = 1}^m \alpha(n)\bigg) > 0.\]
By $(\PMS_4)$ we conclude that there is some $x \in X$ such that 
$x \in \bigcap_{n \in \Nat}\lambda_0^1 (\alpha (n))$.
\end{proof}

\begin{corollary}\label{cor: measurespaceofdetachables2}
Let $\C M(\B \Delta(X)) := (X, \D F(X, \D2), \mu_{x_0})$ be the pre-measure space of complemented 
detachable subsets of $X$. If $\mu_{x_0}^* : \delta_0 \D F(X, \D 2)(\B X) \sto [0, + \infty)$ is 
defined by $\mu_{x_0}^* \big(\B \delta_0 (f)\big) := \mu_{x_0} (f) := f(x_0),$ for every
$\delta_0 (f) \in \delta_0 \D F(X, \D 2)(\B X)$, then 
$\C M(X) := (X, \delta_0 \D F(X, \D 2)(\B X), \mu^*_{x_0})$\index{$\C M(X)$} is a measure space.
\end{corollary}

Next we formulate in our framework the definition of a complete measure space given by Bishop and
Cheng\footnote{In the definition of Bishop and Cheng the symbol of definitional equality $l := 
\lim_{\mathsmaller{m \to + \infty}} \mu \big(\B \lambda_0 (\beta(m))\big)$ is used, but as this a 
convergence condition, one can use the equality of $\Real$ for the same purpose} (see Note~\ref{not: BCms}). 

\begin{definition}\label{def: completems}
A measure space $\C M := (X, \lambda_0 I(\B X), \mu)$ is called \textit{complete}\index{complete measure space},
if the following conditions hold:
\vspace{-2mm}
\[ (\CM_1) \ \ \ \ \  \ \ \ \ \ \ \ \ \ \   \forall_{i \in I}\forall_{\B A(X) \in \Fam(\D 1, \B X)}\bigg( \lambda_0^1 (i) \subseteq  \alpha_0^1(0) \ \& \ \lambda_0^0 (i) \subseteq  \alpha_0^0(0) \To \exists_{k \in I}\big(\B \alpha_0(0) = \B \lambda_0(k)\big)\bigg). \ \ \ \ \ \ \ \ \ \ \ \ \ \ \ \ \ \ \ \ \ \ \ \ \ \  \]
\[ (\CM_2) \ \ \ \ \ \ \ \ \ \ \ \ \ \ \ \ \ \ \ \ \ \ \ \ \  \forall_{\alpha\in \D F(\Nat, I)}\bigg\{\forall_{\beta \in \D F(\Nat, I)}\forall_{l \in [0, + \infty]}\bigg[\forall_{m \in \Nat}\bigg(\bigcup_{n = 1}^m \B \lambda_0 (\alpha(n)) = \B \lambda_0 (\beta(m))\bigg) \ \& \ \ \ \ \ \ \ \ \ \ \ \ \ \ \ \ \ \ \  \ \ \ \ \ \ \  \ \ \ \ \  \]
\[  \ \ \ \ \ \ \ \ \ \ \  \ \ \ \ \ \ \ \ \ \ \ \ \ \ \ \ \ \ \ \ \ \ \ \ \ \exists  \lim_{\mathsmaller{m \to + \infty}} \mu \big(\B \lambda_0 (\beta(m))\big) \  \& \ \lim_{\mathsmaller{m \to + \infty}} \mu \big(\B \lambda_0 (\beta(m))\big) = l \ \  \ \ \ \ \ \ \ \ \ \]
\[ \  \ \ \ \ \ \ \ \ \ \ \ \ \ \ \ \ \ \ \ \ \ \ \ \ \ \ \ \ \  \ \To \exists_{k \in I}\bigg(\bigcup_{n \in \Nat} \B \lambda_0 (\alpha(n)) = \B \lambda_0 (k)  \ \& \ \mu\big(\B \lambda_0 (k)\big) = l\bigg)\bigg]\bigg\}. \]
\[ (\CM_3) \ \  \ \ \ \ \ \ \ \ \ \ \ \ \ \ \ \ \  \ \ \ \ \ \ \   \forall_{i, j \in I}\forall_{\B A(X) \in \Fam(\D 1, \B X)}\bigg(\B \lambda_0(i) \subseteq \B \alpha_0(0) \subseteq \B \lambda_0(j) \ \& \ \mu \big(\B \lambda_0 (i)\big) = \mu \big(\B \lambda_0 (j)\big) \ \ \ \ \ \ \ \ \  \]
\[ \ \ \ \ \ \ \ \ \ \ \ \ \ \ \ \ \ \ \ \ \ \ \ \ \ \ \ \ \ \ \ \ \  \ \ \ \ \ \ \ \ \ \ \ \ \ \ \ \ \ \ \ \ \ \ \ \ \To \exists_{k \in I}\big(\B \alpha_0(0) = \B \lambda_0(k)\big)\bigg). \ \ \ \ \ \ \ \ \ \ \ \ \ \ \ \ \ \ \ \ \ \ \ \ \ \ \ \ \ \ \]
\end{definition}

Regarding the completeness conditions and the space $\C M(X)$, we show the following.

\begin{proposition}\label{prp: detachablecomplete}
Let $\C M(X) := (X, \delta_0 \D F(X, \D 2)(\B X), \mu^*_{x_0})$ be the measure space of 
complemented detachable subsets of $X$.\\[1mm]
\normalfont (i)
\itshape $\C M(X)$ satisfies condition $(\CM_1)$.\\[1mm]
\normalfont (ii)
\itshape The limited principle of omniscience $(\LPO)$ implies that $\C M(X)$ satisfies condition $(\CM_2)$.\\[1mm]
\normalfont (iii)
\itshape In general, $\C M(X)$ does not satisfy condition $(\CM_3)$.

\end{proposition}

\begin{proof}
(i) Let $f \in \D F (X, \D 2)$, let $\B B := (B^1, B^0)$ be a given complemented subset of
$X$ with $\B \alpha_0 (0) := \B B$, and let $\delta_0^1 (f) \subseteq B^1$ and $\delta_0^0 (f) 
\subseteq B^0$. Since $X = \delta_0^1 (f) \cup \delta_0^0 (f) \subseteq (B^1 \cup B^0) \subseteq X$, 
we get $B^1 \cup B^0 = X$, and hence $\B B = \B \delta_0 (\chi_{\B B})$.\\
(ii) Let $\alpha, \beta : \Nat \to \D F(X, \D 2)$, and $l \in (0, + \infty)$ such that 
\[ \forall_{m \in \Nat}\bigg(\bigcup_{n = 1}^m \B \delta_0 (\alpha_n) = \B \delta_0 (\beta_m)\bigg) 
\ \& \  \exists  \lim_{\mathsmaller{m \to + \infty}} \beta_m (x_0) \ \& \ \lim_{\mathsmaller{m \to + \infty}} 
\beta_m (x_0) = l. \]
The last conjunct is equivalent to $\exists_{m_0 \in \Nat}\forall_{m \geq m_0}\big(\beta_m (x_0) = l\big)$, and
since $\beta_m (x_0) \in \D 2$, we get $l \in \D 2$. For every $x \in X$ the sequence $n \mapsto 
\alpha_n (x)$ is in $\D F(\Nat, \D 2)$, hence by $(\LPO)$ we define the function $f$ from $X$ to $\D 2$ by the rule
\[ f(x) := \left\{ \begin{array}{ll}
                 1    &\mbox{, $\exists_{n \in \Nat}\big(\alpha_n (x) = 1\big)$}\\
                 0             &\mbox{, $\forall_{n \in \Nat}\big(\alpha_n (x) = 0\big)$.}
                 \end{array}
          \right.\] 
By the definition of interior union and intersection it is immediate to show that 
\[ \bigcup_{n \in \Nat} \B \delta_0 (\alpha_n) = \B \delta_0 (f)  \TOT
 \bigcup_{n \in \Nat} \delta_0^1 (\alpha_n) = \delta_0^1 (f) \ \& \ \bigcap_{n \in \Nat} 
 \delta_0^0 (\alpha_n) = \delta_0^0 (f).\]
It remains to show that $f(x_0) = l$. If $l = 0$, then 
$\exists_{m_0 \in \Nat}\forall_{m \geq m_0}\big(\beta_m (x_0) = 0\big)$, which implies that 
$\forall_{n \in \Nat}\big(\alpha_n (x_0) = 0\big) : \TOT f(x_0) := 0$. If $l = 1$, then
$\exists_{m_0 \in \Nat}\forall_{m \geq m_0}\big(\beta_m (x_0) = 1\big)$, which implies that 
$\exists_{n \in \{1, \ldots, m_0\}}\big(\alpha_n (x_0) = 1\big)  \To f(x_0) := 1$.\\
(iii) If $X := \D 3$, let $f : \D 3 \to \D 2$ be defined by $f(0) := 1, f(1) := 0 =: f(2)$ and
let $g \colon \D 3 \to \D 2$ be the constant function with value $1$. If 
$\B B := \big(\{0\}, \{1\}\big)$, then $\delta_0^1 (f) := \{0\} = B^1 \subseteq 
\delta_0^1 (g)$ and $\emptyset = \delta_0^0 (g) \subseteq B^0 \subseteq \delta_0^0 (f)$. If $x_0 := 0$, 
then $\mu_0 (f) = \mu_0 (g)$, but $\B B$ cannot ``pseudo-belong'' to $\C D(\D 3)$, since $B^1 \cup B^0$ 
is a proper subset of $\D 3$.
\end{proof}


\section{Real-valued partial functions}
\label{sec: realpartial}

We present here all facts on real-valued partial functions necessary to the definition of an
integration space within $\BST$ (Definition~\ref{def: intspace}).

\begin{definition}\label{def: realpfunction}
If $(X, =_X, \neq_X)$ is an inhabited set, we denote by $f_A := (A, i_A^X, f_A^{\Real})$ a 
real-valued partial function on $X$
\begin{center}
\begin{tikzpicture}

\node (E) at (0,0) {$A$};
\node[right=of E] (F) {$X$};
\node[below=of F] (A) {$\Real$.};

\draw[right hook->] (E)--(F) node [midway,above] {$i_A^X$};
\draw[->] (E)--(A) node [midway,left] {$f_A^{\Real} \ \ $};

\end{tikzpicture}
\end{center}
We say that $f_A$ is strongly extensional\index{strongly extensional partial function}, if 
$f_A^{\Real}$ is strongly extensional, where $A$ is equipped with its canonical inequality as a 
subset of $X$ i.e., for every $a, a{'} \in A$
\[f_A^{\Real}(a) \neq_{\Real} f_A^{\Real}(a{'}) \To i_A^X(a) \neq_X i_A^X(a{'}).\]
Let $\C F(X) := \C F(X, \Real)$\index{$\C F(X)$} be the class of partial functions 
from $X$ to $\Real$, and $\C F^{\se}(X)$\index{$\C F^{\se}(X)$} the class of strongly extensional 
partial functions from $(X =_X, \neq_X)$ to $(\Real, =_{\Real}, \neq_{\Real})$.
\end{definition}


\begin{definition}\label{def: rpfunoper}
Let $f_A := (A, i_A^X, f_A^{\Real}), f_B := (B, i_B^X, f_B^{\Real})$ in $\C F(X)$
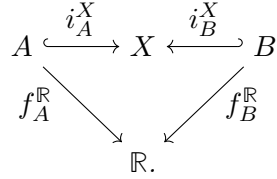
\begin{center}
\begin{tikzpicture}

\node (E) at (0,0) {$A$};
\node[right=of E] (B) {$X$};
\node[right=of B] (F) {$B$};
\node[below=of B] (C) {$\Real$.};

\draw[right hook->] (E)--(B) node [midway,above] {$i_A^X$};
\draw[left hook->] (F)--(B) node [midway,above] {$i_B^X$};
\draw[->] (E)--(C) node [midway,left] {$f_A^{\Real} \ \ $};
\draw[->] (F)--(C) node [midway,right] {$ \ f_B^{\Real}$};

\end{tikzpicture}
\end{center} 
If $\lambda \in \Real$, let $\lambda 
f_A := (A, i_A^X, \lambda f_A^{\Real}) \in \C F(X)$\index{$\lambda f_A$},\index{$f_A \ \square \ f_B$} and 
$ f_A \ \square \ f_B := \big(A \cap B, i_{A \cap B}^X, \big(f_A^{\Real} \ \square 
\ f_B^{\Real}\big)_{A \cap B}^{\Real}\big)$, 
\[ (f_A^{\Real} \ \square \ f_B^{\Real}\big)_{A \cap B}^{\Real} := f_A^{\Real}(a) \ 
\square \ f_B^{\Real}(b); \ \ \ \ (a, b) \in A \cap B, \  \ \square \in \{+, \cdot, \wedge, \vee\}.\]
\end{definition}

The operation $(f_A^{\Real} \ \square \ f_B^{\Real}\big)_{A \cap B}^{\Real} \colon A \cap B \sto \Real$ 
is a function, as if $(a, b) =_{A \cap B} (a{'}, b{'}) :\TOT i_A^X =_X i_A^X(a{'})$,
and since $i_B^X(b) =_X i_A^X(a)$ and $i_B^X(b{'}) =_X i_A^X(a{'})$, we get $a =_A a{'}$, hence 
$f_A^{\Real}(a) =_{\Real} f_A^{\Real}(a{'})$, and $b =_B b{'}$, hence $f_B^{\Real}(b) =_{\Real} f_B^{\Real}(b{'})$. 
If $\lambda$ denotes also the constant function $\lambda \in \Real$ on $X$
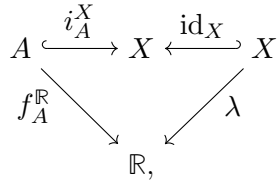
\begin{center}
\begin{tikzpicture}

\node (E) at (0,0) {$A$};
\node[right=of E] (B) {$X$};
\node[right=of B] (F) {$X$};
\node[below=of B] (C) {$\Real$,};

\draw[right hook->] (E)--(B) node [midway,above] {$i_A^X$};
\draw[left hook->] (F)--(B) node [midway,above] {$\id_X$};
\draw[->] (E)--(C) node [midway,left] {$f_A^{\Real} \ \ $};
\draw[->] (F)--(C) node [midway,right] {$ \ \lambda$};

\end{tikzpicture}
\end{center} 
we get as a special case the partial function $f_A \wedge \lambda := \big(A \cap X, i_{A \cap X}^X, 
\big(f_A^{\Real} \wedge \lambda\big)_{A \cap X}^{\Real}\big)$, where $A \cap X := \{ (a, x) \in A \times X 
\mid i_A^X(a) =_X x\}$, $i_{A \cap X}(a, x) := i_A^X(a)$, and $\big(f_A^{\Real} \wedge \lambda\big)_{A \cap X}^{\Real}(a,x) := f_A^{\Real}(a) \wedge \lambda(x) := f_A^{\Real}(a) \wedge \lambda$, for every $(a, x) \in A \cap X$.
By Definition~\ref{def: famofpartial}, if $\Lambda(X,\Real) := (\lambda_0, \C E^X, \lambda_1, \C P^{\Real})
\in \Fam(I, X, \Real)$, 
if $(i, j) \in D(I)$, the following diagram commutes
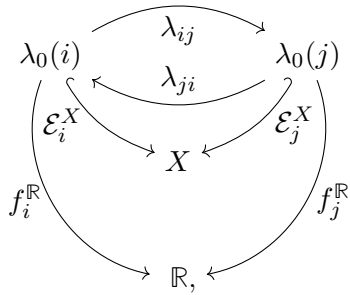
\begin{center}
\begin{tikzpicture}

\node (E) at (0,0) {$\lambda_0(i)$};
\node[right=of E] (B) {};
\node[right=of B] (F) {$\lambda_0(j)$};
\node[below=of B] (A) {$X$};
\node[below=of A] (C) {$ \ \Real$,};

\draw[->,bend left] (E) to node [midway,below] {$\lambda_{ij}$} (F);
\draw[->,bend left] (F) to node [midway,above] {$\lambda_{ji}$} (E);
\draw[->,bend right=50] (E) to node [midway,left] {$f_i^{\Real}$} (C);
\draw[->,bend left=50] (F) to node [midway,right] {$f_j^{\Real}$} (C);
\draw[right hook->,bend right=20] (E) to node [midway,left] {$ \ \C E_i^X \ $} (A);
\draw[left hook->,bend left=20] (F) to node [midway,right] {$\ \C E_j^X \ $} (A);

\end{tikzpicture}
\end{center} 
\[ f_i := \big(\lambda_0(i), \C E_i^X, f_i^{\Real}\big) \in \C F(X), \ \ \ \
f_i^{\Real} := \C P^{\Real}_i \colon \lambda_0(i) \to \Real; \ \ \ \ i \in I.\]
If $f_i^{\Real}$ is strongly extensional, then, for every $u, w \in \lambda_0(i)$, we have that
$f_i^{\Real}(u) \neq_{\Real} f_i^{\Real}(w) \To \C E_i^X(u) \neq_X \C E_i^X(w)$.
As in Definition~\ref{def: subfamsubsets}, If $\kappa \colon \Nat^+ \to I$, the family 
$\Lambda(X, \Real) \circ \kappa := \big(\lambda_0 \circ \kappa, \C E^X \circ \kappa, \lambda_1
\circ \kappa, \C P^{\Real} \circ \kappa\big) \in \Fam(\Nat^+, X, \Real)$ is the $\kappa$-subsequence
of $\Lambda(X, \Real)$, where 
\[ (\lambda_0 \circ \kappa)(n) := \lambda_0(\kappa(n)), \ \ \ \ \big(\C E^X \circ \kappa\big)_n :=
\C E_{\kappa(n)}^X, \ \  \ \ (\lambda_1 \circ \kappa)(n, n) := \lambda_{\kappa(n)\kappa(n)} :=
\id_{\lambda_0(\kappa(n))}, \]
\[  \big(\C P^{\Real} \circ \kappa\big)_n := \C P_{\kappa(n)}^{\Real} := f_{\kappa(n)}^{\Real}; \ \ \ \ n \in \Nat^+.\]
Let $\lambda_0I (X, \Real)$\index{$\lambda_0I (X, \Real)$} be the totality $I$, and we write $f_i
\in \lambda_0I(X, \Real)$, instead of $i \in I$, as we define 
\[ i =_{\lambda_0I (X, \Real)} j :\TOT f_i =_{\C F(X)} f_j.\]
If we consider the intersection $\bigcap_{n \in \Nat^+}(\lambda_0 \circ \kappa)(n) := 
\bigcap_{n \in \Nat^+}\lambda_0(\kappa(n))$, 
by Definition~\ref{def: intfamilyofsubsets}
\[ \Phi \colon \bigcap_{n \in \Nat^+}\lambda_0(\kappa(n))  :\TOT \Phi \colon 
\bigcurlywedge_{n \in \Nat^+}\lambda_0(\kappa(n)) \ \ \& \ \ 
\forall_{n, m \in \Nat^+}\big(\C E_{\kappa(n)}^X(\Phi_n) =_X \C E_{\kappa(m)}^X(\Phi_m)\big), \]
\[ \Phi =_{\mathsmaller{\bigcap_{n \in \Nat^+}\lambda_0(\kappa(n))}} \Theta 
:\TOT \C E_{\kappa(1)}^X(\Phi_1) =_X \C E_{\kappa(1)}^X(\Theta_1), \]
\[ e_{\mathsmaller{\bigcap}}^{\Lambda(X, \Real) \circ \kappa} \colon 
\bigcap_{n \in \Nat^+}\lambda_0(\kappa(n)) \eto X, \ \ \ \  
e_{\mathsmaller{\bigcap}}^{\Lambda(X, \Real) \circ \kappa}(\Phi) := 
\big(\C E^X \circ \kappa)_1 (\Phi_1) := \C E_{\kappa(1)}^X(\Phi_1).\]

\begin{definition}\label{def: sigmapfun}
Let $\Lambda(X,\Real) := (\lambda_0, \C E^X, \lambda_1, \C P^{\Real}) \in 
\Fam(I, X, \Real)$, $\kappa \colon \Nat^+ \to I$, and 
$\Lambda(X, \Real) \circ \kappa$ the $\kappa$-subsequence of $\Lambda(X, \Real)$. 
If $(A, i_A^{\mathsmaller{\bigcap}}) \subseteq \bigcap_{n \in \Nat^+}\lambda_0(\kappa(n))$, we define the function
\[ \sum_{n \in \Nat^+}^Af_{\kappa(n)}^{\Real} \colon A \to \Real, \ \ \ \
\bigg(\sum_{n \in \Nat^+}^Af_{\kappa(n)}^{\Real}\bigg)(a)
:= \sum_{n \in \Nat^+}f_{\kappa(n)}^{\Real}\bigg(\big[i_A^{\mathsmaller{\bigcap}}(a)\big]_n\bigg); \ \ \ \ a \in A,\]
under the assumption that the series on the right converge in $\Real$, for every $a \in \Real$. 
\end{definition}

In the special case 
$\bigcap_{n \in \Nat^+}\lambda_0(\kappa(n)), id_{\mathsmaller{\bigcap_{n \in \Nat^+}\lambda_0(\kappa(n))}})
\subseteq \bigcap_{n \in \Nat^+}\lambda_0(\kappa(n))$, we get the function
\[ \sum_{n \in \Nat^+}^{\mathsmaller{\bigcap}}f_{\kappa(n)}^{\Real} \colon 
\bigcap_{n \in \Nat^+}\lambda_0(\kappa(n)) \to \Real, \ \ \ \ 
\bigg(\sum_{n \in \Nat^+}^{\mathsmaller{\bigcap}}f_{\kappa(n)}^{\Real}\bigg)(\Phi)
:= \sum_{n \in \Nat^+}f_{\kappa(n)}^{\Real}\big(\Phi_n\big); \ \ \ \ 
\Phi \in \bigcap_{n \in \Nat^+}\lambda_0(\kappa(n)),\]
under the assumption on the convergence of the corresponding series.

\begin{proposition}\label{prp: sextpf1}
If in Definition~\ref{def: sigmapfun} the partial functions $f_{\kappa(n)} := 
\big(\lambda_0(\kappa(n)), \C E_{\kappa(n)}^X, f_{\kappa(n)}^{\Real}\big)$ are strongly extensional, 
for every $n \in \Nat^+$, then the real-valued partial function
\begin{center}
\begin{tikzpicture}

\node (E) at (0,0) {$A$};
\node[right=of E] (G) {$\bigcap_{n \in \Nat^+}\lambda_0(\kappa(n))$};
\node[right=of G] (H) {};
\node[right=of H] (F) {$X$};
\node[below=of F] (A) {$\Real$};

\draw[right hook->] (E)--(G) node [midway,above] {$i_A^{\mathsmaller{\bigcap}}$};
\draw[right hook->] (G)--(F) node [midway,above] {$e_{\mathsmaller{\bigcap}}^{\Lambda(X, \Real) \circ \kappa} $};
\draw[right hook->,bend right=20] (E) to node [midway,above] {$ \ \ \ \ \ 
\sum_{n \in \Nat^+}^Af_{\kappa(n)}^{\Real}  $} (A);

\end{tikzpicture}
\end{center}
\[ f_A := \bigg(A, \ e_{\mathsmaller{\bigcap}}^{\Lambda(X, \Real) \circ \kappa}
\circ i_A^{\mathsmaller{\bigcap}}, \ \sum_{n \in \Nat^+}^Af_{\kappa(n)}^{\Real} \bigg)\]
is strongly extensional.
\end{proposition}

\begin{proof}
Let $a, a{'} \in A$ such that 
\[ \bigg(\sum_{n \in \Nat^+}^Af_{\kappa(n)}^{\Real}\bigg)(a) := l \neq_{\Real} l{'} := 
\bigg(\sum_{n \in \Nat^+}^Af_{\kappa(n)}^{\Real}\bigg)(a{'}).\]
There is $N \in \Nat^+$ such that, if $\varepsilon := |l - l{'}| >0$, then 
\begin{align*}
\varepsilon & \leq \bigg| l - \sum_{n = 1}^Nf_{\kappa(n)}^{\Real}\bigg(\big[i_A^{\mathsmaller{\bigcap}}(a)\big]_n\bigg) 
\bigg| + 
\bigg| \sum_{n = 1}^Nf_{\kappa(n)}^{\Real}\bigg(\big[i_A^{\mathsmaller{\bigcap}}(a)\big]_n\bigg) -
\sum_{n = 1}^Nf_{\kappa(n)}^{\Real}\bigg(\big[i_A^{\mathsmaller{\bigcap}}(a{'})\big]_n\bigg) \bigg| \\
& \ \ \ \ \  \ \ \ \ \ \ \ \ \ \ \ \ \ \ \ \ \ \ \ \ \ \ \ \ \ \ \  +
\bigg| \sum_{n = 1}^Nf_{\kappa(n)}^{\Real}\bigg(\big[i_A^{\mathsmaller{\bigcap}}(a{'})\big]_n\bigg) - l{'}\bigg|\\
& \leq \frac{\varepsilon}{4} + 
\bigg| \sum_{n = 1}^Nf_{\kappa(n)}^{\Real}\bigg(\big[i_A^{\mathsmaller{\bigcap}}(a)\big]_n\bigg) -
\sum_{n = 1}^Nf_{\kappa(n)}^{\Real}\bigg(\big[i_A^{\mathsmaller{\bigcap}}(a{'})\big]_n\bigg) \bigg| +
\frac{\varepsilon}{4} \To
\end{align*}
\begin{align*}
0 & < \bigg| \sum_{n = 1}^Nf_{\kappa(n)}^{\Real}\bigg(\big[i_A^{\mathsmaller{\bigcap}}(a)\big]_n\bigg) -
\sum_{n = 1}^Nf_{\kappa(n)}^{\Real}\bigg(\big[i_A^{\mathsmaller{\bigcap}}(a{'})\big]_n\bigg) \bigg|\\
& = \bigg| \sum_{n = 1}^N \bigg[f_{\kappa(n)}^{\Real}\bigg(\big[i_A^{\mathsmaller{\bigcap}}(a)\big]_n\bigg) -
f_{\kappa(n)}^{\Real}\bigg(\big[i_A^{\mathsmaller{\bigcap}}(a{'})\big]_n\bigg)\bigg] \bigg|\\
& \leq \sum_{n = 1}^N \bigg| f_{\kappa(n)}^{\Real}\bigg(\big[i_A^{\mathsmaller{\bigcap}}(a)\big]_n\bigg) -
f_{\kappa(n)}^{\Real}\bigg(\big[i_A^{\mathsmaller{\bigcap}}(a{'})\big]_n\bigg)\bigg|.
\end{align*}
By the property of positive real numbers $x + y > 0 \To [x > 0 \vee y > 0$ (see~\cite{BR87}, p.~13), then
\[ 0 < \bigg| f_{\kappa(n)}^{\Real}\bigg(\big[i_A^{\mathsmaller{\bigcap}}(a)\big]_n\bigg) -
f_{\kappa(n)}^{\Real}\bigg(\big[i_A^{\mathsmaller{\bigcap}}(a{'})\big]_n\bigg)\bigg|,\]
for some $n \in \Nat^+$ with $1 \leq n \leq N$. Since $f_{\kappa(n)}^{\Real}$ is strongly extensional,
by Definition~\ref{def: realpfunction}
\[ \C E_{\kappa(1)}^X\bigg(\big[i_A^{\mathsmaller{\bigcap}}(a)\big]_1\bigg) =_X 
\C E_{\kappa(n)}^X\bigg(\big[i_A^{\mathsmaller{\bigcap}}(a)\big]_n\bigg) \neq_X
\C E_{\kappa(n)}^X\bigg(\big[i_A^{\mathsmaller{\bigcap}}(a{'})\big]_n\bigg) =_X 
\C E_{\kappa(1)}^X\bigg(\big[i_A^{\mathsmaller{\bigcap}}(a{'})\big]_1\bigg),\]
which is the required canonical equality of $a, a{'}$ in $A$ as a subset of $X$.
\end{proof}

Clearly, we have that (see Definition~\ref{def: partialfunction})
\begin{center}
\begin{tikzpicture}

\node (E) at (0,0) {$A$};
\node[right=of E] (B) {};
\node[right=of B] (F) {$\bigcap_{n \in \Nat^+}\lambda_0(\kappa(n))$};
\node[below=of B] (A) {$X$};
\node[below=of A] (C) {$ \ \Real$};

\draw[right hook->] (E) to node [midway,above] {$i_A^{\mathsmaller{\bigcap}}$} (F);
\draw[->,bend right=70] (E) to node [midway,left] {$\sum_{n \in \Nat^+}^Af_{\kappa(n)}^{\Real} $} (C);
\draw[->,bend left=70] (F) to node [midway,left] {$\sum_{n \in \Nat^+}^{\mathsmaller{\bigcap}}f_{\kappa(n)}^{\Real} \ \ $} (C);
\draw[right hook->,bend right=40] (E) to node [midway,right] {$ \ \mathsmaller{e_{\mathsmaller{\bigcap}}^{\Lambda(X, \Real) \circ \kappa} \circ i_A^{\mathsmaller{\bigcap}}}$} (A);
\draw[left hook->,bend left=40] (F) to node [midway,right] {$\mathsmaller{e_{\mathsmaller{\bigcap}}^{\Lambda(X, \Real) \circ \kappa}}$} (A);

\end{tikzpicture}
\end{center} 

\[ \bigg(A, \ e_{\mathsmaller{\bigcap}}^{\Lambda(X, \Real) \circ \kappa} \circ
i_A^{\mathsmaller{\bigcap}}, \ \sum_{n \in \Nat^+}^Af_{\kappa(n)}^{\Real} \bigg) 
\leq \bigg(\bigcap_{n \in \Nat^+}\lambda_0(\kappa(n)), 
\ e_{\mathsmaller{\bigcap}}^{\Lambda(X, \Real) \circ \kappa}, \
\sum_{n \in \Nat^+}^{\mathsmaller{\bigcap}}f_{\kappa(n)}^{\Real} \bigg).\]

\section{Integration and pre-integration spaces}
\label{sec: preintegration}

Next we reformulate predicatively the Bishop-Cheng definition of an integration space 
(see Note~\ref{not: BCis} to compare it with the original definition).

\begin{definition}[Integration space within $\BST$]\label{def: intspace}
Let $(X, =_X, \neq_X)$ be an inhabited set, $\Lambda(X,\Real) := 
(\lambda_0, \C E^X, \lambda_1, \C P^{\Real}) \in \Fam(I, X, \Real)$, such that $f_i := 
\big(\lambda_0(i), \C E_i^X, f_i^{\Real}\big)$ is strongly extensional, for every $i \in I$,
let $\lambda_0I (X, \Real)$ be the totality $I$, equipped with the equality $i =_{\lambda_0I (X, \Real)} j 
:\TOT f_i =_{\C F(X)} f_j$, for every $i \in I$, and let a mapping 
\[ \int \colon \lambda_0I (X, \Real) \to \Real, \ \ \ \ f_i \mapsto \int f_i; \ \ \ \ i \in I,\]
such that the following conditions hold:
\vspace{-2mm}
\[(\IS_1) \ \  \ \ \ \ \ \ \ \ \ \ \ \ \ \ \ \ \ \ \ \  \ \ \forall_{i \in I}\forall_{a \in \Real}\exists_{j \in I}\bigg(a f_i =_{\C F(X)} f_j \ \ \& \ \ \int f_j =_{\Real} a \int f_i\bigg). \ \ \ \ \ \ \ \ \ \ \ \ \ \ \ \ \ \ \ \ \ \ \ \ \ \ \ \  \]
\[(\IS_2) \ \  \ \ \ \ \ \ \ \ \ \ \ \ \ \ \ \ \ \ \ \  \ \ \forall_{i, j \in I}\exists_{k \in I}\bigg(f_i + f_j =_{\C F(X)} f_k \  \ \& \ \ \int f_k =_{\Real} \int f_i\ + \int f_j\bigg). \ \ \ \ \ \ \ \ \ \ \ \ \ \ \ \ \ \ \ \ \ \ \ \ \ \ \ \  \]
\[(\IS_3) \ \  \ \ \ \ \ \ \ \ \ \ \ \ \ \ \ \ \ \ \ \ \ \ \ \ \ \ \ \ \ \ \ \ \ \ \ \  \ \ \forall_{i \in I}\exists_{j \in I}\big( |f_i| =_{\C F(X)} f_j \big). \ \ \ \ \ \ \ \ \ \ \ \ \ \ \ \ \ \ \ \ \ \ \ \ \ \ \ \ \ \ \ \ \ \ \ \ \ \ \ \ \ \ \ \ \ \ \   \]
\[(\IS_4) \ \  \ \ \ \ \ \ \ \ \ \ \ \ \ \ \ \ \ \ \ \ \ \ \ \ \ \ \ \ \ \ \ \ \ \ \ \  \ \ \forall_{i \in I}\exists_{j \in I}\big( f_i \wedge 1 =_{\C F(X)} f_j \big). \ \ \ \ \ \ \ \ \ \ \ \ \ \ \ \ \ \ \ \ \ \ \ \ \ \ \ \ \ \ \ \ \ \ \ \ \ \ \ \ \ \ \ \ \ \ \   \] 
\[(\IS_5) \  \ \ \ \  \ \ \ \ \ \ \ \ \  \ \ \ \ \ \ \ \  \forall_{i \in I}\forall_{\kappa \in \D F(\Nat^+, I)}\bigg\{\bigg[\sum_{n \in \Nat^+}\int f_{\kappa(n)} \in \Real \  \ \& \ \sum_{n \in \Nat^+}\int f_{\kappa(n)} < \int f_i\bigg] \To  \ \ \ \ \ \ \ \ \ \ \ \ \ \ \ \ \ \ \ \ \ \]
\[ \exists_{(\Phi, u) \ \in \ \big(\bigcap_{n \in \Nat^+}\lambda_0(\kappa(n)) \big)\ \cap \ \lambda_0(i)}\bigg( \bigg(\sum_{n \in \Nat^+}^{\mathsmaller{\bigcap}}f_{\kappa(n)}\bigg)(\Phi) := \sum_{n \in \Nat^+}f_{\kappa(n)}^{\Real}(\Phi_n) \in \Real \  \ \& \ \]
\[ \ \ \ \ \ \ \ \ \ \ \ \ \ \ \ \ \ \ \  \sum_{n \in \Nat^+}f_{\kappa(n)}^{\Real}(\Phi_n) < f_i^{\Real}(u). \bigg)\bigg\}\]
\[(\IS_6) \ \  \ \ \ \ \ \ \ \ \ \ \ \ \ \ \ \ \ \ \ \ \ \ \ \ \ \ \ \ \ \ \ \ \ \ \ \  \ \ \ \ \ \ \  \ \ \exists_{i \in I}\bigg( \int  f_i =_{\Real} 1 \bigg). \ \ \ \ \ \ \  \ \ \ \ \ \ \ \ \ \ \ \ \ \ \ \ \ \ \ \ \ \ \ \ \ \ \ \ \ \ \ \ \ \ \ \ \ \ \ \ \ \ \ \ \ \ \   \] 
\[(\IS_7) \  \ \ \ \ \ \ \ \ \ \ \ \ \ \ \ \ \ \ \ \ \  \ \ \ \  \forall_{i \in I}\forall_{\alpha \in \D F(\Nat^+, I)}\bigg(\forall_{n \in \Nat^+}\bigg(n \bigg(\frac{1}{n} f_i \ \wedge \ 1\bigg) =_{\C F(X)} f_{\alpha(n)}\bigg) \To \ \ \ \ \ \ \ \ \ \ \ \ \ \ \ \ \ \ \ \ \ \ \ \  \]
\[ \ \ \ \ \ \ \ \ \ \ \ \ \ \ \ \ \ \lim_{\mathsmaller{n \longrightarrow +\infty}} \int f_{\alpha(n)} \in \Real \ \ \& \ 
\lim_{\mathsmaller{n \longrightarrow +\infty}} \int f_{\alpha(n)} =_{\Real} \int f_i\bigg). \ \ \ \ \ \ \ \ \ \ \ \ \]
\[(\IS_8) \ \  \ \ \ \ \ \ \ \ \ \ \ \ \ \ \ \ \ \ \ \  \ \ \ \ \ \ \ \ \ \ \  \forall_{i \in I}\forall_{\alpha \in \D F(\Nat^+, I)}\bigg(\forall_{n \in \Nat^+}\bigg(\frac{1}{n}\big(n |f_i| \ \wedge \ 1\big) =_{\C F(X)} f_{\alpha(n)}\bigg) \To \ \ \ \ \ \ \ \ \ \ \ \ \ \ \ \ \ \ \ \ \ \ \ \]
\[ \ \ \ \ \ \ \ \ \ \ \ \ \ \ \ \ \ \ \ \ \ \ \ \ \ \ \ \ \ \ \ \ \ \ \ \ \ \ \lim_{\mathsmaller{n \longrightarrow +\infty}} \int f_{\alpha(n)} \in \Real  \ \ \& \  \lim_{\mathsmaller{n \longrightarrow +\infty}} \int f_{\alpha(n)} =_{\Real} 0\bigg). \ \ \ \ \ \ \ \ \ \ \ \ \ \ \ \ \ \ \ \ \ \ \ \ \ \ \ \  \]
We call the triplet $\C L := (X, \lambda_0I (X, \Real), \int)$ an integration space\index{integration space}.

\end{definition}

In the formulation of $(\IS_5)$ we have that 
\[ \bigg(\bigcap_{n \in \Nat^+}\lambda_0(\kappa(n))\bigg) \cap  \lambda_0(i) 
:= \bigg\{(\Phi, u) \in \bigg(\bigcap_{n \in \Nat^+}\lambda_0(\kappa(n))\bigg) \times
\lambda_0(i) \ \mathlarger{\mathlarger{\mid}} \ \C E^X_{\kappa(1)}(\Phi) =_X \C E_i^X(u)\bigg\},\] 
\begin{center}
\begin{tikzpicture}

\node (E) at (0,0) {$\lambda_0(\kappa(n))$};
\node[right=of E] (B) {$X$};
\node[right=of B] (F) {$\lambda_0(i)$};
\node[below=of B] (C) {$\Real$.};

\draw[right hook->] (E)--(B) node [midway,above] {$\C E_{\kappa(n)}^X$};
\draw[left hook->] (F)--(B) node [midway,above] {$\C E_i^X$};
\draw[->] (E)--(C) node [midway,left] {$f_{\kappa(n)}^{\Real} \ \ $};
\draw[->] (F)--(C) node [midway,right] {$ \ f_i^{\Real}$};

\end{tikzpicture}
\end{center} 
If, for every $a \in \Real$, such that $a > 0$, and every $i \in I$, we define
\[ f_i \ \wedge \ a := a \bigg(\frac{1}{a} f_i \ \wedge \ 1\bigg),\]
the formulation of $(\IS_7)$ and $(\IS_8)$ becomes, respectively,
\[ \forall_{i \in I}\forall_{\alpha 
\in \D F(\Nat^+, I)}\bigg(\forall_{n \in \Nat^+}\big(f_i \ \wedge \ n =_{\C F(X)} 
f_{\alpha(n)}\big) \To \lim_{\mathsmaller{n \longrightarrow +\infty}} \int f_{\alpha(n)}
=_{\Real} \int f_i\bigg),  \]
\[ \forall_{i \in I}\forall_{\alpha \in \D F(\Nat^+, I)}\bigg(\forall_{n \in \Nat^+}\bigg(|f_i| \ \wedge \
\frac{1}{n}\bigg) =_{\C F(X)} f_{\alpha(n)}\bigg) \To \lim_{\mathsmaller{n \longrightarrow +\infty}}
\int f_{\alpha(n)} =_{\Real} 0\bigg),  \]
where, for simplicity, we skip to mention he existence of the corresponding limits in $\Real$.
We also quantify over $\D F(\Nat^+, I)$, in order to avoid the use of countable choice. If we had written 
in its premise the formula $\forall_{n \in \Nat^+}\exists_{j \in I}\big(f_i \ \wedge \ n =_{\C F(X)} f_j\big)$, 
we would need countable choice ($\Nat{-}I$) to generate a sequence in $I$ to describe the limit of the
corresponding integrals. Moreover, by $(\IS_1)$, $(\IS_3)$ and $(\IS_4)$, and the definition of $f \wedge a $ 
above we get  $\forall_{i \in I}\forall_{a \in \Real^+}\exists_{j \in I}\big(f_i \wedge  a =_{\C F(X)} f_j\big)$.

\begin{definition}[Pre-integration space within $\BST$]\label{def: preintspace}
Let $(X, =_X, \neq_X)$ be an inhabited set, and let the set $(I, =_I)$ be equipped with 
operations $\cdot_a \colon I \sto I$, for every $a \in \Real$, $+ \colon I \times I \sto I$, $|.| \colon I \sto I$, 
and $\wedge_1 \colon I \sto I$, where 
\[ \cdot_a(i) := a \cdot i, \ \ \ \ +(i,j) := i + j, \ \ \ \ |.|(i) := |i|; \ \ \ \ i \in I, \ a \in \Real.\]  
Let also the operation $\wedge_a \colon I \sto I$, defined by the previous operations with the rule
\[ \wedge_a := \cdot_a \circ \wedge_1 \circ \cdot_{a^{-1}}; \ \ \ \ a \in \Real \ \& \ a > 0.\]
Let $\Lambda(X,\Real) := (\lambda_0, \C E^X, \lambda_1, \C P^{\Real}) \in \Set(I, X, \Real)$ i.e., 
$f_i =_{\C F(X)} f_j \To i =_I j$, for every $i, j \in I$, and $f_i := \big(\lambda_0(i),
\C E_i^X, f_i^{\Real}\big)$ is strongly extensional, for every $i \in I$. Let also a mapping 
\[ \int \colon I \to \Real, \ \ \ \ i \mapsto \int i; \ \ \ \ i \in I,\]
such that the following conditions hold:
\vspace{-2mm}
\[(\PIS_1) \ \  \ \ \ \ \ \ \ \ \ \ \ \ \ \ \ \ \ \ \ \  \ \ \forall_{i \in I}\forall_{a \in \Real}\bigg(a f_i =_{\C F(X)} f_{a \cdot i} \ \ \& \ \ \int a \cdot i =_{\Real} a \int i\bigg). \ \ \ \ \ \ \ \ \ \ \ \ \ \ \ \ \ \ \ \ \ \ \ \ \ \ \ \  \]
\[(\PIS_2) \ \  \ \ \ \ \ \ \ \ \ \ \ \ \ \ \ \ \ \ \ \  \ \ \forall_{i, j \in I}\bigg(f_i + f_j =_{\C F(X)} f_{i + j} \  \ \& \ \ \int (i + j) =_{\Real} \int i \ + \int j\bigg). \ \ \ \ \ \ \ \ \ \ \ \ \ \ \ \ \ \ \ \ \ \ \ \ \ \ \ \  \]
\[(\PIS_3) \ \  \ \ \ \ \ \ \ \ \ \ \ \ \ \ \ \ \ \ \ \ \ \ \ \ \ \ \ \ \ \ \ \ \ \ \ \  \ \ \forall_{i \in I}\big( |f_i| =_{\C F(X)} f_{|i|} \big). \ \ \ \ \ \ \ \ \ \ \ \ \ \ \ \ \ \ \ \ \ \ \ \ \ \ \ \ \ \ \ \ \ \ \ \ \ \ \ \ \ \ \ \ \ \ \   \]
\[(\PIS_4) \ \  \ \ \ \ \ \ \ \ \ \ \ \ \ \ \ \ \ \ \ \ \ \ \ \ \ \ \ \ \ \ \ \ \ \ \ \  \ \ \forall_{i \in I}\big( f_i \wedge 1 =_{\C F(X)} f_{\wedge_1(i)} \big). \ \ \ \ \ \ \ \ \ \ \ \ \ \ \ \ \ \ \ \ \ \ \ \ \ \ \ \ \ \ \ \ \ \ \ \ \ \ \ \ \ \ \ \ \ \ \   \] 
\[(\PIS_5) \  \ \ \ \  \ \ \ \ \ \ \ \ \  \ \ \ \ \ \ \ \  \forall_{i \in I}\forall_{\kappa \in \D F(\Nat^+, I)}\bigg\{\bigg[\sum_{n \in \Nat^+}\int \kappa(n) \in \Real \  \ \& \ \sum_{n \in \Nat^+}\int \kappa(n) < \int i\bigg] \To  \ \ \ \ \ \ \ \ \ \ \ \ \ \ \ \ \ \ \ \ \ \]
\[ \exists_{(\Phi, u) \ \in \ \big(\bigcap_{n \in \Nat^+}\lambda_0(\kappa(n)) \big)\ \cap \ \lambda_0(i)}\bigg( \bigg(\sum_{n \in \Nat^+}^{\mathsmaller{\bigcap}}f_{\kappa(n)}\bigg)(\Phi) := \sum_{n \in \Nat^+}f_{\kappa(n)}^{\Real}(\Phi_n) \in \Real \  \ \& \ \]
\[ \ \ \ \ \ \ \ \ \ \ \ \ \ \ \ \ \ \ \  \sum_{n \in \Nat^+}f_{\kappa(n)}^{\Real}(\Phi_n) < f_i^{\Real}(u). \bigg)\bigg\}\]
\[(\PIS_6) \ \  \ \ \ \ \ \ \ \ \ \ \ \ \ \ \ \ \ \ \ \ \ \ \ \ \ \ \ \ \ \ \ \ \ \ \ \  \ \ \ \ \ \ \  \ \ \exists_{i \in I}\bigg( \int  i =_{\Real} 1 \bigg). \ \ \ \ \ \ \  \ \ \ \ \ \ \ \ \ \ \ \ \ \ \ \ \ \ \ \ \ \ \ \ \ \ \ \ \ \ \ \ \ \ \ \ \ \ \ \ \ \ \ \ \ \ \   \] 
\[(\PIS_7) \ \  \ \ \ \ \ \ \ \ \ \ \ \ \ \ \ \ \ \ \ \  \ \   \forall_{i \in I}\bigg(\lim_{\mathsmaller{n \longrightarrow +\infty}} \int \wedge_n(i) \in \Real \ \ \& \ \lim_{\mathsmaller{n \longrightarrow +\infty}} \int \wedge_n(i) =_{\Real} \int i\bigg). \  \ \ \ \ \ \ \ \ \ \ \  \ \ \ \ \ \ \ \ \ \ \]
\[(\PIS_8) \ \  \ \ \ \ \ \ \ \ \ \ \ \ \ \ \ \ \ \ \ \  \ \   \forall_{i \in I}\bigg(\lim_{\mathsmaller{n \longrightarrow +\infty}} \int \wedge_{\frac{1}{n}}(|i|) \in \Real \ \ \& \ \lim_{\mathsmaller{n \longrightarrow +\infty}} \int \wedge_{\frac{1}{n}}(|i|) =_{\Real} 0\bigg). \ \ \ \ \ \ \ \ \ \ \ \ \ \ \ \ \ \ \ \   \]
We call the triplet $\C L_0 := (X, I, \int)$ a pre-integration space\index{pre-integration space}.

\end{definition}

All the operations on $I$ defined above are functions. E.g., since $\Lambda(X, \Real)  \in \Set(I, X, \Real)$, 
\[ i =_I i{'} \To f_i =_{\C F(X)} f_{i{'}} \To a f_i =_{\C F(X)} a f_{i{'}} \To f_{a \cdot i} 
=_{\C F(X)} f_{a{'} \cdot i} \To a \cdot i =_I a{'} \cdot i.\] 
It is immediate to see that a pre-integration space induces an integration space, if 
\[ \forall_{i \in I}\forall_{a \in \Real^+}\big(f_i \ \wedge \ a  =_{\C F(X)} f_j \To \wedge_a(i) =_I j\big),\]
and hence $(\PIS_7)$ and $(\PIS_8)$ imply $(\IS_7)$ and $(\IS_8)$, respectively, with the integral 
\[ \int^* f_i := \int i; \ \ \ \ i \in I.\]
The notion of a pre-integration space is simpler than that of an integration space, and also closer
to the Bishop-Cheng notion of an integration space. One could say that a pre-integration space is the 
``right'' notion of integration space within $\BST$. In~\cite{BC72}, p.~52 Bishop and Cheng  formulate
the non-trivial theorem that a measure space induces the integration space of the corresponding simple 
functions (see also~\cite{BB85}, p.~285). In~\cite{Ze19} Zeuner interpreting the various constructions 
of Bishop and Cheng into the framework of pre-measure and pre-integration
spaces\footnote{The notion of pre-measure space used in~\cite{Ze19}, which is a bit 
different from the one included here, is an appropriate copy of Bishop's definition of a measure space given 
in~\cite{Bi67} (see Note~\ref{not: Bms}).} gave a proof of this theorem within $\BST$. Here we only
sketch this construction.

Let $\B \Lambda(X) := \big(\lambda_0^1, \C E^{1,X}, \lambda_1^1, \lambda_0^0, \C E^{0,X}, 
\lambda_1^0) \in \Fam(I, \B X)$
and $i_0 \in I$. To each $i \in I$ corresponds the real-valued partial function
\[ \chi_i := \big(\lambda_0^1(i) \cup \lambda_0^0(i), \C E_i^X, \chi_i^{\Real}\big) \in \C F(X), \]
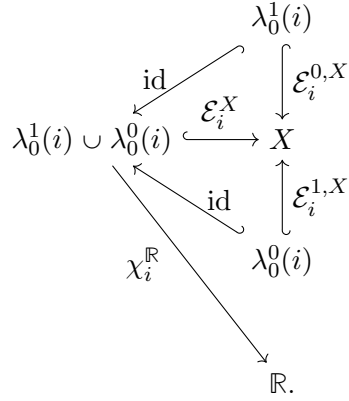
\begin{center}
\begin{tikzpicture}

\node (E) at (0,0) {$ \lambda_0^1(i) \cup \lambda_0^0(i)$};
\node[right=of E] (F) {$X$};
\node[below=of F] (B) {$\lambda_0^0(i)$};
\node[above=of F] (D) {$\lambda_0^1(i)$};
\node[below=of B] (C) {$\Real$.};

\draw[right hook->] (E)--(F) node [midway,above] {$\C E_i^X$};
\draw[->] (E)--(C) node [midway,left] {$\chi_i^{\Real} \ \ $};
\draw[right hook->] (D)--(F) node [midway,right] {$\C E_i^{0,X}$};
\draw[left hook->] (B)--(F) node [midway,right] {$\C E_i^{1,X}$};
\draw[left hook->] (B)--(E) node [midway,right] {$ \ \id$};
\draw[right hook->] (D)--(E) node [midway,left] {$\id \ $};

\end{tikzpicture}
\end{center}
where $\C E_i^X$ is the canonical embedding of $\lambda_0^1(i) \cup \lambda_0^0(i)$ into $X$, 
and $\chi_i^{\Real}$ is given by the rule of the partial function $\chi_{\B \lambda_0(i)}$.
The symbol $\id$ in the above diagram denotes the corresponding function defined by the identity-map rule. 
If $m, n \in \Nat^+$, $i_1, \ldots, i_n, j_1, \ldots, j_m \in I$, and $a_{1_i}, \ldots a_{i_n}, b_1, 
\ldots, b_m \in \Real$, the equality of the following real-valued partial functions is given by the commutativity 
of the following diagram
\[ \sum_{k = 1}^n a_{i_k} \chi_{i_k} := \bigg( \bigcap_{k = 1}^n \big(\lambda_0^1(i_k) \cup 
\lambda_0^0(i_k)\big), i^X_{\mathsmaller{\bigcap_{k = 1}^n (\lambda_0^1(i_k) \cup \lambda_0^0(i_k))}},
\ \sum_{k = 1}^n a_{i_k}\chi_{i_k}^{\Real}\bigg) \in \C F(X),\]
\[ \sum_{l = 1}^m b_{j_l} \chi_{j_l} := \bigg( \bigcap_{l = 1}^m \big(\lambda_0^1(j_l) \cup 
\lambda_0^0(j_l)\big), i^X_{\mathsmaller{\bigcap_{l = 1}^m (\lambda_0^1(j_l) \cup \lambda_0^0(j_l)}},
\ \sum_{l = 1}^m b_{j_l}\chi_{j_l}^{\Real}\bigg) \in \C F(X),\]
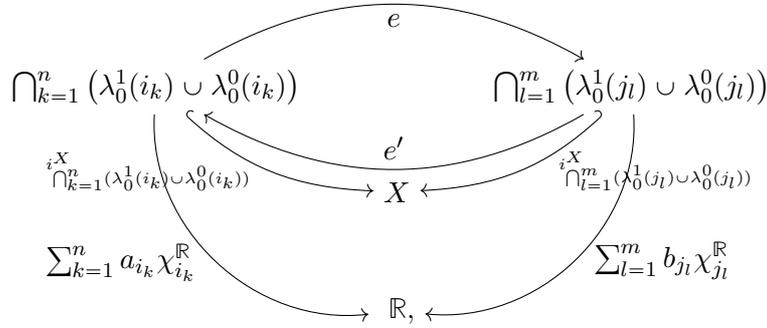
\begin{center}
\begin{tikzpicture}

\node (E) at (0,0) {$\bigcap_{k = 1}^n \big(\lambda_0^1(i_k) \cup \lambda_0^0(i_k)\big)$};
\node[right=of E] (B) {};
\node[right=of B] (F) {$\bigcap_{l = 1}^m \big(\lambda_0^1(j_l) \cup \lambda_0^0(j_l)\big)$};
\node[below=of B] (A) {$X$};
\node[below=of A] (C) {$ \ \Real$,};

\draw[->,bend left] (E) to node [midway,below] {$e$} (F);
\draw[->,bend left] (F) to node [midway,above] {$e{'}$} (E);
\draw[->,bend right=50] (E) to node [midway,left] {$\sum_{k = 1}^n a_{i_k}\chi_{i_k}^{\Real} \ $} (C);
\draw[->,bend left=50] (F) to node [midway,right] {$ \ \sum_{l = 1}^m b_{j_l}\chi_{j_l}^{\Real}$} (C);
\draw[right hook->,bend right=20] (E) to node [midway,left] {$ \mathsmaller{\mathsmaller{i^X_{\mathsmaller{\bigcap_{k = 1}^n (\lambda_0^1(i_k) \cup \lambda_0^0(i_k))}}}} \ $} (A);
\draw[left hook->,bend left=20] (F) to node [midway,right] {$\ \ \ \mathsmaller{\mathsmaller{ i^X_{\mathsmaller{\bigcap_{l = 1}^m (\lambda_0^1(j_l) \cup \lambda_0^0(j_l))}}}} $} (A);

\end{tikzpicture}
\end{center} 
for some (unique up to equality) functions $e$ and $e{'}$. If $\mu_0 \colon \Nat^+ \sto \D V_0$ is defined by 
the rule $\mu_0(n) := (\Real \times I)^n$, for every $n \in \Nat^+$, and if the corresponding dependent 
operation $\mu_1$ is defined in the obvious way, let the totality 
\[ S(I, \B \Lambda(X)) := \sum_{n \in \Nat^+}(\Real \times I)^n, \]
\[(n, u) =_{S(I, \B \Lambda(X))} (m, w) :\TOT \sum_{k = 1}^n a_{i_k} \chi_{i_k} 
=_{\C F(X)} \sum_{l = 1}^m b_{j_l} \chi_{j_l};,\]
where $u := \big((a_1, i_1), \ldots, (a_n, i_n)\big)$ and $w := \big((b_1, j_1), \ldots, (b_m, j_m)\big)$. 
The family of simple functions\index{family of simple functions} generated by the family $\B \Lambda (I, \B X)$ 
of complemented subsets of $X$ is the structure
$\Delta (X, \Real) := \big(\dom_0, \C Z^X, \dom_1, \C P^{\Real}\big) \in \Fam\big(S(I, \B \Lambda(X), X, \Real \big)$, 
where $\dom_0 \colon S(I, \B \Lambda(X)) \sto \D V_0$ is defined by the rule
\[ \dom_0 (n, u) := \bigcap_{k = 1}^n \big(\lambda_0^1(i_k) \cup \lambda_0^0(i_k)\big);  \ \ \ \ u 
:= \big((a_1, i_1), \ldots, (a_n, i_n)\big), \ n \in \Nat^+, \]
the embedding $\C Z^X_{n, u} \colon \dom_0(n, u) \eto X$ is defined in a canonical way through the embeddings 
$\C E_{i_k}^{1,X}$ and $\C E_{i_k}^{0,X}$, where $k \in \{1, \ldots, n\}$. 

If $(n, u) =_{S(I, \B \Lambda(X))} (m, w)$,
the mapping $\dom_{(n, u)(m, w)} \colon \dom_0(n, u) \to \dom_0(m, w)$ is defined, in order to avoid choice, as 
the mapping $E_{(n,u)(m,w)}$, where $E$ is a modulus of equality for $D(S(I, \B \Lambda(X)))$ 
with $E_{(n,u)(n,u)} := \id_{\dom_0(n, u)}$, for every $(n,u) \in S(I, \B \Lambda(X))$. The fact that 
$\Delta (X, \Real) \in \Set\big(S(I, \B \Lambda(X), X, \Real \big)$ is immediate to show. Hence, 
to every $(n, u) \in  S(I, \B \Lambda(X))$ corresponds the partial function
\[ s_{(n,u)} := \bigg(\dom_0(n, u), \C Z^X_{(n, u)},  \sum_{k = 1}^n a_{i_k}\chi_{i_k}^{\Real}\bigg) 
\in \C F(X); \ \ \ \ u := \big((a_1, i_1), \ldots, (a_n, i_n)\big). \]
If $\C M (\B \Lambda(X)) := (X, I, \mu)$ is a pre-measure space, then $\C M (\B \Lambda(X))$ induces 
the pre-integration space $\C L(\B \Lambda(X)) := (X, S(I, \B \Lambda(X)), \int_{\mu})$, where 
\[ \int_{\mu} \colon S(I, \B \Lambda(X)) \to \Real, \ \ \ \ (n, u) \mapsto \int_{\mu} (n, u),\]
\[ \int_{\mu} (n, u) := \sum_{k = 1}^n a_k \mu(i_k); \ \ \ \ u := \big((a_1, i_1), \ldots, (a_n, i_n)\big).\]
The many steps of this involved proof of Bishop and Cheng, appropriately translated into the
predicative framework of $\BST$, are found in~\cite{Ze19}, pp.~34--45.

\section{Notes}
\label{sec: notes6}

\begin{note}\label{not: historyborel}
\normalfont
The set of Borel sets generated by a given family of complemented subsets of a set $X$, with respect 
to a set $\Phi$ of real-valued functions on $X$, was introduced in~\cite{Bi67}, p.~68. This set is 
inductively defined and plays a crucial role in providing important examples of measure spaces in Bishop's 
measure theory developed in~\cite{Bi67}. As this measure theory was replaced in~\cite{BB85} by the 
Bishop-Cheng measure theory, an enriched version of~\cite{BC72} that made no use of Borel sets, the Borel 
sets were somehow ``forgotten'' in the constructive literature. 
In the introduction of~\cite{BC72}, Bishop and Cheng explained why they consider their new measure 
theory ``much more natural and powerful theory''. They do admit though, that some results are harder 
to prove (see~\cite{BC72}, p.~v). As it is also noted in~\cite{Sp02}, p.~25, the Bishop-Cheng measure 
theory is highly impredicative, while Bishop's measure theory in~\cite{Bi67} is highly predicative. 
This fact makes the original Bishop-Cheng measure theory hard to implement in some functional-programming 
language, a serious disandvantage from the computational point of view. This is maybe why, later attempts 
to develop constructive measure theory were done within an abstract algebraic framework 
(see~\cite{CP02},~\cite{CS09} and~\cite{Sp06}). 
Despite the above history of measure theory within Bishop-style constructive mathematics the set of
Borel sets is interestingly connected to the theory of Bishop spaces.
 
\end{note}

\begin{note}\label{not: onextensionalityborel}
\normalfont
The definition of $\Borel(\B \Lambda(X))$ is given by Bishop in~\cite{Bi67}, p.~68, although a
rough notion of a family of complemented subsets is used, condition $(\Borel_3)$ is not mentioned, 
and $F$ is an arbitrary subset of $\D F(X)$, and not necessarily a Bishop topology.
If we want to avoid the extensionality of $\Borel(\B \Lambda(X))$, we need to introduce a ``pseudo''-membership 
condition \[ \B A \ \dot {\in} \ \Borel(\B \Lambda(X)) 
:\TOT \exists_{\B B \in \Borel(\B \Lambda(X))}\big(\B A =_{\C P^{\Disj_F}(X)} \B B\big). \]
A similar condition is necessary, if we want to avoid extensionality in the definition the 
least Bishop topology $\bigvee F_0$. Such an approach though, is not practical, and not compatible 
to the standard practice to study extensional subsets of sets. The quantification over $\Fam(\D 1, F, \B X)$
is not equivalent to the quantification over the class $\C P^{\Disj_F}(X)$, as in order to define a family
in $\Fam(\D 1, F, \B X)$, we need to have \textit{already constructed} an $F$-complemented subset of $X$. 
I.e., an element of $\Fam(\D 1, F, \B X)$ is generated by an already constructed, or given element of
$\C P^{\Disj_F}(X)$, and not from an abstract element of it. Recall that we never define an assignment
routine from a class, like $\C P^{\Disj_F}(X)$, to a set like $\Fam(\D 1, F, \B X)$.

\end{note}

\begin{note}\label{not: oninddefinitionborel}
\normalfont
The notion of a least Bishop topology generated by a given set of function from $X$ to $\Real$, 
together with the set of Borel sets generated by a family of complemented subsets of $X$, are the main
two inductively defined concepts found in~\cite{Bi67}. The difference between the two inductive
definitions is non-trivial. The first is the inductive definition of a subset of $\D F(X)$, while 
the second is the inductive definition of a subset of the class $\C P^{\Disj_F}(X)$.

\end{note}

\begin{note}\label{not: onminus}
\normalfont
As Bishop remarks in~\cite{Bi67}, p.~69, the proof of Proposition~\ref{prp: borel1}(iii) rests on the
property of $F$ that $\big(\frac{1}{n} - f\big) \in F$, for every $f \in F$ and $n \geq 1$. If we 
define similarly the Borel sets generated by any set of real-valued functions $\Theta$ on $X$, then we 
can find $\Theta$ such that $\Borel(\Theta)$ is closed under complements without satisfying the condition 
$f \in \Theta \To \big(\frac{1}{n} - f \big)\in \Theta$. Such a set is $\D F (X, \D 2)$.  In this case we have that
\[ \B o_{\D F (X, \D 2)} (f) := \big([f = 1], [f = 0]\big) \ \ \ \& \ \ - \B o_{\D F (X, \D 2)} (f) = 
\B o_{\D F (X, \D 2)}(1 - f). \]
Hence, the property mentioned by Bishop is sufficient, but not necessary.

\end{note}

%
%
%
%
%

\begin{note}\label{not: BCms}
\normalfont
A measure space is defined in~\cite{BB85}, p.~282, and a complete measure space in~\cite{BB85}, p.~289.
These definitions appeared first in~\cite{BC72} p.~47 and p.~55, respectively\footnote{In~\cite{BC72}, p.~55,
condition $(\BCM_1)$ appears in the equivalent form: if $B$ is an element of $M$ such that 
$B^1 \subseteq A^1$ and $B^0 \subseteq A^0$, then $A \in M$, where we have used the terminology that 
corresponds to the formulation of $(\BCM_1)$ in the definition of Bishop-Cheng.}.\\[1mm]
\textbf{Bishop-Cheng definition of a measure space.}
A measure space is a triplet $(X, M, \mu)$ consisting of a nonvoid\index{Bishop-Cheng measure space}
set $X$ with an inequality $\neq$, a set $M$ of complemented sets in
$X$, and a mapping $\mu$ of $M$ into $\Real^{0+}$, such that the following properties
hold.\\[1mm]
$(\BCMS_1)$ If $A$ and $B$ belong to $M$, then so do $A \vee B$ and $A \wedge B$, and
$\mu(A) + \mu(B) = \mu(A \vee B) + \mu(A \wedge B)$.\\[1mm]
$(\BCMS_2)$ If $A$ and $A \wedge B$ belong to $M$, then so does $A - B$, and $\mu(A)
= \mu(A \wedge B) + \mu(A - B)$.\\[1mm]
$(\BCMS_3)$ There exists $A$ in $M$ such that $\mu(A) >0$.\\[1mm]
$(\BCMS_4)$ If $(A_n)$ is a sequence of elements of $M$ such that $\lim_{k \to \infty} 
\mu\big(\bigwedge_{n = 1}^k A_n\big)$ exists and is positive, then $\bigcap_{n}A_n^1$ is nonvoid.\\[1mm]
We then call $\mu$ the \textit{measure}, and the elements of $M$ the \textit{integrable sets}, of the
measure space $(X, M, \mu)$. For each $A$ in $M$ the nonnegative number $\mu(A)$ is called
the \textit{measure} of $A$.\\[2mm]
\textbf{Bishop-Cheng definition of a complete measure space.} A measure space $(X, M, \mu)$ is 
\textit{complete} if the following three conditions hold.\\[1mm]
$(\BCCMS_1)$ If $A$ is a complemented set, and $B$ is an element of $M$ such\index{Bishop-Cheng complete measure space}
that $\chi_A = \chi_B$ on $B^1 \cup B^0$, then $A \in M$.\\[1mm]
$(\BCCMS_2)$ If $(A_n)$ is a sequence of elements of $M$ such that 
$$l := \lim_{N \to \infty}\mu \bigg(\bigvee_{n = 1}^N A_n\bigg)$$
exists, then $\bigvee_n A_n$ belongs to $M$ and has measure $l$.\\[1mm]
$(\BCCMS_3)$ If $A$ is a complemented set, and if $B, C$ are elements of $M$ 
such that $B < A < C$ and $\mu(B) = \mu(C)$, then $A \in M$.\\[2mm]
As there is no indication of indexing in the description of $M$, the Bishop-Cheng definition of
a measure space seems to employ the powerset axiom in the formulation of $M$. The powerset axiom is
clearly used in $(\BMS_1)$ and $\BCM_3$. 

\end{note}

\begin{note}\label{not: Bms}
\normalfont
The following definition of Bishop is given in~\cite{Bi67}, p.~183.\\[1mm]
\textbf{Bishop definition of a measure space.} Let $F$ be a nonvoid family of real-valued functions on a
set $X$, such that $\epsilon - f \in F$ whenever $\epsilon > 0$ and $f \in F$. Let $\C F$ be any
family of complemented subsets of $X$ $($relative to $F)$, closed with respect
to countable unions, countable intersections, and complementation.
Let $\C M$ be a subfamily of $\C F$ closed under finite unions, intersections,
and differences. Let the function $\mu : \C M \to \Real^{0+}$ satisfy the following
conditions:\\[1mm]
$(\BMS_1)$ There exists a sequence $S_1 \subset S_2 \subset \ldots $ of elements of $\C M$ such 
that\footnote{Bishop here means $\bigcup_{n = 1}^{\infty} S_n = X$.} $\bigcup_{n = 1}^{\infty} S_n = X_0$ 
and $\lim_{n \to \infty} \mu(A \cap S_n) = \mu (A)$ for all $A$ in $\C M$.\\[1mm]
$(\BMS_2)$ If $A \in \C F$, and if there exist $B$ and $N$ in $\C M$ such that $(i)$ $\mu(N) = 0$,
$(ii)$ $x \in A$ whenever $x \in B - N$, and $(iii)$ $x \in -A$ whenever $x \in -B - N$, then $A \in \C M$
and $\mu(A) = \mu(B)$.\\[1mm]
$(\BMS_3)$ If $A \in \C M$, if $B \in \C F$, and if $A \cap B \in \C M$, then $A - B \in \C M$ and
$\mu(A) = \mu(A - B) + \mu(A \cap B)$.\\[1mm]
$(\BMS_4)$ We have $\mu(A) + \mu(B) = \mu(A \cup B) + \mu(A \cap B)$ for all $A$ and $B$ in $\C M$.\\[1mm]
$(\BMS_5)$ For each sequence $\{An\}$ of sets in $\C M$ such that $c := \lim_{n \to \infty} 
\mu\big(\bigcup_{k = 1}^n A_k\big)$ $[$respectively, $c := \lim_{n \to \infty} \mu\big(\bigcap_{k = 1}^n A_k\big)]$
exists, the set $A := \bigcup_{k = 1}^{\infty}A_k$ $($respectively, $A := \bigcap_{k = 1}^{\infty}A_k)$ 
is in $\C M$, and $\mu(A) = c$.\\[1mm]
$(\BMS_6)$ Each $A$ in $\C M$ with $\mu(A) > 0$ is nonvoid.\\[1mm]
Then the quintuple $(X, F, \C F, \C M, \mu)$ is called a \textit{measure space}, $\mu$ is the
\textit{measure}, $\C F$ is the class of \textit{Borel sets}, and $\C M$ is the class of \textit{integrable sets}.\\[2mm]
If in Bishop's definition we understand the families of complemented subsets $\C F$ and $\C M$ 
as indexed families $(\B A)_{i \in I}, (\B A)_{j \in J}$ over some sets $I$ and $J$, respectively, with 
$J \subseteq I$, then the quantifications involved in the clauses of Bishop's definition are over $I$ and $J$, 
and not over some class. Since in~\cite{Bi67}, p.~65, a family of subsets of $X$ is defined as an 
appropriate set-indexed family of sets, Bishop's first definition of measure space is predicative.

\end{note}

\begin{note}\label{not: numerical}
\normalfont
Regarding the exact definition of a measure space within the formal system $\Sigma$ introduced by Bishop
in~\cite{Bi70},
Bishop writes in~\cite{Bi70}, p.~67, the following:
\begin{quote}
To formalize in $\Sigma$ the notion of an abstract measure space, definition 1 of chapter 7 
of~\cite{Bi67} must be 
rewritten as follows. A \textit{measure space} is a family $\C M \equiv \{A_t\}_{t \in T}$ of 
complemented subsets
of a set $X$ relative to a certain family $\C F$ of real-valued function on $X$, a map $\mu : T \to \Real^{0+}$,
and an additional structure as follows: The void set $\emptyset$ is an element $A_{t_0}$ of $\C M$, and
$\mu(t_0) = 0$.
If $s$ and $t$ are in $T$, there exists an element $s \mathsmaller{\vee} t$ of $T$ such that
$A_{s \mathsmaller{\vee} t} < A_s \cup A_t$.
Similarly,
there exist operations $\mathsmaller{\wedge}$ and $\mathsmaller{\sim}$ on $T$, corresponding 
to the set theoretic operations $\cap$ and $-$. 
The usual algebraic axioms are assumed, such as $\mathsmaller{\sim}(s \ \mathsmaller{\vee} \ t) = 
\mathsmaller{\sim} s \ \mathsmaller{\wedge} \ \mathsmaller{\sim} t$. Certain measure-theoretic 
axioms, such as $\mu(s \ \mathsmaller{\vee} \ t) + \mu(s \ \mathsmaller{\wedge} \ t) = \mu(s) + \mu(t)$, 
are also assumed. Finally, there exist operations $\vee$ and $\wedge$. If, for example, $\{t_n\}$ 
is a sequence such that $C \equiv \lim_{k \to \infty} \mu(t_1 \ \mathsmaller{\vee} \ldots 
\mathsmaller{\vee} \ t_k)$
exists, then $\vee \{t_n\}$ is an element of $T$ with measure $C$. Certain axioms for $\vee$ and $\wedge$
are assumed. If $T$ is the family of measurable sets of a compact space relative to a measure $\mu$, and the 
set-theoretic function $\mu : T \to \Real^{0+}$ and the associated operations are defined as indicated above, 
the result is a measure space in the sense just described.

Considerations such as the above indicate that essentially all of the material in~\cite{Bi67}, 
appropriately modified, can be comfortably formalised in $\Sigma$.
\end{quote}

The expression $A_{s \mathsmaller{\vee} t} < A_s \cup A_t$ is probably a typo (it is the 
writing $A_{s \mathsmaller{\vee} t} = A_s \cup A_t$, which expresses the ``weak belongs to'' 
relation for $\lambda_0 I$). Bishop does not mention that $\C M$ is a set of complemented subsets of $X$, 
he only says that it is a family of such sets. This is not the case in~\cite{BB85}, p.~282. This explanation
given by Bishop regarding the explicit and unfolded writing of many of the definitions in constructive mathematics
refer to~\cite{Bi67}. I have found no similar comment of Bishop with respect to his later measure theory,
developed with Cheng. Moreover, I have found no such comment in the extensive work of Chan on Bishop-Cheng 
measure and probability theory.
 
\end{note}

\begin{note}\label{not: onMy75}
\normalfont
In~\cite{My75}, p.~354, Myhill criticised Bishop for using a set of subsets $M$ in the definition of 
a measure space, hence, according to Myhill, Bishop used the powerset axiom. Since $M$ is an $I$-set of
subsets of $X$, in the sense described in section~\ref{sec: setofsubsets}, Myhill's critique is not 
correct. Bishop's exaplanation in the previous extract is also a clear reply to a critique like Myhill's. 
Notice that Myhill's paper~\cite{My75} refers only to~\cite{Bi67}, and it does not mention~\cite{Bi70}, 
which includes Bishop's clear explanation. This is quite surprising, as Myhill's paper, received in 
January 1974, was surely written after the publication of~\cite{KMV70}, in which Bishop's paper~\cite{Bi70} is included 
and Myhill is one of its three editors!
Myhill's critique would be correct, if he was referring to the Bishop-Cheng measure space defined
in~\cite{BC72}, a work published quite some time before Myhill submit~\cite{My75}. Myhill though, 
does not refer to~\cite{BC72} in~\cite{My75}.

\end{note}

\begin{note}\label{not: ondefsigmapfun}
\normalfont
Definition~\ref{def: sigmapfun} is the explicit writing within $\BST$ of the corresponding 
definition in~\cite{BB85}, pp.~216--217.
\end{note}

\begin{note}\label{not: BCis}
\normalfont
The following definition is given in~\cite{BC72}, p.~2, and it is repeated in~\cite{BB85}, p.~217.\\[1mm]
\textbf{Bishop-Cheng definition of an integration space.}
A triplet $(X, L, I)$ is an integration space if $X$ is a nonvoid set with an inequality $\neq$, $L$ 
is a subset of $\C F(X)$ (this set is $\C F^{se}(X)$ in our terminology), and $I$ is a mapping of $L$ 
into $\Real$ such that the following properties hold\index{Bishop-Cheng definition of an integration space}.\\[2mm]
$(\BCIS_1)$ If $f, g \in L$ and $\alpha, \beta \in \Real$, then $\alpha f + \beta g $, $|f|$, and $f 
\wedge 1$ belong to $L$, and $I(\alpha f + \beta g) = \alpha I(f) + \beta I(g)$.\\[1mm]
$(\BCIS_2)$ If $f \in L$ and $(f_n)$ is a sequence of nonnegative functions in $L$ such that $\sum_{n}I(f_n)$
converges and $\sum_nI(f_n) < I(f)$, then there exists $x \in X$ such that $\sum_nf_n(x)$ converges 
and $\sum_nf_n(x) < f(x)$.\\[1mm]
$(\BCIS_3)$ There exists a function $p$ in $L$ with $I(p) = 1$.\\[1mm]
$(\BCIS_4)$ For each $f$ in $L$, $\lim_{n \to \infty} I(f \wedge n) = I(f)$ and $\lim_{n \to \infty} I(|f| \wedge n^{-1}) = 0$.\\[2mm]
The notion of an integration space is a constructive version of the Daniell integral, 
introduced in~\cite{Da18}. The Bishop-Cheng definition of an integration space is impredicative, as 
the class $\C F(X)$ is treated as a set. The notion of a subset is defined only for sets, and $L$ is
considered a subset of $\C F(X)$. The extensional character of $L$ is also not addressed. This
impredicative approach to $L$ is behind the simplicity of the Bishop-Cheng definition. E.g., in $(\BCIS_4)$ 
the formulation of the limit is immediate as the terms $f \wedge n \in L$ and $I$ is defined on $L$. In
Definition~\ref{def: intspace} though, we need to use an element $\alpha(n)$ of the index-set $I$ such 
that $f \wedge n =_{\C F(X)} f_{\alpha(n)}$, in order to express the limit.
\end{note}

\begin{note}\label{not: onZeuner}
\normalfont
The Bishop-Cheng definition of the ``set'' $L^1$\index{$L^1$} (or $L^p$, \index{$L^p$}where $p \geq 1$) 
of \textit{integrable functions}\index{integrable function} is also impredicative, as it rests on
the use of the totality $\C F^{se}(X)$ as a set (see Definition (2.1) in~\cite{BB85}, p.~222). 
In~\cite{Ze19}, pp.~49--60, the pre-integration space $L^1$ of \textit{canonically} integrable functions
is studied instead within $\BST$, as the completion of an integration space. The set $L^1$ is 
predicatively defined in~\cite{Bi67}, p.~190, as an integrable function is an appropriate measurable 
function, which is defined using quantification over the set-indexed family $\C M$ of integrable sets
in a Bishop measure space (see Note~\ref{not: Bms}).
\end{note}

\chapter{Epilogue}
\label{chapter: epilogue}

\section{$\BST$ between dependent type theory and category theory}
\label{sec: typescats}

Here we tried to show how the elaboration of the notion of a set-indexed family of
sets within $\BST$ expands the range of $\BISH$ both in its foundation and its practice. 
Chapters~\ref{chapter: BST}-\ref{chapter: proofrelevance} are concerned with the foundations of
$\BISH$, and chapters~\ref{chapter: bspaces} and~\ref{chapter: measure} with the practice of $\BISH$. 

Chapter~\ref{chapter: BST} presents the set-like objects, the families of which are studied later: 
sets, subsets, partial functions, and complemented subsets. Operations between these objects generate
corresponding operations between their families and family-maps. Chapter~\ref{chapter: familiesofsets}
includes the fundamental notions and results about set-indexed families of sets. A family of 
sets $\Lambda \in \Fam(I)$, together with its $\sum$- and $\prod$-set, and a family map $\Psi \colon \Lambda \To M$,
are examples of notions with a strong type-theoretic, or categorical flavour, depending on the point of 
observation view. This is not accident, as $\MLTT$ was motivated by Bishop's book~\cite{Bi67}. 
Moreover, $\BST$ can roughly be described as a fundamental informal theory of totalities and assignment
routines, and (informal) category theory as a fundamental  (informal) theory of objects and arrows. A 
fundamental similarity between $\BST$ and $\MLTT$ is the explicit use of dependency, which is suppressed in category theory. The fundamental categorical concepts of a functor and a natural transformation,
which are translated within $\BST$ as an $I$-family of sets and a family-map between $I$-families of sets,
have an immediate and explicit formulation within dependent type theory or within $\BST$ 
(see Note~\ref{not: categories}). The formulation of dependency though, within category theory is much
more involved (see e.g.,~\cite{Pa14}). On the other hand, a fundamental similarity between $\BST$ and 
category theory is the use of definitions that do not ``force''  facts and results, as in the case of 
$\MLTT$ and its recent extension $\HoTT$. While the language of $\MLTT$ is clearly closer to $\BST$, a
large part of pure category theory, the size of the totalities involved excluded, follows the ``pattern'' 
of doing constructive mathematics in the style of $\BISH$:  all notions are defined, no powerful axioms
are used, and despite the generality in the categorical formulations, most results have a concrete 
algorithmic\footnote{The question of the constructive character of general category theory is addressed in~\cite{Mc06}. 
There constructivism in mathematics is identified with Brouwer's intuitionism. The inclusion of Bishop-style constructivism 
and of type-theoretic constructivism in the interpretation of mathematical constructivism is necessary and sheds more
light on the original question.} meaning. 

The interconnections between category theory and dependent type theory is a standard theme behind 
foundational studies on mathematics and theoretical computer science the last forty years. The recent explosion of univalent foundations, spearheaded 
by the Fields medalist Vladimir Voevodsky, regenerated the study of these interconnections. The appropriate categorical understanding of the univalence axiom brought category-theorists and type-theorists even closer. $\BST$ seems to be 
in some kind of common territory between dependent type theory and category theory. It also features 
simultaneously the proof-irrelevance of category theory and classical mathematics,  and the proof-relevance of
$\MLTT$. In contrast to $\HoTT$, where a type $A$ has a rich space-structure due to the induction principle
corresponding to the introduction of the identity family $=_A \colon A \to A \to \C U$ on $A$, the notion 
of space in $\BISH$, as in classical set-based mathematics, is not identical to that of a set. This is also
captured in category theory, where the category of sets behaves differently from the category of topological spaces.
We need to add, by definition, extra structure to a set $X$, in order to acquire a non-trivial space structure.
In this work the concept of space considered was that of a Bishop space. This is one option, which is shown to
be very fruitful, if we work within $\BISH^*$, but it is not the only one.



As in the case of $\MLTT$ or $\HoTT$, a non-trivial part of category theory can be studied \textit{within} $\BST$. We gave a glimpse of that in Note~\ref{not: categories}. Working in a similar fashion, most of the theory of small categories can
appropriately be translated into $\BST$. This modelling of pure category theory ``suffers'', as any modelling, 
from the inclusion of features, like conditions $(\Cat_3), (\Cat_4)$, and $(\Funct_3)$, that depend on the 
system $\BST$ itself and are not part of the original theory. In any event, such a translation is not meant
to be an attempt to replace pure category theory, but to embed into $\BISH$ concepts and facts from category
theory useful to the practice of $\BISH$. For example, all categorical notions and facts of constructive algebra 
presented in~\cite{MRR88}, within a category theory irrelevant to the version of Bishop's theory of sets
underlying~\cite{MRR88}, can, in principle, be approached within $\BST$ and the corresponding category
theory within $\BST$.  
Unfolding proof-relevance in $\BISH$ through $\BST$, categorical facts, like the Yoneda lemma 
for $\Fam(\widehat{I})$ 
can be translated from $\MLTT + \FunExt$ to $\BISH$.
It remains to find though, interesting applications of such results to $\BISH$.


Inductive definitions bring the language of $\BISH^*$ closer to dependent type theory. The induction
principles that accommodate inductive definitions in the latter correspond to universal properties 
in category theory. The formalisation of $\BST$, and its possible extension $\BST^*$ with inductive 
definitions with rules of countably many premises, is an important open problem. The natural requirement
for a faithful and adequate formal system for $\BST$ and $\BST^*$ makes the choice of the formal framework
even more difficult. It seems that a version of a formal version of extensional Martin-L\"of type theory,
and the corresponding theory of setoids within it, is a formal system very close to the informal system $\BST$.
As we have explained in Note~\ref{not: onpresentationaxiom}, a formal version of intensional $\MLTT$ does
not seem to be a faithful formal system for the informal theory $\BISH$. The logical framework of an extensional version of dependent type theory though, and the identification of propositions with types, is quite far from the usual practice of $\BISH$, which is, in this respect, close to the standard practice of  classical mathematics. It is natural to search for a formal 
system of $\BST$ where logic is not built in, as in $\MLTT$, and which reflects the way sets are defined in $\BST$. 
We hope that the presentation of $\BST$ in this work will be helpful to the construction of such a formal counterpart.   


Category theory can also be very helpful to the formulation of the properties of Bishop sets and functions. 
The work of Palmgren~\cite{Pa12b} on the categorical properties of the category of setoids and setoid maps 
within intensional $\MLTT$ is expected to be very useful to this. A similar formulation of the categorical
properties of the theory of setoids and setoid maps within extensional $\MLTT$ could be even closer to 
the formulation of the categorical properties of Bishop sets and functions.

%
%

\vspace{-6mm}

\section{Further open questions and future tasks}
\label{sec: open}

We collect here some further open questions and future tasks stemming form this work.

\begin{enumerate}

\item To develop the theory of neighbourhood spaces using the notion of a neighbourhood family of
subsets of  a set $X$ that covers $X$ (see Note~\ref{not: exfamofsubsets}). 

\item Is it possible to use families of complemented subsets to describe a neighbourhood space? The 
starting idea is to assign to each $i \in I$ a complemented subset $\B \nu_0(i) := \big(\nu_0^1(i), 
\nu_0^0(i)\big)$ of $X$, such that $\nu_0^1(i)$ is open and
$\nu_0^0(i)$ is closed. The benefit of such an approach to constructive topology is that the classical 
duality between open and closed sets is captured constructively.  E.g., the $1$-component of the 
complement $- \B \nu_0(i) := \big(\nu_0^0(i), \nu_0^1(i)\big)$ of $\B \nu_0(i)$ is a closed set and
its $0$-component is an open set.

\item Can we use complemented subsets of $\Nat$ in a constructive reconstruction of recursion theory, 
instead of just subsets of $\Nat$? This question is inspired from the work of Nemoto on recursion theory 
within intuitionistic logic.

\item To explore further the notion of an impredicative set, and the hierarchy mentioned in
Note~\ref{not: onfamsoffams}.

\item To find interesting purely mathematical applications of set-relevant families of sets and
of families of families of sets.
\item To investigate the possibility of a $\BHK$-interpretation of a negated formula
(see Note~\ref{not: onnegationbhk}).

\item To develop a (predicative) theory of ordinals within $\BST$.

\item To study families of sets with a proof-relevant equality over an index-set with a proof-relevant equality.

\item To translate more notions and results from $\MLTT$ and $\HoTT$ to $\BISH$ through $\BST$. As a special case,
to translate higher inductive types (HITs) into $\BISH$, other than the truncation $||A||$ of $A$. If we work
directly with a space in $\BST$ i.e., with a Bishop space $\C F := (X, F)$, and not with an arbitrary type, 
as in $\HoTT$, we can define within $\BST$ notions like the cone\index{cone of a Bishop space} and the 
suspension\index{suspension of a Bishop space} of $\C F$.
If $\oi := [0, 1]$, we call $\D I_{\mathsmaller{01}}^{\mcup} := \{0\} \cup (0, 1) \cup \{1\}$ the 
\textit{pseudo-interval} $[0, 1]$\index{pseudo-interval $[0, 1]$}. To $\D I_{\mathsmaller{01}}^{\mcup}$ 
we can associate the least Bishop topology generated by the restriction of the identity map to it. The 
relation $\sim_X^{\tau}$ on $X \times \oic$, defined by
\[ (x, i) \eqt (x{'}, i{'})  : \TOT \big( i, i{'} \in \{0\} \cup (0, 1) \ \& \ i =_{\oic} i{'} \ \& \ x =_X x{'} \big) 
\ \  \mbox{or} \ \   i =_{\oic} 1 =_{\oic} i{'},\]
is an extensional equivalence relation. If $Y := X \times \oic$ and $\tau_0^X \colon Y \to \tau_0^X Y$ is the function that maps $(x, i)$ to its equivalence class (see section~\ref{sec: fameqclass}), then $\tau_0^X Y(Y)$, equipped with an an appropriate Bishop topology, is the cone\index{cone of a Bishop space} of $\C F$. For the suspension of $\C F$ we\index{suspension of a Bishop space} work similarly.

\item To find interesting mathematical applications of $(-2), (-1)${-} and $0${-}sets in $\BISH$. 

\item To elaborate the study of category theory within $\BST$. So far we have formulated within $\BST$ most of the category theory formulated within the Calculus of Inductive Constructions in~\cite{HS98}.





\item  To develop along the lines of Chapter~\ref{chapter: bspaces} the theory of spectra of other structures,
like groups, rings, modules etc.

\item To develop further the theory of Borel sets of a Bishop topology. E.g., to find the exact 
relation between the Borel sets $\Borel(\C F)$ and $\Borel(\C G)$ and the Borel sets $\Borel(\C F \times \C G)$
of the product Bishop space $\C F \times \C G$. And similarly for all important constructions of new Bishop 
spaces from given ones.

\item To formulate various parts of the constructive algebra developed in~\cite{MRR88} and~\cite{LQ15} within $\BST$.

\item To elaborate the theory of (pre-)measure spaces and (pre-)integrations spaces. The past work~\cite{Ze19} 
and the forthcoming work~\cite{PZ20} are in this direction. 

\item To approach Chan's probability theory in~\cite{Ch19}, which is within $\BCMT$, 
through a predicative reconstruction of $\BCMT$ within $\BST$.

\end{enumerate}

\chapter{Appendix}
\label{chapter: appendix}

\section{Bishop spaces}
\label{sec: bishop}

We present the basic notions and facts on Bishop spaces that are used in the previous sections.
For all concepts and results from constructive real analysis that we use here without
further explanation we refer to~\cite{BB85}. For all proofs that are not included in this section we 
refer to~\cite{Pe15}. We work within the extension $\BISH^*$ of $\BISH$ with inductive definitions
with rules of countably many premises. A Bishop space is a constructive, function-theoretic 
alternative to the classical notion of a topological space, and a Bishop morphism is the
corresponding function-theoretic notion of ``continuous function'' between Bishop spaces.

\begin{definition}\label{def: cont1}
If $X$ is a set and $\Real$ is the set of real numbers, we denote by $\mathbb{F}(X)$
the set of functions from $X$ to $\Real$, and by $\Const(X)$\index{$\Const(X)$} the 
subset of $\mathbb{F}(X)$ of all constant functions on $X$. 
If $a \in \Real$, we denote by $\overline{a}^X$\index{$\overline{a}^X$}
the constant function on $X$ with value $a$. We denote by $\Nat^+$ the set of non-zero natural numbers.  
A function $\phi: \Real \rightarrow \Real$ is called \textit{Bishop continuous}, or simply continuous,
if for every $n \in \Nat^+$ there is a function\index{$\omega_{\phi,n}$} $\omega_{\phi, n}:
\mathbb{R}^{+} \rightarrow \mathbb{R}^{+}$,
$\epsilon \mapsto \omega_{\phi, n}(\epsilon)$, which is called a \textit{modulus 
of continuity}\index{modulus of (uniform) continuity} of $\phi$ on $[-n, n]$, such that the following 
condition is satisfied
\[\forall_{x, y \in [-n, n]}(|x - y| < \omega_{\phi, n}(\epsilon) \Rightarrow 
|\phi(x) - \phi(y)| \leq \epsilon),\]
for every $\epsilon > 0$ and every $n \in \Nat^+$. We denote by $\BR$\index{$\BR$} the set of continuous functions 
from $\Real$ to $\Real$, which is equipped with the equality inherited from $\D F(\Real)$.

\end{definition}

Note that we could have defined the modulus of continuity $\omega_{\phi, n}$ as a function from $\Nat^+$ 
to $\Nat^+$. Clearly, a continuous function $\phi \colon \Real \to \Real$ is uniformly continuous on every bounded 
subset of $\Real$. The latter is an impredicative formulation of uniform continuity, since it 
requires quantification over the class of all subsets of $\Real$. The formulation of uniform continuity 
in the Definition~\ref{def: cont1} though, is predicative, since it requires quantification over the sets
$\Nat^+, \D F(\Real^+, \Real^+)$ and $[-n, n]$.

\begin{definition}\label{def: notation1}
If $X$ is a set, $f, g \in \mathbb{F}(X)$, $\epsilon > 0$, and $\Phi \subseteq \mathbb{F}(X)$, 
let\index{$U(X; \Phi, g, \epsilon)$}  
\[ U(X; f, g, \epsilon) :\TOT \forall_{x \in X}\big(|g(x) - f(x)| \leq \epsilon\big),\]
\[ U(X; \Phi, f) :\TOT \forall_{\epsilon > 0}\exists_{g \in \Phi}\big(U(f, g, \epsilon)\big).\]
If the set $X$ is clear from the context, we write simpler $U(f, g, \epsilon)$\index{$U(f,g, \epsilon)$} 
and\index{$U(\Phi, f)$} $U(\Phi, f)$, respectively.\index{$U(X; \Phi, f)$}
We denote by $\Phi^*$ the bounded elements of $\Phi$, and its uniform 
closure\index{uniform closure} $\overline{\Phi}$ is defined by\index{$\overline{\Phi}$} 
\[ \overline{\Phi} := \{f \in \D F(X) \mid U(\Phi, f)\}.\]
\end{definition}

A Bishop topology on $X$ is a certain subset of $\D F(X)$. Since the Bishop topologies considered here
are all extensional subsets of $\D F(X)$, we do not mention the embedding $i_F^{\D F(X)} \colon F \eto \D F(X)$, 
which is given in all cases by the identity map-rule.

\begin{definition}\label{def: bishopspace}
A \textit{Bishop space}\index{Bishop space} is a pair $\C F := (X, F)$, where $F$ is an extensional subset of $\D F(X)$,
which is called a \textit{Bishop topology}, or simply a \textit{topology}\index{Bishop topology}
of functions on $X$, that satisfies the following conditions:\\[1mm]
$(\BS_1)$ If $a \in \Real$, then $\overline{a}^X \in F$.\\[1mm]
$(\BS_2)$ If $f, g \in F$, then $f + g \in F$.\\[1mm]
$(\BS_3)$ If $f \in F$ and $\phi \in \Bic(\Real)$, then $\phi \circ f \in F$
\begin{center}
\begin{tikzpicture}

\node (E) at (0,0) {$X$};
\node[right=of E] (F) {$\Real$};
\node[below=of F] (A) {$\Real$.};

\draw[->] (E)--(F) node [midway,above] {$f$};
\draw[->] (E)--(A) node [midway,left] {$F \ni \phi \circ f \ $};
\draw[->] (F)--(A) node [midway,right] {$\phi \in \BR$};

\end{tikzpicture}
\end{center}
$(\BS_4)$ $\overline{F} = F$.
\end{definition}

If $\C F := (X, F)$ is a Bishop space, then $\C F^* := (X, F^*)$ is the Bishop space of 
bounded\index{$\C F^*$}\index{$F^*$} elements of $F$. The constant functions
$\Const(X)$ is the \textit{trivial}
topology on $X$, while $\D F(X)$ is the \textit{discrete} topology on $X$. Clearly, if $F$ is a topology on $X$,
then $\Const(X) \subseteq F \subseteq \D F(X)$, and the set of its bounded elements
$F^{*}$ is also a topology on $X$. It is straightforward to see that the pair $\C R := (\Real, \BR)$ is a
Bishop space, which we call the \textit{Bishop space of reals}.
A Bishop topology $F$ is a ring and a lattice; since $|\id_{\Real}| \in \Bic(\Real)$, where $\id_{\Real}$ is the 
identity function on $\Real$, by BS$_{3}$ we get that if $f \in F$ then $|f| \in F$. 
By BS$_{2}$ and BS$_{3}$, and using the following equalities  
\[ f{\cdot}g = \frac{(f + g)^{2} - f^{2} - g^{2}}{2} \in F,\]
\[ f \vee g = \max\{f, g\} = \frac{f + g + |f - g|}{2} \in F, \]  
\[ f \wedge g = \min\{f, g\} = \frac{f + g - |f - g|}{2} \in F,\]
we get similarly that if $f, g \in F$, then $f{\cdot}g, f \vee g, f \wedge g \in F$.
Turning the definitional clauses of a Bishop topology into inductive rules, Bishop defined in~\cite{Bi67}, p.~72,
the least topology including a given subbase $F_{0}$. This inductive definition, which 
is also found in~\cite{BB85}, p.~78, is crucial to the definition of new Bishop topologies from given ones.

\begin{definition}\label{def: Least}
The \textit{Bishop closure} of $F_{0}$, or the \textit{least topology}\index{least Bishop topology} 
$\bigvee{F_{0}}$\index{$\bigvee F_0$} 
generated by some $F_{0} \subseteq \mathbb{F}(X)$, is defined by the following inductive rules:
\[ \frac{f_{0} \in F_{0}}{f_{0} \in \bigvee F_{0}}, \ \ \ \frac{a \in \Real}{\overline{a}^{X} \in
\bigvee F_{0}},
\ \ \ \frac{f, g \in \bigvee F_{0}}{f + g \in \bigvee F_{0}}, \ \ \  \frac{f \in \bigvee F_{0}, \ 
g =_{\D F(X)} f}{g \in \bigvee F_{0}},\]
\[ \frac{f \in \bigvee F_{0}, \ \phi \in \Bic(\Real)}{\phi \circ f \in \bigvee F_{0}},
\ \ \ \ \ \frac{\big(g \in \bigvee F_{0}, \ U(f, g, \epsilon)\big)_{\epsilon > 0}}{f \in \bigvee F_{0}}.\]
We call $\bigvee F_{0}$ the \textit{Bishop closure} of $F_{0}$, and $F_{0}$ a \textit{subbase}\index{subbase} of 
$\bigvee F_{0}$.
\end{definition}

If $F_{0}$ is inhabited, then $(\BS_1)$ is provable by $(\BS_3)$. The last, most complex rule
above can be replaced by the rule  
\[ \frac{g_{1} \in \bigvee F_{0} \ \wedge \ U\big(f, g_{1}, \frac{1}{2}\big), 
\ \  g_{2} \in \bigvee F_{0} \ \wedge \ 
      U\big(f, g_{2}, \frac{1}{2^{2}}\big), \   \ldots}{f \in \bigvee F_{0}},\]
a rule with countably many premisses. 
The corresponding induction principle $\Ind_{\bigvee F_{0}}$\index{$\Ind_{\bigvee F_0}$} is 
\[\bigg[\forall_{f_{0} \in F_{0}}\big(P(f_{0})\big) \ \& \
\forall_{a \in \mathbb{R}}\big(P(\overline{a}^{X})\big) \ 
\&  \ \forall_{f, g \in \bigvee F_{0}}\big(P(f) \ \& \ P(g) \Rightarrow P(f + g) \]
\[ \& \ \forall_{f \in \bigvee F_0}\forall_{g \in \D F(X)}\big(g =_{\D F(X)} f \To P(g)\big)\]
\[ \& \ \forall_{f \in \bigvee F_{0}}\forall_{\phi \in \Bic(\mathbb{R})}\big(P(f)
\Rightarrow P(\phi \circ f)\big) \]
\[\ \ \ \ \ \ \ \ \ \ \ \ \ \ \ \ \ \ \ \ \ \& \  \forall_{f \in \bigvee F_{0}}\big( 
 \forall_{\epsilon > 0}\exists_{g \in \bigvee F_{0}}(P(g) \ \& \ U(f, g, \epsilon))
 \Rightarrow P(f)\big)\bigg]\]
 \[\Rightarrow \forall_{f \in \bigvee F_{0}}\big(P(f)\big),\]
where $P$ is any bounded formula.
Next we define the notion of a Bishop morphism\index{Bishop morphism} between Bishop spaces. 
The Bishop morphisms are the arrows in the category of Bishop\index{$\Bis$} spaces $\Bis$.

\begin{definition}\label{def: bmorphism}
If $\C F := (X, F)$ and $\C G = (Y, G)$ are Bishop spaces, a function $h: X \rightarrow Y$ is called
a \textit{Bishop morphism}, if $\forall_{g \in G}(g \circ h \in F)$ 
\begin{center}
\begin{tikzpicture}

\node (E) at (0,0) {$X$};
\node[right=of E] (F) {$Y$};
\node[below=of F] (A) {$\Real$.};

\draw[->] (E)--(F) node [midway,above] {$h$};
\draw[->] (E)--(A) node [midway,left] {$F \ni g \circ h \ $};
\draw[->] (F)--(A) node [midway,right] {$g \in G$};

\end{tikzpicture}
\end{center}
We denote by $\Mor(\C F, \C G)$\index{$\Mor(\C F, \C G)$} the set of Bishop morphisms 
from $\C F$ to $\C G$. As $F$ is an extensional subset of $\D F(X)$, $\Mor(\C F, \C G)$ is an
extensional subset of $\D F(X,Y)$.  If $h \in \Mor(\C F, \C G)$, the \textit{induced mapping} 
$h^* \colon G \to F$ from $h$\index{induced mapping from a Bishop morphism}\index{$h^*$}
is defined by the rule
\[ h^*(g) := g \circ h; \ \ \ \ g \in G.\]
\end{definition}

If $\mathcal{F} := (X, F)$ is a Bishop space, then $F = \Mor(\mathcal{F}, \mathcal{R})$,
and one can show inductively that if $\C G := (Y, \bigvee G_0)$, then 
$h: X \rightarrow Y \in \Mor(\mathcal{F}, \C G)$ if and only if
$\forall_{g_{0} \in G_{0}}(g_{0} \circ h \in F)$
\begin{center}
\begin{tikzpicture}

\node (E) at (0,0) {$X$};
\node[right=of E] (F) {$Y$};
\node[below=of F] (A) {$\Real$.};

\draw[->] (E)--(F) node [midway,above] {$h$};
\draw[->] (E)--(A) node [midway,left] {$F \ni g_0 \circ h \ \ $};
\draw[->] (F)--(A) node [midway,right] {$g_0 \in G_0$};

\end{tikzpicture}
\end{center}
We call this fundamental fact the\index{$\bigvee$-lifting of morphisms} $\bigvee$-\textit{lifting of morphisms}. 
A Bishop morphism is a \textit{Bishop isomorphism}\index{Bishop isomorphism}, if it is an isomorphism in
the category $\Bis$. We write $\C F \simeq \C G$\index{$\C F \simeq \C G$} to denote that $\C F$ and $\C G$ 
are Bishop isomorphic. If $h \in \Mor(\C F, \C G)$ is a bijection, then 
$h$ is a Bishop isomorphism if and only if it is \textit{open}\index{open morphism} i.e., 
$\forall_{f \in F}\exists_{g \in G}\big(f = g \circ h\big)$. 
%

\begin{definition}\label{def: new} Let $\C F := (X, F), \C G := (Y, G)$ be Bishop spaces, and $(A, i_A) \subseteq X$ 
inhabited.
The \textit{product} Bishop space\index{product Bishop space}
$\C F \times \C G := (X \times Y, F \times G)$ of $\C F$\index{$F \times G$}
and $\C G$, the \textit{relative}\index{relative Bishop space} Bishop space $\Fii_{|A} := (A, F_{|A})$ on $A$, and 
the \textit{pointwise exponential Bishop space}\index{$F_{|A}$} \index{$F \to G$}
$\mathcal{F} \rightarrow \mathcal{G} = (\Mor(\mathcal{F}, \mathcal{G}), 
F \rightarrow G)$ are defined, respectively, by 
\[ F \times G := \bigvee \left[\{f \circ \pr_{X}, \mid f \in F\} \cup \{g \circ \pr_{Y} \mid g \in G\}\right] 
=: \bigvee_{f \in F}^{g \in G}f \circ \pr_{X}, g \circ \pr_{Y},\]
\[ F_{|A} = \bigvee \{f_{|A} \mid f \in F\} =: \bigvee_{f \in F}f_{|A}\]
\begin{center}
\begin{tikzpicture}

\node (E) at (0,0) {$A$};
\node[right=of E] (F) {$X$};
\node [right=of F] (B) {$\Real$,};

\draw[right hook->] (E)--(F) node [midway,above] {$i_A$};
\draw[->] (F)--(B) node [midway,above] {$f$};
\draw[->,bend right] (E) to node [midway,below] {$f_{|A}$} (B);

\end{tikzpicture}
\end{center}
\[ F \rightarrow G := \bigvee \big\{\phi_{x, g} \mid x \in X, g \in G \big\} :=
\bigvee_{x \in X}^{g \in G}\phi_{x, g},\]
\[ \phi_{x, g}: \Mor(\mathcal{F}, \mathcal{G}) \rightarrow \mathbb{R}, \ \ \ \ \phi_{x, g}(h) = g(h(x));
\ \ \ \ x \in X, \ g \in G.\]
\end{definition}

One can show inductively the following $\bigvee$-liftings
\begin{align*}
\bigvee F_{0} \times \bigvee G_{0} & := \bigvee \left[\{f_{0} \circ \pr_{X}, \mid f_{0} \in F_{0}\} \cup
\{g_{0} \circ \pr_{Y} \mid g_{0} \in G_{0}\}\right]\\
& =: \bigvee_{f_{0} \in F_{0}}^{g_{0} \in G_{0}}f_{0} \circ \pr_{X}, g_{0} \circ \pr_{Y},
\end{align*}
\[ \big(\bigvee F_{0}\big)_{|A} = \bigvee \{{f_{0}}_{|A} \mid f_{0} \in F_{0}\} =: \bigvee_{f_{0}
\in F_{0}}{f_{0}}_{|A},\]
\[ F \rightarrow \bigvee G_0 = \bigvee \big\{\phi_{x, g_0} \mid x \in X, g_0 \in G_0 \big\} :=
\bigvee_{x \in X}^{g_0 \in G_0}\phi_{x, g_0}.\]
The relative topology $F_A$ is the least topology on $A$ that makes $i_A$ a Bishop morphism, and 
the product topology $F \times G$ is the least topology on $X \times Y$ 
that makes the projections $\pr_X$ and $\pr_Y$ Bishop morphisms. The term pointwise exponential Bishop
topology is due to the fact that $F \to G$ behaves like the the classical topology of the pointwise convergence 
on $C(X, Y)$, the set of continuous functions from the topological space $X$ to the topological space $Y$.

%
%
%
%
%
%
%
%
%

\section{Directed sets}
\label{sec: directed}

\begin{definition}\label{def: dirset}
Let $I$ be a set and $i \lt_I j$ a binary extensional relation\index{binary extensional relation} on $I$ i.e.,
\[\forall_{i, j, i{'}, j{'} \in I}\big(i =_I i{'} \ \& \ j =_I j{'} \ \& \ i \lt_I j \To i{'} \lt_I j{'}\big).\]
If $i \lt_I j$ is reflexive and transitive, then $(I, \lt_I)$ is called a 
a preorder\index{preorder}. We call a preorder $(I, \lt_I)$ a \textit{directed set}\index{directed set}, and 
\textit{inverse-directed}\index{inverse-directed set}, respectively, if
\[ \forall_{i, j \in I}\exists_{k \in I}\big(i \lt_I k \ \& \ j \lt_I k\big), \]
\[\forall_{i, j \in I}\exists_{k \in I}\big(i \mt_I k \ \& \ j \mt_I  k\big) .\]
The covariant \textit{covariant diagonal}\index{covariant diagonal} $D^{\lt}(I)$\index{$D^{\lt}(I)$} of $\lt_I$, the 
contravariant diagonal\index{contravariant diagonal} $D^{\mt}(I)$\index{$D^{\mt}(I)$} of $\lt_I$, and
the $\lt_I$-\textit{upper set} $I_{ij}^{\lt}$\index{$I_{ij}^{\lt}$}
of $i, j \in I$\index{$\lt_I$-upper set} are defined, respectively, by 
\[ D^{\lt}(I) := \big\{(i, j) \in I \times I \mid i \lt_I j \big\}, \]
\[ D^{\mt}(I) := \big\{(j, i) \in I \times I \mid j \mt_I i \big\}, \]
\[ I_{ij}^{\lt} := \{k \in I \mid  i \lt_I k \ \& \ j \lt_I k\}. \]

\end{definition}

Since $i \lt_I j$ is extensional, $D^{\lt}(I), D^{\mt}(I)$, and $I_{ij}^{\lt}$ are extensional 
subsets of $I \times I$.

\begin{definition}\label{def: modulus}
Let $(I, \lt_I)$ be a poset i.e., a preorder such that $\big[i \lt_I j \ \& \ j \lt_I i \big]\To i =_I j$, for every 
$i, j, \in I$. A \textit{modulus of directedness}\index{modulus of directedness} for $I$ is a 
function $\delta \colon I \times I \to I$, such
that for every $i, j, k \in I$ the following conditions are satisfied:\\[1mm]
$(\delta_1)$ $i \lt_I \delta(i, j)$ and $j \lt_I \delta(i, j)$.\\[1mm]
$(\delta_2)$ If $i \lt_I j$, then $\delta(i, j) =_I \delta(j, i) =_I j$.\\[1mm]
$(\delta_3)$ $\delta\big(\delta(i, j), k\big) =_I \delta\big(i, \delta(j, k)\big)$.
\end{definition}

In what follows we avoid for simplicity the use of subscripts on the relation symbols.
If $(I, \lt)$ is a preordered set and $(J, e) \subseteq I$, where $e : J \hookrightarrow I$, 
and using for simplicity the same symbol $\lt$, if we define
$j \lt j{'} : \TOT e(j) \lt e(j{'}),$
for every $j, j{'} \in J$, then $(J, \lt)$ is only a preordered set. If $J$ is a cofinal subset of $I$, 
which classically it is defined by the condition
$\forall_{i \in I}\exists_{j \in J}\big(i \lt j \big)$, then $(J, \lt)$ becomes a
directed set. To avoid the use of dependent choice, we add in the definition of a cofinal subset $J$
of $I$ a modulus of cofinality for $J$.

\begin{definition}\label{def: cofinal}
Let $(I, \lt)$ be a directed set and $(J, e) \subseteq I$, and let 
$j \lt j{'} : \TOT e(j) \lt e(j{'})$, for every $j, j{'} \in J$. We say that $J$ is 
\textit{cofinal in} $I$, if there is a function $\cof_J : I \to J$, which we call a \textit{modulus of
cofinality} of $J$ in $I$, that satisfies the following conditions:\\[1mm]
$(\Cf_1)$ $\forall_{j \in J}\big(\cof_J(e(j)) =_J j\big)$.
\begin{center}
\begin{tikzpicture}

\node (E) at (0,0) {$J$};
\node[right=of E] (F) {$I$};
\node[right=of F] (A) {$J$.};

\draw[right hook->] (E)--(F) node [midway,below] {$e$};
\draw[->] (F)--(A) node [midway,below] {$\cof_J$};
\draw[->,bend left] (E) to node [midway,above] {$\id_J$} (A) ;

\end{tikzpicture}
\end{center}
$(\Cf_2)$  $\forall_{i, i{'} \in I}\big(i \lt i{'} \To \cof_J(i) \lt \cof_J(i{'})\big)$.\\[1mm] 
$(\Cf_3)$ $\forall_{i \in I}\big(i \lt e(\cof_J (i))\big)$.\\[1mm]
We denote the fact that $J$ is cofinal in $I$ by $(J, e, \cof_J) \subseteq^{\cof} I$, or, simpler, by\index{$(J, e, \cof_J) 
\subseteq^{\cof} I$}\index{$J \subseteq^{\cof} I$}$J \subseteq^{\cof} I$. 
\end{definition}

Taking into account the embedding $e$ of $J$ into $I$, the condition (iii) is the exact writing of
the classical defining condition $\forall_{i \in I}\exists_{j \in J}\big(i \lt j \big)$. 
To add the condition (i) is harmless, since $\lt$ is reflexive. If we consider the condition (iii) on $e(j)$,
for some $j \in J$, then by the condition (i) we get the transitivity $e(j) \lt e(\cof_J (e(j)))  = e(j)$.
The condition (ii) is also harmless to add. In the classical setting if $i \lt i{'}$, and $j, j{'} \in J$
such that $i \lt j$ and $i{'} \lt j{'}$, then there is some $i{''} \in I$ such that $j{'} \lt i{''}$ and
$j \lt i{''}$. If $i{''} \lt j{''}$, for some $j{''} \in J$, 
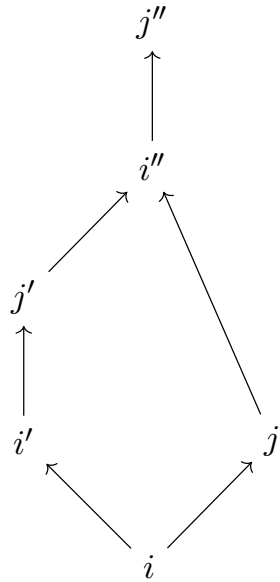
\begin{center}
\resizebox{4cm}{!}{%
\begin{tikzpicture}

\node (E) at (0,0) {$i$};
\node[above=of E] (F) {};
\node[right=of F] (A) {$j$};
\node[left=of F] (B) {$i{'}$};
\node[above=of B] (C) {$j{'}$};
\node[above=of C] (D) {};
\node[right=of D] (H) {$i{''}$};
\node[above=of H] (G) {$j{''}$};

\draw[->] (E)--(B) node [midway,below] {};
\draw[->] (E)--(A) node [midway,left] {};
\draw[->] (B)--(C) node [midway,below] {};
\draw[->] (C)--(H) node [midway,right] {};
\draw[->] (A)--(H) node [midway,right] {};
\draw[->] (H)--(G) node [midway,right] {};

\end{tikzpicture}
}
\end{center}
then $j \lt j{''}$. Since $i{'} \lt j{''}$ too,
the condition (ii) is justified. The added conditions (i) and (ii) are used in the proofs of 
Theorem~\ref{thm: cofinal2} and Lemma~\ref{lem: cofinallemma}(ii), respectively. Moreover, they are used
in the proof of Theorem~\ref{thm: cofinal3}.
The extensionality of $\lt$
is also used in the proofs of Theorem~\ref{thm: cofinal2} and Theorem~\ref{thm: cofinal3}.

E.g., if $\Even$ and $\Odd$ denote the sets of even and odd natural numbers, respectively, 
let $e \colon \Even \eto \Nat$, defined by the identity map-rule, and $\cof_{\Even} \colon \Nat \to 2 \Nat$, 
defined by the rule
\[ \cof_{2 \Nat}(n) := \left\{ \begin{array}{ll}
                 n    &\mbox{, $n \in \Even$}\\
                 n+1           &\mbox{, $n \in \Odd$.}
                 \end{array}
          \right.\]           
Then $(\Even, e, \cof_{\Even}) \subseteq \Nat$.

\begin{remark}\label{rem: cofinal1}
If $(I, \lt)$ is a directed set and $(J, e, \cof_J) \subseteq^{\cof} I$, then $(J, \lt)$ is directed.
\end{remark}

\begin{proof}
 Let $j, j{'} \in J$ and let $i \in I$ such that $e(j) \lt i$ and $e(j{'}) \lt i$. Since $i \lt e(\cof_J (i))$, 
 we get $e(j) \lt e(\cof_J (i))$ and $e(j{'}) \lt e(\cof_J (i))$ i.e., $j \lt \cof_J(i)$ and $j{'} \lt \cof_J(i)$.
\end{proof}


\begin{thebibliography}{50}
\bibitem{AR10} P.~Aczel, M. Rathjen: \textit{Constructive Set Theory}, book draft, 2010.
\bibitem{Aw95} S.~Awodey: Axiom of Choice and Excluded Middle in Categorical Logic, Bulletin of Symbolic Logic,
volume 1, 1995, 344.
\bibitem{Aw10} S.~Awodey: \textit{Category Theory}, Oxford University Press, 2010.
\bibitem{BS18} A.~Bauer, A.~Swan: Every metric space is separable in function realizability, arXiv:1804.\\
00427, 2018.
\bibitem{Be81} M.~J.~Beeson: Formalizing constructive mathematics: why and how, in~\cite{Ri81}, 1981, 146--190.
\bibitem{Be82} M.~J.~Beeson: Problematic principles in constructive mathematics, in D.~van Dalen, D.~Lascar, 
J.~Smiley (Eds.) \textit{Logic Colloquium '80}, North-Holland, 1982, 11--55.
\bibitem{Be85} M.~J.~Beeson: \textit{Foundations of Constructive Mathematics}, Ergebnisse der
Mathematik und ihrer Grenzgebiete, Springer Verlag, 1985.
\bibitem{BAV12} S.~Bhat, A.~Agarwal, R.~Vuduc: A type theory for probability density
functions, POPL'12, ACM, 2012, 545--556.
\bibitem{Bi67} E.~Bishop: \textit{Foundations of Constructive Analysis}, McGraw-Hill, 1967.
\bibitem{Bi68a} E.~Bishop: A General Language, unpublished manuscript, 1968(9)?
\bibitem{Bi68b} E.~Bishop: How to Compile Mathematics into Algol, unpublished manuscript, 1968(9)?
\bibitem{Bi70} E.~Bishop: Mathematics as a Numerical Language, in~\cite{KMV70}, 1970, 53--71. 
\bibitem{Bi71} E.~Bishop: The Neat Category of Stratified Spaces, unpublished manuscript, University of 
California, San Diego, 1971.
\bibitem{Bi72} E.~Bishop: Aspects of Constructivism, Notes on the lectures delivered at the Tenth Holiday 
Mathematics Symposium held at New Mexico State University, Las Cruses, during the period December 27-31, 1972.
\bibitem{Bi73} E.~Bishop: Schizophrenia in Contemporary Mathematics, American Mathematical
Society Colloquium Lectures, Missoula University of Montana 1973, and in~\cite{Ro85}. 
\bibitem{Bi75} E.~Bishop: The crisis in contemporary mathematics, Historia Mathematics 2, 1975, 507--517.
\bibitem{Bi86} \textit{Selected Papers, Errett Bishop}, edited by J.~Wermer, World Scientific Publishing,
Singapore and Philadelphia, 1986. 
\bibitem{BC72} E.~Bishop and H.~Cheng: \textit{Constructive Measure Theory}, Mem. Amer. Math. Soc. 116, 1972.
\bibitem{BB85} E.~Bishop and D.~S.~Bridges: \textit{Constructive Analysis}, Grundlehren der math.
Wissenschaften 279, Springer-Verlag, Heidelberg-Berlin-New York, 1985.
\bibitem{BR87} D.~S.~Bridges and F.~Richman: \textit{Varieties of Constructive Mathematics}, 
Cambridge University Press, 1987.
\bibitem{BD91} D.~Bridges, O.~Demuth: On the Lebesgue measurability of continuous functions in 
constructive analysis, Bulletin of the American Mathematical Society, Vol. 24, No. 2, 1991, 259--276.
\bibitem{BR99} D.~S.~Bridges and S.~Reeves: Constructive Mathematics in Theory and Programming Practice, 
Philosophia Mathematica (3), 1999, 65--104.
\bibitem{Br99} D.~S.~Bridges: Constructive mathematics; a foundation for computable analysis, 
Theoretical Computer Science 219, 1999, 95--109.
\bibitem{BV06} D.~S.~Bridges and L.~S.~V\^{\i}\c{t}\u{a}: \textit{Techniques of Constructive Analysis},
in: Universitext, Springer, New York, 2006.
\bibitem{BV11} D.~S.~Bridges, L.~S.~V\^{\i}\c{t}\u{a}: \textit{Apartness and Uniformity: A Constructive 
Development}, in: CiE series ``Theory and Applications of Computability'', Springer Verlag, Berlin Heidelberg, 2011.
\bibitem{Br12} D.~S.~Bridges: Reflections on function spaces, Annals of Pure and Applied Logic 163, 2012, 101--110.
\bibitem{Ca82} G.~Cantor: \"Uber unendliche, lineare Punktmannichfaltigkeiten, Nummer 3. Mathematische
Annalen, 20, 1882, 113--121.
\bibitem{BIRS21} D.~S.~Bridges, H.~Ishihara, M.~Rathjen, H.~Schwichtenberg (Eds.): \textit{Handbook of Bishop Constructive
Mathematics}, Cambridge University Press, to appear, 2021.
\bibitem{Ca56} C.~Carath\'eodory: \textit{Mass und Integral und ihre algebraisierung}, Springer Basel AG, 1956.
\bibitem{Ch72} Y.~K.~Chan: A constructive approach to the theory of stochastic processes, Transactions
of the American Mathematical Society, Vol. 65, 1972, 37--44.
\bibitem{Ch72b} Y.~K.~Chan: A constructive study of measure theory, Pacific Journal of Mathematics, Vol. 41, No. 1,
1972, 63--79.
\bibitem{Ch72c} Y.~K.~Chan: A constructive approach to the theory of stochastic processes, 
Transactions of the American Mathematical Society, Vol. 165, 1972, 37--44.
\bibitem{Ch74} Y.~K.~Chan: Notes on constructive probability theory, The Annals of Probability, Vol. 2, 
No. 1, 1974, 51--75.
\bibitem{Ch75} Y.~K.~Chan: A short proof of an existence theorem in constructive measure theory,
Proceedings of the American Mathematical Society, Vol. 48, No. 2, 1975, 435--436.
\bibitem{Ch19} Y. K. Chan: \textit{Foundations of Constructive Probability Theory}, arXiv:1906.01803v2, 2019.
\bibitem{Co86} R.~L.~Constable et al. \textit{Implementing Mathematics with the Nuprl Proof Development
373 System}, Prentice-Hall, Inc., Upper Saddle River, NJ, USA, 1986.
\bibitem{CP98} T.~Coquand, H.~Persson: Integrated Development of Algebra in Type Theory, preprint, 1998.
\bibitem{CP02} T.~Coquand, E.~Palmgren: Metric Boolean algebras and constructive measure theory, Arch. Math. Logic
41, 2002, 687--704.
\bibitem{CDPS05} T.~Coquand, P.~Dybjer, E.~Palmgren, A.~Setzer: \textit{Type-theoretic Foundations of 
Constructive Mathematics}, book-draft, 2005.
\bibitem{CS07} T.~Coquand, A.~Spiwack: Towards Constructive Homological Algebra in Type Theory, in LNCS 4573, 
2007, 40--54.
\bibitem{CS09} T.~Coquand, B.~Spitters: Integrals and valuations, Journal of Logic $\&$ Analysis, 1:3, 2009, 1--22.
\bibitem{Co14} T.~Coquand: A remark on singleton types, manuscript, 2014, available at
http://www.cse.chalmers.se/$\sim$coquand/singl.pdf, 2014.
\bibitem{Da18} P.~J.~Daniell: A general form of integral, Annals of mathematics, Second Series, 19 (4),
1918, 279--294.
\bibitem{De88} R.~Dedekind: \textit{Was sind und was sollen die Zahlen?}, 1. Auflage, Vieweg, Braunschweig, 1888.
\bibitem{Du66} J.~Dugundji: \textit{Topology}, Allyn and Bacon, 1966.
\bibitem{Ed09} A.~Edalat: A computable approach to measure and integration theory, Information and Computation
207, 2009, 642--659.
\bibitem{FS16} F.~Faissole, B.~Spitters: Synthetic topology in Homotopy Type Theory for probabilistic 
programming, preprint, 2016.
\bibitem{Fe75} S.~Feferman: A language and axioms for explicit mathematics, in
J.~N.~Crossley (Ed.) \textit{Algebra and Logic}, Springer Lecture Notes 450, 1975, 87--139.
\bibitem{Fe79} S.~Feferman: Constructive theories of functions and classes, in Boffa et al. (Eds.) 
\textit{Logic Colloquium 78}, North-Holland, 1979, 159--224.
\bibitem{Fr73} H.~Friedman: The consistency of classical set theory relative to a set theory with
intuitionistic logic, J. Symbolic Logic, 38, 1973, 315--319.
\bibitem{Fr77} H.~Friedman: Set theoretic foundations for constructive analysis, Annals of Math. 105,
1977, 1--28.
\bibitem{GJ60} L.~Gillman, M.~Jerison: \textit{Rings of Continuous Functions}, Van Nostrand, 1960.
\bibitem{Go73} R.~Godement: \textit{Topologie alg\'ebrique et th\'eorie de faisceaux}, Hermann, Paris, 1973.
\bibitem{Go84} R.~Goldblatt: \textit{Topoi, The Categorical Analysis of Logic}, North-Holland, 1984.
\bibitem{Gr81} N.~Greenleaf: Liberal constructive set theory, in~\cite{Ri81}, 1981, 213--240.
\bibitem{GV87} N.~D.~Goodman, R.~E.~Vesley: Obituary: John R.~Myhill (1923-1987), History and Philosophy of Logic,
8, 1987, 243--244.
\bibitem{Ha74} P.~R.~Halmos: \textit{Measure theory}, Springer-Verlag New York Inc.~1974.
\bibitem{Ha14} F.~Hausdorff: Grundz\"uge der Mengenlehre, Veit $\&$ Comp., Leipzig 1914 
(reproduced in Srishti D. Chatterji et al.~(Eds.): Felix Hausdorff. Gesammelte Werke. Band II: 
Grundz\"uge der Mengenlehre, Springer, Berlin 2002).
\bibitem{He56} A.~Heyting: \textit{Intuitionism, An introduction}, North-Holland, 1956.
\bibitem{He50} E.~Hewitt: Linear functionals on spaces of continuous functions, Fund. Math. 37, 1950, 161--189.
\bibitem{HS98} G.~Huet, A.~Sa\"ibi: Constructive Category Theory, in Proceedings of the joint CLICS-TYPES workshop on categories and type theory, Goteborg, MIT Press, 1998.
\bibitem{IP06} H.~Ishihara, E.~Palmgren: Quotient topologies in constructive set theory and type theory, Annals of 
Pure and Applied Logic 141, 2006, 257--265.
\bibitem{Is13} H.~Ishihara: Relating Bishop's function spaces to neighborhood spaces, Annals of Pure and 
Applied Logic 164, 2013, 482--490.
\bibitem{Ja19} L.~Jaun: \textit{Category Theory in Explicit Mathematics}, PhD Thesis, University of Bern, 2019.
\bibitem{Jo02} P.~T.~Johnstone: \textit{Sketches of an Elefant A Topos Theory Compendium}, Volume 1, 
Oxford University Press, 2002.
\bibitem{KMV70} A.~Kino, J.~Myhill, R.~E.~Vesley (Eds.): \textit{Intuitionism and Proof Theory}, 
North-Holland, 1970.
\bibitem{Ko48} A.~N.~Kolmogoroff: Alg\`ebres de Boole m\'etrique compl\`etes, VI Zjazd Mathematyk\'ow Polskich, 
Warwaw, 1948, 21--30.
\bibitem{LQ15} H. Lombardi, C. Quitt\'e: \textit{Commutative Algebra: Constructive Methods}, Springer, 2015.
\bibitem{LR03}  R.~S.~Lubarsky, M.~Rathjen: On the regular extension axiom and its variants, 
Mathematical Logic Quarterly 49 (5), 511--518, 2003.
\bibitem{Ma09} M.~E.~Maietti: A minimalist two-level foundation for constructive mathematics, 
Annals of Pure and Applied Logic 160(3), 2009, 319--354. 
\bibitem{ML68} P.~Martin-L\"{o}f: \textit{Notes on Constructive Mathematics}, Almqvist and Wiksell, 1970.
\bibitem{ML75} P.~Martin-L\"{o}f: An intuitionistic theory of types: predicative part, in H. E. Rose and
J.~C.~Shepherdson (Eds.) \textit{Logic Colloquium'73, Proceedings of the Logic Colloquium}, volume 80 of
Studies in Logic
and the Foundations of Mathematics,North-Holland, 1975, 73--118.
\bibitem{ML84} P.~Martin-L\"{o}f: \textit{Intuitionistic type theory: Notes by Giovanni Sambin on a series
of lectures given in Padua, June 1980}, Napoli: Bibliopolis, 1984.
\bibitem{ML98} P.~Martin-L\"{o}f: An intuitionistic theory of types, in~\cite{SS98}, 127--172. 
\bibitem{Mc06} C.~McLarty: Two Constructicist Aspects of Category Theory, Philosophia Scienti$\oe$, Cahier sp\'ecial 6, 2006, 95--114.
\bibitem{MRR88} R.~Mines, F.~Richman, W.~Ruitenburg: \textit{A course in constructive algebra}, Springer
Science$+$Business Media New York, 1988.
\bibitem{My68} J.~Myhill: Formal systems of intuitionistic analysis I, in van Rootselaar and Stall (Eds.)
\textit{Logic, methodology and philsosphy of science III}, North-Holland, Amsterdam, 1968, 161--178.
\bibitem{My70} J.~Myhill: Formal systems of intuitionistic analysis II, in~\cite{KMV70}, 1970, 151--162.
\bibitem{My73} J.~Myhill: Some properties of intuitionistic Zermelo-Fraenkel set theory, in A. Matthias , 
H.~Rogers (Eds.) \textit{Cambridge Summer School in Mathematical Logic}, LNM 337, Springer, 1972, 206--231.
\bibitem{My75} J.~Myhill: Constructive Set Theory, J.~Symbolic Logic 40, 1975, 347-382.
\bibitem{Pa05} E.~Palmgren: Bishop's set theory, Slides from TYPES Summer School 2005, Gothenburg, 
in http://staff.math.su.se/palmgren/, 2005.
\bibitem{Pa12a} E.~Palmgren: Proof-relevance of families of setoids and identity in type theory, Arch. Math. Logic, 51,
2012, 35-47. 
\bibitem{Pa12b} E.~Palmgren: Constructivist and structuralist foundations: Bishop's and Lawvere's theories
of sets, Annals of Pure and Applied Logic 163, 2012, 1384--1399.
\bibitem{Pa13} E.~Palmgren: Bishop-style constructive mathematics in type theory - A tutorial, Slides, in 
http://staff.math.su.se/palmgren/, 2013. 
\bibitem{Pa14} E.~Palmgren: \textit{Lecture Notes on Type Theory}, 2014.
\bibitem{Pa17} E.~Palmgren: On Equality of Objects in Categories in Constructive Type Theory, TYPES 2017, A.~Abel 
et.~al.~(Eds.), Article No.~7; pp.~7:1-7:7.
\bibitem{PW14} E.~Palmgren, O.~Wilander: Constructing categories and setoids of setoids in type theory, 
Logical Methods in Computer Science. 10 (2014), Issue 3, paper 25. 
\bibitem{Pe15} I.~Petrakis: \textit{Constructive Topology of Bishop Spaces}, PhD Thesis,
Ludwig-Maximilians-Universit\"{a}t, M\"{u}nchen, 2015. 
\bibitem{Pe15a} I. Petrakis: Completely Regular Bishop Spaces, in A. Beckmann, V. Mitrana and M. Soskova (Eds.): 
\textit{Evolving Computability}, CiE 2015, LNCS 9136, 2015, 302--312.
\bibitem{Pe16a} I.~Petrakis: The Urysohn Extension Theorem for Bishop Spaces, in 
S. Artemov and A. Nerode (Eds.) \textit{Symposium on Logical Foundations of Computer Science 2016},
LNCS 9537, Springer, 2016, 299--316.
\bibitem{Pe16b} I.~Petrakis: A constructive function-theoretic approach to topological compactness, 
Proceedings of the 31st Annual ACM-IEEEE Symposium on Logic in Computer Science (LICS 2016), 
July 5-8, 2016, NYC, USA, 605--614.
\bibitem{Pe19a} I.~Petrakis: Borel and Baire sets in Bishop Spaces, in F. Manea et. al. (Eds): 
Computing with Foresight and Industry, CiE 2019, LNCS 11558, Springer, 2019, 240--252.
\bibitem{Pe19b} I.~Petrakis: Constructive uniformities of pseudometrics and Bishop topologies,
Journal of Logic and Analysis 11:FT2, 2019, 1--44.
\bibitem{Pe19d} I.~Petrakis: Direct spectra of Bishop spaces and their limits, Logical Methods in Computer Science,
Volume 17, Issue 2, 2021, pp. 4:1-4:50. 
\bibitem{Pe19c} I.~Petrakis: Dependent Sums and Dependent Products in Bishop's Set Theory,
in P. Dybjer et. al. (Eds) TYPES 2018, LIPIcs, Vol. 130, Article No. 3, 2019.
\bibitem{Pe19d} I.~Petrakis: A Yoneda lemma-formulation of the univalence axiom, manuscript, 2019, available at 
http://www.mathematik.uni-muenchen.de/$\sim$petrakis/content//Preprints.php.
\bibitem{Pe20a} I.~Petrakis: McShane-Whitney extensions in constructive analysis, Logical Methods in Computer
Science, Volume 16, Issue  2020, 18:1--18:23
\bibitem{Pe20b} I.~Petrakis: Embeddings of Bishop spaces, Journal of Logic and Computation, 
 exaa015, 2020, https://doi.org/10.1093/logcom/exaa015.
\bibitem{Pe20c} I. Petrakis: Functions of Baire-class one over a Bishop topology, in M. Anselmo et.~al.~(Eds.)
\textit{Beyond the Horizon of Computability} CiE 2020, Springer, LNCS 12098, 2020, 215-227.
\bibitem{Pe20e} I.~Petrakis: Proof-relevance in Bishop-style constructive mathematics, submitted, 2020.
\bibitem{Pe20d} I.~Petrakis: Bases of pseudocompact Bishop spaces, invited chapter in~\cite{BIRS21}, to appear,
2021.
\bibitem{Pe21a} I.~Petrakis: Closed subsets in Bishop topological groups, submitted, 2021.
\bibitem{Pe21c} I.~Petrakis: Chu representations of categories related to constructive mathematics,
arXiv:2106.01878v1, 2021. 
\bibitem{PZ20} I.~Petrakis, M.~Zeuner: Pre-measure spaces and pre-integration spaces in predicative Bishop-Cheng
measure theory, in preparation, 2021.
\bibitem{Ra06} M.~Rathjen: Choice principles in constructive and classical set theories, Lecture 
Notes in Logic 27, 2006, 299--326.
\bibitem{Ra09} M.~Rathjen: The constructive Hilbert program and the limits of Martin-L\"{o}f Type Theory,
in S.~Lindstr\"om et al. (Eds.) \textit{Logicism, Intuitionism, and Formalism: What has become of them?}
Synthese Library Volume 341, Springer, 2009, 397--433.
\bibitem{Ri81} F.~Richman: \textit{Constructive Mathematics}, LNM 873, Springer-Verlag, 1981.
\bibitem{Ri01} F.~Richman: Constructive mathematics without choice, in~\cite{Sc01}, pp.199--205.
\bibitem{Ri12} E. Rijke: \textit{Homotopy Type Theory}, Master Thesis, Utrecht University, 2012.
\bibitem{Ro85} M.~Rosenblatt (Ed.): \textit{Errett Bishop:
Reflections on Him and His Research}, Contemporary Mathematics Volume 39, American Mathematical Society, 1985.
\bibitem{RS65} K.~A.~Ross, K.~Stromberg: Baire sets and Baire measures, Arkiv f\"or Matematik Band 6 nr 8, 1965, 151--160.
\bibitem{SS98} G.~Sambin, J.~M.~Smith (Eds.): \textit{Twenty-five years of constructive type theory},
Oxford University Press, 1998.
\bibitem{Sa19} G.~Sambin: The Basic Picture: Structures for Constructive Topology, Oxford University Press, 2020.
\bibitem{Sa68} N.~\v Sanin: \textit{Constructive Real Numbers and Function Spaces}, 
Translations of Mathematical Monographs, vol. 21, AMS, Providence Rhode Island, 1968.
\bibitem{Sh18} M.~Shulman: Linear Logic for Constructive Mathematics, arXiv:1805.07518v1, 2018.
\bibitem{Sc01} P.~Schuster, U.~Berger, H.~Osswald (Eds.): \textit{Reuniting the Antipodes Constructive 
and Nonstandard Views of the Continuum}, Proc. 1999 Venice Symposium, Dordrecht: Kluwer, 2001.
\bibitem{Sc04} P.~Schuster: Countable Choice as a Questionable Uniformity Principle, Philosophia Mathematica (3) Vol.~12, pp.~106--134, 2004.
\bibitem{SW12} H.~Schwichtenberg, S.~Wainer: \textit{Proofs and Computations},  Perspectives in Logic, 
Association for Symbolic Logic and Cambridge University Press, 2012. 
\bibitem{Sc19} H.~Schwichtenberg: \textit{Constructive analysis with witnesses}, Lecture notes, LMU, 2019. 
\bibitem{Se54} I.~Segal: Abstract probability spaces and a theorem of Kolmogoroff, Amer. J. Math. 76, 1954,
721--732.
\bibitem{Se65} I.~Segal: Algebraic integration theory, Bull. Amer. Math. Soc. 71, 1965, 419--489.
\bibitem{Sp02} B.~Spitters: Constructive and intuitionistic integration theory and functional analysis.
PhD Thesis, University of Nijmegen, 2002.
\bibitem{Sp06} B.~Spitters: Constructive algebraic integration theory. Ann. Pure Appl. Logic 137(1-3), 2006,
380--390.
\bibitem{St18} T.~Streicher: \textit{Realizability}, Lecture Notes, TU Darmstadt, 2018.
\bibitem{Ta73} S.~J.~Taylor: \textit{Introduction to Measure and Integration}, Cambridge University Press, 1973.
\bibitem{HoTT13} The Univalent Foundations Program: \textit{Homotopy Type Theory:
Univalent Foundations of Mathematics}, Institute for Advanced Study, Princeton, 2013.
\bibitem{TD88I} A.~S.~Troelstra and D.~van Dalen: \textit{Constructivism in Mathematics}, Volume I, North-Holland, 1988.
\bibitem{TD88II} A.~S.~Troelstra and D.~van Dalen: \textit{Constructivism in Mathematics}, Volume II, North-Holland, 1988.
\bibitem{We40} A.~Weil: Calcul des probabilit\'es, m\'ethode axiomatique, int\'egration, Revue Sci. 
(Rev. Rose Illus.) 78, 1940, 201--208.
\bibitem{Wi13} T.~Wiklund: Locally cartesian closed categories, coalgebras, and containers,
Uppsala Universitet, U.U.D.M Project Report, 2013:5.
\bibitem{Ze19} M.~Zeuner: \textit{Families of Sets in Constructive Measure Theory}, Master's Thesis, LMU, 2019.
\end{thebibliography}
\end{document}